\documentclass[letterpaper,reqno,10pt,twoside]{amsart}
\usepackage{amsmath,amsthm,amsfonts,amssymb,euscript,mathrsfs,graphics,color,latexsym,marginnote}
\usepackage{graphicx}
\usepackage[left=.75 in, right=.75 in,top=.75 in, bottom=.75 in]{geometry}
\usepackage{hyperref}
\usepackage[toc,page]{appendix}
\usepackage{relsize}
\usepackage[shortlabels]{enumitem}
\usepackage[all]{xy}
\usepackage{float}
\usepackage{tikz}
\usepackage{xargs}
\usepackage[dvipsnames]{xcolor}
\usepackage{hyperref}
\usepackage[colorinlistoftodos,prependcaption,textsize=small]{todonotes}

\AtBeginDocument{\setlength{\abovedisplayskip}{3pt plus 2pt minus 1pt}\setlength{\belowdisplayskip}{\abovedisplayskip}\setlength{\abovedisplayshortskip}{0pt plus 1.5pt}\setlength{\belowdisplayshortskip}{1.5pt plus 1.5pt minus 0.5pt}}

\hypersetup{colorlinks=true, pdfstartview=FitV,linkcolor=blue!70!black,citecolor=red!70!black, urlcolor=green!60!black}
\definecolor{labelkey}{rgb}{0.6,0,0}
\newcommandx{\change}[2][1=]{\todo[#1]{#2}}
\newcommandx{\unsure}[2][1=]{\todo[linecolor=red,backgroundcolor=red!25,bordercolor=red,#1]{#2}}
\newcommandx{\rmk}[2][1=]{\todo[linecolor=blue,backgroundcolor=blue!25,bordercolor=blue,#1]{#2}}
\newcommandx{\info}[2][1=]{\todo[linecolor=OliveGreen,backgroundcolor=OliveGreen!25,bordercolor=OliveGreen,#1]{#2}}
\newcommandx{\improvement}[2][1=]{\todo[linecolor=Plum,backgroundcolor=Plum!25,bordercolor=Plum,#1]{#2}}
\newcommandx{\thiswillnotshow}[2][1=]{\todo[disable,#1]{#2}}
\allowdisplaybreaks
\makeatletter
\renewcommand \theequation {\ifnum \c@section>\z@ \@arabic\c@section.\fi\@arabic\c@equation} \@addtoreset{equation}{section}
\@namedef{subjclassname@2020}{2020 Mathematics Subject Classification}
\makeatother
\newtheorem{theorem}{Theorem}[section]
\newtheorem{corollary}{Corollary}[theorem]
\newtheorem{lemma}[theorem]{Lemma}
\newtheorem{proposition}[theorem]{Proposition}

\theoremstyle{definition}

\theoremstyle{remark}
\newtheorem{remark}{Remark}[section]

\def\XXint#1#2#3{{\setbox0=\hbox{$#1{#2#3}{\int}$ }
\vcenter{\hbox{$#2#3$ }}\kern-.6\wd0}}

\def\p{\partial}

\def\no{\nonumber}

\let\caronaccent\v \def\v{\mathcal{V}}

\def\a{\mathcal{A}}
\def\ep{\varepsilon}

\def\ep{\varrho}

\newcommand{\cA}{\mathcal{A}}

\newcommand{\cE}{\mathcal{E}}
\newcommand{\cD}{\mathcal{D}}

\newcommand{\cU}{\mathcal{U}}

\newcommand{\bbD}{\mathbb{D}}

\newcommand{\nat}{\natural}
\newcommand{\Ft}{\widetilde F}
\newcommand{\fM}{\mathfrak{M}}

\newcommand{\sE}{\mathscr{E}}
\newcommand{\chis}{\chi_{\mathrm{bet}+}}

\newcommand{\Lt}{L^2_\theta}

\theoremstyle{plain}
\newtheorem*{maintheorem}{Main Theorem}

\providecommand{\func}[1]{\operatorname{#1}}

\begin{document}
\title[Asymptotic stability of capillary surfaces in 3D]{Asymptotic stability of capillary surfaces in three dimensions}
\author{Yan Guo}
\author{Ian Tice}
\author{Xiaoding Yang}

\subjclass[2020]{Primary: 35Q30, 35R35, 76D45; Secondary: 35B40, 76E17, 47A60}

\keywords{Contact lines, Navier-Stokes equations, free boundary problems}

\begin{abstract}
Contact line occurs wherever a free surface meets a fixed solid boundary. Steady capillary surfaces were studied by Gauss, Laplace and Young, and are characterized by the Young--Laplace equation together with the Young relation, which prescribes the contact angle $\theta_{\mathrm{eq}}$, measured through the fluid, at which the surface meets the container:  $\cos\theta_{\mathrm{eq}}=\frac{\gamma_{sv}-\gamma_{sf}}{\sigma}$. We establish the asymptotic stability of such capillary surfaces for incompressible Navier--Stokes flow in a three-dimensional vessel, with the Navier slip condition on the wall and the dynamic contact-line law
\begin{equation*}
-\kappa\partial_{t}\eta-\hat{\mathcal{W}}(\partial_{t}\eta)=\sigma\Big[\tfrac{\partial_{\nu}(\zeta+\eta)}{\sqrt{1+|\nabla(\zeta+\eta)|^{2}}}-\tfrac{\partial_{\nu}\zeta}{\sqrt{1+|\nabla\zeta|^{2}}}\Big].
\end{equation*}
The basic difficulty is the low regularity at the contact line, which in three dimensions is a curve of triple junctions. To overcome it, in addition to the global temporal vector field $\partial_{t}$ of \cite{GT}, we introduce the localized tangential vector field $\partial_{\theta}$ along the contact line, together with slicing Sobolev norms in the transverse variables and a decomposition of the velocity into its component along the contact line and its transverse part. We develop energy estimates, the gain of $3/2$ derivatives for the surface, and elliptic estimates for these slicing norms of $\partial_{\theta}$, which control the nonlinearity within the energy framework built in \cite{GT} for the temporal field $\partial_{t}$. At the highest temporal level, the surface satisfies $\partial_{t}^{m+1}\eta\in H^{3/2-\alpha}$, as in \cite{GT}, with a loss of $\alpha$ and with no estimate available for $\partial_{t}^{m+1}p$. At the highest angular level, by contrast, we obtain $\partial_{\theta}^{m+1}\eta\in H^{3/2}$ and $\partial_{\theta}^{m+1}p\in L^{2}$, with no loss. This subtle but critical gain is what closes the nonlinearity.
\end{abstract}

\maketitle
\tableofcontents

\setcounter{tocdepth}{2}

\section{Introduction}\label{sec:intro}

\subsection{The vessel and the equations of motion}

We study the motion of a viscous incompressible fluid in a three-dimensional
open-top vessel. We prove an a~priori exponential decay estimate for
sufficiently regular solutions that stay close to a capillary equilibrium.
Write $x'=(x,y)\in\mathbb R^2$ for the horizontal
coordinates and $z$ for height. The vessel is a bounded connected open
set $\mathcal V\subset\mathbb R^3$. Its upper part is a vertical channel
over a bounded domain $D\subset\mathbb R^2$ with smooth connected
boundary. We use the level $h_0$ of Subsection~\ref{subsec:cutoff_extension}
for the bottom of the channel. The analysis uses the wall up to a
neighbourhood of the contact line. We represent its dry continuation by
ending the channel at a horizontal rim of height $h_{\mathrm{rim}}>h_0$.
Thus
\begin{equation}\label{eq:intro-vessel}
\begin{aligned}
  \mathcal V_{\mathrm{top}}
    &:=\mathcal V\cap\{z\geq h_0\}
      =D\times[h_0,h_{\mathrm{rim}}),\\
  \mathcal V_{\mathrm{btm}}
    &:=\mathcal V\cap\{z\leq h_0\},\qquad
  \mathcal V=\mathcal V_{\mathrm{btm}}\cup\mathcal V_{\mathrm{top}}.
\end{aligned}
\end{equation}
The two parts meet along the interior cross-section $D\times\{h_0\}$.
The solid side wall in the upper channel is
$\partial D\times(h_0,h_{\mathrm{rim}})$, and the top is open over $D$.
The lower part of the vessel may have any shape for which its solid wall,
including the bottom, is smooth and joins the vertical side wall smoothly.

The lower part $\mathcal V_{\mathrm{btm}}$ is filled with fluid. At
equilibrium the free surface is the graph of a height $\zeta$ over $D$.
At time $t$ its height is $\zeta+\eta(\cdot,t)$. We assume
$h_0<\min_{\overline D}\zeta$ and
$\max_{\partial D}\zeta<h_{\mathrm{rim}}$.
For the small perturbations considered here the free surface remains
above $h_0$ and its boundary remains below the rim. The height condition
at the rim concerns the trace on $\partial D$. The equilibrium fluid
domain and the moving fluid domain are, respectively,
\begin{equation}\label{eq:intro-fluid-domains}
\begin{aligned}
  \Omega
    &:=\mathcal V_{\mathrm{btm}}\cup
       \{(x',z):x'\in D,\ h_0<z<\zeta(x')\},\\
  \Omega(t)
    &:=\mathcal V_{\mathrm{btm}}\cup
       \{(x',z):x'\in D,\ h_0<z<\zeta(x')+\eta(x',t)\}.
\end{aligned}
\end{equation}
As in \eqref{eq:intro-vessel}, the slice at $z=h_0$ belongs to
$\mathcal V_{\mathrm{btm}}$. The free surfaces and the contact lines
are
\begin{equation*}
\begin{aligned}
  \Sigma&=\{(x',\zeta(x')):x'\in D\},
  &\Gamma&=\{(x',\zeta(x')):x'\in\partial D\},\\
  \Sigma(t)&=\{(x',\zeta(x')+\eta(x',t)):x'\in D\},
  &\Gamma(t)&=\{(x',\zeta(x')+\eta(x',t)):x'\in\partial D\}.
\end{aligned}
\end{equation*}
The wetted solid boundaries are
$\Sigma_s=\partial\Omega\setminus\overline\Sigma$ and
$\Sigma_s(t)=\partial\Omega(t)\setminus\overline{\Sigma(t)}$.
They include the bottom of the vessel. On the vertical side wall the fluid
at time $t$ wets the points with $x'\in\partial D$ and
$h_0<z<\zeta(x')+\eta(x',t)$. The free and solid faces meet along
$\Gamma(t)=\overline{\Sigma(t)}\cap\overline{\Sigma_s(t)}$.
Figure~\ref{fig:contact_line} shows the equilibrium configuration.

The unknowns are the velocity $u$ and the pressure $p$ on $\Omega(t)$
and the height perturbation $\eta$ on $D$. We write $g>0$ for the
gravitational acceleration, $\mu>0$ for the viscosity, $\sigma>0$ for
the surface tension and $\beta_0>0$ for the coefficient in the Navier
slip condition. The stress is $S(p,u)=pI-\mu(\nabla u+(\nabla u)^T)$.
The upward non-unit normal to the free surface is
$\mathcal N=(-\nabla_{x,y}(\zeta+\eta),1)^T$. We write $n_0$ for the
outward unit normal on the solid wall and $\nu$ for the outward unit
normal to $\partial D$, so $n_0=(\nu,0)$ on the vertical side wall.
The subscript $\tan$ denotes projection onto the tangent plane of the wall.
The contact-line friction function $\mathcal W$ is smooth and increasing,
with $\mathcal W(0)=0$, $\kappa=\mathcal W'(0)>0$, and
$\hat{\mathcal W}(z)=\mathcal W(z)-\kappa z$.

With the stress written relative to equilibrium, the equations of motion are
\begin{equation}  \label{eq:moving-domain-full-system_intro} \left\{ \begin{alignedat}{2} &\partial_tu + u\cdot\nabla u + \operatorname{div} S(p,u) =0 &&\qquad \text{in }\Omega(t), \\
&\operatorname{div}u =0 &&\qquad \text{in }\Omega(t), \\
&S(p,u)\mathcal N = \left\{ g\eta -\sigma \operatorname{div}_{x,y} \left[ \frac{\nabla_{x,y}(\zeta+\eta)} {\sqrt{1+|\nabla_{x,y}(\zeta+\eta)|^2}} - \frac{\nabla_{x,y}\zeta} {\sqrt{1+|\nabla_{x,y}\zeta|^2}} \right] \right\} \mathcal N &&\qquad \text{on }\Sigma(t), \\
&\partial_t\eta = u\cdot\mathcal N &&\qquad \text{on }\Sigma(t), \\
&u\cdot n_0 =0 &&\qquad \text{on }\Sigma_s(t), \\
&\bigl( S(p,u)n_0 \bigr)_{\tan} -\beta_0u_{\tan} =0 &&\qquad \text{on }\Sigma_s(t), \\
&-{\kappa \partial_t\eta}-\hat{\mathcal{W}}(\partial_{t}\eta) = \sigma \left[ \frac{\partial_\nu(\zeta+\eta)} {\sqrt{1+|\nabla_{x,y}(\zeta+\eta)|^2}} - \frac{\partial_\nu\zeta} {\sqrt{1+|\nabla_{x,y}\zeta|^2}} \right] &&\qquad \text{on }\Gamma(t), \\
&u(\cdot,0)=u_0^{\sharp} &&\qquad\text{in }\Omega(0), \\
&\eta(\cdot,0)=\eta_0 &&\qquad\text{in }D.
\end{alignedat} \right.
\end{equation}

The first two equations give the incompressible Navier--Stokes system.
The next two are the balance of stress and the kinematic condition on
the free surface. On the wetted wall we impose impermeability and the
Navier slip condition. The contact-line equation is written in the
horizontal coordinates of $\partial D$, which parameterize $\Gamma(t)$.
Its physical origins are discussed in Subsection~\ref{subsec:model_origins}.
The initial height $\eta_0$ determines $\Omega(0)$ by
\eqref{eq:intro-fluid-domains}, and the initial velocity $u_0^\sharp$ is
prescribed on that domain. The same formula gives
$|\Omega(t)|-|\Omega|=\int_D\eta(x',t)\,dx'$.
Incompressibility, the kinematic condition and impermeability of the
solid wall imply conservation of fluid volume. Thus the condition
$\int_D\eta_0=0$ gives $\int_D\eta(\cdot,t)=0$ throughout the evolution,
as shown in Subsection~\ref{subsec:volume}.

\subsection{Equilibria}\label{subsec:steady_state}

\begin{figure}[ht]
\centering
\begin{tikzpicture}[x=1.18cm,y=1.18cm,>=stealth,font=\small,line cap=round,line join=round]
\definecolor{clsurface}{RGB}{47,101,145}
\definecolor{clcurve}{RGB}{167,40,49}
\definecolor{cltangent}{RGB}{23,115,89}
\definecolor{clnormal}{RGB}{176,103,29}
\definecolor{cltransverse}{RGB}{108,70,139}
\node[anchor=west,font=\small\bfseries] at (-2.65,2.85) {(a) Vessel and contact frame};
\path[fill=black!2] (-2,1.8) -- (-2,-1.25) .. controls (-2,-2.45) and (-0.95,-2.8) .. (0,-2.8) .. controls (0.95,-2.8) and (2,-2.45) .. (2,-1.25) -- (2,1.8) arc[start angle=0,end angle=180,x radius=2,y radius=0.7] -- cycle;
\path[fill=clsurface!9] (-2,0) -- (-2,-1.25) .. controls (-2,-2.45) and (-0.95,-2.8) .. (0,-2.8) .. controls (0.95,-2.8) and (2,-2.45) .. (2,-1.25) -- (2,0) -- plot[domain=0:-180,samples=101,variable=\clt] ({2*cos(\clt)},{0.7*sin(\clt)+0.25*sin(2*\clt)}) -- cycle;
\path[fill=clsurface!19] plot[domain=0:360,samples=181,variable=\clt] ({2*cos(\clt)},{0.7*sin(\clt)+0.25*sin(2*\clt)}) -- cycle;
\draw[black!65,line width=0.85pt] (-2,1.8) -- (-2,-1.25) .. controls (-2,-2.45) and (-0.95,-2.8) .. (0,-2.8) .. controls (0.95,-2.8) and (2,-2.45) .. (2,-1.25) -- (2,1.8);
\draw[black!65,line width=0.85pt] (0,1.8) ellipse[x radius=2,y radius=0.7];
\draw[black!60,thin] (1.55,2.24) -- (2.28,2.60);
\node[anchor=west] at (2.28,2.60) {rim};
\node[text=black!65] at (-0.30,0.98) {wall above the fluid};
\draw[clcurve,line width=1.25pt,dashed] plot[domain=0:180,samples=101,variable=\clt] ({2*cos(\clt)},{0.7*sin(\clt)+0.25*sin(2*\clt)});
\draw[clcurve,line width=1.25pt] plot[domain=180:360,samples=101,variable=\clt] ({2*cos(\clt)},{0.7*sin(\clt)+0.25*sin(2*\clt)});
\draw[clcurve,thin] ({2*cos(225)},{0.7*sin(225)+0.25*sin(450)}) -- (-2.48,-0.48);
\node[anchor=east,text=clcurve] at (-2.48,-0.48) {$\Gamma$};
\node[text=clsurface] at (-0.55,0.08) {$\Sigma$};
\node at (-0.65,-1.55) {$\Omega$};
\node[anchor=east] at (-2.10,-1.45) {$\Sigma_s$};
\pgfmathsetmacro{\clphi}{300}
\pgfmathsetmacro{\clzp}{0.5*cos(2*\clphi)}
\pgfmathsetmacro{\clnorm}{sqrt(4+\clzp*\clzp)}
\pgfmathsetmacro{\clpx}{2*cos(\clphi)}
\pgfmathsetmacro{\clpy}{0.7*sin(\clphi)+0.25*sin(2*\clphi)}
\coordinate (clP) at (\clpx,\clpy);
\draw[->,line width=1.0pt,cltangent] (clP) -- ++({1.35*(-2*sin(\clphi))/\clnorm},{1.35*(0.7*cos(\clphi)+\clzp)/\clnorm}) node[anchor=west,xshift=2pt] {$Qe_\theta$};
\draw[->,line width=1.0pt,clnormal] (clP) -- ++({1.25*(-cos(\clphi))},{1.25*(-0.35*sin(\clphi))}) node[anchor=south east,yshift=1pt] {$Qe_r$};
\draw[->,line width=1.0pt,cltransverse] (clP) -- ++({1.25*\clzp*sin(\clphi)/\clnorm},{1.25*(2-0.35*\clzp*cos(\clphi))/\clnorm}) node[anchor=east,xshift=-3pt] {$Qe_h$};
\fill[clcurve] (clP) circle[radius=1.7pt];
\node[below left,xshift=-1pt,yshift=-2pt] at (clP) {$P$};
\node[font=\footnotesize,text=black!65] at (0,-3.13) {Front wall shown transparently.};
\begin{scope}[shift={(6.8,-0.40)}]
\node[anchor=west,font=\small\bfseries] at (-0.70,3.25) {(b) Transverse chart at $P$};
\path[fill=clsurface!12] (0,0) -- (2.85,{2.85*tan(25)}) -- (2.85,-1.70) -- (0,-1.70) -- cycle;
\draw[black!60,line width=0.9pt] (0,-1.70) -- (0,1.95);
\draw[clsurface,line width=1.15pt] (0,0) -- (2.85,{2.85*tan(25)}) node[above right] {$h=-r\cot\theta_{\mathrm{eq}}$};
\draw[->,black!45,thin] (0,0) -- (3.10,0) node[right,text=black] {$r$};
\draw[->,black!55,thin] (0,1.95) -- (0,2.30) node[above,text=black] {$h$};
\draw[->,black!80,thin] (0,-0.68) arc[start angle=-90,end angle=25,radius=0.68];
\node at (0.94,-0.53) {$\theta_{\mathrm{eq}}$};
\fill[clcurve] (0,0) circle[radius=1.7pt];
\node[anchor=east,text=clcurve] at (-0.13,-0.03) {$\Gamma_0$};
\node[anchor=east] at (-0.13,-1.27) {wall};
\node[text=clsurface] at (2.12,-1.22) {fluid};
\node[text=black!65,font=\footnotesize,align=center] at (1.35,-2.73) {$Qe_r,Qe_h$ span $T(P)^\perp$.\\$\theta_{\mathrm{eq}}$ is measured through the fluid.};
\end{scope}
\end{tikzpicture}
\caption{Left: the fluid domain $\Omega$, the free surface $\Sigma$, the solid wall $\Sigma_s$, and the contact line $\Gamma$ where the free surface meets the wall. The solid vessel wall is smooth and
vertical near $\Gamma$, and extends above the contact line. At $P\in\Gamma$, the vectors $Qe_\theta$, $Qe_r$, and $Qe_h$ form an orthonormal frame: $Qe_\theta$ is tangent to $\Gamma$, while $Qe_r$ and $Qe_h$ span its normal plane. Right: the corresponding transverse chart, in which the wall and
free surface are straightened and meet at the angle $\protect\theta_{\mathrm{eq}}$ prescribed by the Young relation \eqref{eq:young}, measured through the fluid. This corner limits
the elliptic regularity.}
\label{fig:contact_line}
\end{figure}
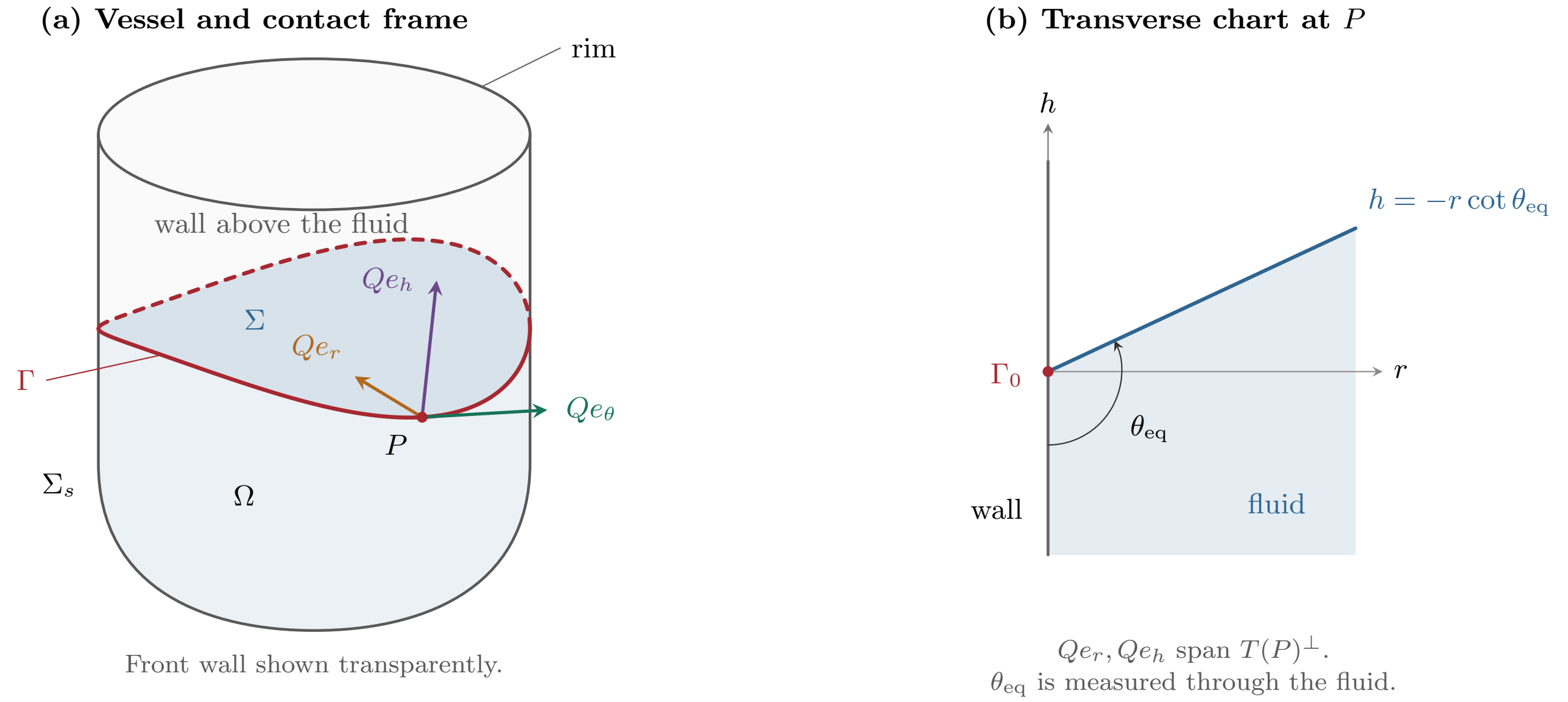

The equilibrium shapes of a liquid at rest in a container were studied by Gauss, Laplace and Young at the beginning of the nineteenth century. Young and Laplace found that the pressure jump across a fluid interface is proportional to its mean curvature, and Gauss derived the equilibrium from a variational principle: the equilibrium surface is a critical point of the energy
\begin{equation*}
\mathcal{F}(S)=\sigma\,|S|+\big(\gamma_{sf}-\gamma_{sv}\big)\,|W(S)|+\int g\,z,
\end{equation*}
in which $\sigma$ is the surface tension of the free surface $S$, $\gamma_{sf}$ and $\gamma_{sv}$ are the tensions of the solid--fluid and solid--vapour interfaces, $W(S)$ is the wetted part of the wall, and the last term is the gravitational energy, all at fixed volume. Varying the surface gives the Young--Laplace equation for the equilibrium height $\zeta$ over $D$,
\begin{equation}  \label{eq:young-laplace} g\zeta-\sigma\func{div}\!\left(\frac{\nabla\zeta}{\sqrt{1+|\nabla\zeta|^{2}}}\right)=\text{const}\quad\text{in }D,
\end{equation}
and varying the wetted region gives the boundary condition at the contact line, the Young relation
\begin{equation}  \label{eq:young} \cos\theta_{\mathrm{eq}}=\frac{\gamma_{sv}-\gamma_{sf}}{\sigma},
\end{equation}
which fixes the angle $\theta_{\mathrm{eq}}$ between the free surface and the wall, measured through the fluid~\cite{Young}. The angle depends only on the three materials in contact, and not on the shape of the vessel or on the amount of liquid. Throughout, $\theta_{\mathrm{eq}}$ denotes this angle (Figure~\ref{fig:contact_line}); it is the opening angle of the fluid at the contact line. On the vertical wall the graph satisfies $\cos\theta_{\mathrm{eq}}=\partial_\nu\zeta/\sqrt{1+|\nabla_{x,y}\zeta|^2}$ on $\partial D$. We assume that $\zeta$ solves \eqref{eq:young-laplace}--\eqref{eq:young}, so that $\theta_{\mathrm{eq}}$ is constant along the contact line.

In two dimensions the contact set is a pair of points. In three dimensions it is a closed curve $\Gamma$, a curve of triple junctions where the liquid, the solid and the vapour meet, and the whole difficulty of this paper lies in the analysis along it. Our result is that this equilibrium is asymptotically stable: a small perturbation of the surface and of the fluid decays, and the contact line returns to $\Gamma$.

\subsection{Previous work and the origin of the model}\label{subsec:model_origins}

Contact lines, and contact points in two dimensions, form at triple junctions between three phases, here fluid, solid and vapour, and have been studied since the work of Young~\cite{Young} in 1805. For an exhaustive overview we refer to de Gennes~\cite{deGennes}; what follows is a brief review. The story begins with the equilibrium configurations studied by Young~\cite{Young}, Laplace~\cite{Laplace} and Gauss~\cite{Gauss}, who found the variational principle of Subsection~\ref{subsec:steady_state} and the relation \eqref{eq:young} between the equilibrium contact angle and the free-energy parameters $\gamma_{sf}$, $\gamma_{sv}$ and $\sigma$. For the equilibrium theory and the solvability of \eqref{eq:young-laplace}--\eqref{eq:young} we refer to Finn~\cite{Finn}.

The dynamics of a contact line are more delicate. The free boundary kinematics and the no-slip condition $u=0$ on the solid are incompatible: together they force $\partial_{t}\eta=0$ at the contact line, that is, the fluid is pinned where it started, which everyday experience contradicts. The no-slip condition must therefore be abandoned. The surveys of Dussan~\cite{Dussan} and Blake~\cite{Blake} discuss the efforts of physicists and chemists to determine the correct law. Write $V_{\mathrm{cl}}=\partial_t\eta$ for the vertical speed of the contact line at fixed $x'\in\partial D$, taken positive upwards. With the dynamic angle $\theta_{\mathrm{dyn}}$ and the equilibrium angle $\theta_{\mathrm{eq}}$ both measured through the fluid, the response law has the form
\begin{equation}  \label{eq:contact_law_heuristic} V_{\mathrm{cl}}=F\big(\cos\theta_{\mathrm{eq}}-\cos\theta_{\mathrm{dyn}}\big),
\end{equation}
where $F$ is increasing with $F(0)=0$. In the graph coordinates of
\eqref{eq:moving-domain-full-system_intro},
$\cos\theta_{\mathrm{dyn}}=\partial_\nu(\zeta+\eta)/\sqrt{1+|\nabla_{x,y}(\zeta+\eta)|^2}$
on $\partial D$, so the contact-line condition reads
$\mathcal W(V_{\mathrm{cl}})=\sigma(\cos\theta_{\mathrm{eq}}-\cos\theta_{\mathrm{dyn}})$.
In particular, $F(s)=\mathcal W^{-1}(\sigma s)$ near $s=0$.
In the angle convention of~\cite{GT}, measured through the vapour, both cosines have the opposite sign. Laws of this form have been derived in several ways: by Blake and Haynes~\cite{BlakeHaynes} from thermodynamic and molecular kinetics arguments, by Cox~\cite{Cox} by matched asymptotics, and by Ren and E~\cite{RenE11} from thermodynamic principles; the molecular dynamics simulations of Ren and E~\cite{RenE} also produced a law of this form, and indicated that the fluid slips along the solid according to the Navier condition. The system studied here combines the Navier slip condition with the general form \eqref{eq:contact_law_heuristic}, written through the friction function $\mathcal{W}$.

For this model in two dimensions, Guo and Tice proved the asymptotic stability of the equilibria for the Stokes equations in~\cite{GT18}, and for the Navier--Stokes equations in~\cite{GT}; Tice and Wu~\cite{TiceWu} proved the corresponding results for the Stokes droplet problem, and Tice and Zheng~\cite{TiceZheng} established local existence. Global well-posedness in two dimensions, that is, the construction of solutions satisfying the dynamic contact line law together with all the other boundary conditions, was obtained in~\cite{GTWYZ}. The present paper proves the a~priori estimates in three dimensions, where the contact set is a closed curve rather than a pair of points; the construction of solutions in three dimensions is not treated here.

We are not aware of other results in which the full fluid mechanics is considered together with a dynamic contact line and a dynamic contact angle, but there are results with a subset of these features. Schweizer~\cite{Schweizer97,Schweizer01} studied a two-dimensional Navier--Stokes problem with the contact angle fixed at $\pi/2$, and Bodea~\cite{Bodea} a similar problem in a three-dimensional channel, periodic in one direction. Bravin et~al.~\cite{BGKMRS} proved well-posedness and regularity in weighted Sobolev spaces for the stationary Stokes equations in a fixed planar wedge with Navier slip on both sides. Gnann, Kn\"upfer, Lam,
Masmoudi and Sauer~\cite{GKLMS} proved local well-posedness for two-dimensional Stokes flow with a moving contact point and a fixed contact angle, for sufficiently small angles and small compatible initial data.

Kn\"upfer and Masmoudi~\cite{KM13,KnupferMasmoudi} studied the dynamics of a two-dimensional drop with fixed contact angle for a fluid governed by Darcy's law. For the stationary Navier--Stokes system with a free but unmoving boundary, the analysis was carried out in two dimensions by Solonnikov~\cite{Solonnikov95,Solonnikov98}, in three dimensions with angle $\pi/2$ by Jin~\cite{Jin}, by Socolowsky~\cite{Socolowsky} for two-dimensional coating problems, and by Kr\"oner~\cite{Kroner} for dynamic contact angles. For inviscid fluids there has been progress on models with floating or emerging solids, without a dynamic law for the contact angle: see Lannes~\cite{Lannes}, de Poyferr\'e~\cite{dePoyferre} and Bocchi~\cite{Bocchi}. A simplified droplet model without fluid coupling was studied by Feldman and Kim~\cite{FeldmanKim}, who proved asymptotic stability by gradient flow techniques. Much work has also been done on contact points in simplified thin-film models; we refer to the survey of Bertozzi~\cite{Bertozzi}.  For a thin-film model in the complete-wetting regime, Giacomelli, Gnann and Peschka~\cite{GGP} used formal asymptotics to determine how
contact-line friction affects the long-time spreading rate of a drop.

\subsection{Reformulation on the equilibrium domain}

We pull \eqref{eq:moving-domain-full-system_intro} back to the equilibrium
fluid domain $\Omega$ of \eqref{eq:intro-fluid-domains}.
Let $\bar\eta$ be the extension of $\eta$ constructed in
Section~\ref{sec:extension}. Choose the weight $W$ of
\eqref{eq:W-definition} so that $W=1$ on $\Sigma$ and $W=0$ below $h_0$.
The change of variables is
\begin{equation*}
 \Phi_\eta(x',z,t)=(x',z+W(x',z)\bar\eta(x',z,t)),\qquad
 \Phi_\eta(\cdot,t):\Omega\longrightarrow\Omega(t).
\end{equation*}
For the small perturbations considered here this map is a diffeomorphism.
It is the identity in the lower portion of the vessel and sends
$(x',\zeta(x'))$ to $(x',\zeta(x')+\eta(x',t))$. In particular it
carries $\Sigma$, $\Sigma_s$ and $\Gamma$ onto $\Sigma(t)$, $\Sigma_s(t)$
and $\Gamma(t)$, respectively.

While taking the pullback, write $(\widetilde u,\widetilde p)$ for the
physical fields $(u,p)$ in \eqref{eq:moving-domain-full-system_intro}.
For $X\in\Omega$ we set
\begin{equation*}
 v(X,t)=\widetilde u(\Phi_\eta(X,t),t),\qquad
 p(X,t)=\widetilde p(\Phi_\eta(X,t),t).
\end{equation*}
The height perturbation $\eta$ is already defined on $D$.
Write $K=(\det D\Phi_\eta)^{-1}$ and $\mathcal A_\eta=(D\Phi_\eta)^{-T}$.
In $\nabla_{\mathcal A_\eta}$, $\operatorname{div}_{\mathcal A_\eta}$
and $S_{\mathcal A_\eta}$, each Cartesian derivative $\partial_i$ is
replaced by $(\mathcal A_\eta)_{ij}\partial_j$.
For a graph of height $f$ we use the curvature convention
\begin{equation*}
 H(f)=\operatorname{div}_{x,y}\left(
 \frac{\nabla_{x,y} f}{\sqrt{1+|\nabla_{x,y} f|^2}}\right).
\end{equation*}
The fixed-domain problem for $(v,p,\eta)$ is
\begin{equation}\label{eq:intro-fixed-system}
\left\{\begin{alignedat}{2}
 &\partial_t v-WK\partial_t\bar\eta\,\partial_3 v
    +v\cdot\nabla_{\mathcal A_\eta}v
    +\operatorname{div}_{\mathcal A_\eta}S_{\mathcal A_\eta}(p,v)=0
       &&\qquad\text{in }\Omega,\\
 &\operatorname{div}_{\mathcal A_\eta}v=0
       &&\qquad\text{in }\Omega,\\
 &S_{\mathcal A_\eta}(p,v)\mathcal N
       =\bigl(g\eta-\sigma[H(\zeta+\eta)-H(\zeta)]\bigr)\mathcal N
       &&\qquad\text{on }\Sigma,\\
 &\partial_t\eta=v\cdot\mathcal N
       &&\qquad\text{on }\Sigma,\\
 &v\cdot n_0=0
       &&\qquad\text{on }\Sigma_s,\\
 &\bigl(S_{\mathcal A_\eta}(p,v)n_0\bigr)_{\tan}-\beta_0v_{\tan}=0
       &&\qquad\text{on }\Sigma_s,\\
 &-\mathcal W(\partial_t\eta)=\sigma\left[
    \frac{\partial_\nu(\zeta+\eta)}{\sqrt{1+|\nabla_{x,y}(\zeta+\eta)|^2}}
   -\frac{\partial_\nu\zeta}{\sqrt{1+|\nabla_{x,y}\zeta|^2}}
   \right] &&\qquad\text{on }\Gamma,\\
 &v(\cdot,0)=u_0^\sharp\circ\Phi_\eta(\cdot,0)
       &&\qquad\text{in }\Omega,\\
 &\eta(\cdot,0)=\eta_0 &&\qquad\text{in }D.
\end{alignedat}\right.
\end{equation}
Here $\mathcal N$ is evaluated at the same horizontal point $x'$ as in
the moving-domain problem. Section~\ref{sec:whole_space} derives
\eqref{eq:intro-fixed-system} starting from \eqref{eq:moving-domain-full-system}.
Separating the linear and nonlinear parts of the capillary and contact-line
terms gives the form \eqref{eq:moving-domain-elliptic-system}, with the
nonlinear terms specified in \eqref{eq:F05}. There $\Gamma_0$ denotes the
contact line in the coordinates of the local chart.

\subsection{Main result}

The energy $\mathcal E$ and dissipation $\mathcal D$ of
Subsection~\ref{subsec:energy_def} collect the norms of the time
derivatives, the additional surface and pressure norms, and the norms in
the interior, contact-line and intermediate regions. In the notation of
that subsection we have
\begin{equation*}
\begin{aligned}
 \mathcal{E}&=\mathcal{E}_{||}+\mathcal{E}_{\mathrm{enhanced}}
     +\mathcal{E}_{\mathrm{int}}+\mathcal{E}_{\mathrm{cl}}+\mathcal{E}_{\mathrm{bet}},\\
 \mathcal{D}&=\mathcal{D}_{||}+\mathcal{D}_{\mathrm{enhanced}}
     +\mathcal{D}_{\mathrm{int}}+\mathcal{D}_{\mathrm{cl}}+\mathcal{D}_{\mathrm{bet}}.
\end{aligned}
\end{equation*}
All estimates are a~priori estimates for sufficiently regular solutions.

\smallskip\noindent\textbf{Hypotheses.}

\begin{enumerate}[label=(\arabic*)]

\item $\theta_{\mathrm{eq}}\in(0,\pi)$, and $\zeta$ satisfies Young's law, so that $\theta_{\mathrm{eq}}$ is constant along $\Gamma$.

\item $\partial D$ is of class $C^{\infty}$ and is connected, and the equilibrium height is smooth up to the boundary, $\zeta\in C^{\infty}(\overline{D})$; the latter follows from the former by elliptic regularity for the capillarity problem \eqref{eq:young-laplace}--\eqref{eq:young}. The side wall is vertical on a neighbourhood of the free surface that contains the level $h_{0}$ of Subsection~\ref{subsec:cutoff_extension}.

\item The closure of the solid boundary $\Sigma_{s}=\partial\Omega\setminus\overline{\Sigma}$ is a $C^{\infty}$ surface with boundary $\Gamma$, without edges or corners; in particular the bottom of the vessel is smooth. Thus the only singular set of $\partial\Omega$ is the contact line $\Gamma$, where $\Sigma$ meets $\Sigma_{s}$ at the angle $\theta_{\mathrm{eq}}$, and $\Omega$ is a bounded Lipschitz domain.

\item $\beta_{0}>0$ and $\kappa>0$.

\item The perturbation has zero mean, $\int_{D}\eta=0$; see Subsection~\ref{subsec:volume}.

\item $m\geq3$, and the parameters satisfy the conditions of Subsection~\ref{subsec:constants}.
\end{enumerate}

\begin{maintheorem}
Assume (1)--(6), and let $\mathcal{E}$ and $\mathcal{D}$ be the energy and the dissipation of Subsection~\ref{subsec:energy_def}. There are constants $\delta_{*}>0$, $C>0$ and $\lambda>0$ with the following property. If $(v,p,\eta)$ is a sufficiently regular solution on $[0,T)$ of \eqref{eq:intro-fixed-system} (equivalently, of \eqref{eq:moving-domain-elliptic-system} with \eqref{eq:F05}) satisfying $\mathcal{D}\in L^1_{\mathrm{loc}}([0,T))$ and $\sup_{0\leq s<T}\mathcal{E}(s)\leq\delta_{*}$, then
\begin{equation*}
\mathcal{E}(t)+\int_{s}^{t}\mathcal{D}\leq C\,\mathcal{E}(s) \qquad\text{and}\qquad \mathcal{E}(t)\leq C\,e^{-\lambda t}\,\mathcal{E}(0) \qquad\text{for }0\leq s\leq t<T.
\end{equation*}
\end{maintheorem}

The proof is completed at the end of Subsection~\ref{subsec:cl_local_energy}. Near the contact line, at the purely angular top level, it uses Theorem~\ref{thm:top_level} of Subsection~\ref{subsec:top_level}. Its hypotheses hold under (1)--(6): (H1), (H2) and (H4) by the construction of the cutoffs and of the extension (Subsections~\ref{subsec:cutoff}--\ref{subsec:cutoff_extension} and Section~\ref{sec:extension}) and by the definition of $\mathcal{D}_{\mathrm{int}}$, and (H3) because the solution is sufficiently regular. Its remainder bounds are proved in Lemma~\ref{lem:tl_rem}, and its right-hand side is controlled by $\mathcal{D}$ itself, with no further hypothesis.  

In forthcoming work we plan to address the construction and continuation of solutions in the functional setting of the Main Theorem. Together with the a priori estimate proved here, such a theory would yield global existence and exponential convergence to equilibrium for a class of compatible initial data with sufficiently small energy. This would then establish the  asymptotic stability of the equilibrium within that class.

\subsection{Technical overview}

The two-dimensional problem was treated by Guo and Tice~\cite{GT}, and the basic architecture of the proof follows theirs: an energy--dissipation structure for time derivatives, an enhancement of the regularity of the free surface to $H^{\frac32-\alpha}$, and elliptic estimates whose integrability index near the contact line is limited by the corner exponent of Orlt and S\"andig~\cite{OS}.  However, the three-dimensional domain creates numerous fundamental new difficulties.

The corner at the contact line limits the elliptic regularity: for the Stokes problem with our boundary conditions, the exponent of Orlt and S\"andig~\cite{OS} allows only $W^{2,q}$ with $q=q_{\varepsilon}=2/(2-\varepsilon)$ close to $1$ when the contact angle is large. In two dimensions this is enough, because the slice is the whole domain; in three dimensions a norm with $q$ near $1$ cannot control the nonlinear terms, as the Sobolev exponents do not close. Three features are new in three dimensions.

\begin{itemize}
\item \emph{Tangential derivatives along the contact line.} Near the contact line $\Gamma$ the problem is written in a tubular chart $(\theta,r,h)$ (Section~\ref{sec:flattening}) and estimated on the two-dimensional slices $\{\theta=\mathrm{const}\}$ (Sections~\ref{sec:preparation} and~\ref{sec:local_charts}). The angular derivatives $\partial_{\theta}$ are estimated by an energy method, by induction on their number (Section~\ref{sec:contact_energy}).

\item \emph{Transport of the contact line along itself.} The velocity has a component along $\Gamma$, so the kinematic forcing on the contact line contains $v_{\gamma}\partial_{\gamma}\eta$, which vanishes in two dimensions. At the purely angular top level this term requires $\partial_{\theta}^{m+2}\eta\in L^{2}(\Gamma_{0})$, one tangential derivative more than the $H^{\frac32-\alpha}$ bound provides. It is controlled by the estimate without loss of Subsection~\ref{subsec:top_level}, which gives $\partial_{\theta}^{m+1}\eta\in H^{\frac32}$ near $\Gamma$.

\item \emph{Anisotropy and thin cutoffs.} The extension of $\eta$ near $\Gamma$ is a slice-wise Poisson extension (Section~\ref{sec:extension}), whose second $\theta$-derivatives are not controlled by the $H^{\frac32}$ norm of the data; top-order terms are therefore integrated by parts in all variables. The cutoffs near $\Gamma$ have width $\varepsilon$, and their derivatives cost $\varepsilon^{-1}$; the elliptic and Ne\caronaccent{c}as constants are taken on domains whose size does not depend on $\varepsilon$ (Subsection~\ref{subsec:artificial}).
\end{itemize}

Our scheme rests on two vector fields and three regions. The two vector fields are the temporal field $\partial_{t}$, as in \cite{GT}, and the tangential field $\partial_{\theta}$ along the contact line, which we introduce here. Both are tangent to the geometry: $\partial_{\theta}$ generates a flow that preserves $\Gamma$, so applying it to the system produces no boundary term of top order. We work with the mixed fields $\partial_{t}^{k_{1}}\partial_{\theta}^{k_{2}}$ with $k_{1}+k_{2}\leq m+1$, $m\geq3$.

The three regions are the interior, away from the contact line, where the elliptic theory is the classical one and the exponent may be taken close to $2$; the contact-line region, where the problem is reduced, slice by slice in  $\theta$, to the two-dimensional corner problem of \cite{GT} and only $q$ near $1$ is available; and an intermediate region between them, where the integrability is bootstrapped from $q$ up to an exponent close to $2$. In the slicing norms used near the contact line, the first space acts in $\theta$ and the second in the transverse variables $(r,h)$.

The proof is a single induction. Within it, an energy estimate for a field needs the elliptic estimates one level lower, the elliptic estimates near the contact line need the interior ones, and the interior ones need a bootstrap of integrability across the intermediate region. The sections below carry this out in order, and the introduction follows that order.

\subsection{The main estimates}

We state here the main theorem of each section, in the order in which they are proved; each is proved in the section indicated. The energy $\mathcal{E}$, the dissipation $\mathcal{D}$ and their pieces are defined in Subsection~\ref{subsec:energy_def}, the parameters in Subsection~\ref{subsec:constants}, and all estimates are a~priori estimates under the smallness assumption $\sup_{0\leq t<T}\mathcal{E}(t)\leq\gamma^{2}$.

\smallskip\noindent\textbf{Geometry (Sections~\ref{sec:flattening}--\ref{sec:extension}).} Near $\Gamma$ we use the chart $\Phi$ of Section~\ref{sec:flattening}, in which $\theta$ is equivalent to the arclength parameter along the contact line and $(r,h)$ are transverse. Theorem~\ref{thm:localized-surface-coordinate-estimate} and Theorem~\ref{cor:localized-norm-equivalence} compare the norms of $\eta$ with those of its chart representation $\mathfrak{h}_{\eta}$, so that the two descriptions may be used interchangeably. The free surface is then flattened by an extension $\bar{\eta}$ of $\eta$ into the bulk, built in Section~\ref{sec:extension} from the half-space extender of Theorem~\ref{thm:half-space-extender} and a Poisson extension; Theorem~\ref{thm:extension_trace} says that $\bar{\eta}=\eta$ on $\Sigma$, and Theorem~\ref{thm:extension-commutator-estimate} is the commutator estimate that lets a cutoff pass through the extension. The extension near the contact line is the \emph{slice-wise}, two-dimensional Poisson extension, taken at each fixed $\theta$. This is what makes $\partial_{\theta}$ commute with the extension, and it is why the angular derivatives of $\bar{\eta}$ obey slicing bounds rather than isotropic ones; the price is paid in Theorem~\ref{thm:top_level}.

\smallskip\noindent\textbf{The temporal field (Section~\ref{sec:whole_space}).} The starting point is the energy--dissipation relation for the fields $\partial_{t}^{j}$, $j\leq m+1$.

\begin{theorem}[Energy--dissipation relation for $\partial_{t}$]
\label{thm:whole-space-energy} We have the energy--dissipation relation
\begin{align}
\mathcal{E}_{||}(t)+\int_{0}^{t}\mathcal{D}_{||}\lesssim \mathcal{E}_{||}(0)+\mathcal{E}^{\frac{3}{2}}(t)+\mathcal{E}^{\frac{3}{2}}(0)+\int_{0}^{t}\sqrt{\mathcal{E}}\mathcal{D}.
\end{align}
\end{theorem}

\noindent The nonlinear terms here are out of reach of isotropic three-dimensional estimates when the contact angle is large. Near the contact line they are projected to the two-dimensional estimates of \cite{GT} by Lemma~\ref{lem:derivative_replacement}, which trades transverse derivatives for angular ones, and by Lemma~\ref{thm:contact_line_forcing_estimate}, which collects the resulting bounds on the forcing. The relation above is not enough to close the nonlinearity: it controls $\partial_{t}^{m+1}\eta$ only in $H^{1}$, and it says nothing about the pressure. Both are supplied by the two enhancements, the gain of $3/2$ derivatives for the surface (Subsection~\ref{subsec:enhanced_surface}, through the harmonic field of Proposition~\ref{thm:whole_space_threehalf}) and the pressure bound of Theorem~\ref{thm:pressure_enhancement} (Subsection~\ref{subsec:enhanced_pressure}, through the Bogovski\u{\i} field of Proposition~\ref{prop:B_d}, constructed as in \cite{Galdi,CM}, and the average identity of Theorem~\ref{thm:energy_pressure}).

We note that, up to the nonlinearity, the globally defined $D_t$ energy estimate is self-contained; it is the foundation of the paper. To circumvent the basic difficulty of low regularity near the contact line, we introduce the angular derivative in the collar (the ordinary $\partial_\theta$ on scalars) and control the nonlinearity there by the anisotropic slicing norm, schematically
\begin{align}
\|fg\|_{H^1} \lesssim \|fg\|_{L^2_\theta H^1_{r,h}} +\|\partial_\theta(fg)\|_{L^2_\theta L^2_{r,h}}.
\end{align}
It suffices to control the two-dimensional slicing norms $L^2_{r,h}$ and $H_{r,h}^{1}$ for fixed $\theta$, which are controlled as in \cite{GT}.
\begin{theorem}[Enhanced energy estimate]
\label{thm:whole_space_estimate_final}

The following inequality holds
\begin{align}
\mathcal{E}_{\mathrm{enhanced}}(t)+\int_{0}^{t}\mathcal{D}_{\mathrm{enhanced}}\lesssim \mathcal{E}(0)+\mathcal{E}^{\frac{3}{2}}(t)+\mathcal{E}^{\frac{3}{2}}(0)+\int_{0}^{t}\sqrt{\mathcal{E}}\mathcal{D}.
\end{align}
\end{theorem}

\noindent Here $\mathcal{E}_{\mathrm{enhanced}}$ contains $\partial_{t}^{j}\eta$ in $H^{3/2-\alpha}$ and $\partial_{t}^{j-1}p$ in $L^{2} $. The loss of $\alpha$ is intrinsic to the temporal field, and at the top level $j=m+1$ there is no estimate for $\partial_{t}^{m+1}p$ at all; this is the limitation that the angular field removes. We note that, up to the nonlinearity, the enhancement is self-contained.

\smallskip\noindent\textbf{The interior (Section~\ref{sec:interior}).} Away from the contact line the elliptic theory is the classical one, and the exponent may be taken close to $2$.

\begin{theorem}[Interior elliptic estimate]
\label{thm:interior_space_final}  The following estimate holds
\begin{align}
{\label{est:interior_final}} \mathcal{E}_{\mathrm{enhanced}}+\mathcal{E}_{\mathrm{int}}+\int_{0}^{t}(\mathcal{D}_{\mathrm{\mathrm{enhanced}}}+\mathcal{D_{\mathrm{int}}})\lesssim C(\varepsilon)\varepsilon^{-m-1}(\mathcal{E}(0)+\mathcal{E}^{\frac{3}{2}}(t)+\mathcal{E}^{\frac{3}{2}}(0))+\int_{0}^{t}\sqrt{\mathcal{E}}\mathcal{D}.
\end{align}
\end{theorem}

\noindent The difficulty is the term $\nabla^{2}\bar{\eta}\,\nabla u$, in which $\bar{\eta}$ is non-local: a cutoff in the interior still sees the surface near the contact line, where only the weak slicing norms are available. This is why the estimate is run with nested cutoffs, each step gaining one derivative, and why it needs the intermediate region below. The power $\varepsilon^{-m-1}$ is harmless; only its finiteness for fixed $\varepsilon$ is used.

We note that, up to the ($\epsilon$-dependent) nonlinearity, the interior estimate is self-contained. Controlling the nonlinearity in the strong $L^p$ norm ($p$ close to $2$) is delicate near the contact line, where only the weak $L^{q_\delta}$ integrability with $q_\delta$ near $1$ is available. To control the non-local nonlinearity $\chi D^2 \bar \eta \nabla v$, we must use the crucial extra gain in the commutator estimate of Theorem \ref{thm:enhanced_commutator} for the non-local $\bar\eta$. On the other hand, to control $||\chi v\nabla \eta||_{W^{2-1/p,p}(\Sigma)}$ on the surface $\Sigma$, we group $\chi$ with $\nabla\eta$ and then control $v$ via the slicing norm used in \eqref{eq:intrinsic_slicing_norm}.

\smallskip\noindent\textbf{The contact line, slice by slice (Sections~\ref{sec:preparation}--\ref{sec:local_charts}).} Section~\ref{sec:preparation} localizes the system near $\Gamma$ and flattens it. In the chart we decompose the velocity into its component $v_{\theta}$ along the contact line and its transverse part $v_{\perp}=(v_{r},v_{h})$. At each fixed $\theta $ the transverse part solves a two-dimensional Stokes problem in a corner, so the $W^{2,q}$--$W^{1,q}$--$W^{3-1/q,q}$ framework of \cite{GT} applies with $q=q_{\delta}$ close to $1$, while $v_{\theta}$ solves a scalar problem whose forcing involves the pressure. Theorem~\ref{thm:local-equation-h-eta} identifies the capillary operator in the chart, and the elliptic estimate at the top level is the following.

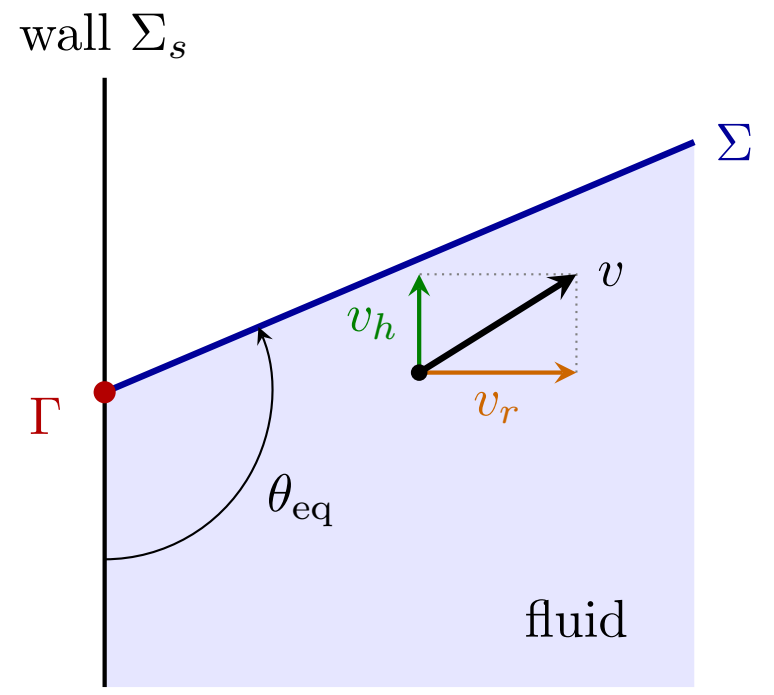
\begin{figure}[ht]
\centering
\begin{tikzpicture}[scale=1.3,>=stealth]
  \path[fill=blue!10] (0,0) -- (3.0,1.273) -- (3.0,-1.5) -- (0,-1.5) -- cycle;
  \draw[thick] (0,-1.5) -- (0,1.6) node[above] {wall $\Sigma_{s}$};
  \draw[very thick,blue!60!black] (0,0) -- (3.0,1.273) node[right] {$\Sigma$};
  \draw[->] (0,-0.85) arc (-90:23:0.85);
  \node at (1.0,-0.55) {$\theta_{\mathrm{eq}}$};
  \fill[red!70!black] (0,0) circle (1.6pt);
  \node[red!70!black,left] at (-0.12,-0.12) {$\Gamma$};
  \node at (2.4,-1.15) {fluid};
  \coordinate (P)  at (1.6,0.1);
  \coordinate (Vr) at (2.4,0.1);   
  \coordinate (Vt) at (1.6,0.6);   
  \coordinate (V)  at (2.4,0.6);   
  \draw[densely dotted,gray] (Vr) -- (V) -- (Vt);
  \draw[->,thick,orange!80!black] (P) -- (Vr) node[below,pos=0.5] {$v_r$};
  \draw[->,thick,green!50!black]  (P) -- (Vt) node[left,pos=0.5] {$v_h$};
  \draw[->,very thick] (P) -- (V) node[right,pos=1] {$v$};
  \fill (P) circle (1.2pt);
\end{tikzpicture}
\caption{A transverse slice in the local charts at a fixed point of the contact line $\Gamma$, with a flat free surface $\Sigma$ meeting the wall $\Sigma_{s}$ at the contact angle $\theta_{\mathrm{eq}}$, measured through the fluid. The velocity $v$ at a point of the fluid is decomposed as $v = v_r e_r + v_h e_h$, with $v_r$ horizontal and $v_h$ vertical.} \label{fig:contact_angle_flat}
\end{figure}

\begin{theorem}[Contact-line elliptic estimate]
\label{thm:Higher-order contact line elliptic estimate}

Suppose $v,p,\eta$ solve equation \eqref{eq:final-localized-flattened-system-j}. We then have the following elliptic estimate when $k_{1}+k_{2}=m$
\begin{align}
\begin{aligned} \mathscr E_{k_{1},k_{2},el,cl} \lesssim \mathcal{E}\mathcal{D}+\frac{C(\varepsilon)}{\varepsilon} \sum_{0\leq j_{1}\leq k_{1}}\sum_{0\leq j_{2}\leq k_{2}}\mathscr{E}_{j_{1},j_{2},\mathrm{energy},\mathrm{cl}}+\sum_{j=0}^{k_{2}+1}\mathscr{E}_{k_{1},j,\mathrm{energy},\mathrm{cl}}+\sum_{j=0}^{k_{2}}\mathscr{E}_{k_{1}+1,j,\mathrm{energy},cl}, \end{aligned}
\end{align}
and when $k_{1}+k_{2}\leq m-1$
\begin{align}
\begin{aligned} \mathscr E_{k_{1},k_{2},el,cl} \lesssim \mathcal{E}^{2}+\frac{C(\varepsilon)}{\varepsilon} \sum_{0\leq j_{1}\leq k_{1}}\sum_{0\leq j_{2}\leq k_{2}}\mathscr{E}_{j_{1},j_{2},\mathrm{energy},\mathrm{cl}}+\sum_{j=0}^{k_{2}+1}\mathscr{E}_{k_{1},j,\mathrm{energy},\mathrm{cl}}+\sum_{j=0}^{k_{2}}\mathscr{E}_{k_{1}+1,j,\mathrm{energy},cl}. \end{aligned}
\end{align}
\end{theorem}

\noindent The slicing elliptic estimate (see Figure~\ref{fig:contact_angle_flat}) is established by a two-dimensional argument as in \cite{GT}, after projecting the velocity field onto the transverse slice. To estimate $v_{\theta}$, we construct an elliptic problem with a Neumann boundary condition and bound its elliptic norm.  The exponent $q_{\delta}$ here is too small to control $\nabla^{2}\bar{\eta}\,\nabla u$ in three dimensions. On an intermediate region, away from the contact line but still inside the chart, the integrability is therefore raised from $q_{\delta}$ to the exponent $p$ close to $2$. This is done by the nested interior estimates of Section~\ref{sec:interior}, whose cutoffs have their transitions there; Theorem~\ref{thm:enhanced_commutator} is the commutator estimate that makes this possible. On the intermediate region itself, the estimate at the exponent $q_{\delta}$ is the following.

\begin{theorem}[Intermediate elliptic estimate]
\label{thm:intermediate_estimate}

The following intermediate elliptic estimate holds given that $k_{1}+k_{2}=m$
\begin{align}
\mathscr{E}_{k_{1},k_{2},el,bet}\lesssim \mathcal{E}\mathcal{D}+\frac{C(\varepsilon)}{\varepsilon} \sum_{0\leq j_{1}\leq k_{1}}\sum_{0\leq j_{2}\leq k_{2}}\mathscr{E}_{j_{1},j_{2},\mathrm{energy},\mathrm{bet}}+\sum_{j=0}^{k_{2}+1}\mathscr{E}_{k_{1},j,\mathrm{energy},\mathrm{bet}}+\sum_{j=0}^{k_{2}}\mathscr{E}_{k_{1}+1,j,\mathrm{energy},bet}+\frac{C(\varepsilon)}{\varepsilon} \mathcal{D}_{\mathrm{int}},
\end{align}
and when $k_{1}+k_{2}\leq m-1$
\begin{align}
\mathscr{E}_{k_{1},k_{2},el,bet}\lesssim \mathcal{E}^{2}+\frac{C(\varepsilon)}{\varepsilon} \sum_{0\leq j_{1}\leq k_{1}}\sum_{0\leq j_{2}\leq k_{2}}\mathscr{E}_{j_{1},j_{2},\mathrm{energy},\mathrm{bet}}+\sum_{j=0}^{k_{2}+1}\mathscr{E}_{k_{1},j,\mathrm{energy},\mathrm{bet}}+\sum_{j=0}^{k_{2}}\mathscr{E}_{k_{1}+1,j,\mathrm{energy},bet}+\frac{C(\varepsilon)}{\varepsilon} \mathcal{E}_{\mathrm{int}}.
\end{align}
\end{theorem}

\noindent It has the same shape as Theorem~\ref{thm:Higher-order contact line elliptic estimate}, with $\chi_{\mathrm{bet}}$ in place of $\chi_{\mathrm{cl}}$ and the interior dissipation on the right, and it is proved at the end of Section~\ref{sec:local_charts}. Its norms enter the interior estimate of Theorem~\ref{thm:interior_space_final} only through nonlinear terms, which is how the three regions are tied together without a loss of constants.

\smallskip\noindent\textbf{The angular energy estimate (Section~\ref{sec:contact_energy}).} We can now run the energy estimate for the mixed fields $\partial_{t}^{k_{1}}\partial_{\theta}^{k_{2}}$ near the contact line, by induction on $k_{2}$; the case $k_{2}=0$ is Theorem~\ref{thm:whole-space-energy}.

\begin{theorem}[Angular energy estimate near the contact line]
\label{thm:contact_line_energy}Suppose $(\partial_{\theta}^{k_{2}}\partial_{t}^{k_{1}}v_{cl},\partial_{\theta}^{k_{2}}\partial_{t}^{k_{1}}\eta,\partial_{t}^{k_{1}}\partial_{\theta}^{k_{2}}p)$ satisfies the equation system \eqref{eq:localized-after-flattening_original_chart} for $k_{1}+k_{2}\leq m+1$. Then, for $\varepsilon$ small enough and with a constant $C(\varepsilon)$ as in Subsection~\ref{subsec:smallness},
\begin{align}
\mathcal{E}_{cl}(t)+\mathcal{E}_{bet}(t)+\int_{0}^{t}(\mathcal{D}_{cl}+\mathcal{D}_{bet})\leq C(\varepsilon)\Big(\mathcal{E}(0)+\mathcal{E}^{\frac{3}{2}}(t)+\int_0^{t}\sqrt{\mathcal{E}}\mathcal{D}\Big).
\end{align}
The proof gives $C(\varepsilon)=C\varepsilon^{-10(m+1)}C^{\prime }(\varepsilon)$, where $C^{\prime }(\varepsilon)$ collects the constants of the cutoff commutators and of Theorem~\ref{thm:top_level}; only its finiteness is used. The interior dissipation $\int_{0}^{t}\mathcal{D}_{\mathrm{int}}$ produced by the cutoff commutators is bounded by the same right-hand side, through Theorem~\ref{thm:interior_space_final}; see \eqref{eq:Dint_E0}.
\end{theorem}

\noindent Each step uses the surface gain of Lemma~\ref{lem:local_threehalf} and the pressure bound of Lemma~\ref{lem:enhance_pressure}, and each produces a factor $\varepsilon^{\iota}$ from the small support of the cutoffs near $\Gamma$. Moreover, to compensate for and absorb the singular $\ep$-dependent constant in Theorem 1.4, we use the key compactness of the linear contributions to create a small parameter:
\begin{align}
\|\p_{\theta}v\|_{L^2_\theta L^{2/\nu}_{r,h}} \lesssim \mathcal{M}(\varepsilon)^{\ell_{1}}\|\p_{\theta}v\|_{H^{1}}+\frac{1}{\mathcal{M}(\varepsilon)^{\ell_{2}}}\|v\|_{H^{1}},
\end{align}
where $\nu=\varepsilon_-/2$ is fixed as in Section~\ref{sec:contact_energy}, $\mathcal{M}(\varepsilon)>0$ is chosen sufficiently small, and $\ell_1,\ell_2>0$. This allows the top-level terms to be absorbed.

One term resists. At the purely angular top level $(k_{1},k_{2})=(0,m+1)$, the kinematic condition contributes the contact-line transport term
\begin{equation*}
\kappa\langle\partial_{t}\eta_{0,m+1},\,v_{\gamma}\partial_{\gamma}\eta_{0,m+1}\rangle_{\Gamma_{0}}, \qquad v_{\gamma}=v\cdot\dot{\gamma},
\end{equation*}
which has no two-dimensional counterpart. At a contact point on the vertical wall in~\cite{GT}, impermeability gives $u_1=0$ and the kinematic condition allows vertical motion. In three dimensions a tangential velocity along the contact curve contributes the transport term above. This term needs $v_{\gamma}\p_{\gamma}\eta_{0,m+1}\in L^{2}(\Gamma_{0})$, that is, $\partial_{\theta}^{m+2}\eta\in L^{2}(\Gamma_{0})$, since $v_{\gamma}\in L^{\infty}(\Gamma_{0})$, while the dissipation gives $\partial_{\theta}^{m+1}\eta\in H^{3/2-\alpha}(\Sigma)$, whose trace is short by exactly $\alpha$.

\smallskip\noindent\textbf{The top angular level (Subsection~\ref{subsec:top_level}).} The following theorem closes that gap. Its hypotheses (H1)--(H4), stated in Subsection~\ref{subsec:top_level}, are the structural properties of the localized system that the previous sections establish.

\begin{theorem}[Top $\protect\theta$-level: $\partial_\protect\theta^{m+1}p\in L^2$ and $\partial_\protect\theta^{m+1}\protect\eta\in H^{3/2}$, without loss]
\label{thm:top_level} Assume (H1)--(H4). There is a constant $C_0\ge1$, depending only on $\zeta$, $\Omega$, $\mu,\sigma,g,\kappa,\beta_0$, $\alpha,\varepsilon_\pm$, the interior exponent $p$, $m$ and the profiles of the cutoffs, but not on $\varepsilon$, and for each $\varepsilon$ there are constants $\delta_0(\varepsilon)\in(0,1]$ and $C(\varepsilon)$, with the following property. At every time $t$ with $\mathcal{E}(t)\le\delta_0(\varepsilon)$,
\begin{align}
\big\|\chi_{\mathrm{bet}+}\,\partial_\theta^{m+1}p\big\|_{L^2(\Omega)}&\le C_0\Big(\big\|\chi_{\mathrm{bet}+}\,\partial_\theta^{m+1}v\big\|_{H^1(\Omega)}+\mathscr{E}_{1,m,\mathrm{energy,cl}}^{1/2}\Big)+C(\varepsilon)\Big(\mathcal{D}_{\mathrm{int}}^{1/2}+\sqrt{\mathcal{E}}\,\mathcal{D}^{1/2}\Big),  \label{tl:pmain} \\
\big\|\chi_{\mathrm{cl},\Sigma}\,\partial_\theta^{m+1}\eta\big\|_{H^{3/2}(\Sigma)}&\le C(\varepsilon)\,\mathcal{D}^{1/2},  \label{tl:main}
\end{align}
all quantities being evaluated at $t$. The pressure bound \eqref{tl:pmain} does not use $\chi_{\mathrm{cl},\Sigma}\partial_\theta^{m+1}\eta(t)\in H^{3/2}(\Sigma)$. In particular, since $\mathscr{E}_{1,m,\mathrm{energy,cl}}\le\mathcal{D}$,
\begin{equation*}
\big\|\chi_{\mathrm{bet}+}\,\partial_\theta^{m+1}p(t)\big\|^2_{L^2(\Omega)}+\big\|\chi_{cl,\Sigma}\,\partial_\theta^{m+1}\eta(t)\big\|^2_{H^{3/2}(\Sigma)}\le C(\varepsilon)\,\mathcal{D}(t).
\end{equation*}
\end{theorem}

\begin{corollary}
\label{cor:top_level} Under the assumptions of Theorem~\ref{thm:top_level}, the following hold.

\begin{enumerate}[label=(\roman*)]

\item For every $t$ with $\mathcal{E}(t)\le\delta_0(\varepsilon)$,
\begin{equation*}
\|\partial_\theta^{m+2}\eta(t)\|_{L^2(\Gamma_0)}+\|\partial_\theta^{m+1}\eta(t)\|_{H^1(\Gamma_0)}\le C(\varepsilon)\,\mathcal{D}(t)^{1/2}.
\end{equation*}

\item If $\sup_{[s,t]}\mathcal{E}\le\delta_0(\varepsilon)$, then
\begin{equation*}
\int_s^t\Big(\|\chi_{ \mathrm{cl},\Sigma}\partial_\theta^{m+1}\eta\|^2_{H^{3/2}(\Sigma)}+\|\partial_\theta^{m+2}\eta\|^2_{L^2(\Gamma_0)}\Big)\le C(\varepsilon)\int_s^t\mathcal{D}.
\end{equation*}

\item The contact-line transport term in \eqref{eq:local-domain-energy-1-2} at $(k_1,k_2)=(0,m+1)$ satisfies
\begin{equation*}
\big|\kappa\langle\partial_t\eta_{0,m+1},\,v_\gamma\,\partial_\gamma\eta_{0,m+1}\rangle_{\Gamma_0}\big|\le C(\varepsilon)\sqrt{\mathcal{E}}\,\mathcal{D},
\end{equation*}
where $v_\gamma=v\cdot\dot\gamma$. The same $C(\varepsilon)\sqrt{\mathcal{E}}\,\mathcal{D}$ bound holds for the contribution of the top transport component of $[\widetilde F^4]_{\Gamma_0}$ of Step~7 of the proof of Theorem~\ref{thm:contact_line_energy} when paired with the contact test in the energy identity. Its $L^2(\Gamma_0)$ norm is bounded by $C(\varepsilon)\sqrt{\mathcal{E}\mathcal{D}}$, as shown below.
\end{enumerate}
\end{corollary}

\noindent Part~(iii) of the corollary is the missing bound for the transport term above, and the induction of Theorem~\ref{thm:contact_line_energy} closes. The reason no $\alpha$ is lost is that the top level is purely angular: the time-derivative term, the only source of the loss, is $\partial_{t}\partial_{\theta}^{m+1}v=\partial_{\theta}(\partial_{t}\partial_{\theta}^{m}v)$, an angular derivative of a quantity the dissipation controls, so the derivative can be moved onto the test function instead of being paid for by interpolation in time. This is the gain announced in the abstract, and it has no analogue for $\partial_{t}$.

The proof has two halves. The pressure is obtained as in Lemma~\ref{lem:enhance_pressure}, but on a domain whose size is fixed by $\zeta$, with the localized pressure extended by zero, so that the constant in the Ne\caronaccent{c}as inequality does not degenerate as the support shrinks. For the surface we first test the differentiated stress against $H^{1}$ fields, which places the normal stress in $H^{-1/2}(\Sigma)$, and then invert the capillary operator with its contact-line condition, which is Lemma~\ref{lem:tl_ell}. The one delicate point is the top-order term $\mu\func{div}_{\mathcal{A}}(\mathbb{D}_{\mathfrak{M}}v)$, which carries two derivatives of $\partial_{\theta}^{m+1}\bar{\eta}$. Because the extension is slice-wise, its $\theta\theta$ derivatives are \emph{not} controlled by $\|\partial_{\theta}^{m+1}\eta\|_{H^{3/2}}$: for $G=e^{ik\theta}\varphi(r)$ one has $\|G\|_{H^{3/2}}\simeq k^{3/2}$, while $\|\partial_\theta^2P_{2-D}EG\|_{L^2}\simeq k^2$ near $\Gamma$. We therefore never differentiate $\nabla\partial_{\theta}^{m+1}\bar{\eta}$: the whole term is integrated by parts in all variables, its face terms on $\Sigma$ and $\Sigma_{s}$ cancel against the matching boundary forcing, and the remaining bulk term needs only $\partial_{\theta}^{m+1}\eta\in H^{3/2-\alpha}$. The same device repairs the top-level forcing bound \eqref{eq:F01_dual}, which fails in slicing norms but holds against $H^{1}(\Omega)$ test functions; see the proof of Lemma~\ref{thm:contact_line_forcing_estimate}.

\smallskip\noindent\textbf{Closing.} Combining Theorem~\ref{thm:contact_line_energy} with Theorem~\ref{thm:interior_space_final} and choosing $\varepsilon$ small, the induction closes and yields
\begin{equation*}
\mathcal{E}(t)+\int_{0}^{t}\mathcal{D}\lesssim\mathcal{E}(0)+\int_{0}^{t}\sqrt{\mathcal{E}}\mathcal{D}+\mathcal{E}^{3/2}(t).
\end{equation*}
Under the smallness assumption the last two terms are absorbed. Since the system is autonomous, the same estimate started at time $s$ gives $\mathcal{E}(t)+\int_{s}^{t}\mathcal{D}\lesssim\mathcal{E}(s)$, and since $\mathcal{E}\lesssim\mathcal{D}$ term by term, the Main Theorem follows with exponential decay. The argument is completed at the end of Subsection~\ref{subsec:cl_local_energy}.

\section{Notation}

\label{sec:notation_guide}

The symbols used throughout are collected here, grouped by role, each with the place where it is defined. Nothing is defined here; the references point to the definition in the text.

\begingroup
\setlength{\parskip}{0pt}
\renewcommand{\addcontentsline}[3]{}

\subsection*{Geometry}

\begin{description}[leftmargin=5.6em,style=sameline,font=\normalfont,itemsep=1pt]

\item[$D$, $\partial D$] the base domain of the equilibrium surface and its boundary curve; \eqref{eq:moving-domain-full-system}.

\item[$\protect\zeta$] the equilibrium height over $D$, a solution of the Young--Laplace equation; \eqref{eq:young-laplace}--\eqref{eq:young}.

\item[$\Omega$, $\Sigma$, $\Sigma_{s}$] the equilibrium domain, the free surface and the side wall; \eqref{eq:moving-domain-full-system} and Figure~\ref{fig:contact_line}.

\item[$\Gamma$] the contact line, where $\Sigma$ meets $\Sigma_{s}$; Figure~\ref{fig:contact_line}.

\item[$\Gamma_{0}$] the contact line in the chart, $\{r=0\}$; \eqref{eq:contact_chart}.

\item[$\protect\theta_{\mathrm{eq}}$] the equilibrium contact angle, measured through the fluid (the opening angle of the fluid at $\Gamma$), constant along $\Gamma$; \eqref{eq:young} and \eqref{eq:s1-cot}.

\item[$\protect\gamma(\protect\theta)$, $\dot\protect\gamma$] the arclength parametrization of $\partial D$ and its unit tangent; Section~\ref{sec:flattening}.

\item[$\Phi$, $(\protect\theta,r,h)$] the contact-line chart and its variables: $\theta$ arclength along $\Gamma$, $(r,h)$ transverse; \eqref{eq:contact_chart}.

\item[$\protect\nu$] the outward unit normal of $\partial D$; the remark after \eqref{eq:s1-24}.

\item[$\mathcal{N}$, $n$, $n_{0}$] the unnormalized and unit normals of the free surface and of the wall; \eqref{eq:transformed-normal}--\eqref{eq:N-equality}.

\item[$H_{0}=-r\cot\protect\theta_{\mathrm{eq}}$] the equilibrium surface in chart variables; \eqref{eq:equilibrium-H}.

\item[$Q$, $\mathcal{A}_{\Phi}$] the rotation and the chart matrix, with $Q^{T}\mathcal{A}_{\Phi}=I$ on $\Gamma$; \eqref{eq:Phi_exact_decomposition}--\eqref{eq:exact_QTA_identity_boxed}.
\end{description}

\subsection*{Unknowns and transformed quantities}

\begin{description}[leftmargin=5.6em,style=sameline,font=\normalfont,itemsep=1pt]

\item[$u$, $P$, $\protect\eta$] the velocity, pressure and surface perturbation on the moving domain; \eqref{eq:moving-domain-full-system}.

\item[$v$, $p$] their pullbacks to the equilibrium domain; \eqref{eq:pullback-unknowns}--\eqref{eq:moving-domain-elliptic-system}.

\item[$\mathfrak{h}_{\protect\eta}$] the surface perturbation in the chart; \eqref{eq:local-eta-definition}--\eqref{eq:local-eta-H-relation}.

\item[$\bar\protect\eta$] the extension of $\eta$ into the bulk, slice-wise near $\Gamma$; \eqref{eq:eta_bar_def}.

\item[$\Phi_{\protect\eta}$, $\mathcal{A}_{\protect\eta}$, $J$, $K$] the flattening map, its inverse-gradient matrix, and the Jacobian factors; \eqref{eq:JKA-definition}--\eqref{eq:A-explicit}.

\item[$S_{\mathcal{A}}(p,v)$, $\mathbb{D}_{\mathcal{A}}$] the geometric stress and symmetric gradient; \eqref{eq:geometric-stress}.

\item[$\mathcal{L}_{\Sigma,0}$, $a_{\protect\zeta}$] the linearized capillary operator and its coefficient; \eqref{def:linearized_operator}.

\item[$\mathcal{B}_{\mathrm{cl},0}$] the linearized contact-line operator; \eqref{def:contact_operator}.

\item[$\hat{\mathcal{W}}$] the contact-line friction, less its linear part; \eqref{eq:W_hat}.

\item[$\mathcal{R}(\protect\nabla\protect\eta)$] the quadratic remainder of the mean curvature; \eqref{eq:F05}.

\item[$F^{i,j}_{0}$, $i=1,\dots,7$] the forcing terms at time level $j$: bulk, divergence, dynamic, kinematic, slip, wall and contact-line; \eqref{eq:F01j}--\eqref{eq:F07j}.

\item[$\tilde F^{i,k_{1},k_{2}}_{\mathrm{bet}+}$] their localized, flattened counterparts at a mixed level; \eqref{eq:final-localized-flattened-system-j}.

\item[$\mathfrak{M}$, $\widetilde F^{i}$] the top-level coefficient matrix and forcing; Subsection~\ref{subsec:top_level}.
\end{description}

\subsection*{Cutoff functions and regions}

\begin{description}[leftmargin=5.6em,style=sameline,font=\normalfont,itemsep=1pt]

\item[$\protect\chi_{\mathrm{cl}}$, $\protect\chi_{\mathrm{bet}}$, $\protect\chi_{\mathrm{int}}$] the cutoffs at the contact line, on the intermediate region and in the interior, on $\Sigma$; Subsection~\ref{subsec:cutoff}.

\item[$\protect\chi_{\mathrm{bet}+}$, $\protect\chi_{\mathrm{cl}+}$, $\protect\chi_{\mathrm{int}+}$] the sums $\chi_{\mathrm{cl}}+\chi_{\mathrm{bet}}$ and the corresponding enlarged cutoffs; Subsection~\ref{subsec:cutoff}.

\item[$\protect\chi_{j,\mathrm{int}}$] the nested interior cutoffs, one for each derivative level; Subsection~\ref{subsec:cutoff}.

\item[$\protect\chi_{\mathrm{bot}}$] the bottom cutoff, where $\mathcal{A}_{\eta}=I$; Subsection~\ref{subsec:cutoff_extension}.

\item[$h_{0}$] the level below which the interior cutoffs vanish; Subsection~\ref{subsec:cutoff_extension}.

\item[$\protect\varepsilon$] the width of the supports of the cutoffs near $\Gamma$; \eqref{eq:cutoff_eps}, and the gains \eqref{eq:QA_small}--\eqref{eq:QA_theta_small}.

\item[$\Omega_{\mathrm{cl}}$, $\Omega_{\mathrm{bet}}$, $\Omega_{\mathrm{bet}+}$, $\Omega_{\mathrm{int}}$] the artificial domains, of a size fixed by $\zeta$, on which the elliptic and Ne\caronaccent{c}as constants are taken; Subsection~\ref{subsec:artificial}.
\end{description}

\subsection*{Derivatives, norms and spaces}

\begin{description}[leftmargin=5.6em,style=sameline,font=\normalfont,itemsep=1pt]

\item[$\partial_{t}$, $\partial_{\protect\theta}$] the two vector fields: time, and arclength along $\Gamma$ in the chart; Subsection~\ref{subsec:cl_local_energy}.

\item[$v_{k_{1},k_{2}}$, etc.] the localized $\partial_{t}^{k_{1}}\partial_{\theta}^{k_{2}}$ of a quantity; \eqref{eq:localized-after-flattening_original_chart}.

\item[$X_{\protect\theta}Y_{rh}$] the slicing convention: the first space acts in $\theta$, the second in the transverse variables; Subsection~\ref{subsec:notation}.

\item[{$[\,\cdot\,]_{\Gamma_{0}}$}] the $L^{2}$ norm, or pairing, on the contact line; Subsection~\ref{subsec:notation}.

\item[$A\lesssim B$] inequality up to a constant independent of the solution and of $\varepsilon$; \eqref{eq:lesssim}.

\item[$E$, $P_{2-D}$] the half-space extender and the two-dimensional Poisson extension; \eqref{eq:half-space-extension} and \eqref{eq:contact-line-extension}.
\end{description}

\subsection*{Parameters}

\begin{description}[leftmargin=5.6em,style=sameline,font=\normalfont,itemsep=1pt]

\item[$m$] the top derivative level; $m\geq3$, Subsection~\ref{subsec:constants}.

\item[$\protect\mu$, $\protect\sigma$, $g$, $\protect\kappa$, $\protect\beta_{0}$] viscosity, surface tension, gravity, contact-line friction, Navier slip; \eqref{eq:moving-domain-full-system}.

\item[$\protect\varepsilon_{\max}$] the corner exponent, $\min\{1,\pi/\theta_{\mathrm{eq}}-1\}$; Subsection~\ref{subsec:constants}.

\item[$\protect\varepsilon_{\pm}$, $\protect\alpha$, $\protect\delta$] the regularity parameters, with $0<\alpha<\varepsilon_{-}<\varepsilon_{+}<\varepsilon_{\max}$ and $\delta=\varepsilon_{-}$; \eqref{eq:parameters}.

\item[$q_{\protect\varepsilon}$, $q_{\pm}$, $q_{\protect\delta}$] the slice exponents $2/(2-\varepsilon)$, close to $1$; \eqref{eq:q_eps} and \eqref{eq:q_pm}.

\item[$p$] the interior and intermediate exponent, close to $2$; Subsection~\ref{subsec:constants}.

\item[$\protect\iota$] the small power of $\varepsilon$ gained near the contact line; the proof of Theorem~\ref{thm:contact_line_energy}.

\item[$\protect\gamma$, $\protect\delta_{0}$] the smallness thresholds for $\mathcal{E}$; Proposition~\ref{prop:B_d} and Theorem~\ref{thm:top_level}.
\end{description}

\subsection*{Energies and dissipations}

\begin{description}[leftmargin=5.6em,style=sameline,font=\normalfont,itemsep=1pt]

\item[$\mathcal{E}$, $\mathcal{D}$] the total energy and dissipation, the sums of the pieces below; Subsection~\ref{subsec:energy_def}.

\item[$\mathcal{E}_{||}$, $\mathcal{D}_{||}$] the temporal pieces, for the fields $\partial_{t}^{j}$; Subsection~\ref{subsec:energy_def}.

\item[$\mathcal{E}_{\mathrm{enhanced}}$, $\mathcal{D}_{\mathrm{enhanced}}$] the enhancements: $H^{3/2-\alpha}$ for the surface and $L^{2}$ for the pressure; Theorem~\ref{thm:whole_space_estimate_final}.

\item[$\mathcal{E}_{\mathrm{int}}$, $\mathcal{D}_{\mathrm{int}}$] the interior elliptic pieces; \eqref{def:interior_energy}--\eqref{def:interior_dissipation}.

\item[$\mathcal{E}_{\mathrm{cl}}$, $\mathcal{E}_{\mathrm{bet}}$, and $\mathcal{D}_{\mathrm{cl}}$, $\mathcal{D}_{\mathrm{bet}}$] the pieces at the contact line and on the intermediate region; \eqref{def:split_energy}--\eqref{def:energy_cl}.

\item[$\mathscr{E}_{k_{1},k_{2},\bullet}$] the piece at the single level $(k_{1},k_{2})$, with $\bullet$ one of $\mathrm{energy}$, $\mathrm{el}$, $\mathrm{cl}$, $\mathrm{bet}$, $\mathrm{int}$; Subsection~\ref{subsec:energy_def}.
\end{description}

\endgroup

\section{Flattening map and local charts}

\label{sec:flattening}

\subsection{Implicit function argument}

We first record an identity for the contact angle $\theta_{\mathrm{eq}}$ of \eqref{eq:young}:
\begin{equation}  \label{eq:s1-cot} \cot\theta_{\mathrm{eq}} = \frac{\cos\theta_{\mathrm{eq}}}{\sin\theta_{\mathrm{eq}}} = \frac{ -\nabla\zeta\circ\gamma\cdot\dot{\gamma}^{\perp} }{ \sqrt{ 1+\left|\nabla\zeta\circ\gamma\cdot\dot{\gamma}\right|^{2} } }.
\end{equation}
Here $\dot\gamma^{\perp}$ is the inward unit normal of $\partial D$ (see the remark after \eqref{eq:s1-24}), so $-\nabla\zeta\cdot\dot\gamma^{\perp}=\partial_{\nu}\zeta$ is the outward normal derivative of $\zeta$. The angle $\theta_{\mathrm{eq}}$, measured through the fluid, satisfies $\cos\theta_{\mathrm{eq}}=\partial_{\nu}\zeta/\sqrt{1+|\nabla\zeta|^{2}}$ on $\partial D$, and $|\nabla\zeta|^{2}=|\partial_{\nu}\zeta|^{2}+|\nabla\zeta\cdot\dot\gamma|^{2}$ there; this gives \eqref{eq:s1-cot}. By Young's law, $\theta_{\mathrm{eq}}$ is constant along $\Gamma$.

Define the smooth map $f:\mathbb{T}_{L}\times\mathbb{R}^{2}\times\mathbb{T}_{L} \to\mathbb{R}$ via
\begin{equation}
\begin{aligned}
&f(\theta,r,h,\tau) ={}  \Bigg( \gamma(\tau)+r\dot{\gamma}^{\perp}(\tau), \zeta\bigl( \gamma(\tau)+r\dot{\gamma}^{\perp}(\tau) \bigr) + \sqrt{ 1+ \left| \nabla\zeta\circ\gamma \cdot\dot{\gamma}(\theta) \right|^{2} } \bigl[h+\cot\theta_{\mathrm{eq}}\,r\bigr] \Bigg) \cdot T(\theta) \\
&\quad - \bigl( \gamma(\theta), \zeta\circ\gamma(\theta) \bigr)\cdot T(\theta).  \label{eq:s1-7}
\end{aligned}
\end{equation}
Note that $f(\theta,0,0,\theta)=0$ for every $\theta\in\mathbb{T}_{L}$. We compute
\begin{equation}
\begin{aligned}
&\partial_{\theta}f(\theta,0,0,\theta) ={}  \Big[ \bigl( \gamma(\theta), \zeta\circ\gamma(\theta) \bigr) - \bigl( \gamma(\theta), \zeta\circ\gamma(\theta) \bigr) \Big]\cdot T^{\prime }(\theta) \\
&\quad - \bigl( \dot{\gamma}(\theta), \nabla\zeta\circ\gamma(\theta) \cdot\dot{\gamma}(\theta) \bigr)\cdot T(\theta) ={}  - \sqrt{ 1+ \left| \nabla\zeta\circ\gamma \cdot\dot{\gamma}(\theta) \right|^{2} }.  \label{eq:s1-8}
\end{aligned}
\end{equation}
Using \eqref{eq:s1-cot}, we also have
\begin{equation}
\partial_{r}f(\theta,0,0,\theta) ={} \Bigg( \dot{\gamma}^{\perp}(\theta), \nabla\zeta\circ\gamma(\theta) \cdot\dot{\gamma}^{\perp}(\theta) + \cot\theta_{\mathrm{eq}} \sqrt{ 1+ \left| \nabla\zeta\circ\gamma \cdot\dot{\gamma}(\theta) \right|^{2} } \Bigg)\cdot T(\theta) ={} \bigl( \dot{\gamma}^{\perp}(\theta),0 \bigr)\cdot T(\theta) =0. \label{eq:s1-9}
\end{equation}
Moreover,
\begin{equation}
\partial_{h}f(\theta,0,0,\theta) = \left( 0, \sqrt{ 1+ \left| \nabla\zeta\circ\gamma \cdot\dot{\gamma}(\theta) \right|^{2} } \right)\cdot T(\theta)    = \nabla\zeta\circ\gamma(\theta) \cdot\dot{\gamma}(\theta), \label{eq:s1-10}
\end{equation}
and, most importantly,
\begin{equation}
\partial_{\tau}f(\theta,0,0,\theta) = \bigl( \dot{\gamma}(\theta), \nabla\zeta\circ\gamma(\theta) \cdot\dot{\gamma}(\theta) \bigr)\cdot T(\theta)    = \sqrt{ 1+ \left| \nabla\zeta\circ\gamma \cdot\dot{\gamma}(\theta) \right|^{2} } \geq 1. \label{eq:s1-11}
\end{equation}

Using the implicit function theorem at each point $(\rho,0,0,\rho), \; \rho\in\mathbb{T}_{L},$ we can construct open sets $W_{\rho}\subset\mathbb{T}_{L}\times\mathbb{R}^{2}, \; V_{\rho}\subset\mathbb{T}_{L},$ and a smooth function $\tau_{\rho}:W_{\rho}\to V_{\rho}$ such that:

\begin{enumerate}
\item $(\rho,0,0)\in W_{\rho}, \; \rho\in V_{\rho};$ \item $\tau_{\rho}(\rho,0,0)=\rho;$ \item $f( \theta,r,h,\tau_{\rho}(\theta,r,h) )=0$ for every $(\theta,r,h)\in W_{\rho}$, and if $f(\theta,r,h,\sigma)=0$ for $(\theta,r,h,\sigma)\in W_{\rho}\times V_{\rho},  $  then $\sigma=\tau_{\rho}(\theta,r,h).$
\end{enumerate}

Then $\{W_{\rho}\}_{\rho\in\mathbb{T}_{L}}$ is an open cover of $\mathbb{T}_{L}\times\{0\}\times\{0\},$ which is compact. Hence, we can choose a finite subcover $\{W_{\rho_{1}},\dots,W_{\rho_{n}}\}.$ We construct a single branch on a uniform collar as follows. Set $G(\theta,r,h,\gamma):=f(\theta,r,h,\theta+\gamma),$ where the last argument is read modulo $L$. Since $G(\theta,0,0,0)=0$ and $\partial_\gamma G(\theta,0,0,0)\ge1$, compactness gives $0<d<L/4$ and $a>0$ such that $\partial_\gamma G(\theta,r,h,\gamma)\ge\frac12 \; (|r|+|h|<a,\ |\gamma|\le d).$ After decreasing $a$, we also have $|G(\theta,r,h,0)|<d/4$. Consequently $G(\theta,r,h,-d)<0<G(\theta,r,h,d).$ For each $(\theta,r,h)$ in this fixed collar there is therefore exactly one $\gamma(\theta,r,h)\in(-d,d)$ for which $G=0$. The implicit function theorem and uniqueness show that $\gamma$ is smooth and periodic in $\theta$. Define $\tau(\theta,r,h):=\theta+\gamma(\theta,r,h)\pmod L.$ Set $U:=\mathbb T_L\times\{(r,h):|r|+|h|<a\}.$ The branch $\tau:U\to\mathbb T_L$ constructed above is smooth and satisfies $\tau(\theta,0,0)=\theta,\; f(\theta,r,h,\tau(\theta,r,h))=0 \;\text{for every }(\theta,r,h)\in U.$ Differentiating the latter identity, we obtain
\begin{equation}
\bigl( \partial_{\theta}\tau, \partial_{r}\tau, \partial_{h}\tau \bigr)(\rho,0,0)  = - \frac{ \bigl( \partial_{\theta}f, \partial_{r}f, \partial_{h}f \bigr)(\rho,0,0,\rho) }{ \partial_{\tau}f(\rho,0,0,\rho) } = \left( 1, 0, - \frac{ \nabla\zeta\circ\gamma(\rho) \cdot\dot{\gamma}(\rho) }{ \sqrt{ 1+ \left| \nabla\zeta\circ\gamma(\rho) \cdot\dot{\gamma}(\rho) \right|^{2} } } \right).  \label{eq:s1-15}
\end{equation}

\subsection{The map}

Restricting the size of $U$, if necessary, we may assume that $U=\mathbb{T}_{L}\times\widetilde{U}$ for some open set $\widetilde{U}\subset\mathbb{R}^{2}.$ This will be convenient in the next part of the construction. Define ${\Phi}:U\to\mathbb{R}^{3}$ via
\begin{equation}
\begin{aligned}
&{\Phi}(\theta,r,h) = \Bigg(  \gamma\bigl(\tau(\theta,r,h)\bigr) + r\dot{\gamma}^{\perp} \bigl(\tau(\theta,r,h)\bigr), \zeta\left( \gamma\bigl(\tau(\theta,r,h)\bigr) + r\dot{\gamma}^{\perp} \bigl(\tau(\theta,r,h)\bigr) \right) \\
&\quad + \sqrt{ 1+ \left| \nabla\zeta\circ\gamma \cdot\dot{\gamma}(\theta) \right|^{2} } \bigl[h+\cot\theta_{\mathrm{eq}}\,r\bigr] \Bigg). \label{eq:s1-16}
\end{aligned}
\end{equation}
Note that
\begin{gather}
{\Phi} \Bigl( \bigl\{ (\theta,r,h)\in U \mid \sin\theta_{\mathrm{eq}}\,h+\cos\theta_{\mathrm{eq}}\,r=0,\ r>0 \bigr\} \Bigr) \subseteq \bigl\{ x\in\mathbb{R}^{3} \mid x_{3}=\zeta(x^{\prime }) \bigr\},  \label{eq:s1-17} \\
{\Phi} \Bigl( \bigl\{ (\theta,r,h)\in U \mid \sin\theta_{\mathrm{eq}}\,h+\cos\theta_{\mathrm{eq}}\,r>0,\ r>0 \bigr\} \Bigr) \subseteq \bigl\{ x\in\mathbb{R}^{3} \mid x_{3}>\zeta(x^{\prime }) \bigr\},  \label{eq:s1-18}
\end{gather}
and
\begin{align}
{\Phi} \Bigl( \bigl\{ (\theta,r,h)\in U \mid \sin\theta_{\mathrm{eq}}\,h+\cos\theta_{\mathrm{eq}}\,r<0,\ r>0 \bigr\} \Bigr) \subseteq \bigl\{ x\in\mathbb{R}^{3} \mid x_{3}<\zeta(x^{\prime }) \bigr\}.  \label{eq:s1-19}
\end{align}
Together, these identities mean that the top of the wedge maps into the free surface, while the volumes above and below the top of the wedge map to the corresponding volumes above and below the free surface.

Similarly, ${\Phi} ( \{ (\theta,r,h)\in U \mid r=0 \} ) \subseteq \partial D\times\mathbb{R},$ \begin{align} {\Phi} \Bigl( \bigl\{ (\theta,r,h)\in U \mid r=0,\ \sin\theta_{\mathrm{eq}}\,h+\cos\theta_{\mathrm{eq}}\,r>0 \bigr\} \Bigr) \subseteq \bigl\{ x\in\mathbb{R}^{3} \mid x^{\prime }\in\partial D,\ x_{3}>\zeta(x^{\prime }) \bigr\},  \label{eq:s1-21}
\end{align}
and
\begin{align}
{\Phi} \Bigl( \bigl\{ (\theta,r,h)\in U \mid r=0,\ \sin\theta_{\mathrm{eq}}\,h+\cos\theta_{\mathrm{eq}}\,r<0 \bigr\} \Bigr) \subseteq \bigl\{ x\in\mathbb{R}^{3} \mid x^{\prime }\in\partial D,\ x_{3}<\zeta(x^{\prime }) \bigr\}.  \label{eq:s1-22}
\end{align}
Together, these identities mean that the appropriate parts of the wedge side wall map into the corresponding parts of the container wall. Finally, ${\Phi}(\theta,0,0) = ( \gamma(\theta), \zeta\circ\gamma(\theta) ),$ which shows that the line $\mathbb{T}_{L}\times\{0\}\times\{0\}$ maps to the contact line.

Next, we use the previously established identities to compute
\begin{equation}
\partial_{\theta}{\Phi}(\theta,0,0)  = \Bigl( \dot{\gamma}(\theta) \partial_{\theta}\tau(\theta,0,0), \nabla\zeta\circ\gamma(\theta) \cdot\dot{\gamma}(\theta) \partial_{\theta}\tau(\theta,0,0) \Bigr) = \bigl( \dot{\gamma}(\theta), \nabla\zeta\circ\gamma(\theta) \cdot\dot{\gamma}(\theta) \bigr)   = \sqrt{ 1+ \left| \nabla\zeta\circ\gamma \cdot\dot{\gamma}(\theta) \right|^{2} } \,T(\theta).  \label{eq:s1-23}
\end{equation}
Also,
\begin{equation}
\begin{aligned}
&\partial_{r}{\Phi}(\theta,0,0) ={}  \Bigg( \dot{\gamma}(\theta) \partial_{r}\tau(\theta,0,0) + \dot{\gamma}^{\perp}(\theta), \nabla\zeta\circ\gamma(\theta) \cdot\dot{\gamma}(\theta) \partial_{r}\tau(\theta,0,0) \\
&\quad + \nabla\zeta\circ\gamma(\theta) \cdot\dot{\gamma}^{\perp}(\theta) + \cot\theta_{\mathrm{eq}} \sqrt{ 1+ \left| \nabla\zeta\circ\gamma \cdot\dot{\gamma}(\theta) \right|^{2} } \Bigg) ={}  \bigl( \dot{\gamma}^{\perp}(\theta),0 \bigr) = \nu(\theta). \label{eq:s1-24}
\end{aligned}
\end{equation}
Here $\nu(\theta)=(\dot{\gamma}^{\perp}(\theta),0)$ points into $D$; it is the inward unit normal of $\partial D$. In \eqref{eq:moving-domain-full-system} and \eqref{def:contact_operator}, $\nu$ denotes the outward unit normal $-\nu(\theta)$; see also the remark after \eqref{eq:local-system-h-eta}.

Furthermore,
\begin{equation}
\begin{aligned}
&\partial_{h}{\Phi}(\theta,0,0) ={}  \Bigg( \dot{\gamma}(\theta) \partial_{h}\tau(\theta,0,0), \nabla\zeta\circ\gamma(\theta) \cdot\dot{\gamma}(\theta) \partial_{h}\tau(\theta,0,0) + \sqrt{ 1+ \left| \nabla\zeta\circ\gamma \cdot\dot{\gamma}(\theta) \right|^{2} } \Bigg) \\
&\quad ={}  \frac{1}{ \sqrt{ 1+ \left| \nabla\zeta\circ\gamma \cdot\dot{\gamma}(\theta) \right|^{2} } } \Bigl( - \dot{\gamma}(\theta) \, \nabla\zeta\circ\gamma(\theta) \cdot\dot{\gamma}(\theta), 1 \Bigr)   ={} T(\theta)\wedge\nu(\theta).  \label{eq:s1-25}
\end{aligned}
\end{equation}

From these identities, it is clear that $\{ \partial_{\theta}{\Phi}, \partial_{r}{\Phi}, \partial_{h}{\Phi} \}$ forms an orthogonal system at every point $(\theta,0,0)$. There is, however, a problem with the normalization of $\partial_{\theta}\overline{\Phi}$, which can be dealt with by reparameterizing the $\theta$-variable, since
\begin{equation}
\sqrt{ 1+ \left| \nabla\zeta\circ\gamma \cdot\dot{\gamma}(\theta) \right|^{2} } \geq 1>0.  \label{eq:s1-26}
\end{equation}
Indeed, set
\begin{equation}
M = \int_{0}^{L} \sqrt{ 1+ \left| \nabla\zeta\circ\gamma \cdot\dot{\gamma}(\vartheta) \right|^{2} } \,d\vartheta \geq L>0,  \label{eq:s1-27}
\end{equation}
and define $\chi:\mathbb{R}\to\mathbb{R}$ via
\begin{equation}
\chi(R) = \int_{0}^{R} \sqrt{ 1+ \left| \nabla\zeta\circ\gamma \cdot\dot{\gamma}(\vartheta) \right|^{2} } \,d\vartheta.  \label{eq:s1-28}
\end{equation}
Here, we regard the map
\begin{equation*}
\mathbb{R}\ni\vartheta \longmapsto \sqrt{ 1+ \left| \nabla\zeta\circ\gamma \cdot\dot{\gamma}(\vartheta) \right|^{2} } \in\mathbb{R}
\end{equation*}
as $L$-periodic. Then $\chi$ is smooth and $\chi^{\prime }(R) = \sqrt{ 1+ | \nabla\zeta\circ\gamma \cdot\dot{\gamma}(R) |^{2} },$ which means that $\chi$ is strictly increasing.

Moreover, due to the $L$-periodicity of $\gamma$, if $m\in\mathbb{Z}$ and $0\leq\rho<L$, then $\chi(mL+\rho) = \chi(\rho)+mM.$ We therefore deduce that $\chi$ is a diffeomorphism from $\mathbb{R}$ onto itself. Moreover, $\chi$ naturally descends to a smooth diffeomorphism $\chi:\mathbb{T}_{L}\to\mathbb{T}_{M}.$ Its inverse $\chi^{-1}:\mathbb{R}\to\mathbb{R}$ is also a diffeomorphism and descends to a smooth diffeomorphism $\chi^{-1}:\mathbb{T}_{M}\to\mathbb{T}_{L}.$ We denote the map defined above by $\bar{\Phi}$ and then define $V=\mathbb{T}_{M}\times\widetilde{U},$ and set $\Phi:V\to\mathbb{R}^{3}$ via
\begin{equation}
\Phi(\theta,r,h) = \overline{\Phi} \bigl( \chi^{-1}(\theta),r,h \bigr). \label{eq:contact_chart}
\end{equation}
Then
\begin{equation}
\partial_{\theta}\Phi(\theta,0,0)  = \partial_{\theta}\overline{\Phi} \bigl( \chi^{-1}(\theta),0,0 \bigr) (\chi^{-1})^{\prime }(\theta) = \partial_{\theta}\overline{\Phi} \bigl( \chi^{-1}(\theta),0,0 \bigr) \frac{1}{ \sqrt{ 1+ \left| \nabla\zeta\circ\gamma \cdot \dot{\gamma}\bigl(\chi^{-1}(\theta)\bigr) \right|^{2} } }   = T\circ\chi^{-1}(\theta).  \label{eq:contact_chart_tangent}
\end{equation}
Similarly,
\begin{equation}
\partial_{r}\Phi(\theta,0,0)  = \partial_{r}\overline{\Phi} \bigl( \chi^{-1}(\theta),0,0 \bigr) = \nu\circ\chi^{-1}(\theta), \qquad \partial_{h}\Phi(\theta,0,0)  = \partial_{h}\overline{\Phi} \bigl( \chi^{-1}(\theta),0,0 \bigr) = (T\wedge\nu)\circ\chi^{-1}(\theta). \label{eq:s1-32}
\end{equation}
We therefore deduce that $\{ \partial_{\theta}\Phi(\theta,0,0), \partial_{r}\Phi(\theta,0,0), \partial_{h}\Phi(\theta,0,0) \}$ is an orthonormal basis. Hence, $D\Phi(\theta,0,0)\in O(3).$ Additionally, reparameterizing in the $\theta$-variable as above does not change any of the mapping properties \eqref{eq:s1-17}--\eqref{eq:s1-22}. These properties continue to hold with $\Phi$ and $V$ in place of $\overline{\Phi}$ and $U$.

It remains to show that $\Phi$ is a diffeomorphism on a possibly smaller set. We know that $D\Phi(\theta,0,0)\in O(3)$ for every $(\theta,0,0)\in V$. By continuity, we may shrink $\widetilde{U}$, if necessary, so that $\Phi$ is a local diffeomorphism around every point of $\mathbb{T}_{M}\times\{0\}\times\{0\}.$ Let $\Gamma$ denote the contact curve and set $\Gamma_{0} = \mathbb{T}_{M}\times\{0\}\times\{0\} = \Phi^{-1}(\Gamma) \subset V.$ By construction, the restriction $\left.\Phi\right|_{\Gamma_{0}}:\Gamma_{0}\to\Gamma$ is a diffeomorphism.

We claim that there exists an open set $\widetilde{W}\subset\widetilde{U}$ such that $\Phi$ is a diffeomorphism onto its image when restricted to $W = \mathbb{T}_{M}\times\widetilde{W} \subseteq V.$ In light of the local diffeomorphism condition, it suffices to show that there exists such an open set on which $\Phi$ is injective. Suppose, by way of contradiction, that no such open set exists. Then, for every $n\in\mathbb{N}$, we may choose $\widetilde{W}_{n} = B(0,2^{-n})\cap\widetilde{U}$ and points $x_{n},y_{n} \in \mathbb{T}_{M}\times\widetilde{W}_{n}$ such that $\Phi(x_{n})=\Phi(y_{n}), \; x_{n}\neq y_{n}.$ By compactness, after extracting a subsequence, we may assume that $x_{n}\to x\in\Gamma_{0}, \; y_{n}\to y\in\Gamma_{0}$ as $n\to\infty$. Then $\Phi(x) = \lim_{n\to\infty}\Phi(x_{n}) = \lim_{n\to\infty}\Phi(y_{n}) = \Phi(y).$ Since $\Phi$ is injective on $\Gamma_{0}$, it follows that $x=y.$ However, for sufficiently large $n$, this contradicts the fact that $\Phi$ is injective when restricted to an open neighborhood of $x=y$, since both $x_{n}$ and $y_{n}$ must belong to this neighborhood.

This contradiction provides the desired set $\widetilde{W}$, and we obtain the desired diffeomorphism
\begin{equation}
\Phi:W\to\Phi(W),  \label{eq:s1-34}
\end{equation}
with the property that $\Gamma_{0}\subset W$ and $\Phi$ maps the wedge onto a neighborhood of the contact curve.

\section{Some notations}

\subsection{The definition of some important constants}

\label{subsec:constants}

\begin{itemize}
\item We define $m$ to be the maximum order of derivative we can take for the elliptic estimate.

\item The integrability exponents $q_{\varepsilon}$, $q_{\pm}$ near the contact line are defined below in terms of the contact angle.

We recall that the equilibrium contact angle $\theta_{\mathrm{eq}}\in(0,\pi)$ of \eqref{eq:young} is measured through the fluid. It is the opening angle of the fluid domain along the contact line $\Gamma$ at equilibrium: in every normal slice $\{\theta=\mathrm{const}\}$ of the chart of Section~\ref{sec:flattening}, the fluid occupies a corner of angle $\theta_{\mathrm{eq}}$. This angle, which can take on any value between $0$ and $\pi$ depending on the choice of the capillary parameters $\sigma$, $\gamma_{sv}$, and $\gamma_{sf}$, plays an important role in the elliptic regularity theory associated to $\Omega$, as it determines the possible regularity gain. For the Stokes problem with the boundary conditions of \eqref{eq:moving-domain-elliptic-system} (stress condition on $\Sigma$, Navier slip and no penetration on $\Sigma_{s}$), the regularity is related to the following parameter, computed by Orlt--S\"andig~\cite{OS}:
\begin{equation}
\varepsilon_{\max} =\varepsilon_{\max}(\theta_{\mathrm{eq}}) :=\min\left\{1,-1+\frac{\pi}{\theta_{\mathrm{eq}}}\right\} \in(0,1].
\end{equation}
For a parameter $0<\varepsilon<\varepsilon_{\max}\leq1$, we set
\begin{equation}
q_\varepsilon:=\frac{2}{2-\varepsilon}\in(1,2).  \label{eq:q_eps}
\end{equation}
Note that $0<\varepsilon_-<\varepsilon_+<\varepsilon_{\max}(\theta_{\mathrm{eq}})$ implies that $q_{\varepsilon_-} <q_{\varepsilon_+} <q_{\varepsilon_{\max}}.$ With $\varepsilon_{\max}$ in hand, we fix three parameters $\alpha$, $\varepsilon_-$, and $\varepsilon_+$ such that
\begin{equation}
0<\alpha<\varepsilon_-<\varepsilon_+<\varepsilon_{\max}, \qquad \alpha< \min\left\{ \frac{\varepsilon_-}{3}, \frac{\varepsilon_+-\varepsilon_-}{2}, \frac{1-\varepsilon_+}{6}, \frac{1}{12} \right\}, \qquad \varepsilon_+\leq\frac{\varepsilon_-+1}{2}.  \label{eq:parameters}
\end{equation}
For brevity, we also write
\begin{equation}
q_+:=q_{\varepsilon_+} \qquad\text{and}\qquad q_-:=q_{\varepsilon_-} \label{eq:q_pm}
\end{equation}
in the notation established in (2.3). We will crucially use these parameters to track regularity in this paper.

\item We choose $\alpha$ to be a constant sufficiently close to $0$.

\item We choose $p$ to be the index for the elliptic estimate in the interior part. It is chosen to be sufficiently close to $2$ such that $p<\frac{2}{1+2\alpha}$.

\item The elliptic norms of $\mathcal{E}$ and $\mathcal{D}$ near the contact line use the single exponent $q_{\delta}:=q_{-}$, that is, $\delta=\varepsilon_{-}$, in both $\mathcal{E}$ and $\mathcal{D}$. In one dimension $W^{3-1/q_{\delta},q_{\delta}}\subset H^{3/2+\delta}$; the estimates \eqref{eq:E_local_enhance}--\eqref{eq:mixed_enhance} need $\delta\geq\varepsilon_{-}$, and the products of Subsection~8.4 need a single $\delta$.

\item Besides \eqref{eq:parameters}, the proofs use the following conditions.

\begin{itemize}
\item $\delta\geq2\alpha$, used in Lemma~\ref{lem:hard_nonlinearity}, (5.142) and (8.51)--(8.53), and $3\alpha\leq\varepsilon_{-}$, used in \eqref{eq:Fint1-transport}. Both follow from $\alpha<\varepsilon_{-}/3$.

\item $\alpha<(1-\varepsilon_{+})/6$, used in (9.4) through the exponent $\beta\in[3/p-3/2,(1-\varepsilon_{+})/2]$ of (3.29) and (9.4).

\item Lower bounds on $p$: $p\geq12/7$ for the $L^{3}$ bound \eqref{eq:A_d}; $p\geq3/(2-2\alpha)$ for (5.146)--(5.147); $p>3/(2-\alpha)$ for (5.90)--(5.94), (5.110)--(5.113) and (9.21); $p>2/(2-\alpha)$ for \eqref{eq:interior_e_1}; $p\geq6/(5-2\varepsilon_{-})$ for (6.115); $p\geq q_{+}$ for the factors on the intermediate region in Subsection~8.4 and Section~9; $p\geq3/(2-3\alpha)$ for \eqref{eq:full-stress-commutator-estimate}; $p>2/(1+\varepsilon_{-})$ for the kinematic term in the proof of the basic interior estimate, whose exponent $\beta^{\prime }$ needs $0<\beta^{\prime }<\frac12+\frac{\varepsilon_{-}}{2}-\frac1p$; and $p\geq6/(4-\varepsilon_{+})$, so that the range of $\beta$ is not empty. Finally $p\geq6/5$, so that $W^{1,p}\subset L^{2}$ in three dimensions; this is used in Subsection~\ref{subsec:top_level} and follows from $p\geq12/7$.
\end{itemize}

Under \eqref{eq:parameters} each of these lower bounds is smaller than $2/(1+2\alpha)$, so the whole system can be satisfied for every $\theta_{\mathrm{eq}}\in(0,\pi)$: take $\varepsilon_{-}=\varepsilon_{\max}/3$ and $\varepsilon_{+}=2\varepsilon_{\max}/3$, then $\alpha$ as in \eqref{eq:parameters}, and $p$ just below $2/(1+2\alpha)$. For example, for $\theta_{\mathrm{eq}}=2\pi/3$, so that $\varepsilon_{\max}=1/2$, one may take $\varepsilon_{-}=0.2$, $\varepsilon_{+}=0.3$, $\delta=0.2$, $\alpha=0.04$, $p=1.8$ and $\beta=0.2$.

\item \emph{Notation.} The letter $\varepsilon$ without subscript denotes the width of the cutoff functions near the contact line (Subsection~\ref{subsec:cutoff}); it is unrelated to the exponents $\varepsilon_{\pm}$. The letter $\delta$ in $q_{\delta}$ always stands for $\varepsilon_{-}$.
\end{itemize}

\subsection{The choice of cutoff functions}

\label{subsec:cutoff}

To avoid confusion, we clarify the choice of cutoff functions. We first choose $\chi_{0,\Sigma}$ to be the background cutoff function for the local chart. It satisfies the following properties:

\begin{itemize}
\item In a small neighborhood of the contact line $\chi_{\mathrm{local},\Sigma}=1.$ \item $0\leq \chi_{\mathrm{local},\Sigma}\leq 1$.

\item $\chi_{\mathrm{local},\Sigma}=0$ on the interior of the equilibrium surface $\Sigma$.
\end{itemize}

We choose $\chi_{cl,\Sigma}$ to be the background cutoff function that we use to define the extension of the perturbation function $\eta$. It has the following properties:

\begin{itemize}
\item In a neighborhood of the contact line $\chi_{\mathrm{cl},\Sigma}=1.$ \item $0\leq \chi_{\mathrm{cl},\Sigma}\leq 1$.

\item The support set of $\chi_{cl,\Sigma}$ is a subset of the set where $\chi_{\mathrm{local},\Sigma}=1$. Equivalently $\func{supp}\{\chi_{cl,\Sigma}\}\subseteq \{x\in \Sigma|\chi_{\mathrm{local},\Sigma}=1\}.$ \end{itemize}

We choose $\chi_{int,\Sigma}$ to be the cutoff function for the interior domain. It has to satisfy the following properties:

\begin{itemize}
\item In the interior of the equilibrium surface $\chi_{int,\Sigma}=1.$ \item $0\leq \chi_{int,\Sigma}\leq 1$.

\item $\{x\in \Sigma |\chi_{\mathrm{local},\Sigma}\neq 1\}\subseteq \{x\in \Sigma|\chi_{int,\Sigma}=1\}.$ \item $\func{supp}\{\chi_{cl,\Sigma}\}\cap \func{supp}\{\chi_{int,\Sigma}\}=\emptyset.$ \end{itemize}

We choose $\chi_{bet,\Sigma}$ to be the cutoff function for the intermediate domain. It can be written as $\chi_{bet,\Sigma}=1-\chi_{cl,\Sigma}-\chi_{int,\Sigma}.$ Moreover, we define $\chi_{m+1,\mathrm{int},\Sigma}:=1-\chi_{\mathrm{cl},\Sigma}$.

We choose $\chi_{j,\mathrm{int},\Sigma}$, $\chi_{j-,\mathrm{int},\Sigma}$ and $\chi_{j+,\mathrm{int},\Sigma}$ to be the cutoff functions used for the local elliptic estimate near the contact line. They must satisfy the following properties. Suppose $m$ is the highest-order temporal derivative we can take for the elliptic estimate and suppose that $\func{supp}\{\chi_{\mathrm{int},\Sigma}\}\subset \{x\in \Sigma|\chi_{m,int,\Sigma}=1\}$ and that $\func{supp}\{\nabla \chi_{m,\mathrm{int},\Sigma}\}\subset \{x\in \Sigma|\chi_{\mathrm{bet},\Sigma}=1\}$. For any nonnegative integer $j\leq m$

\begin{itemize}
\item $\func{supp}\{\chi_{\mathrm{int},\Sigma}\}\subseteq \{x\in \Sigma|\chi_{j,\mathrm{int},\Sigma}=1\},$  $\func{supp}\{\nabla \chi_{j,\mathrm{int},\Sigma}\}\subseteq \{x\in \Sigma|\chi_{\mathrm{bet},\Sigma}=1\},$ for any $j\leq m$.

\item $0\leq \chi_{j,\Sigma}\leq 1$, $0\leq \chi_{j\pm,\Sigma}\leq 1$.

\item $\func{supp}\{\chi_{j-1+,int,\Sigma}\}\subseteq \{x\in\Sigma|\chi_{j-,int,\Sigma}=1\},$  $\func{supp}\{\chi_{j-,int,\Sigma}\}\subseteq \{x\in\Sigma|\chi_{j,int,\Sigma}=1\},$  $\func{supp}\{\chi_{j,int,\Sigma}\}\subseteq \{x\in\Sigma|\chi_{j+,int,\Sigma}=1\},$  $\func{supp}\{\chi_{j+,int,\Sigma}\}\subseteq \{x\in \Sigma|\chi_{j+1-,int,\Sigma}=1\}.$ \end{itemize}

Moreover, we let $\chi_{\mathrm{int}+,\Sigma}$ denote the cutoff function defined by the following properties:

\begin{itemize}
\item $0\leq \chi_{\mathrm{int}+,\Sigma}\leq 1$.

\item $\func{supp}\{\chi_{\mathrm{int},\Sigma}\}\subseteq \{x\in\Sigma|\chi_{\mathrm{int}+,\Sigma}=1\},$  $\func{supp}\{\chi_{\mathrm{int}+,\Sigma}\}\subseteq \{x\in\Sigma|\chi_{0-,int,\Sigma}=1\}.$ \end{itemize}

We choose $\chi_{\mathrm{cl}+,\Sigma}$ to be $\chi_{\mathrm{cl}+,\Sigma}=1-\chi_{\mathrm{int}+}.$ For smooth cutoffs, write $a\prec b$ when $b=1$ on a neighborhood of $\operatorname{supp}a$. We choose the endpoint cutoffs so that
\begin{equation}
\chi_{m+,\mathrm{int}}\prec\chi_{m+1-,\mathrm{int}} \prec\chi_{m+1,\mathrm{int}}\prec\chi_{m+1+,\mathrm{int}}, \qquad \chi_{m+1,\mathrm{int}}:=1-\chi_{\mathrm{cl}}.
\end{equation}
The surface cutoffs are their traces. All inclusions have positive gaps in the fixed profiles before the scaling by $\varepsilon$. In particular, $\chi_{m+1-,\mathrm{int},\Sigma}\chi_{\mathrm{cl},\Sigma}=0.$ Since $\chi_{\mathrm{local},\Sigma}$ is independent of $\theta$, its support can be written as $\mathcal{W}\times\mathbb{T}$ in the chart. The cross-section $\mathcal{W}$ is fixed by $\zeta$ and does not depend on $\varepsilon$, and $\Phi$ is a diffeomorphism on the collar $\Phi(\mathbb{T}\times\mathcal{W})$. Only the supports of the cutoff functions near the contact line are $\varepsilon$-thin: we choose $\chi_{\mathrm{cl},\Sigma}$ and $\chi_{\mathrm{bet},\Sigma}$, and their extensions, so that every point $(r,h)$ of their supports satisfies
\begin{align}
|r|+|h|<\varepsilon.  \label{eq:cutoff_eps}
\end{align}

By hypotheses (2) and (3) of the Main result, the lateral wall in the upper collar and the equilibrium graph are smooth up to the contact line. In particular $\zeta$ has a smooth extension across $\partial D$ in this collar, and the chart $\Phi$ of Section~\ref{sec:flattening} is smooth. Set $Z=Q^T(D\Phi)^{-T}-I$. Since $Z(\theta,0,0)=0$ for every $\theta$, for each angular order $j$ used below we have $|\partial_\theta^j Z(\theta,r,h)| \leq C_{j,\partial D,\zeta} (|r|+|h|).$ This follows by integrating the transverse gradient of $\partial_\theta^j Z$ from $(0,0)$ to $(r,h)$ in a fixed collar. Since $Q^{T}\mathcal{A}_{\Phi}=I$ on $\Gamma=\{r=h=0\}$ for every $\theta$ (Subsection~2.11), on these supports
\begin{align}
|(Q^{T}\mathcal{A}_{\Phi}-I)(r,\theta,h)|\lesssim C_{j,\partial D,\zeta}(|r|+|h|)\lesssim C_{j,\partial D,\zeta}\varepsilon, \label{eq:QA_small}
\end{align}
and
\begin{align}
|\partial_{\theta}^{j}(Q^{T}\mathcal{A}_{\Phi})(r,\theta,h)|\lesssim C_{j,\partial D,\zeta}(|r|+|h|)\lesssim C_{j,\partial D,\zeta}\varepsilon \label{eq:QA_theta_small}
\end{align}
for every $\theta\in\mathbb{T}$ and $j\geq1$. The domains on which the elliptic and the Ne\caronaccent{c}as constants are taken have a size fixed by $\zeta$ and the boundary regularity; see Subsection~\ref{subsec:artificial}.

When a derivative falls on one of these cutoff functions, it produces a factor $\varepsilon^{-1}$. Precisely, we fix profiles $\tilde{\chi}_{\mathrm{cl}}$, $\tilde{\chi}_{\mathrm{bet}}$ that do not depend on $\varepsilon$ and set $\chi_{\mathrm{cl}}(r,h)=\tilde{\chi}_{\mathrm{cl}}((r,h)/\varepsilon)$ and $\chi_{\mathrm{bet}}(r,h)=\tilde{\chi}_{\mathrm{bet}}((r,h)/\varepsilon)$ near $\Gamma$, and the transitions of the nested cutoffs $\chi_{j,\mathrm{int}}$, $\chi_{j\pm,\mathrm{int}}$, $\chi_{\mathrm{int}\pm}$ near $\Gamma$ are scaled in the same way. Hence
\begin{align}
|\nabla^{k}\chi|\leq C_{k}\,\varepsilon^{-k}\qquad\text{near }\Gamma, \label{eq:cutoff_scaling}
\end{align}
with $C_{k}$ independent of $\varepsilon$. Since $\Gamma\subset\partial\Omega $, every point $x=\Phi(\theta,r,h)$ of the supports in \eqref{eq:cutoff_eps} satisfies $\func{dist}(x,\partial\Omega)\leq C(|r|+|h|)\leq C\varepsilon$.

\subsection{The extension of the cutoff functions}

\label{subsec:cutoff_extension}

In the last subsection, we defined cutoff functions on the free surface $\Sigma$. We now extend them to the three-dimensional domain $\Omega$. We fix a level $h_{0}$ below the range of $\zeta$ such that the side wall is vertical on $\{z\geq h_{0}\}$, and we choose the function $\varphi$ in \eqref{eq:W-definition} so that $\varphi(z)=0$ for $z\leq h_{0}$.

\begin{itemize}
\item \emph{Tube cutoffs.} The cutoffs $\chi_{\mathrm{cl}}$, $\chi_{\mathrm{bet}}$, $\chi_{\mathrm{bet}+}=\chi_{\mathrm{cl}}+\chi_{\mathrm{bet}}$ and the shells $\chi_{j,\mathrm{int}}-\chi_{j-1,\mathrm{int}}$, together with their $\pm$ versions, are defined in the chart $\Phi$ of Section~\ref{sec:flattening}. Each of them is $\tilde{\chi}\circ\Phi^{-1}$ for a smooth function $\tilde{\chi}=\tilde{\chi}(r,h)$ that is compactly supported in the fixed cross-section $\mathcal{W}$ of Subsection~\ref{subsec:cutoff}. Hence these cutoffs do not depend on $\theta$, $[\partial_{\theta},\chi]=0$, and their supports and their nesting are those of the corresponding surface cutoffs.

\item \emph{Interior cutoffs.} Each interior cutoff is the complement of tube cutoffs. In particular $\chi_{\mathrm{int}}+\chi_{\mathrm{bet}}+\chi_{\mathrm{cl}}=1,$ and $\chi_{j,\mathrm{int}}$, $\chi_{j\pm,\mathrm{int}}$, $\chi_{\mathrm{int}\pm}$ are defined by the summation of the tube extension of $\chi_{j+1,\mathrm{int}}-\chi_{j,\mathrm{int}}$ and $\chi_{\mathrm{int}}$, with their vertical transitions near $z=h_{0}$ nested in the same order as their transitions in the tube. For example, $\func{supp}\{\nabla(\chi_{\mathrm{bet}}+\chi_{\mathrm{cl}})\neq 0\}\subset\{x\in\Omega|\chi_{j,\mathrm{int}}=1\}$ for any $j\geq 0$. Thus the interior cutoffs equal $1$ for $z\leq h_{0}-c_{1}$ for some constant $c_{1}$, and their complements $1-\chi_{j,\mathrm{int}}$ are tube cutoffs.
\end{itemize}

On $\Sigma$, each extension equals the corresponding surface cutoff. Near the bottom only global norms are used: the forcing terms carried by $\chi_{\mathrm{int}+}$ in Theorem~\ref{thm:whole-space-energy} are bounded through the $H^{1}$ norms of $\partial_{t}^{j}v$ in $\mathcal{E}_{\mathrm{enhanced}}$ and $\mathcal{D}_{||}$; see the paragraph after \eqref{eq:forcing_split}.

\subsection{The artificial domain}

\label{subsec:artificial}

In this paper, we use the notation $\Sigma_{}\mathrm{art}$ for the artificial boundary of the domains used for the interior elliptic estimate, the intermediate elliptic estimate, the contact-line elliptic estimate and the contact-line energy estimate.

For example, for the contact line elliptic estimate, we construct $\Omega_{\mathrm{cl}}$ bounded by $\Sigma_{s},\Sigma$ and $\Sigma_{\mathrm{art}}$ such that $\partial \Omega_{\mathrm{art}}$ is smooth except on $\Sigma_{s}\cap\Sigma$ which is the contact line. Moreover, $\Sigma_{\mathrm{art}}$ is nicely chosen such that $\func{supp}\{\chi_{\mathrm{cl}}\}\subset \Omega_{\mathrm{cl}}\subset \text{Tubular neighborhood of the contact line}.$ In the same way, $\Omega_{\mathrm{bet}}$ and $\Omega_{\mathrm{bet}+}$ are domains bounded by $\Sigma$, $\Sigma_{s}$ and an artificial boundary $\Sigma_{\mathrm{art}}$, contained in the collar $\Phi(\mathbb{T}\times\mathcal{W})$, with $\func{supp}\chi_{\mathrm{bet}}\subset\Omega_{\mathrm{bet}}$ and $\func{supp}\chi_{\mathrm{bet}+}\subset\Omega_{\mathrm{bet}+}$, where $\chi_{\mathrm{bet}+}=\chi_{\mathrm{cl}}+\chi_{\mathrm{bet}}$. They are the domains used in Section~\ref{sec:contact_energy}. Lemma~\ref{lem:enhance_pressure} uses zero-trace tests on $\Omega_{\mathrm{bet}+}$. In Lemma~\ref{lem:local_threehalf} the unknowns are supported in that collar, but the harmonic test is constructed on $\Omega$ and restricted to their support. The size of the domain $\tilde{\Omega}_{\mathrm{cl}}$ depends only on the size of the neighborhood, which depends only on the regularity of the steady state $\zeta$, and is independent of the size of the supports of the cutoff functions. Therefore, in this paper, the constant depending on $\Omega,\Omega_{\mathrm{int}},\Omega_{\mathrm{bet}},\Omega_{\mathrm{cl}}$ does not depend on the smallness of the supports of the cutoff functions we choose.

The size of the supports of the cutoff functions matters in two places only. First, in the commutators between spatial derivatives and cutoff functions: each derivative that falls on a cutoff costs a factor $\varepsilon^{-1}$, so these terms, such as $p_{k}\nabla\chi$, are bounded through $\mathcal{D}_{\mathrm{int}}$ and are not absorbed; see the cutoff forcing terms in the proof of Theorem~\ref{thm:contact_line_energy}. Second, in \eqref{eq:QA_small}--\eqref{eq:QA_theta_small}, which give the factors $\varepsilon^{\iota}$ in Section~\ref{sec:contact_energy}. This is one of the major differences between the three-dimensional and the two-dimensional problems.

\subsection{Some notation conventions}

\label{subsec:notation}

$~$

If the domain of a norm is not specified, the norm is taken over $\Omega$.  If a norm is written as $\|\|_{\mathcal{B}_{1}\mathcal{B}_{2}}$ with $\mathcal{B}_{1},\mathcal{B}_{2}$ two different Sobolev spaces, then the first Sobolev space is in $\theta\in \mathbb{T}:=[0,2\pi]$.  The notation $a\lesssim b$ means
\begin{align}
a\leq Cb,  \label{eq:lesssim}
\end{align}
where the constant $C$ depends only on the steady state $\zeta$, $\alpha$, $p$, $q_{\pm}$, the domain $\Omega$ and the solid wall $\Sigma_{s}$.

\subsection{The definition of energy and dissipation functionals}

\label{subsec:energy_def}

$~$

We introduce the energy functional as follows
\begin{align}
\mathcal{E}:=\mathcal{E}_{||}+\mathcal{E}_{\mathrm{enhanced}}+\mathcal{E}_{\mathrm{int}}+\mathcal{E}_{\mathrm{cl}}+\mathcal{E}_{\mathrm{bet}}.
\end{align}
We introduce the dissipation functional as follows
\begin{align}
\mathcal{D}:=\mathcal{D}_{||}+\mathcal{D}_{\mathrm{enhanced}}+\mathcal{D}_{\mathrm{int}}+\mathcal{D}_{\mathrm{cl}}+\mathcal{D}_{\mathrm{bet}}.
\end{align}

We first define the energy and dissipation functionals in the whole space.
\begin{gather}
\mathcal{E}_{||}:=\sum_{j=0}^{m+1}(\|\partial_{t}^{j}v\|_{L^{2}}^{2}+\|\partial_{t}^{j}\eta\|_{H^{1}(\Sigma)}^{2}). \qquad \mathcal{D}_{||}:=\sum_{j=0}^{m+1}(\|\partial_{t}^{j}v\|_{H^{1}}^{2}+\|\partial_{t}^{j+1}\eta\|^{2}_{L^{2}(\Gamma_{0})}). \\
\mathcal{E}_{\mathrm{enhanced}}:=\sum_{j=0}^{m}(\|\partial_{t}^{j}\eta\|_{H^{\frac{3}{2}-\alpha}(\Sigma)}^{2}+\|\partial_{t}^{j}p\|_{L^{2}}^{2}+\|\partial_{t}^{j}v\|_{H^{1}}^{2}+\|\partial_{t}^{j+1}\eta\|^{2}_{L^{2}(\Gamma_{0})}+\|\partial_{t}^{j+1}\eta\|_{H^{\frac{1}{2}-\alpha}}^{2}). \\
\mathcal{D}_{\mathrm{enhanced}}:=\sum_{j=0}^{m+1}(\|\partial_{t}^{j}\eta\|_{H^{\frac{3}{2}-\alpha}(\Sigma)}^{2}).
\end{gather}

We now define the interior energy and dissipation functionals
\begin{equation}
\begin{aligned}
&{\label{def:interior_energy}} \mathcal{E}_{\mathrm{int}}:=\mathscr{E}_{m,\mathrm{el},\mathrm{int}}:= \sum_{k=0}^{m-1}\sum_{n=0}^{k}(\|\chi_{k+1,int}\partial_{t}^{n}v\|^{2}_{W^{2+(m-1-k),p}}+\|\chi_{k+1,int}\partial_{t}^{n}p\|^{2}_{W^{1+(m-1-k),p}} \\
&\quad +\|\chi_{k+1,int,\Sigma}\partial_{t}^{n}\eta\|^{2}_{W^{3+m-1-k-\frac{1}{p},p}}) +\sum_{k=0}^{m}(\|\chi_{0,\mathrm{int}}\partial_{t}^{k}v\|^{2}_{H^{1+\frac{\epsilon_{-}}{2}}}+\|\chi_{0,\mathrm{int},\Sigma}\partial_{t}^{k}\eta\|^{2}_{H^{\frac{3}{2}+\frac{\epsilon_{-}-\alpha}{2}}}),
\end{aligned}
\end{equation}
\begin{align}
{\label{def:interior_dissipation}} \mathcal{D}_{\mathrm{int}}:=\mathscr{D}_{m,\mathrm{el},\mathrm{int}}:=\sum_{k=0}^{m}\sum_{n=0}^{k}(\|\chi_{k,int}\partial_{t}^{n}v\|^{2}_{W^{2+(m-k),p}}+\|\chi_{k,int}\partial_{t}^{n}p\|^{2}_{W^{1+(m-k),p}}+\|\chi_{k,int,\Sigma}\partial_{t}^{n}\eta\|^{2}_{W^{3+m-k-\frac{1}{p},p}}).
\end{align}

We now define the energy and dissipation near the contact line and in the intermediate part. We define
\begin{equation}
\begin{aligned}
&\mathcal{E}_{\mathrm{cl}}:=\sum_{k_{1}+k_{2}=m-1}\mathscr{E}_{k_{1},k_{2},\mathrm{el},\mathrm{cl}}+\sum_{k_{1}+k_{2}=m}\mathscr{E}_{k_{1},k_{2},\mathrm{energy},\mathrm{cl}} +\sum_{k_{1}+k_{2}=m}(\|\chi_{\mathrm{cl}}\partial_{t}^{k_{1}}\partial_{\theta}^{k_{2}}v\|^{2}_{L^{2}H^{1+\frac{\epsilon_{-}}{2}}} \\
&\quad +\|\chi_{\mathrm{cl},\Sigma}\partial_{t}^{k_{1}}\partial_{\theta}^{k_{2}}\eta\|^{2}_{L^{2}H^{\frac{3}{2}+\frac{\epsilon_{-}-\alpha}{2}}}), +\sum_{k_{1}+k_{2}=m+1}\sum_{j_{1}=0}^{k_{1}}\sum_{j_{2}=0}^{k_{2}}(\|(\chi_{cl}+\chi_{\mathrm{bet}})\partial_{t}^{j_{1}}\partial_{\theta}^{j_{2}}\eta\|^{2}_{H^{1}}+\|(\chi_{cl}+\chi_{\mathrm{bet}})\partial_{t}^{j_{1}}\partial_{\theta}^{j_{2}}v\|_{L^{2}}^{2}), \\
&\quad \mathcal{E}_{\mathrm{bet}}:=\sum_{k_{1}+k_{2}=m-1}\mathscr{E}_{k_{1},k_{2},\mathrm{el},\mathrm{bet}}+\sum_{k_{1}+k_{2}=m}\mathscr{E}_{k_{1},k_{2},\mathrm{energy},\mathrm{bet}} \\
&\quad +\sum_{k_{1}+k_{2}=m}(\|\chi_{\mathrm{bet}}\partial_{t}^{k_{1}}\partial_{\theta}^{k_{2}}v\|_{L^{2}H^{1+\frac{\epsilon_{-}}{2}}}^{2}+\|\chi_{\mathrm{bet},\Sigma}\partial_{t}^{k_{1}}\partial_{\theta}^{k_{2}}\eta\|_{L^{2}H^{\frac{3}{2}+\frac{\epsilon_{-}-\alpha}{2}}}^{2}),
\end{aligned}
\end{equation}
and
\begin{equation}
\begin{aligned}
&\mathcal{D}_{\mathrm{cl}}:=\sum_{k_{1}+k_{2}=m}\mathscr{E}_{k_{1},k_{2},\mathrm{el},\mathrm{cl}}+\sum_{k_{1}+k_{2}=m+1}\mathscr{E}_{k_{1},k_{2},\mathrm{energy},\mathrm{cl}}\\
&\quad +\sum_{k_{1}+k_{2}=m}(\|\chi_{\mathrm{cl},\Sigma}\partial_{t}^{k_{1}}\partial_{\theta}^{k_{2}}\eta\|^{2}_{L^{\infty}H^{\frac{3}{2}+\frac{\epsilon_{-}-\alpha}{2}}}+\|\chi_{\mathrm{cl}}\partial_{t}^{k_{1}}\partial_{\theta}^{k_{2}}v\|^{2}_{L^{\infty}H^{1+\frac{\epsilon_{-}}{2}}}),   \\
&\quad \mathcal{D}_{\mathrm{bet}}:=\sum_{k_{1}+k_{2}=m}\mathscr{E}_{k_{1},k_{2},\mathrm{el},\mathrm{bet}}+\sum_{k_{1}+k_{2}=m+1}\mathscr{E}_{k_{1},k_{2},\mathrm{energy},\mathrm{bet}}   \\
&\quad +\sum_{k_{1}+k_{2}=m}(\|\chi_{\mathrm{bet},\Sigma}\partial_{t}^{k_{1}}\partial_{\theta}^{k_{2}}\eta\|_{L^{\infty}H^{\frac{3}{2}+\frac{\epsilon_{-}-\alpha}{2}}}^{2}+\|\chi_{\mathrm{bet}}\partial_{t}^{k_{1}}\partial_{\theta}^{k_{2}}v\|_{L^{\infty}H^{1+\frac{\epsilon_{-}}{2}}}^{2}),
\end{aligned}
\end{equation}
where
\begin{align}
{\label{def:split_energy}} \begin{aligned} \mathscr{E}_{k_{1},k_{2},el,cl}:= &\sum_{j_{1}=0}^{k_{1}}\sum_{j_{2}=0}^{k_{2}}(\|\chi_{cl}\partial_{t}^{j_{1}}\partial_{\theta}^{j_{2}}v\|^{2}_{L^{2}W^{2,q_{\delta}}}+\|\chi_{cl,\Sigma}\partial_{t}^{j_{1}}\partial_{\theta}^{j_{2}}\eta\|^{2}_{L^{2}W^{3-\frac{1}{q_{\delta}},q_{\delta}}}+\|\chi_{cl}\partial_{t}^{j_{1}}\partial_{\theta}^{j_{2}}p\|_{L^{2}W^{1,q_{\delta}}}^{2}),\\
\mathscr{E}_{k_{1},k_{2},el,bet}:=&\sum_{j_{1}=0}^{k_{1}}\sum_{j_{2}=0}^{k_{2}}(\|\chi_{bet}\partial_{t}^{j_{1}}\partial_{\theta}^{j_{2}}v\|^{2}_{L^{2}W^{2,q_{\delta}}}+\|\chi_{\mathrm{bet},\Sigma}\partial_{t}^{j_{1}}\partial_{\theta}^{j_{2}}\eta\|^{2}_{L^{2}W^{3-\frac{1}{q_{\delta}},q_{\delta}}}+\|\chi_{\mathrm{bet}}\partial_{t}^{j_{1}}\partial_{\theta}^{j_{2}}p\|_{L^{2}W^{1,q_{\delta}}}^{2}), \end{aligned}
\end{align}
and when $k_{1}\neq m+1$
\begin{equation}
\begin{aligned}
&\mathscr{E}_{k_{1},k_{2},energy,cl}:=\sum_{j_{1}=0}^{k_{1}}\sum_{j_{2}=0}^{k_{2}}(\|(\chi_{cl}+\chi_{\mathrm{bet}})\partial_{t}^{j_{1}}\partial_{\theta}^{j_{2}}v\|^{2}_{H^{1}} +\|(\chi_{cl,\Sigma}+\chi_{\mathrm{bet},\Sigma})\partial_{t}^{j_{1}}\partial_{\theta}^{j_{2}}\eta\|^{2}_{H^{\frac{3}{2}-\alpha}} \\
&\quad +\|(\chi_{cl,\Sigma}+\chi_{\mathrm{bet},\Sigma})\partial_{t}^{j_{1}+1}\partial_{\theta}^{j_{2}}\eta\|^{2}_{H^{\frac{1}{2}-\alpha}}) +\sum_{j_{1}=0}^{k_{1}}\sum_{j_{2}=0}^{k_{2}}(\|\partial_{t}^{j_{1}+1}\partial_{\theta}^{j_{2}}\eta\|^{2}_{L^{2}(\Gamma_{0})} \\
&\quad +\|(\chi_{cl,\Sigma}+\chi_{\mathrm{bet},\Sigma})\partial_{t}^{j_{1}}\partial_{\theta}^{j_{2}}\eta\|^{2}_{H^{1}}) +\sum_{j_{1}=0}^{k_{1}}\sum_{j_{2}=0}^{k_{2}}\|(\chi_{\mathrm{cl}}+\chi_{\mathrm{bet}})\partial_{t}^{j_{1}}\partial_{\theta}^{j_{2}}p\|^{2}_{L^{2}}, \label{def:energy_cl}
\end{aligned}
\end{equation}
when $k_{1}=m+1$
\begin{equation}
\begin{aligned}
&\mathscr{E}_{k_{1},k_{2},energy,cl}:=\sum_{j_{1}=0}^{k_{1}}\sum_{j_{2}=0}^{k_{2}}(\|(\chi_{\mathrm{cl}}+\chi_{\mathrm{bet}})\partial_{t}^{j_{1}}\partial_{\theta}^{j_{2}}v\|^{2}_{H^{1}} +\|(\chi_{\mathrm{cl},\Sigma}+\chi_{\mathrm{bet},\Sigma})\partial_{t}^{j_{1}}\partial_{\theta}^{j_{2}}\eta\|^{2}_{H^{\frac{3}{2}-\alpha}} \\
&\quad +\|(\chi_{\mathrm{cl},\Sigma}+\chi_{\mathrm{bet},\Sigma })\partial_{t}^{j_{1}+1}\partial_{\theta}^{j_{2}}\eta\|^{2}_{H^{\frac{1}{2}-\alpha}}) +\sum_{j_{1}=0}^{k_{1}}\sum_{j_{2}=0}^{k_{2}}(\|\partial_{t}^{j_{1}+1}\partial_{\theta}^{j_{2}}\eta\|^{2}_{L^{2}(\Gamma_{0})} \\
&\quad +\|(\chi_{\mathrm{cl},\Sigma}+\chi_{\mathrm{bet},\Sigma})\partial_{t}^{j_{1}}\partial_{\theta}^{j_{2}}\eta\|^{2}_{H^{1}}) +\sum_{j_{1}=0}^{m}\|(\chi_{\mathrm{cl}}+\chi_{\mathrm{bet}})\partial_{t}^{j_{1}}p\|^{2}_{L^{2}}.
\end{aligned}
\end{equation}
In both cases we set $\mathscr{E}_{k_{1},k_{2},\mathrm{energy},\mathrm{bet}}:=\mathscr{E}_{k_{1},k_{2},\mathrm{energy},\mathrm{cl}}$. The energy estimates of Section~\ref{sec:contact_energy} are carried out for the unknowns localized by $\chi_{\mathrm{bet}+}=\chi_{\mathrm{cl}}+\chi_{\mathrm{bet}}$, so the energy pieces near the contact line and on the intermediate region carry the same weight $\chi_{\mathrm{cl}}+\chi_{\mathrm{bet}}$.

\subsection{Constants and the order in which they are chosen}

\label{subsec:smallness}

Constants written $C$, or hidden in $\lesssim$, depend only on the quantities listed after \eqref{eq:lesssim}; in particular they do not depend on $\varepsilon$. Constants that may depend on $\varepsilon$ are written $C(\varepsilon)$; each $C(\varepsilon)$ is a finite power of $\varepsilon^{-1} $ times a constant of the first kind, but only its finiteness for fixed $\varepsilon$ is used, and its value may change from line to line. The width $\varepsilon$ is fixed first, small enough for the absorptions in the proof of Theorem~\ref{thm:contact_line_energy}, which use only constants of the first kind. The smallness of $\mathcal{E}$ is chosen afterwards, depending on $\varepsilon$. Thus a nonlinear term such as
\begin{align}
C(\varepsilon)\sqrt{\mathcal{E}}\sqrt{\mathcal{D}}\sqrt{\mathcal{D}}\quad\text{or}\quad C(\varepsilon)\sqrt{\mathcal{E}}\,\mathcal{E}
\end{align}
is absorbed once $C(\varepsilon)\sqrt{\mathcal{E}}$ is small. When the factor $C(\varepsilon)$ in front of a nonlinear term plays no role, we suppress it and write $\sqrt{\mathcal{E}}\mathcal{D}$ or $\mathcal{E}^{3/2}$.

\subsection{The conservation of total volume}

\label{subsec:volume}

We choose the initial height so that the fluid has its equilibrium volume:
\begin{equation}\label{eq:volume_initial} \int_D (\zeta+\eta(0))\,dx=\int_D\zeta\,dx .
\end{equation}
Since the perturbed surface is a graph over the fixed horizontal domain $D$, $|\Omega(t)|-|\Omega|=\int_D\eta(x,t)\,dx .$ On the graph, the outward unit normal is $n=\mathcal N/|\mathcal N|$ and $dS=|\mathcal N|\,dx$. The kinematic condition therefore gives
\begin{equation}
\frac{d}{dt}\int_D\eta(x,t)\,dx =\int_D u(x,\zeta(x)+\eta(x,t),t)\cdot\mathcal N(x,t)\,dx =\int_{\Sigma(t)}u\cdot n\,dS .
\end{equation}
Since $u\cdot n_0=0$ on the wetted solid wall, the divergence theorem and incompressibility in the moving domain give
\begin{equation}
\int_{\Sigma(t)}u\cdot n\,dS =\int_{\Omega(t)}\operatorname{div}u\,dX=0 .
\end{equation}
Thus $\int_D\eta(x,t)\,dx=0$ for all $t$ by \eqref{eq:volume_initial}.

\subsection{The computation for the geometric matrix}

Write $q=W\bar\eta$ as in \eqref{eq:A-explicit}. Then $\mathcal A_\eta=\mathfrak A(\nabla q)$, where $\mathfrak A(0)=I$ and $\mathfrak A$ and its derivatives are bounded on $|\nabla q|\leq1/2$. The first time derivative is linear in $\nabla\partial_t q$, with a coefficient depending on $\nabla q$. For $j\geq1$ the chain rule gives
\begin{equation}\label{eq:A_equi} \|\partial_t^j\mathcal A_\eta\|_{L^\infty} \leq C_j\sum_{r=1}^{j}\ \sum_{\substack{i_1+\cdots+i_r=j\\i_1,\ldots,i_r\geq1}} \prod_{k=1}^{r}\|\nabla\partial_t^{i_k}q\|_{L^\infty}.
\end{equation}
The same inequality holds with all norms on the same subset of $\Omega$. The term $r=1$ is the only term containing the $j$th time derivative. We will show in the following Section that
\begin{align}
\|\p_{t}^{m}(W\bar{\eta})\|_{W^{1,\infty}}\lesssim \sqrt{\mathcal{D}}, \quad\operatorname{and } \quad\sum_{j=0}^{m-1}\|\p_{t}^{j}(W\bar{\eta})\|_{W^{1,\infty}}\lesssim\sqrt{\mathcal{E}}
\end{align}
For $r\geq2$ we shall bound one factor by dissipation and the remaining factors by energy when a dissipation bound is needed; see \eqref{eq:interior_e_1} for the localized estimate. Small energy allows higher powers of $\sqrt{\mathcal E}$ in the resulting sum to be absorbed into its bound. The same argument applies to $K=(1+\partial_3q)^{-1}$.

\subsection{Exact decomposition of $Q^{T}\mathcal{A}_{\Phi}-I$}

Recall that the horizontal part of the local transformation before reparametrization: $X_{\mathrm{old}}(s,r,h) := \gamma(\tau(s,r,h)) +r\dot{\gamma}^{\perp}(\tau(s,r,h))$ for $s\in\mathbb T_L$.  Define
\begin{equation*}
\Lambda_{\mathrm{old}}(s):= \sqrt{1+|\nabla\zeta(\gamma(s))\cdot\dot\gamma(s)|^2}, \qquad \mathscr H:=h+r\cot\theta_{\mathrm{eq}}.
\end{equation*}
For $\theta\in\mathbb T_M$, we define
\begin{equation}
\Lambda_*(\theta) :=\Lambda_{\mathrm{old}}(\chi^{-1}(\theta)), \qquad X_*(\theta,r,h) :=X_{\mathrm{old}}(\chi^{-1}(\theta),r,h).
\end{equation}
Then, since $\chi'=\Lambda_{\mathrm{old}}$,
\begin{equation}
\Phi=(X_*,\zeta(X_*)+\Lambda_*\mathscr H), \qquad \Lambda_*'(\theta)= \frac{\Lambda_{\mathrm{old}}'(\chi^{-1}(\theta))} {\Lambda_{\mathrm{old}}(\chi^{-1}(\theta))}.
\end{equation}
Write $X_{*,\Gamma}=X_*(\theta,0,0)$, $DX_{*,\Gamma}=DX_*(\theta,0,0)$, $D\widetilde X_{*,\Gamma}=D\widetilde X_*(\theta,0,0)$, and $\widetilde X_*=(X_*,0)$. With $Q=D\Phi(\theta,0,0)$ and $\mathcal A_\Phi=(D\Phi)^{-T}$,
\begin{gather}
{\label{eq:exact_QTA_identity_boxed}} D\Phi-Q={} D\widetilde X_*-D\widetilde X_{*,\Gamma} +e_h\otimes\bigl[(DX_*)^T\nabla\zeta(X_*) -(DX_{*,\Gamma})^T\nabla\zeta(X_{*,\Gamma})\bigr] +\Lambda_*'\mathscr H\,e_h\otimes e_\theta, \\
Q^T\mathcal A_\Phi-I={} -(D\Phi-Q)^T\mathcal A_\Phi.
\end{gather}
Thus the three contributions are
\begin{gather}
{\label{eq:Phi_exact_decomposition}} Q_X :=-(D\widetilde X_*-D\widetilde X_{*,\Gamma})^T\mathcal A_\Phi, \\
Q_\zeta :=-\bigl[(DX_*)^T\nabla\zeta(X_*) -(DX_{*,\Gamma})^T\nabla\zeta(X_{*,\Gamma})\bigr] \otimes(\mathcal A_\Phi^Te_h), \qquad Q_{\mathrm{ang}} :=-\Lambda_*'\mathscr H\, e_\theta\otimes(\mathcal A_\Phi^Te_h).
\end{gather}

Substituting \eqref{eq:Phi_exact_decomposition} into \eqref{eq:exact_QTA_identity_boxed}, we find $Q^{T}A-I = Q_{X}+Q_{\zeta}+Q_{\mathrm{ang}}.$ Let $\mathcal{B}:=Q^{T}\mathcal{A}_{\Phi}=I+\mathcal{Z}, \; \mathcal{Z}=Q_X+Q_\zeta+Q_{\mathrm{ang}},$ and introduce the skew-symmetric connection matrix $\mathscr C:=Q^{T}\partial_{\theta}Q.$ Since $Q^{T}Q=I$, we have $(\partial_{\theta}Q^{T})Q=-\mathscr C, \; \mathscr C^{T}=-\mathscr C.$ The coefficient bounds depend on the full fixed geometry, including the curvature of $\partial D$. For every angular order $i$ used below, we have $\|\partial_\theta^i\mathcal Z\|_{L^\infty(\mathbb T_M\times\mathcal W)} \le C_{i,\partial D,\zeta}.$ Since $\partial_\theta^i\mathcal Z(\theta,0,0)=0$, integration in the transverse variables gives the following estimates for the difference matrix $\mathcal{Z}$
\begin{equation}
|\partial_\theta^i\mathcal Z(\theta,r,h)| \le C_{i,\partial D,\zeta}(|r|+|h|), \qquad \|\partial_\theta^i\mathcal Z\|_{L^\infty(\operatorname{supp}\chi_{\mathrm{bet}+})} \le C_{i,\partial D,\zeta}\varepsilon.
\end{equation}

\section{Localized surface function}

\label{sec:localized_surface}

To define the extension of the surface perturbation function $\eta$, we first express $\eta$ in local coordinates. In this section and all subsequent uses of the contact-line chart, write $X:=X_*$ and $\Lambda:=\Lambda_*$. Thus the angular variable is the final parameter $\theta\in\mathbb T_M$, and
\begin{equation*}
\Phi(\theta,r,h) = \left( X(\theta,r,h), \zeta\bigl(X(\theta,r,h)\bigr) + \Lambda(\theta)\bigl(h+r\cot\theta_{\mathrm{eq}}\bigr) \right).
\end{equation*}

Let the steady surface in the physical coordinates be $\Sigma = \left\{ (x,z)\in\mathbb{R}^{2}\times\mathbb{R }\;\middle|\; z=\zeta(x) \right\}.$ A point $(\theta,r,h)$ belongs to the representation of $\Sigma$ in the local chart if and only if $\Phi(\theta,r,h)\in\Sigma.$ We note that the equilibrium surface in the local chart can be represented by
\begin{equation}  \label{eq:equilibrium-H} H_0(\theta,r) = -r\cot\theta_{\mathrm{eq}}.
\end{equation}
Suppose $X=X(\theta,r,h)$ is induced by the diffeomorphism $\Phi$ and is smooth in a neighborhood of the equilibrium surface. Then any point on the equilibrium surface near the contact line can be represented by the surface map $\widehat X(\theta,r):=X(\theta,r,H_{0}(\theta,r)).$ Therefore, the representation of the original perturbation $\eta$ in the local surface coordinates is given by
\begin{equation}  \label{eq:local-eta-definition} \mathfrak{h}_{\eta}(\theta,r,t) := \eta\bigl(\widehat X(\theta,r),t\bigr).
\end{equation}
By \eqref{eq:equilibrium-H}, this is equivalently given by
\begin{equation}  \label{eq:local-eta-H-relation} \mathfrak{h}_{\eta}(\theta,r,t) := \eta\left( X\bigl(\theta,r,-r\cot\theta_{\mathrm{eq}}\bigr), t \right).
\end{equation}
We shall use $\mathfrak{h}_{\eta}$ to represent $\eta$ in the local chart.

\emph{Remark.} We now have the representation of $\eta$ in the local chart. The remaining problem is that, when we carry out estimates near the contact line or on the intermediate region, the local charts near the contact line always overlap with the interior part. Therefore, we have to establish a theorem comparing the norms of $\eta$ in the interior region with the norms of $\mathfrak{h}_{\eta}$ in the intermediate part. Specifically, we need this equivalence in two cases. The first case is the intermediate elliptic estimate, where the commutator between the cutoff functions and spatial derivatives induces a forcing term in the interior part. The second case is the pure time energy estimate, where we have to derive a norm equivalence near the contact line to complete the estimate. We establish three theorems regarding this norm equivalence.

\begin{theorem}[Localized estimates in the surface coordinates]
\label{thm:localized-surface-coordinate-estimate} Let $H_0(\theta,r)=-r\cot\theta_{\mathrm{eq}},$ and define $\widehat X(\theta,r) = X(\theta,r,H_{0}(\theta,r)),$ and $\mathfrak{h}_\eta(\theta,r,t) = \eta(\widehat X(\theta,r),t).$ Assume that $\eta(\cdot,t)$ is supported in the part $\Sigma\cap\Phi(W)$ of $\Sigma$ covered by the chart \eqref{eq:s1-34}. Then the localized coordinate displacement satisfies
\begin{gather}
\label{eq:localized-heta-H1}
\|\eta\|_{W^{a,b}}\lesssim \|\mathfrak{h}_{\eta}\|_{W^{a,b}} \lesssim \|\eta\|_{W^{a,b}}, \\
\label{eq:localized-heta-H_high_1}
\|\mathfrak{h}_{\eta}\|_{L^{\infty}W^{a,2}}\lesssim \|\eta\|_{W^{a+1,2}}, \\
\label{eq:localized-heta-H_high_2}
\|\mathfrak{h}_{\eta}\|_{L^{2}W^{a,2}}\lesssim \|\eta\|_{W^{a,2}},
\end{gather}
for every $a\geq0$ and $1<b<\infty$, and for every integer $a\geq0$ and $b=\infty$; the constants depend on $a$ and $b$. In particular,
\begin{gather}
\label{eq:localized-heta-H_high_3}
\|\mathfrak{h}_{\eta}\|_{L^{\infty}H^{\frac{5}{2}-\beta}}\lesssim\|\eta\|_{H^{\frac{7}{2}-\beta}}. \\
\label{eq:localized-heta-H_high_4}
\|\mathfrak{h}_{\eta}\|_{L^{2}H^{\frac{5}{2}-\beta}}\lesssim\|\eta\|_{H^{\frac{5}{2}-\beta}}. \\
\label{eq:localized-heta-Hp}
\|\eta\|_{W^{3-\frac{1}{p},p}}\lesssim \|\mathfrak{h}_{\eta}\|_{W^{3-\frac{1}{p},p}} \lesssim \|\mathfrak{h}_{\eta}\|_{L^{2}W^{3-\frac{1}{p},p}}+\|\p_{\theta}\mathfrak{h}_{\eta}\|_{W^{2-\frac{1}{p},p}}.
\end{gather}
\end{theorem}

\begin{proof}
\textbf{Step 1. Low-order control of $\mathfrak{h}_\eta$ and $H^{1}$ estimate.}

Recall the notation $\widehat X(\theta,r) = X(\theta,r,H_{0}(\theta,r)).$ Then $\mathfrak{h}_\eta(r,\theta) = \eta(\widehat X(\theta,r),t).$ Differentiating once in $r$,
\begin{equation}  \label{eq:heta-first-derivative} \partial_r\mathfrak{h}_\eta = \nabla\eta(Y)\cdot \partial_rY.
\end{equation}
Moreover, $\partial_rY = X_r(\theta,r,H_{0}) + X_h(\theta,r,H_{0})\partial_rH_{0}.$ Since $\partial_rH_{0} = -\cot\theta_{\mathrm{eq}},$ substituting this into \eqref{eq:heta-first-derivative}, we obtain
\begin{align}  \label{eq:first-order-implicit_r} \partial_r\mathfrak{h}_\eta = \nabla\eta(Y) \cdot \left( X_r-\cot\theta_{\mathrm{eq}}\,X_h \right).
\end{align}
Similarly, the $\theta$-derivative of $\mathfrak{h}_{\eta}$ can be written as
\begin{equation}  \label{eq:first-order-implicit_theta} \partial_\theta\mathfrak{h}_\eta = \nabla\eta(Y) \cdot (X_\theta).
\end{equation}
Since $X$ is induced by $\Phi$, which depends only on the first spatial derivatives of the steady state function $\zeta$, we obtain the following result from equation \eqref{eq:first-order-implicit_r} and equation \eqref{eq:first-order-implicit_theta} $\|\eta\|_{H^{1}}\lesssim \|\mathfrak{h}_{\eta}\|_{H^{1}} \lesssim \|\eta\|_{H^{1}},$ where the reverse inequality is derived by considering the inverse map.

\textbf{Step 2. Higher-order spatial derivative equivalence.}

Recall that $H_0=-r\cot\theta_{\mathrm{eq}}$ is a linear function depending on $r$. Hence $\partial_r^jH_{0} =0, \; j=2,3,4.$ This implies that $|\partial_{r}^{j}\partial_{\theta}^{j_{2}}\hat{X}|\lesssim |(1+\nabla^{j_{1}+j_{2}}) \zeta|$ for any nonnegative numbers $j_{1}$ and $j_{2}$. We now apply Fa\`{a} di Bruno's formula to $\mathfrak{h}_\eta=\eta\circ \hat{X}.$ For the second derivative,
\begin{equation}  \label{eq:Faa-second} \partial_r^2\mathfrak{h}_\eta = D^2\eta(\hat{X}) [\partial_r\hat{X},\partial_r\hat{X}] + \nabla\eta(Y)\cdot\partial_r^2\hat{X}.
\end{equation}
Similarly, we have
\begin{equation}  \label{eq:Faa-second_theta} \partial_\theta^2\mathfrak{h}_\eta = D^2\eta(\hat{X}) [\partial_\theta\hat{X},\partial_\theta\hat{X}] + \nabla\eta(Y)\cdot\partial_\theta^2\hat{X},
\end{equation}
and
\begin{align}
\partial_\theta\p_{r}\mathfrak{h}_\eta = D^2\eta(\hat{X}) [\partial_\theta\hat{X},\partial_r\hat{X}] + \nabla\eta(Y)\cdot\partial_\theta\partial_{r}\hat{X}.
\end{align}
Therefore, we directly obtain
\begin{align}
{\label{eq:second_order_equivalence}} \|\partial_{r}^{2}\mathfrak{h}_{\eta}\|_{L_{\theta}^{\infty}L^{q}}+\|\partial_{\theta}^{2}\mathfrak{h}_{\eta}\|_{L_{\theta}^{\infty}L^{q}}+\|\partial_{r}\p_{\theta}\mathfrak{h}_{\eta}\|_{L_{\theta}^{\infty}L^{q}} \lesssim \|\eta\|_{W^{3,q}},
\end{align}
for any $q\in[1,\infty]$. When $q\lesssim 2$, we have
\begin{align}
\|\partial_{r}^{2}\mathfrak{h}_{\eta}\|_{L_{\theta}^{\infty}L^{q}}+\|\partial_{\theta}^{2}\mathfrak{h}_{\eta}\|_{L_{\theta}^{\infty}L^{q}}+\|\partial_{r}\p_{\theta}\mathfrak{h}_{\eta}\|_{L_{\theta}^{\infty}L^{q}} \lesssim& \|\eta\|_{H^{3-\beta}}.
\end{align}
Moreover, we obtain $\|\eta\|_{W^{2,b}}\lesssim\|\mathfrak{h}_{\eta}\|_{W^{2,b}}\lesssim \|\eta\|_{W^{2,b}}$ for any $b\in(0,\infty]$.

For the third and higher derivatives, by Fa\`{a} di Bruno's formula and the fact that $X_{r},X_{\theta}$ depend only on the steady state $\zeta$, we obtain $\|\partial_{r}^{3}\mathfrak{h}_{\eta}\|_{L_{\theta}^{\infty}L^{q}} \lesssim \|\eta\|_{W^{3+1,q}},$ and $\|\eta\|_{W^{3,b}}\lesssim\|\mathfrak{h}_{\eta}\|_{W^{3,b}}\lesssim \|\eta\|_{W^{3,b}}.$ By interpolation theory, we derive \eqref{eq:localized-heta-H_high_1}, \eqref{eq:localized-heta-H_high_2} and \eqref{eq:localized-heta-H1}.

\textbf{Step 3. Recovery of the full surface regularity.}

We next prove
\begin{equation}  \label{eq:eta-anisotropic-estimate} \|\eta\|_{W^{3-\frac1p,p}} \lesssim \|\partial_\theta\mathfrak h_\eta\|_ {W_{r,\theta}^{2-\frac1p,p}} + \|\mathfrak h_\eta\|_ {L_\theta^pW_r^{3-\frac1p,p}}.
\end{equation}
We first establish the corresponding anisotropic estimate for $\mathfrak{h}_\eta$. Set $s=3-\frac1p.$ We claim that
\begin{equation}  \label{eq:anisotropic-Sobolev} |\mathfrak{h}_\eta|_{W^{s,p}_{r,\theta}} \lesssim |\partial_\theta\mathfrak{h}_\eta|_ {W^{s-1,p}_{r,\theta}} + |\mathfrak{h}_\eta|_{L_\theta^pW_r^{s,p}}.
\end{equation}

Since the functions under consideration are localized in $r$, we may first extend them from the fixed $r$-interval to $\mathbb{R}$ by means of a bounded extension operator $E$. Hence it suffices to prove \eqref{eq:anisotropic-Sobolev} on $\mathbb{T}_L\times\mathbb{R}.$ For noninteger $s>0$, we use the standard equivalence $W^{s,p} \simeq B^s_{p,p}.$ Let $\Delta_j$ denote the isotropic Littlewood--Paley projection onto the frequency region $|(k,\xi)|\sim 2^j,$ where $k$ is the Fourier frequency corresponding to $\theta$ and $\xi$ is the Fourier frequency corresponding to $r$. On each dyadic annulus, introduce a smooth decomposition $1 = \chi_\theta(k,\xi) + \chi_r(k,\xi),$ such that $\func{supp}\chi_\theta \subset \{ |k|\gtrsim |(k,\xi)| \},$ and $\func{supp}\chi_r \subset \{ |\xi|\gtrsim |(k,\xi)| \}.$ Correspondingly, write $\Delta_j\mathfrak{h}_\eta = \Delta_j^\theta\mathfrak{h}_\eta + \Delta_j^r\mathfrak{h}_\eta.$ We estimate the two pieces separately.

\medskip

\noindent \textit{The $\theta$-dominant region.} On the support of $\Delta_j^\theta$, $|k|\simeq 2^j.$ Hence we may invert one $\theta$-derivative and obtain  $\|\Delta_j^\theta\mathfrak{h}_\eta\|_{L^p} \lesssim 2^{-j} \| \Delta_j^\theta \partial_\theta\mathfrak{h}_\eta \|_{L^p}$.  Multiplying by $2^{js}$ and summing,
\begin{equation}
\sum_{j\geq0} 2^{jsp} \|\Delta_j^\theta\mathfrak{h}_\eta\|_{L^p}^p  \lesssim \sum_{j\geq0} 2^{j(s-1)p} \| \Delta_j^\theta \partial_\theta\mathfrak{h}_\eta \|_{L^p}^p \lesssim \| \partial_\theta\mathfrak{h}_\eta \|_{W^{s-1,p}_{r,\theta}}^p. \label{eq:theta-dominant}
\end{equation}

\medskip

\noindent \textit{The $r$-dominant region.} On the support of $\Delta_j^r$, $|\xi|\simeq 2^j.$ Let $\Delta_j^{(r)}$ denote the one-dimensional Littlewood--Paley projection in the $r$-variable. Since the multiplier defining $\Delta_j^r$ is supported where the $r$-frequency has size $2^j$, the $L^p$-boundedness of the corresponding multipliers gives $\|\Delta_j^r\mathfrak{h}_\eta\|_ {L^p_{\theta,r}} \lesssim \|\Delta_j^{(r)} \mathfrak{h}_\eta\|_ {L^p_{\theta,r}}$. Consequently,
\begin{equation}
\sum_{j\geq0} 2^{jsp} \|\Delta_j^r\mathfrak{h}_\eta\|_{L^p_{\theta,r}}^p \lesssim \sum_{j\geq0} 2^{jsp} \|\Delta_j^{(r)} \mathfrak{h}_\eta\|_{L^p_{\theta,r}}^p \ = \int_{\mathbb{T}_L} \sum_{j\geq0} 2^{jsp} | \Delta_j^{(r)} \mathfrak{h}_\eta(\theta,\cdot) |_{L_r^p}^p ,d\theta \lesssim \|\mathfrak{h}_\eta\|_ {L_\theta^pW_r^{s,p}}^p. \label{eq:r-dominant}
\end{equation}

Combining \eqref{eq:theta-dominant} and \eqref{eq:r-dominant}, and including the low-frequency part, yields
\begin{equation*}
\|\mathfrak{h}_\eta\|_{W^{s,p}_{r,\theta}} \lesssim \|\partial_\theta\mathfrak{h}_\eta\|_ {W^{s-1,p}_{r,\theta}} + \|\mathfrak{h}_\eta\|_ {L_\theta^pW_r^{s,p}}.
\end{equation*}
Since $s=3-\frac1p,$ we obtain
\begin{equation}  \label{eq:heta-anisotropic-final} \|\mathfrak h_\eta\|_ {W^{3-\frac1p,p}_{r,\theta}} \lesssim \|\partial_\theta\mathfrak h_\eta\|_ {W^{2-\frac1p,p}_{r,\theta}} + \|\mathfrak h_\eta\|_ {L_\theta^pW_r^{3-\frac1p,p}}.
\end{equation}
Finally, using \eqref{eq:localized-heta-H1}, we have $\|\eta\|_{W^{3-\frac1p,p}} \lesssim \|\mathfrak{h}_\eta\|_ {W^{3-\frac1p,p}_{r,\theta}}.$ This finishes the proof of the theorem.
\end{proof}

The proof of the previous theorem shows the following norm separation lemma.

\begin{lemma}
{\label{lem:norm_separation}}  Suppose $f\in W^{3-\frac{1}{p},p}$ on the free surface $\Sigma$ in the local chart. We have
\begin{align}
\|f\|_{W^{3-\frac{1}{p},p}}\lesssim\|\partial_{\theta}f\|_{W_{r,\theta}^{2-\frac{1}{p},p}}+\|f\|_{L_{\theta}^{2}W_{r}^{3-\frac{1}{p},p}}.
\end{align}
\end{lemma}

Now we use the basic equivalence theorem to prove the following theorem.

\begin{theorem}
{\label{cor:localized-norm-equivalence}}

Suppose $\mathfrak{h}_{\eta}$ is the surface perturbation in the local chart, and that $\eta$ is supported in $\Sigma\cap\Phi(W)$, as in Theorem~\ref{thm:localized-surface-coordinate-estimate}. The following inequalities for $\mathfrak{h}_{\eta}$ hold
\begin{align}
{\label{eq:cut_off_eq_estimate}} \sum_{j=0}^{k}\|\partial_{\theta}^{m-k}\partial_{t}^{j}\mathfrak{h}_{\eta}\|_{L_{\theta}^{2}H_{r}^{\frac{5}{2}-\beta}}\lesssim\sum_{j=0}^{m}\|\partial_{t}^{j}\eta\|_{H^{\frac{5}{2}+(m-j)-\beta}}\lesssim\sum_{j=0}^{m}\|\partial_{t}^{j}\eta\|_{W^{3+(m-j)-\frac{1}{p},p}},
\end{align}
\begin{align}
{\label{eq:norm_equivalence}} \begin{aligned} \|\partial_{t}^{j}\eta\|_{W^{a,b}}\lesssim \|\partial_{t}^{j}\mathfrak{h}_{\eta}\|_{W^{a,b}}\lesssim\ \|\partial_{t}^{j}\eta\|_{W^{a,b}}, \end{aligned}
\end{align}
where $a\in(0,3)$, $b\in [1,2]$ and $0\leq j\lesssim m+2$.
\end{theorem}

\begin{proof}
Note that $\partial_{t}\eta(X(\theta,r,H_{0}(\theta,r)),t)=\partial_{t}\mathfrak{h}_{\eta}(\theta,r,t)$. The inequality \eqref{eq:norm_equivalence} is proved by replacing $\eta$ and $\mathfrak{h}_{\eta}$ in the proof of the previous theorem by $\partial_{t}^{j}\eta$ and $\partial_{t}^{j}\mathfrak{h}_{\eta}$. Inequality \eqref{eq:cut_off_eq_estimate} is proved in the same way. We omit the details here.
\end{proof}

Now, we have defined the surface function $\mathfrak{h}_{\eta}$ in the local charts and established the corresponding norm equivalence theorem. For simplicity, in the following estimates, we still use $\eta$ instead of $\mathfrak{h}_{\eta}$ even in the local charts to keep the notation unified.

\section{Extension of the surface function}

\label{sec:extension}

\subsection{Half-space extender}

Now we have the diffeomorphism and the corresponding local charts. With this flattening map, we can define the extension along one direction, as in \cite{GT}. We give the following definition.

\begin{theorem}
\label{thm:half-space-extender} Let $\Omega=\mathbb{T}\times\mathbb{R}, \; \Omega_{+}=\mathbb{T}\times(0,\infty).$ Fix $1\leq k\in\mathbb{N}$. Then there exists a bounded linear map $E:H^{0}(\Omega_{+})\to H^{0}(\Omega)$ satisfying the following properties.

\begin{enumerate}
\item For $0\leq s\leq k$, we have $E\bigl(H^{s}(\Omega_{+})\bigr)\subseteq H^{s}(\Omega),$  and the induced restriction $E:H^{s}(\Omega_{+})\to H^{s}(\Omega)$  is bounded and linear. More generally, $E:W^{s,p}(\Omega_{+})\to W^{s,p}(\Omega)$ is bounded for $0\leq s\leq k$ and $1<p<\infty$.

\item $\left.Ef\right|_{\Omega_{+}}=f$  for every $f\in H^{0}(\Omega_{+})$.

\item $[E,\partial_{1}]f=0$  for every $f\in H^{1}(\Omega_{+})$.
\end{enumerate}
\end{theorem}

\begin{proof}
Since Vandermonde matrices are invertible, we can choose $a_{1},\dots,a_{k+1}\in\mathbb{R}$ such that
\begin{equation}  \label{eq:Vandermonde-condition} \sum_{j=1}^{k+1} \left(-\frac{1}{j}\right)^{i}a_{j} =1, \qquad i=0,\dots,k.
\end{equation}
Given $f\in C^{k}(\Omega_{+})$, define $Ef:\Omega\to\mathbb{R}$ via
\begin{equation}  \label{eq:half-space-extension} Ef(x) =
\begin{cases}
f(x), & x_{2}>0, \\[1ex]
\displaystyle \sum_{j=1}^{k+1} a_{j} \chi_{\alpha }(x_{1},x_{2})f\left(x_{1},-\frac{x_{2}}{j}\right), & x_{2}<0,\end{cases}\end{equation} where $\chi_{\alpha}=\chi_{\alpha}(x_{2})$ is a smooth cutoff function of $x_{2}$ only, which equals $1$ in a small neighborhood of $\{x_{2}=0\}$ and satisfies, for any $j$, $\chi_{\mathrm{cl}}(x_{1},-\frac{x_{2}}{j})=1,$ for any $(x_{1},x_{2})\in \func{supp}\{\chi_{\alpha}\}$.

Routine computations then show that, for $x\in\Omega_{-}$ and $\lvert\alpha\rvert\leq k$,
\begin{equation}  \label{eq:derivative-half-space-extension} \partial^{\alpha}Ef(x) = \sum_{m=1}^{k+1} \left(-\frac{1}{m}\right)^{\alpha_{2}} a_{m} \partial^{\alpha}f \left(x_{1},-\frac{x_{2}}{m}\right),
\end{equation}
where $\Omega_{-}:=\mathbb{T}\times(-\infty,0).$ Formula \eqref{eq:derivative-half-space-extension} holds where $\chi_{\alpha}=1$; elsewhere the derivatives of $\chi_{\alpha}$ produce additional terms of lower order in $f$. Since $\chi_{\alpha}$ depends on $x_{2}$ only, it commutes with $\partial_{1}$, which gives property (3). Moreover,
\begin{align}
\partial^{\alpha}Ef(x_{1},0) &= \partial^{\alpha}f(x_{1},0) \sum_{m=1}^{k+1} \left(-\frac{1}{m}\right)^{\alpha_{2}}a_{m}= \partial^{\alpha}f(x_{1},0) \label{eq:matching-derivatives}
\end{align}
for $x_{1}\in\mathbb{T}$. We readily deduce from this that $Ef\in C^{k}(\Omega).$ Moreover, by construction, we have $Ef=f$ in $\Omega_{+}$,  $\partial_{1}Ef=E\partial_{1}f$, and $\|Ef\|_{H^{s}(\Omega)} \lesssim C(\varepsilon)\|f\|_{H^{s}(\Omega_{+})}$ for every $s\in[0,k]\cap\mathbb{N}$, and in the same way $\|Ef\|_{W^{s,p}(\Omega)}\lesssim C(\varepsilon)\|f\|_{W^{s,p}(\Omega_{+})}$ for $s\in[0,k]\cap\mathbb{N}$ and $1<p<\infty$, by \eqref{eq:derivative-half-space-extension}. The fractional orders follow by interpolation. We may then use continuous linear extension to obtain the desired operator.
\end{proof}

Throughout we fix $k:=m+3$ in Theorem~\ref{thm:half-space-extender}. This covers every order used below; the highest one is $W^{m+3-\frac1p,p}$, in $\mathcal{D}_{\mathrm{int}}$.

\subsection{Free surface function extension}

\label{subsec:extension_global}

We view $\eta:\Sigma\to\mathbb{R},$ where $\Sigma$ is the equilibrium free surface. We decompose $\Sigma=\Sigma_{\mathrm{int}}\cup\Sigma_{\mathrm{cl}},$ where $\Sigma_{\mathrm{cl}}$ is the portion of $\Sigma$ near the contact line corresponding to the region that can be reached with our global diffeomorphism. Using this decomposition and smooth cutoffs, we can split $\eta=\chi_{\mathrm{int},\Sigma}\eta+(1-\chi_{\mathrm{int},\Sigma})\eta,$ where $\chi_{\mathrm{int},\Sigma}$ is a cutoff function in the 2-D $(\theta,r)-$plane that vanishes in a neighborhood of the contact line and depends only on the variable $r$ defined above.

We first handle the interior part in a very simple way. We view a set $D_{\mathrm{int}}\subset D$, where $D\subset\mathbb{R}^{2}$ is the parameterization domain for the equilibrium, as lying above the entirety of the fluid body; that is, it is translated above the cylinder in which the free surface lies. We define
\begin{align}
\bar{\eta}_{\mathrm{int}}(x,y,z,t):=\chi_{\mathrm{int}+}\mathcal{P}\circ E_{D_{\mathrm{int}}}[\eta_{\mathrm{int}}\circ \pi](x,y,z-E\zeta(x,y),t),
\end{align}
where:

\begin{enumerate}
\item $\pi$ is the projection from $D_{\mathrm{int}}$ down to  $\Sigma_{\mathrm{int}}$;

\item $E_{D_{\mathrm{int}}}$ is a standard Sobolev extension  operator; Stein's extension operator will work, although it is  more than is needed here;

\item $P$ is the standard downward Poisson extension.

\item On the boundary $\Sigma$, we have $\bar{\eta}_{j,\mathrm{int}}=\chi_{j,\mathrm{int},\Sigma}\eta.$ \item $\chi_{\mathrm{dif,j}}$, $\chi_{\mathrm{dif,j+}}$ and $\chi_{\mathrm{0,int}}$ are cutoff functions in 3-D defined in Subsection~\ref{subsec:cutoff_extension}.
\end{enumerate}

Thus, this is the easy part. The more complicated part uses our diffeomorphism $\Phi:W\to\Phi(W)\subset\Omega.$ We view $\eta_{\mathrm{cl}}\circ\Phi$ as being defined on the flat top portion of the canonical dihedron. This function is compactly supported, so we can regard it as a function on the half-space $\mathbb{T}\times(0,\infty)$, as in the previous subsection. We then define
\begin{equation}  \label{eq:contact-line-extension} \overline{\eta}_{\mathrm{cl}} = \left[ (1-\chi_{m+1-,\mathrm{int}}) P_{2-D}\circ E\bigl(\eta_{\mathrm{cl}}\circ\Phi\bigr)(\theta,r,h+r\cot\theta_{\mathrm{eq}},t) \right]\circ\Phi^{-1},
\end{equation}
where $E$ is the extension operator constructed above, interpreted with the top boundary lying at a slant; $P_{2-D}$ is the downward Poisson extension for any fixed $\theta$, normal to the extended top $\mathbb{T}\times\mathbb{R}$. Using the definition of the cutoff functions, we have $\left.\bar\eta_{\mathrm{cl}}\right|_\Sigma =(1-\chi_{m+1-,\mathrm{int},\Sigma}) \chi_{\mathrm{cl},\Sigma}\eta =\chi_{\mathrm{cl},\Sigma}\eta$.

Similarly, we define $\chi_{\mathrm{bet}+}=\chi_{\mathrm{cl}}+\chi_{\mathrm{bet}}$ and $\chi_{m+1,\mathrm{int},\Sigma}=1-\chi_{\mathrm{cl},\Sigma}$. We note that
\begin{align}
(1-\chi_{cl,\Sigma}-\chi_{\mathrm{int},\Sigma})\eta=\sum_{j=0}^{m}(\chi_{j+1,int,\Sigma}\eta -\chi_{j,int,\Sigma}\eta)+\chi_{\mathrm{0},\mathrm{int},\Sigma}\eta-\chi_{\mathrm{int},\Sigma}\eta .
\end{align}
We then introduce the notation $\chi_{\func{dif},j,\Sigma}=\chi_{j+1,int,\Sigma}-\chi_{j,int,\Sigma},$ and $\chi_{\func{dif},j+,\Sigma}=\chi_{j+1+,int,\Sigma}-\chi_{j-,int,\Sigma}.$ The key point is that we can define for any $0\leq k\leq m$
\begin{equation}
\begin{aligned}
&\label{eq:interior-extension} \overline{\eta}_{k+1,\mathrm{int}}(x,y,z,t) = \sum_{j=0}^{k}\bigl[\chi_{\mathrm{diff},j+}P_{2-D}\circ E \bigl(\eta_{\mathrm{dif},j,\Sigma}\circ\Phi\bigr)(\theta,r,h+r\cot\theta_{\mathrm{eq}},t)\bigr]\circ{\Phi}^{-1} \\
&\quad +\bigl[(\chi_{0+,int}-\chi_{\mathrm{int}-})P_{2-D}\circ E\bigl[(\eta_{0,\mathrm{int},\Sigma}-\eta_{{\mathrm{int},\Sigma}})\circ\Phi\bigr](\theta,r,h+r\cot\theta_{\mathrm{eq}},t)\bigr]\circ \Phi^{-1}+\bar{\eta}_{\mathrm{int}},
\end{aligned}
\end{equation}
and $\bar{\eta}_{j,\mathrm{int},d}=\overline{\eta_{{j,\mathrm{int}}}-\eta_{\mathrm{int}}}:=\bar{\eta}_{j,\mathrm{int}}-\bar{\eta}_{\mathrm{int}}.$ \emph{Remark.} The definition above implies a crucial property of this sequence of cutoff functions. We have $\chi_{j_{1}-}\chi_{\mathrm{dif},j_{2}+}=0$ if $j_{2}-j_{1}\geq 1$. Specifically, $\chi_{cl}\chi_{j,\mathrm{int}}=0$ if $j\leq m$. This implies that the interior elliptic estimate does not depend on the elliptic estimate in the contact line part. In other words $\chi_{j,int}\mathcal{A}_{\eta}$ does not depend on $\chi_{cl,\Sigma}\eta$.

Finally, we define $\bar{\eta}$ as follows
\begin{align}
\bar{\eta}:=\bar{\eta}_{m+1,\mathrm{int}}+\bar{\eta}_{\mathrm{cl}}. \label{eq:eta_bar_def}
\end{align}
We now prove the key property that $\eta$ must satisfy.

\begin{theorem}
\label{thm:extension_trace}Suppose $\bar{\eta}$ is the extension function of $\eta$. Then $\bar{\eta}=\eta$ on $\Sigma$.
\end{theorem}

\begin{proof}
By definition $\bar{\eta}_{\mathrm{int}}=\chi_{\mathrm{int},\Sigma}\eta$ on the surface $\Sigma$. Then it suffices to prove that $\bar{\eta}_{\mathrm{cl}}=\chi_{\mathrm{cl},\Sigma}\eta$ on the equilibrium surface $\Sigma$. The proof of $\bar{\eta}_{m+1,\mathrm{int}}-\chi_{\mathrm{int}}=\eta_{m+1,\mathrm{int}}-\chi_{\mathrm{int}}$ follows in the same way. For any fixed time $t$ and any point $x\in\Sigma$ near the contact line, there exist $\theta,r$ such that $x=X(\theta,r,H_{0}(\theta,r))$. Therefore $\bar{\eta}_{\mathrm{cl}}(x)|_{x\in\Sigma}=\bar{\mathfrak{h}} _{\eta}(\theta_{x},r_{x},H_{0}(\theta,r))=\mathfrak{h}_{\eta}( \theta_{x},r_{x})=\chi_{\mathrm{cl}}\eta(X(\theta_{x},r_{x},H_{0}))=\chi_{ \mathrm{cl}}\eta(x)$. This finishes the proof.
\end{proof}

\subsection{The commutator and the bootstrap of elliptic regularity.}

Let $\bar{E}:=P_{2-D}\circ E.$ We assume that $\chi_{j,int,\Sigma}\prec \chi_{j+1,int,\Sigma},$ in the sense that $\chi_{j+1,\mathrm{int},\Sigma}\equiv 1 \; \text{in a neighborhood of } \func{supp}\chi_{j,\mathrm{int},\Sigma},$ and we also assume that $\chi_{j,int}\prec \chi_{j+1,int}.$ We then state the following theorem.

\begin{theorem}
{Enhanced regularity of the commutator in the intermediate area.} Fix the cutoff width $\varepsilon>0$. Assume
\begin{equation}
\|\nabla(W\bar\eta)\|_{L^\infty(\Omega)}\le\frac12.
\end{equation}\label{thm:extension-commutator-estimate} The extension function $\bar{\eta}_{int}$ of $\eta_{int}$ satisfies the following estimate
\begin{equation}
\begin{aligned}
&\label{eq:extension-commutator-estimate} \|\chi_{j,int}\bar{\eta}_{m,int}\|_{W^{m-j+2,p}}+\|\chi_{j+,int}\bar{\eta}_{m,int}\|_{W^{m-j+2,p}} \lesssim C(\varepsilon)\sum_{k=0}^{j}\|P_{2-D}\circ E(\chi_{\mathrm{diff},k,\Sigma}\eta)\|_{W^{m-(j)+2,p}} \\
&\quad +C(\varepsilon)\|P\circ E_{D_{\mathrm{int}}}(\chi_{\mathrm{int},\Sigma}\eta)\|_{W^{m-(j+1)+3,p}}+C(\varepsilon)\|P_{2-D}\circ E((\chi_{0,\mathrm{int},\Sigma}-\chi_{\mathrm{int},\Sigma})\eta)\|_{W^{m-j+2,p}},
\end{aligned}
\end{equation}
for any $j\leq m-1$. This directly implies the following estimate
\begin{equation}
\begin{aligned}
&\label{eq:extension-commutator-estimate_0} \|\chi_{j+,int}(I-\mathcal{A}_{\eta})\|_{W^{m-j+1,p}} \lesssim C(\varepsilon)(\sum_{k=0}^{j}\|P_{2-D}\circ E(\chi_{\mathrm{diff},k,\Sigma}\eta)\|_{W^{m-(j)+2,p}} \\
&\quad +\|P\circ E_{D_{\mathrm{int}}}(\chi_{\mathrm{int},\Sigma}\eta)\|_{W^{m-(j+1)+3,p}} +\|P_{2-D}\circ E((\chi_{0,\mathrm{int},\Sigma}-\chi_{\mathrm{int},\Sigma})\eta)\|_{W^{m-j+2,p}}),
\end{aligned}
\end{equation}
for any positive integer $j\leq m-1$.
\end{theorem}

\begin{proof}
It suffices to prove that
\begin{equation}
\begin{aligned}
&\|\chi_{j+,int}\bar{\eta}_{m,int}\|_{W^{m-j+2,p}} \lesssim C(\varepsilon)(\sum_{k=0}^{j}\|P_{2-D}\circ E(\chi_{\mathrm{diff},k,\Sigma}\eta)\|_{W^{m-(j)+2,p}} \\
&\quad +\|P\circ E_{D_{\mathrm{int}}}(\chi_{\mathrm{int},\Sigma}\eta)\|_{W^{m-(j+1)+3,p}} +\|P_{2-D}\circ E((\chi_{0,\mathrm{int},\Sigma}-\chi_{\mathrm{int},\Sigma})\eta)\|_{W^{m-j+2,p}}).
\end{aligned}
\end{equation}
For simplicity, in the following computation, we use $P_{2-D}\circ E$ to denote the local extension of a surface function $f$ given by
\begin{equation*}
\left[(P_{2-D}\circ E)(f\circ\widehat X) \bigl(\theta,r,h+r\cot\theta_{\mathrm{eq}},t\bigr)\right]\circ\Phi^{-1}.
\end{equation*}
The bulk cutoff multiplying this extension is retained explicitly. By the definition of $\bar{\eta}_{m,int}$ given by \eqref{eq:interior-extension}, we have
\begin{equation*}
\overline{\eta}_{m,\mathrm{int}} = \sum_{k=0}^{m-1}\chi_{\mathrm{diff},k+}P_{2-D}\circ E \bigl[\eta_{\mathrm{dif},k,\Sigma}\bigr]+(\chi_{0+,int}-\chi_{\mathrm{int}-})P_{2-D}\circ E \bigl[(\eta_{0,\mathrm{int},\Sigma}-\eta_{\mathrm{int},\Sigma})\bigr]+\chi_{j+,int}\bar{\eta}_{\mathrm{int}}.
\end{equation*}
Therefore, multiplying $\bar{\eta}_{int}$ by $\chi_{j+,\mathrm{int}}$, we obtain
\begin{align}
\chi_{j+,\mathrm{int}} \overline{\eta}_{m,\mathrm{int}}=& \sum_{k=0}^{m-1}\chi_{j+,int}\chi_{\mathrm{diff},k+}P_{2-D}\circ E\bigl[\eta_{\mathrm{dif},k,\Sigma}\bigr] \\
&+\chi_{j+,int}(\chi_{0+,int}-\chi_{\mathrm{int}-})P_{2-D}\circ E\bigl[(\eta_{0,\mathrm{int},\Sigma}-\eta_{\mathrm{int},\Sigma})\bigr]+\chi_{j+,int}\bar{\eta}_{\mathrm{int}}.
\end{align}

For smooth cutoffs $a,b$, write $a\prec b$ if $b=1$ on a neighborhood of $\func{supp}a$. By the definition of $\chi_{j+,int}$, we have $\chi_{j+,int}\prec \chi_{j+1-,int}.$ Therefore $\chi_{j+,int}\chi_{k+,\mathrm{diff},int}=0$ for any integer $k\geq j+1$. Hence,
\begin{equation}
\begin{aligned}
&\chi_{j+,\mathrm{int}} \overline{\eta}_{m,\mathrm{int}}= \sum_{k=0}^{j}\chi_{j+,int}\chi_{\mathrm{diff},k+}P_{2-D}\circ E\bigl[\eta_{\mathrm{dif},k,\Sigma}\bigr] \\
&\quad +\chi_{j+,int}(\chi_{0+,int}-\chi_{\mathrm{int}-})P_{2-D}\circ E \bigl[(\eta_{0,\mathrm{int},\Sigma}-\eta_{\mathrm{int},\Sigma})\bigr] +\chi_{j+,int}\bar{\eta}_{\mathrm{int}}.
\end{aligned}
\end{equation}
Therefore, we have
\begin{align}
\begin{aligned}
\|\chi_{j+,int}\bar{\eta}_{m,int}\|_{W^{m-j+2,p}}\lesssim&C(\varepsilon)\sum_{k=0}^{j}\|P_{2-D}\circ E(\chi_{\mathrm{diff},k,\Sigma}\eta)\|_{W^{m-(j)+2,p}}+C(\varepsilon)\|P\circ E_{D_{\mathrm{int}}}(\chi_{\mathrm{int},\Sigma}\eta)\|_{W^{m-j+2,p}}\\ &+C(\varepsilon)\|P_{2-D}\circ E(\chi_{0,\mathrm{int},\Sigma}\eta-\chi_{\mathrm{int},\Sigma}\eta)\|_{W^{m-j+2,p}}. \end{aligned}
\end{align}
This finishes the proof.
\end{proof}

\emph{Remark.} This is one of the key lemmas of the paper. We have to exchange the extension operator and the cutoff function. This theorem tells us that when we test $\bar{\eta}$ by $\chi_{j,\mathrm{int}}$, its value does not depend on the value of $(1-\chi_{j+1,\mathrm{int},\Sigma})\eta$. Moreover, since we did not use any property of the Sobolev norm $W^{{m-j+2,p}}$, we can change it to any Sobolev norm of any $t$-derivative of $\eta$ since $P\circ{E}_{D_{\mathrm{int}}},P_{2-D}\circ E$ and cutoff functions are all independent of $t$.

We now state another enhanced commutator estimate.

\begin{theorem}
\label{thm:enhanced_commutator} Let $\bar{E}_{\rho}:=P_{2-D}\circ E$, acting slice by slice as in \eqref{eq:contact-line-extension}. Let $\psi=\psi(r,h)$ be smooth and supported in the chart, and let $f$ be supported in the part of $\Sigma$ covered by the chart and vanish where $\chi_{\mathrm{cl},\Sigma}=1$. Then, for every $1<q<\infty$,
\begin{align}
\big\|\psi\,\bar{E}_{\rho}f-\bar{E}_{\rho}(\psi_{\Sigma}f)\big\|_{L^{2}W^{4,q}}\leq C_{\psi}\,\|f\|_{L^{2}W^{3-\frac{1}{q},q}}, \label{eq:enhanced_commutator}
\end{align}
where $\psi_{\Sigma}$ is the restriction of $\psi$ to $\Sigma$ and $C_{\psi}$ depends on finitely many derivatives of $\psi$; for the cutoffs of Subsection~\ref{subsec:cutoff}, $C_{\psi}=C(\varepsilon)$. Thus the commutator needs one spatial derivative less than $\psi\bar{E}_{\rho}f$ itself. In particular, for $f:=(\chi_{m+1,\mathrm{int},\Sigma}-\chi_{\mathrm{int},\Sigma})\eta$ and $f^{\prime }:=\chi_{m,\mathrm{int},\Sigma}f=(\chi_{m,\mathrm{int},\Sigma}-\chi_{\mathrm{int},\Sigma})\eta$,
\begin{align}
\big\|\chi_{m,\mathrm{int}}\bar{E}_{\rho}f-\chi_{m+,\mathrm{int}}\bar{E}_{\rho}f^{\prime }\big\|_{L^{2}W^{4,q}}\leq C(\varepsilon)\,\|f\|_{L^{2}W^{3-\frac{1}{q},q}}.  \label{eq:enhanced_commutator_m}
\end{align}
\end{theorem}

\begin{proof}
Fix $\theta$ and work on the slice in the variables $(r,\mathscr{H})$, $\mathscr{H}=h+r\cot\theta_{\mathrm{eq}}$. There $\bar{E}_{\rho}g$ is the harmonic extension to the half-plane $\mathbb{H}:=\{\mathscr{H}<0\}$ of the function $Eg$ on $\{\mathscr{H}=0\}$, and $\Delta$ is the Laplacian in $(r,\mathscr{H})$. The function $\psi$ is smooth in $(r,\mathscr{H})$, and its restriction to $\{\mathscr{H}=0\}$ is $\psi_{\Sigma}$. Since $f$ vanishes where $\chi_{\mathrm{cl},\Sigma}=1$, which contains the reflected points used by $E$ (see \eqref{eq:half-space-extension}), $Ef$ and $E(\psi_{\Sigma}f)$ are the extensions of $f$ and $\psi_{\Sigma}f$ by zero, and $E(\psi_{\Sigma}f)=\psi_{\Sigma}Ef$.

\emph{The commutator.} Put $\mathfrak{M}_{1}:=\psi\bar{E}_{\rho}f-\bar{E}_{\rho}(\psi_{\Sigma}f)$. On $\{\mathscr{H}=0\}$ both terms equal $\psi_{\Sigma}Ef$, so $\mathfrak{M}_{1}=0$ there. Since $\bar{E}_{\rho}f$ and $\bar{E}_{\rho}(\psi_{\Sigma}f)$ are harmonic in $\mathbb{H}$, $\Delta\mathfrak{M}_{1}=(\Delta\psi)\,\bar{E}_{\rho}f+2\nabla\psi\cdot\nabla \bar{E}_{\rho}f$  in $\mathbb{H}$.  The right-hand side is supported in $\func{supp}\nabla\psi$, a bounded set, and the Poisson extension maps $W^{3-\frac1q,q}(\mathbb{R})$ into $W^{3,q}$ of bounded subsets of $\mathbb{H}$. Hence $\|\Delta\mathfrak{M}_{1}\|_{W^{2,q}}\leq C_{\psi}\|f\|_{W^{3-\frac1q,q}(\mathbb{R})}$. The boundary $\{\mathscr{H}=0\}$ is flat, so local $W^{4,q}$ regularity for the Dirichlet problem, on half-balls covering the bounded region where the left-hand side of \eqref{eq:enhanced_commutator} is measured, gives $\|\mathfrak{M}_{1}\|_{W^{4,q}}\leq C\big(\|\Delta\mathfrak{M}_{1}\|_{W^{2,q}}+\|\mathfrak{M}_{1}\|_{L^{q}}\big)\leq C_{\psi}\|f\|_{W^{3-\frac1q,q}(\mathbb{R})}$. No boundary condition on the wall is needed: after the extension $E$ the slice problem lives on the whole half-plane. The constants do not depend on $\theta$; squaring and integrating in $\theta$ proves \eqref{eq:enhanced_commutator}.

\emph{The particular case.} Since $\chi_{m+1,\mathrm{int},\Sigma}=1$ on $\func{supp}\chi_{m,\mathrm{int},\Sigma}$ and $\chi_{m,\mathrm{int},\Sigma}=1$ on $\func{supp}\chi_{\mathrm{int},\Sigma}$, we have $\chi_{m,\mathrm{int},\Sigma}f=f^{\prime }$, and $f$ vanishes where $\chi_{\mathrm{cl},\Sigma}=1$. Hence
\begin{equation*}
\chi_{m,\mathrm{int}}\bar{E}_{\rho}f-\chi_{m+,\mathrm{int}}\bar{E}_{\rho}f^{\prime }=\big(\chi_{m,\mathrm{int}}\bar{E}_{\rho}f-\bar{E}_{\rho}(\chi_{m,\mathrm{int},\Sigma}f)\big)+(1-\chi_{m+,\mathrm{int}})\bar{E}_{\rho}f^{\prime }.
\end{equation*}
The first term is bounded by \eqref{eq:enhanced_commutator} with $\psi=\chi_{m,\mathrm{int}}$. For the second, $f^{\prime }$ is supported in $\func{supp}\chi_{m,\mathrm{int},\Sigma}$, and $\chi_{m+,\mathrm{int}}=1$ on a neighbourhood of $\func{supp}\chi_{m,\mathrm{int}}$. So the distance $d$ from $\func{supp}(1-\chi_{m+,\mathrm{int}})$ to $\func{supp}f^{\prime }\subset\{\mathscr{H}=0\}$ is positive, and $d\geq c\,\varepsilon$ by \eqref{eq:cutoff_scaling}. The Poisson kernel of $\mathbb{H}$ is smooth away from the support of the data: on $\func{supp}(1-\chi_{m+,\mathrm{int}})$ every derivative of order $k$ of $\bar{E}_{\rho}f^{\prime }$ is bounded by $C_{k}d^{-k-1}\|f^{\prime }\|_{L^{1}}$. Hence $\|(1-\chi_{m+,\mathrm{int}})\bar{E}_{\rho}f^{\prime }\|_{W^{4,q}}\leq C(\varepsilon)\|f\|_{L^{q}}$, which proves \eqref{eq:enhanced_commutator_m}.
\end{proof}

\section{Whole space energy estimate}

\label{sec:whole_space}

\subsection{Transformation from the moving domain to the equilibrium domain}

We first formulate the problem in the original Cartesian coordinates. Let  $D\subset \mathbb{R}^{2}$ be the horizontal cross-section of the vessel. The equilibrium free surface is given by
\begin{equation*}
\Sigma = \left\{ (x,z)\in \mathbb{R}^{2}\times\mathbb{R }\;\middle|\; x\in D,\quad z=\zeta(x) \right\},
\end{equation*}
while the perturbed free surface is
\begin{equation*}
\Sigma(t) = \left\{ (x,z)\in \mathbb{R}^{2}\times\mathbb{R }\;\middle|\; x\in D,\quad z=\zeta(x)+\eta(x,t) \right\}.
\end{equation*}
We denote the equilibrium fluid domain by $\Omega$ and the moving fluid domain by $\Omega(t)$. The corresponding equilibrium and moving contact lines are  $\Gamma = \left\{ (x,\zeta(x)) \;\middle|\; x\in\partial D \right\}$, and $\Gamma(t) = \left\{ (x,\zeta(x)+\eta(x,t)) \;\middle|\; x\in\partial D \right\}$.

We use the convention
\begin{equation*}
H(f) := \func{div}_{x} \left( \frac{\nabla_x f} {\sqrt{1+|\nabla_x f|^{2}}} \right),
\end{equation*}
and define  $\mathbb{D }u := \nabla u+(\nabla u)^{T}, \; S(p,u) := pI-\mu\mathbb{D }u$. The upward-pointing non-unit normal to the moving free surface is
\begin{equation*}
N_{\zeta+\eta} =
\begin{pmatrix}
-\nabla_x(\zeta+\eta) \\
1\end{pmatrix}.
\end{equation*}

The perturbation formulation in the moving domain is
\begin{equation}  \label{eq:moving-domain-full-system} \left\{ \begin{alignedat}{2} &\partial_tu + u\cdot\nabla u + \operatorname{div} S(p,u) =0 &&\qquad \text{in }\Omega(t), \\
&\operatorname{div}u =0 &&\qquad \text{in }\Omega(t), \\ &S(p,u)\mathcal N = \left\{ g\eta -\sigma \operatorname{div}_{x,y} \left[ \frac{\nabla_{x,y}(\zeta+\eta)} {\sqrt{1+|\nabla_{x,y}(\zeta+\eta)|^2}} - \frac{\nabla_{x,y}\zeta} {\sqrt{1+|\nabla_{x,y}\zeta|^2}} \right] \right\} \mathcal N &&\qquad \text{on }\Sigma(t), \\ &\partial_t\eta = u\cdot\mathcal N &&\qquad \text{on }\Sigma(t), \\ &u\cdot n_0 =0 &&\qquad \text{on }\Sigma_s, \\ &\bigl( S(p,u)n_0 \bigr)_{\tan} -\beta_0u_{\tan} =0 &&\qquad \text{on }\Sigma_s, \\ &-{\kappa \partial_t\eta}-\hat{\mathcal{W}}(\partial_{t}\eta) = \sigma \left[ \frac{\partial_\nu(\zeta+\eta)} {\sqrt{1+|\nabla_{x,y}(\zeta+\eta)|^2}} - \frac{\partial_\nu\zeta} {\sqrt{1+|\nabla_{x,y}\zeta|^2}} \right] &&\qquad \text{on }\Gamma, \\ &u(\cdot,0)=u_0^{\sharp}, \\ &\eta(\cdot,0)=\eta_0.
\end{alignedat} \right.
\end{equation}
In the equation above, $\mathcal{W}(z)$ is an increasing smooth function of $z$ such that $\mathcal{W}(0)=0$ and $\mathcal{W}^{\prime }(0)>0$. We define
\begin{align}
\hat{\mathcal{W}}(z)=\mathcal{W}(z)-\kappa z,  \label{eq:W_hat}
\end{align}
where $\kappa=\mathcal{W}^{\prime }(0)>0$. Therefore $\hat{\mathcal{W}}(\partial_{t}\eta)\sim O((\partial_{t}\eta)^{2})$.

\subsubsection{Construction of the global diffeomorphism}

We now construct a diffeomorphism from the equilibrium domain $\Omega$ onto the moving domain $\Omega(t)$. For simplicity, suppose that the lower part of the vessel has been normalized so that the equilibrium free surface lies strictly above it. Choose a fixed smooth function  $\phi\in C^\infty(\mathbb{R})$ such that  $\phi(z)=0 \;\text{near the fixed lower boundary}$, and $\phi(z)=z$ for all $z$  in a neighborhood containing the range of $\zeta$. We choose coordinates such that $\zeta(x)>c>0$ for a constant $c$ and any $x$. Define
\begin{equation}  \label{eq:W-definition} W(x,z) := \frac{\phi(z)}{\zeta(x)}.
\end{equation}
Then  $W(x,\zeta(x))=1$. More generally, one may take any $W\in C^\infty(\overline\Omega)$ so that $W|_\Sigma=1, \; W=0\;\text{on }\overline\Omega\cap\{z\le h_b\}$, where the wall is vertical above $h_b$. We define
\begin{equation}  \label{eq:global-diffeomorphism} \Phi_\eta(x,z,t) :=
\begin{pmatrix}
x \\
z+W(x,z)\bar\eta(x,z,t)\end{pmatrix}.
\end{equation}
Thus  $\Phi_\eta(t): \Omega \longrightarrow \Omega(t)$. Its inverse  $\Psi_\eta(t) := \Phi_\eta(t)^{-1} : \Omega(t) \longrightarrow \Omega$ is the desired map from the moving domain to the equilibrium domain.

On the equilibrium free surface $z=\zeta(x)$, we have  $q(x,\zeta(x),t) = W(x,\zeta(x)) \bar\eta(x,\zeta(x),t) = \eta(x,t)$, and hence  $\Phi_\eta(x,\zeta(x),t) = \left( x,\zeta(x)+\eta(x,t) \right)$. Therefore,  $\Phi_\eta(\Sigma)=\Sigma(t)$. Since the horizontal coordinate $x$ is unchanged, $\Phi_\eta$ maps the solid side wall into itself. Moreover, since $W=0$ near the fixed lower boundary,  $\Phi_\eta=\func{Id}$ near the lower boundary. If $\eta$ is sufficiently small so that, for example,  $\|\nabla q\|_{L^\infty(\Omega)} \ll1$, then $\Phi_\eta$ is a diffeomorphism from $\Omega$ onto $\Omega(t)$.

\subsubsection{Geometric quantities associated with the diffeomorphism}

The deformation gradient is
\begin{equation*}
\nabla\Phi_\eta =
\begin{pmatrix}
1 & 0 & 0 \\
0 & 1 & 0 \\
\partial_1q & \partial_2q & 1+\partial_3q\end{pmatrix}.
\end{equation*}
Define
\begin{equation}  \label{eq:JKA-definition} J_{\eta} := \det\nabla\Phi_\eta = 1+\partial_3q, \qquad K_{\eta} := J_{\eta}^{-1}, \qquad \mathcal{A}_\eta := (\nabla\Phi_\eta)^{-T}.
\end{equation}
Then
\begin{equation}  \label{eq:A-explicit} \mathcal A_\eta = \begin{pmatrix} 1&0&-K\partial_1q\\
0&1&-K\partial_2q\\ 0&0&K \end{pmatrix}.
\end{equation}

For a scalar function $f$, we define  $(\nabla_{\mathcal{A}_\eta}f)_i := (\mathcal{A}_\eta)_{ij}\partial_j f$. For a vector field $v$, we define  $(\nabla_{\mathcal{A}_\eta}v)_{ij} := (\mathcal{A}_\eta)_{jk}\partial_k v_i$, $\func{div}_{\mathcal{A}_\eta}v := (\mathcal{A}_\eta)_{ij}\partial_jv_i$, and, for a matrix field $M$,  $(\func{div}_{\mathcal{A}_\eta}M)_i := (\mathcal{A}_\eta)_{jk}\partial_kM_{ij}$. We further define  $\mathbb{D}_{\mathcal{A}_\eta}v := \nabla_{\mathcal{A}_\eta}v + (\nabla_{\mathcal{A}_\eta}v)^T$, and
\begin{equation}  \label{eq:geometric-stress} S_{\mathcal{A}_\eta}(p,v) := pI-\mu\mathbb{D}_{\mathcal{A}_\eta}v.
\end{equation}

The non-unit normal to the equilibrium free surface is
\begin{equation*}
N_\zeta =
\begin{pmatrix}
-\nabla_x\zeta \\
1\end{pmatrix}.
\end{equation*}
Define the transformed free-surface normal by
\begin{equation}  \label{eq:transformed-normal} \mathcal N := J\mathcal A_\eta N_\zeta.
\end{equation}
Since  $q(x,\zeta(x),t)=\eta(x,t)$, differentiation along the equilibrium free surface gives  $\nabla_x\eta = \nabla_xq + (\partial_3q)\nabla_x\zeta \;\text{on }\Sigma$. Consequently,
\begin{equation}  \label{eq:N-equality} \mathcal N = J_{\eta}\mathcal A_\eta N_\zeta = \begin{pmatrix} -\nabla_x(\zeta+\eta)\\ 1 \end{pmatrix} = N_{\zeta+\eta} \qquad\text{on }\Sigma.
\end{equation}

\subsubsection{Pullback of the unknowns}

To distinguish the physical fields from their pullbacks, write
$(\widetilde u,\widetilde p)$ for the velocity and pressure $(u,p)$ in
\eqref{eq:moving-domain-full-system}. Define the fields on $\Omega$ by
\begin{equation}\label{eq:pullback-unknowns}
  v:=\widetilde u\circ\Phi_\eta,\qquad
  p:=\widetilde p\circ\Phi_\eta .
\end{equation}
For $f=f^{\mathrm{phys}}\circ\Phi_\eta$ the chain rule gives
$(\nabla f^{\mathrm{phys}})\circ\Phi_\eta=\nabla_{\mathcal A_\eta} f$.
Since $\partial_t\Phi_\eta=W\partial_t\bar\eta\,e_3$, we have
\begin{equation*}
 (\partial_t\widetilde u)\circ\Phi_\eta
 =\partial_t v-WK\partial_t\bar\eta\,\partial_3v .
\end{equation*}
The convection, divergence and stress terms transform according to
\begin{equation*}
\begin{aligned}
 (\widetilde u\cdot\nabla\widetilde u)\circ\Phi_\eta
   &=v\cdot\nabla_{\mathcal A_\eta}v,\\
 (\operatorname{div}\widetilde u)\circ\Phi_\eta
   &=\operatorname{div}_{\mathcal A_\eta}v,\\
 \bigl(\operatorname{div}S(\widetilde p,\widetilde u)\bigr)\circ\Phi_\eta
   &=\operatorname{div}_{\mathcal A_\eta}S_{\mathcal A_\eta}(p,v).
\end{aligned}
\end{equation*}

\subsection{Whole domain energy estimate}

Applying the transformation established in the previous subsection, we derive from \eqref{eq:moving-domain-full-system} the following system on the fixed equilibrium domain:
\begin{equation}  \label{eq:moving-domain-elliptic-system} \left\{ \begin{alignedat}{2} &\partial_tv + \operatorname{div}_{\mathcal{A}_{\eta}}S_{\mathcal{A}_{\eta}}(p,v) = F_{0}^1 &&\qquad \text{in } \Omega, \\ &\operatorname{div}_{\mathcal{A}_{\eta}}v = F_{0}^2 &&\qquad \text{in } \Omega, \\
&S_{\mathcal{A}_{\eta}}(p,v)\mathcal{N} = \left( g\eta - \sigma\operatorname{div}_{\Sigma}\mathcal L_{\Sigma,0}\eta-\sigma \operatorname{div}_{\Sigma}F_0^{5} \right)\mathcal{N} + F_{0}^3 &&\qquad \text{on }\Sigma, \\ &\partial_t\eta-v\cdot \mathcal{N} = F_{0}^4 &&\qquad \text{on }\Sigma, \\ &v\cdot n_{0} = 0 &&\qquad \text{on }\Sigma_{s}, \\
&U_{0}S_{\mathcal{A}_{\eta}}(p,v)n_{0} - \beta_{0} U_{0}v = F_{0}^6 &&\qquad \text{on }\Sigma_{s}, \\ &(\sigma\mathcal B_{\mathrm{cl},0}\eta+\sigma F_{0}^{5}\cdot n_{0}) =-\kappa \partial_{t}\eta-F_{0}^7 &&\qquad \text{on }\Gamma_0, \end{alignedat} \right.
\end{equation}
where
\begin{align}
&F_0^{1} :={}  WK\partial_t\bar\eta\,\partial_3v - v\cdot\nabla_{\mathcal{A}_\eta}v,\quad F_0^{2} :={} 0,\quad F_0^{3} :={} 0,\quad F_0^{4} :={} 0, \notag \\
&F_{0}^{5} :={} \big(\mathcal{R}(\nabla \eta)\big),\quad F_0^{6}:={} 0,\quad F_0^{7}:={} \hat{W}(\partial_{t}\eta). \label{eq:F05}
\end{align}
Moreover, $\mathcal{L}_{\Sigma}, \mathcal{B}_{\mathrm{cl},0},\mathcal{R}(\nabla \eta)$ are defined as follows
\begin{align}
\mathcal{L}_{\Sigma,0}:=&\frac{I}{\sqrt{1+|\nabla_{x}\zeta|^{2}}}-\frac{\nabla \zeta \otimes \nabla \zeta}{(\sqrt{1+|\nabla_{x}\zeta|^{2}})^{3}}, \label{def:linearized_operator} \\
\mathcal{B}_{\mathrm{cl},0}(\eta):=&\mathcal{L}_{\Sigma,0}\eta\cdot \nu . \label{def:contact_operator}
\end{align}
The constant $\beta_{0}$ in the Navier slip boundary condition is positive.

Let $F_{\zeta}(z)=z(1+|z|^2)^{-1/2}$ and define
\begin{equation*}
\mathcal{R}(\nabla \eta)=R_\zeta(\eta):=F_{\zeta}(\nabla\zeta+\nabla\eta) -F_{\zeta}(\nabla\zeta)-\mathcal{L}_{\Sigma_{0},0}(\nabla\eta).
\end{equation*}
Then
\begin{equation*}
\left[ \frac{\nabla_{x,y}(\zeta+\eta)} {\sqrt{1+|\nabla_{x,y}(\zeta+\eta)|^2}} - \frac{\nabla_{x,y}\zeta} {\sqrt{1+|\nabla_{x,y}\zeta|^2}} \right]  =(\mathcal{L}_{\Sigma,0}\nabla\eta) +R_\zeta(\eta).
\end{equation*}
Pointwise,  $|\func{div}R_\zeta(\eta)| \le C_\zeta\bigl(|\nabla\eta|\,|D^2\eta| +|\nabla\eta|^2\bigr)$ for a fixed bounded range of $\nabla\eta$. Under a fixed chart $Y$ with Jacobian $J_Y$,  $(\func{div}R_\zeta(\eta))\circ Y =J_Y^{-1}\func{div} \bigl(J_Y(DY)^{-1}(R_\zeta(\eta)\circ Y)\bigr)$. All temporal and angular derivatives are taken on this exact expression using the full product and chain rules.

\begin{remark}
\label{rem:R_x_dependence} The remainder depends on $x$ through $\nabla\zeta$: $\mathcal{R}(\nabla\eta)=R(x,\nabla\eta)$ with $R$ smooth in $(x,\xi)$ and $R(x,\xi)=O(|\xi|^{2})$. When a spatial or angular derivative falls on the $x $-dependence, the resulting term has the same structure as the term in which that derivative falls on $\eta$, with one derivative of $\eta$ fewer, and it is bounded by the same estimates. We do not write these terms separately.
\end{remark}

Moreover, for simplicity, we use $\mathcal{L}_{\Sigma,0}\eta$ to denote $\mathcal{L}_{\Sigma,0}\nabla \eta$.

We now apply $\partial_{t}^{j}$ for $0\leq j\leq m+1$ to the system \eqref{eq:moving-domain-elliptic-system}. Using the notation $(v_{j},p_{j},\eta_{j}):=(\partial_{t}^{j}v,\partial_{t}^{j}p,\partial_{t}^{j}\eta)$, we obtain
\begin{equation}  \label{eq:whole-domain-elliptic-system-j} \left\{ \begin{alignedat}{2} &\partial_tv_{j} + \operatorname{div}_{\mathcal{A}_{\eta}}S_{\mathcal{A}_{\eta}}(p_{j},v_{j}) = F_{0}^{1,j} &&\qquad \text{in } \Omega, \\
&\operatorname{div}_{\mathcal{A}_{\eta}}v_j = F_{0}^{2,j} &&\qquad \text{in } \Omega, \\ &S_{\mathcal{A}_{\eta}}(p_j,v_j)\mathcal{N} = \left( g\eta_j - \sigma\operatorname{div}_{\Sigma}\mathcal L_{\Sigma,0}\eta_j -\sigma \operatorname{div}_{\Sigma}F_0^{5,j}\right)\mathcal{N} + F_{0}^{3,j} &&\qquad \text{on }\Sigma, \\ &\partial_t\eta_{j}-v_{j}\cdot \mathcal{N} = F_{0}^{4,j} &&\qquad \text{on }\Sigma, \\ &v_{j}\cdot n_{0} = 0 &&\qquad \text{on }\Sigma_{s}, \\ &U_{0}S_{\mathcal{A}_{\eta}}(p_{j},v_{j})n_{0} - \beta_0U_{0}v_{j} = F_{0}^{6,j} &&\qquad \text{on }\Sigma_{s}, \\
&(\sigma\mathcal B_{\mathrm{cl},0}\eta_{j}+\sigma F_{0}^{5,j}\cdot n_{0}) =-\kappa \partial_{t}\eta_{j}- F_{0}^{7,j} &&\qquad \text{on }\Gamma_0,
\end{alignedat} \right.
\end{equation}
where
\begin{align}
&F_0^{1,j} :={}  \sum_{\substack{ a+b+c=j \\ a,b,c\geq0}} \frac{j!}{a!\,b!\,c!} \partial_t^a(WK)\, \partial_t^{b+1}\bar\eta\, \partial_3\partial_t^c v - \sum_{\substack{ a+b+c=j \\ a,b,c\geq0}} \frac{j!}{a!\,b!\,c!} \partial_t^a v\cdot \nabla_{\partial_t^b\mathcal{A}_\eta}\partial_t^c v \notag \\
&+ \mu\sum_{r=1}^{j}\binom{j}{r} \func{div}_{\mathcal{A}_\eta} \left( \mathbb{D}_{\partial_t^r\mathcal{A}_\eta}v_{j-r} \right) - \sum_{r=1}^{j}\binom{j}{r} \func{div}_{\partial_t^r\mathcal{A}_\eta} S_{\mathcal{A}_\eta}(p_{j-r},v_{j-r}) \notag \\
&+ \mu\sum_{r=1}^{j}\sum_{s=1}^{j-r} \binom{j}{r}\binom{j-r}{s} \func{div}_{\partial_t^r\mathcal{A}_\eta} \left( \mathbb{D}_{\partial_t^s\mathcal{A}_\eta}v_{j-r-s} \right), \label{eq:F01j} \\
&F_0^{2,j} :={}  -\sum_{r=1}^{j}\binom{j}{r} \func{div}_{\partial_t^r\mathcal{A}_\eta}v_{j-r},\quad F_0^{3,j} :={}  -\sigma\sum_{r=1}^{j}\binom{j}{r} \func{div}_{\Sigma}\left[ \partial_t^{j-r}\bigl(\mathcal{R}(\nabla\eta)\bigr) \right] \partial_t^r\mathcal{N} \label{eq:F02j} \\
&+ \mu\sum_{r=1}^{j}\binom{j}{r} \left( \mathbb{D}_{\partial_t^r\mathcal{A}_\eta}v_{j-r} \right)\mathcal{N} + \sum_{r=1}^{j}\binom{j}{r} \left( g\eta_{j-r} -\sigma\func{div}_{\Sigma}\mathcal{L}_{\Sigma,0}\eta_{j-r} \right)\partial_t^r\mathcal{N}  - \sum_{r=1}^{j}\binom{j}{r} S_{\mathcal{A}_\eta}(p_{j-r},v_{j-r}) \partial_t^r\mathcal{N} \notag \\
&+ \mu\sum_{r=1}^{j}\sum_{s=1}^{j-r} \binom{j}{r}\binom{j-r}{s} \left( \mathbb{D}_{\partial_t^s\mathcal{A}_\eta}v_{j-r-s} \right)\partial_t^r\mathcal{N},\quad F_0^{4,j} :={}  \sum_{r=1}^{j}\binom{j}{r} v_{j-r}\cdot\partial_t^r\mathcal{N}, \notag \\
&F_0^{5,j} :={} \partial_t^{j}\bigl[\mathcal{R}(\nabla\eta)\bigr], \qquad 0\leq j\leq m+1,\quad F_0^{6,j} :={}  \mu\sum_{r=1}^{j}\binom{j}{r} U_0 \left( \mathbb{D}_{\partial_t^r\mathcal{A}_\eta}v_{j-r} \right)n_0, \label{eq:F05j} \\
&F_0^{7,j} :={}  \partial_t^j\bigl[\hat W(\partial_t\eta)\bigr]. \label{eq:F07j}
\end{align}

For $1\leq j\leq m+1$, the Fa\`a di Bruno expansion of \eqref{eq:F05j} is
\begin{equation*}
F_0^{5,j}=\sum_{r=1}^{j}\frac1{r!} \sum_{\substack{i_1+\cdots+i_r=j\\i_1,\ldots,i_r\geq1}} \frac{j!}{i_1!\cdots i_r!} \mathcal R^{(r)}(\nabla\eta) [\nabla\eta_{i_1},\ldots,\nabla\eta_{i_r}].
\end{equation*}
At $j=0$, instead, $F_0^{4,0}=0$ and $F_0^{5,0}=\mathcal R(\nabla\eta)$. Here all sums with an empty index range are understood to be zero. The reference quantities $n_0$, $U_0$, $\beta_0$, $\kappa$, $\mathcal{N}_0$, $\mathcal{L}_{\Sigma,0}$, and $\mathcal{B}_{\mathrm{cl},0}$ are assumed to be independent of time.

Let $\chi_{\mathrm{int}+},\chi_{\mathrm{cl}+}$ be time-independent smooth cutoffs, and set $\chi_{\mathrm{int}+}+\chi_{\mathrm{cl}+}=1$.Moreover, set
\begin{align}
F_{\mathrm{cl}+}^{i,j}:=\chi_{\mathrm{cl}+}F_{0}^{i,j},\qquad F_{\mathrm{int}+}^{i,j}:=\chi_{\mathrm{int}+}F_{0}^{i,j} .  \label{eq:forcing_split}
\end{align}

We first prove the following lemma for projecting nonlinear terms near the contact line into 2-D.

\begin{lemma}
{\label{lem:derivative_replacement}}

We define $\mathcal{E}_{\mathrm{cl}-},\mathcal{E}_{\mathrm{bet}-},\mathcal{D}_{\mathrm{cl}-},\mathcal{D}_{\mathrm{bet}-}$ as follows
\begin{align}
\mathcal{E}_{\mathrm{cl}-}:=&\sum_{k_{1}+k_{2}=m-1}\mathscr{E}_{k_{1},k_{2},\mathrm{el},\mathrm{cl}}+\sum_{k_{1}+k_{2}=m}\mathscr{E}_{k_{1},k_{2},\mathrm{energy},\mathrm{cl}}+\sum_{k_{1}+k_{2}=m}(\|\chi_{\mathrm{cl},\Sigma}\partial_{t}^{k_{1}}\partial_{\theta}^{k_{2}}\eta\|^{2}_{L^{2}H^{\frac{3}{2}+\frac{\epsilon_{-}-\alpha}{2}}}),  \notag \\
&\quad+\sum_{k_{1}+k_{2}=m+1}\sum_{j_{1}=0}^{k_{1}}\sum_{j_{2}=0}^{k_{2}}(\|(\chi_{cl}+\chi_{\mathrm{bet}})\partial_{t}^{j_{1}}\partial_{\theta}^{j_{2}}\eta\|^{2}_{H^{1}}+\|(\chi_{cl}+\chi_{\mathrm{bet}})\partial_{t}^{j_{1}}\partial_{\theta}^{j_{2}}v\|_{L^{2}}^{2}), \\
\mathcal{E}_{\mathrm{bet}-}:=&\sum_{k_{1}+k_{2}=m}\mathscr{E}_{k_{1},k_{2},\mathrm{energy},\mathrm{cl}}+\sum_{k_{1}+k_{2}=m}(\|\chi_{\mathrm{bet},\Sigma}\partial_{t}^{k_{1}}\partial_{\theta}^{k_{2}}\eta\|_{L^{2}H^{\frac{3}{2}+\frac{\epsilon_{-}-\alpha}{2}}}^{2}) \\
&+\sum_{k_{1}+k_{2}=m-1}\sum_{j_{1}=0}^{k_{1}}\sum_{j_{2}=0}^{k_{2}}(\|\chi_{bet}\partial_{t}^{j_{1}}\partial_{\theta}^{j_{2}}v\|^{2}_{L^{2}W^{2,q_{\delta}}}+\|\chi_{\mathrm{bet},\Sigma}\partial_{t}^{j_{1}}\partial_{\theta}^{j_{2}}\eta\|^{2}_{L^{2}W^{3-\frac{1}{q_{\delta}},q_{\delta}}}+\|\chi_{\mathrm{bet}}\partial_{t}^{j_{1}}\partial_{\theta}^{j_{2}}p\|_{L^{2}W^{1,q_{\delta}}}^{2}), \\
\mathcal{D}_{\mathrm{cl}-}:=&\sum_{k_{1}+k_{2}=m}\mathscr{E}_{k_{1},k_{2},\mathrm{el},\mathrm{cl}}+\sum_{k_{1}+k_{2}=m+1}\mathscr{E}_{k_{1},k_{2},\mathrm{energy},\mathrm{cl}}, \\
\mathcal{D}_{\mathrm{bet}-}:=&\sum_{k_{1}+k_{2}=m}\mathscr{E}_{k_{1},k_{2},\mathrm{el},\mathrm{bet}}+\sum_{k_{1}+k_{2}=m+1}\mathscr{E}_{k_{1},k_{2},\mathrm{energy},\mathrm{bet}} \\
&+\sum_{k_{1}+k_{2}=m}\sum_{j_{1}=0}^{k_{1}}\sum_{j_{2}=0}^{k_{2}}(\|\chi_{bet}\partial_{t}^{j_{1}}\partial_{\theta}^{j_{2}}v\|^{2}_{L^{2}W^{2,q_{\delta}}}+\|\chi_{\mathrm{bet},\Sigma}\partial_{t}^{j_{1}}\partial_{\theta}^{j_{2}}\eta\|^{2}_{L^{2}W^{3-\frac{1}{q_{\delta}},q_{\delta}}}+\|\chi_{\mathrm{bet}}\partial_{t}^{j_{1}}\partial_{\theta}^{j_{2}}p\|_{L^{2}W^{1,q_{\delta}}}^{2}).
\end{align}
Suppose $\ell\in\{1,2\}$, and $\ell=1$ at the top level $k_{1}+k_{2}=m+1$. If we replace an $\ell$th-order $r,h$ derivative in the terms included in $\mathcal{E}_{\mathrm{cl}-},\mathcal{D}_{\mathrm{cl}-},\mathcal{E}_{\mathrm{bet}-},\mathcal{D}_{\mathrm{bet}-}$ by an $\ell$th-order $\theta$ derivative, the newly obtained terms are still included in $\mathcal{E}_{\mathrm{cl}-},\mathcal{D}_{\mathrm{cl}-},\mathcal{E}_{\mathrm{bet}-},\mathcal{D}_{\mathrm{bet}-}$, respectively.
\end{lemma}

\begin{proof}
In this proof, we use $(v,\eta,p)$ to denote $(\chi_{\mathrm{cl}}v,\chi_{\mathrm{cl},\Sigma}\eta,\chi_{\mathrm{cl}}p)$ or $(\chi_{\mathrm{bet}}v,\chi_{\mathrm{bet},\Sigma}\eta,\chi_{\mathrm{bet}}p)$, respectively.

For the energy part, by replacing an $r,h$ derivative by a $\theta$ derivative, using the fact that $\frac{\varepsilon_{-}-\alpha}{2}<\frac{1}{2}$, we obtain
\begin{align}
&\sum_{j=0}^{m}\sum_{k_{1}+k_{2}=j}(\|\partial_{t}^{k_{1}}\partial_{\theta}^{k_{2}}\partial_{\theta}v\|^{2}_{L^{2}L^{2}}+\|\partial_{t}^{k_{1}}\partial_{\theta}^{k_{2}}\partial_{\theta}\eta\|^{2}_{H^{\frac{1}{2}-\alpha}}+\|\partial_{t}^{k_{1}}\partial_{\theta}^{k_{2}}\partial_{\theta}\eta\|^{2}_{L^{2}H^{\frac{1}{2}+\frac{\epsilon_{-}-\alpha}{2}}})+\sum_{j=0}^{m+1}\sum_{k_{1}+k_{2}=j}(\|\partial_{t}^{k_{1}}\partial_{\theta}^{k_{2}}\partial_{\theta}\eta\|^{2}_{L^{2}L^{2}}) \notag \\
&+\sum_{j=0}^{m-1}\sum_{k_{1}+k_{2}=j}\Bigg(\|\partial_{t}^{k_{1}}\partial_{\theta}^{k_{2}}\partial_{\theta}v\|^{2}_{L^{2}W^{1,q_{\delta}}}+\|\partial_{t}^{k_{1}}\partial_{\theta}^{k_{2}}\partial_{\theta}\eta\|^{2}_{L^{2}W^{2-\frac{1}{q_{\delta}},q_{\delta}}}+\|\partial_{t}^{k_{1}}\partial_{\theta}^{k_{2}}\partial_{\theta}p\|^{2}_{L^{2}L^{q_{\delta}}} \notag \\
&+\|\partial_{t}^{k_{1}}\partial_{\theta}^{k_{2}}\partial_{\theta}^{2}v\|^{2}_{L^{2}L^{q_{\delta}}}+\|\partial_{t}^{k_{1}}\partial_{\theta}^{k_{2}}\partial_{\theta}^{2}\eta\|^{2}_{L^{2}W^{1-\frac{1}{q_{\delta}},q_{\delta}}}\Bigg) \notag \\
&\lesssim \sum_{j=0}^{m}\sum_{k_{1}+k_{2}=j}(\|\partial_{t}^{k_{1}}\partial_{\theta}^{k_{2}}v\|^{2}_{H^{1}}+\|\partial_{t}^{k_{1}}\partial_{\theta}^{k_{2}}\eta\|^{2}_{H^{\frac{3}{2}-\alpha}}+\|\partial_{t}^{k_{1}}\partial_{\theta}^{k_{2}+1}\eta\|^{2}_{H^{1}}) \notag \\
&+\sum_{j=0}^{m+1}\sum_{k_{1}+k_{2}=j}(\|\partial_{t}^{k_{1}}\partial_{\theta}^{k_{2}}\eta\|^{2}_{H^{1}})+\sum_{j=0}^{m}\sum_{k_{1}+k_{2}=j}(\|\partial_{t}^{k_{1}}\partial_{\theta}^{k_{2}}v\|^{2}_{H^{1}}+\|\partial_{t}^{k_{1}}\partial_{\theta}^{k_{2}}\eta\|^{2}_{H^{\frac{3}{2}-\alpha}}+\|\partial_{t}^{k_{1}}\partial_{\theta}^{k_{2}}p\|^{2}_{L^{2}}) \lesssim  \mathcal{E}_{\mathrm{cl}-}+\mathcal{E}_{\mathrm{bet-}}.
\end{align}

We then do the same for the dissipation; by definition, we have
\begin{align}
&\sum_{j=0}^{m+1}\sum_{k_{1}+k_{2}=j}(\|\partial_{t}^{k_{1}}\partial_{\theta}^{k_{2}}\partial_{\theta}v\|^{2}_{L^{2}L^{2}}+\|\partial_{t}^{k_{1}}\partial_{\theta}^{k_{2}}\partial_{\theta}\eta\|^{2}_{H^{\frac{1}{2}-\alpha}}) \notag \\
&+\sum_{j=0}^{m}\sum_{k_{1}+k_{2}=j}\Bigg(\|\partial_{t}^{k_{1}}\partial_{\theta}^{k_{2}}\partial_{\theta}v\|^{2}_{L^{2}W^{1,q_{\delta}}}+\|\partial_{t}^{k_{1}}\partial_{\theta}^{k_{2}}\partial_{\theta}\eta\|^{2}_{L^{2}W^{2-\frac{1}{q_{\delta}},q_{\delta}}}+\|\partial_{t}^{k_{1}}\partial_{\theta}^{k_{2}}\partial_{\theta}p\|^{2}_{L^{2}L^{q_{\delta}}} \notag \\
&+\|\partial_{t}^{k_{1}}\partial_{\theta}^{k_{2}}\partial_{\theta}^{2}v\|^{2}_{L^{2}L^{q_{\delta}}}+\|\partial_{t}^{k_{1}}\partial_{\theta}^{k_{2}}\partial_{\theta}^{2}\eta\|^{2}_{L^{2}W^{1-\frac{1}{q_{\delta}},q_{\delta}}}\Bigg) \notag \\
&\lesssim \sum_{j=0}^{m+1}\sum_{k_{1}+k_{2}=j}(\|\partial_{t}^{k_{1}}\partial_{\theta}^{k_{2}}v\|^{2}_{H^{1}}+\|\partial_{t}^{k_{1}}\partial_{\theta}^{k_{2}}\eta\|^{2}_{H^{\frac{3}{2}-\alpha}})+\sum_{j=0}^{m}\sum_{k_{1}+k_{2}=j}(\|\partial_{t}^{k_{1}}\partial_{\theta}^{k_{2}+1}v\|^{2}_{H^{1}}+\|\partial_{t}^{k_{1}}\partial_{\theta}^{k_{2}+1}\eta\|^{2}_{H^{\frac{3}{2}-\alpha}}) \notag \\
&+\sum_{k_{1}+k_{2}\leq m+1, k_{2}\geq 1}\|\partial_{t}^{k_{1}}\partial_{\theta}^{k_{2}}p\|_{L^{2}}^{2} \lesssim  \mathcal{D}_{\mathrm{cl}-}+\mathcal{D}_{\mathrm{bet}-}.
\end{align}
\end{proof}

\begin{proposition}
For the terms $\|\chi_{\mathrm{cl}}\partial_{t}^{k_{1}}\partial_{\theta}^{k_{2}}v\|^{2}_{L^{2}H^{1+\frac{\epsilon_{-}}{2}}}$ and $\|\chi_{\mathrm{bet}}\partial_{t}^{k_{1}}\partial_{\theta}^{k_{2}}v\|^{2}_{L^{2}H^{1+\frac{\epsilon_{-}}{2}}}$ included in $\mathcal{E}_{\mathrm{cl}},\mathcal{E}_{\mathrm{bet}}$ but excluded from $\mathcal{E}_{\mathrm{cl}}-,\mathcal{E}_{\mathrm{bet}-}$, we do not need to replace their derivatives while establishing the nonlinear estimates.
\end{proposition}

\begin{proof}
Consider the term
\begin{align}
\sum_{j=0}^{m}\|\partial_{t}^{j}v\|_{L^{2}H^{1+\frac{\epsilon_{-}}{2}}}.
\end{align}
If we replace one $r,h$ derivative by a $\theta$ derivative, then the following norm of the velocity,
\begin{equation}
\sum_{j=0}^{m} \bigl\| \partial_t^{j} \partial_\theta v \bigr\|_{L^2 H^{\epsilon_-/2}},
\end{equation}
cannot be controlled by the energy functional. In the nonlinear estimates of this paper, the norm $\|\cdot\|_{L^{2}H^{1+\frac{\epsilon_{-}}{2}}}$ enters in the following two ways, and in neither of them is the derivative replaced.

\begin{itemize}
\item Case 1: In 2-D, we use the Sobolev embedding   $H^{1+\frac{\epsilon_{-}}{2}}\hookrightarrow L^{\infty}$ to bound $\|u\|_{L^{\infty}}$. Therefore, in this case, we do not need to replace the derivative.

\item Case 2: To estimate the Bogovskii operator, we have to derive the following boundedness   $\|F^{2}\|_{W_{r,h}^{1,\frac{4}{3-2\epsilon_{+}}}}\lesssim \mathcal{E}$. But actually, what we really care about is   $\|F^{2}\|_{W_{r,h}^{1,\min\{\frac{4}{3-2\epsilon_{+}},2\}}}\lesssim \mathcal{E}$. Then it suffices to modify the estimate   $\|\func{div}_{\partial_{t}\mathcal{A}}\partial_{t}u\|_{L^{2}}\lesssim\|\partial_{t}\eta\|_{W^{1,\infty}}\|\partial_{t}u\|_{H^{1}}\lesssim\mathcal{E}_{2-D}$. \end{itemize}

Therefore, for the $H^{1+\frac{\epsilon_{-}}{2}}$ norm, we do not need to do derivative replacement.
\end{proof}

Before we state the estimate for forcing terms near the contact line, we first separate the forcing terms $F_{0}^{4,j}, F_{0}^{5,j}$. For $1\leq j\leq m+1$, we have  $F_{0}^{4,j}=v \cdot \partial_{t}^{j}\mathcal{N}+\mathcal{R}_{\mathrm{low}}^{4,j}$, and
\begin{equation}
{} F_{0}^{4,j}= v \cdot \partial_{t}^{j}\mathcal{N}+\mathcal{R}_{\mathrm{low}}^{4,j},\quad \mathcal{R}^{4,j}_{\mathrm{low}}:= \sum_{r=1}^{j-1}\binom{j}{r}\partial_{t}^{r}v\cdot\partial_{t}^{j-r}\mathcal{N}.
\end{equation}
\begin{equation}
{} F_{0}^{5,j}= \mathcal{R}^{\prime }(\nabla \eta)\nabla \eta_{j}+\mathcal{R}_{\mathrm{low}}^{5,j},\quad \mathcal{R}^{5,j}_{\mathrm{low}}:= \sum_{r=2}^{j}\frac{1}{r!}\sum_{i_{1}+\ldots+i_{r}=j, i_{1},\ldots,i_{r}>0}\frac{j!}{\prod_{k=1}^{r} i_{k}!}\mathcal{R}^{(r)}(\nabla \eta)[\nabla\eta_{i_{1}},\ldots,\nabla\eta_{i_{r}}].
\end{equation}
We put $\mathcal R_{\mathrm{low}}^{4,0} =\mathcal R_{\mathrm{low}}^{5,0}=0$ for bookkeeping. These conventions do not replace $F_0^{5,0}=\mathcal R(\nabla\eta)$; its contribution is treated separately below. For $a\in\{4,5\}$, we use $\mathcal R_{\mathrm{cl}+}^{a,j}=\chi_{\mathrm{cl}+}\mathcal R_{\mathrm{low}}^{a,j}$ and $\mathcal R_{\mathrm{int}+}^{a,j}=\chi_{\mathrm{int}+}\mathcal R_{\mathrm{low}}^{a,j}$. Before the main theorem, we state another lemma establishing bounds for the local forcing terms.

\begin{lemma}
{\label{lem:hard_nonlinearity}} For $1\le j\le m$,
\begin{align}
&\left|g\int_\Sigma\eta_j v\cdot\partial_t^j\mathcal N +\sigma\int_\Sigma(\mathcal L_{\Sigma,0}\nabla_\Sigma\eta_j) \cdot\nabla_\Sigma(v\cdot\partial_t^j\mathcal N)\right| \lesssim\mathcal E^{3/2},\quad \left|\sigma\int_\Sigma (\mathcal R'(\nabla\eta)\nabla_\Sigma\eta_j) \cdot\nabla_\Sigma\partial_t\eta_j\right| \lesssim\mathcal E^{3/2}, \notag \\
&\left|\sigma\int_\Sigma (\mathcal R'(\nabla\eta)\nabla_\Sigma\eta_j) \cdot\nabla_\Sigma(v\cdot\partial_t^j\mathcal N)\right| \lesssim\mathcal E^{3/2}.
\end{align}
Moreover, for $1\le j\le m+1$,
\begin{align}
&\left|g\int_{\Sigma}\eta_{j}v\cdot\partial_{t}^{j}\mathcal N +\sigma\int_{\Sigma}(\mathcal L_{\Sigma,0}\nabla_{\Sigma}\eta_j) \cdot\nabla_{\Sigma}(v\cdot\partial_t^j\mathcal N)\right| \lesssim \sqrt{\mathcal E}\,\mathcal D, \label{eq:4_e} \\
&\left|\sigma\int_{\Sigma}(\mathcal R'(\nabla\eta)\nabla_{\Sigma}\eta_j) \cdot\nabla_{\Sigma}\partial_t\eta_j -\frac{\sigma}{2}\frac{d}{dt}\int_{\Sigma} (\mathcal R'(\nabla\eta)\nabla_{\Sigma}\eta_j)\cdot\nabla_{\Sigma}\eta_j\right| \lesssim \sqrt{\mathcal E}\,\mathcal D, \label{eq:5_e} \\
&\left|\int_{\Sigma}(\mathcal R'(\nabla\eta)\nabla_{\Sigma}\eta_j) \cdot\nabla_{\Sigma}\eta_j\right| \lesssim \mathcal E^{3/2},\quad \left|\sigma\int_{\Sigma}(\mathcal R'(\nabla\eta)\nabla_{\Sigma}\eta_j) \cdot\nabla_{\Sigma}(v\cdot\partial_t^j\mathcal N)\right| \lesssim \sqrt{\mathcal E}\,\mathcal D. \label{eq:mix_e}
\end{align}
\end{lemma}

\begin{proof}
It suffices to prove the case $j=m+1$. When $j\leq m$, all of the terms can be bounded by the energy functional $\mathcal{E}$ with one fewer temporal derivative.

We first prove \eqref{eq:4_e}. We have
\begin{align}
{\label{eq:separation}} &g\int_{\Sigma}\eta_{j}v\cdot \partial_{t}^{j}\mathcal{N}+\sigma\int_{\Sigma}\mathcal{L}_{\Sigma,0}\nabla_{\Sigma}\eta_{j}\cdot \nabla_{\Sigma}(v\cdot \partial_{t}^{j}\mathcal{N})  \notag \\
=&g\int_{\Sigma}\eta_{j}v\cdot \partial_{t}^{j}\mathcal{N}-\sigma\int_{\Sigma}\mathcal{L}_{\Sigma,0}\nabla_{\Sigma}\eta_{j}\cdot (\nabla_{\Sigma}v)^{T} \nabla_{\Sigma}\eta_{j}-\sigma\int_{\Sigma}\mathcal{L}_{\Sigma,0}\nabla_{\Sigma}\eta_{j}\cdot (\nabla_{\Sigma}^{2}\eta_j{} v_{||}),
\end{align}
where
\begin{align}
\nabla_{\Sigma}v:=
\begin{pmatrix}
\partial_{1}v_{1},\partial_{2}v_{1} \\
\partial_{1}v_{2},\partial_{2}v_{2}\end{pmatrix},
\end{align}
and $v_{\Sigma}=(v_{1},v_{2})$.  For the first two terms on the right-hand side of the inequality, we have
\begin{align*}
&g\int_{\Sigma}\eta_{j}v\cdot \partial_{t}^{j}\mathcal{N}-\sigma\int_{\Sigma}\mathcal{L}_{\Sigma,0}\nabla_{\Sigma}\eta_{j}\cdot (\nabla_{\Sigma}v)^{T} \nabla_{\Sigma}\eta_{j} \\
&\lesssim \|\eta_{j}\|_{L^{+\infty}}\|v\|_{L^{2}(\Sigma)}\|\nabla \eta_{j}\|_{L^{\frac{4}{1+2\alpha}}}+\sigma\int_{\Sigma}\mathcal{L}_{\Sigma,0}\nabla_{\Sigma}\eta_{j}\cdot (\nabla_{\Sigma}(v_{\mathrm{cl}+}+v_{\mathrm{int}+}))^{T} \nabla_{\Sigma}\eta_{j} \\
&\lesssim \|v\|_{H^{1}}\|\eta_{j}\|_{H^{\frac{3}{2}-\alpha}}^{2}+\sigma \|\nabla_{\Sigma}\eta_{j}\|^{2}_{L^{\frac{4}{1+2\alpha}}}\|\nabla_{\Sigma}v_{\mathrm{int}+}\|_{L^{\frac{2p}{3-p}}(\Sigma)}+\sigma\|\nabla_{\Sigma}\eta_{j}\|^{2}_{L^{2}L^{\frac{1}{\alpha}}}\|\nabla_{\Sigma}v_{\mathrm{cl}+}\|_{L^{\infty}L^{\frac{1}{1-2\alpha}}} \\
&\lesssim \|v\|_{H^{1}}\|\eta_{j}\|_{H^{\frac{3}{2}-\alpha}}^{2}+\sigma \|\eta_{j}\|^{2}_{H^{\frac{3}{2}-\alpha}}\|v_{\mathrm{int}+}\|_{W^{2,p}}+\sigma\|\eta_{j}\|^{2}_{H^{\frac{3}{2}-\alpha}}(\|v_{\mathrm{cl}+}\|_{L^{\infty}W^{2,q_{\delta}}}+\|\partial_{\theta}v_{\mathrm{cl}+}\|_{L^{\infty}H^{1}}) \lesssim \sqrt{\mathcal{E}}\mathcal{D}.
\end{align*}
For the last term on the right-hand side of \eqref{eq:separation}, we write the surface calculation in the horizontal coordinates $x\in D$. Every gradient and divergence in the following display is horizontal. Integration by parts, the symmetry of $\mathcal L_{\Sigma,0}$, and $v_\parallel\cdot\nu_D=v\cdot n_0=0$ on $\partial D$ give
\begin{align}
&\int_D \left(\mathcal{L}_{\Sigma,0}\nabla\partial_t^j\eta\right) \cdot\left(\nabla^2\partial_t^j\eta\,v_\parallel\right)\,dx = \frac12\int_{\partial D} (v_\parallel\cdot\nu_D) \left( \nabla\partial_t^j\eta\cdot \mathcal{L}_{\Sigma,0}\nabla\partial_t^j\eta \right)\,d\ell \notag \\
&- \frac12\int_D \nabla\partial_t^j\eta\cdot \left[ (\func{div}v_\parallel)\mathcal{L}_{\Sigma,0} + (v_\parallel\cdot\nabla)\mathcal{L}_{\Sigma,0} \right] \nabla\partial_t^j\eta\,dx \notag \\
&=-\frac12\int_D \nabla\partial_t^j\eta\cdot \left[ (\func{div}v_\parallel)\mathcal{L}_{\Sigma,0} + (v_\parallel\cdot\nabla)\mathcal{L}_{\Sigma,0} \right] \nabla\partial_t^j\eta\,dx.
\end{align}
As in the estimate for the first two terms of \eqref{eq:separation}, we have
\begin{align}
&\frac12\int_D \nabla\partial_t^j\eta\cdot \left[ (\func{div}v_\parallel)\mathcal{L}_{\Sigma,0} + (v_\parallel\cdot\nabla)\mathcal{L}_{\Sigma,0} \right] \nabla\partial_t^j\eta\,dx \\
&\lesssim \|\eta_{j}\|_{H^{\frac{3}{2}-\alpha}}^{2}(\|v_{\mathrm{int}+}\|_{W^{2,p}}+\|v_{\mathrm{cl}+}\|_{L^{\infty}W^{2,q_{\delta}}}+\|\partial_{\theta}v_{\mathrm{cl}+}\|_{L^{\infty}H^{1}})\lesssim\sqrt{\mathcal{E}}\mathcal{D}.
\end{align}

For \eqref{eq:5_e}, put $A=\mathcal R'(\nabla\eta)$ and $w=\nabla_{\Sigma}\eta_j$. The matrix $A$ is symmetric, so
\begin{equation}
\sigma\int_{\Sigma}(Aw)\cdot\partial_t w -\frac{\sigma}{2}\frac{d}{dt}\int_{\Sigma}(Aw)\cdot w =-\frac{\sigma}{2}\int_{\Sigma}((\partial_t A)w)\cdot w .
\end{equation}
By Remark~\ref{rem:R_x_dependence}, $|A|\lesssim|\nabla\eta|$ and $|\partial_t A|\lesssim |\nabla\partial_t\eta|$. Both gradients are bounded in $L^\infty(\Sigma)$ by $C(\varepsilon)\sqrt{\mathcal E}$. Indeed, for $\eta,\partial_t\eta$, the interior bound follows from $W^{3-1/p,p}\subset W^{1,\infty}$ on the two-dimensional surface. In the collar the radial derivatives of $u$ and $\partial_\theta u$ are bounded in $L^2_\theta L^\infty_r$ by the radial elliptic norms through mixed order two. For the angular part, $\partial_\theta \eta$ and $\partial_\theta^2 \eta$, $\partial_\theta \p_{t}\eta$ and $\partial_\theta^2 \p_{t}\eta$  are bounded in $L^2_\theta H^1_r$ by the surface energy. The one-dimensional angular Sobolev bound then controls both components in $L^\infty_{\theta,r}$. Here $m\ge3$, and the fixed chart and the cutoffs only change the constants. Finally, $\|w\|_{L^2}^2\le\mathcal E$ for the capillary correction itself, and $\|w\|_{L^2}^2\le\mathcal D$ for its time derivative.

For \eqref{eq:mix_e}, use $v\cdot\partial_t^j\mathcal N =-v_\parallel\cdot\nabla\eta_j$. Integrate the term with $\nabla^2\eta_j$ in the direction $v_\parallel$, as in \eqref{eq:4_e}. The boundary flux vanishes since $v_\parallel\cdot\nu_D=0$. By symmetry of $A$, the absolute value on the left of \eqref{eq:mix_e} is bounded by
\begin{equation}
C\int_{\Sigma} \bigl(|A|\,|\nabla v_\parallel| +|v_\parallel|\,|\nabla A|\bigr)|\nabla\eta_j|^2\,dx .
\end{equation}
Specifically, after integration by part, we separate this nonlinear term into local parts by estimating it via slicing norms and interior part by estimating it by interior norms.  We have,
\begin{align}
&\int_{\Sigma}\mathcal{R}^{\prime }(\nabla \eta)\nabla_{\Sigma}\eta_{j}\cdot \nabla_{\Sigma}v\nabla_{\Sigma}\eta_{j} -\frac12\int_{\Sigma} \nabla\partial_t^j\eta\cdot \left[ (\func{div}v_\parallel)\mathcal{R}^{\prime }(\nabla \eta) + (v_\parallel\cdot\nabla)\mathcal{R}^{\prime }(\nabla \eta) \right] \nabla\partial_t^j\eta\,dx \notag \\
&\lesssim \|\nabla(\chi_{\mathrm{bet}+}\eta_{j})\|^{2}_{L^{2}L^{\frac{1}{\alpha}}}(\|v_{\mathrm{bet}+}\|_{L^{\infty}W^{2,q_{\delta}}}+\|\p_{\theta}v_{\mathrm{bet}+}\|_{L^{\infty}H^{1}}) \notag \\
&(\|\eta_{\mathrm{bet}+}\|_{L^{\infty}W^{3-\frac{1}{q_{\delta}},q_{\delta}}} +\|\partial_{\theta}\eta_{\mathrm{bet}+}\|_{L^{\infty}W^{3-\frac{1}{q_{\delta}},q_{\delta}}}+\|\p_{\theta}^{2}\eta_{\mathrm{bet}+}\|_{L^{\infty}H^{1}}) \notag \\
&+\|\nabla(\chi_{\mathrm{int}+}\eta_{j})\|^{2}_{L^{\frac{4}{1+2\alpha}}}(\|v_{\mathrm{int}+}\|_{W^{1,\frac{2p}{3-p}}})(\|\eta_{\mathrm{int}+}\|_{W^{2,\frac{2p}{3-p}}} ) \notag \\
&\lesssim\sigma\|\eta_{j}\|^{2}_{H^{\frac{3}{2}-\alpha}} \Bigl(\|v_{\mathrm{int}+}\|_{W^{2,p}} +\|v_{\mathrm{bet}+}\|_{L^{\infty}W^{2,q_{\delta}}} +\|\partial_{\theta}v_{\mathrm{bet}+}\|_{L^{\infty}H^{1}}\Bigr) \notag \\
&\times\Bigl(\|\eta_{\mathrm{int}+}\|_{W^{3-\frac{1}{p},p}} +\|\eta_{\mathrm{bet}+}\|_{L^{\infty}W^{3-\frac{1}{q_{\delta}},q_{\delta}}} +\|\partial_{\theta}\eta_{\mathrm{bet}+}\|_{L^{\infty}W^{3-\frac{1}{q_{\delta}},q_{\delta}}}+\|\p_{\theta}^{2}\eta_{\mathrm{bet}+}\|_{L^{\infty}H^{1}}\Bigr) \lesssim\sqrt{\mathcal{E}}\mathcal{D},
\end{align}
where we used the Sobolev embedding $W^{2-\frac{1}{p},p}\hookrightarrow L^{\infty}$ $W^{3-\frac{1}{p},p}\hookrightarrow W^{2,\frac{2p}{3-p}}$ in two-dimensional space, and $W^{2-\frac{1}{q_{\delta}},q_{\delta}}\hookrightarrow L^{\infty} $, $W^{3-\frac{1}{q_{\delta}},q_{\delta}}\hookrightarrow W^{2,\frac{1}{1-\epsilon_{-}}}$ in one-dimensional space in slicing norms estimate, and the H\"older index
\begin{align}
\frac{3-p}{2p}+\frac{1+2\alpha}{4}\times 2+\frac{1}{\infty}<1,\qquad 1-\varepsilon_{-}+2\alpha+\frac{1}{\infty}<1.
\end{align}
The derivatives in this display are horizontal derivatives of the surface fields. In particular, $|\nabla A|\lesssim |\nabla\eta|+|\nabla^2\eta|$: the first term allows for the fixed $x$-dependence in Remark~\ref{rem:R_x_dependence}.

\end{proof}

We now state the lemma estimating the forcing terms near the contact line.

\begin{lemma}
{\label{thm:contact_line_forcing_estimate}} Suppose $m\geq 3$. Then
\begin{align}
&\sum_{j=0}^{m}(\|F_{\mathrm{cl}+}^{1,j}\|_{L^{2}(H^{1})^{*}(\Omega)} +\|F_{\mathrm{cl}+}^{2,j}\|_{L^{2}L^{2}(\Omega)}+\|F_{\mathrm{cl}+}^{3,j}\|_{L^{2}(H^{\frac{1}{2}})^{*}(\Sigma)}+\|\nabla_{\Sigma}\mathcal{R}_{\mathrm{cl}+}^{4,j}\|_{L^{2}L^{\frac{1}{1-\alpha}}} \notag \\
&+\|\partial_{t}\mathcal{R}_{\mathrm{cl}+}^{5,j}\|_{L^{2}L^{\frac{1}{1-\alpha}}}+\|\mathcal{R}_{\mathrm{cl}+}^{5,j+1}\|_{L^{2}L^{2}}+\|\nabla_{\Sigma}\mathcal{R}_{\mathrm{cl}+}^{5,j}\|_{L^{2}L^{\frac{1}{1-\alpha}}} +\|F_{\mathrm{cl}+}^{6,j}\|_{L^{2}(H^{\frac{1}{2}})^{*}(\Sigma_{s})}+[F_{\mathrm{cl}+}^{7,j}]_{\Gamma_{0}})\lesssim\mathcal{E}, \label{eq:contact_line_forcing_estimate_e} \\
&\|F_{\mathrm{cl}+}^{1,m+1}\|_{(H^{1}(\Omega))^{*}}+\sum_{j=0}^{m}\|F_{\mathrm{cl}+}^{1,j}\|_{L^{2}(H^{1})^{*}(\Omega)} +\sum_{j=0}^{m+1}(\|F_{\mathrm{cl}+}^{2,j}\|_{L^{2}L^{2}(\Omega)}+\|F_{\mathrm{cl}+}^{3,j}\|_{L^{2}(H^{\frac{1}{2}})^{*}(\Sigma)}+\|\nabla_{\Sigma}\mathcal{R}_{\mathrm{cl}+}^{4,j}\|_{L^{2}L^{\frac{1}{1-\alpha}}} \notag \\
&+\|\partial_{t}\mathcal{R}_{\mathrm{cl}+}^{5,j}\|_{L^{2}L^{\frac{1}{1-\alpha}}}+\|\mathcal{R}_{\mathrm{cl}+}^{5,j}\|_{L^{2}L^{2}}+\|\nabla_{\Sigma}\mathcal{R}_{\mathrm{cl}+}^{5,j}\|_{L^{2}L^{\frac{1}{1-\alpha}}} +\|F_{\mathrm{cl}+}^{6,j}\|_{L^{2}(H^{\frac{1}{2}})^{*}(\Sigma_{s})}+[F_{\mathrm{cl}+}^{7,j}]_{\Gamma_{0}})\lesssim\sqrt{\mathcal{E}\mathcal{D}}. \label{eq:contact_line_forcing_estimate_d}
\end{align}
\begin{equation}
\sum_{j=0}^{m}\int_{\Sigma}\mathcal{R}_{\mathrm{cl}+}^{5,j}\cdot \nabla_{\Sigma}\mathcal{R}_{\mathrm{cl}+}^{4,j}\lesssim \mathcal{E}^{2},\quad \sum_{j=0}^{m+1}\int_{\Sigma}\mathcal{R}_{\mathrm{cl}+}^{5,j}\cdot \nabla_{\Sigma}\mathcal{R}_{\mathrm{cl}+}^{4,j}\lesssim \mathcal{E}\mathcal{D}.
\end{equation}
\end{lemma}

\begin{proof}
\emph{Step 1: reduction to slices.} Near $\Gamma$ we write every forcing term in the chart $(\theta,r,h)$. By \eqref{eq:F01j}--\eqref{eq:F07j}, each $F^{i,j}_{\mathrm{cl}+}$ is a finite sum of terms $\chi_{\mathrm{cl}+}\,c\,\prod_{k=1}^{\ell}\partial_{t}^{s_{k}}D_{k}F_{k}$, for $\ell\geq2$, $\sum_{k}s_{k}\leq j$. Here $c$ is a smooth function of $y$, $\nabla\bar\eta$ and $\nabla\zeta$, bounded with its derivatives for $\mathcal{E}\leq1$. Each $F_{k}$ is one of $v$, $p$, $\bar\eta$, $\eta$, and $D_{k}$ is a Cartesian derivative of order at most two. A factor $\partial_{t}^{b+1}\bar\eta$, which occurs in $F^{1,j}_{0}$, is counted with $s_{k}=b$: by the kinematic equation its norms are those of $\partial_{t}^{b}v$ on $\Sigma$, extended. A Cartesian derivative is $e^{\theta}\partial_{\theta}+e^{r}\partial_{r}+e^{h}\partial_{h}$ with smooth coefficients. The terms that contain only $\partial_{r}$ and $\partial_{h}$ are, slice by slice, the forcing terms of the two-dimensional problem, and the slice norms in which \cite{GT} estimates them are the slice norms in $\mathcal{E}_{\mathrm{cl}-}$, $\mathcal{D}_{\mathrm{cl}-}$, $\mathcal{E}_{\mathrm{bet}-}$, $\mathcal{D}_{\mathrm{bet}-}$. The terms in which some $\partial_{r}$ or $\partial_{h}$ is replaced by $\partial_{\theta}$ are estimated in the same way, by Lemma~\ref{lem:derivative_replacement}: the replaced norms are still part of these functionals, with $\ell=1$ replacement at the top level $j=m+1$.

\emph{Step 2: integrability in $\theta$.} Fix one product and let $k_{0}$ be an index with the largest $s_{k}$. If $m\geq4 $, then $s_{k}\leq m-2$ for every $k\neq k_{0}$, since two indices with $s_{k}\geq m-1$ would give $2(m-1)\leq m+1$. For such $k$ the slice norm of $\partial_{t}^{s_{k}}D_{k}F_{k}$ is controlled by $\mathcal{E}$ together with one more $\theta$-derivative, at a level $\leq m-1$. By $H^{1}(\mathbb{T})\subset L^{\infty}(\mathbb{T})$ it is bounded in $L^{\infty}_{\theta}$ by $\sqrt{\mathcal{E}}$. The factor $k_{0}$ is estimated in $L^{2}_{\theta}$ of its slice norm. H\"older's inequality in $\theta$, with exponents $(2,\infty,\dots,\infty)$, and the two-dimensional estimates of \cite{GT} in each slice give the bounds \eqref{eq:contact_line_forcing_estimate_e} and \eqref{eq:contact_line_forcing_estimate_d}. One gets $\mathcal{E}$ or $\sqrt{\mathcal{E}\mathcal{D}}$ according to whether the slice norm of the factor $k_{0}$ belongs to $\mathcal{E}$ or to $\mathcal{D}$. The integrals of products of forcing terms in the last two lines of the lemma are bounded in the same way. If $m=3$, in this case, the temporal derivative of $\mathcal{R}_{\mathrm{cl}+}^{5,j}$ and $\mathcal{R}_{\mathrm{cl}+}^{5,j+1}$ will be controlled by the end of the proof. Then it suffices to estimate the forcing terms in inequality \eqref{eq:contact_line_forcing_estimate_d}. For those forcing terms, two indices with $s_{k}=2$ can occur only when $j=4=m+1$ and every other $s_{k}$ vanishes. Now suppose $s_{1}=s_{2}=2$ and other indices vanish. Then the factor $k=k_{1}$ is estimated in $L_{\theta}^{2}$ of its slice norm which is bounded by $\sqrt{\mathcal{E}}$; while the factor $k_{2}$ is estimated in $L_{\theta}^{\infty}$ of its sliced norm which is bounded by $\sqrt{\mathcal{D}}$.  This proves \eqref{eq:contact_line_forcing_estimate_d} except for the terms $\mathcal{R}_{\mathrm{cl}+}^{5,j}$ and $\mathcal{R}_{\mathrm{cl}+}^{5,j+1}$.

At $j=m+1$ the $F^{1}$ part of \eqref{eq:contact_line_forcing_estimate_d} is measured in $(H^{1}(\Omega))^{*}$, not in $L^{2}(H^{1})^{*}$, which by the convention of Subsection~\ref{subsec:notation} is $L^{2}$ in $\theta$ of a slice dual. Indeed, $F^{1,m+1}_{\mathrm{cl}+}$ contains $\mu\func{div}_{\mathcal{A}_{\eta}}(\mathbb{D}_{\partial_{t}^{m+1}\mathcal{A}_{\eta}}v)$, whose $\theta\theta$ part contains $\partial_{\theta}^{2}\bar{\eta}_{m+1}\,\partial_{3}v$ with a smooth coefficient, while $\mathcal{D}$ controls $\eta_{m+1}$ only in $H^{\frac{3}{2}-\alpha}(\Sigma)$. For $\eta_{m+1}=f(\theta)\phi(r)$ with $\hat{f}_{k}=|k|^{-5/2}$ one has $\|\eta_{m+1}\|^{2}_{H^{3/2-\alpha}}\sim\sum_{k}|k|^{-2-2\alpha}<\infty$, while the square of the $L^{2}_{\theta}$ slice-dual norm of this part is $\sim\sum_{k}|k|^{-1}=\infty$: the extension acts slice by slice, and a slice dual cannot integrate by parts in $\theta$. Against $w\in H^{1}(\Omega)$ the bound holds. Integrating the second-order term by parts in all variables, as in Step~A2 and Step~B1(c) of the proof of Theorem~\ref{thm:top_level} with $\partial_{\theta}^{m+1}\bar{\eta}$ replaced by $\partial_{t}^{m+1}\bar{\eta}$, the bulk term needs only $\eta_{m+1}\in H^{\frac{3}{2}-\alpha}(\Sigma)$ and the face terms on $\Sigma$ and $\Sigma_{s}$ are bounded as in Lemma~\ref{lem:tl_first}; this gives $\lesssim\sqrt{\mathcal{E}}\sqrt{\mathcal{D}}\,\|w\|_{H^{1}}$. Every use of $F^{1}$ pairs it with such a $w$. At $m=3$, $\partial_{t}\mathcal{R}_{\mathrm{cl}+}^{5,j}$ in \eqref{eq:contact_line_forcing_estimate_e} contains $\mathcal{R}^{\prime \prime }[\nabla\eta_{2},\nabla\eta_{2}]$, with two factors at time level $m-1 $; they are bounded with the energy norms $\partial_{\theta}\eta_{2}\in H^{\frac{3}{2}-\alpha}(\Sigma)\subset L^{\infty}$ and $\partial_{r}\eta_{2}\in H^{1}_{\theta}L^{1/\alpha}_{r}\subset L^{\infty}_{\theta}L^{4}_{r}$. At $m=4$, $\partial_{t}\mathcal{R}_{\mathrm{cl}+}^{5,j}$ in \eqref{eq:contact_line_forcing_estimate_d} contains $\mathcal{R}^{\prime \prime }[\nabla\eta_{3},\nabla\eta_{3}]$ with two factors at time level $m-1$; they are bounded with the energy norms $\partial_{\theta}\eta_{3}\in H^{\frac{3}{2}-\alpha}(\Sigma)\subset L^{\infty}$ and $\partial_{r}\eta_{3}\in H^{1}_{\theta}L^{1/\alpha}_{r}\subset L^{\infty}_{\theta}L^{4}_{r}$. At $m=3$, $\partial_{t}\mathcal{R}_{\mathrm{cl}+}^{5,j}$ in \eqref{eq:contact_line_forcing_estimate_d} contains $\mathcal{R}^{\prime \prime }[\nabla\eta_{2},\nabla\eta_{3}]$ with two factors at time level $m,m-1$; $\nabla \eta_{2}$ is bounded with the energy norms $\partial_{\theta}\eta_{2}\in H^{\frac{3}{2}-\alpha}(\Sigma)\subset L^{\infty} $ and $\partial_{r}\eta_{2}\in H^{1}_{\theta}L^{1/\alpha}_{r}\subset L^{\infty}_{\theta}L^{4}_{r}$, and $\nabla \eta_{3}$ is bounded with the dissipation norms $\partial_{\theta}\eta_{3}\in H^{\frac{3}{2}-\alpha}(\Sigma)\subset L^{\infty}$ and $\partial_{r}\eta_{3}\in H^{1}_{\theta}L^{1/\alpha}_{r}\subset L^{\infty}_{\theta}L^{4}_{r}$.
\end{proof}

The essential boundary condition is  $v_{j}\cdot n_0=0 \;\text{on }\Sigma_{\mathrm{art}}$. For every smooth vector test function $u$ satisfying  $u\cdot n_0=0 \;\text{on }\Sigma_{\mathrm{s}}$, and every smooth scalar test function $q$ on $\Omega$, by multiplying both sides of the first equation of \eqref{eq:whole-domain-elliptic-system-j} by $Ju$ and integrating by parts as in \cite{GT}, we obtain the weak formulation as follows. Assume that  $S_{\mathcal{A}_\eta}(p,v) =pI-\mu\mathbb{D}_{\mathcal{A}_\eta}v, \; U_0=I-n_0\otimes n_0$, and that the geometric coefficients satisfy the Piola identity  $\partial_k\bigl(J(\mathcal{A}_\eta)_{ik}\bigr)=0$. We use surface measures compatible with $J\mathcal{A}_\eta\nu=\mathcal{N } / |\mathcal{N}_{0}|$ on $\Sigma$ and $J\mathcal{A} _\eta n_0=J_s n_0$ on $\Sigma_s$. Here $J_s$ is the wall surface Jacobian; if the flattening map is the identity on the wall, then $J_s=1$. For every smooth vector test function $u$ satisfying $u\cdot n_0=0$ on $\Sigma_s$, the weak momentum equation obtained by testing equation \eqref{eq:whole-domain-elliptic-system-j} by $Ju$ is
\begin{align}
&\int_\Omega J\,\partial_t v_j\cdot u + \frac{\mu}{2}\int_\Omega J\,\mathbb{D}_{\mathcal{A}_\eta}v_j : \mathbb{D}_{\mathcal{A}_\eta}u - \int_\Omega Jp_j\func{div}_{\mathcal{A}_\eta}u  \notag \\
&\quad+ \beta_0\int_{\Sigma_s} J_s(U_0v_j)\cdot(U_0u) + g\int_\Sigma\eta_j(u\cdot\mathcal{N}) - \sigma\left\langle \func{div}_{\Sigma} \mathcal{L}_{\Sigma,0}\eta_j,u\cdot\mathcal{N }\right\rangle_\Sigma \notag \\
&\qquad= \int_\Omega JF_0^{1,j}\cdot u - \left\langle F_0^{3,j},u\right\rangle_\Sigma - \left\langle J_sF_0^{6,j},U_0u\right\rangle_{\Sigma_s} + \sigma\left\langle \func{div}_{\Sigma}F_0^{5,j},u\cdot\mathcal{N }\right\rangle_\Sigma. \label{eq:whole-domain-weak-momentum-j}
\end{align}

When $j\leq m$, taking $u=v_j$ and using the divergence and kinematic equations, we obtain
\begin{align}
&\int_\Omega J\,\partial_t v_j\cdot v_j + \frac{\mu}{2}\int_\Omega J|\mathbb{D}_{\mathcal{A}_\eta}v_j|^2 + \beta_0\int_{\Sigma_s}J_s|U_0v_j|^2 + g\int_\Sigma\eta_j \left(\partial_t\eta_j-F_0^{4,j}\right) - \sigma\left\langle \func{div}_{\Sigma}\mathcal{L}_{\Sigma,0}\eta_j, \partial_t\eta_j-F_0^{4,j} \right\rangle_\Sigma \notag \\
&= \int_\Omega Jp_jF_0^{2,j} + \int_\Omega JF_0^{1,j}\cdot v_j - \left\langle F_0^{3,j},v_j\right\rangle_\Sigma - \left\langle J_sF_0^{6,j},U_0v_j\right\rangle_{\Sigma_s} + \sigma\left\langle \func{div}_{\Sigma}F_0^{5,j}, \partial_t\eta_j-F_0^{4,j} \right\rangle_\Sigma.
\end{align}

For the kinetic and gravitational terms,
\begin{equation}
\int_\Omega J\,\partial_t v_j\cdot v_j  = \frac{d}{dt}\frac12\int_\Omega J|v_j|^2 - \frac12\int_\Omega(\partial_tJ)|v_j|^2,\quad g\int_\Sigma\eta_j\partial_t\eta_j  = \frac{d}{dt}\frac g2\|\eta_j\|_{L^2(\Sigma)}^2.
\end{equation}

Integration by parts on $\Sigma$ gives
\begin{align}
&-\sigma\left\langle \func{div}_{\Sigma }\mathcal{L}_{\Sigma,0}\eta_j+\func{div}_{\Sigma}F_{0}^{5,j}, \partial_t\eta_j-F_0^{4,j} \right\rangle_\Sigma = \sigma\int_\Sigma a_\zeta\nabla_\Sigma\eta_j \cdot\nabla_\Sigma\partial_t\eta_j - \sigma\int_\Sigma a_\zeta\nabla_\Sigma\eta_j \cdot\nabla_\Sigma F_0^{4,j} \notag \\
&-\sigma \left\langle \mathcal{B}_{\mathrm{cl},0}\eta_j+F_{0}^{5,j}\cdot n_{0}, \partial_t\eta_j-F_0^{4,j} \right\rangle_{\Gamma_0}+\sigma\langle F_{0}^{5,j},\nabla_{\Sigma}(\partial_{t}\eta_{j}-F_{0}^{4,j}) \rangle_{\Sigma},
\end{align}
where we use $a_{\zeta}$ to denote $\mathcal{L}_{\Sigma,0}$. Since $a_\zeta$ is symmetric and the reference surface gradient and measure are independent of time, we have
\begin{align}
\sigma\int_\Sigma a_\zeta\nabla_\Sigma\eta_j \cdot\nabla_\Sigma\partial_t\eta_j ={}& \frac{d}{dt} \left( \frac{\sigma}{2}\int_\Sigma a_\zeta\nabla_\Sigma\eta_j\cdot\nabla_\Sigma\eta_j \right).
\end{align}
Consequently, using $\sigma\mathcal{B}_{\mathrm{cl},0}\eta_j+\sigma F_{0}^{5,j}\cdot n_{0} =-\kappa\partial_t\eta_j-F_0^{7,j}$,
\begin{align}
&-\sigma\left\langle \func{div}_{\Sigma}\mathcal{L}_{\Sigma,0}\eta_j+\func{div}_{\Sigma}F_{0}^{5,j}, \partial_t\eta_j-F_0^{4,j} \right\rangle_\Sigma = \frac{d}{dt} \left( \frac{\sigma}{2}\int_\Sigma a_\zeta\nabla_\Sigma\eta_j\cdot\nabla_\Sigma\eta_j \right) - \frac{\sigma}{2}\int_\Sigma (\partial_ta_\zeta) \nabla_\Sigma\eta_j\cdot\nabla_\Sigma\eta_j \notag \\
&- \sigma\int_\Sigma a_\zeta\nabla_\Sigma\eta_j\cdot\nabla_\Sigma F_0^{4,j} + \kappa\|\partial_t\eta_j\|_{L^2(\Gamma_0)}^2 \notag \\
&- \kappa\left\langle \partial_t\eta_j,F_0^{4,j} \right\rangle_{\Gamma_0} + \left\langle F_0^{7,j},\partial_t\eta_j-F_0^{4,j} \right\rangle_{\Gamma_0}+\sigma\int_{\Sigma}F_{0}^{5,j}\cdot \nabla_{\Sigma}(\partial_{t}\eta_{j}-F_{0}^{4,j}).
\end{align}

Combining the preceding identities, we obtain
\begin{align}
&\frac{d}{dt} \Bigg[ \frac12\int_\Omega J|v_j|^2 + \frac g2\|\eta_j\|_{L^2(\Sigma)}^2 + \frac{\sigma}{2}\int_\Sigma a_\zeta\nabla_\Sigma\eta_j\cdot\nabla_\Sigma\eta_j \Bigg] + \frac{\mu}{2}\int_\Omega J|\mathbb{D}_{\mathcal{A}_\eta}v_j|^2 + \beta_0\int_{\Sigma_s}J_s|U_0v_j|^2 + \kappa\|\partial_t\eta_j\|_{L^2(\Gamma_0)}^2 \notag \\
&= \frac12\int_\Omega(\partial_tJ)|v_j|^2 + \frac{\sigma}{2}\int_\Sigma (\partial_ta_\zeta) \nabla_\Sigma\eta_j\cdot\nabla_\Sigma\eta_j \notag \\
&+ \int_\Omega Jp_jF_0^{2,j} + \int_\Omega JF_0^{1,j}\cdot v_j - \left\langle F_0^{3,j},v_j\right\rangle_\Sigma - \left\langle J_sF_0^{6,j},U_0v_j\right\rangle_{\Sigma_s} \notag \\
&+ g\int_\Sigma\eta_jF_0^{4,j} + \sigma\int_\Sigma a_\zeta\nabla_\Sigma\eta_j\cdot\nabla_\Sigma F_0^{4,j} + \kappa\left\langle \partial_t\eta_j,F_0^{4,j} \right\rangle_{\Gamma_0} \notag \\
&- \left\langle F_0^{7,j},\partial_t\eta_j-F_0^{4,j} \right\rangle_{\Gamma_0} - \sigma\int_\Sigma F_0^{5,j}\cdot\nabla_\Sigma\partial_t\eta_j + \sigma\int_\Sigma F_0^{5,j}\cdot\nabla_\Sigma F_0^{4,j}. \label{eq:whole-domain-energy-j}
\end{align}

When $j=m+1$, taking the test function $u=(v_{j}-\omega)$, where $\omega$ is defined by Proposition~\ref{prop:B_d}, and using a computation similar to that for $j\leq m$, we obtain
\begin{align}
&\frac{d}{dt} \Bigg[ \frac12\int_\Omega J|v_j|^2 + \frac g2\|\eta_j\|_{L^2(\Sigma)}^2 + \frac{\sigma}{2}\int_\Sigma a_\zeta\nabla_\Sigma\eta_j\cdot\nabla_\Sigma\eta_j -\int_{\Omega}Jv_{j}\omega-\int_{\Omega}p_{j-1}\langle JF_{0}^{2,j}\rangle\Bigg]  \notag \\
&\quad+ \frac{\mu}{2}\int_\Omega J|\mathbb{D}_{\mathcal{A}_\eta}v_j|^2 + \beta_0\int_{\Sigma_s}J_s|U_0v_j|^2 + \kappa\|\partial_t\eta_j\|_{L^2(\Gamma_0)}^2  \notag \\
&\qquad= \frac12\int_\Omega(\partial_tJ)|v_j|^2 + \frac{\sigma}{2}\int_\Sigma (\partial_ta_\zeta) \nabla_\Sigma\eta_j\cdot\nabla_\Sigma\eta_j-\int_{\Omega}Jv_{j}\partial_{t}\omega+((v_{j},\omega))-\int_{\Omega}\partial_{t}Jv_{j}\omega  \notag \\
&\qquad\quad- \int_\Omega p_{j-1}\partial_{t}\langle JF_{0}^{2,j}\rangle + \int_\Omega JF_0^{1,j}\cdot (v_j-\omega) - \left\langle F_0^{3,j},v_j\right\rangle_\Sigma - \left\langle J_sF_0^{6,j},U_0v_j\right\rangle_{\Sigma_s}  \notag \\
&\qquad\quad+ g\int_\Sigma\eta_jF_0^{4,j} + \sigma\int_\Sigma a_\zeta\nabla_\Sigma\eta_j\cdot\nabla_\Sigma F_0^{4,j} + \kappa\left\langle \partial_t\eta_j,F_0^{4,j} \right\rangle_{\Gamma_0}  \notag \\
&\qquad\quad- \left\langle F_0^{7,j},\partial_t\eta_j-F_0^{4,j} \right\rangle_{\Gamma_0} - \sigma\int_\Sigma F_0^{5,j}\cdot\nabla_\Sigma\partial_t\eta_j + \sigma\int_\Sigma F_0^{5,j}\cdot\nabla_\Sigma F_0^{4,j}.
\end{align}

Define
\begin{equation}
\mathcal{E}_{j,||} :={}  \frac12\int_\Omega J|v_{j}|^2 +\|\eta_{j}\|_{H^{1}(\Sigma)}^{2},\quad \mathcal{D}_{j,||} :={}  \frac{\mu}{2}\int_\Omega J|\mathbb{D}_{\mathcal{A}_\eta}v_{j}|^2 +\beta_0 \|U_0v_{j}\|_{L^2(\Sigma_s)}^2+\frac{\kappa}{2}\int_{\Gamma}|\partial_{t}\eta_{j}|^{2}.
\end{equation}
Here the wall surface Jacobian $J$ is close to one. In particular,  $\mathcal{E}_{j,||} \asymp \|v_{j}\|_{L^2(\Omega)}^2 +\|\eta_{j}\|_{H^1(\Sigma)}^2$. For the terms on the left-hand side, we have
\begin{align}
\frac12\int_\Omega J|v_j|^2 + \frac g2\|\eta_j\|_{L^2(\Sigma)}^2 + \frac{\sigma}{2}\int_\Sigma a_\zeta\nabla_\Sigma\eta_j\cdot\nabla_\Sigma\eta_j\lesssim \mathcal{E}_{j,||},
\end{align}
and, conversely, $\mathcal{D}_{j,||}$ controls $\|v_{j}\|_{H^{1}}$. Indeed, by Korn's inequality with the slip term, for $v_{j}\cdot n_{0}=0$ on $\Sigma_{s}$ and $\eta$ small,
\begin{align}
\|v_{j}\|_{H^{1}(\Omega)}^{2}\lesssim\frac{\mu}{2}\int_\Omega J|\mathbb{D}_{\mathcal{A}_\eta}v_{j}|^2 +\beta_0 \|U_0v_{j}\|_{L^2(\Sigma_s)}^2\lesssim\mathcal{D}_{j,||},
\end{align}
with a constant that depends on $1/\beta_{0}$; this is where $\beta_{0}>0$ is used. The inequality is Lemma~\ref{lem:korn} below, combined with $|\mathcal{A}_{\eta}-I|+|J-1|\lesssim\sqrt{\mathcal{E}}$, which gives $\|\mathbb{D}v_{j}\|_{L^{2}}\leq\|\mathbb{D}_{\mathcal{A}_{\eta}}v_{j}\|_{L^{2}}+C\sqrt{\mathcal{E}}\|v_{j}\|_{H^{1}}$; the last term is absorbed for $\mathcal{E}$ small.

\begin{lemma}[Korn's inequality with the slip term]
\label{lem:korn} Let $\Omega$ be a bounded Lipschitz domain and let $\Sigma_{s}\subset\partial\Omega$ be a relatively open part of the boundary that is not contained in a line. There is $C_{K}$, depending only on $\Omega$ and $\Sigma_{s}$, such that
\begin{equation*}
\|u\|_{H^{1}(\Omega)}^{2}\leq C_{K}\Big(\|\mathbb{D}u\|_{L^{2}(\Omega)}^{2}+\|U_{0}u\|_{L^{2}(\Sigma_{s})}^{2}\Big)
\end{equation*}
for every $u\in H^{1}(\Omega;\mathbb{R}^{3})$ with $u\cdot n_{0}=0$ on $\Sigma_{s}$.
\end{lemma}

\begin{proof}
Suppose not. Then there are $u_{k}$ with $u_{k}\cdot n_{0}=0$ on $\Sigma_{s}$, $\|u_{k}\|_{H^{1}}=1$ and $\|\mathbb{D}u_{k}\|_{L^{2}}+\|U_{0}u_{k}\|_{L^{2}(\Sigma_{s})}\to0$. By Rellich's theorem a subsequence converges in $L^{2}(\Omega)$. Korn's second inequality on Lipschitz domains, $\|u\|_{H^{1}}\lesssim\|\mathbb{D}u\|_{L^{2}}+\|u\|_{L^{2}}$ \cite{Nitsche}, shows that the subsequence is Cauchy in $H^{1}(\Omega)$. Its limit $u$ has $\|u\|_{H^{1}}=1$ and $\mathbb{D}u=0$, so $u(x)=a+b\times x$ is a rigid motion. By continuity of the trace, $u\cdot n_{0}=0$ and $U_{0}u=0$ on $\Sigma_{s}$, so $u=0$ on $\Sigma_{s}$. A rigid motion that vanishes on a set not contained in a line is zero, a contradiction.
\end{proof}

The wall $\Sigma_{s}$ contains a vertical cylinder over $\partial D$, so the lemma applies, and its constant does not depend on $\varepsilon$. It also applies, up to lower-order terms, to the localized fields $v_{k_{1},k_{2}}$ of Section~\ref{sec:contact_energy}, extended by zero: their normal component on $\Sigma_{s}$ equals that of the correction $\omega$ in \eqref{eq:local-domain-energy-1-2}, up to sign, which involves fewer angular derivatives.

We now state the main theorem.

\begin{proof}[Proof of Theorem~\protect\ref{thm:whole-space-energy}]
We now estimate every term on the right-hand side. For the additional terms when $j=m+1$, we give their estimates in Proposition~\ref{prop:B_d}. The Jacobian term satisfies
\begin{equation}
\left| \frac12\int_\Omega (\partial_tJ)|v_j|^2 \right| \leq \frac12\|\partial_tJ\|_{L^\infty(\Omega)} \|v_{j}\|_{L^2(\Omega)}^2 \leq C\|\partial_tJ\|_{L^\infty(\Omega)} \mathcal{E}_{j,||},
\end{equation}
where
\begin{align}
\|\partial_{t}J\|_{L^{\infty}}\lesssim \|\partial_{t}\nabla \bar{\eta}\|_{L^{\infty}}&\lesssim \|\chi_{1,\mathrm{int}}\partial_{t}\eta\|_{W^{3-\frac{1}{p}+(m-1),p}}+\|\chi_{\mathrm{cl},\Sigma}\partial_{t}\eta\|_{L^{\infty}W^{3-\frac{1}{q_{\delta}},p}}+\|\chi_{\mathrm{cl},\Sigma}\partial_{t}\partial_{\theta}\eta\|_{L^{\infty}H^{\frac{3}{2}-\alpha}}  \notag \\
&\lesssim\sqrt{\mathcal{D}_{\mathrm{int}}}+\sqrt{\mathcal{D}_{\mathrm{cl}}}, \label{eq:dtJ_Linfty}
\end{align}
given that $m\geq 1$. Since $\|v_j\|_{L^2}^2\lesssim\mathcal{E}$, the term is bounded by $\mathcal{E}\sqrt{\mathcal{D}}$, and $\mathcal{E}\sqrt{\mathcal{D}}\lesssim\sqrt{\mathcal{E}}\,\mathcal{D}$ by Lemma~\ref{lem:E_le_D}. Therefore,
\begin{align}
\left| \frac12\int_\Omega (\partial_tJ)|v_j|^2 \right|\lesssim \sqrt{\mathcal{E}}\mathcal{D}.
\end{align}

\textbf{Estimate for the term including $F_{0}^{2,j}$}

For every $\varepsilon>0$, and $j\leq m$, the pressure term satisfies
\begin{align}
\left| \int_\Omega Jp_{j}F_{0}^{2,j} \right| &\leq C\|p_{j}\|_{L^2(\Omega)} \|F_{0}^{2,j}\|_{L^2(\Omega)}.
\end{align}
For $F_{0}^{2,j}$, we separate it by applying $\chi_{cl+}$ and $\chi_{int+}$. We have  $F_{0}^{2,j}=F_{\mathrm{cl}+}^{2,j}+F_{\mathrm{int}+}^{2,j}$, where
\begin{equation*}
F_{\mathrm{int}+}^{2,j} = -\chi_{\mathrm{int}+} \sum_{k=1}^{j}\binom{j}{k} \func{div}_{\partial_t^k\mathcal{A}_\eta}v_{j-k} ,
\end{equation*}
and
\begin{equation*}
F_{\mathrm{cl}+}^{2,j} = -\chi_{\mathrm{cl}+} \sum_{k=1}^{j}\binom{j}{k} \func{div}_{\partial_t^k\mathcal{A}_\eta}v_{j-k}.
\end{equation*}
For $\|F_{\mathrm{cl}+}^{2,j}\|_{L^{2}(\Omega)}$, we estimate it by Lemma \ref{thm:contact_line_forcing_estimate} from the fact that $\func{supp}\{\chi_{\mathrm{cl}+,\Sigma}\}\subseteq\{x\in \Sigma| \chi_{\mathrm{cl},\Sigma}+\chi_{\mathrm{bet},\Sigma}=1\}$.

For $\|F_{\mathrm{int}+}^{2,j}\|_{L^{2}(\Omega)}$ with any $j\leq m+1$, by H\"older's inequality, we have
\begin{align}
&\|F_{\mathrm{int}+}^{2,j}\|_{L^2(\Omega)} \lesssim{} \sum_{k=1}^{j-1}\binom{j}{k} \left\| \chi_{\mathrm{int}+} (\partial_t^k\mathcal{A}_\eta)\nabla v_{j-k} \right\|_{L^2(\Omega)} + \left\| \chi_{\mathrm{int}+} (\partial_t^{j}\mathcal{A}_\eta)\nabla v \right\|_{L^2(\Omega)} \notag \\
&\lesssim{} \sum_{k=1}^{j-1}\binom{j}{k} \|\chi_{0,\mathrm{int}}\partial_t^k\mathcal{A}_\eta\|_{L^{\frac{6}{1+2\alpha}}(\Omega)} \|\nabla (\chi_{\mathrm{0,\mathrm{int}}}v_{j-k})\|_{L^{\frac{3p}{3-p}}(\Omega)}+\sum_{k=1}^{j-1}\binom{j}{k} \|\chi_{0,\mathrm{int}}\partial_t^k\mathcal{A}_\eta\|_{L^3(\Omega)} \| \chi_{0,\mathrm{int}}v_{j-k}\|_{L^6(\Omega)} \notag \\
&+\left\| (\chi_{0,\mathrm{int}}\partial_t^{j}\mathcal{A}_\eta)\right\|_{L^{2}(\Omega)}\left\|\nabla (\chi_{0,\mathrm{int}}v) \right\|_{L^{\infty}(\Omega)}+\left\| (\chi_{0,\mathrm{int}}\partial_t^{j}\mathcal{A}_\eta)\right\|_{L^{2}(\Omega)}\left\|\chi_{0,\mathrm{int}}v \right\|_{L^{\infty}(\Omega)} \notag \\
&\lesssim{} \sum_{k=1}^{j-1}\binom{j}{k} (\|\partial_t^k(\nabla\bar{\eta}_{0,\mathrm{int},d})\|_{L^{\infty}L^{\frac{6}{1+2\alpha}}(\Omega)}+\|\partial_{t}^{k}\bar{\eta}_{\mathrm{int}}\|_{L^{\frac{6}{1+2\alpha}}}) \|\nabla (\chi_{0,\mathrm{int}}v_{j-k})\|_{L^{\frac{3p}{3-p}}(\Omega)} \notag \\
&+\sum_{k=1}^{j-1}\binom{j}{k} (\|\partial_t^k(\nabla\bar{\eta}_{0,\mathrm{int},d})\|_{L^{\infty}L^{\frac{6}{1+2\alpha}}(\Omega)}+\|\partial_{t}^{k}\bar{\eta}_{\mathrm{int}}\|_{L^{\frac{6}{1+2\alpha}}}) \| \chi_{0,\mathrm{int}}v_{j-k}\|_{L^6(\Omega)} \notag \\
&+\left\| (\partial_t^{j}((\nabla \bar{\eta}_{\mathrm{int}}+\nabla\bar{\eta}_{0,\mathrm{int},d})))\right\|_{L^{2}(\Omega)}\left\|\nabla (\chi_{0,\mathrm{int}}v) \right\|_{L^{\infty}(\Omega)} +\left\| (\partial_t^{j}((\nabla \bar{\eta}_{\mathrm{int}}+\nabla\bar{\eta}_{0,\mathrm{int},d})))\right\|_{L^{2}(\Omega)}\left\|\chi_{0,\mathrm{int}}v \right\|_{L^{\infty}(\Omega)} \notag \\
&\lesssim{} \sum_{k=1}^{j-1}\binom{j}{k} (\|\partial_t^k\partial_{\theta}((\chi_\mathrm{bet}){\eta})\|_{H^{\frac{3}{2}-\alpha}(\Sigma)}+\|\partial_t^k((\chi_\mathrm{bet}){\eta})\|_{L^{2}H^{\frac{3}{2}-\alpha}(\Sigma)}+\|\partial_{t}^{k}\eta_{\mathrm{int}}\|_{H^{\frac{3}{2}-\alpha}(\Sigma)}) \| \chi_{0,\mathrm{int}}v_{j-k}\|_{W^{2,p}(\Omega)} \notag \\
&+\sum_{k=1}^{j-1}\binom{j}{k} (\|\partial_t^k\partial_{\theta}((\chi_\mathrm{bet}){\eta})\|_{H^{\frac{3}{2}-\alpha}(\Sigma)}+\|\partial_t^k((\chi_\mathrm{bet}){\eta})\|_{L^{2}H^{\frac{3}{2}-\alpha}(\Sigma)}+\|\partial_{t}^{k}\eta_{\mathrm{int}}\|_{H^{\frac{3}{2}-\alpha}(\Sigma)}) \| \chi_{1,\mathrm{int}}v_{j-k}\|_{W^{2,p}(\Omega)} \notag \\
&+\left\| (\partial_t^{j}((\chi_{\mathrm{int},\Sigma}+\chi_{\mathrm{bet},\Sigma}){\eta}))\right\|_{H^{1}(\Sigma)}\left\|v_{0,\mathrm{int}} \right\|_{W^{3,p}(\Omega)} +\left\| (\partial_t^{j}((\chi_{\mathrm{int},\Sigma}+\chi_{\mathrm{bet},\Sigma}){\eta}))\right\|_{H^{1}(\Sigma)}\left\|\chi_{1,\mathrm{int}}v \right\|_{W^{2,p}(\Omega)} \notag \\
&\lesssim \sum_{k=1}^{j-1}(\sqrt{\mathscr{E}_{k,\mathrm{energy},\mathrm{int}}}+\sqrt{\mathscr{E}_{k,\mathrm{energy},\mathrm{bet}}})\sqrt{\mathscr{D}_{\mathrm{int}}} +(\sqrt{\mathscr{E}_{j,\mathrm{energy},\mathrm{int}}}+\sqrt{\mathscr{E}_{j,\mathrm{energy},\mathrm{bet}}})\sqrt{\mathscr{D}_{\mathrm{int}}} \lesssim \sqrt{\mathcal{E}}\sqrt{\mathcal{D}},
\end{align}
where we used the fact that $H^{1-\alpha}\hookrightarrow L^{\frac{6}{1+2\alpha}}$ in two-dimensional space.

When $j\leq m$, the terms above, with one fewer derivative, can all be bounded by the energy functional $\sqrt{\mathcal{E}}$. Therefore, when $j\leq m$, we have
\begin{align}
\left| \int_\Omega Jp_{j}F_{0}^{2,j} \right| &\leq C\|p_{j}\|_{L^2(\Omega)} \|F_{0}^{2,j}\|_{L^2(\Omega)}\lesssim \sqrt{\mathcal{E}}^{3}.
\end{align}

When $j=m$, we deal with this term in a different way.

\textbf{Estimate for $F_{\mathrm{int}+}^{1,j}$}

We now estimate the interior part of $F_{\mathrm{int}+}^{1,j}$. The contact line part can be estimated by Lemma \ref{thm:contact_line_forcing_estimate}. Recall that
\begin{align}
&F_{\mathrm{int}+}^{1,j} ={}  \chi_{\mathrm{int}+} \partial_t^j \left( WK\partial_t\bar\eta\,\partial_3v -v\cdot\nabla_{\mathcal{A}_\eta}v \right) + \mu\chi_{\mathrm{int}+} \sum_{k=1}^{j}\binom{j}{k} \func{div}_{\mathcal{A}_\eta} \left( \mathbb{D}_{\partial_t^k\mathcal{A}_\eta}v_{j-k} \right) \notag \\
&- \chi_{\mathrm{int}+} \sum_{k=1}^{j}\binom{j}{k} \func{div}_{\partial_t^k\mathcal{A}_\eta} \left[ \partial_t^{j-k}S_{\mathcal{A}_\eta}(p,v) \right].
\end{align}

First, by $H^1(\Omega)\hookrightarrow L^6(\Omega)$ and H\"older's inequality,
\begin{align}
&\left\| J\chi_{\mathrm{int}+} \partial_t^j \left( WK\partial_t\bar\eta\,\partial_3v \right) \right\|_{(H^1)^*} \lesssim \sum_{\substack{ a+b+c=j \\ a,b,c\geq0 \\ a,b,c\neq j}} \frac{j!}{a!\,b!\,c!} \left\| \chi_{0-,\mathrm{int}}\, \partial_t^a(WK)\, \partial_t^{b+1}\bar\eta \right\|_{L^3} \left\| \chi_{\mathrm{int}+}\nabla\partial_t^c v \right\|_{L^2} \notag \\
&+ \left\| \chi_{0-,\mathrm{int}}\, \partial_{t}^{j}(WK)\, \partial_t^{1}\bar\eta \right\|_{L^3} \left\| \chi_{\mathrm{int}}\nabla v \right\|_{L^2}+ \left\| \chi_{0-,\mathrm{int}}\, (WK)\, \partial_t^{1}\bar\eta \right\|_{L^3} \left\| \chi_{\mathrm{int}}\nabla\partial_t^j v \right\|_{L^2} \notag \\
&+\left\| \chi_{0-,\mathrm{int}}\, (WK)\, \partial_t^{j+1}\bar\eta \right\|_{L^2} \left\| \chi_{\mathrm{int}}\nabla v\right\|_{L^{\infty}} \notag \\
&\lesssim \sum_{\substack{ a+b+c=j \\ a,b,c\geq0 \\ a,b,c\neq j}} \frac{j!}{a!\,b!\,c!} (1+\sum_{d=0}^{a}\left\| \chi_{0,\mathrm{int}}\, \partial_t^d\bar{\eta}\,\right\|_{W^{1,\infty}}) \left\|(\chi_{0,\mathrm{int}})\partial_t^{b+1}\bar\eta \right\|_{L^3} \left\| (\chi_{\mathrm{int}}+\chi_{\mathrm{bet}})\partial_t^c v \right\|_{W^{3,p}} \notag \\
&+(1+\sum_{d=0}^{j}\left\| \chi_{0-,\mathrm{int}}\, \partial_t^{d}\nabla \bar\eta \right\|_{L^{3}}) \left\| \chi_{0-,\mathrm{int}}\, \partial_t^{1}\bar\eta \right\|_{L^{\infty}} \left\| \chi_{\mathrm{int}}\nabla v \right\|_{L^2}+ \left\| \chi_{0+,\mathrm{int}}\, \partial_t^{1}\bar\eta \right\|_{L^3} \left\| (\chi_{\mathrm{int}}+\chi_{\mathrm{bet}})\partial_t^j v \right\|_{H^1} \notag \\
&+\left\| \chi_{0+,\mathrm{int}}\, \partial_t^{j+1}\bar\eta \right\|_{L^2} \left\| (\chi_{1,\mathrm{int}}) v\right\|_{W^{3,p}} \notag \\
&\lesssim \sum_{d=0}^{j-1} (1+\sum_{k=0}^{j-1}(\|\partial_{t}^{d}\bar{\eta}_{1,\mathrm{int},d}\|_{L^{\infty}W^{1,+\infty}}+\|\partial_{t}^{d}\partial_{\theta}\bar{\eta}_{1,\mathrm{int},d}\|_{L^{\infty}L^{+\infty}}+\|\partial_{t}^{d}\eta_{\mathrm{int}}\|_{W^{3-\frac{1}{p},p}})) \notag \\
&\times \left\|(\chi_{\mathrm{int}}+\chi_{\mathrm{bet}})\partial_t^{b+1}\eta \right\|_{H^{1}}(\mathscr{E}_{\mathrm{int}}^{\frac{1}{2}}+\mathscr{E}_{m,\mathrm{energy},\mathrm{bet}}^{\frac{1}{2}}) \notag \\
&+(1+\sum_{d=0}^{j}(\left\| \, \partial_t^{d} \eta_{int} \right\|_{H^{1}}+\|\partial_{t}^{d}\bar{\eta}_{1,\mathrm{int},d}\|_{L^{3}H^{\frac{3}{2}}}+\|\partial_{t}^{d}\partial_{\theta}\bar{\eta}_{1,\mathrm{int},d}\|_{L^{3}H^{\frac{1}{2}}})) \left\| \chi_{0,\mathrm{int}}\, \partial_t^{1}\eta \right\|_{W^{3-\frac{1}{p},p}} (\mathscr{E}^{\frac{1}{2}}_{\mathrm{int}}+\mathscr{E}^{\frac{1}{2}}_{m,\mathrm{energy},\mathrm{bet}}) \notag \\
&+ \left\| (\chi_{\mathrm{int}}+\chi_{\mathrm{bet}})\, \partial_t^{1}\eta \right\|_{H^1} (\mathscr{D}^{\frac{1}{2}}_{\mathrm{int}}+\mathscr{D}^{\frac{1}{2}}_{m,\mathrm{energy},\mathrm{bet}}) +\left\| (\chi_{\mathrm{int}}+\chi_{\mathrm{bet}})\, \partial_t^{j+1}\eta \right\|_{H^{\frac{1}{2}-\alpha}} (\mathscr{E}^{\frac{1}{2}}_{\mathrm{int}}+\mathscr{E}^{\frac{1}{2}}_{m,\mathrm{energy},\mathrm{bet}}) \notag \\
&\lesssim \sqrt{\mathcal{E}\mathcal{D}}+(1+\sum_{d=0}^{j}(\left\| \, \partial_t^{d} \eta_{int} \right\|_{H^{1}}+\|\partial_{t}^{d}{\eta}_{\mathrm{bet}}\|_{H^{\frac{3}{2}-\alpha}}+\|\partial_{t}^{d}\partial_{\theta}{\eta}_{\mathrm{bet}}\|_{H^{\frac{1}{2}-\alpha}})\mathcal{E} \notag \\
&+\sum_{d=0}^{j-1}(1+\left\| \, \partial_t^{d} \eta_{int} \right\|_{W^{3-\frac{1}{p},p}}+\|\partial_{t}^{d}{\eta}_{\mathrm{bet}}\|_{L^{\infty}H^{\frac{3}{2}+\alpha}}+\|\partial_{t}^{d}\partial_{\theta}\bar{\eta}_{\mathrm{bet}}\|_{L^{\infty}H^{\frac{1}{2}+\alpha}}))\mathcal{E} \notag \\
&\lesssim\sqrt{\mathcal{E}\mathcal{D}}+\sum_{d=0}^{j-1}\mathcal{E}+\sum_{d=0}^{j-1}(1+\left\| \partial_t^{d} \eta_{int} \right\|_{W^{3-\frac{1}{p},p}}+\|\partial_{t}^{d}{\eta}_{\mathrm{bet}}\|_{L^{\infty}H^{\frac{3}{2}+\frac{\varepsilon_{-}-\alpha}{2}}}+\|\partial_{t}^{d}\partial_{\theta}\bar{\eta}_{\mathrm{bet}}\|_{H^{1+2\alpha}}))\mathcal{E} \lesssim \sqrt{\mathcal{E}\mathcal{D}}, \label{eq:Fint1-transport}
\end{align}
where $\chi_{1,\mathrm{int}}=1$ on $\func{supp}\chi_{\mathrm{int}+}$ and we used the fact that $H^{\frac{1}{2}}\hookrightarrow L^{3}$ in three-dimensional space, and the fact that $H_{\theta}^{\frac{1}{2}+\alpha}\hookrightarrow L_{\theta}^{\infty}$ in one-dimensional space, and the interpolation of $L^{\infty}$. In the last step, $L^{\infty}_{\theta}H^{\frac{3}{2}+\alpha}$ is replaced by $L^{\infty}_{\theta}H^{\frac{3}{2}+\frac{\varepsilon_{-}-\alpha}{2}}$, which uses $3\alpha\leq\varepsilon_{-}$ (Subsection~\ref{subsec:constants}), and
\begin{equation*}
\|\partial_{t}^{d}\partial_{\theta}\bar{\eta}_{\mathrm{bet}}\|_{L^{\infty}_{\theta}H^{\frac{1}{2}+\alpha}}\lesssim\|\partial_{t}^{d}\partial_{\theta}\bar{\eta}_{\mathrm{bet}}\|_{H^{1+2\alpha}}
\end{equation*}
, since the supremum in $\theta$ costs more than half a derivative.

We next estimate the convection term. Assume that $\chi_{0-,\mathrm{int}}=1$ on $\func{supp}\chi_{\mathrm{int}+}$ and that $J$ and $\mathcal{A}_\eta$ are uniformly bounded in $L^\infty$. By the Leibniz rule, H\"older's inequality, and the embedding $H^1(\Omega)\hookrightarrow L^6(\Omega)$, we obtain, for $1\leq j\leq m$,
\begin{align}
&\left\| J\chi_{\mathrm{int}+} \partial_t^j \left(v\cdot\nabla_{\mathcal{A}_\eta}v\right) \right\|_{(H^1)^*} \lesssim \sum_{\substack{ a+b+c=j \\ 0\leq a,b,c<j}} \frac{j!}{a!\,b!\,c!} \left\| \chi_{0-,\mathrm{int}}\partial_t^a v \right\|_{L^3} \left\| \chi_{0-,\mathrm{int}}\partial_t^b\mathcal{A}_\eta \right\|_{L^\infty} \left\| \chi_{\mathrm{int+}}\nabla\partial_t^c v \right\|_{L^2} \notag \\
&+ \left\| \chi_{0-,\mathrm{int}}\partial_t^jv \right\|_{L^3} \|\chi_{0-,\mathrm{int}}\mathcal{A}_\eta\|_{L^\infty} \left\| \chi_{\mathrm{int}+}\nabla v \right\|_{L^2} + \left\| \chi_{0-,\mathrm{int}}v \right\|_{L^\infty} \left\| \chi_{0-,\mathrm{int}}\partial_t^j\mathcal{A}_\eta \right\|_{L^3} \left\| \chi_{\mathrm{int}+}\nabla v \right\|_{L^2} \notag \\
&+ \left\| \chi_{0-,\mathrm{int}}v \right\|_{L^3} \|\chi_{0-,\mathrm{int}}\mathcal{A}_\eta\|_{L^\infty} \left\| \chi_{\mathrm{int}+}\nabla\partial_t^jv \right\|_{L^2}.
\end{align}
Using the estimate for the previous term, we know that $\|\partial_{t}^{k}\nabla \chi_{0+,\mathrm{int}}\bar{\eta}\|_{L^{\infty}}$ is bounded by $1+\sqrt{\mathcal{E}}$ for any $0\leq k\leq j-1$ with $0\leq j\leq m$, and that $\|\partial_{t}^{k}\chi_{0+,\mathrm{int}}\nabla\bar{\eta}\|_{L^{3}}$ is bounded by $1+\sqrt{\mathcal{D}}$ for any $0\leq k\leq j$ with $0\leq j\leq k $. Therefore, we have
\begin{align}
&\sum_{\substack{ a+b+c=j \\ 0\leq a,b,c<j}} \frac{j!}{a!\,b!\,c!} \left\| \chi_{0-,\mathrm{int}}\partial_t^a v \right\|_{L^3} \left\| \chi_{0-,\mathrm{int}}\partial_t^b\mathcal{A}_\eta \right\|_{L^\infty} \left\| \chi_{\mathrm{int}+}\nabla\partial_t^c v \right\|_{L^2} \notag \\
&+ \left\| \chi_{0-,\mathrm{int}}\partial_t^jv \right\|_{L^3} \|\chi_{0-,\mathrm{int}}\mathcal{A}_\eta\|_{L^\infty} \left\| \chi_{\mathrm{int}+}\nabla v \right\|_{L^2} + \left\| \chi_{0-,\mathrm{int}}v \right\|_{L^\infty} \left\| \chi_{0-,\mathrm{int}}\partial_t^j\mathcal{A}_\eta \right\|_{L^3} \left\| \chi_{\mathrm{int}+}\nabla v \right\|_{L^2} \notag \\
&+ \left\| \chi_{0-,\mathrm{int}}v \right\|_{L^3} \|\chi_{0-,\mathrm{int}}\mathcal{A}_\eta\|_{L^\infty} \left\| \chi_{\mathrm{int}+}\nabla\partial_t^jv \right\|_{L^2} \notag \\
&\lesssim \sum_{\substack{ n_{1}+n_{2}\leq j \\ 0\leq n_{1},n_{2}<j}} \left\| (\chi_{\mathrm{int}}+\chi_{\mathrm{bet}})\partial_t^{n_{1}} v \right\|_{H^{1}} \left\| (\chi_{\mathrm{int}}+\chi_{\mathrm{bet}})\partial_t^{n_{2}} v \right\|_{H^{1}}(1+\sqrt{\mathcal{E}}) \notag \\
&+ \left\| (\chi_{\mathrm{int}}+\chi_{\mathrm{bet}})\partial_t^jv \right\|_{H^{1}} \left\| (\chi_{\mathrm{int}}+\chi_{\mathrm{bet}})v \right\|_{H^1} \notag \\
&+ \left\| \chi_{0,\mathrm{int}}v \right\|_{W^{2,p}} \left( \sum_{r=1}^{j} \sum_{\substack{ \ell_1+\cdots+\ell_r=j \\ \ell_i\geq1,\ell_{1}=\max\{\ell_{i}\}}} \left\| (\chi_{0,\mathrm{int}})\partial_t^{\ell_1}\bar{\eta} \right\|_{W^{1,3}} \prod_{i=2}^{r} \left\| (\chi_{0,\mathrm{int}})\partial_t^{\ell_i}\bar\eta \right\|_{W^{1,\infty}} \right) \left\| (\chi_{\mathrm{int}}+\chi_{\mathrm{bet}})v \right\|_{H^1} \lesssim \sqrt{\mathcal{E}\mathcal{D}},
\end{align}
where we omit the $\sqrt{\mathcal{E}}^{i}\sqrt{\mathcal{D}}$ for $i\geq 2$ since $\sqrt{\mathcal{E}}$ is small so that this can be controlled by $\sqrt{\mathcal{E}\mathcal{D}}$.

For the remaining temporal commutator, by $H^1(\Omega)\hookrightarrow L^6(\Omega)$ and H\"older's inequality, for $j\geq1$, the Leibniz rule gives
\begin{align}
&\sum_{k=0}^{j}\binom{j}{k} \func{div}_{\partial_t^k\mathcal{A}_\eta} \left[ \partial_t^{j-k}S_{\mathcal{A}_\eta}(p,v) \right] = \func{div}_{\mathcal{A}_\eta} S_{\mathcal{A}_\eta}(p_j,v_j) + \func{div}_{\partial_t^j\mathcal{A}_\eta} S_{\mathcal{A}_\eta}(p,v) - \mu\func{div}_{\mathcal{A}_\eta} \left( \mathbb{D}_{\partial_t^j\mathcal{A}_\eta}v \right) \notag \\
&+ \sum_{k=1}^{j-1}\binom{j}{k} \nabla_{\partial_t^k\mathcal{A}_\eta}p_{j-k} - \mu \sum_{\substack{ a+b+c=j \\ 0\leq a,b,c<j}} \frac{j!}{a!\,b!\,c!} \func{div}_{\partial_t^a\mathcal{A}_\eta} \left( \mathbb{D}_{\partial_t^b\mathcal{A}_\eta}v_c \right).
\end{align}
The term $\func{div}_{\mathcal{A}_\eta} S_{\mathcal{A}_\eta}(p_j,v_j)$ is included in the principal Stokes operator. We estimate the remaining terms directly, without integrating by parts. Using $L^{6/5}(\Omega)\hookrightarrow(H^1(\Omega))^*$, H\"older's inequality, $W^{3,p}\hookrightarrow H^2(\Omega)\hookrightarrow W^{1,6}(\Omega)$, and the previous bound for $\|\partial_{t}^{k}\nabla\bar{\eta}\|_{L^{\infty}}, \|\partial_{t}^{j}\nabla \bar{\eta}\|_{L^{3}}$ with $k\leq j-1$, we obtain
\begin{align}
&\left\| J\chi_{\mathrm{int}+} \left\{ \sum_{k=0}^{j}\binom{j}{k} \func{div}_{\partial_t^k\mathcal{A}_\eta} \left[ \partial_t^{j-k}S_{\mathcal{A}_\eta}(p,v) \right] - \func{div}_{\mathcal{A}_\eta} S_{\mathcal{A}_\eta}(p_j,v_j) + \mu\func{div}_{\mathcal{A}_\eta} \left( \mathbb{D}_{\partial_t^j\mathcal{A}_\eta}v \right)\right\} \right\|_{(H^1)^*} \notag \\
&\lesssim \left\| \chi_{0-,\mathrm{int}}\partial_t^j\mathcal{A}_\eta \right\|_{L^3} \Bigg( \left\|\chi_{\mathrm{int}+}\nabla p\right\|_{L^2} + \mu\left\|\chi_{0-,\mathrm{int}}\mathcal{A}_\eta\right\|_{L^\infty} \left\|\chi_{\mathrm{int}+}\nabla^2v\right\|_{L^2} + \mu\left\|\chi_{0-,\mathrm{int}}\nabla\mathcal{A}_\eta\right\|_{L^3} \left\|\chi_{\mathrm{int}+}\nabla v\right\|_{L^6} \Bigg) \notag \\
&+ \sum_{k=2}^{j-1}\binom{j}{k} \left\|\chi_{0-,\mathrm{int}}\partial_t^k\mathcal{A}_\eta\right\|_{L^3} \left\|\chi_{\mathrm{int}+}\nabla p_{j-k}\right\|_{L^2}+\left\|\chi_{0-,\mathrm{int}}\partial_t\mathcal{A}_\eta\right\|_{L^{\infty}} \left\|\chi_{\mathrm{int}+}\nabla p_{j-1}\right\|_{L^p} \notag \\
&+ \mu\sum_{\substack{ a+b+c=j \\ 0\leq a,b,c<j}} \frac{j!}{a!\,b!\,c!} \left\|\chi_{0-,\mathrm{int}}\partial_t^a\mathcal{A}_\eta\right\|_{L^\infty} \Bigg( \left\|\chi_{0-,\mathrm{int}}\partial_t^b\mathcal{A}_\eta\right\|_{L^{\frac{6p}{5p-6}}} \left\|\chi_{\mathrm{int}+}\nabla^2v_c\right\|_{L^{p}} \notag \\
&+ \left\|\chi_{0-,\mathrm{int}}\nabla\partial_t^b\mathcal{A}_\eta\right\|_{L^{\frac{2}{1+2\alpha}}} \left\|\chi_{\mathrm{int}+}\nabla v_c\right\|_{L^{\frac{3}{1-3\alpha}}} \Bigg) \lesssim \left\|\chi_{0-,\mathrm{int}}\partial_t^j\mathcal{A}_\eta\right\|_{L^3} \left( \|\chi_{1,\mathrm{int}}p\|_{W^{2,p}} + \mu\|\chi_{1,\mathrm{int}}v\|_{W^{3,p}} \right) \notag \\
&+ \sum_{k=2}^{j-1}\binom{j}{k} \left\|\chi_{0-,\mathrm{int}}\partial_t^k\mathcal{A}_\eta\right\|_{L^3} \|\chi_{1,\mathrm{int}}p_{j-k}\|_{W^{2,p}}+\left\|\chi_{0-,\mathrm{int}}\partial_t\mathcal{A}_\eta\right\|_{L^{\infty}} \|\chi_{1,\mathrm{int}}p_{j-1}\|_{W^{1,p}} \notag \\
&+ \mu\sum_{\substack{ a+b+c=j \\ 0\leq a,b,c<j}} \frac{j!}{a!\,b!\,c!} \left\|\chi_{0-,\mathrm{int}}\partial_t^a\mathcal{A}_\eta\right\|_{L^\infty} \left( \left\|\chi_{0-,\mathrm{int}}\partial_t^b\mathcal{A}_\eta\right\|_{L^{4}} + \left\|\chi_{0-,\mathrm{int}}\nabla\partial_t^b\mathcal{A}_\eta\right\|_{L^{\frac{2}{1+2\alpha}}} \right) \|\chi_{0,\mathrm{int}}v_{c}\|_{W^{2,p}} \lesssim \sqrt{\mathcal{E}\mathcal{D}}. \label{eq:full-stress-commutator-estimate}
\end{align}
In the last sum $b\leq m$ and $c\leq m$. To estimate the data for $\nabla\partial_t^b\mathcal A_\eta$ on $\func{supp}\chi_{0-,\mathrm{int}}$, we distinguish the energy and dissipation bounds. For $b\leq m-1$, the term $k=n=b$ of \eqref{def:interior_energy} controls $\chi_{b+1,\mathrm{int},\Sigma}\partial_t^b\eta \;\hbox{in }W^{2+m-b-1/p,p}(\Sigma)$. For $b\leq m$, the term $k=n=b$ of \eqref{def:interior_dissipation} controls $\chi_{b,\mathrm{int},\Sigma}\partial_t^b\eta \;\hbox{in }W^{3+m-b-1/p,p}(\Sigma)$. By nesting, the respective enlarged cutoffs equal one on the data supports under consideration; multiplication by the fixed cutoff profiles costs at most $C(\varepsilon)$. In their respective index ranges, both orders are at least $3-1/p$. Since $p>3/2$, $W^{3-1/p,p}(\Sigma)\hookrightarrow H^{4-3/p}(\Sigma) \hookrightarrow H^2(\Sigma)$. Thus the required $H^2$ data norm is bounded by $C(\varepsilon)\sqrt{\mathcal E}$ for $b\leq m-1$ and by $C(\varepsilon)\sqrt{\mathcal D}$ for $b\leq m$. The three-dimensional Poisson extension maps $H^{2}(\Sigma) $ into $H^{5/2}$ near $\Sigma$. The slice-wise extension maps it into $H^{2}$ near $\Sigma$, by \eqref{tl:P2} and $|k|^{4}\langle\xi\rangle^{-1}\leq\langle(k,\xi)\rangle^{4}$. Hence $\|\chi_{0-,\mathrm{int}}\nabla\partial_{t}^{b}\mathcal{A}_{\eta}\|_{L^{2/(1+2\alpha)}}\lesssim\|\chi_{0-,\mathrm{int}}\nabla\partial_{t}^{b}\mathcal{A}_{\eta}\|_{L^{2}}$ is bounded by $\sqrt{\mathcal{D}}$, and by $\sqrt{\mathcal{E}}$ if $b\leq m-1$. Also $W^{2,p}\subset W^{1,3/(1-3\alpha)}$ in three dimensions, since $p\geq3/(2-3\alpha)$, and $\frac{1+2\alpha}{2}+\frac{1-3\alpha}{3}=\frac56$, so the product lies in $L^{6/5}\subset(H^{1})^{*}$. In each product one factor is at an energy level, which gives $\sqrt{\mathcal{E}\mathcal{D}}$.

For the term $- \mu\func{div}_{\mathcal{A}_\eta} \left( \mathbb{D}_{\partial_t^j\mathcal{A}_\eta}v \right)$, we have
\begin{align}
&{}  |\int_{\Omega}J\chi_{\mathrm{int}+}\mu\func{div}_{\mathcal{A}_\eta} \left( \mathbb{D}_{\partial_t^j\mathcal{A}_\eta}v \right)\cdot \partial_{t}^{j}v| \lesssim  |\int_{\partial\Omega} \mu \chi_{\mathrm{int}+}\mathbb{D}_{\partial_{t}^{j}\mathcal{A}_{\eta}}v\mathcal{N}\cdot\partial_{t}^{j}v|+|\int_{\Omega}J\mu \left( \mathbb{D}_{\partial_t^j\mathcal{A}_\eta}v \right): \mathbb{D}_{\mathcal{A}}\partial_{t}^{j}(\chi_{\mathrm{int}+}v)| \label{eq:highest_IBP} \\
&\lesssim \|\partial_{t}^{j}\nabla \eta_{1,\mathrm{int}}\|_{L^{2}(\Sigma)}\|v_{0,\mathrm{int}}\|_{W^{1,\infty}}\|\partial_{t}^{j}v\|_{H^{1}}+\|\partial_{t}^{j}\nabla \bar{\eta}_{1,\mathrm{int}}\|_{L^{2}}\|v_{0,\mathrm{int}}\|_{W^{1,\infty}}\|\partial_{t}^{j}v\|_{H^{1}} \notag \\
&\lesssim \|\partial_{t}^{j}\eta_{1,\mathrm{int}}\|_{H^{1}(\Sigma)}\|v_{0,\mathrm{int}}\|_{W^{3,p}}\|\partial_{t}^{j}v\|_{H^{1}}+\|\partial_{t}^{j}{\eta}_{1,\mathrm{int}}\|_{H^{1}(\Sigma)}\|v_{0,\mathrm{int}}\|_{W^{3,p}}\|\partial_{t}^{j}v\|_{H^{1}} \lesssim  \sqrt{\mathcal{E}}\mathcal{D}. \notag
\end{align}
The integration by parts in the last display does not use the choice of the test function. For every $w\in H^{1}(\Omega)$ with $w\cdot n_{0}=0$ on $\Sigma_{s}$ it gives
\begin{equation*}
\Big|\int_{\Omega}J\chi_{\mathrm{int}+}\mu\func{div}_{\mathcal{A}_\eta}\big(\mathbb{D}_{\partial_t^j\mathcal{A}_\eta}v\big)\cdot w\Big|\lesssim\sqrt{\mathcal{E}\mathcal{D}}\,\|w\|_{H^{1}},
\end{equation*}
which is an $(H^{1})^{*}$ bound for this part of $F^{1,j}_{0}$. This is the form used in Proposition~\ref{prop:B_d}, in Proposition~\ref{thm:whole_space_threehalf} against $w=M\nabla\psi$, in \eqref{eq:F01_dual}, and in Step~9 of the proof of Theorem~\ref{thm:contact_line_energy}.

Combining all of the estimates above, we obtain
\begin{align}
\int_{\Omega}JF_{0}^{1,j}\cdot v_{j}\lesssim \|v_{j}\|_{H^{1}}\sqrt{\mathcal{E}\mathcal{D}}.
\end{align}

\textbf{Term including $F^{3,j}$}

It suffices to estimate the interior part of $F_{0}^{3,j}$ denoted by $F_{\mathrm{int}}^{3,j}$. Since $\Omega$ is three-dimensional, the trace theorem and the Sobolev embedding on $\Sigma$ give  $H^1(\Omega)\longrightarrow H^{1/2}(\Sigma) \hookrightarrow L^4(\Sigma)$. Consequently,
\begin{align}
\left| \int_\Sigma F_{\mathrm{int}+}^{3,j}\cdot v_j \right| &\lesssim \|F_{\mathrm{int}+}^{3,j}\|_{L^{4/3}(\Sigma)} \|v_j\|_{H^1(\Omega)}. \label{eq:Fint3-trace}
\end{align}

Choose a time-independent smooth cutoff $\widetilde\chi$ by assuming that $\tilde{\chi}=1$ in a neighborhood of $\func{supp}(\chi_{\mathrm{int}+})$ and $\operatorname{supp}\widetilde{\chi} \Subset \operatorname{int}\{\chi_{0-,\mathrm{int}} = 1\}$. We use this cut off function in the estimate below. The localized forcing is $F_{\mathrm{int}+}^{3,j} =\chi_{\mathrm{int}+}F_0^{3,j}$, so its complete expansion is
\begin{align}
&F_{\mathrm{int}+}^{3,j} ={} -\sigma\chi_{\mathrm{int}+}\sum_{r=1}^{j}\binom jr \Bigl(\operatorname{div}_\Sigma [\partial_t^{j-r}\mathcal R(\nabla\eta)]\Bigr)\partial_t^r\mathcal N +\mu\chi_{\mathrm{int}+}\sum_{r=1}^{j}\binom jr \bigl(\mathbb D_{\partial_t^r\mathcal A_\eta}v_{j-r}\bigr)\mathcal N \notag \\
&+\chi_{\mathrm{int}+}\sum_{r=1}^{j}\binom jr \bigl(g\eta_{j-r}-\sigma\operatorname{div}_\Sigma \mathcal L_{\Sigma,0}\eta_{j-r}\bigr)\partial_t^r\mathcal N -\chi_{\mathrm{int}+}\sum_{r=1}^{j}\binom jr S_{\mathcal A_\eta}(p_{j-r},v_{j-r})\partial_t^r\mathcal N \notag \\
&+\mu\chi_{\mathrm{int}+}\sum_{r=1}^{j}\sum_{s=1}^{j-r} \binom jr\binom{j-r}s \bigl(\mathbb D_{\partial_t^s\mathcal A_\eta}v_{j-r-s}\bigr) \partial_t^r\mathcal N.
\end{align}

Fix $0<\alpha<1/2$. The exponents satisfy  $\frac{1+2\alpha}{4} +\frac{1-\alpha}{2} +\frac14=1$. Separating the terms with $r=j$, applying H\"older's inequality, using the uniform $L^\infty$ bounds for $\mathcal{A}_\eta$, $\mathcal{N}$, and $|\mathcal{N}_0|^{-1}$, we obtain
\begin{align}
&\|F_{\mathrm{int}+}^{3,j}\|_{L^{4/3}(\Sigma)} \lesssim \sum_{r=1}^{j-1}\binom jr \left\| \chi_{0-,\mathrm{int},\Sigma}\partial_t^r\mathcal{N }\right\|_{L^{\frac4{1+2\alpha}}(\Sigma)} \left\| \chi_{0-,\mathrm{int},\Sigma} \func{div}_{\Sigma}\partial_t^{j-r}\bigl(\mathcal{R}(\nabla\eta)\bigr) \right\|_{L^{\frac2{1-\alpha}}(\Sigma)} \notag \\
&+ \mu\sum_{r=1}^{j-1}\binom jr \left\| \tilde{\chi}_{\Sigma}\partial_t^r\mathcal{A}_\eta \right\|_{L^{\frac4{1+2\alpha}}(\Sigma)} \left\| \chi_{0-,\mathrm{int},\Sigma}\nabla v_{j-r} \right\|_{L^{\frac2{1-\alpha}}(\Sigma)} \notag \\
&+ \sum_{r=1}^{j-1}\binom jr \left\| \chi_{0-,\mathrm{int},\Sigma}\partial_t^r\mathcal{N }\right\|_{L^{\frac4{1+2\alpha}}(\Sigma)} \Bigg( g\left\| \chi_{0-,\mathrm{int},\Sigma}\eta_{j-r} \right\|_{L^{\frac2{1-\alpha}}(\Sigma)} + \sigma\left\| \chi_{0-,\mathrm{int},\Sigma}\func{div}_{\Sigma}\mathcal{L}_{\Sigma,0}\eta_{j-r} \right\|_{L^{\frac2{1-\alpha}}(\Sigma)} \notag \\
&+ \left\| \chi_{0-,\mathrm{int},\Sigma}p_{j-r} \right\|_{L^{\frac2{1-\alpha}}(\Sigma)} + \mu\left\| \chi_{0-,\mathrm{int},\Sigma}\nabla v_{j-r} \right\|_{L^{\frac2{1-\alpha}}(\Sigma)} \Bigg) \notag \\
&+ \mu\sum_{r=1}^{j-1}\sum_{s=1}^{j-r} \binom jr\binom{j-r}{s} \left\| \chi_{0-,\mathrm{int},\Sigma}\partial_t^r\mathcal{N }\right\|_{L^{\frac4{1+2\alpha}}(\Sigma)} \left\| \tilde{\chi}_{\Sigma}\partial_t^s\mathcal{A}_\eta\right\|_{L^{\infty}} \left\|\chi_{0-,\mathrm{int},\Sigma} \nabla v_{j-r-s} \right\|_{L^{\frac2{1-\alpha}}(\Sigma)} \notag \\
&+ \mu\left\| \tilde{\chi}_{\Sigma}\partial_t^j\mathcal{A}_\eta \right\|_{L^2(\Sigma)} \left\| \chi_{0-,\mathrm{int},\Sigma}\nabla v \right\|_{L^4(\Sigma)} \notag \\
&+ \left\| \chi_{0-,\mathrm{int},\Sigma}\partial_t^j\mathcal{N }\right\|_{L^2(\Sigma)} \Bigg( \left\| \chi_{0-,\mathrm{int},\Sigma}\mathcal{\mathop{\rm div}\nolimits}_{\Sigma}R(\nabla\eta) \right\|_{L^4(\Sigma)} + g\left\| \chi_{0-,\mathrm{int},\Sigma}\eta \right\|_{L^4(\Sigma)} + \sigma\left\| \chi_{0-,\mathrm{int},\Sigma}\func{div}_{\Sigma}\mathcal{L}_{\Sigma,0}\eta \right\|_{L^4(\Sigma)} \notag \\
&+ \left\| \chi_{0-,\mathrm{int},\Sigma}p \right\|_{L^4(\Sigma)} + \mu\left\| \chi_{0-,\mathrm{int},\Sigma}\nabla v \right\|_{L^4(\Sigma)} \Bigg). \label{eq:Fint3-products}
\end{align}
All sums with an empty index range are understood to be zero.

For the double sum, we apply the Sobolev embedding and trace theorems. The observation that $j-s\leq m$ implies  $\min\{j-r-s,r\}\leq \frac m2<m-1, \; m\geq3$ allows us to distinguish cases according to which of the factors $\nabla v_{j-r-s}$ and $\partial_t^r\mathcal{N}$ has fewer temporal derivatives. The factor with fewer temporal derivatives can be bounded by the energy functional, and the factor with more temporal derivatives is estimated by the dissipation functional. Therefore, when $1\leq s\leq j-1$, by \eqref{eq:A_equi}, we have
\begin{align}
&\left\| \tilde{\chi}_{\Sigma}\partial_t^s\mathcal{A}_\eta \right\|_{L^\infty(\Sigma)} \left\| \chi_{0-,\mathrm{int},\Sigma}\nabla v_{j-r-s} \right\|_{L^{\frac2{1-\alpha}}(\Sigma)}  \notag \\
&\quad\lesssim (\left\| \chi_{0,\mathrm{int},\Sigma} \partial_t^s\eta \right\|_{W^{3-\frac{1}{p},p}(\Sigma)}+\sqrt{\mathcal{E}\mathcal{D}}\mathbf{1}_{s=j-1}+{\mathcal{E}}\mathbf{1}_{0\leq s<j-1}) \left\| \chi_{0,\mathrm{int}}v_{j-r-s} \right\|_{W^{2,p}(\Omega)}  \notag \\
&\quad\lesssim \sqrt{ \mathscr E_{m-1,\mathrm{high},\mathrm{int}} \mathscr D_{m,\mathrm{high},\mathrm{int}}}+{\mathcal{E}}\sqrt{\mathcal{D}}.  \label{eq:Fint3-double-product}
\end{align}

We next estimate the boundary norms by the spatial norms appearing in the energy and dissipation. The surface embedding  $H^{\frac12-\alpha}(\Sigma) \hookrightarrow L^{\frac4{1+2\alpha}}(\Sigma)$ implies, when the non-unit normal depends affinely on $\eta$ and its first surface derivatives,
\begin{equation}
\left\|\chi_{0-,\mathrm{int}}\partial_t^r\mathcal{N }\right\|_{L^{\frac4{1+2\alpha}}(\Sigma)}  \lesssim \|\chi_{0,\mathrm{int}}\eta_r \|_{H^{\frac32-\alpha}(\Sigma)},\quad \left\|\chi_{0-,\mathrm{int}}\partial_t^j\mathcal{N }\right\|_{L^2(\Sigma)} \lesssim \|\chi_{0,\mathrm{int}}\eta_j\|_{H^1(\Sigma)}.
\end{equation}
For nonlinear dependence of $\mathcal{N}$ on these variables, the corresponding lower-order tame products must also be included. The derivatives of $\mathcal{A}_\eta$ are estimated by the geometric coefficient estimates, retaining all products generated by the chain rule.

Moreover, the trace theorem gives
\begin{align}
&\|\chi_{0-,\mathrm{int}}\nabla v_\ell\|_{L^{\frac2{1-\alpha}}(\Sigma)} \lesssim \|\chi_{0,\mathrm{int}} v_\ell\|_{W^{2,p}(\Omega)},\quad \|\chi_{0-,\mathrm{int}}p_\ell\|_{L^{\frac2{1-\alpha}}(\Sigma)}  \lesssim \|\chi_{0,\mathrm{int}} p_\ell\|_{W^{1,p}(\Omega)}, \notag \\
&\|\chi_{0-,\mathrm{int}}\eta_\ell\|_{L^{\frac2{1-\alpha}}(\Sigma)} + \|\chi_{0-,\mathrm{int}}\func{div}_{\Sigma}\mathcal{L}_{\Sigma,0}\eta_\ell \|_{L^{\frac2{1-\alpha}}(\Sigma)}  \lesssim \|\chi_{0,\mathrm{int}}\eta_\ell \|_{W^{3-\frac1p,p}(\Sigma)}. \notag \\
&\|\chi_{0-,\mathrm{int}}\eta_\ell\|_{L^{4}(\Sigma)} + \|\chi_{0-,\mathrm{int}}\func{div}_{\Sigma}\mathcal{L}_{\Sigma,0}\eta_\ell \|_{L^{4}(\Sigma)} \lesssim \|\chi_{0,\mathrm{int}}\eta_\ell \|_{W^{4-\frac1p,p}(\Sigma)}.
\end{align}
In the last line we take $\ell=0$ for the undifferentiated height terms in \eqref{eq:Fint3-products}. Substituting these estimates, we find
\begin{align}
&\|F_{\mathrm{int}+}^{3,j}\|_{L^{4/3}(\Sigma)}  \notag \\
&\quad\lesssim_{\varepsilon} \sum_{r=1}^{j-1}\binom jr \|\chi_{0,\mathrm{int},\Sigma}\eta_r \|_{H^{\frac{3}{2}-\alpha}(\Sigma)} \Bigg( \|\chi_{0,\mathrm{int},\Sigma}\eta_{j-r} \|_{W^{3-\frac1p,p}(\Sigma)} + \|\chi_{0,\mathrm{int}} p_{j-r}\|_{W^{1,p}(\Omega)} + \mu\|\chi_{0,\mathrm{int}} v_{j-r}\|_{W^{2,p}(\Omega)} \Bigg)  \notag \\
&\qquad+ \mu\sum_{r=1}^{j-1}\sum_{s=1}^{j-r} \binom jr\binom{j-r}{s} \|\chi_{0,\mathrm{int},\Sigma}\eta_r \|_{H^{\frac32-\alpha}(\Sigma)}\sqrt{ \mathscr E_{m-1,\mathrm{high},\mathrm{int}} \mathscr D_{m,\mathrm{high},\mathrm{int}}}  \notag \\
&\qquad+ \mu\left\|\chi_{0,\mathrm{int},\Sigma}\partial_t^j\eta \right\|_{H^1(\Sigma)} \|\chi_{0,\mathrm{int}} v\|_{W^{3,p}(\Omega)}  \notag \\
&\qquad+ \left\|\chi_{0,\mathrm{int},\Sigma}\partial_t^j\eta \right\|_{H^1(\Sigma)} \Bigg( \|\chi_{0,\mathrm{int},\Sigma}\eta\|_{W^{4-\frac1p,p}(\Sigma)} + \|\chi_{0,\mathrm{int}} p\|_{W^{2,p}(\Omega)} + \mu\|\chi_{0,\mathrm{int}} v\|_{W^{3,p}(\Omega)} \Bigg)+{\mathcal{E}}\sqrt{\mathcal{D}}  \notag \\
&\quad\lesssim C(\varepsilon)(\mathscr{E}^{\frac{1}{2}}_{m+1,\mathrm{energy},\mathrm{bet}}+\mathscr{E}^{\frac{1}{2}}_{\mathrm{int}})(\mathscr{D}^{\frac{1}{2}}_{\mathrm{int}}+C(\varepsilon)\mathscr{D}^{\frac{1}{2}}_{m+1,\mathrm{energy},\mathrm{bet}})+{\mathcal{E}}\sqrt{\mathcal{D}}. \label{eq:Fint3-spatial-products}
\end{align}
Combining this estimate with \eqref{eq:Fint3-trace} yields
\begin{align}
\left| \int_\Sigma F_{\mathrm{int}+}^{3,j}\cdot v_j \right| \lesssim C(\varepsilon)\sqrt{\mathcal{E}}\mathcal{D}\,.  \label{eq:Fint3-energy}
\end{align}

\textbf{Estimate for the interior contribution involving $F_0^{4,j}$.}

Since $F_0^{4,0}=0$, the terms containing $F_0^{4,j}$ vanish at $j=0$. For $1\leq j\leq m+1$, we now estimate the term
\begin{align}
g\int_\Sigma\eta_jF_0^{4,j} + \sigma\int_\Sigma a_\zeta\nabla_\Sigma\eta_j\cdot\nabla_\Sigma F_0^{4,j}.
\end{align}

For the contact line part, by Lemma \ref{lem:hard_nonlinearity} and Lemma \ref{thm:contact_line_forcing_estimate}, we have
\begin{align}
&g\int_\Sigma\eta_j\mathcal{R}_{\mathrm{cl}+}^{4,j} + \sigma\int_\Sigma a_\zeta\nabla_\Sigma\eta_j\cdot\nabla_\Sigma \mathcal{R}_{\mathrm{cl}+}^{4,j}+g\int_\Sigma\eta_jv\cdot \partial_{t}^{j}\mathcal{N} + \sigma\int_\Sigma a_\zeta\nabla_\Sigma\eta_j\cdot\nabla_\Sigma v\cdot \partial_{t}^{j}\mathcal{N}  \notag \\
\lesssim& \sqrt{\mathcal{E}}\mathcal{D}+\|\eta_{j}\|_{L^{2}H^{\frac{3}{2}-\alpha}}\|\nabla_{\Sigma}\mathcal{R}_{\mathrm{cl}+}^{4,j}\|_{L^{2}L^{\frac{1}{1-\alpha}}}\lesssim\sqrt{\mathcal{E}}\mathcal{D}.
\end{align}

For the interior part, recall that
\begin{equation*}
\mathcal{R}_{\mathrm{int}+}^{4,j} = \chi_{\mathrm{int}+} \sum_{r=1}^{j-1}\binom jr v_{j-r}\cdot\partial_t^r\mathcal{N}.
\end{equation*}
We assume that $\chi_{\mathrm{int}+}|_\Sigma$ is supported away from $\partial\Sigma$ and that $\chi_{0,\mathrm{int},\Sigma}=1$ on a neighborhood of its support. All cutoffs are independent of time. For $1\leq r\leq j-1$, the preceding product estimates give
\begin{align}
&\left| g\int_\Sigma\eta_j\chi_{\mathrm{int}+} \sum_{r=1}^{j-1}\binom jr v_{j-r}\cdot\partial_t^r\mathcal{N }\right.\quad \left. + \sigma\int_\Sigma a_\zeta\nabla_\Sigma\eta_j\cdot \nabla_\Sigma\left( \chi_{\mathrm{int}+} \sum_{r=1}^{j-1}\binom jr v_{j-r}\cdot\partial_t^r\mathcal{N }\right) \right| \notag \\
&\lesssim \|\eta_j\|_{H^1(\Sigma)} \sum_{r=1}^{j-1}\binom jr \Bigg[ \left\| \chi_{0-,\mathrm{int},\Sigma}v_{j-r} \right\|_{L^\infty(\Sigma)} \left\| \chi_{0-,\mathrm{int},\Sigma}\partial_t^r\mathcal{N }\right\|_{H^1(\Sigma)} \notag \\
&+ \left\| \chi_{0-,\mathrm{int},\Sigma}\nabla_\Sigma v_{j-r} \right\|_{L^{\frac{2p}{3-p}}(\Sigma)} \left\| \chi_{0-,\mathrm{int},\Sigma}\partial_t^r\mathcal{N }\right\|_{L^{\infty}(\Sigma)} \Bigg] \notag \\
&\lesssim \|\eta_j\|_{H^1(\Sigma)} \sum_{r=1}^{j-1}\binom jr \Bigg[ \left\| \chi_{0,\mathrm{int}}v_{j-r} \right\|_{W^{2,p}} \left\| \chi_{0,\mathrm{int},\Sigma}\partial_t^{r}\eta \right\|_{W^{3-\frac{1}{p},p}(\Sigma)} \Bigg] \lesssim \sqrt{\mathcal{E}}\mathcal{D}. \label{eq:Fint4-lower}
\end{align}
This finishes the estimate for this term.

\textbf{Estimate for the interior contribution involving $F_0^{5,j}$.}

Let  $F_0^{5,j} = \partial_t^j\bigl(\mathcal{R}(\nabla_\Sigma\eta)\bigr), \; F_{\mathrm{int}+}^{5,j} = \chi_{\mathrm{int}+}F_0^{5,j}$, where $\chi_{\mathrm{int}+}$ is independent of time. We assume that $\nabla_\Sigma$ and the reference surface measure are also independent of time.

At $j=0$, we keep the full curvature remainder. In the horizontal coordinates $x\in D$, define
\begin{equation}
\mathscr R_0^5:=\sigma\int_D\Bigg( \sqrt{1+|\nabla\zeta+\nabla\eta|^2} -\sqrt{1+|\nabla\zeta|^2} -\frac{\nabla\zeta\cdot\nabla\eta}{\sqrt{1+|\nabla\zeta|^2}} -\frac12\mathcal L_{\Sigma,0}\nabla\eta\cdot\nabla\eta \Bigg)\,dx.
\end{equation}
The derivative of the integrand with respect to $\nabla\eta$ is $\mathcal R(\nabla\eta)$. Thus the capillary term in \eqref{eq:whole-domain-energy-j}, written in these coordinates, satisfies
\begin{equation*}
-\sigma\int_D F_0^{5,0}\cdot\nabla\partial_t\eta\,dx =-\frac{d}{dt}\mathscr R_0^5.
\end{equation*}
Taylor's formula at $\nabla\eta=0$ and the fixed equilibrium coefficients give $|\mathscr R_0^5| \leq C\|\nabla\eta\|_{L^\infty(D)}\|\nabla\eta\|_{L^2(D)}^2 \leq C(\varepsilon)\mathcal E^{3/2}$. There are no mixed terms at this level, since $F_0^{4,0}=0$. For the remainder of the curvature calculation let $1\leq j\leq m+1$. The chain rule gives $F_0^{5,j} = \mathcal{R}^{\prime }(\nabla_\Sigma\eta)\nabla_\Sigma\eta_j + \mathcal{R}_{\mathrm{low}}^{\,j}$, where
\begin{align}
\mathcal{R}_{\mathrm{low}}^{\,j} :={}& \sum_{\ell=2}^{j}\frac1{\ell!} \sum _{\substack{ i_1+\cdots+i_\ell=j \\ i_k\geq1}} \frac{j!}{i_1!\cdots i_\ell!} \mathcal{R}^{(\ell)}(\nabla_\Sigma\eta) \left[ \nabla_\Sigma\eta_{i_1},\ldots, \nabla_\Sigma\eta_{i_\ell} \right]. \label{eq:Fint5-lower-definition}
\end{align}
For $j=1$, this remainder is zero.

For the highest-order term, by Lemma \ref{lem:hard_nonlinearity}, we have
\begin{align}
&-\sigma\int_\Sigma \left( \mathcal{R}^{\prime }(\nabla_\Sigma\eta)\nabla_\Sigma\eta_j \right)\cdot\nabla_\Sigma\partial_t\eta_j = -\frac{d}{dt} \left[ \frac{\sigma}{2}\int_\Sigma \left( \mathcal{R}^{\prime }(\nabla_\Sigma\eta)\nabla_\Sigma\eta_j \right)\cdot\nabla_\Sigma\eta_j \right] \notag \\
&+ \frac{\sigma}{2}\int_\Sigma \left( \partial_t\bigl[\mathcal{R}^{\prime }(\nabla_\Sigma\eta)\bigr] \nabla_\Sigma\eta_j \right)\cdot\nabla_\Sigma\eta_j,
\end{align}
where
\begin{align}
\frac{\sigma}{2}\int_\Sigma \left( \mathcal{R}^{\prime }(\nabla_\Sigma\eta)\nabla_\Sigma\eta_j \right)\cdot\nabla_\Sigma\eta_j\lesssim \sqrt{\mathcal{E}}^{3},
\end{align}
and
\begin{align}
\frac{\sigma}{2}\int_\Sigma \left( \partial_t\bigl[\mathcal{R}^{\prime }(\nabla_\Sigma\eta)\bigr] \nabla_\Sigma\eta_j \right)\cdot\nabla_\Sigma\eta_j\lesssim\sqrt{\mathcal{E}}\mathcal{D}.
\end{align}

For the lower-order remainder, we again separate the forcing term into the contact line part and the interior part. For the contact line part, temporal integration by parts gives
\begin{equation}
-\sigma\int_\Sigma \mathcal{R}_{\mathrm{cl}+}^{5,j} \cdot\nabla_\Sigma\partial_t\eta_j = -\frac{d}{dt} \left[ \sigma\int_\Sigma R_{\mathrm{cl}+}^{5,j} \cdot\nabla_\Sigma\eta_j \right] + \sigma\int_\Sigma \partial_{t}R_{\mathrm{cl}+}^{5,j} \cdot\nabla_\Sigma\eta_j,
\end{equation}
where
\begin{align}
\left[ \sigma\int_\Sigma \mathcal{R}_{\mathrm{cl}+}^{5,j} \cdot\nabla_\Sigma\eta_j \right]\lesssim \|\mathcal{R}_{\mathrm{cl}+}^{5,j}\|_{L^{2}L^{2}} \|\eta_{j}\|_{H^{1}}\lesssim \sqrt{\mathcal{E}}^{3},
\end{align}
and
\begin{align}
\sigma\int_\Sigma \partial_{t}\mathcal{R}_{\mathrm{cl}+}^{5,j} \cdot\nabla_\Sigma\eta_j\lesssim\|\nabla_{\Sigma}\eta_{j}\|_{L^{2}L^{\frac{1}{\alpha}}}\|\partial_{t}\mathcal{R}_{\mathrm{cl}+}^{5,j}\|_{L^{2}L^{\frac{1}{1-\alpha}}}\lesssim\sqrt{\mathcal{E}}\mathcal{D}.
\end{align}

For the interior part, we have
\begin{equation}
-\sigma\int_\Sigma \chi_{\mathrm{int}+} \mathcal{R}_{\mathrm{low}}^{\,j} \cdot\nabla_\Sigma\partial_t\eta_j = -\frac{d}{dt} \left[ \sigma\int_\Sigma \chi_{\mathrm{int}+} \mathcal{R}_{\mathrm{low}}^{\,j} \cdot\nabla_\Sigma\eta_j \right] + \sigma\int_\Sigma \chi_{\mathrm{int}+} \partial_t\mathcal{R}_{\mathrm{low}}^{\,j} \cdot\nabla_\Sigma\eta_j.
\end{equation}

Define
\begin{equation}
\mathscr R_j^5 :={}  \frac{\sigma}{2}\int_\Sigma \left( \mathcal{R}^{\prime }(\nabla_\Sigma\eta)\nabla_\Sigma\eta_j \right)\cdot\nabla_\Sigma\eta_j + \sigma\int_\Sigma \mathcal{R}_{\mathrm{low}}^{\,j} \cdot\nabla_\Sigma\eta_j. \label{eq:Fint5-correction}
\end{equation}
Combining the preceding identities yields
\begin{equation}
-\sigma\int_\Sigma F_{0}^{5,j} \cdot\nabla_\Sigma\partial_t\eta_j = -\frac{d}{dt}\mathcal{R}_j^5 + \frac{\sigma}{2}\int_\Sigma \left( \partial_t\bigl[\mathcal{R}^{\prime }(\nabla_\Sigma\eta)\bigr] \nabla_\Sigma\eta_j \right)\cdot\nabla_\Sigma\eta_j + \sigma\int_\Sigma \chi_{\mathrm{int}+} \partial_t\mathcal{R}_{\mathrm{low}}^{\,j} \cdot\nabla_\Sigma\eta_j+\mathcal{O}_{\mathbf{contact}},
\end{equation}
where $\mathcal{O}_{\mathrm{contact}}$ contains the contact line part and $\mathcal{O}_{\mathbf{contact}}\lesssim \sqrt{\mathcal{E}}\mathcal{D}$.

We now estimate the only remainder integral. Assume that $\chi_{0,\mathrm{int},\Sigma}=1$ on a neighborhood of $\func{supp}(\chi_{\mathrm{int}+}|_\Sigma)$. For $0<\alpha<1/2$, by H\"older's inequality and  $H^{\frac12-\alpha}(\Sigma) \hookrightarrow L^{\frac4{1+2\alpha}}(\Sigma)$, we have
\begin{align}
&\left| \sigma\int_\Sigma \chi_{\mathrm{int}+} \partial_t\mathcal{R}_{\mathrm{low}}^{\,j} \cdot\nabla_\Sigma\eta_j \right| \lesssim \left\| \chi_{0,\mathrm{int},\Sigma} \partial_t\mathcal{R}_{\mathrm{low}}^{\,j} \right\|_{L^{\frac4{3-2\alpha}}(\Sigma)} \left\| \chi_{0-,\mathrm{int},\Sigma}\eta_j \right\|_{H^{\frac32-\alpha}(\Sigma)} \notag \\
&\lesssim \|\chi_{1,\mathrm{int},\Sigma}\partial_{t}\eta\|_{W^{1,\infty}}\left\| \chi_{0,\mathrm{int},\Sigma}\eta_j \right\|_{H^{\frac32-\alpha}(\Sigma)}^{2}+\left\| \chi_{0,\mathrm{int},\Sigma}\eta_j \right\|_{H^{\frac32-\alpha}(\Sigma)}\sum_{\ell=2}^{j+1}\sum _{\substack{ i_1+\cdots+i_\ell=j+1 \\ j-1\geq i_k\geq1}}\prod_{\ell}\|\chi_{0,\mathrm{int},\Sigma}\partial_{t}^{i_{\ell}}\eta\|_{W^{1,\infty}} \notag \\
&\lesssim\sqrt{\mathcal{E}}(\mathcal{D}).
\end{align}
This estimate holds since every factor $\eta_{i_k}$ in $\mathcal{R}_{\mathrm{low}}^{\,j}$ has $i_k\leq j-1$. Thus, after one temporal derivative, the largest temporal order appearing in $\partial_t\mathcal{R}_{\mathrm{low}}^{\,j}$ is $j$. In particular, this remainder contains no $\eta_{j+1}$. Consequently,
\begin{equation}
\left| -\sigma\int_\Sigma F_{0}^{5,j} \cdot\nabla_\Sigma\partial_t\eta_j + \frac{d}{dt}\mathcal{R}_j^5 \right| \quad\lesssim\sqrt{\mathcal{E}}(\mathcal{D}) . \label{eq:Fint5-final}
\end{equation}

Finally, the energy correction satisfies
\begin{align}
|\mathscr R_j^5| \lesssim{}&\sqrt{\mathcal{E}}^{3} + \left\| \chi_{0,\mathrm{int},\Sigma} \mathcal{R}_{\mathrm{low}}^{\,j} \right\|_{L^2(\Sigma)} \left\| \chi_{0,\mathrm{int},\Sigma}\eta_j \right\|_{H^1(\Sigma)}. \label{eq:Fint5-correction-bound}
\end{align}
For $j\geq2$, recall that
\begin{align}
\mathcal{R}_{\mathrm{low}}^{\,j} = \sum_{\ell=2}^{j}\frac1{\ell!} \sum _{\substack{ i_1+\cdots+i_\ell=j \\ i_k\geq1}} \frac{j!}{i_1!\cdots i_\ell!} \mathcal{R}^{(\ell)}(\nabla_\Sigma\eta) \left[ \nabla_\Sigma\eta_{i_1},\ldots, \nabla_\Sigma\eta_{i_\ell} \right].
\end{align}
In each summand, choose $k_*$ such that  $i_{k_*}=\max_{1\leq k\leq\ell}i_k$. Since $\ell\geq2$ and each $i_k\geq1$, we have $i_{k_*}\leq j-1$. We place this factor in $L^2(\Sigma)$ and all remaining factors in $L^\infty(\Sigma)$. Choose $\chi_{1,\mathrm{int},\Sigma}=1$ on a neighborhood of $\func{supp}\chi_{0,\mathrm{int},\Sigma}$. Assuming that the derivatives of $\mathcal{R}$ through order $j$ are uniformly bounded on the range of $\nabla_\Sigma\eta$, H\"older's inequality gives
\begin{align}
&\left\| \chi_{0,\mathrm{int},\Sigma} \mathcal{R}_{\mathrm{low}}^{\,j} \right\|_{L^2(\Sigma)} \lesssim \sum_{\ell=2}^{j} \sum_{\substack{ i_1+\cdots+i_\ell=j \\ i_k\geq1}} \left\| \chi_{1,\mathrm{int},\Sigma} \nabla_\Sigma\eta_{i_{k_*}} \right\|_{L^2(\Sigma)} \prod_{\substack{ k=1 \\ k\neq k_*}}^\ell \left\| \chi_{1,\mathrm{int},\Sigma} \nabla_\Sigma\eta_{i_k} \right\|_{L^\infty(\Sigma)} \notag \\
&\lesssim \sum_{\ell=2}^{j} \sum_{\substack{ i_1+\cdots+i_\ell=j \\ i_k\geq1}} \left\| \chi_{1,\mathrm{int},\Sigma}\eta_{i_{k_*}} \right\|_{H^1(\Sigma)} \prod_{\substack{ k=1 \\ k\neq k_*}}^\ell \left\| \chi_{1,\mathrm{int},\Sigma}\eta_{i_k} \right\|_{W^{1,\infty}(\Sigma)}. \label{eq:Fint5-lower-L2-Linfty}
\end{align}
The implicit constant may depend on $j$. In particular, for the terms containing $\eta_{j-1}$, the derivative distribution is $(j-1,1)$, and their contribution is bounded by
\begin{equation*}
C_j \left\| \chi_{1,\mathrm{int},\Sigma}\eta_{j-1} \right\|_{H^1(\Sigma)} \left\| \chi_{1,\mathrm{int},\Sigma}\eta_1 \right\|_{W^{1,\infty}(\Sigma)}.
\end{equation*}
For the remaining terms, the largest temporal order is at most $j-2$.

Consequently, the second term in \eqref{eq:Fint5-correction-bound} satisfies
\begin{align}
&\left\| \chi_{0,\mathrm{int},\Sigma} \mathcal{R}_{\mathrm{low}}^{\,j} \right\|_{L^2(\Sigma)} \left\| \chi_{0,\mathrm{int},\Sigma}\eta_j \right\|_{H^1(\Sigma)}  \notag \\
&\quad\lesssim \left\| \chi_{0,\mathrm{int},\Sigma}\eta_j \right\|_{H^1(\Sigma)} \sum_{\ell=2}^{j} \sum_{\substack{ i_1+\cdots+i_\ell=j \\ i_k\geq1}} \left\| \chi_{1,\mathrm{int},\Sigma}\eta_{i_{k_*}} \right\|_{H^1(\Sigma)} \prod_{\substack{ k=1 \\ k\neq k_*}}^\ell \left\| \chi_{1,\mathrm{int},\Sigma}\eta_{i_k} \right\|_{W^{1,\infty}(\Sigma)}. \label{eq:Fint5-correction-lower-products}
\end{align}
For $p>3/2$, the remaining factors can be estimated using  $W^{3-\frac1p,p}(\Sigma)\hookrightarrow W^{1,\infty}(\Sigma)$. Therefore
\begin{align}
&\left\| \chi_{0,\mathrm{int},\Sigma} \mathcal{R}_{\mathrm{low}}^{\,j} \right\|_{L^2(\Sigma)} \left\| \chi_{0,\mathrm{int},\Sigma}\eta_j \right\|_{H^1(\Sigma)} \lesssim \mathcal{E}^{3/2},
\end{align}
and hence $\|\mathscr{R}_{j}^{5}\|_{L^{2}(\Sigma)}\lesssim \mathcal{E}^{\frac{3}{2}}$.

\textbf{Estimate for the interior contribution involving $F_0^{5,j}\cdot\nabla_\Sigma F_0^{4,j}$.}

We consider
\begin{equation*}
\sigma\int_\Sigma F_0^{5,j} \cdot\nabla_\Sigma F_0^{4,j},
\end{equation*}
where $\chi_{\mathrm{int}+}|_\Sigma$ is supported away from $\partial\Sigma$. Note that $F_{0}^{5,j}$ can be rewritten as follows for some symmetric tensor $\mathcal{R}^{\prime }(\nabla\eta)$  $F_0^{5,j} = \mathcal{R}^{\prime }(\nabla_\Sigma\eta)\nabla \partial_{t}^{j}\eta + \mathcal{R}_{\mathrm{low}}^{5,j}$. Also note that $\partial_t^j\mathcal{N}=(-\nabla_\Sigma\eta_j,0)$. We have  $F_0^{4,j} = -v_\parallel\cdot\nabla_\Sigma\eta_j +\mathcal{R}_{\mathrm{low}}^{4,j}$. Consequently,
\begin{align*}
&\sigma\int_\Sigma F_0^{5,j} \cdot\nabla_\Sigma F_0^{4,j} =\sigma\int_{\Sigma}(\mathcal R'(\nabla_\Sigma\eta)\nabla_\Sigma\eta_j + \mathcal R_{\mathrm{cl}+}^{5,j}+\mathcal R_{\mathrm{int}+}^{5,j}) \cdot \nabla_{\Sigma}(-v_\parallel\cdot\nabla_\Sigma\eta_j +\mathcal{R}_{\mathrm{cl}+}^{4,j}+\mathcal{R}_{\mathrm{int}+}^{4,j}) \\
&=-\sigma \int_{\Sigma}\mathcal R'(\nabla_\Sigma\eta)\cdot \nabla_\Sigma\eta_j\nabla_{\Sigma}(v_{||}\cdot \nabla_{\Sigma}\eta_{j}) +\sigma \int_{\Sigma}\mathcal R'(\nabla_\Sigma\eta)\nabla_\Sigma\eta_j\cdot \nabla_{\Sigma}\mathcal{R}_{\mathrm{cl}+}^{4,j} +\sigma \int_{\Sigma}\mathcal R'(\nabla_\Sigma\eta)\nabla_\Sigma\eta_j\cdot \nabla_{\Sigma}\mathcal{R}_{\mathrm{int}+}^{4,j} \\
&-\sigma \int_{\Sigma}\mathcal{R}_{\mathrm{cl}+}^{5,j}\cdot \nabla_{\Sigma}(v_{||}\cdot \nabla_{\Sigma}\eta_{j}) +\sigma \int_{\Sigma}\mathcal{R}_{\mathrm{cl}+}^{5,j}\cdot \nabla_{\Sigma}\mathcal{R}_{\mathrm{cl}+}^{4,j} +\sigma \int_{\Sigma}\mathcal{R}_{\mathrm{cl}+}^{5,j}\cdot\nabla_{\Sigma}\mathcal{R}_{\mathrm{int}+}^{4,j} -\sigma \int_{\Sigma}\mathcal{R}_{\mathrm{int}+}^{5,j}\cdot \nabla_{\Sigma}(v_{||}\cdot \nabla_{\Sigma}\eta_{j}) \\
&+\sigma \int_{\Sigma}\mathcal{R}_{\mathrm{int}+}^{5,j}\cdot \nabla_{\Sigma}\mathcal{R}_{\mathrm{cl}+}^{4,j} +\sigma \int_{\Sigma}\mathcal{R}_{\mathrm{int}+}^{5,j}\cdot\nabla_{\Sigma}\mathcal{R}_{\mathrm{int}+}^{4,j} =\sum_{i=1}^{9}I_{i}.
\end{align*}

For the term $I_{1}$, using Lemma \ref{lem:hard_nonlinearity}, we have  $(I_{1})\lesssim\sqrt{\mathcal{E}}\mathcal{D}$. For the term $I_{2}$, by Lemma \ref{thm:contact_line_forcing_estimate}, we have
\begin{align}
I_{2}\lesssim \|(\chi_{\mathrm{cl}}+\chi_{\mathrm{bet}})\eta\|_{L^{\infty}W^{3-\frac{1}{q_{\delta}},q_{\delta}}}\|\nabla_{\Sigma}((\chi_{\mathrm{cl}}+\chi_{\mathrm{bet}})\eta_{j})\|_{L^{2}L^{\frac{1}{\alpha}}}\|\nabla_{\Sigma}\mathcal{R}^{4,j}_{\mathrm{cl}+}\|_{L^{2}L^{\frac{1}{1-\alpha}}}\lesssim \sqrt{\mathcal{E}}\mathcal{D}.
\end{align}
For the term $I_{4}$, using Lemma~\ref{thm:contact_line_forcing_estimate}, we have
\begin{align}
I_{4}=-\sigma\int_{\Sigma}\mathcal{R}_{\mathrm{cl}+}^{5,j}\cdot ((\nabla_{\Sigma}((\chi_{\mathrm{cl},\Sigma}+\chi_{\mathrm{bet},\Sigma})v_{||}))^{T}\nabla_{\Sigma}\eta_{j})-\sigma\int_{\Sigma}(\mathcal{R}_{\mathrm{cl}+}^{5,j})\cdot (\nabla_{\Sigma}^{2}((\chi_{\mathrm{cl},\Sigma}+\chi_{\mathrm{bet},\Sigma})\eta_{j})((\chi_{\mathrm{cl},\Sigma}+\chi_{\mathrm{bet},\Sigma})v_{||})).
\end{align}
Using integration by parts and the fact that $v_{||}\cdot n_{0}=0$ and Lemma \ref{thm:contact_line_forcing_estimate}, we have
\begin{align}
&I_{4}\lesssim  \|\mathcal{R}_{\mathrm{cl}+}^{5,j}\|_{L^{2}L^{\frac{1}{\alpha}}}\|(\chi_{\mathrm{cl},\Sigma}+\chi_{\mathrm{bet},\Sigma})\nabla_{\Sigma}v\|_{L^{\infty}L^{\frac{1}{1-2\alpha}}(\Sigma)}\|(\chi_{\mathrm{cl},\Sigma}+\chi_{\mathrm{bet},\Sigma})\eta_{j}\|_{L^{2}H^{\frac{3}{2}-\alpha}} \notag \\
&+\|\nabla_{\Sigma}\mathcal{R}_{\mathrm{cl}+}^{5,j}\|_{L^{2}L^{\frac{1}{1-\alpha}}}\|(\chi_{\mathrm{cl},\Sigma}+\chi_{\mathrm{bet},\Sigma})v\|_{L^{\infty}(\Sigma)}\|(\chi_{\mathrm{cl},\Sigma}+\chi_{\mathrm{bet},\Sigma})\eta_{j}\|_{L^{2}H^{\frac{3}{2}-\alpha}} \notag \\
&\lesssim  \|\mathcal{R}_{\mathrm{cl}+}^{5,j}\|_{L^{2}L^{\frac{1}{\alpha}}}\|(\chi_{\mathrm{cl}}+\chi_{\mathrm{bet}})v\|_{L^{\infty}W^{2,q_{\delta}}(\Sigma)}\|(\chi_{\mathrm{cl},\Sigma}+\chi_{\mathrm{bet},\Sigma})\eta_{j}\|_{L^{2}H^{\frac{3}{2}-\alpha}} \notag \\
&+\|\nabla_{\Sigma}\mathcal{R}_{\mathrm{cl}+}^{5,j}\|_{L^{2}L^{\frac{1}{1-\alpha}}}\|(\chi_{\mathrm{cl}}+\chi_{\mathrm{bet}})v\|_{L^{\infty}W^{2,q_{\delta}}}\|(\chi_{\mathrm{cl},\Sigma}+\chi_{\mathrm{bet},\Sigma})\eta_{j}\|_{L^{2}H^{\frac{3}{2}-\alpha}} \lesssim   \sqrt{\mathcal{E}}\mathcal{D}.
\end{align}
For the term $I_{5}$, by Lemma \ref{thm:contact_line_forcing_estimate}, we have $I_{5}\lesssim \sqrt{\mathcal{E}}\mathcal{D}$.

We now estimate the terms $I_{j}$ that contain the interior part of the forcing. Note that there are some mixed contact-line and interior terms, such as $\int_{\Sigma}\mathcal{R}_{\mathrm{cl}+}^{5,j}\cdot \nabla_{\Sigma}\mathcal{R}_{\mathrm{int}+}^{4,j}$. We use the fact that $\func{supp}\{\chi_{\mathrm{cl}+}\chi_{\mathrm{int}+}\}\subseteq\{x\in \Omega| \chi_{0,\mathrm{int}}(x)=1\}$, which means that we can use the interior estimate to bound $\mathcal{R}_{\mathrm{cl}+}^{5,j}$. Therefore, we use the estimate for $I_{9}$ to bound all of these mixed terms. It then suffices to estimate the following terms. Recall that $\mathcal R_{\mathrm{int}+}^{a,j} =\chi_{\mathrm{int}+}\mathcal R_{\mathrm{low}}^{a,j}$ for $a\in\{4,5\}$. The exact interior sum is
\begin{equation}
I_3+I_7+I_9 ={} -\sigma\int_\Sigma\mathcal R_{\mathrm{int}+}^{5,j} \cdot\nabla_\Sigma(v_\parallel\cdot\nabla_\Sigma\eta_j) +\sigma\int_\Sigma \left(\mathcal R'(\nabla_\Sigma\eta)\nabla_\Sigma\eta_j +\mathcal R_{\mathrm{int}+}^{5,j}\right) \cdot\nabla_\Sigma\mathcal R_{\mathrm{int}+}^{4,j}.
\end{equation}
For the first integral, we transfer the entire surface derivative onto the localized lower-order remainder:
\begin{equation}
I_7 =-\sigma\int_\Sigma\mathcal R_{\mathrm{int}+}^{5,j} \cdot\nabla_\Sigma(v_\parallel\cdot\nabla_\Sigma\eta_j) =\sigma\int_\Sigma\func{div}_\Sigma\mathcal R_{\mathrm{int}+}^{5,j} (v_\parallel\cdot\nabla_\Sigma\eta_j). \label{eq:Fint5-Fint4-low-ibp}
\end{equation}
There is no boundary term because $\chi_{\mathrm{int}+}|_\Sigma$ vanishes near $\partial\Sigma$. In the other integral we retain the full product rule $\nabla_\Sigma\mathcal R_{\mathrm{int}+}^{4,j} =\chi_{\mathrm{int}+}\nabla_\Sigma\mathcal R_{\mathrm{low}}^{4,j} +\mathcal R_{\mathrm{low}}^{4,j}\nabla_\Sigma\chi_{\mathrm{int}+}$. In particular, the second term contributes
\begin{equation*}
\mathcal C_j:=\sigma\int_\Sigma \left(\mathcal R'(\nabla_\Sigma\eta)\nabla_\Sigma\eta_j +\mathcal R_{\mathrm{int}+}^{5,j}\right) \cdot\nabla_\Sigma\chi_{\mathrm{int}+}\, \mathcal R_{\mathrm{low}}^{4,j}.
\end{equation*}
On the support of $\nabla_\Sigma\chi_{\mathrm{int}+}$, the interior product bounds already used in \eqref{eq:Fint4-lower} and \eqref{eq:Fint5-lower-L2-Linfty} give
\begin{equation*}
\|\chi_{0,\mathrm{int},\Sigma}\mathcal R_{\mathrm{low}}^{4,j}\|_{L^2} \leq C(\varepsilon)\sqrt{\mathcal E\mathcal D},\quad \|\mathcal R_{\mathrm{int}+}^{5,j}\|_{L^2} +\|\chi_{0,\mathrm{int},\Sigma}\mathcal R'(\nabla_\Sigma\eta)\|_{L^\infty} \|\chi_{0,\mathrm{int},\Sigma}\eta_j\|_{H^1} \leq C(\varepsilon)\mathcal E.
\end{equation*}
Thus $|\mathcal C_j|\leq C(\varepsilon)\mathcal E^{3/2}\mathcal D^{1/2} \leq C(\varepsilon)\sqrt{\mathcal E}\mathcal D$, where we used Lemma~\ref{lem:E_le_D} and $\mathcal E\leq1$.

Assume that $\chi_{0,\mathrm{int},\Sigma}=1$ on a neighborhood of $\func{supp}(\chi_{\mathrm{int}+}|_\Sigma)$. By H\"older's inequality and $H^{1/2-\alpha}(\Sigma)\hookrightarrow L^{4/(1+2\alpha)}(\Sigma)$, \begin{align} &|I_3+I_7+I_9| \leq{}  C\|\func{div}_\Sigma\mathcal R_{\mathrm{int}+}^{5,j} \|_{L^{\frac4{3-2\alpha}}(\Sigma)} \|\chi_{0,\mathrm{int},\Sigma}v\|_{L^\infty(\Sigma)} \|\chi_{0,\mathrm{int},\Sigma}\eta_j\|_{H^{\frac32-\alpha}(\Sigma)} \notag \\
&+C(\varepsilon)\Bigl[ \|\chi_{0,\mathrm{int},\Sigma}\mathcal R'(\nabla_\Sigma\eta)\|_{L^\infty} \|\chi_{0,\mathrm{int},\Sigma}\eta_j\|_{H^1} +\|\mathcal R_{\mathrm{int}+}^{5,j}\|_{L^2}\Bigr] \times\sum_{r=1}^{j-1}\binom jr \|\chi_{0,\mathrm{int},\Sigma}(v_{j-r}\cdot\partial_t^r\mathcal N) \|_{H^1(\Sigma)}.
\end{align}
The second line includes both product-rule terms, including $\mathcal C_j$.

Recall that
\begin{align}
\mathcal{R}_{\mathrm{int}+}^{5,j} = \chi_{\mathrm{int}+,\Sigma}\sum_{\ell=2}^{j}\frac1{\ell!} \sum_{\substack{ i_1+\cdots+i_\ell=j \\ i_k\geq1}} \frac{j!}{i_1!\cdots i_\ell!} \mathcal{R}^{(\ell)}(\nabla_\Sigma\eta) \left[ \nabla_\Sigma\eta_{i_1},\ldots, \nabla_\Sigma\eta_{i_\ell} \right].
\end{align}
A spatial derivative can fall on the cutoff, on the coefficient $\mathcal{R}^{(\ell)}(\nabla_\Sigma\eta)$, or on one of the differentiated factors. Assuming that the derivatives of $\mathcal{R}$ through order $j+1$ are uniformly bounded on the range of $\nabla_\Sigma\eta$, using the Leibniz rule and H\"older's inequality, we have
\begin{align}
&\left\| \func{div}_\Sigma \mathcal{R}_{\mathrm{int}+}^{5,j} \right\|_{L^{\frac4{3-2\alpha}}(\Sigma)} \leq C(\varepsilon) \sum_{\ell=2}^{j} \sum_{\substack{ i_1+\cdots+i_\ell=j \\ i_k\geq1}} \Bigg[ \left( 1+ \left\| \chi_{1,\mathrm{int},\Sigma}\nabla_\Sigma^2\eta \right\|_{L^{\frac4{3-2\alpha}}(\Sigma)} \right) \prod_{k=1}^{\ell} \left\| \chi_{1,\mathrm{int},\Sigma}\nabla_\Sigma\eta_{i_k} \right\|_{L^\infty(\Sigma)} \notag \\
&+ \sum_{a=1}^{\ell} \left\| \chi_{1,\mathrm{int},\Sigma}\nabla_\Sigma^2\eta_{i_a} \right\|_{L^{\frac4{3-2\alpha}}(\Sigma)} \prod_{\substack{ k=1 \\ k\neq a}}^{\ell} \left\| \chi_{1,\mathrm{int},\Sigma}\nabla_\Sigma\eta_{i_k} \right\|_{L^\infty(\Sigma)} \Bigg] \notag \\
&\leq C(\varepsilon) \left( 1+ \left\| \chi_{1,\mathrm{int},\Sigma}\eta \right\|_{W^{3-\frac1p,p}(\Sigma)} \right) \sum_{\ell=2}^{j} \sum_{\substack{ i_1+\cdots+i_\ell=j \\ i_k\geq1}} \prod_{k=1}^{\ell} \left\| \chi_{1,\mathrm{int},\Sigma}\eta_{i_k} \right\|_{W^{3-\frac1p,p}(\Sigma)} \leq C(\varepsilon)\sqrt{\mathcal{D}_{\mathrm{int}}\mathcal{E}_{\mathrm{int}}}.
\end{align}
Here the cutoff is compactly supported, and $C(\varepsilon)$ may depend on $j$ and the fixed scaled cutoff profiles. The derivative of the cutoff is included in this bound. Since $\ell\geq2$, every index satisfies $i_k\leq j-1$. Therefore, the corresponding contribution to the energy estimate satisfies
\begin{align}
&\left| \sigma\int_\Sigma \func{div}_\Sigma \mathcal{R}_{\mathrm{int}+}^{5,j} \left(v_\parallel\cdot\nabla_\Sigma\eta_j\right) \right| \leq C(\varepsilon) \left\|\chi_{1,\mathrm{int}}v\right\|_{W^{2,p}(\Omega)} \left\| \chi_{0,\mathrm{int},\Sigma}\eta_j \right\|_{H^{\frac32-\alpha}(\Sigma)} \left( 1+ \left\| \chi_{1,\mathrm{int},\Sigma}\eta \right\|_{W^{3-\frac1p,p}(\Sigma)} \right) \notag \\
&\times \sum_{\ell=2}^{j} \sum_{\substack{ i_1+\cdots+i_\ell=j \\ i_k\geq1}} \prod_{k=1}^{\ell} \left\| \chi_{1,\mathrm{int},\Sigma}\eta_{i_k} \right\|_{W^{3-\frac1p,p}(\Sigma)} \leq C(\varepsilon)\sqrt{\mathcal{E}}\mathcal{D}.
\end{align}

For the term $\|\chi_{0,\mathrm{int},\Sigma}(v_{j-r}\cdot \partial_{t}\mathcal{N})\|_{H^{1}(\Sigma)}$, we have
\begin{align}
&\|\chi_{0,\mathrm{int},\Sigma}(v_{j-r}\cdot \partial_{t}^{r}\mathcal{N})\|_{H^{1}(\Sigma)}  \notag \\
\leq& C(\varepsilon) \|\nabla (\chi_{1,\mathrm{int},\Sigma}v_{j-r})\|_{L^{2}}\|\chi_{1,\mathrm{int},\Sigma}\eta_{r}\|_{W^{1,\infty}(\Sigma)}+C(\varepsilon)\| (\chi_{1,\mathrm{int},\Sigma}v_{j-r})\|_{L^{\infty}}\|\chi_{1,\mathrm{int},\Sigma}\eta_{r}\|_{H^{2}(\Sigma)}  \notag \\
\leq& C(\varepsilon)\|(\chi_{1,\mathrm{int},\Sigma}v_{j-r})\|_{W^{2,p}}\|\chi_{1,\mathrm{int},\Sigma}\eta_{r}\|_{W^{3-\frac{1}{p},p}(\Sigma)}+C(\varepsilon)\|(\chi_{1,\mathrm{int}}v_{j-r})\|_{W^{2,p}}\|\chi_{1,\mathrm{int},\Sigma}\eta_{r}\|_{W^{3-\frac{1}{p},p}(\Sigma)}  \notag \\
\leq& C(\varepsilon) \sqrt{\mathcal{E}\mathcal{D}}.
\end{align}

Combining the spatial integration by parts, the estimate for $\func{div}_\Sigma\mathcal R_{\mathrm{int}+}^{5,j}$, and the full product-rule estimate for $\nabla_\Sigma\mathcal R_{\mathrm{int}+}^{4,j}$, we obtain
\begin{equation}
\left|\sigma\int_\Sigma\mathcal R_{\mathrm{int}+}^{5,j} \cdot\nabla_\Sigma(v_\parallel\cdot\nabla_\Sigma\eta_j)\right| +\left|\sigma\int_\Sigma \left(\mathcal R'(\nabla_\Sigma\eta)\nabla_\Sigma\eta_j +\mathcal R_{\mathrm{int}+}^{5,j}\right) \cdot\nabla_\Sigma\mathcal R_{\mathrm{int}+}^{4,j}\right| \leq C(\varepsilon)\sqrt{\mathcal E}\mathcal D.
\end{equation}
Therefore, we have
\begin{align}
\sigma \int_{\Sigma}F_{0}^{5,j}\cdot \nabla_{\Sigma}F_{0}^{4,j}\leq C(\varepsilon)\sqrt{\mathcal{E}}\mathcal{D}.
\end{align}

\textbf{Boundary term $\langle \partial_{t}\eta_{j},F_{0}^{4,j}\rangle_{\Gamma_{0}}:=\langle \partial_{t}\eta_{j},F_{0}^{4,j}\rangle_{L^{2}(\Gamma_{0})}$}

$~$

By the definition of $F_0^{4,j}$ and the impermeability condition $v_k\cdot\nu=0$, the normal component of each $v_k$ vanishes on the wall. Hence,
\begin{align}
\kappa\left\langle \partial_t\eta_j,F_0^{4,j} \right\rangle_{\Gamma_0} = \kappa\left\langle \partial_t\eta_j, -\sum_{k=0}^{j-1}\binom jk v_{k,\gamma}\, \partial_t^{j-k}\partial_\gamma\eta \right\rangle_{\Gamma_0}.
\end{align}
Here $\gamma$ denotes the boundary curve, $v_{k,\gamma}$ is the tangential component of $v_k$ along $\gamma$, and $\partial_\gamma$ denotes differentiation with respect to arclength along $\gamma$. By the definition of the local diffeomorphism, we have the following computation on the fixed boundary
\begin{equation}
\partial_{\theta}\mathfrak{h}_{\eta} =\partial_{\theta}\eta(X(\theta,r,H_{0}(\theta,r)))=\nabla \eta \cdot \partial_{\theta}(X(\theta,r,H_{0}(\theta,r))) =\nabla \eta\cdot {P}_{x,y}\nabla (T(\theta))=\nabla \eta \cdot \dot{\gamma}(\theta)=\nabla_{\gamma}\eta.
\end{equation}
Therefore  $|\partial_{\gamma}\eta|\leq C |\partial_{\theta}\mathfrak{h}_{\eta}|$, where $C$ only depends on the steady state $\zeta$. Therefore,
\begin{align}
&\kappa\left\langle \partial_t\eta_j,F_0^{4,j} \right\rangle_{\Gamma_0} \lesssim [\partial_{t}\eta_{j}]_{\Gamma_{0}}[\partial_{\theta}\eta_{j}]_{\Gamma_{0}} \|v\|_{L^{\infty}(\Gamma_{0})}+[\partial_{t}\eta_{j}]_{\Gamma_{0}}[\partial_{\theta}\eta_{j-1}]_{\Gamma_{0}} \|\chi_{\mathrm{cl}}\partial_{t}v\|_{L^{\infty}W^{2,q_{\delta}}} \notag \\
&+\sum_{k=2}^{j-1}[\partial_{t}\eta_{j}]_{\Gamma_{0}}\|v_{k}\|_{L^{2}(\Gamma_{0})} \|\partial_{t}^{j-k}\partial_{\theta}\eta\|_{L^{\infty}(\Gamma_{0})} \lesssim\sqrt{\mathcal{E}}\mathcal{D}.
\end{align}
In the term $k=0$ the factor $[\partial_{\theta}\eta_{j}]_{\Gamma_{0}}=\|\partial_{t}^{j}\partial_{\theta}\eta\|_{L^{2}(\Gamma_{0})}$ is a contact-line dissipation term: for $j=m+1$ it belongs to $\mathscr{E}_{m,1,\mathrm{energy},\mathrm{cl}}\subset\mathcal{D}_{\mathrm{cl}}$, and for $j\leq m$ to the corresponding terms of lower order. Moreover $\|v\|_{L^{\infty}(\Gamma_{0})}\lesssim\sqrt{\mathcal{E}}$ by Lemma~\ref{lem:tl_coef}. No norm of $\partial_{\theta}\eta_{m+1}$ in $H^{3/2-\alpha}(\Sigma)$ is needed; that norm is not part of $\mathcal{D}$. In the term $k=1$ the three factors are bounded by $\sqrt{\mathcal{D}}$, $\sqrt{\mathcal{E}}$ and $\sqrt{\mathcal{E}}$, and $\mathcal{E}\sqrt{\mathcal{D}}\lesssim\sqrt{\mathcal{E}}\,\mathcal{D}$ by Lemma~\ref{lem:E_le_D}. In the terms $k\geq2$, $\|v_{k}\|_{L^{2}(\Gamma_{0})}\lesssim\|\chi_{\mathrm{cl}}\partial_{t}^{k}v\|_{L^{2}W^{2,q_{\delta}}}\lesssim\sqrt{\mathcal{D}}$, because $W^{2,q_{\delta}}\subset C^{0}$ on the slices and $k\leq m$. Moreover $\|\partial_{t}^{j-k}\partial_{\theta}\eta\|_{L^{\infty}(\Gamma_{0})}\lesssim\|\partial_{t}^{j-k}\partial_{\theta}\eta\|_{H^{1}(\Gamma_{0})}\lesssim\sqrt{\mathcal{E}}$, because $j-k\leq m-1$ and the contact-line terms $\|\partial_{t}^{j-k}\partial_{\theta}^{i}\eta\|_{L^{2}(\Gamma_{0})}$, $i\leq2$, belong to $\mathcal{E}_{\mathrm{cl}}$. The factor $\|\chi_{\mathrm{cl}}\partial_{t}^{k}v\|_{L^{\infty}W^{2,q_{\delta}}}$ is not used, since it is not controlled for $k\geq m-1$.

From the computation above we observe that $[F_{0}^{4,j}]_{\Gamma_{0}}\lesssim\sqrt{\mathcal{E}\mathcal{D}}$.

\textbf{Contact line term $\left\langle F_{0}^{7,j},\partial_{t}\eta_{j}-F_{0}^{4,j}\right\rangle_{\Gamma_{0}}$}

From the previous computation, we have
\begin{align}
\left\langle F_{0}^{7,j},\partial_{t}\eta_{j}-F_{0}^{4,j}\right\rangle_{\Gamma_{0}}\lesssim \left\langle F_{0}^{7,j}\right\rangle_{\Gamma_{0}}(\left\langle\partial_{t}\eta_{j}\right\rangle_{\Gamma_{0}}+\left\langle F_{0}^{4,j}\right\rangle_{\Gamma_{0}})\lesssim \sqrt{\mathcal{E}}\mathcal{D}.
\end{align}

\textbf{Energy dissipation relation when $j= m+1$}

We apply the Bogovskii operator, which is stated in the proposition after this main theorem. Taking the test function to be $u=(\partial_{t}^{m+1}v-\omega)$ in equation \eqref{eq:whole-domain-weak-momentum-j}, we obtain the following relation for $j=m+1$:
\begin{align}
&\frac{d}{dt} \Bigg[ \frac12\int_\Omega J|v_j|^2 + \frac g2\|\eta_j\|_{L^2(\Sigma)}^2 + \frac{\sigma}{2}\int_\Sigma a_\zeta\nabla_\Sigma\eta_j\cdot\nabla_\Sigma\eta_j -\int_{\Omega}Jv_{j}\omega-\int_{\Omega}p_{j-1}\langle JF_{0}^{2,j}\rangle\Bigg]  \notag \\
&\quad+ \frac{\mu}{2}\int_\Omega J|\mathbb{D}_{\mathcal{A}_\eta}v_j|^2 + \beta_0\int_{\Sigma_s}J_s|U_0v_j|^2 + \kappa\|\partial_t\eta_j\|_{L^2(\Gamma_0)}^2  \notag \\
&\qquad= \frac12\int_\Omega(\partial_tJ)|v_j|^2 + \frac{\sigma}{2}\int_\Sigma (\partial_ta_\zeta) \nabla_\Sigma\eta_j\cdot\nabla_\Sigma\eta_j-\int_{\Omega}Jv_{j}\partial_{t}\omega+((v_{j},\omega))-\int_{\Omega}\partial_{t}Jv_{j}\omega  \notag \\
&\qquad\quad- \int_\Omega p_{j-1}\partial_{t}\langle JF_{0}^{2,j}\rangle + \int_\Omega JF_0^{1,j}\cdot (v_j-\omega) - \left\langle F_0^{3,j},v_j\right\rangle_\Sigma - \left\langle J_sF_0^{6,j},U_0v_j\right\rangle_{\Sigma_s}  \notag \\
&\qquad\quad+ g\int_\Sigma\eta_jF_0^{4,j} + \sigma\int_\Sigma a_\zeta\nabla_\Sigma\eta_j\cdot\nabla_\Sigma F_0^{4,j} + \kappa\left\langle \partial_t\eta_j,F_0^{4,j} \right\rangle_{\Gamma_0}  \notag \\
&\qquad\quad- \left\langle F_0^{7,j},\partial_t\eta_j-F_0^{4,j} \right\rangle_{\Gamma_0} - \sigma\int_\Sigma F_0^{5,j}\cdot\nabla_\Sigma\partial_t\eta_j + \sigma\int_\Sigma F_0^{5,j}\cdot\nabla_\Sigma F_0^{4,j}.  \label{eq:whole-domain-energy-m+1}
\end{align}
Except for the terms with $\omega$, the nonlinear terms are estimated as above. The estimates for the terms with $\omega$ are proved in Proposition~\ref{prop:B_d} and Theorem \ref{thm:energy_pressure} after the proof of this theorem. We have
\begin{equation}
\left| \int_{\Omega}\partial_t^{m+1}v\,J\omega \right|+\left|\int_{\Omega}p^{m}\langle JF_{0}^{2,m+1}\rangle\right| \lesssim \mathcal{E}^{3/2},
\end{equation}
and
\begin{equation}
\begin{aligned} &\left| \int_{\Omega} \partial_t^{m+1}v\,\partial_t(J\omega) \right| +\left| \bigl(\!\bigl(\partial_t^{m+1}v,\omega\bigr)\!\bigr) \right| +\left| \int_{\Omega}JF_{0}^{1,m+1}\cdot\omega \right|+\left|\int_{\Omega}p^{m}\partial_{t}\langle JF_{0}^{2,m+1}\rangle\right| \lesssim \bigl(\sqrt{\mathcal{E}}+\mathcal{E}\bigr)\mathcal{D}. \end{aligned}
\end{equation}

\textbf{Main result in this theorem}

We add the identity at $j=0$, with the correction $\mathscr R_0^5$, to those for $1\leq j\leq m+1$, with the corrections $\mathscr R_j^5$ from \eqref{eq:Fint5-correction}. Each correction has absolute value at most $C(\varepsilon)\mathcal E^{3/2}$. Combining all the estimates in this proof, we obtain
\begin{align}
\mathcal{E}_{||}(t)-\mathcal{E}^{\frac{3}{2}}(t)+\int_{0}^{t}\mathcal{D}_{||}\lesssim \mathcal{E}_{||}(0)+\int_{0}^{t}\sqrt{\mathcal{E}}\mathcal{D}+\mathcal{E}^{\frac{3}{2}}(0).
\end{align}
This finishes the proof.
\end{proof}

Using the estimates above, we derive the following corollary.
\begin{corollary}{\label{energy_whole_enhance}} Suppose $0\leq j\leq m$. Then we have the following estimate
\begin{align}
\sum_{j=0}^{m}(\|\p_{t}^{j}v\|_{H^{1}}^{2}+\|\p_{t}^{j+1}\eta\|_{L^{2}(\Gamma_{0})}^{2})(t)\lesssim \mathcal{E}(0)+\int_{0}^{t}\mathcal{D}_{||}(s)ds.
\end{align}
\end{corollary}
\begin{proof}
By the fundamental theorem of calculus, we have
\begin{align}
\sum_{j=0}^{m}(\|\p_{t}^{j}v\|_{H^{1}}^{2}+\|\p_{t}^{j+1}\eta\|_{L^{2}(\Gamma_{0})}^{2})(t) \lesssim \sum_{j=0}^{m}(\|\p_{t}^{j}v\|_{H^{1}}^{2}+\|\p_{t}^{j+1}\eta\|_{L^{2}(\Gamma_{0})}^{2})(t)(0)+\int_{0}^{t}(\sum_{j=0}^{m+1}(\|\p_{t}^{j}v\|_{H^{1}}^{2}+\|\p_{t}^{j+1}\eta\|_{L^{2}(\Gamma_{0})}^{2}))(s)ds.
\end{align}
This finishes the proof.
\end{proof}

\textbf{Bogovskii operator at the highest order}

When $j=m+1$, we have to change the test function from $Jv_{j }$ to $J(v_{k}-\omega)$, where $\omega$ satisfies the following properties.

\begin{proposition}
{\label{prop:B_d}} Let $0<T\leq\infty$ and suppose that  $\sup_{0\leq t<T}\mathcal{E}(t)\leq\gamma^2$, where $\gamma\in(0,1)$ is small enough that $\|J-1\|_{L^{\infty}}+\|K-1\|_{L^{\infty}}+\|\mathcal{A}_{\eta}-I\|_{L^{\infty}}\leq\frac{1}{2}$, so that $\Phi_{\eta}$ is a diffeomorphism (Subsection~5.1.1). Let $F^{2}:=F^{2,m+1}_{0}$ be given by \eqref{eq:F02j}, and let $\langle\cdot\rangle_{\Omega}$ denote the spatial average on $\Omega$, i.e.,
\begin{equation}
\langle g\rangle_{\Omega} =\frac{1}{|\Omega|}\int_{\Omega}g.
\end{equation}
Then there exists $\omega:\Omega\times[0,T)\to\mathbb{R}^3$ satisfying the following:

\begin{enumerate}
\item $\omega(\cdot,t)\in H_0^1(\Omega;\mathbb{R}^3)$  for $0\leq t<T$, and $J\func{div}_{\mathcal{A}}\omega =JF_{0}^{2,m+1}-\langle JF_{0}^{2,m+1}\rangle_{\Omega}$.

\item $\omega$ obeys the estimates
\begin{equation}
\|\omega\|_{W_0^{1,\min\{4/(3-2\varepsilon_-),p\}}} \lesssim \mathcal{E}, \qquad \|\omega\|_{W_0^{1,2}} +\|\partial_t\omega\|_{L^{2}} \lesssim \bigl(\sqrt{\mathcal{E}}+\mathcal{E}\bigr) \sqrt{\mathcal{D}}.
\end{equation}

\item We have the interaction estimates
\begin{equation}
\left| \int_{\Omega}\partial_t^{m+1}v\,J\omega \right| \lesssim \mathcal{E}^{3/2},
\end{equation}
and
\begin{equation}
\begin{aligned} &\left| \int_{\Omega} \partial_t^{m+1}v\,\partial_t(J\omega) \right| +\left| \bigl(\!\bigl(\partial_t^{m+1}v,\omega\bigr)\!\bigr) \right| +\left| \int_{\Omega}JF_{0}^{1,m+1}\cdot\omega \right| \lesssim \bigl(\sqrt{\mathcal{E}}+\mathcal{E}\bigr)\mathcal{D}. \end{aligned}
\end{equation}
\end{enumerate}
\end{proposition}

\begin{proof}
\emph{Construction.} We use one right inverse $B$ of the divergence for all the bounds below. On the bounded Lipschitz domain $\Omega$, \cite[Section~4]{CM} provides linear operators $B_{0}$ and $K$ with $\func{div}B_{0}=I-K$. Here $B_{0}$ is bounded from $H^{s}_{p,\overline\Omega}(\mathbb{R}^{3})$ to $H^{s+1}_{p,\overline\Omega}(\mathbb{R}^{3};\mathbb{R}^{3})$ for every $s\in\mathbb{R}$ and $1<p<\infty$, where $H^{s}_{p,\overline\Omega}$ denotes the distributions in the Bessel potential space $H^{s}_{p}(\mathbb{R}^{3})$ supported in $\overline\Omega$, and $K$ is smoothing. Let $B_{1}$ be the Bogovski\u{\i} operator of $\Omega$ \cite[Section~III.3]{Galdi}, a right inverse of the divergence that is bounded from $L^{q}_{0}(\Omega):=\{g\in L^{q}(\Omega):\int_{\Omega}g=0\}$ to $W^{1,q}_{0}(\Omega)$, $1<q<\infty$. Put $B:=B_{0}+B_{1}K$. If $\langle g,1\rangle=0$, then $\langle Kg,1\rangle=\langle g,1\rangle-\langle\func{div}B_{0}g,1\rangle=0$, since $B_{0}g$ is supported in $\overline\Omega$. Hence $\func{div}Bg=g$. Moreover, $B$ is bounded from $L^{q}_{0}(\Omega)$ to $W^{1,q}_{0}(\Omega)$, and from $\{h\in(H^{1}(\Omega))^{*}:h(1)=0\}=\{h\in H^{-1}_{2,\overline\Omega}:\langle h,1\rangle=0\}$ to $L^{2}(\Omega)$. Put $g:=JF^{2,m+1}_{0}-\langle JF^{2,m+1}_{0}\rangle_{\Omega}$, $\bar\omega:=Bg$, $\omega:=M\bar\omega$, and $M:=J^{-1}\mathcal{A} _{\eta}^{-T}$.  By the Piola identity, $J\func{div}_{\mathcal{A}}(M\bar\omega)=\func{div}(J\mathcal{A}^{T}M\bar\omega)=\func{div}\bar\omega=g$, which is property (1), and $\omega\in H^{1}_{0}$.

\emph{The bounds on $\omega$.} $\|\omega\|_{W^{1,q}}\lesssim\|g\|_{L^{q}}$, since $|M-I|\lesssim\sqrt{\mathcal{E}}$ and $\nabla M\in L^{\infty}_{\theta}L^{r_{-}}_{rh}$ near $\Gamma$ by Lemma~\ref{lem:tl_coef}. By \eqref{eq:F02j}, $F^{2,m+1}_{0}$ is a sum of products $\nabla\partial_{t}^{r}\bar\eta\,\nabla v_{m+1-r}$, $1\leq r\leq m+1$, with smooth coefficients. Near $\Gamma$ we use mixed norms. The worst factor is $\nabla\partial_{t}^{m+1}\bar\eta$. For the slice-wise pieces, \eqref{tl:P2} and \eqref{tl:aniso} give $\nabla\partial_{t}^{m+1}\bar\eta\in L^{2}_{\theta}H^{1/2}_{rh}\subset L^{2}_{\theta}L^{4}_{rh}$ with norm $\lesssim\|\partial_{t}^{m+1}\eta\|_{H^{1}(\Sigma)}\leq\sqrt{\mathcal{E}}$. The same bounds give $\nabla\partial_{t}^{m+1}\bar\eta\in L^{2}_{\theta}H^{1-\alpha}_{rh}\subset L^{2}_{\theta}L^{2/\alpha}_{rh}$ with norm $\lesssim\|\partial_{t}^{m+1}\eta\|_{H^{3/2-\alpha}(\Sigma)}\leq\sqrt{\mathcal{D}}$. The $\theta$-derivative of the slice-wise extension of an $H^{1}$ datum $G$ is in $L^{2}_{\theta}H^{1/2}_{rh}$, by \eqref{tl:P2} applied to $\partial_{\theta}G$ and $\|\partial_{\theta}G\|_{L^{2}}\leq\|G\|_{H^{1}}$. The other factor is $\nabla v\in L^{\infty}_{\theta}L^{r_{-}}_{rh}$, with norm $\lesssim\sqrt{\mathcal{E}}$. The terms with $r\leq m$ have a factor bounded in $L^{\infty}$, as in the interior computation below. Set $q_*:=\frac4{3-2\varepsilon_-},\; q_E:=\min\{q_*,p\}$. Near the contact line, the energy bounds, the fact that $\frac14+\frac{1-\varepsilon_{-}}{2}=\frac{3-2\varepsilon_{-}}{4}$ and $\frac{\alpha}{2}+\frac{1-\varepsilon_{-}}{2}<\frac12$, and the above argument give the mixed norm $\|\nabla \bar{\eta}_{m+1}\,\nabla v\|_{L^2_\theta L^{q_*}_{rh}} \le C(\varepsilon)\mathcal E$. Since $q_E\le p<2$ and $q_E\le q_*$, the bounded chart implies $L^2_\theta L^{q_*}_{rh}\hookrightarrow L^{q_E}(\mathcal U), \; \|g\|_{L^{q_E}(\Omega)}\le C(\varepsilon)\mathcal E$. Use the same exponent in the lift and its application. As $q_E>6/5$,
\begin{equation}
\|\omega\|_{L^2} \le C\|\omega\|_{W^{1,q_E}_0} \le C(\varepsilon)\mathcal E, \qquad \left|\int_\Omega Jv_{m+1}\cdot\omega\right| \le C(\varepsilon)\mathcal E^{3/2}.
\end{equation}
The separate dissipative bound uses $\alpha<\varepsilon_-$ and is unchanged: $\|\nabla Q_{m+1}\,\nabla v\|_{L^2} \le C(\varepsilon)\sqrt{\mathcal E\mathcal D}$. Away from $\Gamma$, the computation is the following.
\begin{align}
&\|F_{\mathrm{int}+}^{2,m+1}\|_{L^{\min\{\frac{4}{3-2\epsilon_{+}},p\}}}\lesssim  \sum_{k=1}^{m+1}\|\chi_{\mathrm{int}+}\func{div}_{\partial_{t}^{k}\mathcal{A}}v_{m+1-k}\|_{L^{\min\{\frac{4}{3-2\epsilon_{+}},p\}}} \notag \\
&\lesssim  \sum_{k=2}^{m-1}\|\chi_{\mathrm{int}+}\func{div}_{\partial_{t}^{k}\mathcal{A}}v_{m+1-k}\|_{L^{\frac{4}{3-2\epsilon_{+}}}}+\|\chi_{\mathrm{int}+}\func{div}_{\partial_{t}\mathcal{A}}v_{m}\|_{L^{p}}+\|\chi_{\mathrm{int}+}\func{div}_{\partial_{t}^{m}\mathcal{A}}\partial_{t}v\|_{L^{p}} +\|\chi_{\mathrm{int}+}\func{div}_{\partial_{t}^{m+1}\mathcal{A}}v\|_{L^{p}} \notag \\
&\lesssim \sum_{k=2}^{m-1}\|\chi_{0-,\mathrm{int}}\nabla {\partial_{t}^{k}\bar{\eta}}\|_{L^{\infty}}\|\chi_{0-,\mathrm{int}}\nabla v_{m+1-k}\|_{L^{\frac{4}{3-2\epsilon_{+}}}}+\|\chi_{0-,\mathrm{int}}\nabla {\partial_{t}\bar{\eta}}\|_{L^{\infty}}\|\chi_{0-,\mathrm{int}}\nabla v_{m}\|_{L^{p}} \notag \\
&+\|\chi_{0-,\mathrm{int}}\nabla {\partial_{t}^{m}\bar{\eta}}\|_{L^{p}}\|\chi_{0-,\mathrm{int}}\nabla \partial_{t}v\|_{L^{\infty}}+\|\chi_{0-,\mathrm{int}}\nabla {\partial_{t}^{m+1}\bar{\eta}}\|_{L^{p}}\|\chi_{0-,\mathrm{int}}\nabla v\|_{L^{\infty}} \lesssim \mathcal{E}.
\end{align}
The $L^{2}$ bound of the interior part is the same computation with $L^{2}$ norms, and gives $\sqrt{\mathcal{E}\mathcal{D}}$. Hence the first two bounds of (2) hold. The exponent $4/(3-2\varepsilon_{-})$ in (2) comes from $\nabla v\in L^{\infty}_{\theta}L^{r_{-}}_{rh}$; the energy controls only the exponent $q_{\delta}=q_{-}$ near $\Gamma$ (Subsection~\ref{subsec:constants}).

\emph{The time derivative: divergence form.} The bound $\|\partial_{t}\omega\|_{L^{2}}$ cannot be obtained from $\|\partial_{t}g\|_{L^{q}}$: $\partial_{t}g$ contains $\nabla\partial_{t}^{m+2}\bar\eta\,\nabla v$, and $\partial_{\theta}$ of the slice-wise extension of the datum $\partial_{t}^{m+2}\eta\in H^{1/2-\alpha}(\Sigma)$ lies in no $L^{q}$; for $e^{ik\theta}\varphi(r)$ its norm grows like $k$, against $k^{1/2-\alpha}$ for the datum. We use instead the structure of $g$. By \eqref{eq:JKA-definition}--\eqref{eq:A-explicit}, $J\mathcal{A}=(1+\partial_{3}q)I-\nabla q\otimes e_{3}$ is affine in $\nabla q$, $q=W\bar\eta$, so  $J\mathcal{A}^{T}u=(1+\partial_{3}q)\,u-(\nabla q\cdot u)\,e_{3}$. Since $\func{div}(J\mathcal{A}^{T}v)=J\func{div}_{\mathcal{A}}v=0$, differentiating $j$ times in $t$ gives
\begin{equation}
\begin{aligned}  \label{eq:F2_divergence} JF^{2,j}_{0}&=J\func{div}_{\mathcal{A}}v_{j}=\func{div}(J\mathcal{A}^{T}v_{j})=-\func{div}G_{j},\qquad G_{j}:=\sum_{r=1}^{j}\binom{j}{r}\Big((\partial_{3}Q_{r})\,v_{j-r}-(\nabla Q_{r}\cdot v_{j-r})\,e_{3}\Big),
\end{aligned}
\end{equation}
with $Q_{r}:=\partial_{t}^{r}q$. On the vertical part of $\Sigma_{s}$, $v_{l}\cdot n_{0}=0$ and $e_{3}\cdot n_{0}=0$, and $Q_{r}=0$ below the level $h_{0}$; hence $G_{j}\cdot n_{0}=0$ on $\Sigma_{s}$.

\emph{The time derivative: the null structure.} For a vector field $u$ and a function $Q$,  $\func{div}\big((\partial_{3}Q)\,u-(\nabla Q\cdot u)\,e_{3}\big)=(\func{div}u)\,\partial_{3}Q-\partial_{3}u\cdot\nabla Q$: the second derivatives of $Q$ cancel. Integrating by parts once more, for $\varphi\in H^{1}(\Omega)$,
\begin{equation}
\int_{\Omega}\func{div}\big((\partial_{3}Q)u-(\nabla Q\cdot u)e_{3}\big)\varphi =-\int_{\Omega}Q\big((\func{div}u)\,\partial_{3}\varphi-\partial_{3}u\cdot\nabla\varphi\big) +\int_{\Sigma}Q\,\varphi\,\big((\func{div}u)\,\nu_{3}-\partial_{3}u\cdot\nu\big)\,dS . \label{eq:null_structure}
\end{equation}
The terms on $\Sigma_{s}$ vanish when $Q=0$ below $h_{0}$ and $u\cdot n_{0}=0 $ on the vertical wall: there $n_{3}=0$, and $\partial_{3}u\cdot n_{0}=\partial_{3}(u\cdot n_{0})=0$ because $\partial_{3}$ is tangent to the wall and $n_{0}$ does not depend on $z$. No derivative of $Q$ remains in \eqref{eq:null_structure}.

\emph{The time derivative: the bound.} By \eqref{eq:F2_divergence}, $\partial_{t}g=-\func{div}\partial_{t}G_{m+1}-\partial_{t}\langle JF^{2,m+1}_{0}\rangle_{\Omega}$. We regard $\partial_{t}g$ as a functional on $H^{1}(\Omega)$, and we have $\langle\partial_{t}g,1\rangle=0$. In $\partial_{t}G_{m+1}$, the only term with $m+2$ time derivatives on $\eta$ is $(\partial_{3}Q)v-(\nabla Q\cdot v)e_{3}$ with $Q=Q_{m+2}=W\partial_{t}^{m+2}\bar\eta$, and \eqref{eq:null_structure} applies to it with $u=v$. Near $\Gamma$, $Q\in L^{2}_{\theta}H^{1-\alpha}_{rh}\subset L^{2}_{\theta}L^{2/\alpha}_{rh}$ and $Q|_{\Sigma}=\partial_{t}^{m+2}\eta\in L^{2}_{\theta}H^{1/2-\alpha}_{r}\subset L^{2}_{\theta}L^{1/\alpha}_{r}$, both with norm $\lesssim\|\partial_{t}^{m+2}\eta\|_{H^{1/2-\alpha}(\Sigma)}$. Also $\nabla v\in L^{\infty}_{\theta}L^{r_{-}}_{rh}$ and $\nabla v|_{\Gamma_{\pm}}\in L^{\infty}_{\theta}L^{s_{-}}$, with norm $\lesssim\sqrt{\mathcal{E}}$ by Lemma~\ref{lem:tl_coef}. Finally $\varphi\in L^{2}_{\theta}H^{1}_{rh}$, and its trace $\varphi|_{\Sigma}\in L^{2}_{\theta}H^{1/2}_{r}\subset L^{2}_{\theta}L^{s}_{r}$ for every $s<\infty $. H\"older's inequality applies in both integrals: in $(r,h)$ because $\frac{\alpha}{2}+\frac{1-\varepsilon_{-}}{2}+\frac12<1$, and on $\Sigma$, with $s=1/(\varepsilon_{-}-\alpha)$, because $\alpha+(1-\varepsilon_{-})+\frac1s=1$. Away from $\Gamma$, $\nabla v\in L^{\infty}$. By the kinematic equation, $\partial_{t}^{m+2}\eta=v_{m+1}\cdot\mathcal{N}+F^{4,m+1}_{0}$, and $\|\partial_{t}^{m+2}\eta\|_{H^{1/2-\alpha}(\Sigma)}\lesssim\|v_{m+1}\|_{H^{1}}+\sqrt{\mathcal{E}\mathcal{D}}\lesssim\sqrt{\mathcal{D}}$. Hence this term is $\lesssim\sqrt{\mathcal{E}\mathcal{D}}\,\|\varphi\|_{H^{1}}$. Every other term of $\partial_{t}G_{m+1}$ is a product $(\nabla Q_{k})\,v_{l} $ with $k+l=m+2$ and $1\leq k,l\leq m+1$. These products lie in $L^{2}(\Omega)$, and their traces lie in $L^{2}(\Sigma)$, with norm $\lesssim(\sqrt{\mathcal{E}}+\mathcal{E})\sqrt{\mathcal{D}}$. For $k=m+1$ we use $\nabla Q_{m+1}\in L^{2}$ and $v_{1}\in L^{\infty}$; for $k=1$ we use $\nabla Q_{1}\in L^{\infty}$, by the proof of Lemma~\ref{lem:tl_coef}, and $v_{m+1}\in H^{1}$. For these terms we use $\langle-\func{div}H,\varphi\rangle=\int_{\Omega}H\cdot\nabla\varphi-\int_{\Sigma}(H\cdot\nu)\varphi\,dS$, since $H\cdot n_{0}=0$ on $\Sigma_{s}$. Finally $|\partial_{t}\langle JF^{2,m+1}_{0}\rangle_{\Omega}|\lesssim\sqrt{\mathcal{E}\mathcal{D}}$ by Theorem~\ref{thm:energy_pressure}. Altogether,
\begin{equation}  \label{eq:dtg_dual} \|\partial_{t}g\|_{(H^{1}(\Omega))^{*}}\lesssim(\sqrt{\mathcal{E}}+\mathcal{E})\sqrt{\mathcal{D}}.
\end{equation}

\emph{The time derivative: the right inverse in negative order.} By the construction, $B$ is bounded from $\{h\in(H^{1}(\Omega))^{*}:h(1)=0\}$ to $L^{2}(\Omega)$. Since $B$ is linear and does not depend on $t$, $\partial_{t}\bar\omega=B\,\partial_{t}g$, and by \eqref{eq:dtg_dual},
\begin{equation*}
\|\partial_{t}\bar\omega\|_{L^{2}} \lesssim(\sqrt{\mathcal{E}}+\mathcal{E})\sqrt{\mathcal{D}},\quad \|\partial_{t}\omega\|_{L^{2}} \leq\|\partial_{t}M\|_{L^{\infty}}\|\bar\omega\|_{L^{2}}+\|M\|_{L^{\infty}}\|\partial_{t}\bar\omega\|_{L^{2}}\lesssim(\sqrt{\mathcal{E}}+\mathcal{E})\sqrt{\mathcal{D}} .
\end{equation*}
This completes property (2).

For property (3), using property (2), we have
\begin{equation}
\left| \int_{\Omega}v_{m+1}\,J\omega \right| \lesssim \|v_{m+1}\|_{L^{2}}\|J\omega\|_{L^{2}}\lesssim\|v_{m+1}\|_{L^{2}}\|\omega\|_{W_{0}^{1,\min\{\frac{4}{3-2\epsilon_{-}},p\}}}\lesssim \mathcal{E}^{\frac{3}{2}},
\end{equation}
since $W^{1,4/(3-2\varepsilon_{-})}_{0}\subset L^{2}$ in three dimensions. Moreover, we have
\begin{align}
&\left| \int_{\Omega} v_{m+1}\,\partial_t(J\omega) \right| +\left| \bigl(\!\bigl(v_{m+1},\omega\bigr)\!\bigr) \right| +\left| \int_{\Omega}JF_{0}^{1,m+1}\cdot\omega \right| \notag \\
&\lesssim \|v_{m+1}\|_{L^{2}}\|\partial_{t}J\|_{L^{\infty}}\|\omega\|_{L^{2}}+\|v_{m+1}\|_{L^{2}}\|\partial_{t}\omega\|_{L^{2}}+\|\omega\|_{H^{1}}(\|F_{0}^{1,m+1}\|_{(H^{1})^{*}}+\|v_{m+1}\|_{H^{1}}) \lesssim \bigl(\sqrt{\mathcal{E}}+\mathcal{E}\bigr)\mathcal{D}.
\end{align}
\end{proof}

Using the computation in Proposition~\ref{prop:B_d}, we have the following theorem.

\begin{theorem}
{\label{thm:energy_pressure}} The following estimates hold for $F_{0}^{2,m+1} $
\begin{equation}
\partial_{t}\langle JF_{0}^{2,m+1}\rangle\lesssim\sqrt{\mathcal{E}}\sqrt{\mathcal{D}}, \langle JF_{0}^{2,m+1}\rangle\lesssim\mathcal{E}.
\end{equation}
Moreover, the following estimates hold for the pressure $\partial_{t}^{m}p$
\begin{equation}
\int_{\Omega}p_{m}\langle JF_{0}^{2,m+1}\rangle\lesssim \|p_{m}\|_{L^{2}} \langle JF_{0}^{2,m+1}\rangle\lesssim \sqrt{\mathcal{E}}^{3},\quad \int_{\Omega}p_{m}\partial_{t}\langle JF_{0}^{2,m+1}\rangle\lesssim \|p_{m}\|_{L^{2}} \partial_{t}\langle JF_{0}^{2,m+1}\rangle\lesssim \sqrt{\mathcal{E}}\mathcal{D}.
\end{equation}
\end{theorem}

\begin{proof}
The estimates for the pressure follow from those for the average: $|\int_{\Omega}p_{m}c|\leq|\Omega|^{1/2}\|p_{m}\|_{L^{2}}|c|$ for a constant $c$, with $\|p_{m}\|_{L^{2}}\lesssim\sqrt{\mathcal{E}}$ (it is part of $\mathcal{E}_{\mathrm{enhanced}}$) for the first, and $\|p_{m}\|_{L^{2}}\lesssim\sqrt{\mathcal{D}}$ for the second. The latter is the proof of Theorem~\ref{thm:pressure_enhancement}, with the forcing at the level $m$ bounded by $\sqrt{\mathcal{E}\mathcal{D}}$, as in \eqref{eq:contact_line_forcing_estimate_d}, instead of by $\mathcal{E}$.

For the average we use the Piola identity, not the $L^{q}$ norms of $F^{2,m+1}_{0}$ and its time derivative. By \eqref{eq:F2_divergence} and the divergence theorem,
\begin{equation*}
|\Omega|\,\langle JF^{2,m+1}_{0}\rangle_{\Omega}=\int_{\Omega}\func{div}(J\mathcal{A}^{T}v_{m+1})=\int_{\partial\Omega}v_{m+1}\cdot J\mathcal{A}n\,dS=\int_{D}v_{m+1}\cdot\mathcal{N}\,dx .
\end{equation*}
On $\Sigma_{s}$ we have $J\mathcal{A}n_{0}=Jn_{0}$ (or $J\mathcal{A}=I$ below $h_{0}$) and $v_{m+1}\cdot n_{0}=0$; on $\Sigma$, $J\mathcal{A}\nu\,dS=\mathcal{N}\,dx$. By the kinematic equation and $\int_{D}\partial_{t}^{m+2}\eta=0$ (Subsection~\ref{subsec:volume}),
\begin{equation*}
|\Omega|\,\langle JF^{2,m+1}_{0}\rangle_{\Omega}=-\int_{D}F^{4,m+1}_{0}\,dx=\sum_{r=1}^{m+1}\binom{m+1}{r}\int_{D}v_{m+1-r,x}\cdot\nabla_{x}\partial_{t}^{r}\eta\,dx,
\end{equation*}
since $\partial_{t}^{r}\mathcal{N}=-(\nabla_{x}\partial_{t}^{r}\eta,0)$ for $r\geq1$; here $v_{x}$ denotes the horizontal components. Each term is
\begin{equation*}
\lesssim\|v_{m+1-r}\|_{L^{2}(\Sigma)}\|\nabla\partial_{t}^{r}\eta\|_{L^{2}(\Sigma)}\lesssim\mathcal{E}
\end{equation*}
. In the time derivative, the only term with $m+2 $ time derivatives on $\eta$ is
\begin{equation*}
\int_{D}v_{x}\cdot\nabla_{x}\partial_{t}^{m+2}\eta\,dx=-\int_{D}(\func{div}_{x}v_{x})\,\partial_{t}^{m+2}\eta\,dx,
\end{equation*}
with no boundary term, because $v_{x}\cdot\nu=v\cdot n_{0}=0$ on $\partial D$. Near $\Gamma$ we bound it with mixed norms: $\nabla v|_{\Sigma}\in L^{\infty}_{\theta}L^{s_{-}}_{r}$ with norm $\lesssim\sqrt{\mathcal{E}}$ (Lemma~\ref{lem:tl_coef}), and $\partial_{t}^{m+2}\eta\in L^{2}_{\theta}H^{1/2-\alpha}_{r}\subset L^{2}_{\theta}L^{1/\alpha}_{r}$ with norm $\lesssim\|\partial_{t}^{m+2}\eta\|_{H^{1/2-\alpha}(\Sigma)}\lesssim\sqrt{\mathcal{D}}$, by the kinematic equation. Since $(1-\varepsilon_{-})+\alpha\leq1$, H\"older's inequality on the bounded slices applies. Away from $\Gamma$, $\nabla v\in L^{\infty}$. So the term is $\lesssim\sqrt{\mathcal{E}\mathcal{D}}$. The other terms are products of $v_{l}|_{\Sigma}$, $l\leq m+1$, with $\nabla\partial_{t}^{r}\eta$, $r\leq m+1$, and are $\lesssim\sqrt{\mathcal{E}\mathcal{D}}$. Hence $|\partial_{t}\langle JF^{2,m+1}_{0}\rangle_{\Omega}|\lesssim\sqrt{\mathcal{E}\mathcal{D}}$ and $|\langle JF^{2,m+1}_{0}\rangle_{\Omega}|\lesssim\mathcal{E}$.
\end{proof}

\begin{remark}
The time derivative of $\langle JF^{2,m+1}_{0}\rangle$ cannot be bounded through $\|\partial_{t}F^{2,m+1}_{0}\|_{L^{q_{-}}}$. Indeed $\partial_{t}F^{2,m+1}_{0}$ contains $\nabla\partial_{t}^{m+2}\bar\eta$, which lies in no $L^{q}$ near $\Gamma$, because the extension there acts slice by slice. The proof above avoids this by the divergence theorem. The proof of Proposition~\ref{prop:B_d} avoids it by the null structure \eqref{eq:null_structure}.
\end{remark}

\subsection{Enhanced regularity for the surface function}

\label{subsec:enhanced_surface}

In this subsection, we derive the $\frac{3}{2}-\alpha$ estimate for the surface function $\eta_{j}$. We use the functional calculus of the capillary operator $K=g-\sigma\func{div}(a_{\zeta}\nabla\,\cdot\,)$ on $D$ with the conormal boundary condition, as in \cite[Section~8]{GT}; the operators $D^{s}_{j}$ and the space $H^{1/2}_{K}$ used below are those of \cite[Section~8]{GT}. The facts needed in three dimensions are the following. The operator $K$ is self-adjoint on $L^{2}(D)$ with form domain $H^{1}(D)$, and the constants form its ground state, with eigenvalue $g$; in particular $\int_{D}D^{s}_{j}f=g^{s/2}\int_{D}f$. The identities of \cite[Proposition~8.26 and Lemma~8.28]{GT} are spectral and hold verbatim; \cite[Lemma~8.27]{GT} follows from the ground state; \cite[Theorem~8.23]{GT} is used only for $t\leq1$, where it reads $[L^{2},H^{1}]_{t}=H^{t}$; and \cite[Proposition~8.30]{GT} uses $H^{r}_{0}=H^{r}$ for $r\leq1/2$. For the $H^2$ regularity of the harmonic test on the three-dimensional vessel we use the edge regularity theorem in the proof below.

For the positive conormal realization of $K=g-\sigma\operatorname{div}(a_\zeta\nabla\cdot)$ on the smooth domain $D$, the fractional-domain theorem shows that $\mathcal H_K^t:=D(K^{t/2})=H^t(D),\; 0\le t<\frac32$, with equivalent norms. Let $P_j$ denote the projection onto the first $j+1$ eigenfunctions, including the constant eigenfunction, and set $D_j^s=K^{s/2}P_j$. For $0\le s<1$, $\|D_j^sf\|_{H^{1-s/2}} \lesssim\|P_jf\|_{\mathcal H_K^{1+s/2}} \le\|f\|_{\mathcal H_K^{1+s/2}} \lesssim\|f\|_{H^{1+s/2}}$, uniformly in $j$. Moreover, $(f,D_j^sf)_{1,\Sigma} =\|D_j^{s/2}f\|_{\mathcal H_K^1}^2 =\|P_jf\|_{\mathcal H_K^{1+s/2}}^2$. Thus the limit of this pairing is equivalent to $\|f\|_{H^{1+s/2}(D)}^2$. This statement does not include the critical case $s=1$. With the facts above, we construct $\psi$ via the following proposition.

\begin{proposition}
{\label{thm:whole_space_threehalf}\label{rem:psi_existence}} Let $0\leq k\leq m+1$, $j,k\in\mathbb N$, and $0\leq s<1$. There exists a unique $\psi\in H^2(\Omega)$ with $\int_\Omega\psi=0$ solving
\begin{equation}
\begin{cases}
-\Delta\psi=0 & \text{in }\Omega, \\
\partial_{\nu}\psi =\dfrac{D_j^s\partial_t^k\eta}{|\mathcal{N}_0|} & \text{on }\Sigma, \\
\partial_{\nu}\psi=0 & \text{on }\Sigma_s,\end{cases}\end{equation} where $\nu$ is the outward unit normal and $\nu=\mathcal N_0/|\mathcal N_0|$ on $\Sigma$. With constants independent of $j$,
\begin{equation}
\|\psi\|_{H^1(\Omega)} \lesssim\|\partial_t^k\eta\|_{\mathcal H_K^{s-1/2}}, \|\psi\|_{H^2(\Omega)} \lesssim\|D_j^s\partial_t^k\eta\|_{H^{1/2}(D)},\quad \|\partial_t\psi\|_{H^1(\Omega)} \lesssim\|\partial_t^{k+1}\eta\|_{\mathcal H_K^{s-1/2}} .
\end{equation}
Here the scale at negative index is the dual scale under the $L^2(D)$ pairing, as in \eqref{functional:negative}.
\end{proposition}

\begin{proof}
Write $f=\partial_t^k\eta$ and $G=D_j^sf$. The compatibility condition is $\int_D G\,dx=g^{s/2}\int_D f\,dx=0$ by conservation of volume (Subsection~\ref{subsec:volume}): on $\Sigma$ the surface measure is $dS=|\mathcal N_0|\,dx$. More generally, for any $G\in H^{1/2}(D)$ of zero horizontal mean, the weak problem is
\begin{equation}
\int_\Omega\nabla\psi\cdot\nabla\phi =\langle G,\phi|_\Sigma\rangle_D,\qquad \phi\in H^1(\Omega).
\end{equation}
The trace on the free face is bounded into $H^{1/2}(D)$. Lax--Milgram and Poincar\'e on $\{\phi\in H^1(\Omega):\int_\Omega\phi=0\}$ therefore give a unique normalized solution with $\|\psi\|_{H^1(\Omega)}\lesssim\|G\|_{(H^{1/2}(D))'}$.

By hypotheses (1)--(3), the two smooth faces of $\Omega$ meet along one smooth closed edge, with no vertices. Apply the Neumann case of \cite[Theorem~3, p.~17]{DaugeEdge} at regularity index $2$: the condition on the opening is $2-1<\pi/\theta_{\mathrm{eq}}$, since $0<\theta_{\mathrm{eq}}<\pi$. The bulk datum is zero, and the two face data are $G/|\mathcal N_0|$ and zero. They belong to the required facewise $H^{1/2}$ space, with norm bounded by $C\|G\|_{H^{1/2}(D)}$. Thus $\psi\in H^2(\Omega)$. The normalized solution operator on zero-mean $H^{1/2}(D)$ data is bounded into $H^2(\Omega)$ by the closed graph theorem: its graph is closed by the weak identity and uniqueness. Both lifting constants depend only on the fixed vessel.

For $G=D_j^s f$, the spectral expansion and \eqref{functional:negative} give $\|G\|_{(H^{1/2}(D))'}\lesssim \|P_j f\|_{\mathcal H_K^{s-1/2}}\leq\|f\|_{\mathcal H_K^{s-1/2}}$. The normalized lifting operator and $D_j^s$ are independent of time. Applying the weak estimate to $D_j^s\partial_t f$ gives the last bound.
\end{proof}

Set $G=D_j^s\partial_t^k\eta$, $M=J^{-1}\mathcal A_\eta^{-T}$, and take $w=M\nabla\psi$. To see that $M$ is an $H^1$ multiplier, use $M\in L^\infty$ and $\nabla M\in L^\infty_\theta L^{r_-}_{r,h}$ with $r_->2$ in the collar (Lemma~\ref{lem:tl_coef}). On each fixed transverse slice, $H^1\hookrightarrow L^{2r_-/(r_--2)}$ bounds $(\nabla M)\nabla\psi$ in $L^2$ after integration in $\theta$. Away from the edge use the interior coefficient bound $\nabla M\in L^3$ and $H^1\hookrightarrow L^6$. Hence $\|w\|_{H^1}\lesssim\|\psi\|_{H^2}$ under the smallness assumption. The Piola identity gives $J\operatorname{div}_{\mathcal A_\eta}w=\Delta\psi=0$ and $w\cdot\mathcal N=\nabla\psi\cdot\mathcal N_0=G$ on $\Sigma$. The fixed wall condition also holds: the vertical change of variables preserves the wall tangent plane in the collar and is the identity outside the support of the deformation. The finite spectral sum $G$ belongs to $H^1(D)$, so $w$ is an admissible test. We then have the following. (For simplicity, we write down the estimate for the highest order. For the lower-order estimates, there is always one temporal derivative fewer, so we can replace the dissipation used in the estimate by the energy functional.)
\begin{align}
&\left\langle \partial_t v_{m+1}, Jw \right\rangle +\left(\!\left(v_{m+1},w\right)\!\right) +\left(\eta_{m+1},w\cdot\mathcal{N}\right)_{1,\Sigma} +\kappa\left\langle\partial_t\eta_{m+1},w\cdot\mathcal{N}\right\rangle_{\Gamma_{0}}  \notag \\
&= \int_\Omega F_{0}^{1,m+1}\cdot wJ -\int_{\Sigma_s}JU_{0}(w)F_0^{6,m+1} -\int_{\Sigma} \left[ \sigma F_0^{5,m+1}\nabla (w\cdot\mathcal{N})+F_{0}^{3,m+1}\cdot w \right] -\left[w\cdot\mathcal{N},F_0^{7,m+1}\right]_{\Gamma_{0}}.
\end{align}

For the right-hand side of the equation, using the estimate for forcing terms derived in the previous theorem, we have
\begin{align}
&\int_\Omega F_{0}^{1,m+1}\cdot wJ -\int_{\Sigma_s}JU_{0}(w)F_0^{6,m+1} -\int_{\Sigma} \left[ F_{0}^{3,m+1}\cdot w \right]  \notag \\
&\lesssim \sqrt{\mathcal{E}\mathcal{D}} \|\omega\|_{H^{1}}\lesssim\sqrt{\mathcal{E}\mathcal{D}} \|\psi\|_{H^{2}} \lesssim\sqrt{\mathcal{E}\mathcal{D}} \|D_j^s\partial_t^k\eta\|_{\mathcal{H}_{\mathcal{K}}^{1/2}}.
\end{align}

For the term with $F_{0}^{5,m+1}$, we have
\begin{align}
\int_{\Sigma} \left[ \sigma F_0^{5,m+1}\nabla (w\cdot\mathcal{N}) \right]\lesssim \|F_{\mathrm{cl}+}^{5,m+1}\|_{H^{\frac{s}{2}}}\|D_{j}^{s}\eta_{m+1}\|_{H^{1-\frac{s}{2}}}+\|F_{\mathrm{int}+}^{5,m+1}\|_{H^{\frac{s}{2}}}\|D_{j}^{s}\eta_{m+1}\|_{H^{1-\frac{s}{2}}}.
\end{align}
For the first term on the right-hand side of the inequality, arguing as in \cite{GT}, we have
\begin{align}
&\|F_{\mathrm{cl}+}^{5,m+1}\|_{H^{\frac{s}{2}}}\|D_{j}^{s}\eta_{m+1}\|_{H^{1-\frac{s}{2}}} \lesssim \|\eta_{m+1}\|_{H^{1+\frac{s}{2}}}(\|F_{\mathrm{cl}+}^{5,m+1}\|_{L_{\theta}^{2}H_{r}^{\frac{s}{2}}}+\|F_{\mathrm{cl}+}^{5,m+1}\|_{H_{\theta}^{\frac{s}{2}}L_{r}^{2}}) \notag \\
&\lesssim \|\eta_{m+1}\|_{H^{1+\frac{s}{2}}}\Big(\|\eta_{m+1}\|_{H^{1+\frac{s}{2}}}\|\nabla \eta_{\mathrm{cl}+}\|_{L^{\infty}H^{\frac{1}{2}+}}+\sum_{j=1}^{\frac{(m+1)}{2}}\|\nabla \eta_{\mathrm{cl}+,j}\|_{L^{\infty}H^{\frac{1}{2}+\alpha}}\|\nabla \eta_{\mathrm{cl}+,m+1-j}\|_{L^{2}H^{\frac{1}{2}-\alpha}} \notag \\
&+\|\eta_{m+1}\|_{H^{1+\frac{s}{2}}}\|\nabla \eta_{\mathrm{cl}+}\|_{H^{1}L^{\infty}}+\|\eta_{m}\|_{H^{1+\frac{s}{2}}}\|\nabla \eta_{1,\mathrm{cl}+}\|_{H^{1}L^{\infty}} \notag \\
&+\sum_{j=2}^{\frac{(m+1)}{2}}(\|\partial_{\theta}\nabla\eta_{\mathrm{cl}+,j}\|_{L^{\infty}L^{\frac{1}{\alpha}}}+\|\nabla\eta_{\mathrm{cl}+,j}\|_{L^{\infty}L^{\frac{1}{\alpha}}})\|\nabla\eta_{\mathrm{cl}+,m+1-j}\|_{H^{\frac{s}{2}}L^{\frac{1}{\alpha}}}\Big) \notag \\
&\lesssim \|\eta_{m+1}\|_{H^{1+\frac{s}{2}}}\Big(\|\eta_{m+1}\|_{H^{1+\frac{s}{2}}} (\|\eta_{\mathrm{cl}+}\|_{L^{\infty}W^{3-\frac{1}{q_{\delta}}{q_{\delta}}}} \notag \\
&+\|\partial_{\theta}\eta_{\mathrm{cl}+}\|_{L^{\infty}H^{\frac{3}{2}-\alpha}})+\sum_{j=1}^{\frac{(m+1)}{2}}(\| \eta_{\mathrm{cl}+,j}\|_{L^{\infty}H^{\frac{3}{2}+\frac{\epsilon_{-}-\alpha}{2}}}+\| \p_{\theta}\eta_{\mathrm{cl}+,j}\|_{L^{\infty}H^{1}})\| \eta_{\mathrm{cl}+,m+1-j}\|_{H^{\frac{3}{2}-\alpha}} \notag \\
&+\|\eta_{m+1}\|_{H^{1+\frac{s}{2}}}(\|\eta_{\mathrm{cl}+}\|_{H^{1}H^{\frac{3}{2}+\frac{\varepsilon_{-}-\alpha}{2}}}+\|\partial_{\theta}\eta_{\mathrm{cl}+}\|_{H^{1}H^{\frac{3}{2}-\alpha}}) \notag \\
&+\|\eta_{m}\|_{H^{1+\frac{s}{2}}}(\|\partial_{t}\eta_{\mathrm{cl}+}\|_{H^{1}H^{\frac{3}{2}+\frac{\varepsilon_{-}-\alpha}{2}}}+\|\partial_{t}\p_{\theta}\eta_{\mathrm{cl}+}\|_{H^{1}H^{\frac{3}{2}-\alpha}}) \notag \\
&+\sum_{j=2}^{\frac{(m+1)}{2}}\Big((\|\eta_{j}\|_{H^{\frac{3}{2}-\alpha}}+\|\p_{\theta}\eta_{j}\|_{H^{\frac{3}{2}-\alpha}}+\|\p_{\theta}^{2}\eta_{j}\|_{H^{\frac{3}{2}-\alpha}})(\|\eta_{m+1-j}\|_{H^{\frac{3}{2}-\alpha}}+\|\p_{\theta}\eta_{m+1-j}\|_{H^{\frac{3}{2}-\alpha}})\Big) \notag \\
&\lesssim \|\eta_{m+1}\|_{H^{1+\frac{s}{2}}}^{2}\sqrt{\mathcal{E}}+\|\eta_{m+1}\|_{H^{1+\frac{s}{2}}}\sqrt{\mathcal{E}}\sqrt{\mathcal{D}}+\sqrt{\mathcal{E}}{\mathcal{D}},
\end{align}
where we used the fact that $m+1\geq 4$ and that for every $\rho\in[0,\frac12]$, on $\mathbb{R}^2$ or a bounded product domain
\begin{equation}
\|f\|_{H^\rho_{r,h\theta}} \asymp \left( \|f\|_{H_r^\rho L_\theta^2}^{2} + \|f\|_{L_r^2H_\theta^\rho}^{2} \right)^{1/2}.
\end{equation}
On $\mathbb{R}^2$, this follows from the Fourier weight equivalence  $(1+|\xi_r|^2+|\xi_\theta|^2)^\rho \asymp (1+|\xi_r|^2)^\rho+(1+|\xi_\theta|^2)^\rho$. On a bounded product domain, the same conclusion follows by applying bounded extension operators in each variable. We apply this with $\rho=s/2$. This implies the following norm separation $X=H^1_\theta H_r^{1/2+},\; \|uv\|_{H^{\frac{s}{2}}}\le C\|u\|_X\|v\|_{H^{\frac{s}{2}}}$.

For the second term, we argue similarly. All spatial norms below are taken over $\Sigma$. Recall that $F_{\mathrm{int}+}^{5,m+1} = \chi_{\mathrm{int}+}\partial_t^{m+1} \mathcal{R}(\nabla_\Sigma\eta)$. Separating the highest temporal derivative, we write
\begin{align}
F_{\mathrm{int}+}^{5,m+1} ={}& \chi_{\mathrm{int}+} D\mathcal{R}(\nabla_\Sigma\eta) \nabla_\Sigma\eta_{m+1} +\chi_{\mathrm{int}+}\mathcal{R}_{\mathrm{low}}^{m+1}, \\
\mathcal{R}_{\mathrm{low}}^{m+1} :={}& \sum_{\ell=2}^{m+1} \sum_{\substack{ i_1+\cdots+i_\ell=m+1 \\ i_q\ge1}} \frac{(m+1)!}{\ell!\,i_1!\cdots i_\ell!} D^\ell\mathcal{R}(\nabla_\Sigma\eta) \bigl[ \nabla_\Sigma\eta_{i_1},\ldots, \nabla_\Sigma\eta_{i_\ell} \bigr].
\end{align}
Every temporal index in the remainder satisfies $i_q\le m$. Since $\Sigma$ is two-dimensional and $p>\frac32$, the space $W^{2-\frac1p,p}(\Sigma)$ is an algebra and acts by multiplication on $H^{s/2}(\Sigma)$ in the above range. Hence
\begin{align}
&\left\| \chi_{\mathrm{int}+} D\mathcal{R}(\nabla_\Sigma\eta) \nabla_\Sigma\eta_{m+1} \right\|_{H^{s/2}} \lesssim \left\| \chi_{\mathrm{int}+} D\mathcal{R}(\nabla_\Sigma\eta) \right\|_{W^{2-\frac1p,p}} \|\chi_{0,\mathrm{int},\Sigma}\eta_{m+1}\|_{H^{1+\frac{s}{2}}} \notag \\
&\lesssim \|\chi_{0,\mathrm{int},\Sigma}\eta\|_{W^{3-\frac1p,p}} \|\chi_{0,\mathrm{int},\Sigma}\eta_{m+1}\|_{H^{1+\frac{s}{2}}}.
\end{align}

For the lower-order terms, the algebra and composition estimates yield
\begin{align}
\left\| \chi_{\mathrm{int}+}\mathcal{R}_{\mathrm{low}}^{m+1} \right\|_{H^{s/2}} \lesssim C_M \sum_{\ell=2}^{m+1} \sum_{\substack{ i_1+\cdots+i_\ell=m+1 \\ 1\le i_q\le m}} \prod_{q=1}^{\ell} \|\chi_{0,\mathrm{int},\Sigma}\eta_{i_q}\|_{W^{3-\frac1p,p}}.
\end{align}
Combining these bounds gives
\begin{align}
&\|F_{\mathrm{int}+}^{5,m+1}\|_{H^{s/2}} \lesssim  C_M\Bigg[ \|\chi_{0,\mathrm{int},\Sigma}\eta\|_{W^{3-\frac1p,p}} \|\chi_{0,\mathrm{int},\Sigma}\partial_t^{m+1}\eta\|_{H^{1+\frac{s}{2}}} + \sum_{\ell=2}^{m+1} \sum_{\substack{i_1+\cdots+i_\ell=m+1\\1\le i_q\le m}} \prod_{q=1}^{\ell} \|\chi_{0,\mathrm{int},\Sigma}\partial_t^{i_q}\eta\|_{W^{3-\frac1p,p}} \Bigg] \notag \\
&\lesssim \sqrt{\mathcal{E}}\|\chi_{0,\mathrm{int},\Sigma}\partial_t^{m+1}\eta\|_{H^{1+\frac{s}{2}}}+\sqrt{\mathcal{E}}\sqrt{D}.
\end{align}

We now estimate the terms on the left-hand side. For the first term, we have
\begin{align}
\langle \partial_{t}v_{m+1}, J\omega\rangle=&\partial_{t}\langle v_{m+1}, J\omega\rangle-\langle v_{m+1}, \partial_{t}J\omega\rangle-\langle v_{m+1}, J\partial_{t}\omega\rangle.
\end{align}
For the first term on the right-hand side of the equation above, we have
\begin{align}
\langle v_{m+1}, J\omega\rangle\lesssim\|v_{m+1}\|_{L^{2}}\|\omega\|_{L^{2}}\lesssim \sqrt{\mathcal{E}_{||}}\|\eta_{m+1}\|_{H^{s-\frac{1}{2}}}\lesssim \sqrt{\mathcal{E}_{||}}\|\eta_{m+1}\|_{H^{1}}\lesssim {\mathcal{E}}_{||}.
\end{align}
Similarly, for the second term, we have
\begin{align}
&\langle v_{m+1}, \partial_{t}J\omega\rangle\lesssim {\mathcal{E}}_{||} (\|\partial_{t}\bar{\eta}\|_{W^{1,+\infty}})\lesssim  {\mathcal{E}}_{||} (\|\chi_{0+,\mathrm{int}}\partial_{t}\bar{\eta}\|_{W^{1,+\infty}}+\|\chi_{\mathrm{cl}+}\partial_{t}\bar{\eta}\|_{L^{\infty}W^{1,+\infty}}+\|\chi_{\mathrm{cl}+}\partial_{t}\partial_{\theta}\bar{\eta}\|_{L^{\infty}L^{+\infty}}) \notag \\
&\lesssim {\mathcal{E}}_{||} (\|\chi_{0+,\mathrm{int}}\partial_{t}{\eta}\|_{W^{3-\frac{1}{p},p}}+\|\chi_{\mathrm{cl}+}\partial_{t}{\eta}\|_{L^{\infty}W^{3-\frac{1}{q_{\delta}},q_{\delta}}}+\|\chi_{\mathrm{cl}+}\partial_{t}\partial_{\theta}{\eta}\|_{L^{\infty}H^{\frac{3}{2}-\alpha}}) \lesssim  \mathcal{E}_{||}\sqrt{\mathcal{E}}.
\end{align}
For the third term, we have
\begin{equation}\label{eq:half-_eta}
\begin{aligned} \langle v_{m+1}, J\partial_{t}\omega\rangle&\lesssim \|v_{m+1}\|_{L^{2}}\|\eta_{m+2}\|_{H^{\frac{1}{2}-\alpha }}+\|v_{m+1}\|_{L^{2}}\|\p_{t}\bar{\eta}\|_{W^{1,+\infty}}\|\eta_{m+1}\|_{H^{\frac{1}{2}-\alpha}}\\
&\lesssim \|v_{m+1}\|_{L^{2}}\|\eta_{m+2}\|_{H^{\frac{1}{2}-\alpha }}+\mathcal{E}\|\eta_{m+1}\|_{H^{\frac{1}{2}-\alpha}}.
\end{aligned}
\end{equation}
By the kinematic boundary condition, we have
\begin{align}
&\|\eta_{m+2}\|_{H^{\frac{1}{2}-\alpha }}\lesssim \|v_{m+1}\|_{H^{1}}+\|F_{\mathrm{int}+}^{4,j}\|_{H^{\frac{1}{2}-\alpha}}+\|F_{\mathrm{cl}+}^{4,j}\|_{H^{\frac{1}{2}-\alpha}} \notag \\
&\lesssim \|v_{m+1}\|_{H^{1}}+\|F_{\mathrm{int}+}^{4,j}\|_{H^{\frac{1}{2}-\alpha}}+\|F_{\mathrm{cl}+}^{4,j}\|_{H^{\frac{1}{2}-\alpha}L^{2}}+\|F_{\mathrm{cl}+}^{4,j}\|_{L^{2}H^{\frac{1}{2}-\alpha}} \notag \\
&\lesssim  \sqrt{\mathcal{D}_{||}}+\|\chi_{0,\mathrm{int},\Sigma}v\|_{W^{2-\frac1p,p}(\Sigma)} \|\chi_{\mathrm{int}+,\Sigma}\eta_{m+1}\|_{H^{\frac32-\alpha}(\Sigma)} \notag \\
&+ \sum_{r=1}^{j-1}\binom{j}{r} \|\chi_{0,\mathrm{int},\Sigma}v_{m+1-r}\|_{W^{2-\frac1p,p}(\Sigma)} \|\chi_{\mathrm{int}+,\Sigma}\eta_r\|_{W^{3-\frac1p,p}(\Sigma)} \notag \\
&+(\|\chi_{\mathrm{cl}+,\Sigma}v\|_{L^{\infty}W^{2-\frac1{q_{\delta}},q_{\delta}}(\Sigma)}+\|\chi_{\mathrm{cl}+,\Sigma}\partial_{\theta}v\|_{L^{2}L^{\infty}(\Sigma)}) \|\chi_{\mathrm{bet}+,\Sigma}\eta_{m+1}\|_{H^{\frac32-\alpha}(\Sigma)} \notag \\
&+(\|\chi_{\mathrm{cl}+,\Sigma}v_{1}\|_{L^{\infty}W^{2-\frac1{q_{\delta}},q_{\delta}}(\Sigma)}+\|\chi_{\mathrm{cl}+,\Sigma}\partial_{\theta}v_{1}\|_{L^{2}L^{\infty}(\Sigma)}) \|\chi_{\mathrm{bet}+,\Sigma}\eta_{m}\|_{H^{\frac32-\alpha}(\Sigma)} \notag \\
&+ \sum_{r=1}^{m-1}\binom{j}{r} (\|\chi_{\mathrm{cl}+,\Sigma}v_{m+1-r}\|_{L^{2}W^{2-\frac1{q_{\delta}},q_{\delta}}(\Sigma)}+\|\chi_{\mathrm{cl}+,\Sigma}\partial_{\theta}v_{m+1-r}\|_{L^{2}L^{2}(\Sigma)}) \notag \\
&\times(\|\chi_{\mathrm{bet}+,\Sigma}\eta_{r}\|_{L^{\infty}H^{\frac32-\alpha}(\Sigma)}+\|\chi_{\mathrm{bet}+,\Sigma}\eta_{r}\|_{L^{\infty}W^{3-\frac{1}{q_{\delta}},,q_{\delta}}(\Sigma)}) \notag \\
&\lesssim  \sqrt{\mathcal{D}_{||}}+\|\chi_{0,\mathrm{int},\Sigma}v\|_{W^{2,p}(\Omega)} \|\chi_{0,\mathrm{int},\Sigma}\eta_{m+1}\|_{H^{\frac32-\alpha}(\Sigma)} + \sum_{r=1}^{j-1}\binom{j}{r} \|\chi_{0,\mathrm{int}}v_{m+1-r}\|_{W^{2,p}(\Omega)} \|\chi_{0,\mathrm{int}}\eta_r\|_{W^{3-\frac1p,p}(\Sigma)} \notag \\
&+(\|\chi_{\mathrm{bet}+}v\|_{L^{\infty}W^{2,q_{\delta}}(\Omega)}+\|\chi_{\mathrm{bet}+}\partial_{\theta}v\|_{L^{2}W^{2,q_{\delta}}(\Omega)}) \|\chi_{\mathrm{bet}+,\Sigma}\eta_{m+1}\|_{H^{\frac32-\alpha}(\Sigma)} \notag \\
&+(\|\chi_{\mathrm{bet}+}v_{1}\|_{L^{\infty}W^{2,q_{\delta}}(\Omega )}+\|\chi_{\mathrm{bet}+}\partial_{\theta}v_{1}\|_{L^{2}W^{2,q_{\delta}}(\Omega)}) \|\chi_{\mathrm{bet}+,\Sigma}\eta_{m}\|_{H^{\frac32-\alpha}(\Sigma)} \notag \\
&+ \sum_{r=1}^{m-1}\binom{j}{r} (\|\chi_{\mathrm{bet}+,\Sigma}v_{m+1-r}\|_{L^{2}W^{2,q_{\delta}}(\Omega)}+\|\chi_{\mathrm{bet}+,\Sigma}\p_{\theta}v_{m+1-r}\|_{H^{1}(\Omega)}) \notag \\
&\times(\|\chi_{\mathrm{bet}+,\Sigma}\eta_{r}\|_{L^{\infty}H^{\frac32-\alpha}(\Sigma)}+\|\chi_{\mathrm{bet}+,\Sigma}\eta_{r}\|_{L^{\infty}W^{3-\frac{1}{q_{\delta}},,q_{\delta}}(\Sigma)}) \lesssim \sqrt{\mathcal{D}_{||}}+\sqrt{\mathcal{E}\mathcal{D}}.
\end{align}

For the second term, we have
\begin{align}
|((v_{m+1},\omega))|\lesssim \|v_{m+1}\|_{H^{1}}\|\omega\|_{H^{1}}\lesssim \sqrt{\mathcal{D}_{||}}\|\partial_t^{m+1}\eta\|_{H^{1+\frac{s}{2}}}\lesssim \mathcal{D}_{||}+\varepsilon\|\partial_{t}^{m+1}\eta\|_{H^{1+\frac{s}{2}}}^{2}.
\end{align}

For the third term, the spectral identity above gives $(\eta_{m+1},D_j^s\eta_{m+1})_{1,\Sigma} =\|P_j\eta_{m+1}\|_{\mathcal H_K^{1+s/2}}^2$. As $j\to\infty$, this increases to $\|\eta_{m+1}\|_{\mathcal H_K^{1+s/2}}^2$, which is equivalent to $\|\eta_{m+1}\|_{H^{1+s/2}(D)}^2$. We take $s=1-2\alpha$, so that the mixed-norm exponent is $s/2=1/2-\alpha$.

For the fourth term, combining it with the last term on the right-hand side, we obtain
\begin{align}
[\kappa\partial_{t}\eta_{m+1}+F_{0}^{7,j},\omega \cdot \mathcal{N}]_{\Gamma_{0}}\lesssim ([\partial_{t}\eta_{m+1}]_{\Gamma_{0}}+[F_{0}^{7,j}]_{\Gamma_{0}})\|\partial_t^{m+1}\eta\|_{H^{1+\frac{s}{2}}}\lesssim \mathcal{D}_{||}+\varepsilon\|\partial_{t}^{m+1}\eta\|_{H^{1+\frac{s}{2}}}^{2}+\mathcal{E}\mathcal{D}.
\end{align}

Combining all of the estimates above, we obtain
\begin{align}
\int_{0}^{t}(\|\partial_{t}^{m+1}\eta\|_{H^{\frac{3}{2}-\alpha}}^{2}+\|\partial_{t}^{m+2}\eta\|_{H^{\frac{1}{2}-\alpha}}^{2})\lesssim \mathcal{E}_{||}(t)+\mathcal{E}_{||}(0)+\int_{0}^{t}\mathcal{D}_{||}+\int_{0}^{t}\sqrt{\mathcal{E}}\mathcal{D},  \label{eq:threehalf_top}
\end{align}
where the terms bounded above by $\sqrt{\mathcal{D}_{||}}\,\mathcal{E}_{||}$ are $\lesssim\mathcal{E}\sqrt{\mathcal{D}}\lesssim\sqrt{\mathcal{E}}\,\mathcal{D}$ by Lemma~\ref{lem:E_le_D}, and, at each time,
\begin{align}
\sum_{j=0}^{m}(\|\partial_{t}^{j}\eta\|_{H^{\frac{3}{2}-\alpha}}^{2}+\|\partial_{t}^{j+1}\eta\|_{H^{\frac{1}{2}-\alpha}}^{2})\lesssim \mathcal{E}_{||}+\sum_{j=0}^{m}(\|v_j\|_{H^1}^2 +\|\eta_{j+1}\|_{L^2(\Gamma_0)}^2)+\mathcal{E}^{\frac{3}{2}}.  \label{eq:threehalf_low}
\end{align}
The same argument gives the corresponding inequality for each fixed $0\leq j\leq m$.

Below the top level, $j\leq m$, no integration by parts in time is needed. Indeed $\langle\partial_{t}v_{j},Jw\rangle\leq\|\partial_{t}^{j+1}v\|_{L^{2}}\|\psi\|_{H^{1}}$ with $j+1\leq m+1$, which is part of $\mathcal{D}_{||}$, and the contact-line term $\kappa[\partial_{t}\eta_{j},w\cdot\mathcal{N}]_{\Gamma}$ costs only $H^{\frac{3}{2}-2\alpha+\delta^{\prime }}$ with $0<\delta^{\prime }<\alpha$. Hence $\|\partial_{t}^{j}\eta\|_{H^{\frac{3}{2}-\alpha}}^{2}\lesssim\mathcal{D}_{||}+\mathcal{E}\mathcal{D}$ at each time, and this is integrated in time:
\begin{align}
\sum_{j=0}^{m}\int_{0}^{t}(\|\partial_{t}^{j}\eta\|_{H^{\frac{3}{2}-\alpha}}^{2}+\|\partial_{t}^{j+1}\eta\|_{H^{\frac{1}{2}-\alpha}}^{2})\lesssim\int_{0}^{t}\mathcal{D}_{||}+\int_{0}^{t}\sqrt{\mathcal{E}}\mathcal{D}.  \label{eq:threehalf_low_int}
\end{align}
The level $j=0$ raises nothing new, since $\int_{D}\eta=0$ by Subsection~\ref{subsec:volume}.

\subsection{Enhanced energy estimate for pressure}

\label{subsec:enhanced_pressure}

In this subsection, we estimate $\|p_{m}\|_{L^{2}}$ by the energy functional $\mathcal{E}$.

\begin{theorem}[Energetic enhancement for $\partial_t p$]
\label{thm:pressure_enhancement} Suppose $0<T\le\infty$ and  $\sup_{0\le t<T}\mathcal{E}(t)\le\gamma^2$, with $\gamma\in(0,1)$ chosen as in Proposition~\ref{prop:B_d}. At every $t\in[0,T)$ and for each $0\le j\le m$, \begin{equation} \|p_j\|_{L^2}\lesssim \|v_j\|_{H^1}+\|\partial_t v_j\|_{L^2} +\|\eta_j\|_{H^1}+\mathcal E+\mathcal E^{3/2}.
\end{equation}
\end{theorem}

We give the proof for $j=m$; the same argument applies at each lower order.

\begin{remark}
The pressure is estimated through the weak formulation and the Ne\caronaccent{c}as inequality on $\Omega$, not through $H^{2}$ elliptic regularity, which is not available in the wedge along $\Gamma$. No cutoff is used, so no derivative of a cutoff of width $\varepsilon$ enters.
\end{remark}

\begin{proof}
Apply \eqref{eq:whole-domain-weak-momentum-j} at $j=m$. Both tests used below satisfy $w\cdot n_0=0$ on $\Sigma_s$ and have $w\cdot\mathcal N\in H^1(\Sigma)$ vanishing near $\Gamma_0$. Solving for the pressure pairing on the reference vessel gives
\begin{align}
\label{eq:pressure_estimate}
\int_\Omega Jp_m\func{div}_{\mathcal A_\eta}w ={}&\langle\partial_t v_m,Jw\rangle +\left(\!\left(v_m,w\right)\!\right)-\int_\Omega JF_0^{1,m}\cdot w +\int_{\Sigma_s}J_sF_0^{6,m}\cdot U_0w+\int_\Sigma F_0^{3,m}\cdot w \notag\\
&+\int_\Sigma \bigl[g\eta_m-\sigma\func{div}_\Sigma(a_\zeta\nabla_\Sigma\eta_m+F_0^{5,m})\bigr] (w\cdot\mathcal N).
\end{align}
Using the previous energy estimate
\begin{align}
&\left|\int_\Omega JF_0^{1,m}\cdot w -\int_{\Sigma_s}J_sF_0^{6,m}\cdot U_0w-\int_\Sigma F_0^{3,m}\cdot w \right|\lesssim(\mathcal E+\mathcal E^{3/2})\|w\|_{H^1}.
\end{align}
On the other hand,  $\left| \left(\!\left(v_{m},w\right)\!\right) \right| \lesssim \|v_{m}\|_{H^1}\|w\|_{H^1}$, and $\left|\left\langle\partial_tv_{m},Jw\right\rangle\right| \lesssim \|\partial_t v_{m}\|_{L^2}\|w\|_{L^2}$.

For $w\in H^1_0(\Omega;\mathbb R^3)$ all boundary terms vanish. Put $\overline{p_m}=\int_\Omega Jp_m/\int_\Omega J$, where $\int_\Omega J=|\Omega(t)|=|\Omega|$ by volume conservation. Let $z=\mathcal B(J(p_m-\overline{p_m}))$, with $\mathcal B$ the bounded right inverse of the ordinary divergence on mean-zero data in $\Omega$. As in Proposition~\ref{prop:B_d}, take $w=Mz$ with $M=J^{-1}\mathcal A_\eta^{-T}$. The $H^1$ multiplier bound for $M$ gives $\|w\|_{H^1}\lesssim\|p_m-\overline{p_m}\|_{L^2}$, and the Piola identity gives
\begin{equation}
J\func{div}_{\mathcal A_\eta}w=\func{div}z=J(p_m-\overline{p_m}), \qquad \int_\Omega Jp_m\func{div}_{\mathcal A_\eta}w =\int_\Omega J|p_m-\overline{p_m}|^2 .
\end{equation}
Apply \eqref{eq:pressure_estimate} to this $w$. Since $J$ is bounded above and below by positive constants, the preceding bounds yield
\begin{equation}\label{est:p_1} \|p_m-\overline{p_m}\|_{L^2} \lesssim\|v_m\|_{H^1}+\|\partial_t v_m\|_{L^2} +\mathcal E+\mathcal E^{3/2}.
\end{equation}

It now suffices to estimate the average. To estimate the average of $p$, choose $\varphi\in C_{c}^{\infty}(\Sigma)$, supported in the smooth patch away from the contact line, such that
\begin{equation*}
\int_{\Sigma}\frac{\varphi}{|\mathcal{N}_{0}|}\,dS=|\Omega|=|\Omega(t)|.
\end{equation*}
By the trace extension theorem, there exists $w_{0}\in H^{1}(\Omega;\mathbb{R}^{3})$ satisfying  $w_{0}=0 \;\text{on }\Sigma_{s}, \; w_{0}=\frac{\varphi \mathcal{N}}{|\mathcal{N}|^{2}} \;\text{on }\Sigma$. Denote the push-forward of $\omega_{0}$ to $\Omega(t)$ by $\tilde{\omega_{0}}$. Then the divergence theorem yields
\begin{equation*}
\int_{\Omega(t)}\func{div}\tilde{w}_{0}\,dx =\int_{\Omega}J\func{div}_{\mathcal{A}}{w}_{0}\,dx= \int_{\Sigma}{w}_{0}\cdot J\mathcal{A}_{\eta}\nu \,dS = \int_{\Sigma}\frac{\varphi}{|\mathcal{N}_{0}|}\,dS = |\Omega|,
\end{equation*}
where we used the fact that $J\mathcal{A}_{\eta}\nu=J\mathcal{A}_{\eta}\frac{\mathcal{N}_{0}}{|\mathcal{N}_{0}|}=\frac{\mathcal{N}}{|\mathcal{N}_{0}|}$ on $\Sigma$. Consequently,
\begin{equation*}
\int_{\Omega(t)}\bigl(\func{div}\tilde{w}_{0}-1\bigr)\,dx=\int_{\Omega}J\bigl(\func{div}_{\mathcal{A}_{\eta}}w_{0}-1\bigr)\,dx=0.
\end{equation*}
Let $M=J^{-1}\mathcal A_\eta^{-T}$, and let $\mathcal B$ be a fixed right inverse of the ordinary divergence on zero-average data. Then there exists $f\in H_{0}^{1}(\Omega;\mathbb{R}^{3})$ such that  $G=J\operatorname{div}_{\mathcal A_\eta}w_0-J, \; f=M\mathcal B(G),\; w_*=w_0-f$. The normalization of the top flux and conservation of volume give $\int_\Omega G=0$. The Piola identity therefore yields
\begin{equation}
\operatorname{div}_{\mathcal A_\eta}w_*=1, \qquad w_*|_{\partial\Omega}=w_0|_{\partial\Omega}, \qquad \int_\Omega Jp_m\operatorname{div}_{\mathcal A_\eta}w_* =\int_{\Omega}Jp_{m}=|\Omega|\,\overline{p_m}.
\end{equation}
Here the cutoff $\varphi$ and its extension are fixed independently of $p$. Moreover, $w_*\cdot\mathcal N=\varphi$ on $\Sigma$, and $\varphi$ vanishes near $\Gamma_0$. Thus the contact-line term vanishes for this test.

For the right-hand side of equation \eqref{eq:pressure_estimate}, with $w=w_*$, we have
\begin{equation}
{}  \left\langle \partial_t v_{m},Jw\right\rangle +\left(\!\left(v_{m},w\right)\!\right) - \int_\Omega F_{0}^{1,m}\cdot wJ +\int_{\Sigma_s}J_sU_{0}w\cdot F_{0}^{6,m}+\int_\Sigma F_0^{3,m}\cdot w \lesssim \|\partial_{t}v_{m}\|_{L^{2}}+\|v_{m}\|_{H^{1}}+\mathcal{E}+\mathcal{E}^{\frac{3}{2}}. \label{eq:pressure_hard_1}
\end{equation}

For the remaining surface term use $w=w_*$, so that $w\cdot\mathcal N=\varphi$. Integration by parts on the surface gives
\begin{align}
\label{eq:pressure_hard}
&\int_\Sigma \bigl[g\eta_m-\sigma\func{div}_\Sigma(a_\zeta\nabla_\Sigma\eta_m+F_0^{5,m})\bigr]\varphi =g\int_\Sigma\eta_m\varphi+ \sigma\int_\Sigma(a_\zeta\nabla_\Sigma\eta_m+F_0^{5,m})\cdot\nabla_\Sigma\varphi .
\end{align}
There is no curve integral since $\varphi$ is supported away from $\Gamma_0$. The fixed $H^1$ norm of $\varphi$ therefore gives
\begin{equation}\label{eq:pressure_hard_final} \left|g\int_\Sigma\eta_m\varphi+ \sigma\int_\Sigma(a_\zeta\nabla_\Sigma\eta_m+F_0^{5,m})\cdot\nabla_\Sigma\varphi\right| \lesssim\|\eta_m\|_{H^1}+\|F_0^{5,m}\|_{L^2} \lesssim\|\eta_m\|_{H^1}+\mathcal E .
\end{equation}

Combining equation \eqref{eq:pressure_hard_1} and equation \eqref{eq:pressure_hard_final}, we obtain
\begin{align}
{\label{est:p_2}} \|\overline{p_m}\|_{L^2} \lesssim \|v_m\|_{H^1}+\|\partial_t v_m\|_{L^2} +\|\eta_m\|_{H^1}+\mathcal E+\mathcal E^{3/2}.
\end{align}
Combining \eqref{est:p_1} and \eqref{est:p_2} by the triangle inequality gives
\begin{equation}
\|p_m\|_{L^2}\lesssim \|v_m\|_{H^1}+\|\partial_t v_m\|_{L^2} +\|\eta_m\|_{H^1}+\mathcal E+\mathcal E^{3/2},
\end{equation}
which finishes the proof.
\end{proof}

\begin{lemma}[The energy is controlled by the dissipation]
\label{lem:E_le_D} There is $C_{*}$ such that $\mathcal{E}(t)\leq C_{*}\mathcal{D}(t)$ at every time $t$ at which $\mathcal{E}(t)$ is sufficiently small.
\end{lemma}

\begin{proof}
Every term of $\mathcal{E}$ other than the pressure terms $\|\partial_{t}^{j}p\|_{L^{2}(\Omega)}^{2}$, $j\leq m$, of $\mathcal{E}_{\mathrm{enhanced}}$ is bounded by a term of $\mathcal{D}$:

\begin{itemize}
\item the terms of $\mathcal{E}_{||}$ and the velocity and surface terms of $\mathcal{E}_{\mathrm{enhanced}}$ by $\mathcal{D}_{||}$ and $\mathcal{D}_{\mathrm{enhanced}}$;

\item the terms of $\mathcal{E}_{\mathrm{int}}$ by those of $\mathcal{D}_{\mathrm{int}}$, which carry one more derivative, using $p\geq6/(5-2\varepsilon_{-})$ for the fractional terms;

\item the terms of $\mathcal{E}_{\mathrm{cl}}$ and $\mathcal{E}_{\mathrm{bet}}$ by those of $\mathcal{D}_{\mathrm{cl}}$ and $\mathcal{D}_{\mathrm{bet}}$ at the same or the next level, using $W^{2,q_{\delta}}\subset H^{1+\delta}$ on the slices and $W^{3-1/q_{\delta},q_{\delta}}(\mathbb{R})\subset H^{3/2+\delta}(\mathbb{R})$.
\end{itemize}

The pressure terms also have a direct bound. Write $\chi_+=\chi_{\mathrm{cl}}+\chi_{\mathrm{bet}}$, so that $1=\chi_{\mathrm{int}}+\chi_+$. For $0\leq j\leq m$, the cutoff nesting and $W^{1,p}(\Omega)\subset L^2(\Omega)$ give
\begin{equation}
\|\partial_t^jp\|_{L^2} \leq\|\chi_{\mathrm{int}}\partial_t^jp\|_{L^2} +\|\chi_+\partial_t^jp\|_{L^2} \leq C(\varepsilon)\|\chi_{m,\mathrm{int}}\partial_t^jp\|_{W^{1,p}} +\|\chi_+\partial_t^jp\|_{L^2}.
\end{equation}
The first term occurs in the $k=m$ row of \eqref{def:interior_dissipation}. The second occurs in the $(k_1,k_2)=(m+1,0)$ energy block included in $\mathcal D_{\mathrm{cl}}$. Therefore
\begin{equation}
\sum_{j=0}^{m}\|\partial_t^jp\|_{L^2}^2 \leq C(\varepsilon)\mathcal D.
\end{equation}
Together with the preceding bounds, this gives $\mathcal E\leq C_*\mathcal D$ for a constant $C_*$ depending on the fixed cutoffs.
\end{proof}

The lemma holds at each fixed time. Its proof uses only the definitions of $\mathcal E$ and $\mathcal D$, the cutoff partition and the stated Sobolev embeddings. It uses neither the pressure-enhancement theorem nor an energy inequality, so the lemma may be used in the proofs of Theorems~\ref{thm:whole-space-energy} and~\ref{thm:whole_space_estimate_final}, as above.

We now prove Theorem~\ref{thm:whole_space_estimate_final}.

\begin{proof}[Proof of Theorem~\protect\ref{thm:whole_space_estimate_final}]

Recall that $\mathcal{E}_{\mathrm{enhanced}}$ consists of $\|\partial_{t}^{j}\eta\|_{H^{3/2-\alpha}}^{2}$, $\|\partial_{t}^{j}p\|_{L^{2}}^{2}$, $\|\partial_{t}^{j}v\|_{H^{1}}^{2}$, $\|\partial_{t}^{j+1}\eta\|_{L^{2}(\Gamma_{0})}^{2}$\newline and $\|\partial_{t}^{j+1}\eta\|_{H^{1/2-\alpha}}^{2}$, $j\leq m$, and that
\begin{equation*}
\mathcal{D}_{\mathrm{enhanced}}=\sum_{j\leq m+1}\|\partial_{t}^{j}\eta\|_{H^{3/2-\alpha}}^{2}
\end{equation*}
. We use $\mathcal{E}\leq1$ throughout, so that $\mathcal{E}^{2}\leq\mathcal{E}^{3/2}$ and $\mathcal{E}\mathcal{D}\leq\sqrt{\mathcal{E}}\mathcal{D}$.

\emph{The dissipation.} By \eqref{eq:threehalf_top} and \eqref{eq:threehalf_low_int},
\begin{equation*}
\int_{0}^{t}\mathcal{D}_{\mathrm{enhanced}}\lesssim\mathcal{E}_{||}(t)+\mathcal{E}_{||}(0)+\int_{0}^{t}\mathcal{D}_{||}+\int_{0}^{t}\sqrt{\mathcal{E}}\mathcal{D}.
\end{equation*}

\emph{The energy at time $t$.} We take $j\leq m$.

\begin{itemize}
\item By \eqref{eq:threehalf_low}, $\|\partial_{t}^{j}\eta\|_{H^{3/2-\alpha}}^{2}(t)\lesssim\mathcal{E}_{||}(t)+\|v_j\|_{H^1}^2(t) +\|\eta_{j+1}\|_{L^2(\Gamma_0)}^2(t)+\mathcal{E}^{3/2}(t)$.

\item By the fundamental theorem of calculus and the Cauchy--Schwarz inequality, $\|\partial_{t}^{j}v\|_{H^{1}}^{2}(t)\leq\|\partial_{t}^{j}v\|_{H^{1}}^{2}(0)+\int_{0}^{t}\big(\|\partial_{t}^{j}v\|_{H^{1}}^{2}+\|\partial_{t}^{j+1}v\|_{H^{1}}^{2}\big)\leq\mathcal{E}(0)+2\int_{0}^{t}\mathcal{D}_{||}$, since $j+1\leq m+1$. In the same way, $\|\partial_{t}^{j+1}\eta\|_{L^{2}(\Gamma_{0})}^{2}(t)\leq\mathcal{E}(0)+2\int_{0}^{t}\mathcal{D}_{||}$, since $\mathcal{D}_{||}$ contains $\|\partial_{t}^{i+1}\eta\|_{L^{2}(\Gamma_{0})}^{2}$ for $i\leq m+1$.

\item By the kinematic equation, $\partial_{t}^{j+1}\eta=v_{j}\cdot\mathcal{N}+F_{0}^{4,j}$. The bound following \eqref{eq:half-_eta}, with energy norms in place of dissipation norms since $j\leq m$, gives $\|\partial_{t}^{j+1}\eta\|_{H^{1/2-\alpha}}\lesssim\|v_{j}\|_{H^{1}}+\mathcal{E}$.

\item By Theorem~\ref{thm:pressure_enhancement} and the sentence after it, $\|\partial_{t}^{j}p\|_{L^{2}}\lesssim\|v_{j}\|_{H^{1}}+\|v_{j+1}\|_{L^{2}}+\|\eta_{j}\|_{H^{1}}+\mathcal{E}$, and $\|v_{j+1}\|_{L^{2}}^{2}+\|\eta_{j}\|_{H^{1}}^{2}\lesssim\mathcal{E}_{||}$.
\end{itemize}

Hence $\mathcal{E}_{\mathrm{enhanced}}(t)\lesssim\mathcal{E}(0)+\mathcal{E}_{||}(t)+\mathcal{E}^{3/2}(t)+\int_{0}^{t}\mathcal{D}_{||}$.

\emph{Conclusion.} Adding the two bounds and using Theorem~\ref{thm:whole-space-energy} for $\mathcal{E}_{||}(t)+\int_{0}^{t}\mathcal{D}_{||}$,
\begin{equation*}
\mathcal{E}_{\mathrm{enhanced}}(t)+\int_{0}^{t}\mathcal{D}_{\mathrm{enhanced}}\lesssim\mathcal{E}(0)+\mathcal{E}^{\frac{3}{2}}(t)+\mathcal{E}^{\frac{3}{2}}(0)+\int_{0}^{t}\sqrt{\mathcal{E}}\mathcal{D},
\end{equation*}
which is the claim.
\end{proof}

\section{Interior elliptic estimate}

\label{sec:interior}

\subsection{Interior elliptic system}

Multiplying both sides of equation \eqref{eq:whole-domain-elliptic-system-j} by $\chi_{j,\mathrm{int}}$ for $0\leq j\leq m$, we obtain
\begin{equation}  \label{eq:interior-domain-elliptic-system-j} \left\{ \begin{alignedat}{2} &\partial_tv_{j,\mathrm{int}} + \operatorname{div}S(p_{j,\mathrm{int}},v_{j,\mathrm{int}}) = F_{\mathrm{int}}^{1,j}+F_{s,\mathrm{int}}^{1,j} &&\qquad \text{in } \Omega, \\ &\operatorname{div}v_{j,\mathrm{int}} = F_{\mathrm{int}}^{2,j}+F_{\mathrm{s,int}}^{2,j} &&\qquad \text{in } \Omega, \\ &S(p_{j,\mathrm{int}},v_{j,\mathrm{int}})\frac{1}{|\mathcal{N}_{0}|}\mathcal{N}_{0} = \frac{1}{|\mathcal{N}_{0}|}\left( g\eta_{j,\mathrm{int}} - \sigma\operatorname{div}_{\Sigma}\mathcal L_{\Sigma,0}\eta_{j,\mathrm{int}} -\sigma \operatorname{div}_{\Sigma}F_{\mathrm{int}}^{5,j}\right)\mathcal{N}_{0} +\frac{1}{|\mathcal{N}_{0}|} F_{\mathrm{int}}^{3,j}+\frac{1}{|\mathcal{N}_{0}|}F_{s,\mathrm{int}}^{3,j} &&\qquad \text{on }\Sigma, \\ &\frac{1}{|\mathcal{N}_{0}|}\partial_t\eta_{j,\mathrm{int}}-\frac{1}{|\mathcal{N}_{0}|}v_{j,\mathrm{int}}\cdot \mathcal{N}_{0} = \frac{1}{|\mathcal{N}_{0}|}F_{\mathrm{int}}^{4,j}+\frac{1}{|\mathcal{N}_{0}|}F_{s,\mathrm{int}}^{4,j} &&\qquad \text{on }\Sigma,\\
&U_{0}S(p_{j,\mathrm{int}},v_{j,\mathrm{int}})n_{0}-\beta_{0}U_{0}v_{j,\mathrm{int}} = F_{\mathrm{int}}^{6,j}+F_{s,\mathrm{int}}^{6,j} &&\qquad \text{on }\Sigma_{s}, \\ &v_{j,\mathrm{int}}\cdot n_{0} = 0 &&\qquad \text{on }\Sigma_{s}, \\ &v_{\mathrm{int},j} = 0 &&\qquad \text{near }\Sigma_{\mathrm{art}},\\ &\eta_{\mathrm{int},j}=0&&\qquad \text{near }\Gamma_{\mathrm{art}}, \end{alignedat} \right.
\end{equation}
where
\begin{align}
(v_{j,\mathrm{int}},p_{j,\mathrm{int}},\eta_{j,\mathrm{int}}):=(\chi_{j,\mathrm{int}}\partial_{t}^{j}v,\chi_{j,\mathrm{int}}\partial_{t}^{j}p,\chi_{j,\mathrm{int}}\partial_{t}^{j}\eta),
\end{align}
and
\begin{equation}
\begin{aligned}
&F_{\mathrm{int}}^{1,j}:= \nabla\chi_{j,\mathrm{int}}\cdot S(p_{j},v_{j})-\mu \operatorname{div}(v_{j}\otimes \nabla\chi_{j,\mathrm{int}}+\nabla\chi_{j,\mathrm{int}}\otimes v_{j})+\chi_{j,\mathrm{int}}F_{0}^{1,j},\quad F_{\mathrm{int}}^{2,j}:= \nabla\chi_{j,\mathrm{int}}\cdot v_{j}+\chi_{j,\mathrm{int}}F_{0}^{2,j},\\
&F_{\mathrm{int}}^{3,j}:= -\mu (v_{j}\otimes \nabla\chi_{j,\mathrm{int}}+\nabla\chi_{j,\mathrm{int}}\otimes v_{j})\mathcal{N}_{0}+\sigma [\operatorname{div_{\Sigma}}\mathcal{L}_{\Sigma,0},\chi_{j,\mathrm{int},\Sigma}]\eta_{j}\mathcal{N}_{0}+\chi_{j,\mathrm{int},\Sigma}F_{0}^{3,j}+\sigma(\nabla_{\Sigma }\chi_{j,\mathrm{int},\Sigma}\cdot F_{0}^{5,j})\mathcal{N}_{0},\\
&F_{\mathrm{int}}^{4,j}:= \chi_{j,\mathrm{int},\Sigma}F_{0}^{4,j},\quad F_{\mathrm{int}}^{5,j}:= \chi_{j,\mathrm{int},\Sigma}F_{0}^{5,j},\quad F_{\mathrm{int}}^{6,j}:= -\mu U_{0}(v_{j}\otimes \nabla\chi_{j,\mathrm{int}}+\nabla\chi_{j,\mathrm{int}}\otimes v_{j})n_{0}+\chi_{j,\mathrm{int}}F_{0}^{6,j},
\end{aligned}
\end{equation}
and
\begin{equation}
\begin{aligned}
&F_{s,\mathrm{int}}^{1,j}:= -\chi_{j,\mathrm{int}}\operatorname{div}_{\mathcal{A}_{\eta}-I}S_{\mathcal{A}_{\eta}}(\p_{t}^{j}p,\p_{t}^{j}v)+\mu\chi_{j,\mathrm{int}}\operatorname{div}\mathbb{D}_{\mathcal{A}_{\eta}-I}(v_{j}),\quad F_{s,\mathrm{int}}^{2,j}:= -\chi_{j,\mathrm{int}}\operatorname{div}_{\mathcal{A}_{\eta}-I}v_{j},\\
&F_{s,\mathrm{int}}^{3,j}:= -\mu\chi_{j,\mathrm{int},\Sigma}\mathbb{D}_{I-\mathcal{A}_{\eta}}(v_{j})\mathcal{N}+\chi_{j,\mathrm{int},\Sigma}S(p_{j},v_{j})(\mathcal{N}_{0}-\mathcal{N})\\
&+\chi_{\mathrm{j},\mathrm{int},\Sigma}\left( g\eta_{j} - \sigma\operatorname{div}_{\Sigma}\mathcal L_{\Sigma,0}\eta_{j} -\sigma \operatorname{div}_{\Sigma} F_{0}^{5,j}\right)(\mathcal{N}-\mathcal{N}_{0}),\\
&F_{s,\mathrm{int}}^{4,j}:= \chi_{j,\mathrm{int},\Sigma}v_{j}(\mathcal{N}-\mathcal{N}_{0}),\quad F_{s,\mathrm{int}}^{6,j}:= -\mu\chi_{j,\mathrm{int}}U_{0}\mathbb{D}_{I-\mathcal{A}_{\eta}}(v_{j})n_{0}.
\end{aligned}
\end{equation}

\begin{lemma}[Interior elliptic estimate]
\label{lem:interior_elliptic} Let $1<p<\frac{2}{1+2\alpha}$ and $0\leq j\leq m$. Let $\mathcal{U}$ be the neighbourhood of $\Gamma$ of Lemma~\ref{lem:tl_coef}, and $r_{-}=\frac{2}{1-\varepsilon_{-}}$. Assume that, on $\func{supp}\chi_{j,\mathrm{int}}$, the steady state $\zeta$ is smooth. Then every solution of \eqref{eq:interior-domain-elliptic-system-j} satisfies
\begin{equation}\label{eq:interior_elliptic}
\begin{aligned}
&\|v_{j,\mathrm{int}}\|_{W^{2,p}} +\|\eta_{j,\mathrm{int}}\|_{W^{3-\frac{1}{p},p}} +\|p_{j,\mathrm{int}}\|_{W^{1,p}} \leq C(\varepsilon)\Big(\|F_{\mathrm{int}}^{1,j}\|_{L^{p}} +\|F_{\mathrm{int}}^{2,j}\|_{W^{1,p}} +\|F_{\mathrm{int}}^{3,j}\|_{W^{1-\frac{1}{p},p}}\\
&+\|F_{\mathrm{int}}^{4,j}\|_{W^{2-\frac{1}{p},p}} +\|F_{\mathrm{int}}^{5,j}\|_{W^{2-\frac{1}{p},p}} +\|\partial_{t}\eta_{j,\mathrm{int}}\|_{H^{\frac{3}{2}-\alpha}} +\|\partial_{t}v_{j,\mathrm{int}}\|_{L^{2}} +\|v_{j,\mathrm{int}}\|_{H^{1}} +\|\eta_{j,\mathrm{int}}\|_{H^{\frac{3}{2}-\alpha}} +\|p_{j,\mathrm{int}}\|_{L^{2}}\\
&+\|F_{\mathrm{int}}^{6,j}\|_{W^{1-\frac{1}{p},p}} +\|F_{s,\mathrm{int}}^{1,j}\|_{L^{p}} +\|F_{s,\mathrm{int}}^{2,j}\|_{W^{1,p}} +\|F_{s,\mathrm{int}}^{3,j}\|_{W^{1-\frac{1}{p},p}} +\|F_{s,\mathrm{int}}^{4,j}\|_{W^{2-\frac{1}{p},p}} +\|F_{s,\mathrm{int}}^{6,j}\|_{W^{1-\frac{1}{p},p}}\Big).
\end{aligned}
\end{equation}
The constant depends on $\varepsilon$, since $\func{supp}\chi_{j,\mathrm{int}}$ comes within distance $c\varepsilon$ of $\Gamma$. The smallness assumption holds when $\mathcal{E}$ is small.
\end{lemma}

\begin{proof}
On $\func{supp}\chi_{j,\mathrm{int}}$ the boundary of $\Omega$ consists of a part of $\Sigma$ and a part of $\Sigma_{s}$. Both stay at distance $\geq c\varepsilon$ from $\Gamma$, hence from each other, and both are smooth by hypothesis~(3) of the Main result. The localized unknowns vanish near $\Sigma_{\mathrm{art}}$. So the problem may be posed on a smooth bounded domain that contains $\func{supp}\chi_{j,\mathrm{int}}$, whose boundary consists of these two parts and an artificial part where the unknowns vanish; its curvature is $O(1/\varepsilon)$.

\emph{Constant coefficients.} By the kinematic equation, the normal velocity on $\Sigma$ is $v_{j,\mathrm{int}}\cdot\mathcal{N}_{0}=\partial_{t}\eta_{j,\mathrm{int}}-F^{4,j}_{\mathrm{int}}-F_{s,\mathrm{int}}^{4,j}$, which lies in $W^{2-\frac1p,p}(\Sigma) $, because $H^{\frac32-\alpha}(\Sigma)\subset W^{2-\frac1p,p}(\Sigma)$ for $p<\frac{2}{1+2\alpha}$. The tangential part of the stress condition prescribes the tangential stress on $\Sigma$, in $W^{1-\frac1p,p}$, and on $\Sigma_{s}$ we have $v\cdot n_{0}=0$ and the Navier condition. With the term $\partial_{t}v_{j,\mathrm{int}}$ moved to the right-hand side, this is the Stokes system with prescribed normal velocity and tangential stress. To check the complementing condition, flatten a smooth face, freeze the principal coefficients in orthonormal coordinates, and take a nonzero tangential frequency. A decaying homogeneous half-space mode has zero normal velocity and zero tangential traction on the plane. Testing the momentum equation with the conjugate velocity over one tangential period (per unit length in the other tangential direction) gives zero boundary work. The pressure term vanishes by the divergence equation and the integral of the square of the symmetric gradient is zero. The velocity is therefore rigid and decay forces it to vanish; the momentum equation at this frequency then gives zero pressure. This verifies the complementing condition for the principal system on each physical face. The Navier friction is of lower order: its $W^{1-1/p,p}$ trace norm is controlled by the retained $H^1$ velocity norm since $p<2$. The smooth-boundary systems estimate of Agmon, Douglis and Nirenberg~\cite{ADN1964} now gives the bound for $v_{j,\mathrm{int}}$ in $W^{2,p}$ and for $p_{j,\mathrm{int}}$ minus its mean in $W^{1,p}$. The mean is not determined by these boundary conditions, and it is controlled by $\|p_{j,\mathrm{int}}\|_{L^{2}}$ on the right-hand side. Set $n=\mathcal N_0/|\mathcal N_0|$ and let $a_\zeta$ be the matrix in \eqref{def:linearized_operator}. Taking the normal component gives
\begin{equation}
g\eta_{j,\mathrm{int}}-\sigma\operatorname{div}_{\Sigma} (a_\zeta\nabla_{\Sigma}\eta_{j,\mathrm{int}}) =n\cdot S(p_{j,\mathrm{int}},v_{j,\mathrm{int}})n -\frac{n\cdot(F_{\mathrm{int}}^{3,j}+F_{s,\mathrm{int}}^{3,j})}{|\mathcal N_0|} +\sigma\operatorname{div}_{\Sigma}F_{\mathrm{int}}^{5,j}.
\end{equation}
Every term on the right belongs to $W^{1-1/p,p}$ under the hypotheses of the lemma. This is a uniformly elliptic equation on a smooth part of $D$, with $\eta_{j,\mathrm{int}}=0$ near $\Gamma_{\mathrm{art}}$, and its right-hand side lies in $W^{1-\frac1p,p}(\Sigma)$ by the trace theorem, the capillary flux being in $W^{2-\frac1p,p}$. Hence $\eta_{j,\mathrm{int}}\in W^{3-\frac1p,p}$ with the stated bound.
\end{proof}

This lemma has the following corollary.
\begin{corollary}\label{cor:enhanced elliptic estimate} Let $1<p<\frac{2}{1+2\alpha}$ and $0\leq j\leq m$. Let $\mathcal{U}$ be the neighbourhood of $\Gamma$ of Lemma~\ref{lem:tl_coef}, and $r_{-}=\frac{2}{1-\varepsilon_{-}}$. Assume that, on $\func{supp}\chi_{j,\mathrm{int}}$, the steady state $\zeta$ is smooth. Then, for every nonnegative integer $\ell\geq0$, every solution of \eqref{eq:interior-domain-elliptic-system-j} satisfies
\begin{align}{\label{eq:enhance_el_estimate}} &\|v_{j,\mathrm{int}}\|_{W^{2+\ell,p}}+\|\eta_{j,\mathrm{int}}\|_{W^{3+\ell-\frac{1}{p},p}}+\|p_{j,\mathrm{int}}\|_{W^{1+\ell,p}}  \notag \\
&\leq C(\varepsilon)\Big(\|F_{\mathrm{int}}^{1,j}\|_{W^{\ell,p}}+\|F_{\mathrm{int}}^{2,j}\|_{W^{1+\ell,p}}+\|F_{\mathrm{int}}^{3,j}\|_{W^{1+\ell-\frac{1}{p},p}}+\|F_{\mathrm{int}}^{4,j}\|_{W^{2+\ell-\frac{1}{p},p}}+\|F_{\mathrm{int}}^{5,j}\|_{W^{2+\ell-\frac{1}{p},p}}+\|\partial_{t}\eta_{j,\mathrm{int}}\|_{W^{2-\frac{1}{p}+\ell,p}} \notag \\
&\qquad+\|v_{j,\mathrm{int}}\|_{W^{1+\ell,p}}+\|\eta_{j,\mathrm{int}}\|_{W^{2-\frac{1}{p}+\ell,p}}+\|p_{j,\mathrm{int}}\|_{W^{\ell,p}}+\|F_{\mathrm{int}}^{6,j}\|_{W^{1+\ell-\frac{1}{p},p}}+\|F_{s,\mathrm{int}}^{1,j}\|_{W^{\ell,p}}\notag\\ &\qquad+\|F_{s,\mathrm{int}}^{2,j}\|_{W^{1+\ell,p}}+\|F_{s,\mathrm{int}}^{3,j}\|_{W^{1+\ell-\frac{1}{p},p}}+\|F_{s,\mathrm{int}}^{4,j}\|_{W^{2+\ell-\frac{1}{p},p}}+\|F_{s,\mathrm{int}}^{6,j}\|_{W^{1+\ell-\frac{1}{p},p}}+\|\partial_{t}v_{j,\mathrm{int}}\|_{W^{\ell,p}} \Big),
\end{align}
given that the right-hand side is bounded. The constant depends on $\varepsilon$, since $\func{supp}\chi_{j,\mathrm{int}}$ comes within distance $c\varepsilon$ of $\Gamma$. The smallness assumption holds when $\mathcal{E}$ is small.
\end{corollary}

\begin{proof}
On a smooth bounded domain, higher-order $L^p$ regularity for the Stokes system with prescribed normal velocity and tangential stress follows from the Agmon--Douglis--Nirenberg theory \cite{ADN1964}, since these boundary conditions satisfy the complementing condition. See also \cite[Theorem~5.12]{AACG2018} for a corresponding regularity result with homogeneous normal velocity; nonhomogeneous normal data can be treated by a suitable lifting. Here, we use only integer values of $\ell\ge0$; the case $\ell=0$ is the preceding lemma. The localized boundary pieces are smooth and separated from the contact curve. After lifting the divergence and prescribed normal velocity, the smooth-domain Stokes estimate with normal velocity and tangential stress gives $W^{2+\ell,p}$ regularity for the velocity and $W^{1+\ell,p}$ regularity for the pressure modulo constants. A change of the friction coefficient to the fixed value $1$ is paid for by the lower velocity norm already on the right-hand side. The retained pressure norm controls the pressure constant. The normal stress equation then gives the stated surface estimate by smooth elliptic regularity. The constant may depend on $\ell$, $\varepsilon$, and the corresponding finite collection of smooth geometry norms.
\end{proof}
\subsection{Basic interior estimate}

\begin{theorem}
$\chi_{m,int}v,\chi_{m,int}p,\chi_{m,int,\Sigma}\eta$ satisfy the following bounds
\begin{equation}\label{est:basic_diss_int} \|\chi_{m,int}\partial_{t}^{m}v\|_{W^{2,p}}+\|\chi_{m,int}\partial_{t}^{m}p\|_{W^{1,p}}+\|\chi_{m,int,\Sigma}\partial_{t}^{m}\eta\|_{W^{3-\frac{1}{p},p}} \lesssim{}  \sqrt{\mathcal{E}}\sqrt{\mathcal{D}}+\frac{C(\varepsilon)}{\varepsilon}(\sqrt{\mathcal{D}_{||}}+\sqrt{\mathcal{D}_{\mathrm{enhanced}}}),
\end{equation}
and
\begin{equation}\label{est:basic_energy_int} \sum_{j=0}^{m-1}\bigl(\|\chi_{m,int}\partial_{t}^{j}v\|_{W^{2,p}}+\|\chi_{m,int}\partial_{t}^{j}p\|_{W^{1,p}}+\|\chi_{m,int,\Sigma}\partial_{t}^{j}\eta\|_{W^{3-\frac{1}{p},p}}\bigr) \lesssim{}  {\mathcal{E}}+\frac{C(\varepsilon)}{\varepsilon}(\sqrt{\mathcal{E}_{||}}+\sqrt{\mathcal{E}_{\mathrm{enhanced}}}).
\end{equation}
\end{theorem}

\begin{proof}
We apply Lemma~\ref{lem:interior_elliptic} with $j=m$; its smallness assumption holds for small $\mathcal{E}$, as shown there. It remains to bound the forcing terms on its right-hand side. The major difference between the interior estimate and the intermediate estimate is that we must derive the estimate in 3-D. So we need to be very careful about the exponents in H\"older's inequality. Also, for the nonlinear terms, we do not have difference terms caused by the geometry or flattening. Set $\chi_{**}=\chi_{bet}+\chi_{int}$. Moreover, by Theorem \ref{thm:localized-surface-coordinate-estimate}, we have the norm equivalence between $\mathfrak{h}_{\eta}$ and $\eta$. So we use $\eta$ in this section.

\textbf{Cutoff commutators in $F^{1}$.}

\begin{equation}\label{eq:bet-cut-momentum_int}
\begin{aligned}
&\sum_{j=0}^{m} \left\| \left( S(\partial_{t}^{j}p,\partial_{t}^{j}v)\nabla\chi_{m,\mathrm{int}} +\partial_{t}^{j}v\Delta \chi_{m,int} \right) \right\|_{L^{p}} \lesssim C(\varepsilon)\frac{1}{\varepsilon}(\|S((\chi_{**}\partial_{t}^{j}p),(\chi_{**}\partial_{t}^j{}v))\|_{L^{2}}+\|(\chi_{**}\partial_{t}^{j}v)\|_{L^{2}})\\
&\lesssim C(\varepsilon)\frac{1}{\varepsilon}\sum_{j=0}^{m}(\|(\partial_{t}^{j}\chi_{bet}v)\|_{L_{\theta}^{2}W^{1,2}}+\|\chi_{bet}\partial_{t}^{j}v\|_{L_{\theta}^{2}L^{2}}+\|\partial_{t}^{j}\partial_{\theta}(\chi_{bet}v)\|_{L_{\theta}^{2}L^{2}})+\sum_{j=0}^{m}\|(\partial_{t}^{j}\chi_{bet}v)\|_{L_{\theta}^{2}L^{2}}\\
&+\sum_{j=0}^{m}(\|\partial_{t}^{j}\chi_{int}v)\|_{H^{1}}+\|\chi_{int}\partial_{t}^{j}v\|_{L^{2}}+\|\chi_{bet}\partial_{t}^{j}p\|_{L^{2}}) \lesssim C(\varepsilon)\frac{1}{\varepsilon} (\sqrt{\mathcal{D}_{||}}+\sqrt{\mathcal{D}_{\mathrm{enhanced}}}).
\end{aligned}
\end{equation}

\textbf{Second-order spatial derivative terms in the momentum equation}

We have
\begin{equation}
\begin{aligned}
&\sum_{j=0}^{m}\sum_{a+b+c=j,c\neq j}\left\| \left(\chi_{m,int}\operatorname{div}_{\p_{t}^{a}\mathcal{A}_{\eta}}\mathbb{D}_{\p_{t}^{b}\mathcal{A}_{\eta}}(\partial_{t}^{c}v) \right) \right\|_{L^{p}}+\sum_{j=0}^{m}\sum_{a+b=j,b<j}\left\| \left(\chi_{m,int}\operatorname{div}_{\p_{t}^{a}\mathcal{A}_{\eta}}(\partial_{t}^{b}p I) \right) \right\|_{L^{p}}\\
&\lesssim \sum_{j=0}^{m}\sum_{\substack{a+b+c=j,\ c<j\\a,b,c<m}}\|\chi_{m+,\mathrm{int}}\partial_{t}^{a}\mathcal{A}_{\eta}\|_{L^{\infty}}\|\chi_{m+,\mathrm{int}}{\partial_{t}^{b}\mathcal{A_{\eta}}}\|_{L^{\infty}}\|(\chi_{m,int}\nabla^{2}\partial_{t}^{c}v)\|_{L^{p}}\\
&+\sum_{j=0}^{m}\sum_{\substack{a+b=j,\ b<j\\ a<m}}\|\chi_{m,\mathrm{int}}\nabla\partial_{t}^{b}p\|_{L^{p}}\|\chi_{m+,\mathrm{int}}\partial_{t}^{a}\mathcal{A}_{\eta}\|_{L^{\infty}}\\
&+\|\chi_{m+,\mathrm{int}}\mathcal{A}_{\eta}\|_{L^{\infty}}\|\chi_{m+,\mathrm{int}}{\partial_{t}^{m}\mathcal{A_{\eta}}}\|_{L^{\infty}}\|(\chi_{m,\mathrm{int}}\nabla^{2}v,\chi_{m,\mathrm{int}}\nabla p)\|_{L^{p}}\\
&+\sum_{j=0}^{m}\sum_{\substack{a+b+c=j\\a<m,\ c<j}}\|\chi_{m+,\mathrm{int}}\partial_{t}^{a}\mathcal{A}_{\eta}\|_{L^{\infty}}\|\chi_{m,\mathrm{int}}{\partial_{t}^{b}\nabla \mathcal{A_{\eta}}}(\chi_{m+,\mathrm{int}}\nabla\partial_{t}^{c}v)\|_{L^{p}}\\
&+\|\chi_{m+,int}\partial_{t}^{m}\mathcal{A}_{\eta}\|_{L^{\infty}}\|\chi_{m,\mathrm{int}}\nabla{\mathcal{A_{\eta}}}(\chi_{m+,\mathrm{int}}\nabla v)\|_{L^{p}} =I_{1,1}+I_{1,2}+I_{2}+I_{3}+I_{4}.
\end{aligned}
\end{equation}

\emph{Remark.} The reason why we have to use the commutator here, but not for the intermediate estimate, is that the cutoff function is chosen to commute with the $\theta$ derivative; it cannot commute with all spatial derivatives. Moreover, for the terms without a cutoff function inside the norm, we actually omit the cutoff function $\|(\chi_{\mathrm{int}}+\chi_{\mathrm{bet}})\|$ inside the norm. We will keep this convention throughout this section.

For the term $I_{1}=I_{1,1}+I_{1,2}$, we have to bound $\|\partial_{t}^{b}\mathcal{A}_{\eta}\|_{L^{\alpha}}$ for $0<b\leq m$. By definition, we have
\begin{equation}\label{eq:int_linear_first_derivatives}
\begin{aligned}
&\|\partial_{t}^{b}\bar{\eta}\|_{W^{1,+\infty}} \lesssim (\|\chi_{int}\partial_{t}^{b}{\eta}\|_{W^{3-\frac{1}{p},p}}\\
&+\|\chi_{bet}\partial_{t}^{b}{\eta}\|_{L_{||}^{\infty}H^{\frac{3}{2}+\frac{\varepsilon_{-}-\alpha}{2}}} +\|\chi_{\mathrm{cl}}\partial_{t}^{b}{\eta}\|_{L_{||}^{\infty}H^{\frac{3}{2}+\frac{\varepsilon_{-}-\alpha}{2}}} +\|\chi_{bet}\partial_{t}^{b}\partial_{\theta}{\eta}\|_{H^{\frac{3}{2}-\alpha}} +\|\chi_{\mathrm{cl}}\partial_{t}^{b}\partial_{\theta}{\eta}\|_{H^{\frac{3}{2}-\alpha}}) \lesssim\sqrt{\mathcal{D}},
\end{aligned}
\end{equation}
where we used the Sobolev embedding $H^{2+\frac{\varepsilon_{-}-\alpha}{2}}\hookrightarrow W^{1,\infty}$ in two dimensional space. Moreover, when $0<b\leq m-1$, using the same argument, we have $\|\partial_{t}^{b}\bar{\eta}\|_{L^{\infty}}\lesssim \sqrt{\mathcal{E}}.$ These two estimates above induce the following result by the definition of $W$
\begin{align}
\|\partial_{t}^{m}(W\bar{\eta})\|_{W^{1,+\infty}}\lesssim \sqrt{\mathcal{D}}+\sqrt{\mathcal{E}\mathcal{D}},\quad\operatorname{and}\quad \sum_{j=0}^{m-1}\|\partial_{t}^{j}(W\bar{\eta})\|_{W^{1,+\infty}}\lesssim \sqrt{\mathcal{E}}+{\mathcal{E}}
\end{align}
Use \eqref{eq:int_linear_first_derivatives} and Lemma~\ref{lem:E_le_D} for $\p_{t}^{b}\mathcal{A}_{\eta}$. Thus,
\begin{equation}\label{eq:interior_e_1} \|\partial_t^b\mathcal A_\eta\|_{L^\infty} \leq C(\|\p_{t}^{b}\bar{\eta}\|_{W^{1,+\infty}}+ \sum_{r=2}^{b}\ \sum_{\substack{i_1+\cdots+i_r=b\\i_k\geq1}} \prod_{k=1}^r \|\p_{t}^{i_{k}}(W\bar{\eta})\|_{W^{1,+\infty}}+\sqrt{\mathcal{E}\mathcal{D}} ) \leq C(\varepsilon)(1+\sqrt{\mathcal E})\sqrt{\mathcal D} \leq C(\varepsilon)\sqrt{\mathcal D},
\end{equation}
when $b=m$. And $\|\partial_t^b\mathcal A_\eta\|_{L^\infty} \leq C(\varepsilon)(1+\sqrt{\mathcal E})\sqrt{\mathcal E} \leq C(\varepsilon)\sqrt{\mathcal E},$ when $0<b\leq m-1$.

For the estimates for velocity and pressure, we have
\begin{align}
\|(\chi_{m,\mathrm{int}})\nabla^{2}\partial_{t}^{c}v\|_{L^{p}}\lesssim \frac{1}{\varepsilon}\|(\chi_{\mathrm{int}}+\chi_\mathrm{bet})\partial_{t}^{c}v\|_{H^{1}}+\|\chi_{m,\mathrm{int}}\partial_{t}^{c}v\|_{W^{2,p}}.
\end{align}
Similarly,
\begin{align}
\|(\chi_{m,\mathrm{int}})\nabla\partial_{t}^{b}p\|_{L^{p}}\lesssim \frac{1}{\varepsilon}\|(\chi_{\mathrm{int}}+\chi_\mathrm{bet})\partial_{t}^{b}p\|_{L^{2}}+\|\chi_{m,\mathrm{int}}\partial_{t}^{b}p\|_{W^{1,p}}.
\end{align}
In $I_{1,1}$ we have $c<j$, hence $a+b>0$; in $I_{1,2}$ we have $b<j$, hence $a>0$. At least one matrix factor is therefore bounded by $C(\varepsilon)\sqrt{\mathcal E}$. The indicated velocity and pressure norms of order strictly less than $m$ are also at the energy level, by \eqref{def:interior_energy} and the global $H^1$ velocity and $L^2$ pressure terms. Using small energy and Lemma~\ref{lem:E_le_D} gives $I_1\leq C(\varepsilon)\mathcal E \leq C(\varepsilon)\sqrt{\mathcal E\mathcal D}.$ For $I_2$, apply \eqref{eq:A_equi} on the same set $U$. For the $(v,p)$ term contained in $I_{2}$, we have
\begin{align}
{\label{eq:key_c_1}} \begin{aligned} \|(\chi_{m,\mathrm{int}}\nabla^{2}v,\chi_{m,\mathrm{int}}\nabla p)\|_{L^{p}}&\lesssim \|\chi_{m,int}v\|_{W^{2,p}}+\|\chi_{m,int}p\|_{W^{1,p}}+\|\chi_{bet}v\|_{W^{1,p}}+\|\chi_{bet}p\|_{L^{p}}\\ &\lesssim \|\chi_{m,int}v\|_{W^{2,p}}+\|\chi_{m,int}p\|_{W^{1,p}}+\|\chi_{bet}v\|_{H^{1}}+\|\chi_{bet}p\|_{L^{2}}\lesssim\sqrt{\mathcal{E}}.
\end{aligned}
\end{align}
Therefore, we have $I_{2}\lesssim \sqrt{\mathcal{E}\mathcal{D}}.$ We now estimate the term $I_{3}$. Since $\|\partial_{t}^{a}\mathcal{A}_{\eta}\|_{L^{\infty}}$ has already been estimated, it remains to estimate the other factors. By definition and the commutator estimate \ref{thm:enhanced_commutator}, we have
\begin{equation}
\begin{aligned}
&\sum_{j=0}^{m}\sum_{\substack{a+b+c=j\\a<m,\ c<j}}\|\chi_{m+,\mathrm{int}}\partial_{t}^{a}\mathcal{A}_{\eta}\|_{L^{\infty}}\|\chi_{m,\mathrm{int}}\nabla{\partial_{t}^{b}\mathcal{A_{\eta}}}(\chi_{m+,\mathrm{int}}\nabla\partial_{t}^{c}v)\|_{L^{p}}\\
&\lesssim \sum_{j=0}^{m}\sum_{\substack{a+b+c=j\\a<m,\ c<j}}\|\chi_{m+,\mathrm{int}}\partial_{t}^{a}\mathcal{A}_{\eta}\|_{L^{\infty}}\|\chi_{m-1,\mathrm{int}}\nabla {\partial_{t}^{b}\mathcal{A_{\eta}}}(\nabla\partial_{t}^{c}v)\|_{L^{p}}\\
&+\sum_{j=0}^{m}\sum_{\substack{a+b+c=j\\a<m,\ c<j}}\|\chi_{m+,\mathrm{int}}\partial_{t}^{a}\mathcal{A}_{\eta}\|_{L^{\infty}}\|(\chi_{m,\mathrm{int}}-\chi_{m-1,\mathrm{int}})\nabla {\partial_{t}^{b}\mathcal{A_{\eta}}}(\nabla\partial_{t}^{c}v)\|_{L^{p}} =I_{3,1}+I_{3,2}.
\end{aligned}
\end{equation}
Since $c<j$, we have $a+b>0$. If $a=0$, the undifferentiated coefficient is bounded by $1+\sqrt{\mathcal{E}}$, and the factor with $b>0$ supplies the differentiated geometry. If $b=0$, the positive time derivative lies in $a$. Thus the cutoff split has no additional small factor; the following product estimates give $I_{3}\lesssim C(\varepsilon)\sqrt{\mathcal{E}\mathcal{D}}$.

For the term $I_{3,1}$, we have
\begin{align}
I_{3,1}\lesssim \sum_{j=0}^{m}\sum_{\substack{a+b+c=j\\a<m,\ c<j}}\|\chi_{m+,\mathrm{int}}\partial_{t}^{a}\mathcal{A}_{\eta}\|_{L^{\infty}}\|\chi_{m-1,\mathrm{int}}\nabla {\partial_{t}^{b}\mathcal{A_{\eta}}}\|_{L^{3}}\|(\nabla\partial_{t}^{c}(\chi_{m,\mathrm{int}}v)\|_{L^{\frac{3p}{3-p}}}.
\end{align}

To estimate $\|\chi_{m-1,\mathrm{int}}\nabla {\partial_{t}^{b}\mathcal{A_{\eta}}}\|_{L^{3}}$, we note that it contains linear part and nonlinear part. For the nonlinear part, when spatial derivative hit on $\nabla \bar{\eta}$, then we separate this term into interior part and intermediate part, both parts are bounded in $W^{2,3}$. When spatial derivative hit on other functions independent of $\eta$, the less spatial derivative are applied to the term with $\eta$, it is then bounded by the multiplication of $W^{1,\infty}$ norm of $\bar{\eta}$ with corresponding temporal derivative. More specifically, we have
\begin{equation}\label{eq:A_d}
\begin{aligned}
&\|\chi_{m-1,\mathrm{int}}\partial_{t}^{b}\nabla\mathcal{A}_{\eta}\|_{L^{3}} \lesssim \sqrt{\mathcal{E}}\sqrt{\mathcal{D}}+ \|(\partial_{t}^{b}\overline{(( \chi_{\mathrm{int},\Sigma}){\eta})}\|_{W^{2,3}}+\|(\partial_{t}^{b}\overline{(( \chi_{m,\mathrm{int},\Sigma}-\chi_{\mathrm{int},\Sigma}){\eta})}\|_{W^{2,3}}\\
&+\sum_{r=2}^{b+1}\sum_{\substack{i_{1}+\ldots+i_{\ell}=b\\i_{1}\geq 0, i_{2},\ldots i_{\ell}\geq 1}}\big(\|(\partial_{t}^{i_{1}}\overline{(( \chi_{\mathrm{int},\Sigma}){\eta})}\|_{W^{2,3}}+\|(\partial_{t}^{i_{1}}\overline{(( \chi_{m,\mathrm{int},\Sigma}-\chi_{\mathrm{int},\Sigma}){\eta})}\|_{W^{2,3}}\big)\prod_{k=2}^{r}\|(\partial_{t}^{i_{k}}\overline{{\eta}})\|_{W^{1,\infty}}\\
&\lesssim \sqrt{\mathcal{E}}\sqrt{\mathcal{D}}+ (1+\sqrt{\mathcal{E}})\sum_{k=0}^{b}\Big(\|\chi_{int,\Sigma }\partial_{t}^{k}\eta\|_{W^{3-\frac{1}{p},p}}+\|\partial_{\theta}^{2}\partial_{t}^{k}\overline{ (\chi_{m,\mathrm{int},\Sigma}-\chi_{\mathrm{int},\Sigma})\eta}\|_{H^{\frac{1}{6}}H^{\frac{1}{3}}}\\
&+\|\partial_{\theta}\partial_{t}^{k}\overline{ (\chi_{m,\mathrm{int},\Sigma}-\chi_{\mathrm{int},\Sigma})\eta}\|_{H^{\frac{1}{6}}H^{\frac{4}{3}}}+\|\partial_{t}^{k}\overline{ (\chi_{m,\mathrm{int},\Sigma}-\chi_{\mathrm{int},\Sigma})\eta})\|_{L^{3}W^{\frac{5}{2},p}}\Big)\\
&\lesssim  \sqrt{\mathcal{E}}\sqrt{\mathcal{D}}+(1+\sqrt{\mathcal{E}})\sum_{k=0}^{b}\Big(\|\chi_{int,\Sigma }\partial_{t}^{k}\eta\|_{W^{3-\frac{1}{p},p}}+\|\partial_{\theta}\partial_{t}^{k}{ (\chi_{m,\mathrm{int},\Sigma}-\chi_{\mathrm{int},\Sigma})\eta}\|_{H^{\frac{7}{6}}L^{2}}\\
&+\|\partial_{\theta}\partial_{t}^{k}{ (\chi_{m,\mathrm{int},\Sigma}-\chi_{\mathrm{int},\Sigma})\eta}\|_{H^{\frac{1}{6}}H^{\frac{6}{5}}}+\|\partial_{t}^{k}{ (\chi_{m,\mathrm{int},\Sigma}-\chi_{\mathrm{int},\Sigma})\eta})\|_{L^{3}W^{2,p}}\Big)\\
&\lesssim  \sqrt{\mathcal{E}}\sqrt{\mathcal{D}}+(1+\sqrt{\mathcal{E}})\sum_{k=0}^{b}\Big(\|\chi_{int,\Sigma }\partial_{t}^{k}\eta\|_{W^{3-\frac{1}{p},p}}+\|\partial_{\theta}\partial_{t}^{k}{ (\chi_{m,\mathrm{int},\Sigma}-\chi_{\mathrm{int},\Sigma})\eta}\|_{H^{\frac{3}{2}-\alpha}}\\
&+\|\partial_{t}^{k}{ (\chi_{m,\mathrm{int},\Sigma}-\chi_{\mathrm{int},\Sigma})\eta})\|_{W^{3-\frac{1}{p},p}}+\|\partial_{t}^{k}{ (\chi_{m,\mathrm{int},\Sigma}-\chi_{\mathrm{int},\Sigma})\eta})\|_{L^{\infty}H^{\frac{3}{2}-\alpha}}\Big) \lesssim \sqrt{\mathcal{D}},
\end{aligned}
\end{equation}
where we used interpolation theory for $L^{3}W^{2,p}$ with the assumption that $\alpha<\frac{3}{2}-\frac{2}{p}$. Moreover, $(v,p)$ can be estimated in the same way as \eqref{eq:key_c_1}, as follows:
\begin{align}
\begin{aligned} \|(\nabla \partial_{t}^{c}(\chi_{m,\mathrm{int}}v))\|_{L^{\frac{3p}{3-p}}}\lesssim \|\chi_{m,\mathrm{int}}\partial_{t}^{c}v\|_{W^{2,p}} \lesssim\sqrt{\mathcal{E}} \end{aligned}
\end{align}
for $c\leq m-1$. Therefore, combining the previous estimates for $\|\chi_{m+,\mathrm{int}}\p_{t}^{a}\mathcal{A}_{\eta}\|_{L^{\infty}}$, we have $I_{3,1}\lesssim\sqrt{\mathcal{E}\mathcal{D}}.$ For the term $I_{3,2}$, by the commutator estimate \ref{thm:enhanced_commutator}, we have
\begin{equation}
\begin{aligned}
I_{3,2}\leq C\sum_{j=0}^{m}\sum_{\substack{a+b+c=j\\a<m,\ c<j}} &\|\chi_{m+,\mathrm{int}}\partial_t^a\mathcal A_\eta\|_{L^\infty} \| (\chi_{m,\mathrm{int}}-\chi_{m-1,\mathrm{int}}) \nabla\partial_t^b\mathcal A_\eta\|_{L^p_\theta L^s_{rh}}\\
&\times\big(\|\chi_{\mathrm{bet}}\partial_t^c v\|_{L^\infty_\theta W^{2,q_\delta}_{rh}} +\|\chi_{\mathrm{bet}}\partial_\theta\partial_t^c v\|_{L^\infty_\theta L^{r_-}_{rh}}\big),
\end{aligned}
\end{equation}
where $r_-=2q_\delta/(2-q_\delta)$ and $1/s=1/p-1/r_-$. Indeed, we have $W^{1,q_\delta}_{rh}\subset L^{r_-}_{rh}$ which controls the transverse gradient, and the second norm controls the angular derivative. This estimate holds since $2q_{\delta}+pq_{\delta}-2p>0$. Now it suffices to estimate $\|(\chi_{m,\mathrm{int}}-\chi_{m-1,\mathrm{int}})\nabla \partial_{t}^{b}\mathcal{A}_{\eta}\|_{L^{p}L^{\frac{2pq_{\delta}}{2q_{\delta}+pq_{\delta}-2p}}}$. Use the product rule from \eqref{eq:A_d} in this mixed norm, with the cutoff on its single second-derivative factor. When some $i_\ell=m$, necessarily $b=m$, $r=1$, $k=0$; this contribution is at most $C(\varepsilon)\sqrt{\mathcal E\mathcal D}$ by \eqref{eq:int_linear_first_derivatives} and the $k=0$ linear extension bound below. Write $R_b=\boldsymbol{1}_{\{b=m\}}C(\varepsilon) \sqrt{\mathcal E\mathcal D}$. For the other terms all multipliers have energy bounds. By Theorem~\ref{thm:enhanced_commutator}, for $b\geq1$ we obtain
\begin{align}
&\|(\chi_{m,\mathrm{int}}-\chi_{m-1,\mathrm{int}})\nabla \partial_{t}^{b}\mathcal{A}_{\eta}\|_{L^{p}L^{\frac{2pq_{\delta}}{2q_{\delta}+pq_{\delta}-2p}}}  \notag \\
\lesssim&R_b+C(\varepsilon)(1+\sqrt{\mathcal{E}})\sum_{k=0}^{b}\bigl(\|(\chi_{m,\mathrm{int}}-\chi_{m-1,\mathrm{int}})\partial_{\theta}^{2}\partial_{t}^{k}\bar{\eta}\|_{L^{p}L^{\frac{2pq_{\delta}}{2q_{\delta}+pq_{\delta}-2p}}}+\|(\chi_{m,\mathrm{int}}-\chi_{m-1,\mathrm{int}})\partial_{\theta}\partial_{t}^{k}\bar{\eta}\|_{L^{p}W^{1,\frac{2pq_{\delta}}{2q_{\delta}+pq_{\delta}-2p}}} \\
&+\|(\chi_{m,\mathrm{int}}-\chi_{m-1,\mathrm{int}})\partial_{t}^{k}\bar{\eta}\|_{L^{p}W^{2,\frac{2pq_{\delta}}{2q_{\delta}+pq_{\delta}-2p}}}\big)  \notag \\
\lesssim&R_b+C(\varepsilon)(1+\sqrt{\mathcal{E}})\sum_{k=0}^{b}\bigl(\|(\chi_{\mathrm{bet}})\partial_{\theta}\partial_{t}^{k}\bar{\eta}\|_{H^{1}H^{1-\alpha}}+\|(\chi_{m,\mathrm{int}}-\chi_{m-1,\mathrm{int}})\partial_{\theta}\partial_{t}^{k}\bar{\eta}\|_{L^{2}H^{2-\alpha}}+\|\psi\,\bar{E}_{\rho}f_{k}\|_{L^{p}W^{2,s}}\big),
\end{align}
where $s:=\frac{2pq_{\delta}}{2q_{\delta}+pq_{\delta}-2p}$, $\psi:=\chi_{m,\mathrm{int}}-\chi_{m-1,\mathrm{int}}$ and $f_{k}:=(\chi_{m+1,\mathrm{int},\Sigma}-\chi_{\mathrm{int},\Sigma})\partial_{t}^{k}\eta$. Here we used that on $\func{supp}\psi$ the extension $\partial_{t}^{k}\bar{\eta}$ of Subsection~\ref{subsec:extension_global} equals $\bar{E}_{\rho}f_{k}$, up to the terms $(1-\chi_{\mathrm{dif},j+})P_{2-D}E(\chi_{\mathrm{dif}, j,\Sigma}\partial_t^k\eta)$. The evaluation region is at distance at least $c\varepsilon$ from the support of the extended boundary datum, as in the proof of Theorem~\ref{thm:enhanced_commutator}. Here separation bounds every required \emph{transverse} kernel derivative by $C(\varepsilon)$; an angular derivative passes to the datum. For a datum $g=\chi_{\mathrm{dif},j,\Sigma}\partial_t^k\eta$, all these terms of total spatial order up to two in $L^p_\theta L^s_{rh}$ are therefore bounded by $C(\varepsilon)\sum_{a=0}^2\|\partial_\theta^a g\|_{L^2_\theta L^1_r}$: use $p<2$ and the finite sizes of the collar and of the support of $Eg$. The terms with $a=1,2$ are controlled by $\|\partial_\theta g\|_{H^{3/2-\alpha}}$, since $3/2-\alpha>1$; the term $a=0$ uses the surface $L^2$ norm. These are among the surface bounds below since the difference cutoffs are smooth multipliers of $\chi_{\mathrm{bet},\Sigma}$. The datum $f_k$ has only the slice regularity $W^{3- \frac{1}{q_{\delta}},q_{\delta}}$ on the transition region of $\chi_{m,\mathrm{int}}$, so we split off the commutator of Theorem~\ref{thm:enhanced_commutator}. Since $\chi_{m+1,\mathrm{int},\Sigma}=1$ on $\func{supp}\chi_{m,\mathrm{int},\Sigma}$ and $\psi\chi_{\mathrm{int}}=0$, we have $\psi_{\Sigma}f_{k}=(\chi_{m,\mathrm{int},\Sigma}-\chi_{m-1,\mathrm{int},\Sigma})\partial_{t}^{k}\eta$, which is controlled by the interior norms. With $L^{2}_{\theta}W^{4,q_{\delta}}\subset L^{p}_{\theta}W^{2,s}$ and $L^{p}_{\theta}W^{3,p}\subset L^{p}_{\theta}W^{2,s}$ on the slices, as $s<\frac{2p}{2-p}$, we obtain
\begin{align}
\|\psi\,\bar{E}_{\rho}f_{k}\|_{L^{p}W^{2,s}}\leq{}&\big\|\psi\,\bar{E}_{\rho}f_{k}-\bar{E}_{\rho}(\psi_{\Sigma}f_{k})\big\|_{L^{p}W^{2,s}}+\big\|\bar{E}_{\rho}(\psi_{\Sigma}f_{k})\big\|_{L^{p}W^{2,s}}  \notag \\
\lesssim{}&\big\|\psi\,\bar{E}_{\rho}f_{k}-\bar{E}_{\rho}(\psi_{\Sigma}f_{k})\big\|_{L^{2}W^{4,q_{\delta}}}+\big\|\bar{E}_{\rho}(\psi_{\Sigma}f_{k})\big\|_{L^{p}W^{3,p}}  \notag \\
\leq{}&C(\varepsilon)\Big(\|(\chi_{m+1,\mathrm{int},\Sigma}-\chi_{\mathrm{int},\Sigma})\partial_{t}^{k}\eta\|_{L^{2}W^{3-\frac{1}{q_{\delta}},q_{\delta}}}+\|(\chi_{m,\mathrm{int},\Sigma}-\chi_{m-1,\mathrm{int},\Sigma})\partial_{t}^{k}\eta\|_{W^{3-\frac{1}{p},p}}\Big),  \notag
\end{align}
by \eqref{eq:enhanced_commutator}. Altogether,
\begin{align}
&\|(\chi_{m,\mathrm{int}}-\chi_{m-1,\mathrm{int}})\nabla \partial_{t}^{b}\mathcal{A}_{\eta}\|_{L^{p}L^{s}}   \leq R_b  \notag \\
&+C(\varepsilon)\sum_{k=0}^{b}\bigl(\|(\chi_{\mathrm{bet}})\partial_{\theta}\partial_{t}^{k}{\eta}\|_{H^{\frac{3}{2}-\alpha}}+\|{(\chi_{m,\mathrm{int,\Sigma}}-\chi_{m-1,\mathrm{int},\Sigma})\partial_{t}^{k}\eta}\|_{W^{3-\frac{1}{p},p}}+\|{(\chi_{m+1,\mathrm{int,\Sigma}}-\chi_{\mathrm{int},\Sigma})\partial_{t}^{k}\eta}\|_{L^{2}W^{3-\frac{1}{q_{\delta}},q_{\delta}}}\big).
\end{align}
For $b=0$ the same linear bounds applied to $\nabla\mathcal A_\eta=\mathfrak A'(\nabla q)[\nabla^2q]$ give $\|\psi\nabla\mathcal A_\eta\|_{L^p_\theta L^s_{rh}} \leq C(\varepsilon)\sqrt{\mathcal E}$. More generally the preceding bounds use only energy for $b<m$, and give $C(\varepsilon)\sqrt{\mathcal D}$ for $b=m$. If $c\leq m-2$ the velocity norm $L^\infty_\theta W^{2,q_\delta}$ uses energy: take its $L^2_\theta$ norm and that of one more angular derivative. Their total order is at most $m-1$. If $c=m-1$ the same argument uses dissipation; then $a+b\leq1$ and the matrix gradient is at the energy level. For the angular velocity term take one and two angular derivatives of the localized $\partial_t^c v$ in $L^2_\theta H^1_{rh}$, using $H^1_{rh}\subset L^{r_-}_{rh}$. The required total order $c+2$ is at most $m$ in energy and $m+1$ in dissipation, by \eqref{def:energy_cl}. On the strip $\chi_{\rm bet}=1$, so the $\chi_{\rm bet+}$ localization in that definition suffices. Combining these bounds and Lemma~\ref{lem:E_le_D}, $I_{3,2}\leq C(\varepsilon)\sqrt{\mathcal E\mathcal D}.$ For $I_4$, use \eqref{eq:interior_e_1} in $L^\infty$. On the support of $\chi_{m-1,\mathrm{int}}$, use \eqref{eq:A_d} with $b=0$ and $\nabla(\chi_{m,\mathrm{int}}v)\in L^{3p/(3-p)}$. On the strip use $\psi\nabla\mathcal A_\eta\in L^p_\theta L^s_{rh}$ and $\nabla v\in L^\infty_\theta L^{2q_\delta/(2-q_\delta)}_{rh}$ there, as above with $b=c=0$. These pairs bound the undifferentiated product by $C(\varepsilon)\mathcal E$. Hence, for $\mathcal E\leq1$, $I_4\leq C(\varepsilon)\sqrt{\mathcal D}\,\mathcal E \leq C(\varepsilon)\sqrt{\mathcal E\mathcal D}.$ Finally, combining all of the estimates above, we have
\begin{align}
\sum_{j=0}^{m}\sum_{a+b+c=j, c\neq j} \left\| \left(\chi_{int}\func{div}_{\partial_t^{a}\mathcal{A}_{\eta}}S_{\partial_{t}^{b}\mathcal{A_{\eta}}}(\partial_{t}^{c}v,\partial_{t}^{c}p) \right) \right\|_{L^{p}}\lesssim \sqrt{\mathcal{E}\mathcal{D}}.
\end{align}

In particular, the preceding coefficient estimates give
\begin{equation}
\|\chi_{m+,\mathrm{int}}\partial_t^m\mathcal A_\eta\|_{L^\infty} \leq C(\varepsilon)\sqrt{\mathcal D},\quad \sum_{j=1}^{m-1}\|\chi_{m+,\mathrm{int}}\partial_t^j\mathcal A_\eta\|_{L^\infty} +\|\chi_{m+,\mathrm{int}}(I-\mathcal A_\eta)\|_{L^\infty} \leq C(\varepsilon)\sqrt{\mathcal E},
\end{equation}
and
\begin{align}
\|\chi_{m,\mathrm{int}}\partial_{t}^{m}v\|_{W^{2,p}}\lesssim\sqrt{\mathcal{D}}, \qquad \sum_{j=0}^{m-1}\|\chi_{m,\mathrm{int}}\partial_{t}^{j}v\|_{W^{2,p}}\lesssim\sqrt{\mathcal{E}}.
\end{align}
\textbf{Convection term in the momentum equation}

We now consider the term
\begin{equation}
\begin{aligned}
&\sum_{j=0}^{m}\|\partial_{t}^{j}(\chi_{m,int}v\cdot \nabla_{\mathcal{A}_{\eta}}v)\|_{L^{p}} \leq C_m\Bigl(\|\chi_{m,int}v\cdot \nabla _{\partial_{t}^{m}\mathcal{A}_{\eta}}v\|_{L^{p}}\\
&+\sum_{j=0}^{m}\sum_{a+b+c=j,b\leq m-1}\|\chi_{m,int}\partial_{t}^{a}v\cdot \nabla_{\partial_{t}^{b}\mathcal{A}_{\eta}}\partial_{t}^{c}v\|_{L^{p}}\Bigr)=:C_m(I_{1}+I_{2}).
\end{aligned}
\end{equation}
For $I_{1}$, using the Sobolev inequality,
\begin{align}
I_{1}\lesssim\|\chi_{m,\mathrm{int}}\partial_{t}^{m}\mathcal{A}_{\eta}\|_{L^{\infty}}\|\chi_{m+,\mathrm{int}}v\|_{L^{\infty}}\|\chi_{m+,\mathrm{int}}\nabla v\|_{L^{2}}.
\end{align}
Using the estimate for the previous stress tensor term, we have
\begin{equation}
\begin{aligned}
&\|\chi_{m,int}\partial_{t}^{m}\mathcal{A}_{\eta}\|_{L^{2}}\lesssim  \sqrt{\mathcal{D}},\quad \|\chi_{m+,\mathrm{int}}v\|_{L^{\infty}}\lesssim \|\chi_{int}v\|_{W^{2,p}}+\|(\chi_{bet})v\|_{L_{||}^{\infty}W^{2,q_{\delta}}}\lesssim \sqrt{\mathscr{E}_{0,el,int}}+\sqrt{\mathscr{E}_{0,el,bet}}\lesssim \sqrt{\mathcal{E}},\\
&\|\chi_{m+,\mathrm{int}}\nabla v\|_{L^{2}}\lesssim \|\chi_{m,\mathrm{int}}v\|_{W^{2,p}}+\|\chi_{\mathrm{bet}}v\|_{L^{2}H^{1}}\lesssim\sqrt{\mathcal{E}}.
\end{aligned}
\end{equation}
Therefore, $I_{1}$ satisfies the following boundedness property $I_{1}\lesssim \sqrt{\mathcal{E}\mathcal{D}}.$ For $I_{2}$, using the previous computation again, we have
\begin{equation}
\begin{aligned}
&I_{2}\lesssim \sum_{j=0}^{m}\sum_{a+b+c=j,(b)\neq m}\|\chi_{m+,\mathrm{int}}\partial_{t}^{a}v\|_{L^{\infty}}\|\chi_{m+,int}{\partial_{t}^{b}\mathcal{A}_{\eta}}\|_{L^{+\infty}}\|\chi_{m,\mathrm{int}}\partial_{t}^{c}\nabla v\|_{L^{2}}\\
&\lesssim  \sum_{j=0}^{m}\sum_{a+b+c=j,(b)\neq m}\|\chi_{m+,\mathrm{int}}\partial_{t}^{a}v\|_{W^{2,p}}\|\chi_{m+,int}{\partial_{t}^{b}\mathcal{A}_{\eta}}\|_{L^{+\infty}}(\|\chi_{m,\mathrm{int}}\partial_{t}^{c} v\|_{W^{2,p}}+\|\chi_{\mathrm{bet}}\partial_{t}^{c} v\|_{H^{1}}) \lesssim  \sqrt{\mathcal{E}\mathcal{D}}.
\end{aligned}
\end{equation}
The last inequality holds since $a,c$ cannot both equal $m$. Therefore, we have
\begin{align}
\begin{aligned} \sum_{j=0}^{m}\|\partial_{t}^{j}(\chi_{m,int}v\cdot \nabla_{\mathcal{A}_{\eta}}v)\|_{L^{p}}\lesssim& \sqrt{\mathcal{E}\mathcal{D}}. \end{aligned}
\end{align}

\textbf{The temporal derivative term in the momentum equation}

We now discuss
\begin{align}
&\sum_{j=0}^{m}\|\chi_{m,int}\partial_{t}^{j}(\partial_{t}\bar{\eta}WK\partial_{3}v)\|_{L^{p}}  \notag \\
&\lesssim\|\chi_{m,int}\partial_{t}^{m+1}\bar{\eta}WK\partial_{3}v\|_{L^{p}}+\|\chi_{m,int}\partial_{t}\bar{\eta}\partial_{t}^{m}WK\partial_{3}v\|_{L^{p}}+\|\chi_{m,int}\partial_{t}\bar{\eta}W\partial_{t}^{m}K\partial_{3}v\|_{L^{p}}+\|\chi_{m,int}\partial_{t}\bar{\eta}WK\partial_{3}\partial_{t}^{m}v\|_{L^{p}}  \notag \\
&\qquad+\sum_{j=0}^{m}\sum_{a+b+c+d=j,(a,b,c,d)\neq m}\|\chi_{m,int}(\partial_{t}^{a+1}\bar{\eta}\partial_{t}^{b}W\partial_{t}^{c}K\partial_{3}\partial_{t}^{d}v)\|_{L^{p}}=I_{1}+I_{2}+I_{3}+I_{4}+I_{5}.
\end{align}

For $I_{1}$, using the previous estimate for $\|\nabla \bar{\eta}\|_{W^{1,\infty}}$, we have
\begin{align}
I_{1}\lesssim& \|\chi_{m+,\mathrm{int}}\partial_{t}^{m+1}\bar{\eta}\|_{L^{\frac{1}{\alpha}}}\|\chi_{m+,\mathrm{int}}W\|_{L^{\infty}}\|\chi_{m+,\mathrm{int}}K\|_{L^{\infty}}\|\chi_{m+,\mathrm{int}}\partial_{3}v\|_{L^{2}} \\
\lesssim& (\|\chi_{\mathrm{int},\Sigma}\partial_{t}^{m+1}{\eta}\|_{H^{1}}+\|\partial_{t}^{m+1}\bar{\eta}_{m+1,\mathrm{int},d}\|_{L_{||}^{\infty}H^{1-\alpha}})   \times (\|\chi_{int}v\|_{H^{1}}+\|\chi_{bet}v\|_{H^{1}})   \\
\lesssim& (\|\chi_{int,\Sigma}\partial_{t}^{m+1}{\eta}\|_{H^{1}}+\|\partial_{t}^{m+1}{\eta}_{\mathrm{bet}}\|_{H^{\frac{1}{2}+\alpha}H^{\frac{1}{2}-\alpha}})   \times (\|\chi_{int}v\|_{H^{1}}+\|\chi_{bet}v\|_{H^{1}})   \\
\lesssim&(\sqrt{\mathscr{D}_{\mathrm{enhance}}}+\sqrt{\mathscr{D}_{m+1,energy,bet}})(\sqrt{\mathscr{E}_{0,energy,bet}}+\sqrt{\mathcal{E}_{\mathrm{enhance}}})   \lesssim\sqrt{\mathcal{E}\mathcal{D}}.
\end{align}

For the terms $I_{2}$ and $I_{3}$, using the previous estimate, we have
\begin{equation}
\begin{aligned}
&I_{2}+I_{3}\lesssim  \|\chi_{m+,int}\partial_{t}\bar{\eta}\|_{L^{\infty}}\|\chi_{m,int}\partial_{t}^{m}\nabla \bar{\eta}\|_{L^{+\infty}}\|\chi_{m+,int}\partial_{3}v\|_{L^{2}}\\
&\lesssim  (\|\chi_{int,\Sigma}\partial_{t}{\eta}\|_{W^{3-\frac{1}{p},p}}+\|(\chi_{bet,\Sigma})\partial_{t}\eta\|_{L_{||}^{\infty}W^{3-\frac{1}{q_{\delta}},q_{\delta}}})\\
&\times(\sqrt{\mathscr{D}_{m,el,int}}+\sqrt{\mathscr{D}_{m,el,bet}}+\sqrt{\mathscr{D}_{m+1,energy,bet}}) \times (\|\chi_{int}v\|_{H^{1}}+\|\chi_{bet}v\|_{L_{||}^{2}H^{1}}+\|\chi_{bet}\partial_{\theta}v\|_{L_{||}^{2}L^{2}})\\
&\lesssim (\sqrt{\mathscr{E}_{1,el,int}}+\sqrt{\mathscr{E}_{1,el,bet}}) \times(\sqrt{\mathscr{D}_{m,el,int}}+\sqrt{\mathscr{D}_{m,el,bet}}+\sqrt{\mathscr{D}_{m+1,energy,bet}})\\
&\times (\sqrt{\mathscr{E}_{0,el,int}}+\sqrt{\mathscr{E}_{0,el,bet}}) \lesssim\sqrt{\mathcal{E}\mathcal{D}}.
\end{aligned}
\end{equation}

For the term $I_{4}$, we have
\begin{equation}
\begin{aligned}
&I_{4}\lesssim  \|\chi_{m,\mathrm{int}}\partial_{t}\bar{\eta}\|_{L^{\infty}}\|\chi_{m+,\mathrm{int}}\partial_{3}\partial_{t}^{m}v\|_{L^{2}} \lesssim  (\|\chi_{int,\Sigma}\partial_{t}{\eta}\|_{W^{3-\frac{1}{p},p}}+\|(\chi_{cl,\Sigma}+\chi_{bet,\Sigma})\partial_{t}\eta\|_{L_{||}^{\infty}W^{3-\frac{1}{q_{\delta}},q_{\delta}}})\\
&\times (\|\chi_{int}\partial_{t}^{m}v\|_{H^{1}}+\|\chi_{bet}\partial_{t}^{m}v\|_{L_{||}^{2}H^{1}}+\|\chi_{bet}\partial_{\theta}\partial_{t}^{m}v\|_{L_{||}^{2}L^{2}})\\
&\lesssim (\sqrt{\mathscr{E}_{1,el,int}}+\sqrt{\mathscr{E}_{1,el,bet}}+\sqrt{\mathscr{E}_{1,el,cl}})\times (\sqrt{\mathscr{D}_{m,el,int}}+\sqrt{\mathscr{E}_{m,el,bet}}+\sqrt{\mathscr{E}_{m,el,cl}}) \lesssim \sqrt{\mathcal{E}}\sqrt{\mathcal{D}}.
\end{aligned}
\end{equation}

For the term $I_{5}$, using the previous computations and estimates, we have
\begin{equation}
I_{5}\lesssim  \sum_{j=0}^{m}\sum_{a+b+c+d=j,(a,b,c,d)\neq m} \|\chi_{m+,\mathrm{int}}\partial_{t}^{a+1}\bar{\eta}\|_{L^{\frac{1}{\alpha}}}\|\chi_{m+,\mathrm{int}}\partial_{t}^{b}\nabla \bar{\eta}\|_{L^{\infty}}\|\chi_{m+,\mathrm{int}}\partial_{t}^{c}\nabla \bar{\eta}\|_{L^{\infty}}\|\partial_{t}^{d}\nabla v\|_{L^{2}} \lesssim \sqrt{\mathcal{E}}\sqrt{\mathcal{D}}.
\end{equation}
Therefore,
\begin{align}
\sum_{j=0}^{m}\|\chi_{int}\partial_{t}^{j}(\partial_{t}\bar{\eta}WK\partial_{3}v)\|_{L^{p}}\lesssim \sqrt{\mathcal{E}\mathcal{D}}.
\end{align}

\textbf{Cutoff and nonlinear forcing terms in the divergence equation}

We aim to show the boundedness of
\begin{align}
\sum_{j=0}^{m}\|[\func{div}_{\mathcal{A}},\chi_{m,\mathrm{int}}]v_{j}\|_{W^{1,p}}+\sum_{j_{1}+j_{2}=m,j_{2}\neq m}\|\chi_{m,\mathrm{int}}\func{div}_{\partial_{t}^{j_{1}}\mathcal{A}}\partial_{t}^{j_{2}}v\|_{W^{1,p}}\lesssim\frac{1}{\varepsilon}\sqrt{\mathcal{D}_{||}}+\sqrt{\mathcal{E}\mathcal{D}}.
\end{align}
The estimate of this term is the same as the estimate for the stress tensor and cutoff forcing in the momentum equation. Thus, we omit the details here.

\textbf{Cutoff commutator in the capillary equation}

By definition, we have
\begin{equation}\label{eq:commutator_capillary_el}
\begin{aligned}
&\|\nabla\chi_{m,\mathrm{int},\Sigma}\partial_{t}^{m}v\mathcal{N}_{0}\|_{W^{1-\frac{1}{p},p}(\Sigma)} +\|([\func{div}_{\Sigma}\mathcal{L}_{\Sigma_{0}},\chi_{m,\mathrm{int},\Sigma}]\partial_{t}^{m}\eta)\mathcal{N}_{0}\|_{W^{1-\frac{1}{p},p}(\Sigma)} +\|[\func{div}_{\Sigma},\chi_{j,\mathrm{int},\Sigma}]F_{0}^{5,m}\|_{W^{1-\frac{1}{p},p}}\\
&\lesssim  \frac{C(\varepsilon)}{\varepsilon}\|\partial_{t}^{m}v\|_{H^{1}}+\frac{C(\varepsilon)}{\varepsilon}\|\partial_{t}^{m}\eta\|_{H^{\frac{3}{2}-\alpha}(\Sigma)} +\sum_{\ell=2}^{m+2}\sum_{i_{1}+\ldots i_{\ell}=m}\prod_{j=1}^{\ell}(\|\chi_{m,\mathrm{int},\Sigma}\partial_{t}^{i_{j}}\eta\|_{W^{3-\frac{1}{p},p}}\\
&+\|(\chi_{\mathrm{bet},\Sigma}+\chi_{\mathrm{cl},\Sigma})\partial_{t}^{i_{j}}\eta\|_{L^{2}W^{3-\frac{1}{q_{\delta}},q_{\delta}}} +\|(\chi_{\mathrm{bet},\Sigma}+\chi_{\mathrm{cl},\Sigma})\partial_{t}^{i_{j}}\eta\|_{H^{\frac{3}{2}-\alpha}}) \lesssim  \frac{C(\varepsilon)}{\varepsilon}(\sqrt{\mathcal{D}_{||}}+\sqrt{\mathcal{D}_{\mathrm{enhanced}}})+\sqrt{\mathcal{E}}\sqrt{\mathcal{D}}.
\end{aligned}
\end{equation}
\textbf{Nonlinear pressure term in the capillary equation}

We now estimate the term
\begin{align}
\sum_{a+b+c=j,b\neq j}\|(\chi_{m,\mathrm{int},\Sigma}S_{\partial_{t}^{a}\mathcal{A}_{\eta}}(\partial_{t}^{b}p,\partial_{t}^{b}v)\partial_{t}^{c}\mathcal{N})\|_{W^{1-\frac{1}{p},p}}.
\end{align}
As in the previous computation, we split the estimate of this term into four parts:
\begin{equation}
\begin{aligned}
&\sum_{a+b+c=m,b\neq j}\|(\chi_{m,\mathrm{int},\Sigma}S_{\p_{t}^{a}\mathcal{A}_{\eta}}(\partial_{t}^{b}p,\partial_{t}^{b}v)\partial_{t}^{c}\mathcal{N})\|_{W^{1-\frac{1}{p},p}}= \|(\chi_{m,\mathrm{int},\Sigma}\mathbb{D}_{\p_{t}^{m}\mathcal{A}_{\eta}}(v)\mathcal{N})\|_{W^{1-\frac{1}{p},p}}\\
&+\|(\chi_{m,\mathrm{int},\Sigma}S_{\mathcal{A}_{\eta}}(p,v)\partial_{t}^{m}\mathcal{N})\|_{W^{1-\frac{1}{p},p}} +\sum_{a+b+c=m,(a,b,c\neq m)}\|(\chi_{m,\mathrm{int},\Sigma}S_{\p_{t}^{a}\mathcal{A}_{\eta}}(\partial_{t}^{b}v)\partial_{t}^{c}\mathcal{N})\|_{W^{1-\frac{1}{p},p}}\\
&+\sum_{a+b=m,(a,b,c\neq m)}\|(\chi_{m,\mathrm{int},\Sigma}(\partial_{t}^{a}p)\partial_{t}^{b}\mathcal{N})\|_{W^{1-\frac{1}{p},p}} = I_{1}+I_{2}+I_{3}+I_{4}.
\end{aligned}
\end{equation}

Assume that $m\ge1$ and $1<q_\delta<p<2, \qquad 0<\varepsilon_-<1, \qquad 0<\beta_{1}<1-\frac1p,$ with
\begin{equation*}
0<\alpha< \min\left\{ 2-\frac3p,\, 1-\varepsilon_-,\, 1+\varepsilon_--\frac2p-2\beta_{1},\, \frac{2(q_\delta-1) \left(1-\frac1p-\beta_1\right)} {1+q_\delta\beta_{1}} \right\}.
\end{equation*}
The cutoffs are understood to satisfy the support conditions required by the localized product and composition estimates. All implicit constants may depend on the cutoffs and fixed bounds for the undifferentiated geometric coefficients and the inverse Jacobian.

For the term $I_1$, the two-dimensional interior estimate uses $W^{\frac2p+\alpha,p}$, whereas the mixed estimates use the one-dimensional multiplier space $W^{\frac1p+\beta,p}$. Consequently,
\begin{equation}
\begin{aligned}
&I_1 \lesssim{}  \|\chi_{m-1,\mathrm{int},\Sigma} \mathbb D_{\partial_t^m\mathcal A_\eta}(v) \mathcal N\|_{W^{1-\frac1p,p}} + \|(\chi_{m,\mathrm{int},\Sigma} -\chi_{m-1,\mathrm{int},\Sigma}) \mathbb D_{\partial_t^m\mathcal A_\eta}(v) \mathcal N\|_{W^{1-\frac1p,p}}\\
&\lesssim{}  \|\chi_{m-1,\mathrm{int},\Sigma} \partial_t^m\mathcal A_\eta\|_ {W^{\frac2p+\alpha,p}} \|\chi_{m-1+,\mathrm{int},\Sigma} \nabla v\|_{W^{1-\frac1p,p}} \|\chi_{m-1+,\mathrm{int},\Sigma} \mathcal N\|_{W^{\frac2p+\alpha,p}}\\
&+ \|\chi_{\mathrm{bet},\Sigma} \partial_t^m\mathcal A_\eta\|_ {L_\theta^\infty W_r^{\frac1p+\beta_{1},p}} \|(\chi_{m,\mathrm{int},\Sigma} -\chi_{m-1,\mathrm{int},\Sigma}) \nabla v\|_ {L_\theta^pW_r^{1-\frac1p,p}} \|\chi_{\mathrm{bet},\Sigma} \mathcal N\|_ {L_\theta^\infty W_r^{\frac1p+\beta_{1},p}}\\
&+ \|\chi_{\mathrm{bet},\Sigma} \partial_t^m\mathcal A_\eta\|_ {W_\theta^{\frac1p+\beta_{1},p}L_r^\infty} \|(\chi_{m,\mathrm{int},\Sigma} -\chi_{m-1,\mathrm{int},\Sigma}) \nabla v\|_ {W_\theta^{1-\frac1p,p}L_r^p} \|\chi_{\mathrm{bet},\Sigma} \mathcal N\|_ {W_\theta^{\frac1p+\beta_{1},p}L_r^\infty}.
\end{aligned}
\end{equation}

Accounting for the nonlinear dependence of $\mathcal{A}_\eta$ on $\eta$, the localized tame composition estimates yield
\begin{equation}
\begin{aligned}
&I_1 \lesssim{}  \sqrt{\mathcal{E}} \sum_{\ell=1}^{m} \ \sum_{\substack{ a_1+\cdots+a_\ell=m \\ a_j\ge1}} \ \sum_{i=1}^{\ell} \Bigg( \|\chi_{m,\mathrm{int},\Sigma} \partial_t^{a_i}\eta\|_ {W^{3-\frac1p,p}} + \|\chi_{\mathrm{bet},\Sigma} \partial_t^{a_i}\eta\|_ {L_\theta^\infty H_r^{\frac32+\frac{\varepsilon_--\alpha}{2}}}\\
&+ \|\chi_{\mathrm{bet},\Sigma} \partial_t^{a_i}\eta\|_ {L_\theta^2 W_r^{3-\frac1{q_\delta},q_\delta}} + \|\partial_\theta (\chi_{\mathrm{bet},\Sigma} \partial_t^{a_i}\eta)\|_ {H^{\frac32-\alpha}} \Bigg) \times \prod_{\substack{ 1\le j\le\ell \\ j\ne i}} \Bigg( \|\chi_{m,\mathrm{int},\Sigma} \partial_t^{a_j}\eta\|_{W^{1,\infty}} + \|\chi_{\mathrm{bet},\Sigma} \partial_t^{a_j}\eta\|_{W^{1,\infty}} \Bigg).
\end{aligned}
\end{equation}
Here the product is interpreted as one when $\ell=1$. For $\ell\ge2$, every time index satisfies $a_j\le m-1$. Provided that the higher spatial norms in the first parenthesis are controlled by $\sqrt{\mathcal{D}}$, the lower time-derivative $W^{1,\infty}$ norms are controlled by $\sqrt{\mathcal{E}}$, and the localized velocity norms are controlled by $\sqrt{\mathcal{E}}$, we obtain
\begin{align}
I_1 &\lesssim \sqrt{\mathcal{E}\mathcal{D}} \sum_{\ell=1}^{m} \mathcal{E}^{\frac{\ell-1}{2}} \lesssim \sqrt{\mathcal{E}\mathcal{D}}, \qquad \mathcal{E}\le1.
\end{align}

The embeddings used above are
\begin{align}
\|f\|_{W^{1-\frac1p,p}(\Sigma)} \lesssim{}& \|f\|_{L_\theta^pW_r^{1-\frac1p,p}} + \|f\|_{W_\theta^{1-\frac1p,p}L_r^p},
\end{align}
and
\begin{equation}
H_r^{\frac32+\frac{\varepsilon_--\alpha}{2}}  \hookrightarrow W_r^{1+\frac1p+\beta_{1},p},\quad H^{\frac32-\alpha}(\Sigma)  \hookrightarrow L_\theta^\infty W_r^{\frac1p+\beta_{1},p} \cap W_\theta^{\frac1p+\beta_{1},p}L_r^\infty.
\end{equation}
Moreover, interpolation gives
\begin{align}
\|\partial_t^k\eta\|_ {W_\theta^{\frac1p+\beta_{1},p}W_r^{1,\infty}} \lesssim{}& \|\partial_t^k\eta\|_ {L_\theta^{q_\delta} W_r^{3-\frac1{q_\delta},q_\delta}} + \|\partial_t^k\eta\|_ {W_\theta^{1,p}H_r^{\frac32-\alpha}},
\end{align}
provided that $\alpha< \frac{2(q_\delta-1) \left(1-\frac1p-\beta_{1}\right)} {1+q_\delta\beta_{1}}.$ For the term $I_{2}$, just like the estimate for $I_{1}$, we have $I_{2}\lesssim \sqrt{\mathcal{D}\mathcal{E}},$ where we used the previous estimate for $\|v\|_{W^{2,p}}$ and $\|p\|_{W^{1,p}}$.

Finally, for the term $I_{3}$, we have the following estimate. Similarly to $I_{1}$, we have
\begin{equation}
\begin{aligned}
&I_{3} \lesssim \sum_{j=0}^{m}\sum_{a+b+c=m,(a,b,c)\neq m}\|(\chi_{m-1,\mathrm{int},\Sigma}\mathbb{D}_{\p_{t}^{a}\mathcal{A}_{\eta}}(\partial_{t}^{b}v)\partial_{t}^{c}\mathcal{N})\|_{W^{1-\frac{1}{p},p}}+\|((\chi_{m,\mathrm{int},\Sigma}-\chi_{m-1,\mathrm{int},\Sigma})\mathbb{D}_{\p_{t}^{a}\mathcal{A}_{\eta}}(\partial_{t}^{b}v)\partial_{t}^{c}\mathcal{N})\|_{W^{1-\frac{1}{p},p}}\\
&\lesssim \sum_{j=0}^{m}\sum_{a+b+c=m(a,b,c)\neq m}\big(\big\|\chi_{m-1,int,\Sigma}\partial_{t}^{a}\mathcal{A}_{\eta}\|_{W^{\frac{2}{p}+\beta_{1},p}}\|\chi_{m-1+,int,\Sigma}(\nabla \partial_{t}^{b}v)\|_{W^{1-\frac{1}{p},p}}\|\chi_{m-1+,int,\Sigma}\partial_{t}^{c}\mathcal{N}\|_{W^{\frac{2}{p}+\beta_{1},p}}\\
&+\|\chi_{\mathrm{bet},int,\Sigma}\partial_{t}^{a}\mathcal{A}_{\eta}\|_{L^{\infty}W^{\frac{1}{p}+\beta_{1},p}}\|(\chi_{m,\mathrm{int},\Sigma}-\chi_{m-1+,\mathrm{int},\Sigma})(\nabla \partial_{t}^{b}v)\|_{L^{p}W^{1-\frac{1}{p},p}}\|\chi_{\mathrm{bet},\Sigma}\partial_{t}^{v}\mathcal{N}\|_{L^{\infty}W^{\frac{1}{p}+\beta_{1},p}}\\
&+\|\chi_{\mathrm{bet},int,\Sigma}\partial_{t}^{a}\mathcal{A}_{\eta}\|_{W^{\frac{1}{p}+\beta_{1},p}L^{\infty}}\|(\chi_{m,\mathrm{int},\Sigma}-\chi_{m-1+,\mathrm{int},\Sigma})(\nabla \partial_{t}^{b}v)\|_{W^{1-\frac{1}{p},p}L^{p}}\|\chi_{\mathrm{bet},\Sigma}\partial_{t}^{c}\mathcal{N}\|_{W^{\frac{1}{p}+\beta_{1},p}L^{\infty}}\big) \lesssim \sqrt{\mathcal{E}\mathcal{D}}.
\end{aligned}
\end{equation}

For the term $I_{4}$, similarly, we have $I_{4} \lesssim \sqrt{\mathcal{E}\mathcal{D}}.$ In conclusion
\begin{align}
\begin{aligned}
&\sum_{a+b+c=j,b\neq j}\|(\chi_{m,\mathrm{int},\Sigma}S_{\partial_{t}^{a}\mathcal{A}_{\eta}}(\partial_{t}^{b}p,\partial_{t}^{b}v)\partial_{t}^{c}\mathcal{N})\|_{W^{1-\frac{1}{p},p}}\lesssim \sqrt{\mathcal{E}\mathcal{D}}. \end{aligned}
\end{align}

\textbf{Curvature term in the capillary equation}

\begin{align}
\begin{aligned} \sum_{j=0}^{m}\sum_{a+b=j,a,\neq j}\|\chi_{m,\mathrm{int},\Sigma}((g\partial_{t}^{a}\eta-\sigma\mathcal{K}(\partial_{t}^{a}\eta))\partial_{t}^{b}\mathcal{N})\|_{W^{1-\frac{1}{p},p}},
\end{aligned}
\end{align}
where $\mathcal{K}(\eta):=\func{div}_{\Sigma}\mathcal{L}_{\Sigma,0}\eta$. This estimate follows directly from the estimate for the previous nonlinear stress tensor term and the fact that $\partial_{t}^{j}(\chi_{m,\mathrm{int},\Sigma}\eta)\in W^{3-\frac{1}{p},p}$. We have
\begin{align}
\sum_{j=0}^{m}\sum_{a+b=j,a,\neq j}\|\chi_{m,\mathrm{int}}((g\partial_{t}^{a}\eta-\sigma\mathcal{K}(\partial_{t}^{a}\eta))\partial_{t}^{b}\mathcal{N})\|_{W^{1-\frac{1}{p},p}}\lesssim \sqrt{\mathcal{E}\mathcal{D}}.
\end{align}

\textbf{Residual nonlinear term}

Write $F(z)=z(1+|z|^{2})^{-1/2}$ and $R(z,w)=F(z+w)-F(z)-DF(z)w$, so that $\mathcal{R}(\nabla\eta)=R(\nabla\zeta,\nabla\eta)$. The spatial chain rule is
\begin{align}
\partial_i\mathcal{R}(\nabla\eta) =D_wR(\nabla\zeta,\nabla\eta)\partial_i\nabla\eta +D_zR(\nabla\zeta,\nabla\eta)\partial_i\nabla\zeta.
\end{align}
Thus we estimate the term
\begin{equation}
\begin{aligned}
&\sum_{j=0}^{m}\|\chi_{m,int,\Sigma}\partial_{t}^{j}(\mathcal{R}(\nabla \eta))\|_{W^{2-\frac{1}{p},p}} \lesssim\sum_{j=0}^{m}\|\partial_{t}^{j}(D_wR(\nabla\zeta,\nabla\eta)\chi_{m,int,\Sigma}\nabla^{2}\eta)\|_{W^{1-\frac{1}{p},p}}\\
&+\sum_{j=0}^{m}\|\partial_{t}^{j}(D_zR(\nabla\zeta,\nabla\eta)\chi_{m,int,\Sigma}\nabla^{2}\zeta)\|_{W^{1-\frac{1}{p},p}} +\sum_{j=0}^{m}\|\partial_{t}^{j}(\mathcal{R}(\nabla \eta)\nabla\chi_{m,int,\Sigma})\|_{W^{1-\frac{1}{p},p}}.
\end{aligned}
\end{equation}
Here $D_wR(z,w)=O(|w|)$ and $D_zR(z,w)=O(|w|^{2})$ for $w$ small. The derivatives of the fixed smooth surface $\zeta$ are included in the geometric constant.

Since $D_wR(\nabla\zeta,\nabla\eta)=O(|\nabla\eta|)$, we apply a computation similar to that for the curvature nonlinear term and the nonlinear stress tensor term to $\sum_{j=0}^{m}\|\partial_{t}^{j}(D_wR(\nabla\zeta,\nabla\eta)\chi_{m,int,\Sigma}\nabla^{2}\eta)\|_{W^{1-\frac{1}{p},p}}$ to show that
\begin{align}
\sum_{j=0}^{m}\|\partial_{t}^{j}(D_wR(\nabla\zeta,\nabla\eta)\chi_{m,int,\Sigma}\nabla^{2}\eta)\|_{W^{1-\frac{1}{p},p}}\lesssim\sqrt{\mathcal{E}\mathcal{D}}.
\end{align}
The second term contains the fixed factor $\nabla^{2}\zeta$ and is quadratic in $\nabla\eta$, so it is bounded similarly by $\sqrt{\mathcal{E}\mathcal{D}}$. For the third term, using an estimate similar to that for the cutoff commutator term derived in estimate \eqref{eq:commutator_capillary_el}, we have
\begin{align}
\sum_{j=0}^{m}\|\partial_{t}^{j}(\mathcal{R}(\nabla \eta)\nabla\chi_{m,int,\Sigma})\|_{W^{1-\frac{1}{p},p}}\lesssim\sqrt{\mathcal{E}}\sqrt{\mathcal{D}}.
\end{align}
The difference here is that the curvature residue term is nonlinear. Therefore, we have
\begin{equation}
\begin{aligned}
&\sum_{j=0}^{m}\|\chi_{m,int,\Sigma}\partial_{t}^{j}(\mathcal{R}(\nabla \eta))\|_{W^{2-\frac{1}{p},p}} \lesssim\sum_{j=0}^{m}\|\partial_{t}^{j}(D_wR(\nabla\zeta,\nabla\eta)\chi_{m,int,\Sigma}\nabla^{2}\eta)\|_{W^{1-\frac{1}{p},p}}\\
&+\sum_{j=0}^{m}\|\partial_{t}^{j}(D_zR(\nabla\zeta,\nabla\eta)\chi_{m,int,\Sigma}\nabla^{2}\zeta)\|_{W^{1-\frac{1}{p},p}} +\sum_{j=0}^{m}\|\partial_{t}^{j}(\mathcal{R}(\nabla \eta)\nabla\chi_{m,int,\Sigma})\|_{W^{1-\frac{1}{p},p}}\lesssim C(\varepsilon) \sqrt{\mathcal{E}\mathcal{D}},
\end{aligned}
\end{equation}
and
\begin{align}
\sum_{j=0}^{m}\|\partial_{t}^{j}(\mathcal{R}(\nabla \eta)\chi_{m,int,\Sigma})\|_{W^{1-\frac{1}{p},p}}\lesssim \sqrt{\mathcal{E}}\sqrt{\mathcal{D}}.
\end{align}
We omit the details here.

\textbf{Nonlinear term in the kinematic boundary condition}

For the kinematic boundary condition, there is only one nonlinear term which can be bounded as follows:
\begin{equation}
\sum_{j=0}^{m}\sum_{a+b=j,b\geq 1}(\|\partial_{t}^{a}\chi_{m,int,\Sigma}\nabla v \partial_{t}^{b}\chi_{m+,int,\Sigma}\mathcal{N}\|_{W^{1-\frac{1}{p},p}}+\|\chi_{m+,int,\Sigma}\partial_{t}^{a}v\cdot \chi_{m,int,\Sigma}\partial_{t}^{b}\nabla \mathcal{N}\|_{W^{1-\frac{1}{p},p}}) = I_{1}+I_{2}.
\end{equation}
For the term $I_{1}$, we apply a computation similar to that for the estimate of $\|\partial_{t}^{j}S_{\mathcal{A}}(p,v)\mathcal{N}\|_{W^{1-\frac{1}{p},p}}$ to show that $I_{1}\lesssim\sqrt{\mathcal{E}\mathcal{D}}.$ Assume that $1<q_\delta<p<2, \qquad 0<\varepsilon_-<1, \qquad 0<\beta^{\prime }< \frac12+\frac{\varepsilon_-}{2}-\frac1p.$ On a fixed bounded product chart with coordinates $(r,\theta)$, we have
\begin{equation}\label{eq:intrinsic_slicing_norm}
\begin{aligned}
&\|\chi_{\mathrm{bet},\Sigma}f\|_{W_\theta^{\frac1p+\beta^{\prime },p}L_r^\infty} + \|\chi_{\mathrm{bet},\Sigma}f\|_{L_\theta^\infty W_r^{\frac1p+\beta^{\prime },p}} \lesssim \|\chi_{\mathrm{bet},\Sigma}f\|_{L_\theta^\infty H_r^{\frac12+\frac{\varepsilon_-}{2}}} + \|\chi_{\mathrm{bet},\Sigma}f\|_{W_\theta^{1,2}H_r^{\frac12}}\\
&\le \|\chi_{\mathrm{bet},\Sigma}f\|_{L_\theta^\infty H_r^{\frac12+\frac{\varepsilon_-}{2}}} + \|\chi_{\mathrm{bet},\Sigma}f\|_{L_\theta^2W_r^{2-\frac1{q_\delta},q_\delta}} + \|\chi_{\mathrm{bet},\Sigma}f\|_{W_\theta^{1,2}H_r^{\frac12}}.
\end{aligned}
\end{equation}
If $f=\partial_{t}^{j}v$ for $0\leq j\leq m-1$, then, by the trace theorem,
\begin{equation*}
\|\chi_{\mathrm{bet},\Sigma}f\|_{L_\theta^\infty H_r^{\frac12+\frac{\varepsilon_-}{2}}} + \|\chi_{\mathrm{bet},\Sigma}f\|_{L_\theta^2W_r^{2-\frac1{q_\delta},q_\delta}} + \|\chi_{\mathrm{bet},\Sigma}f\|_{W_\theta^{1,2}H_r^{\frac12}}\lesssim \sqrt{\mathcal{E}}.
\end{equation*}

Indeed, since $\frac1p+\beta^{\prime }< \frac12+\frac{\varepsilon_-}{2} \qquad\text{and}\qquad p<2,$ the one-dimensional Sobolev embedding yields $H_r^{\frac12+\frac{\varepsilon_-}{2}} \hookrightarrow W_r^{\frac1p+\beta^{\prime },p}.$ Consequently,
\begin{align}
\|\chi_{\mathrm{bet},\Sigma}f\|_{L_\theta^\infty W_r^{\frac1p+\beta^{\prime },p}} \lesssim \|\chi_{\mathrm{bet},\Sigma}f\|_{L_\theta^\infty H_r^{\frac12+\frac{\varepsilon_-}{2}}}.
\end{align}

For the other mixed norm, Hilbert-space interpolation gives
\begin{equation}
\|\chi_{\mathrm{bet},\Sigma}f\|_{ H_\theta^{\frac{1+\frac1p+\beta^{\prime }}{2}} H_r^{\frac12+ \frac{\varepsilon_-}{4} \left(1-\frac1p-\beta^{\prime }\right)} } \lesssim \|\chi_{\mathrm{bet},\Sigma}f\|_{L_\theta^2H_r^{\frac12+\frac{\varepsilon_-}{2}}}^{ \frac{1-\frac1p-\beta^{\prime }}{2}} \|\chi_{\mathrm{bet},\Sigma}f\|_{W_\theta^{1,2}H_r^{\frac12}}^{ \frac{1+\frac1p+\beta^{\prime }}{2}}.
\end{equation}
The parameter assumptions imply $\frac1p+\beta^{\prime }<1,$ and hence
\begin{equation*}
\frac{1+\frac1p+\beta^{\prime }}{2}> \frac1p+\beta^{\prime }, \qquad \frac12+ \frac{\varepsilon_-}{4} \left(1-\frac1p-\beta^{\prime }\right)> \frac12.
\end{equation*}
Therefore,
\begin{equation*}
H_\theta^{\frac{1+\frac1p+\beta^{\prime }}{2}} H_r^{\frac12+ \frac{\varepsilon_-}{4} \left(1-\frac1p-\beta^{\prime }\right)} \hookrightarrow W_\theta^{\frac1p+\beta^{\prime },p}L_r^\infty.
\end{equation*}
Using the boundedness of the $\theta$-interval and Young's inequality, we obtain
\begin{equation}
\|\chi_{\mathrm{bet},\Sigma}f\|_{W_\theta^{\frac1p+\beta^{\prime },p}L_r^\infty}  \lesssim \|\chi_{\mathrm{bet},\Sigma}f\|_{L_\theta^\infty H_r^{\frac12+\frac{\varepsilon_-}{2}}}^{ \frac{1-\frac1p-\beta^{\prime }}{2}} \|\chi_{\mathrm{bet},\Sigma}f\|_{W_\theta^{1,2}H_r^{\frac12}}^{ \frac{1+\frac1p+\beta^{\prime }}{2}} \lesssim \|\chi_{\mathrm{bet},\Sigma}f\|_{L_\theta^\infty H_r^{\frac12+\frac{\varepsilon_-}{2}}} + \|\chi_{\mathrm{bet},\Sigma}f\|_{W_\theta^{1,2}H_r^{\frac12}},
\end{equation}
which proves the desired estimate.

Applying the estimates above to $I_{2}$, we have
\begin{equation}
\begin{aligned}
&I_{2}\lesssim  \sum_{j=1}^{m}\sum_{a+b=j,b\geq 1}\big(\|\chi_{m,\mathrm{int},\Sigma}\partial_{t}^{a}v\|_{W^{2-\frac{1}{p},p}}\|\partial_{t}^{b}\chi_{m,\mathrm{int},\Sigma}\eta\|_{W^{3-\frac{1}{p},p}}+\|\chi_{\mathrm{bet},\Sigma}\partial_{t}^{a}v\|_{L^{\infty}W^{\frac{1}{p}+\beta^{\prime },p}}\|\partial_{t}^{b}\chi_{m,\mathrm{int},\Sigma}\eta\|_{W^{3-\frac{1}{p},p}}\\
&+\|\chi_{\mathrm{bet},\Sigma}\partial_{t}^{a}v\|_{W^{\frac{1}{p}+\beta^{\prime },p}L^{\infty}}\|\partial_{t}^{b}\chi_{m,\mathrm{int},\Sigma}\eta\|_{W^{3-\frac{1}{p},p}}\big) \lesssim  \sqrt{\mathcal{E}\mathcal{D}}.
\end{aligned}
\end{equation}

\textbf{Nonlinear terms on the fixed surface}

For $F_{\mathrm{int}}^{6,j}$, the cutoff commutator $-\mu U_{0}(v_{j}\otimes\nabla_{\mathcal{A}}\chi_{j,\mathrm{int}}+\nabla_{\mathcal{A}}\chi_{j,\mathrm{int}}\otimes v_{j})n_{0}$ is linear in $v_{j}$ and is bounded as the capillary commutator above, while the nonlinear part $\chi_{j,\mathrm{int}}F^{6,j}_{0}$ is bounded as $F^{3,j}_{\mathrm{int}}$. Hence
\begin{align}
\|F_{\mathrm{int}}^{6,j}\|_{W^{1-\frac{1}{p},p}(\Sigma_{s})}\lesssim\frac{1}{\varepsilon}\big(\sqrt{\mathcal{D}_{||}}+\sqrt{\mathcal{D}_{\mathrm{enhanced}}}\big)+\sqrt{\mathcal{E}\mathcal{D}}.
\end{align}

\textbf{Constant Variable forcing term $F_{s,\mathrm{int}}^{1,j}$}

It suffices to discuss the case $j=m$. When $j\leq m$, we proceed in the same manner. It is actually even simpler for $j<m$.

Using an estimate similar to that for $F_{\mathrm{int}}^{1,j}$, we have
\begin{align}
\begin{aligned} \|F_{s,\mathrm{int}}^{1,m}\| &\lesssim\|\chi_{m+,\mathrm{int}}(I-\mathcal{A}_{\eta})\|_{L^{\infty}}\|\chi_{m+,\mathrm{int}}{\mathcal{A_{\eta}}}\|_{L^{\infty}}\|(\chi_{m,int}\nabla^{2}\partial_{t}^{m}v)\|_{L^{p}}+\|\chi_{m,int}\nabla\partial_{t}^{m}p\|_{L^{p}}\|\chi_{m+,\mathrm{int}}(I-\mathcal{A}_{\eta})\|_{L^{\infty}}\\
&\quad+\|\chi_{m+,\mathrm{int}}(I-\mathcal{A}_{\eta})\|_{L^{\infty}}\|\chi_{m,\mathrm{int}}{\nabla \mathcal{A_{\eta}}}(\chi_{m+,\mathrm{int}}\nabla\partial_{t}^{m}v)\|_{L^{p}}+\|\chi_{m+,int}\mathcal{A}_{\eta}\|_{L^{\infty}}\|\chi_{m,\mathrm{int}}\nabla{\mathcal{A_{\eta}}}(\chi_{m+,\mathrm{int}}\nabla \partial_{t}^{m}v)\|_{L^{p}}\\ &\qquad=I_{1}+I_{2}+I_{3}+I_{4}. \end{aligned}
\end{align}
For $I_{1}$ and $I_{2}$, using the estimate for $\|\mathcal{A}_{\eta}\|_{L^{\infty}}$ and $\|I-\mathcal{A}_{\eta}\|_{L^{\infty}}$, we have
\begin{align}
I_{1}+I_{2}\lesssim\sqrt{\mathcal{E}}(\|v_{m,\mathrm{int}}\|_{W^{2,p}}+\|p_{m,\mathrm{int}}\|_{W^{1,p}}+\|\partial_{t}^{m}v_{\mathrm{bet}}\|_{H^{1}}+\|\partial_{t}^{m}p_{\mathrm{bet}}\|_{L^{2}}).
\end{align}
For $I_{3}$ and $I_{4}$ we split $\chi_{m,\mathrm{int}}=\chi_{m-1,\mathrm{int}}+\psi$, with $\psi:=\chi_{m,\mathrm{int}}-\chi_{m-1,\mathrm{int}}$. On $\func{supp}\chi_{m-1,\mathrm{int}}$, the bound \eqref{eq:A_d} with $b=0$ controls $\|\chi_{m-1,\mathrm{int}}\nabla\mathcal{A}_{\eta}\|_{L^{3}}$ by norms of $\eta$ without time derivatives. These are at the energy level; for $\|\chi_{m,\mathrm{int},\Sigma}\eta\|_{W^{3-\frac{1}{p},p}}$ this holds because $\mathcal{E}_{\mathrm{int}}$ contains it. Hence $\|\chi_{m-1,\mathrm{int}}\nabla\mathcal{A}_{\eta}\|_{L^{3}}\lesssim\sqrt{\mathcal{E}}$. On $\func{supp}\psi$, which lies in $\{\chi_{\mathrm{bet}}=1\}$, we use H\"older's inequality with $\sigma:=\frac{4p}{4-2p+p\varepsilon_{-}},\qquad\text{so that}\qquad \frac{1}{p}=\frac{1}{\sigma}+\frac{2-\varepsilon_{-}}{4},$ together with the embedding $H^{\varepsilon_{-}/2}\subset L^{\frac{4}{2-\varepsilon_{-}}}$ on the two-dimensional slices. Both exponents are at least $p$, since $\frac{4}{2-\varepsilon_{-}}>2>p$ and $\sigma\geq p$. We obtain
\begin{equation}
I_{3}+I_{4}\lesssim  \|\chi_{m-1,\mathrm{int}}\nabla\mathcal{A}_{\eta}\|_{L^{3}}\|\chi_{m,\mathrm{int}}v_{m}\|_{W^{2,p}} +\|\psi\nabla \mathcal{A}_{\eta}\|_{L^{p}_{\theta}L^{\sigma}_{rh}}\|\chi_{\mathrm{bet}}\nabla\partial_{t}^{m}v\|_{L^{\infty}_{\theta}L^{\frac{4}{2-\varepsilon_{-}}}_{rh}},
\end{equation}
and
\begin{equation}
\begin{aligned}
&\|\chi_{\mathrm{bet}}\nabla\partial_{t}^{m}v\|_{L^{\infty}_{\theta}L^{\frac{4}{2-\varepsilon_{-}}}_{rh}} \leq C(\varepsilon)\Bigl(\|(\chi_{\mathrm{cl}}+\chi_{\mathrm{bet}})\partial_{t}^{m}v\|_{L^{\infty}_{\theta}H_{rh}^{1+\frac{\varepsilon_{-}}{2}}}\\
&+\|(\chi_{\mathrm{cl}}+\chi_{\mathrm{bet}})\partial_{t}^{m}\partial_{\theta}v\|_{H^{1}} +\|\chi_{\mathrm{int}}\partial_{t}^{m}v\|_{W^{2,p}}\Bigr) \leq C(\varepsilon)\sqrt{\mathcal{D}}.
\end{aligned}
\end{equation}
Here $H^{1}$ is the full three-dimensional norm. Since $0<\varepsilon_{-}<1$, choose $1/2<s<1-\varepsilon_{-}/2$. Then $H^{1}(\mathbb{T}\times\mathcal{W})\hookrightarrow H_{\theta}^{s}H_{rh}^{\varepsilon_{-}/2}\hookrightarrow L_{\theta}^{\infty}H_{rh}^{\varepsilon_{-}/2}$, and the angular term is controlled by the dissipation at derivative rank $(m,1)$.

For the shell factor we apply Theorem~\ref{thm:enhanced_commutator} as in the estimate of $I_{3,2}$. On $\func{supp}\psi$ the extension $\bar\eta$ equals $\bar{E}_{\rho}f$, with $f=(\chi_{m+1,\mathrm{int},\Sigma}-\chi_{\mathrm{int},\Sigma})\eta$, up to far-field terms that are smooth there, and
\begin{equation*}
\psi\bar{E}_{\rho}f=\big(\psi\bar{E}_{\rho}f-\bar{E}_{\rho}(\psi_{\Sigma}f)\big)+\bar{E}_{\rho}(\psi_{\Sigma}f),\qquad \psi_{\Sigma}f=(\chi_{m,\mathrm{int},\Sigma}-\chi_{m-1,\mathrm{int},\Sigma})\eta .
\end{equation*}
The commutator lies in $L^{2}_{\theta}W^{4,q_{\delta}}$ by \eqref{eq:enhanced_commutator}, so its second derivatives lie in $L^{2}_{\theta}L^{\infty}\subset L^{p}_{\theta}L^{\sigma}$ on the slices. The second derivatives of $\bar{E}_{\rho}(\psi_{\Sigma}f)$ lie in $L^{p}_{\theta}W^{1,p}_{rh}\subset L^{p}_{\theta}L^{\frac{2p}{2-p}}_{rh}\subset L^{p}_{\theta}L^{\sigma}_{rh}$, because $\sigma\leq\frac{2p}{2-p}$, which is equivalent to $p\varepsilon_{-}\geq0$. Hence
\begin{align}
\|\psi\nabla \mathcal{A}_{\eta}\|_{L^{p}_{\theta}L^{\sigma}_{rh}}\leq C(\varepsilon)\big(\|\chi_{\mathrm{bet},\Sigma}\partial_{\theta}\eta\|_{H^{\frac{3}{2}-\alpha}}+\|\chi_{m,\mathrm{int},\Sigma}\eta\|_{W^{3-\frac{1}{p},p}}+\|\chi_{\mathrm{bet,\Sigma}}\eta\|_{L^{2}W^{3-\frac{1}{q_{\delta}},q_{\delta}}}\big),  \notag
\end{align}
and all three norms on the right are at the energy level. Therefore
\begin{equation}
\begin{aligned}
&I_{3}+I_{4}\lesssim  \sqrt{\mathcal{E}}\,\|\chi_{m,\mathrm{int}}v_{m}\|_{W^{2,p}}\\
&+C(\varepsilon)\big(\|\chi_{\mathrm{bet},\Sigma}\partial_{\theta}\eta\|_{H^{\frac{3}{2}-\alpha}}+\|\chi_{m,\mathrm{int},\Sigma}\eta\|_{W^{3-\frac{1}{p},p}}+\|\chi_{\mathrm{bet,\Sigma}}\eta\|_{L^{2}W^{3-\frac{1}{q_{\delta}},q_{\delta}}}\big)\|\chi_{\mathrm{bet}}\nabla\partial_{t}^{m}v\|_{L^{\infty}_{\theta}L_{rh}^{\frac{4}{2-\varepsilon_{-}}}} \lesssim \sqrt{\mathcal{E}\mathcal{D}}.
\end{aligned}
\end{equation}

\textbf{Constant Variable forcing term $F_{s,\mathrm{int}}^{2,j}$}

The estimate follows from the estimate for $F_{s,\mathrm{int}}^{1,j}$. We write
\begin{equation}
\|F_{s,\mathrm{int}}^{2,m}\|_{W^{1,p}}\lesssim  \|\chi_{m+,\mathrm{int}}(\mathcal{A}_{\eta}-I)\|_{L^{\infty}}\|\chi_{m,\mathrm{int}}\nabla^{2}v_{m}\|_{L^{p}}+\|\chi_{m,\mathrm{int}}\nabla(\mathcal{A}_{\eta}-I)\chi_{m+,\mathrm{int}}\nabla v_{m}\|_{L^{p}} \lesssim \sqrt{\mathcal{E}\mathcal{D}}.
\end{equation}
The first term is bounded as $I_{1}$. The second term has the same structure as $I_{3}+I_{4}$ and is bounded in the same way, with the shell factor in $L^{p}_{\theta}L^{\sigma}_{rh}$ and the velocity factor in $L^{\infty}_{\theta}L^{\frac{4}{2-\varepsilon_{-}}}_{rh}$.

\textbf{Constant Variable forcing term $F_{s,\mathrm{int}}^{3,j}$}

Using the same estimate as for $F_{\mathrm{int}}^{3,m}$, we have
\begin{align}
\|F_{s,\mathrm{int}}^{3,m}\|_{W^{1-\frac{1}{p},p}}\lesssim (\|\chi_{m,\mathrm{int},\Sigma}\partial_{t}^{m}\eta\|_{W^{3-\frac{1}{p},p}}+\|\chi_{m,\mathrm{int}}\partial_{t}^{m}v\|_{W^{2,p}}+\|\chi_{m,\mathrm{int}}\partial_{t}^{m}p\|_{W^{1,p}}+\sqrt{\mathcal{E}\mathcal{D}})\sqrt{\mathcal{E}}.
\end{align}
The only difference is that we take the $m$-th-order derivative of $p,v$ and $\mathcal{K}(\eta)$ with no time derivative on the geometric terms.

\textbf{Constant Variable forcing term $F_{s,\mathrm{int}}^{6,j}$}

Since the term included in $F_{s,\mathrm{int}}^{6,j}$ has a form similar to the stress tensor difference term included in $F_{s,\mathrm{int}}^{3,j}$, its estimate follows directly from the estimate for $F_{\mathrm{int}}^{3,j}$.

\textbf{Constant Variable forcing term $F_{s,\mathrm{int}}^{4,m}$}

By definition, and the estimate for $F_{\mathrm{int}}^{4,m}$, we have
\begin{equation}
\begin{aligned}
&\|F_{s,\mathrm{int}}^{4,m}\|_{W^{2-\frac{1}{p},p}}   \lesssim (\|\nabla(\chi_{m,int,\Sigma}\partial_{t}^{m}v) \chi_{m+,int,\Sigma}(\mathcal{N}-\mathcal{N}_{0})\|_{W^{1-\frac{1}{p},p}}\\
&+\|\chi_{m+,int,\Sigma}\partial_{t}^{m}v\cdot \chi_{m,int,\Sigma}\nabla (\mathcal{N}-\mathcal{N}_{0})\|_{W^{1-\frac{1}{p},p}} +\|\chi_{m,int,\Sigma}\partial_{t}^{m}v\cdot(\mathcal{N}-\mathcal{N}_{0})\|_{W^{1-\frac{1}{p},p}})\\
&\lesssim \sqrt{\mathcal{E}}(\|\chi_{m,\mathrm{int}}\partial_{t}^{m}v\|_{W^{2,p}}+\|\chi_{\mathrm{bet}}\partial_{t}^{m}v\|_{L^{\infty}H^{1+\frac{\varepsilon_{-}}{2}}} +\|\chi_{\mathrm{bet}}\partial_{t}^{m}v\|_{L^{2}W^{2,q_{\delta}}}+\|\chi_{bet}\partial_{t}^{m}\partial_{\theta}v\|_{H^{1}})   \lesssim\sqrt{\mathcal{E}\mathcal{D}}.
\end{aligned}
\end{equation}

This finishes the proof of \eqref{est:basic_diss_int}. The estimate \eqref{est:basic_energy_int} follows by reducing the temporal derivative by one order. Then any dissipation $\mathcal{D}$ can be replaced by $\mathcal{E}$ in the estimates above.
\end{proof}

Now we have the basic elliptic estimate. We then proceed to derive the estimates for other terms included in the interior elliptic energy functional.

\subsection{Enhanced elliptic estimate}

To establish the interior estimate, we employ a bootstrap argument. At each stage, we shrink the interior domain and obtain higher spatial regularity at the expense of one level of temporal differentiation. To carry out this procedure, we repeatedly apply commutator estimates analogous to those established in the previous section. We construct the following sequence of cutoff functions. Suppose $\chi_{j,int}$ for $m\geq j\geq 1$ are cutoff functions defined as in Section 2. Then we have the following theorem.

\begin{theorem}
Define
\begin{align}
\mathscr{D}_{m,\mathrm{basic},\mathrm{int}}:=\sum_{k=0}^{m}\Bigl(\|\chi_{m,\mathrm{int}}\partial_{t}^{k}v\|_{W^{2,p}}^{2}+\|\chi_{m,\mathrm{int}}\partial_{t}^{k}p\|_{W^{1,p}}^{2}+\|\chi_{m,\mathrm{int},\Sigma}\partial_{t}^{k}\eta\|_{W^{3-\frac{1}{p},p}}^{2}\Bigr),
\end{align}
and
\begin{align}
\mathscr{E}_{m,\mathrm{basic},\mathrm{int}}:=\sum_{k=0}^{m-1}\Bigl(\|\chi_{m,\mathrm{int}}\partial_{t}^{k}v\|_{W^{2,p}}^{2}+\|\chi_{m,\mathrm{int}}\partial_{t}^{k}p\|_{W^{1,p}}^{2}+\|\chi_{m,\mathrm{int},\Sigma}\partial_{t}^{k}\eta\|_{W^{3-\frac{1}{p},p}}^{2}\Bigr).
\end{align}

Suppose $(\chi_{j,\mathrm{int}}\p_{t}^{j}v,\chi_{j,\mathrm{int}}\p_{t}^{j}p,\chi_{j,\mathrm{int},\Sigma}\p_{t}^{j}\eta)$ is the solution to the equation system \eqref{eq:interior-domain-elliptic-system-j} for any $0\leq j\leq m-1$. Then the following properties hold for $j\leq m-1$.
\begin{align}
{\label{eq:enhanced_regularity}} \mathscr{D}_{j,high,int} \leq C_m(\varepsilon)\bigl(\mathscr{D}_{m,\mathrm{basic},int} +\mathcal{E}\mathcal{D}\bigr),
\end{align}
where
\begin{align}
\mathscr{D}_{j,\mathrm{high},\mathrm{int}} := \sum_{k=0}^{j} \Big( &\|\chi_{j,\mathrm{int}} \partial_{t}^{k}v\|^{2}_{W^{2+m-j,p}} + \|\chi_{j,\mathrm{int}}\partial_t^k p\|^{2}_{W^{1+m-j,p}} + \|\chi_{j,\mathrm{int},\Sigma}\partial_t^k\eta\|^{2} _{W^{3-\frac{1}{p}+(m-j),p}} \Big).
\end{align}

\begin{align}
{\label{eq:enhanced_regularity_e}} \mathscr{E}_{j,high,int} \leq C_m(\varepsilon)\bigl(\mathscr{E}_{m,basic,int}+\mathcal{E}^{2}\bigr),
\end{align}
and
\begin{align}
\mathscr{E}_{j,\mathrm{high},\mathrm{int}} := \sum_{k=0}^{j} \Big( &\|\chi_{j+1,\mathrm{int}} \partial_{t}^{k}v\|^{2}_{W^{1+m-j,p}} + \|\chi_{j+1,\mathrm{int}}\partial_t^k p\|^{2}_{W^{m-j,p}} + \|\chi_{j+1,\mathrm{int},\Sigma}\partial_t^k\eta\|^{2} _{W^{2-\frac{1}{p}+(m-j),p}} \Big),
\end{align}
for $j\geq 0$.
\end{theorem}

\begin{proof}
We only prove estimate \eqref{eq:enhanced_regularity}. Estimate \eqref{eq:enhanced_regularity_e} follows in a similar way by reducing the temporal derivative by one order. We prove the theorem by a recurrence argument. Suppose $\mathscr{D}_{j,\mathrm{high}\mathrm{int}}$ is bounded by the right-hand side of the inequality \eqref{eq:enhanced_regularity} for $j\geq n$. We then prove that $\mathscr{D}_{j,\mathrm{high}\mathrm{int}}$ for $j=n-1$ is also bounded by the right-hand side of \eqref{eq:enhanced_regularity}. By Corollary~\ref{cor:enhanced elliptic estimate}, it suffices to bound the forcing terms by the right-hand side of \eqref{eq:enhanced_regularity} by setting $\ell=m-j$.

\textbf{Step 0.} We first study a norm that appears in many places in the following estimates. $\|\chi_{j,\mathrm{int}}\partial_{t}^{j}\mathcal{A}\|_{W^{m-j+1,p}}.$ By Theorem \ref{thm:extension-commutator-estimate}, when $j\geq 1$,  we consider the following bound
\begin{equation}
\begin{aligned}
&\sum_{r=1}^{j}\sum_{\substack{i_{1}+\ldots+i_{r}=j\\
i_{k}\geq 1}}\prod_{k=1}^{r} \big(\|\partial_{t}^{i_{k}}\bar{\eta}_{j+1,\mathrm{int},d}\|_{W^{m-j+2,p}}+\|\partial_{t}^{i_{k}}\bar{\eta}_{\mathrm{int}}\|_{W^{m-j+2,p}} \big)\\
&\lesssim \sum_{r=1}^{j}\sum_{\substack{i_{1}+\ldots+i_{r}=j\\
i_{k}\geq 1}}\prod_{k=1}^{r}\big(\|\partial_{t}^{i_{k}}{\eta}_{\mathrm{int}}\|_{W^{m-j+2-\frac{1}{p},p}} +\sum_{a=0}^{m-j+2}\|\partial_{t}^{i_{k}}\bar{\eta}_{j+1,\mathrm{int},d}\|_{W_{\theta}^{a,p}W_{r,h}^{m-j+2-a,p}}  \big)\\
&\lesssim\sum_{r=1}^{j}\sum_{\substack{i_{1}+\ldots+i_{r}=j\\
i_{k}\geq 1}}\prod_{k=1}^{r}\big( \|\partial_{t}^{i_{k}}{\eta}_{\mathrm{int}}\|_{W^{m-j+2-\frac{1}{p},p}} +\sum_{a=0}^{m-j+1}\|\partial_{t}^{i_{k}}\bar{\eta}_{j+1,\mathrm{int},d}\|_{W_{\theta}^{a,p}W_{r,h}^{m-j+2-a,p}}+\|\partial_{t}^{i_{k}}\bar{\eta}_{j+1,\mathrm{int},d}\|_{W_{\theta}^{m-j+2,p}L^{p}} \big)\\
&\lesssim\sum_{r=1}^{j}\sum_{\substack{i_{1}+\ldots+i_{r}=j\\
i_{k}\geq 1}}\prod_{k=1}^{r}\big(\|\partial_{t}^{i_{k}}{\eta}_{\mathrm{int}}\|_{W^{m-j+2-\frac{1}{p},p}} +\sum_{a=0}^{m-j+1}\|\partial_{t}^{i_{k}}{\eta}_{j+1,\mathrm{int}}\|_{W^{a,p}W_{r,h}^{m-j+2-a-\frac{1}{p},p}}+\|\partial_{\theta}^{m-j+2}\partial_{t}^{i_{k}}{\eta}_{\mathrm{bet}}\|_{L^{p}}\big)\\
&\lesssim\sum_{r=1}^{j}\sum_{\substack{i_{1}+\ldots+i_{r}=j\\
i_{k}\geq 1}}\prod_{k=1}^{r}\big(\|\partial_{t}^{i_{k}}{\eta}_{\mathrm{int}}\|_{W^{m-j+2-\frac{1}{p},p}} +\|\partial_{t}^{i_{k}}{\eta}_{j+1,\mathrm{int}}\|_{W^{m-j+2-\frac{1}{p},p}}+\|\partial_{\theta}^{m-j+1}\partial_{t}^{i_{k}}{\eta}_{\mathrm{bet}}\|_{H^{1}} \big)\\
&\lesssim\|\partial_{t}^{j}{\eta}_{\mathrm{int}}\|_{W^{m-j+2-\frac{1}{p},p}} +\|\partial_{t}^{j}{\eta}_{j+1,\mathrm{int}}\|_{W^{m-j+2-\frac{1}{p},p}}+\|\partial_{\theta}^{m-j+1}\partial_{t}^{j}{\eta}_{\mathrm{bet}}\|_{H^{1}} +\sqrt{\mathcal{E}}\sqrt{\mathcal{D}} \lesssim \sqrt{\mathcal{D}}.
\end{aligned}
\end{equation}
Fix $0\leq j\leq m-1$ and set $s=m-j$. We estimate the geometric coefficients on $\operatorname{supp}\chi_{j+,\mathrm{int}}$. Only the wider surface cutoff $\chi_{j+1,\mathrm{int},\Sigma}$ is needed. For $N=s+2$ write
\begin{equation}
B_{i,N}:=\|\chi_{j+1,\mathrm{int},\Sigma}\partial_t^i\eta\|_{W^{N-1/p,p}} +\|(\chi_{\mathrm{cl},\Sigma}+\chi_{\mathrm{bet},\Sigma}) \partial_\theta^{N-1}\partial_t^i\eta\|_{H^1}.
\end{equation}
From the computation above we obtain
\begin{equation}\label{eq:enh-int-linear-extension} \|\partial_t^i\bar\eta_{j+1,\mathrm{int}}\|_{W^{N,p}} \le C_m(\varepsilon) B_{i,N}.
\end{equation}
For $0\leq i\leq j-1$, $B_{i,N}$ is bounded by $C(\varepsilon)\sqrt{\mathcal{E}}$ and when $i=j$, it is bounded by $C(\varepsilon)\sqrt{\mathcal{D}}$.

By \eqref{eq:A-explicit}, write $\mathcal A_\eta=F(\nabla q)$, $q=W\bar\eta$, where $F(0)=I$. In computing the coefficient on the support in question use $q=W\bar\eta_{j+1,\mathrm{int}}$. Outside the admissible range $|\nabla q|\leq 1/2$ we may extend $F$ smoothly with bounded derivatives before estimating this localized expression. Since $s\geq1$ and $p>3/2$, $W^{s+1,p}$ is an algebra in three dimensions. The spatial chain rule and Sobolev composition bound $D^rF(\nabla q)$ in this algebra by $C_m(\varepsilon)$: by \eqref{eq:enh-int-linear-extension} the $W^{s+1,p}$ norm of $\nabla q$ is bounded by $C_m(\varepsilon)B_{0,s+2}$ and the energy is small. Thus for $1\leq k\leq j$ the time chain rule gives
\begin{equation}\label{eq:enh-int-coefficient-products} \|\chi_{j+,\mathrm{int}}\partial_t^k\mathcal A_\eta\|_{W^{s+1,p}} \le C_m(\varepsilon)\sum_{r=1}^k \sum_{\substack{i_1+\cdots+i_r=k\\i_1,\ldots,i_r\geq1}} \prod_{d=1}^r B_{i_d,s+2}.
\end{equation}
At $k=j\geq1$ use dissipation for one factor and energy for the other $r-1$ factors. The term $r=1$ is $O(C_m(\varepsilon)\sqrt{\mathcal D})$; the sum with $r\geq2$ is $O(C_m(\varepsilon)\sqrt{\mathcal E\mathcal D})$. Alternatively all factors are bounded by energy for every $1\leq k\leq j$. At $k=0$ the formula $F(z)-I=\int_0^1DF(\tau z)z\,d\tau$ bounds the localized $\mathcal A_\eta-I$ by $C_m(\varepsilon)\sqrt{\mathcal E}$ also when $j=0$. The identity in $\mathcal A_\eta$ itself is retained. The same proof applies to $K$ with equilibrium value $1$. All these bounds hold also with the smaller cutoff $\chi_{j,\mathrm{int}}$. In subsequent uses in this enhanced interior estimate the cutoff constants are included in $C_m(\varepsilon)$ in \eqref{eq:enhanced_regularity}. Similarly,
\begin{align}
\sum_{k=1}^{j-1}\|\chi_{j,\mathrm{int}}\partial_{t}^{k}\mathcal{A}\|_{W^{m-j+1,p}}+\|\chi_{j,\mathrm{int}}(I-\mathcal{A})\|_{W^{m-j+1,p}}\lesssim \sqrt{\mathcal{E}}.
\end{align}

\textbf{Step 1.} When $j=m-1$, the norm $W^{1,p}$ is not an algebra in $3-D$. In this step, we first derive the elliptic boundedness in this special case. In this case $\ell=m-(m-1)=1$ We now establish the estimates for terms on the right-hand side of \eqref{eq:enhance_el_estimate}.

\textbf{Linear term estimates}

For estimates of linear terms on the right-hand side of \eqref{eq:enhance_el_estimate} which are not included in forcing terms, we have
\begin{align}
\|\eta_{m-1,\mathrm{int}}\|_{W^{3-\frac{1}{p},p}}+ \|\p_{t}\eta_{m-1,\mathrm{int}}\|_{W^{3-\frac{1}{p},p}}&\lesssim C(\varepsilon)(\|\eta_{m,\mathrm{int}}\|_{W^{3-\frac{1}{p},p}}+\|\chi_{m,\mathrm{int},\Sigma}\p_{t}^{m-1}\eta\|_{W^{3-\frac{1}{p},p}})\notag\\ &\lesssim C(\varepsilon)\sqrt{\mathscr{D}_{m,\mathrm{basic},\mathrm{int}}}. \\
\|\p_{t}v_{m-1,\mathrm{int}}\|_{W^{1,p}}+\|v_{m-1,\mathrm{int}}\|_{W^{2,p}}&\lesssim C(\varepsilon) (\|v_{m,\mathrm{int}}\|_{H^{1}}+\|\p_{t}^{m-1}\chi_{m-1,\mathrm{int}}v\|_{W^{2,p}})\lesssim C(\varepsilon)\sqrt{\mathscr{D}_{m,\mathrm{basic},\mathrm{int}}}. \\
\|p_{m-1,\mathrm{int}}\|_{W^{1,p}}&\lesssim C(\varepsilon)\|\p_{t}^{m-1}\chi_{m,\mathrm{int}}p\|_{W^{1,p}}\lesssim C(\varepsilon)\sqrt{\mathscr{D}_{m,\mathrm{basic},\mathrm{int}}}.
\end{align}
\textbf{Cutoff commutators}

In Step 0, setting $j=m-1$, we have
\begin{equation}\label{est:bet-cut-momentum}
\begin{aligned}
&\left\| \partial_t^{m-1} \left( S(p,v)\nabla\chi_{\mathrm{j,int}} +v\Delta \chi_{m-1,int} \right) \right\|_{W^{1,p}} \lesssim \frac{C(\varepsilon)}{\varepsilon}(\|S(\partial_{t}^{m-1}(\chi_{m,int})p,\partial_{t}^{j}(\chi_{m,int})v)\|_{W^{1,p}}\\
&+\|(\chi_{m,int})\partial_{t}^{m-1}v\|_{W^{1,p}}) \lesssim \frac{C(\varepsilon)}{\varepsilon}(\|\chi_{m,int}\partial_{t}^{m-1}v\|_{W^{2,p}} +\|\chi_{m,int}\partial_{t}^{m-1}p\|_{W^{1,p}}) \lesssim\frac{C(\varepsilon)}{\varepsilon}\sqrt{\mathscr{D}_{m,basic,int}}.
\end{aligned}
\end{equation}

\textbf{Stress nonlinear terms in the momentum equation}

We now estimate nonlinear terms induced by $\partial_{t}^{m-1}(\func{div}_{\mathcal{A}}S_{\mathcal{A}}(\chi_{m-1,int}v,\chi_{m-1,int}p))$ or $\partial_{t}^{m-1}(\func{div}_{\mathcal{A}}S_{\chi_{m-1int}\mathcal{A}}(v,p)).$ This is the most delicate term, where we must use the commutator defined before. Since the commutators have already been discussed above, in the following computation, we can commute the cutoff function and the spatial derivatives. By Step 0, we have
\begin{equation}
\begin{aligned}
&\|\chi_{m-1,int}\partial_{t}^{m-1}(\mathop{\rm div}\nolimits_{\mathcal{A}_{\eta}}S_{\mathcal{A}_{\eta}}(p,v))-\chi_{m-1,int}\operatorname{div}_{\mathcal{A_{\eta}}}S_{\mathcal{A}_{\eta}}(\partial_{t}^{m-1}v,\partial_{t}^{m-1}p)\|_{W^{1,p}}\\
&\lesssim \sum_{a+b+c=m-1,c\neq m-1}\|\chi_{m-1,int}(\mathop{\rm div}\nolimits_{\p_{t}^{a}\mathcal{A}_{\eta}}\nabla_{\p_{t}^{b}\mathcal{A}_{\eta}}\partial_{t}^{c} v)\|_{W^{1,p}}+\sum_{b+c=m-1,c\neq m-1}\|\chi_{m-1,int}(\nabla_{\p_{t}^{b}\mathcal{A}_{\eta}}\partial_{t}^{c} p)\|_{W^{1,p}}\\
&\lesssim \sum_{a+b+c=m-1,c\neq m-1}\|\chi_{m-1+,int}\partial_{t}^{a}\mathcal{A}\|_{W^{2,p}}\|\chi_{m-1+,int}\partial_{t}^{b}\mathcal{A}\|_{W^{2,p}}(\|\chi_{m-1,int} \partial_{t}^{c}v\|_{W^{3,p}}+\|\chi_{m,int} \partial_{t}^{c}v\|_{W^{2,p}})\\
&+\sum_{a+b+c=m-1,c\neq m-1}\|\chi_{m-1+,int}\partial_{t}^{a}\mathcal{A}\|_{W^{2,p}}\|\chi_{m-1+,int}\partial_{t}^{b}\nabla \mathcal{A}\|_{W^{1,p}}(\|\chi_{m-1,int} \partial_{t}^{c}v\|_{W^{3,p}}+\|\chi_{m,int} \partial_{t}^{c}v\|_{W^{2,p}})\\
&+\sum_{b+c=m-1,c\neq m-1}\|\chi_{m-1+,int}\partial_{t}^{b}\mathcal{A}\|_{W^{2,p}}(\|\chi_{m-1,int} \partial_{t}^{c}p\|_{W^{2,p}}+\|\chi_{m,int} \partial_{t}^{c}p\|_{W^{1,p}}) \lesssim \sqrt{\mathcal{E}\mathcal{D}},
\end{aligned}
\end{equation}
where we used the fact that $b,c$ cannot both equal $m-1$ given that $b+c=m-1$ and that $m\geq 3$. Similarly, for terms included in $F_{s,\mathrm{int}}^{1,m-1}$, we have
\begin{equation}
\begin{aligned}
&\|\chi_{m-1,\mathrm{int}} \operatorname{div}_{I-\mathcal{A_{\eta}}}S_{\mathcal{A}_{\eta}}(\partial_{t}^{m-1}v,\partial_{t}^{m-1}p)\|_{W^{1,p}} +\|\chi_{m-1,\mathrm{int}} \operatorname{div}\mathbb{D}_{I-\mathcal{A}_{\eta}}(\partial_{t}^{m-1}v)\|_{W^{1,p}}\\
&\lesssim \|\chi_{m-1+,int}(\mathcal{A}-I)\|_{W^{2,p}}(\|\chi_{m-1,int} \partial_{t}^{m-1}v\|_{W^{3,p}} +\|\chi_{m,int} \partial_{t}^{m-1}v\|_{W^{2,p}})\\
&+\|\chi_{m-1+,int}\mathcal{A}\|_{W^{2,p}}\|\chi_{m-1+,int}\nabla \mathcal{A}\|_{W^{1,p}}(\|\chi_{m-1,int} \partial_{t}^{m-1}v\|_{W^{3,p}} +\|\chi_{m,int} \partial_{t}^{m-1}v\|_{W^{2,p}})\\
&+\|\chi_{m-1+,int}(\mathcal{A}-I)\|_{W^{2,p}}(\|\chi_{m-1,int} \p_{t}^{m-1}p\|_{W^{2,p}} +\|\chi_{m,int} \p_{t}^{m-1}p\|_{W^{1,p}}) \lesssim (\|p_{m-1,\mathrm{int}}\|_{W^{2,p}}\\
&+\|v_{m-1,\mathrm{int}}\|_{W^{3,p}} +\|\chi_{m,int} \p_{t}^{m-1}p\|_{W^{1,p}} +\|\chi_{m,\mathrm{int}}\p_{t}^{m-1}v\|_{W^{2,p}})(\|\chi_{m-1+,int}(\mathcal{A}-I)\|_{W^{2,p}}) \lesssim \sqrt{\mathcal{E}\mathcal{D}}.
\end{aligned}
\end{equation}
\textbf{Nonlinear convection term}

We now estimate $\|\partial_{t}^{m-1}(\chi_{m-1int}v\cdot \nabla_{\mathcal{A}}v)\|_{W^{1,p}}.$ We give the calculation for all $0\leq j\leq m-1$, with $s=m-j$; here $j=m-1$. Multiplication maps $W^{s+1,p}\times W^{s,p}$ into $W^{s,p}$. The cutoff $\chi_{j+1,\mathrm{int}}$ equals one near $\operatorname{supp}\chi_{j,\mathrm{int}}$, so the gradient may be taken on $\chi_{j+1,\mathrm{int}}\partial_t^c v$. By Step~0,
\begin{align}\label{eq:enh-int-convection-products} &\|\chi_{j,\mathrm{int}}\partial_t^j(v\cdot\nabla_{\mathcal A}v)\|_{W^{s,p}} \notag\\
&\quad\le C_m(\varepsilon)\sum_{a+b+c=j} \|\chi_{j,\mathrm{int}}\partial_t^a v\|_{W^{s+1,p}} \|\chi_{j+,\mathrm{int}}\partial_t^b\mathcal A\|_{W^{s+1,p}} \|\chi_{j+1,\mathrm{int}}\partial_t^c v\|_{W^{s+1,p}} \le C_m(\varepsilon)\sqrt{\mathcal E\mathcal D}.
\end{align}
Indeed the two velocity factors are bounded by energy at level $j$ and dissipation at level $j+1$, respectively. The middle factor is bounded: for $b>0$ use all the time partitions in \eqref{eq:enh-int-coefficient-products}, and for $b=0$ keep also the identity.

\textbf{Nonlinear temporal derivative term}

We now estimate $\|\partial_{t}^{m-1}\chi_{int}(\partial_{t}\bar{\eta }WK\partial_{3}v)\|_{W^{1,p}}.$ Again take $0\leq j\leq m-1$ and $s=m-j$. The function $W$ is fixed, so time derivatives fall only on $\partial_t\bar\eta$, $K$, and $v$. With the same nested velocity cutoff as above,
\begin{align}\label{eq:enh-int-temporal-products} &\|\chi_{j,\mathrm{int}}\partial_t^j(\partial_t\bar\eta\,WK\partial_3v)\|_{W^{s,p}} \notag\\
&\quad\le C_m(\varepsilon)\sum_{a+c+d=j} \|\chi_{j+,\mathrm{int}}\partial_t^{a+1}\bar\eta\|_{W^{s+1,p}} \|\chi_{j+,\mathrm{int}}\partial_t^cK\|_{W^{s+1,p}} \|\chi_{j+1,\mathrm{int}}\partial_t^d v\|_{W^{s+1,p}} \le C_m(\varepsilon)\sqrt{\mathcal E\mathcal D}.
\end{align}
Use \eqref{eq:enh-int-linear-extension} with $N=s+1$, $i=a+1\leq j+1$, so the first factor is bounded by dissipation. The middle factor is bounded by the product formula of Step~0 for $K$, including its constant term when $c=0$. The last factor is bounded by energy at level $j$.

\textbf{Nonlinear terms in the divergence equation}

Since $W^{2,p}$ is an algebra in 3-D, and again by Step 0, we have
\begin{equation}
\begin{aligned}
&\sum_{a+b=m-1,b\neq m-1}\| \chi_{m-1,int}(\operatorname{div}_{\p_{t}^{a}\mathcal{A}_{\eta}}\partial_{t}^{b}v)\|_{W^{2,p}} \lesssim \sum_{a+b=m-1,b\neq m-1}\| \nabla \chi_{m-1,int}({\partial_{t}^{a}\mathcal{A}_{\eta}}\partial_{t}^{b}v)\|_{W^{2,p}}\\
&+\sum_{a+b=m-1,b\neq m-1}\| (\operatorname{div}_{\p_{t}^{a}\mathcal{A}_{\eta}}\partial_{t}^{b}(\chi_{m-1,int}v))\|_{W^{2,p}} \lesssim \sum_{a+b=m-1,b\neq m-1}\|\chi_{m-1+,int}{\partial_{t}^{a}\mathcal{A}_{\eta}}\|_{W^{2,p}}\|\chi_{m,int}\partial_{t}^{b}v\|_{W^{2,p}}\\
&+\sum_{a+b=m-1,b\neq m-1}\| \chi_{m-1+,int}{\partial_{t}^{a}\mathcal{A}_{\eta}}\|_{W^{2,p}}\|\partial_{t}^{b}(\chi_{m-1,int}v))\|_{W^{3,p}}\lesssim\sqrt{\mathcal{E}\mathcal{D}}.
\end{aligned}
\end{equation}
Similarly,
\begin{equation}
\begin{aligned}
&\| \chi_{m-1,int}(\operatorname{div}_{I-\mathcal{A}_{\eta}}\partial_{t}^{m-1}v)\|_{W^{2,p}} \lesssim \| \nabla \chi_{m-1,int}({(I-\mathcal{A}_{\eta}})\partial_{t}^{m-1}v)\|_{W^{2,p}}+\| \operatorname{div}_{I-\mathcal{A}_{\eta}}(\chi_{m-1,int}\partial_{t}^{m-1}v)\|_{W^{2,p}}\\
&\lesssim \|\chi_{m-1+,int}(I-{\mathcal{A}_{\eta}})\|_{W^{2,p}}(\|\chi_{m,int}\partial_{t}^{m-1}v\|_{W^{2,p}}+\|\chi_{m-1,int}\partial_{t}^{m-1}v\|_{W^{3,p}}) \lesssim \sqrt{\mathcal{E}\mathcal{D}}.
\end{aligned}
\end{equation}
\textbf{Cutoff commutator in the capillary equation}

We have the following estimate for the cutoff commutator
\begin{align}
&\|[\mathbb{D},\chi_{m-1,\mathrm{int},\Sigma}]\partial_{t}^{m-1}v\|_{W^{2-\frac{1}{p},p}}+\|[\func{div}_{\Sigma}\mathcal{L}_{\Sigma,0},\chi_{m-1,\mathrm{int},\Sigma}]\partial_{t}^{m-1}\eta\|_{W^{2-\frac{1}{p},p}} \\
&\lesssim \frac{C(\varepsilon)}{\sqrt{\varepsilon}}(\|\chi_{m,\mathrm{int}}\partial_{t}^{m-1}v\|_{W^{2,p}}+\|\chi_{m,\mathrm{int}}\partial_{t}^{m-1}\eta\|_{W^{3-\frac{1}{p},p}}) \lesssim \frac{C(\varepsilon)}{\sqrt{\varepsilon}}\sqrt{\mathscr{D}_{m,\mathrm{basic},\mathrm{int}}}.
\end{align}

\textbf{Stress tensor term in the capillary equation}

We now estimate the following term
\begin{equation}
\begin{aligned}
&\sum_{a+b+c=m-1,b\neq m-1}\|(\chi_{m-1,int,\Sigma})\mathbb{D}_{\p_{t}^{a}\mathcal{A}}(\partial_{t}^{b}v)\partial_{t}^{c}\mathcal{N}\|_{W^{2-\frac{1}{p},p}}+\sum_{b+c=m-1,b\neq m-1}\|(\chi_{m-1,int,\Sigma})(\partial_{t}^{b}p I)\partial_{t}^{c}\mathcal{N}\|_{W^{2-\frac{1}{p},p}}\\
&\lesssim \sum_{a+b+c=m-1,b\neq m-1}\|\chi_{m-1+,int,\Sigma}\partial_{t}^{a}\mathcal{A}\|_{W^{2-\frac{1}{p},p}}(\|\chi_{m-1,int,\Sigma}\partial_{t}^{b}v\|_{W^{3-\frac{1}{p},p}})\|\chi_{m-1+,int,\Sigma}\partial_{t}^{c}\mathcal{N}\|_{W^{2-\frac{1}{p},p}}\\
&+\sum_{b+c=m-1, b\neq m-1}\|\chi_{m-1,int,\Sigma}\partial_{t}^{b}p\|_{W^{2-\frac{1}{p},p}}\|\chi_{m-1+,int,\Sigma}\partial_{t}^{c}\mathcal{N}\|_{W^{2-\frac{1}{p},p}}\\
&\lesssim \sum_{\substack{a+b+c=m-1\\
a,c>0}}(\|\chi_{m-1+,int,\Sigma}\partial_{t}^{a}\mathcal{A}\|_{W^{2-\frac{1}{p},p}})(\|\chi_{m-1,int}\partial_{t}^{b}v\|_{W^{3,p}})\|\chi_{m,int,\Sigma}\partial_{t}^{c}\eta\|_{W^{3-\frac{1}{p},p}}\\
&+\sum_{\substack{b+c=m-1\\
c>0}}(\|\chi_{m-1,int}\partial_{t}^{b}p\|_{W^{2,p}})\|\chi_{m,int,\Sigma}\partial_{t}^{c}\eta\|_{W^{3-\frac{1}{p},p}}\\
&+\sum_{a_{1}+a_{2}=m-1,a_{1}\neq m-1}(\|\chi_{m-1,int}\partial_{t}^{a_{1}}v\|_{W^{3,p}})\|\chi_{m,int,\Sigma}\partial_{t}^{a_{2}}\eta\|_{W^{3-\frac{1}{p},p}} \lesssim \sqrt{\mathcal{E}\mathcal{D}}.
\end{aligned}
\end{equation}

\textbf{Curvature term in the capillary equation}

We now estimate
\begin{align}
\sum_{a+b=m-1,a\neq m-1}\|(\chi_{m-1,int,\Sigma}(g\partial_{t}^{a}\eta-\sigma \mathcal{K}(\partial_{t}^{a}\eta))\partial_{t}^{b}\mathcal{N})\|_{W^{2-\frac{1}{p},p}}.
\end{align}
Since $W^{2-\frac{1}{p},p}$ is an algebra in 2-D, we have
\begin{equation}
\begin{aligned}
&\sum_{a+b=m-1,a\neq m-1}\|(\chi_{m-1,int,\Sigma}(g\partial_{t}^{a}\eta-\sigma \mathcal{K}(\partial_{t}^{a}\eta))\partial_{t}^{b}\mathcal{N})\|_{W^{2-\frac{1}{p},p}} \lesssim \sum_{a+b=m-1,a\neq m-1}\|\chi_{m-1,int,\Sigma}\partial_{t}^{a}\eta\|_{W^{4-\frac{1}{p},p}}\\
&\|\chi_{m,int,\Sigma}\partial_{t}^{b}\eta\|_{W^{3-\frac{1}{p},p}} +\sum_{a+b=m-1,a\neq m-1}\|\chi_{m,int,\Sigma}\partial_{t}^{a}\eta\|_{W^{3-\frac{1}{p},p}} \|\chi_{m,int,\Sigma}\partial_{t}^{b}\eta\|_{W^{3-\frac{1}{p},p}} \lesssim\sqrt{\mathcal{E}\mathcal{D}}.
\end{aligned}
\end{equation}

\textbf{Fully Taylor nonlinear term in the capillary equation}

We now estimate $\sum_{j=0}^{m-1}\|\chi_{m-1,int,\Sigma}\func{div}_{\Sigma}(\partial_{t}^{j}(\mathcal{R}(\nabla \eta))))\|_{W^{2-\frac{1}{p},p}}.$ \begin{equation}
\begin{aligned}
&\sum_{j=0}^{m-1}\|\chi_{m-1,int,\Sigma}\operatorname{div}_{\Sigma}(\partial_{t}^{j}(\mathcal{R}(\nabla \eta))))\|_{W^{2-\frac{1}{p},p}} \lesssim \sum_{j=0}^{m-1}\sum_{\ell=2}^{j+2}\sum_{i_{1}+\ldots i_{\ell}=j}(\|\chi_{m-1,\mathrm{int},\Sigma}\partial_{t}^{i_{1}}\eta\|_{W^{4-\frac{1}{p},p}})\\
&\prod_{k=2}^{\ell}(\|\chi_{m,\mathrm{int},\Sigma}\partial_{t}^{i_{k}}\eta\|_{W^{3-\frac{1}{p},p}}) +\sum_{j=0}^{m}\sum_{\ell=2}^{j+2}\sum_{i_{1}+\ldots i_{\ell}=j} \prod_{k=1}^{\ell}(\|\chi_{m,\mathrm{int},\Sigma}\partial_{t}^{i_{k}}\eta\|_{W^{3-\frac{1}{p},p}}) \lesssim\sqrt{\mathcal{E}\mathcal{D}}.
\end{aligned}
\end{equation}

\textbf{The estimate for the boundary freezing forcing term $F_{s,\mathrm{int}}^{3,m-1}$}

Since $W^{2-\frac{1}{p},p}$ is an algebra, using an estimate similar to that for the terms included in $F_{\mathrm{int}}^{3,m-1}$, we have
\begin{align}
\|F_{s,\mathrm{int}}^{3,m-1}\|_{W^{2-\frac{1}{p},p}}\lesssim \sqrt{\mathcal{E}}(\|v_{m-1,\mathrm{int}}\|_{W^{3,p}}+\|p_{m-1,\mathrm{int}}\|_{W^{2,p}}+\|\eta_{m-1,\mathrm{int}}\|_{W^{4-\frac{1}{p},p}}).
\end{align}

\textbf{The nonlinear term in the kinematic boundary condition}

We now establish the estimate for $\sum_{a+b=m-1,a\neq m-1}\|\partial_{t}^{a}v\cdot \partial_{t}^{b}\mathcal{N}\|_{W^{3-\frac{1}{p},p}}.$ Since $W^{2-\frac{1}{p},p}$ is an algebra in 2-D, we have
\begin{equation}
\begin{aligned}
&\sum_{a+b=m-1,a\neq m-1}\|\chi_{m-1,int,\Sigma}\partial_{t}^{a}v\cdot \partial_{t}^{b}\mathcal{N}\|_{W^{3-\frac{1}{p},p}} \lesssim \sum_{a+b=m-1,a\neq m-1}\|\partial_{t}^{a}\nabla (\chi_{m-1,int,\Sigma}v)\cdot \partial_{t}^{b}\mathcal{N}\|_{W^{2-\frac{1}{p},p}}\\
&+\sum_{a+b=m-1,a\neq m-1}\|\partial_{t}^{a}v\cdot \partial_{t}^{b}\nabla (\chi_{m-1,int,\Sigma}\mathcal{N})\|_{W^{2-\frac{1}{p},p}} +\sum_{a+b=m-1,a\neq m-1}\|\partial_{t}^{a}v\cdot \partial_{t}^{b} (\chi_{m-1,int,\Sigma}\mathcal{N})\|_{W^{2-\frac{1}{p},p}}\\
&\lesssim \sum_{a+b=m-1,a\neq m-1}\|\partial_{t}^{a} (\chi_{m-1,int,\Sigma}v)\|_{W^{3-\frac{1}{p},p}} \|\chi_{m,int,\Sigma}\partial_{t}^{b}\mathcal{N}\|_{W^{2-\frac{1}{p},p}}\\
&+\sum_{a+b=m-1,a\neq m-1}\|\chi_{m,int,\Sigma}\partial_{t}^{a}v\|_{W^{2-\frac{1}{p},p}}\| \partial_{t}^{b} (\chi_{m-1,int,\Sigma}\eta)\|_{W^{4-\frac{1}{p},p}} \lesssim  \sqrt{\mathcal{E}\mathcal{D}}.
\end{aligned}
\end{equation}
Using a similar estimate, we have
\begin{align}
\begin{aligned} &\|\chi_{m-1,int,\Sigma}\partial_{t}^{m-1}v\cdot (\mathcal{N}-\mathcal{N}_{0})\|_{W^{3-\frac{1}{p},p}}\\
\lesssim& \|\chi_{m,int,\Sigma}\partial_{t}^{m-1}v\|_{W^{2-\frac{1}{p},p}}\| (\chi_{m-1,int,\Sigma}\eta)\|_{W^{4-\frac{1}{p},p}}+\|\chi_{m-1,int,\Sigma}\partial_{t}^{m-1}v\|_{W^{3-\frac{1}{p},p}}\| (\chi_{m,int,\Sigma}\eta)\|_{W^{3-\frac{1}{p},p}}\lesssim\sqrt{\mathcal{E}\mathcal{D}}. \end{aligned}
\end{align}

For the forcing terms on the fixed boundary, we apply a computation similar to that for the forcing terms on the capillary boundary. We omit the details here.  For $k\leq m-1$, the estimate follows the same idea. Therefore, combining all of the estimates in this step, we have
\begin{align}
\mathscr{D}_{m-1,high,int} \leq C_{m-1}(\varepsilon)\bigl(\mathscr{D}_{m,\mathrm{basic},int} +\mathcal{E}\mathcal{D}\bigr).
\end{align}

\textbf{Step 2: The case $j\leq m-1$} We use a recurrence argument. We have already shown the result when $j=m-1$. We then further suppose that the result is true for $n\leq j\leq m-2$ for some positive integer $n<m-1$. We now want to show the result when $j=n-1$. We still show the estimates for nonlinear terms as follows. In this case, $\ell=m-j\geq 2$, and all of the elliptic norms for forcing terms on the right hand side of \eqref{eq:enhance_el_estimate} are algebras in 3-D or 2-D.

\textbf{Linear estimate}

For the linear terms that are not included in the right-hand side of equation \eqref{eq:enhance_el_estimate}, applying a computation similar to that in Step 1, we have
\begin{align}
\|\eta_{j,\mathrm{int}}\|_{W^{2+(m-j)-\frac{1}{p},p}}+ \|\p_{t}\eta_{j,\mathrm{int}}\|_{W^{2+(m-j)-\frac{1}{p},p}}&\lesssim C(\varepsilon)\sqrt{\mathscr{D}_{j+1,\mathrm{high},\mathrm{int}}}. \\
\|\p_{t}v_{j,\mathrm{int}}\|_{W^{(m-j),p}}+\|v_{j,\mathrm{int}}\|_{W^{1+(m-j),p}}&\lesssim C(\varepsilon)\sqrt{\mathscr{D}_{j+1,\mathrm{high},\mathrm{int}}}. \\
\|p_{j,\mathrm{int}}\|_{W^{m-j,p}}&\lesssim\|\p_{t}^{j}\chi_{j+1,\mathrm{int}}p\|_{W^{m-j,p}}\lesssim C(\varepsilon)\sqrt{\mathscr{D}_{j+1,\mathrm{high},\mathrm{int}}}.
\end{align}

\textbf{Cutoff commutators}

Using Step 0, we have
\begin{equation}\label{eq:bet-cut-momentum-t}
\begin{aligned}
&\left\| \partial_t^{j} \left( S_{\mathcal{A}}(p,v)\nabla_{\mathcal{A}}\chi_{\mathrm{j,int}} +v\Delta_{\mathcal{A}} \chi_{j,int} \right) \right\|_{W^{m-j,p}} \lesssim \frac{C(\varepsilon)}{\varepsilon^{\frac{1}{2}}}(\|S_{\mathcal{A}}(\partial_{t}^{j}(\chi_{j,int})p,\partial_{t}^{j}(\chi_{j,int})v)\|_{W^{m-j,p}}\\
&+\|(\chi_{j+1,int})\partial_{t}^{j}v\|_{W^{m-j,p}}) +\sum_{a+b=j,(b\neq j)}(\|S_{\p_{t}^{a}\mathcal{A}}(\partial_{t}^{b}\chi_{j+1,int}p,\partial_{t}^{b}\chi_{j+1,int}v)\|_{W^{m-j,p}})\\
&+\sum_{a+b+c=j,(c\neq j)}\|(\chi_{j+1,int})\partial_{t}^{a}v\nabla \partial_{t}^{c}\bar{\eta}(\nabla^{2}\partial_{t}^{b}\bar{\eta}+\nabla\partial_{t}^{b}\bar{\eta})\|_{W^{m-j,p}}\\
&\lesssim \frac{C(\varepsilon)}{\varepsilon^{\frac{1}{2}}}(\|\chi_{j+1,int}\partial_{t}^{j}v\|_{W^{m-(j+1)+2,p}} +\|\chi_{j+1,int}\partial_{t}^{j}p\|_{W^{m-(j+1)+1,p}})\\
&+\sum_{a+b=j,(b\neq j)}(\|\chi_{j+1,int}\partial_{t}^{b}v\|_{W^{m-(j+1)+2,p}} +\|\chi_{j+1,int}\partial_{t}^{b}p\|_{W^{m-(j+1)+1,p}})\|\chi_{j+,int} \partial_{t}^{a}\bar{\eta}\|_{W^{m-(j+1)+2,p}}\\
&+\sum_{a+b+c=j,(c\neq j)}(\|\chi_{j+1,int}\partial_{t}^{a}v\|_{W^{m-(j+1)+2,p}})\|\chi_{j+,int}\nabla \partial_{t}^{c}\bar{\eta}\|_{W^{m-(j+1)+2,p}}\\
&\|\chi_{j+,int} \partial_{t}^{b}\bar{\eta}\|_{W^{m-(j+1)+3,p}} \lesssim \frac{C(\varepsilon)}{\varepsilon^{\frac{1}{2}}}(\|\chi_{j+1,int}\partial_{t}^{j}v\|_{W^{m-(j+1)+2,p}} +\|\chi_{j+1,int}\partial_{t}^{j}p\|_{W^{m-(j+1)+1,p}})+\sqrt{\mathcal{E}\mathcal{D}}\\
&\lesssim \frac{C(\varepsilon)}{\varepsilon^{\frac{1}{2}}}\sqrt{\mathscr{D}_{j+1,high,int}}+\sqrt{\mathscr{D}_{j+1,high,int}\mathscr{E}_{j+1,high,int}} \lesssim C_j(\varepsilon)\bigl(\sqrt{\mathscr{D}_{j+1,high,int}}+\sqrt{\mathcal{E}\mathcal{D}}\bigr).
\end{aligned}
\end{equation}
For the other commutators, a similar computation shows that they all are bounded by $C_j(\varepsilon)\bigl(\sqrt{\mathscr{D}_{j+1,high,int}}+\sqrt{\mathcal{E}\mathcal{D}}\bigr)$. We omit the details.

\textbf{Stress nonlinear terms in the momentum equation}

We now estimate nonlinear terms induced by $\partial_{t}^{j}(\func{div}_{\mathcal{A}}S_{\mathcal{A}}(\chi_{j,int}u,\chi_{j,int}p))$ or $\partial_{t}^{j}(\func{div}_{\mathcal{A}}S_{\chi_{j,int}\mathcal{A}}(u,p)).$ This is the most delicate term, where we must use the commutator defined before. Since the commutators have already been discussed above, in the following computation, we can commute the cutoff function and the spatial derivatives. By the commutator estimate, we have
\begin{equation}
\begin{aligned}
&\|\chi_{j,int}\partial_{t}^{j}(\mathop{\rm div}\nolimits_{\mathcal{A}_{\eta}}S_{\mathcal{A}_{\eta}}(p,v))-\chi_{j,int}\operatorname{div}_{\mathcal{A_{\eta}}}S_{\mathcal{A}_{\eta}}(\partial_{t}^{j}v,\partial_{t}^{j}p)\|_{W^{m-j,p}}\\
&\lesssim \sum_{a+b+c=j,c\neq j}\|\chi_{j,int}(\mathop{\rm div}\nolimits_{\p_{t}^{a}\mathcal{A}_{\eta}}\mathbb{D}_{\p_{t}^{b}\mathcal{A}_{\eta}}\partial_{t}^{c} v)\|_{W^{m-j,p}}+\sum_{b+c=j,c\neq j}\|\chi_{j,int}(\nabla_{\p_{t}^{b}\mathcal{A}_{\eta}}\partial_{t}^{c} p)\|_{W^{m-j,p}}\\
&\lesssim \sum_{a+b+c=j,c\neq j}\|\chi_{j+,int}\partial_{t}^{a}\mathcal{A}\|_{W^{m-j,p}}\|\chi_{j+,int}\partial_{t}^{b}\mathcal{A}\|_{W^{m-j,p}}(\|\chi_{j,int} \partial_{t}^{c}v\|_{W^{m-j+2,p}}+\|\chi_{j+1,int} \partial_{t}^{c}v\|_{W^{m-j+1,p}})\\
&+\sum_{a+b+c=j,c\neq j}\|\chi_{j+,int}\partial_{t}^{a}\mathcal{A}\|_{W^{m-j,p}}\|\chi_{j+,int}\partial_{t}^{b}\nabla \mathcal{A}\|_{W^{m-j,p}}(\|\chi_{j+1,int} \partial_{t}^{c}v\|_{W^{m-j+1,p}})\\
&+\sum_{b+c=j,c\neq j}\|\chi_{j+,int}\partial_{t}^{b}\mathcal{A}\|_{W^{m-j,p}}(\|\chi_{j,int} \partial_{t}^{c}p\|_{W^{m-j+1,p}}+\|\chi_{j+1,int} \partial_{t}^{c}p\|_{W^{m-j,p}}) \lesssim \sqrt{\mathcal{E}\mathcal{D}}.
\end{aligned}
\end{equation}
Similarly, we have
\begin{align}
\|F_{s,\mathrm{int}}^{1,j}\|_{W^{m-j,p}}\lesssim\sqrt{\mathcal{E}}(\|v_{j,\mathrm{int}}\|_{W^{2+(m-j),p}}+\|p_{j,\mathrm{int}}\|_{W^{1+(m-j),p}}).
\end{align}

\emph{Remark.} Here one sees the key point of the new commutator scheme. After taking $m-j$ derivatives, $\chi_{j}\eta$ has one order more regularity than is required for the nonlinear estimate here. In other words, we gain one more derivative for $\eta$. That is the reason why we only need to choose the cutoff and do not need to enhance the regularity of the commutators.

\textbf{Nonlinear convection term}

We now estimate $\|\partial_{t}^{j}(\chi_{j,int}v\cdot \nabla_{\mathcal{A}}v)\|_{W^{m-j,p}}.$ Apply \eqref{eq:enh-int-convection-products} with the present $j$. The bound is $C_m(\varepsilon)\sqrt{\mathcal E\mathcal D}$.

\textbf{Nonlinear temporal derivative term in the momentum equation}

We now estimate $\|\partial_{t}^{j}\chi_{j,int}(\partial_{t}\bar{\eta }WK\partial_{3}v)\|_{W^{m-j,p}}.$ Estimate \eqref{eq:enh-int-temporal-products} at this $j$ gives $C_m(\varepsilon)\sqrt{\mathcal E\mathcal D}$.

\textbf{Nonlinear terms in the divergence equation}

Since $W^{m-j+1,p}$ is an algebra in 3-D, we have
\begin{equation}
\begin{aligned}
&\sum_{a+b=j,b\neq j}\| \chi_{j,int}(\operatorname{div}_{\p_{t}^{a}\mathcal{A}_{\eta}}\partial_{t}^{b}v)\|_{W^{m-j+1,p}}\\
&\lesssim \sum_{a+b=j,b\neq j}\| \nabla \chi_{j,int}({\partial_{t}^{a}\mathcal{A}_{\eta}}\partial_{t}^{b}v)\|_{W^{m-j+1,p}} +\sum_{a+b=j,b\neq j}\| (\operatorname{div}_{\p_{t}^{a}\mathcal{A}_{\eta}}\partial_{t}^{b}(\chi_{j,int}v))\|_{W^{m-j+1,p}}\\
&\lesssim \sum_{a+b=j,b\neq j}\|\chi_{j+,int}{\partial_{t}^{a}\mathcal{A}_{\eta}}\|_{W^{m-j+1,p}} \|\chi_{j+1,int}\partial_{t}^{b}v\|_{W^{m-(j+1)+2,p}}\\
&+\sum_{a+b=j,b\neq j}\| \chi_{j+,int}{\partial_{t}^{a}\mathcal{A}_{\eta}}\|_{W^{m-j+1,p}} \|\partial_{t}^{b}(\chi_{j,int}v)\|_{W^{m-j+2,p}} \lesssim \sum_{a+b=j,b\neq j}\|\chi_{j+,int}{\partial_{t}^{a}\bar{\eta}}\|_{W^{m-(j+1)+3-\frac{1}{p},p}}\\
&\|\chi_{j+1,int}\partial_{t}^{b}v\|_{W^{m-j+1,p}} +\sum_{a+b=j,b\neq j}\|\chi_{j+,int}{\partial_{t}^{a}\bar{\eta}}\|_{W^{m-(j+1)+3-\frac{1}{p},p}} \|\chi_{j,int}\partial_{t}^{b}v\|_{W^{m-j+2,p}} \lesssim \sqrt{\mathcal{E}\mathcal{D}}.
\end{aligned}
\end{equation}
Similarly, we have $\|F_{s,\mathrm{int}}^{2,j}\|_{W^{1+(m-j),p}}\lesssim\sqrt{\mathcal{E}}\|v_{j,\mathrm{int}}\|_{W^{2+(m-j),p}}.$

\textbf{Stress tensor term in the capillary equation}

Put $s=m-j+1-1/p$. Since $j\leq m-1$ and $p>3/2$, $s>2/p$, so $W^{s,p}(\Sigma)$ is an algebra. For $c\geq1$, $\partial_t^c\mathcal N$ contains one spatial derivative of $\eta_c$; its localized $W^{s,p}$ norm is bounded by $C(\varepsilon)\|\chi_{j+1,int,\Sigma}\eta_c\|_{W^{s+1,p}}$. For $c=0$, the fixed smooth normal $\mathcal N_0$ contributes a bounded coefficient, and its perturbation is estimated by the same surface norm with $c=0$. We now estimate the following term
\begin{equation}
\begin{aligned}
&\|(\chi_{j,int,\Sigma})\partial_{t}^{j}(S_{\mathcal{A}}(p,v)\mathcal{N})-(\chi_{j,int,\Sigma})(S_{\mathcal{A}}(\partial_{t}^{j}p,\partial_{t}^{j}v)\mathcal{N})\|_{W^{m-j+1-\frac{1}{p},p}}\\
&\lesssim \sum_{a+b+c=j,b\neq j}\|\chi_{j+,int,\Sigma}\partial_{t}^{a}\mathcal{A}\|_{W^{m-j+1-\frac{1}{p},p}}(\|\chi_{j,int,\Sigma}\partial_{t}^{b}v\|_{W^{m-j+2-\frac{1}{p},p}})\|\chi_{j+,int,\Sigma}\partial_{t}^{c}\mathcal{N}\|_{W^{m-j+1-\frac{1}{p},p}}\\
&+\sum_{b+c=j,b\neq j}(\|\chi_{j,int,\Sigma}\partial_{t}^{b}p\|_{W^{m-j+1-\frac{1}{p},p}})\|\chi_{j+1,int,\Sigma}\partial_{t}^{c}\eta\|_{W^{m-j+2-\frac{1}{p},p}}\\
&\lesssim \sum_{\substack{a+b+c=j\\a,c> 0}}\|\chi_{j+,int,\Sigma}\partial_{t}^{a}\mathcal{A}\|_{W^{m-j+1-\frac{1}{p},p}}(\|\chi_{j,int}\partial_{t}^{b}v\|_{W^{m-j+2,p}})\|\chi_{j+1,int,\Sigma}\partial_{t}^{c}\eta\|_{W^{m-(j+1)+3-\frac{1}{p},p}}\\
&+ \sum_{\substack{a_{1}+a_{2}=j\\a_{2}> 0}}(\|\chi_{j,int}\partial_{t}^{a_{1}}v\|_{W^{m-j+2,p}})\|\chi_{j+1,int,\Sigma}\partial_{t}^{a_{2}}\eta\|_{W^{m-(j+1)+3-\frac{1}{p},p}}\\
&+\sum_{\substack{b+c=j\\c>0}}\|\chi_{j,int}\partial_{t}^{b}p\|_{W^{m-j+1,p}}\|\chi_{j+1,int,\Sigma}\partial_{t}^{c}\eta\|_{W^{m-(j+1)+3-\frac{1}{p},p}} \lesssim\sqrt{\mathcal{E}\mathcal{D}}.
\end{aligned}
\end{equation}

\textbf{Curvature term in the capillary equation}

We now estimate
\begin{align}
\sum_{a+b=j,a\neq j}\|(\chi_{j,int,\Sigma}(g\partial_{t}^{a}\eta-\sigma \mathcal{K}(\partial_{t}^{a}\eta))\partial_{t}^{b}\mathcal{N})\|_{W^{m-j+1-\frac{1}{p},p}}.
\end{align}
Since $W^{m-j+1-\frac{1}{p},p}$ is an algebra in 2-D, we have
\begin{equation}
\begin{aligned}
&\sum_{a+b=j,a\neq j}\|(\chi_{j,int,\Sigma}(g\partial_{t}^{a}\eta-\sigma \mathcal{K}(\partial_{t}^{a}\eta))\partial_{t}^{b}\mathcal{N})\|_{W^{m-j+1-\frac{1}{p},p}} \lesssim \sum_{a+b=j,a\neq j}\|\chi_{j,int,\Sigma}\partial_{t}^{a}\eta\|_{W^{m-j+3-\frac{1}{p},p}}\\
&\|\chi_{j+1,int,\Sigma}\partial_{t}^{b}\eta\|_{W^{m-(j+1)+3-\frac{1}{p},p}} +\sum_{a+b=j,a\neq j}\|\chi_{j+1,int,\Sigma}\partial_{t}^{a}\eta\|_{W^{m-(j+1)+3-\frac{1}{p},p}}\\
&\|\chi_{j+1,int,\Sigma}\partial_{t}^{b}\eta\|_{W^{m-(j+1)+3-\frac{1}{p},p}} \lesssim \sqrt{\mathcal{E}\mathscr{D}_{j,high,int}}+\mathscr{E}_{j,high,int} \lesssim\sqrt{\mathcal{E}\mathcal{D}}.
\end{aligned}
\end{equation}

\textbf{Fully Taylor nonlinear term in the capillary equation}

We now estimate $\|\chi_{j,int,\Sigma}\func{div}_{\Sigma}(\partial_{t}^{j}(\mathcal{R}(\nabla \eta))))\|_{W^{1+(m-j)-\frac{1}{p},p}}.$ We have
\begin{equation}
\begin{aligned}
&\|\chi_{j,int,\Sigma}\operatorname{div}_{\Sigma}(\partial_{t}^{j}(\mathcal{R}(\nabla \eta))))\|_{W^{1+(m-j)-\frac{1}{p},p}} \lesssim \sum_{\ell=2}^{j+2}\sum_{i_{1}+\ldots i_{\ell}=j}(\|\chi_{j,\mathrm{int},\Sigma}\partial_{t}^{i_{1}}\eta\|_{W^{3+(m-j)-\frac{1}{p},p}})\\
&\prod_{k=2}^{\ell}(\|\chi_{j+1,\mathrm{int},\Sigma}\partial_{t}^{i_{k}}\eta\|_{W^{3+(m-j-1)-\frac{1}{p},p}}) +\sum_{j=0}^{m}\sum_{\ell=2}^{j+2}\sum_{i_{1}+\ldots i_{\ell}=j} \prod_{k=1}^{\ell}(\|\chi_{j+1,\mathrm{int},\Sigma}\partial_{t}^{i_{k}}\eta\|_{W^{3+(m-j-1)-\frac{1}{p},p}}) \lesssim\sqrt{\mathcal{E}\mathcal{D}}.
\end{aligned}
\end{equation}
Similarly, we have
\begin{align}
\|F_{s,\mathrm{int}}^{3,j}\|_{W^{1-\frac{1}{p}+(m-j),p}}\lesssim\sqrt{\mathcal{E}}(\|v_{j,\mathrm{int}}\|_{W^{2+(m-j),p}}+\|p_{j,\mathrm{int}}\|_{W^{1+(m-j),p}}+\|\eta_{j,\mathrm{int}}\|_{W^{3+(m-j)-\frac{1}{p},p}}).
\end{align}
\textbf{The nonlinear term in the kinematic boundary condition}

We now establish the estimate for $\sum_{a+b=j,a\neq j}\|\partial_{t}^{a}v\cdot \partial_{t}^{b}\mathcal{N}\|_{W^{m-j+2-\frac{1}{p},p}}.$ Since $W^{m-j+2-\frac{1}{p},p}$ is an algebra in 2-D, we have
\begin{equation}
\begin{aligned}
&\sum_{a+b=j,a\neq j}\|\chi_{j,int,\Sigma}\partial_{t}^{a}v\cdot \partial_{t}^{b}\mathcal{N}\|_{W^{m-j+2-\frac{1}{p},p}} \lesssim \sum_{a+b=j,a\neq j}\|\partial_{t}^{a}\nabla (\chi_{j+,int,\Sigma}v)\cdot \partial_{t}^{b}\mathcal{N}\|_{W^{m-j+1-\frac{1}{p},p}}\\
&+\sum_{a+b=j,a\neq j}\|\partial_{t}^{a}v\cdot \partial_{t}^{b}\nabla (\chi_{j,int,\Sigma}\mathcal{N})\|_{W^{m-j+1-\frac{1}{p},p}} +\sum_{a+b=j,a\neq j}\|\partial_{t}^{a}v\cdot \partial_{t}^{b} (\chi_{j,int,\Sigma}\mathcal{N})\|_{W^{m-j+1-\frac{1}{p},p}}\\
&\lesssim \sum_{a+b=j,a\neq j}\|\partial_{t}^{a}\nabla (\chi_{j,int,\Sigma}v)\|_{W^{m-j+1-\frac{1}{p},p}} \|\chi_{j+,int,\Sigma}\partial_{t}^{b}\nabla\eta\|_{W^{m-(j+1)+2-\frac{1}{p},p}}\\
&+\sum_{a+b=j,a\neq j}\|\chi_{j+1,int,\Sigma}\partial_{t}^{a}v\|_{W^{m-(j+1)+2-\frac{1}{p},p}}\| \partial_{t}^{b} (\chi_{j,int,\Sigma}\eta)\|_{W^{m-j+3-\frac{1}{p},p}} \lesssim\sqrt{\mathscr{D}_{j,high,int}}\sqrt{\mathcal{E}} \lesssim\sqrt{\mathcal{E}\mathcal{D}}.
\end{aligned}
\end{equation}
Similarly, we have $\|F_{s,\mathrm{int}}^{4,j}\|_{W^{2+(m-j)-\frac{1}{p},p}}\lesssim\sqrt{\mathcal{E}\mathcal{D}}.$ For $k\leq j-1$, the estimate follows the same idea. Therefore, combining all of the estimates in this step, we have
\begin{align}
\mathscr{D}_{j,high,int} \leq C_j(\varepsilon)\bigl(\mathscr{D}_{j+1,\mathrm{high},\mathrm{int}} +\mathcal{E}\mathcal{D}\bigr).
\end{align}
Combining Step 1 and Step 2, the finite recurrence gives
\begin{align}
\sum_{j=0}^{m-1}\mathscr{D}_{j,high,int} \leq C_m(\varepsilon)\bigl(\mathscr{D}_{m,\mathrm{basic},\mathrm{int}} +\mathcal{E}\mathcal{D}\bigr).
\end{align}
The energy recurrence is analogous, with $\mathscr{E}_{m,\mathrm{basic},\mathrm{int}}$ and $\mathcal{E}^{2}$ on the right. This proves \eqref{eq:enhanced_regularity} and \eqref{eq:enhanced_regularity_e}.
\end{proof}

\begin{remark}
The constant $C_m(\varepsilon)$ is finite after the collar width is fixed. We then choose the smallness of $\mathcal{E}$ as in Subsection~\ref{subsec:smallness}, so that $C_m(\varepsilon)\sqrt{\mathcal{E}}\ll1$.
\end{remark}

The resulting interior regularity estimate is summarized in the following theorem.

\begin{proof}[Proof of Theorem~\protect\ref{thm:interior_space_final}]
By Theorems \ref{thm:whole-space-energy} and \ref{thm:whole_space_estimate_final}, we have
\begin{align}
\mathcal{E}_{||}(t)+\mathcal{E}_{\mathrm{enhanced}}(t)+\int_{0}^{t}(\mathcal{D}_{||}+\mathcal{D}_{\mathrm{enhanced}})\lesssim \mathcal{E}(0)+\mathcal{E}^{\frac{3}{2}}(t)+\mathcal{E}^{\frac{3}{2}}(0)+\int_{0}^{t}\sqrt{\mathcal{E}}\mathcal{D}.
\end{align}
To integrate the full $\mathscr D_{m,\mathrm{basic},\mathrm{int}}$, first estimate its summands of time order $j<m$. Theorem~\ref{thm:pressure_enhancement} and the definitions of $\mathcal D_{||}$ and $\mathcal D_{\mathrm{enhanced}}$ give, at each time,
\begin{equation}
\sum_{j=0}^m\|p_j\|_{L^2}^2 \lesssim \mathcal D_{||}+\mathcal D_{\mathrm{enhanced}} +\mathcal E^2+\mathcal E^3 .
\end{equation}
Indeed $v_j$ and $v_{j+1}$ have orders at most $m+1$, and $\eta_j$ is controlled in $H^{3/2-\alpha}(\Sigma)$. The other terms of $\mathcal E_{||}+\mathcal E_{\mathrm{enhanced}}$ are bounded directly by $\mathcal D_{||}+\mathcal D_{\mathrm{enhanced}}$, including the contact traces of $\partial_t\eta_j$. Since $\mathcal E\le1$, Lemma~\ref{lem:E_le_D} gives $\mathcal E^2+\mathcal E^3 \le 2C_*\mathcal E\mathcal D$. Squaring \eqref{est:basic_energy_int} therefore bounds the sum of the squared lower-order norms. Adding the squared top-order estimate \eqref{est:basic_diss_int}, we obtain
\begin{equation}
\mathscr D_{m,\mathrm{basic},\mathrm{int}} \le C(\varepsilon)\bigl(\mathcal D_{||}+\mathcal D_{\mathrm{enhanced}} +\mathcal E\mathcal D\bigr).
\end{equation}
Integrate this inequality in time. The first two dissipations are controlled by the preceding whole-domain estimate, and $\mathcal E\mathcal D\le \sqrt{\mathcal E}\mathcal D$. For the instantaneous energy use \eqref{est:basic_energy_int} at time $t$. Combining these bounds with \eqref{eq:enhanced_regularity} and \eqref{eq:enhanced_regularity_e} gives
\begin{align}
{\label{eq:interior_basic}} &\mathcal{E}_{||}(t)+\mathcal{E}_{\mathrm{enhanced}}(t)+\mathscr{E}_{m,basic,int}(t)+\sum_{j=0}^{m-1}\mathscr{E}_{j,high,int}+\int_{0}^{t}(\mathcal{D}_{\mathrm{enhanced}}+\mathscr{D}_{m,basic,int}+\sum_{j=0}^{m-1}\mathscr{D}_{j,high,int}) \\
\lesssim &C(\varepsilon)\varepsilon^{-m-1}(\mathcal{E}(0)+\int_{0}^{t}\sqrt{\mathcal{E}}\mathcal{D}+\mathcal{E}^{\frac{3}{2}}(0)+\mathcal{E}^{\frac{3}{2}}(t)).
\end{align}

By the Cauchy inequality and the Sobolev embeddings $W^{3-\frac{1}{p},p}\hookrightarrow H^{\frac{3}{2}+\varepsilon_{-}}$ and $W^{2,p}\hookrightarrow H^{1+\varepsilon_{-}}$, we have
\begin{equation}\label{eq:interior_enhanced}
\begin{aligned}
&\sum_{j=0}^{m}(\|\chi_{0,\mathrm{int},\Sigma}\partial_{t}^{j}\eta\|_{H^{\frac{3}{2}+\frac{\varepsilon_{-}-\alpha}{2}}(\Sigma)}^{2}+\|\chi_{0,\mathrm{int}}\partial_{t}^{j}v\|_{H^{1+\frac{\varepsilon_{-}}2}}^{2})(t)\\
&\lesssim  \sum_{j=0}^{m}(\|\chi_{0,\mathrm{int},\Sigma}\partial_{t}^{j}\eta\|_{H^{\frac{3}{2}+\frac{\varepsilon_{-}-\alpha}{2}}(\Sigma)}^{2}+\|\chi_{0,\mathrm{int}}\partial_{t}^{j}v\|_{H^{1+\frac{\varepsilon_{-}}2}}^{2})(0)\\
&+\int_{0}^{t}\sum_{j=0}^{m}(\|\chi_{0,\mathrm{int},\Sigma}\partial_{t}^{j}\eta\|_{W^{3-\frac{1}{p},p}(\Sigma)}^{2}+\|\chi_{0,\mathrm{int}}\partial_{t}^{j}v\|_{W^{2,p}}^{2})\\
&+\int_{0}^{t}\sum_{j=0}^{m}(\|\chi_{0,\mathrm{int},\Sigma}\partial_{t}^{j+1}\eta\|_{H^{\frac{3}{2}-\alpha}(\Sigma)}^{2}+\|\chi_{0,\mathrm{int}}\partial_{t}^{j+1}v\|_{H^{1}}^{2}) \lesssim \mathcal{E}(0)+\int_{0}^{t}(\mathcal{D}_{\mathrm{m,basic,int}}+\mathcal{D}_{||}+\mathcal{D}_{\mathrm{enhance}}).
\end{aligned}
\end{equation}
Therefore, combining equation \eqref{eq:interior_basic} and equation \eqref{eq:interior_enhanced}, we derive estimate \eqref{est:interior_final}.
\end{proof}

\section{Preparation before the main estimate}

\label{sec:preparation}

In this section, we derive the estimate near the contact line. We first derive the equation system in the local charts.

\subsection{Equation system on the fixed equilibrium domain}

\begin{theorem}
Let $\Phi_{\eta}$ be the diffeomorphism defined in the previous subsection. Applying this map $\Phi_{\eta}$ to the equation system \eqref{eq:moving-domain-full-system}, we obtain the following equation system on the fixed domain:
\begin{equation}  \label{eq:fixed-domain-forcing-system} \left\{ \begin{alignedat}{2} &\partial_tu + \operatorname{div}S(p,u) = F^{1} &&\qquad \text{in }\Omega, \\ &\operatorname{div}u = F^{2} &&\qquad \text{in }\Omega, \\ &S(p,u)n_0 = \left\{ g\eta -\sigma \operatorname{div}_{x,y} \left[ \frac{\nabla_{x,y}(\zeta+\eta)} {\sqrt{1+|\nabla_{x,y}(\zeta+\eta)|^2}} - \frac{\nabla_{x,y}\zeta} {\sqrt{1+|\nabla_{x,y}\zeta|^2}} \right] \right\} n_0 + F_{o}^{3} &&\qquad \text{on }\Sigma, \\ &\frac{1}{|N_{0}|}\partial_t\eta = u\cdot n_0 + F^{4} &&\qquad \text{on }\Sigma, \\ &u\cdot n_0 = 0 &&\qquad \text{on }\Sigma_s, \\ &\bigl(S(p,u)n_0\bigr)_{\tan} -\beta_0u_{\tan} = F^{5} &&\qquad \text{on }\Sigma_s, \\ &-{\kappa}\partial_t\eta-\hat{\mathcal{W}}(\partial_{t}\eta) = \sigma \left[ \frac{\partial_\nu(\zeta+\eta)} {\sqrt{1+|\nabla_{x,y}(\zeta+\eta)|^2}} - \frac{\partial_\nu\zeta} {\sqrt{1+|\nabla_{x,y}\zeta|^2}} \right] &&\qquad \text{on }\Gamma, \\
&u(\cdot,0)=u_0^\sharp, \\ &\eta(\cdot,0)=\eta_0, \end{alignedat} \right.
\end{equation}
where the forcing terms are defined by
\begin{align}
F^{1} :={}& WK\partial_t\bar\eta\,\partial_3u - u\cdot\nabla_{\mathcal{A}_\eta}u + \func{div}S(p,u) - \func{div}_{\mathcal{A}_\eta} S_{\mathcal{A}_\eta}(p,u),  \label{eq:F1} \\
F^{2} :={}& \func{div}u - \func{div}_{\mathcal{A}_\eta}u,  \label{eq:F2} \\
F_{o}^{3} :={}&\frac{1}{|N_{0}|}\big( S(p,u)N_0 - S_{\mathcal{A}_\eta}(p,u)\mathcal{N}  \notag \\
&\quad + \left\{ g\eta -\sigma \func{div}_{x,y} \left[ \frac{\nabla_{x,y}(\zeta+\eta)} {\sqrt{1+|\nabla_{x,y}(\zeta+\eta)|^2}} - \frac{\nabla_{x,y}\zeta} {\sqrt{1+|\nabla_{x,y}\zeta|^2}} \right] \right\} (\mathcal{N}-N_0)\big),  \label{eq:F3} \\
F^{4} :={}&\frac{1}{|N_{0}|}\big( u\cdot(\mathcal{N}-N_0)\big), \label{eq:F4} \\
F^{5} :={}& \left[ \left( S(p,u)-S_{\mathcal{A}_\eta}(p,u) \right)n_w \right]_{\tan}.  \label{eq:F5}
\end{align}
\end{theorem}

\begin{proof}
Using the preceding identities for spatial derivatives and the temporal derivative, the momentum equation becomes
\begin{equation}  \label{eq:fixed-momentum} \partial_tu - WK\partial_t\bar\eta\,\partial_3u + u\cdot\nabla_{\mathcal{A}_\eta}u + \func{div}_{\mathcal{A}_\eta} S_{\mathcal{A}_\eta}(p,u) = 0 \qquad\text{in }\Omega.
\end{equation}
The incompressibility condition becomes  $\func{div}_{\mathcal{A}_\eta}u=0 \,\text{in }\Omega.$

The dynamic free-surface condition becomes
\begin{equation}  \label{eq:fixed-dynamic-bc} S_{\mathcal{A}_\eta}(p,u)\mathcal{N }= \left\{ g\eta - \sigma \left[ H(\zeta+\eta)-H(\zeta) \right] \right\} \mathcal{N }\qquad\text{on }\Sigma.
\end{equation}

The kinematic condition can be written in geometric form as  $\partial_t\eta = u\cdot\mathcal{N }\,\text{on }\Sigma.$ Indeed,  $u\cdot\mathcal{N }= u_3-u_h\cdot\nabla_x(\zeta+\eta).$

Since $\Phi_\eta$ preserves the solid wall, the impermeability condition becomes  $u\cdot n_0=0 \,\text{on }\Sigma_s.$ The Navier-slip condition becomes  $\bigl( S_{\mathcal{A}_\eta}(p,u)n_0 \bigr)_{\tan} - \beta_0u_{\tan} = 0 \,\text{on }\Sigma_s.$

Since the horizontal coordinate $x$ is unchanged by $\Phi_\eta$, the moving contact line $\Gamma(t)$ is pulled back to the fixed equilibrium contact line $\Gamma$. Hence the contact-line equation becomes
\begin{equation}  \label{eq:fixed-contact-line} -{\kappa\partial_t\eta}-\hat{\mathcal{W}}(\partial_{t}\eta) = \sigma \left[ \frac{\partial_\nu(\zeta+\eta)} {\sqrt{1+|\nabla_x(\zeta+\eta)|^2}} - \frac{\partial_\nu\zeta} {\sqrt{1+|\nabla_x\zeta|^2}} \right] \qquad\text{on }\Gamma.
\end{equation}

Combining the preceding equations, the transformed problem on the fixed equilibrium domain is \eqref{eq:moving-domain-elliptic-system}.

Let
\begin{equation*}
N_0 :=
\begin{pmatrix}
-\nabla_x\zeta \\
1\end{pmatrix}
\end{equation*}
denote the non-unit normal vector to the equilibrium free surface. Recall that
\begin{equation*}
\mathcal{N }=
\begin{pmatrix}
-\nabla_x(\zeta+\eta) \\
1\end{pmatrix}, \text{ and hence } \mathcal{N}-N_0 =
\begin{pmatrix}
-\nabla_x\eta \\
0\end{pmatrix}.
\end{equation*}

We rewrite the transformed system by extracting the constant-coefficient linear principal part and the equilibrium normal vector. The resulting system takes the form of \eqref{eq:fixed-domain-forcing-system} with the forcing terms defined by \eqref{eq:F1} to \eqref{eq:F5}.
\end{proof}

\begin{remark}
Compared with system \eqref{eq:moving-domain-elliptic-system}, only the differentiated unknowns are kept on the left-hand side; the geometric terms created by the differentiation are moved to the right-hand side as forcing terms.
\end{remark}

\subsection{Localization}

We first localize the equations on the fixed equilibrium domain and only afterwards apply the contact-line flattening map. In this way, all commutators between the cutoff functions and the differential operators are explicitly included in the forcing terms.

\subsubsection{Localization on the fixed equilibrium domain}

Recall that  $\chi_{cl,\Sigma}$ and its extension $\chi_{cl}$ are the contact line cutoff functions chosen to be independent of $t$ and $\theta$. We define  $v_{cl}:=\chi_{cl} u, \, p_{cl}:=\chi_{cl} p, \, \eta_{cl}:=\chi_{cl,\Sigma}\eta.$

Recall that the fixed-domain system has been written in the frozen form
\begin{equation}  \label{eq:fixed-frozen-system-recall} \left\{ \begin{alignedat}{2} &\partial_tu+\operatorname{div}S(p,u) =F^1 &&\qquad \text{in }\Omega, \\ &\operatorname{div}u =F^2 &&\qquad \text{in }\Omega, \\ &S(p,u)n_0 = \left( g\eta-\sigma\operatorname{div}_{\Sigma}\mathcal L_\Sigma\eta \right)n_0 + F^3 &&\qquad \text{on }\Sigma, \\ &\frac{1}{|N_{0}|}\partial_t\eta-u\cdot n_0 =F^4 &&\qquad \text{on }\Sigma, \\ &u\cdot n_0 =0 &&\qquad \text{on }\Sigma_s, \\ &U_0S(p,u)n_0 -\beta_0U_0u =F^5 &&\qquad \text{on }\Sigma_s, \\ &\sigma\mathcal B_{\mathrm{cl},0}\eta =-F^{6}-\kappa\partial_{t}\eta &&\qquad \text{on }\Gamma, \end{alignedat} \right.
\end{equation}
where  $U_0:=I-n_0\otimes n_0,$ is the tangential projection associated with the equilibrium unit normal $n_0 $, and  $F^{6}=\hat{\mathcal{W}}(\partial_{t}\eta)+\sigma \mathcal{R}(\nabla \eta)\cdot n_{0}.$ Moreover, $\mathcal{L}_\Sigma\eta=\mathcal{L}_{\Sigma,0}\eta := ( M_\zeta\nabla_x\eta ),$ with
\begin{equation*}
M_\zeta = \frac{I}{\sqrt{1+|\nabla\zeta|^2}} - \frac{\nabla\zeta\otimes\nabla\zeta} {(1+|\nabla\zeta|^2)^{3/2}}.
\end{equation*}
All nonlinear terms and all differences between the $\mathcal{A}_\eta$-operators and the frozen operators are already included in $F^1,\ldots,F^6$. Moreover,  $F^{3}=F_{o}^{3}-\sigma\func{div}_{\Sigma}\mathcal{R}(\nabla \eta)n_{0}.$

\subsubsection{Cutoff commutators}

For later use, define  $\mathscr{C} _{\chi_{cl}}(u) := u\otimes\nabla\chi_{cl} + \nabla\chi_{cl}\otimes u.$ Since  $S(p,u) = pI-\mu\mathbb{D}u,$ we have
\begin{equation}  \label{eq:stress-cutoff-j} S(p_{cl},v_{cl}) = \chi_{cl}S(p,u) - \mu\mathscr C_{\chi_{cl}}(u).
\end{equation}
Therefore,
\begin{align}
\partial_tv_{cl} + \func{div}S(p_{cl},v_{cl}) ={}& \chi_{cl} \left( \partial_tu+\func{div}S(p,u) \right)   + S(p,u)\nabla\chi_{cl} - \mu \func{div} \mathscr{C}_{\chi_{cl}}(u).
\end{align}
Hence  $\partial_tv_{cl} + \func{div}S(p_{cl},v_{cl}) = F_{cl}^1 \, \text{in }\Omega,$ where
\begin{equation}  \label{eq:Fj1-before-flattening} F_{cl}^1 := \chi_{cl}F^1 + S(p,u)\nabla\chi_{cl} - \mu \operatorname{div} \mathscr C_{\chi_{cl}}(u).
\end{equation}
Equivalently, the last two terms are precisely the commutator  $[ \func{div}S,\chi_{cl} ](p,u) := \func{div}S(p_{cl},v_{cl}) - \chi_{cl}\func{div}S(p,u).$

For the divergence equation,  $\func{div}v_{cl} = \chi_{cl}\func{div}u + u\cdot\nabla\chi_{cl}.$ Thus  $\func{div}v_{cl} = F_{cl}^2 \, \text{in }\Omega,$ where
\begin{equation}  \label{eq:Fj2-before-flattening} F_{cl}^2 := \chi_{cl}F^2 + u\cdot\nabla\chi_{cl}.
\end{equation}
In particular,  $[\func{div},\chi_{cl}]u = u\cdot\nabla\chi_{cl}.$

\subsubsection{Free-surface commutator}

Define
\begin{equation}  \label{eq:surface-commutator} [\func{div}_{\Sigma}\mathcal{L}_\Sigma,\chi_{cl,\Sigma}]\eta := \func{div}_{\Sigma}\mathcal{L}_\Sigma (\chi_{cl,\Sigma}\eta) - \chi_{cl,\Sigma} \operatorname{div}_{\Sigma} \mathcal{L}_\Sigma\eta.
\end{equation}
Since  $\eta_{cl}=\chi_{cl,\Sigma}\eta,$ we have explicitly
\begin{align}
[\func{div}_{\Sigma}\mathcal{L}_\Sigma,\chi_{cl,\Sigma}]\eta ={}& (M_\zeta\nabla\eta) \cdot\nabla\chi_{cl,\Sigma}   + \func{div}_{\Sigma} \left( \eta M_\zeta \nabla\chi_{cl,\Sigma} \right). \label{eq:surface-commutator-explicit}
\end{align}

Multiplying the dynamic boundary condition by $\chi_{j,\Sigma}$ gives
\begin{align}
S(p_{cl},v_{cl})n_0 ={}& \left( g\eta_{cl} - \sigma\func{div}_{\Sigma}\mathcal{L}_\Sigma\eta_{cl} \right)n_0 + F_{cl}^3,
\end{align}
where
\begin{equation}  \label{eq:Fj3-before-flattening}
\begin{aligned} F_{cl}^3 :={}& \chi_{cl,\Sigma}F^3 - \mu\mathscr C_{\chi_{cl}}(u)n_0  + \sigma [\operatorname{div}_{\Sigma}\mathcal L_\Sigma,\chi_{cl,\Sigma}]\eta \,n_0. \end{aligned}
\end{equation}
Hence the surface forcing contains both the stress-cutoff commutator and the capillary-cutoff commutator.

\subsubsection{Kinematic condition}

Since $\chi_{cl,\Sigma}$ is independent of time,  $\partial_t\eta_{cl} = \chi_{cl,\Sigma}\partial_t\eta.$ Also,  $v_{cl}\cdot n_0 = \chi_{cl,\Sigma}u\cdot n_0 \, \text{on }\Sigma.$ Therefore,  $\frac{1}{|N_{0}|}\partial_t\eta_{cl} - v_{cl}\cdot n_0 = F_{cl}^4 \, \text{on }\Sigma,$ where  $F_{cl}^4 := \chi_{cl,\Sigma}F^4.$ There is no cutoff commutator in this equation since the cutoff is independent of time.

\subsubsection{Solid-wall boundary conditions}

The impermeability condition gives  $v_{cl}\cdot n_0=0 \, \text{on }\Sigma_s.$

For the Navier-slip condition, using \eqref{eq:stress-cutoff-j}, we obtain  $U_0S(p_{cl},v_{cl})n_0 - \beta_0U_0v_{cl} = F_{cl}^5 \, \text{on }\Sigma_s,$ where
\begin{equation}  \label{eq:Fj5-before-flattening} F_{cl}^5 := \chi_{cl}F^5 - \mu U_0 \mathscr C_{\chi_{cl}}(u)n_0.
\end{equation}
Thus the notation $(\cdot)_{\tan}$ is completely replaced by the projection operator $U_0$.

\subsubsection{Contact-line condition}

We choose  $\chi_{cl,\Sigma}\equiv1 \, \text{in a neighborhood of }\Gamma.$ Hence  $\nabla\chi_{cl,\Sigma}=0 \, \text{near }\Gamma,$ and therefore  $\eta_{cl}=\eta, \, \partial_t\eta_{cl}=\partial_t\eta \, \text{on }\Gamma.$ Consequently, no cutoff commutator occurs in the contact-line condition, and  $\sigma\mathcal{B}_{\mathrm{cl},0}\eta_{cl} = -F_{cl}^6-\kappa\partial_{t}\eta \, \text{on }\Gamma,$ where  $F_{cl}^6=F^6 \, \text{on }\Gamma.$

Collecting the preceding equations, we obtain
\begin{equation}  \label{eq:localized-before-flattening} \left\{ \begin{alignedat}{2} &\partial_tv_{cl} + \operatorname{div}S(p_{cl},v_{cl}) = F_{cl}^1 &&\qquad \text{in }\Omega_{cl}, \\ &\operatorname{div}v_{cl} = F_{cl}^2 &&\qquad \text{in }\Omega_{cl}, \\ &S(p_{cl},v_{cl})n_0 = \left( g\eta_{cl}-\sigma\operatorname{div}_{\Sigma}\mathcal L_\Sigma\eta_{cl} \right)n_0 + F_{cl}^3 &&\qquad \text{on }\Sigma, \\
&\frac{1}{|N_{0}|}\partial_t\eta_{cl}-v_{cl}\cdot n_0 = F_{cl}^4 &&\qquad \text{on }\Sigma, \\ &v_{cl}\cdot n_0 = 0 &&\qquad \text{on }\Sigma_s, \\
&U_0S(p_{cl},v_{cl})n_0 - \beta_0U_0v_{cl} = F_{cl}^5 &&\qquad \text{on }\Sigma_s, \\ &\sigma \mathcal B_{\mathrm{cl},0}\eta_{\mathrm{cl}} = -F_{cl}^6-\kappa\partial_{t}\eta_{\mathrm{cl}} &&\qquad \text{on }\Gamma,
\end{alignedat} \right.
\end{equation}
where all of the forcing terms are defined by \eqref{eq:Fj1-before-flattening} to \eqref{eq:Fj5-before-flattening}.

\subsubsection{Nonlinear terms and the cutoff commutator}

It is important that the nonlinear terms also generate cutoff commutators. For example,
\begin{equation}
\chi_{cl} \left( u\cdot\nabla_{\mathcal{A}_\eta}u \right) ={}  u\cdot \nabla_{\mathcal{A}_\eta} (\chi_{cl}u) - \left( u\cdot\nabla_{\mathcal{A}_\eta}\chi_{cl} \right)u ={}  u\cdot\nabla_{\mathcal{A}_\eta}v_{cl} - \left( u\cdot\nabla_{\mathcal{A}_\eta}\chi_{cl} \right)u.\label{eq:nonlinear-cutoff-commutator}
\end{equation}
Thus, if $F^1$ contains  $-u\cdot\nabla_{\mathcal{A}_\eta}u,$ then $\chi_jF^1$ contains  $-u\cdot\nabla_{\mathcal A_\eta}v_{cl} + ( u\cdot\nabla_{\mathcal A_\eta}\chi_{cl} )u.$ The second term is another cutoff commutator and is supported in $\func{supp}\nabla\chi_j$.

\emph{Remark.} Although here we choose $\chi_{j}$ and $\chi_{j,\Sigma}$ to vanish in the interior, we can also multiply the equation system by a cutoff function $\chi_{bet}$ which vanishes near the boundary and the interior of the domain. This system serves as an intermediate step for connecting the interior estimate and the contact line estimate. In this case, we obtain the following equation system
\begin{equation}  \label{eq:localized-before-flattening_bet} \left\{ \begin{alignedat}{2} &\partial_tv_{bet} + \operatorname{div}S(p_{\operatorname{bet}},v_{\operatorname{bet}}) = F_{bet}^1 &&\qquad \text{in }\Omega_{bet}, \\ &\operatorname{div}v_{bet} = F_{bet}^2 &&\qquad \text{in }\Omega_{bet}, \\ &S(p_{bet},v_{bet})n_0 = \left( g\eta_{bet}-\sigma\operatorname{div}_{\Sigma}\mathcal L_\Sigma\eta_{bet} \right)n_0 + F_{bet}^3 &&\qquad \text{on }\Sigma_{bet}, \\ &\frac{1}{|N_{0}|}\partial_t\eta_{bet}-v_{bet}\cdot n_0 = F_{bet}^4 &&\qquad \text{on }\Sigma_{bet}, \\ &v_{\mathrm{bet}}\cdot n_0 = 0 &&\qquad \text{on }\Sigma_s, \\ &U_0S(p_{\mathrm{bet}},v_{\mathrm{bet}})n_0 - \beta_0U_0v_{\mathrm{bet}} = F_{\mathrm{bet}}^5 &&\qquad \text{on }\Sigma_s, \\ &v_{bet}\cdot n_0 = 0 &&\qquad \text{on }\Sigma_{art}, \\
&U_0S(p_{bet},v_{bet})n_0 =0 &&\qquad \text{on }\Sigma_{art}, \\ &\mathcal B_{\mathrm{cl},0}\eta_{art} =0 &&\qquad \text{on }\Sigma_{art}\cap\Sigma_{bet}, \end{alignedat} \right.
\end{equation}
where $v_{bet},p_{bet},\eta _{bet}, F_{bet}^{i}$ are defined by replacing $\chi_{j,\Sigma}$ with $\chi_{bet}$.

\subsection{Flattening transformation and the resulting equation system}

\subsubsection{Transformation to the local contact-line chart}

We now apply the equilibrium contact-line flattening map
\begin{equation*}
\Phi(\theta,r,h) = \left( X(\theta,r,h), \zeta(X(\theta,r,h)) + \Lambda(\theta) (h+r\cot\theta_{\mathrm{eq}}) \right).
\end{equation*}
Let  $\mathcal{A}_\Phi:=(D\Phi)^{-T}, \, Q(\theta):=\mathcal{A}_\Phi(\theta,0,0),$ and define  $\mathcal{B} := Q^T\mathcal{A}_\Phi = I+\mathcal{Z}.$ We also set  $\mathfrak{C}:=Q^T\partial_\theta Q.$ \begin{equation*} -R =\mu\big( \left(Q^{T}\mathcal{A}_{\Phi} e_{\theta}\right)\otimes \left(Q^{T}\partial_{1}Q\,u\right) + \left(Q^{T}\partial_{1}Q\,u\right)\otimes \left(Q^{T}\mathcal{A}_{\Phi} e_{\theta}\right)\big).
\end{equation*}

First define the pullbacks  $\widetilde v_{cl} := Q^T(v_{cl}\circ\Phi), \, \widetilde p_{cl} := p_{cl}\circ\Phi.$ For simplicity, after the transformation we relabel  $(\widetilde v_{cl},\widetilde p_{cl},\widetilde\eta_{cl}) \longmapsto (v_{cl},p_{cl},\xi_{cl}).$ In the following computation, we use $u$ to denote $\tilde{v}_{\mathrm{cl}}$ and $v$ to denote $v_{\mathrm{cl}}$.

For a vector field $X$ defined on $\Omega$, we have that  $(\func{div} X)\circ\Phi = \func{div}_{\mathcal{A}_{\Phi}}(X\circ\Phi),$ so if we set $Y = Q^{T} X\circ\Phi$ then
\begin{equation}
\begin{aligned} (\operatorname{div} X)\circ\Phi &= \operatorname{div}_{A}(QY) = \mathcal{A}_{ij}\,\partial_j (Q_{ik} Y_k) = Q_{ik} \mathcal{A}_{ij}\,\partial_j Y_k + \mathcal{A}_{ij}\,\partial_j Q_{ik}\, Y_k \\ &= \operatorname{div}_{Q^{T}\mathcal{A}} Y + \mathcal{A}_{i1}\,\partial_1 Q_{ik}\, Y_k = \operatorname{div}_{Q^{T}\mathcal{A}} Y + \partial_1 Q\, Y \cdot (\mathcal{A}_{\Phi} e_{\theta}). \end{aligned}
\end{equation}
On the other hand, for a symmetric tensor field $G$ on $\Omega$ we have  $(\func{div} G)\circ\Phi = \func{div}_{\mathcal{A}_{\Phi}}(G\circ\Phi),$ and if $S$ is a symmetric tensor field on $D$ then
\begin{equation}
\func{div}_{\mathcal{A}_{\Phi}}(Q S Q^{T}) = Q\left[\func{div}_{Q^{T}\mathcal{A}_{\Phi}} S + S\,\partial_1 Q^{T} \mathcal{A}_{\Phi} e_{\theta} + Q^{T}\partial_1 Q\, S\, Q^{T} \mathcal{A} _{\Phi}e_{\theta}\right].
\end{equation}
An elementary calculation shows that
\begin{equation}
\mathbb{D} v\circ\Phi = Q\left[\mathbb{D}_{Q^{T}\mathcal{A}_{\Phi}} u + (Q^{T} \mathcal{A}_{\Phi} e_{\theta})\otimes(Q^{T}\partial_1 Q\, u) + (Q^{T}\partial_1 Q\, u)\otimes(Q^{T} \mathcal{A}_{\Phi} e_{\theta})\right] Q^{T}.
\end{equation}
For brevity, we write
\begin{equation}
- R = \mu\big((Q^{T} \mathcal{A}_{\Phi} e_1)\otimes(Q^{T}\partial_1 Q\, u) + (Q^{T}\partial_1 Q\, u)\otimes(Q^{T} \mathcal{A}_{\Phi} e_{\theta})\big).
\end{equation}
Then we have that  $S(q,v)\circ\Phi = Q[S_{Q^{T}\mathcal{A}_{\Phi}}(p,u) + R] Q^{T},$ and
\begin{multline}
(\func{div} S(q,v))\circ\Phi = Q\Big[\func{div}_{Q^{T}\mathcal{A}_{\Phi}}\big(S_{Q^{T}\mathcal{A}_{\Phi}}(p,u) + R\big) + \big(S_{Q^{T}\mathcal{A}_{\Phi}}(p,u) + R\big)\,\partial_{\theta } Q^{T} \mathcal{A}_{\Phi} e_{\theta} \\ + Q^{T}\partial_{\theta} Q\,\big(S_{Q^{T}\mathcal{A}_{\Phi}}(p,u) + R\big)\, Q^{T} \mathcal{A}_{\Phi} e_{\theta}\Big].
\end{multline}
Finally, we also know that
\begin{equation}{\label{eq:geometric_normal_vector}} n\circ\Phi = \frac{\mathcal{A}_{\Phi} n_0}{\lvert \mathcal{A}_{\Phi} n_0\rvert},
\end{equation}
where $n_0$ denotes the unit normal in the dihedron. This lets us compute
\begin{equation}
(I - n\otimes n)\circ\Phi = I - \frac{\mathcal{A}_{\Phi} n_0}{\lvert \mathcal{A}_{\Phi} n_0\rvert}\otimes\frac{\mathcal{A}_{\Phi} n_0}{\lvert \mathcal{A}_{\Phi} n_0\rvert} = Q\left(I - \frac{Q^{T} \mathcal{A}_{\Phi} n_0}{\lvert Q^{T} \mathcal{A}_{\Phi} n_0\rvert} \otimes\frac{Q^{T} \mathcal{A}_{\Phi} n_0}{\lvert Q^{T} \mathcal{A}_{\Phi} n_0\rvert}\right) Q^{T}.
\end{equation}

\subsubsection{Transformation of nonlinear derivatives}

The derivatives appearing in the nonlinear forcing terms must also be transformed by the flattening map. For any scalar function $f$,  $(\nabla f)\circ\Phi = \mathcal{A}_\Phi \nabla_{\theta,r,h}(f\circ\Phi).$ Consequently,
\begin{equation}  \label{eq:Aeta-gradient-flattening} \left( \nabla_{\mathcal{A}_\eta}f \right)\circ\Phi = \left( \mathcal{A}_\eta\circ\Phi \right) A_\Phi \nabla_{\theta,r,h}(f\circ\Phi).
\end{equation}
We therefore introduce the combined nonlinear derivative matrix  $\mathcal{A }:= ( \mathcal{A}_\eta\circ\Phi ) \mathcal{A}_\Phi.$ Then  $( \nabla_{\mathcal{A}_\eta}f )\circ\Phi = \nabla_{\mathcal{A}} (f\circ\Phi).$ For vector fields, recalling that  $v=Q^T(u\circ\Phi), \, u\circ\Phi=Qv,$ the transformed nonlinear convection term is
\begin{equation}  \label{eq:nonlinear-flattened-exact} Q^T \left[ \left( u\cdot\nabla_{\mathcal A_\eta}u \right)\circ\Phi \right] = Q^T \left[ (Qv)\cdot \nabla_{\mathcal A}(Qv) \right].
\end{equation}
Expanding the derivative of $Qv$ gives
\begin{equation}  \label{eq:nonlinear-flattened-expanded} Q^T \left[ (Qv)\cdot \nabla_{\mathcal{A}}(Qv) \right] = \mathcal{N}_{\mathcal{A}}(v) + \mathcal{R}_Q^{\mathrm{nl}}(v),
\end{equation}
where $\mathcal{N}_{\mathcal{B}_\eta}(v)$ contains the derivatives of $v$ with respect to the flattening map and $\mathcal{R}_Q^{\mathrm{nl}}(v)$ contains the connection terms generated by $\partial_\theta Q$.  In particular, the localized convection term \eqref{eq:nonlinear-cutoff-commutator} becomes
\begin{equation}
Q^T \left[ \chi_{cl} \left( u\cdot\nabla_{\mathcal{A}_\eta}u \right) \circ\Phi \right] ={}  Q^T \left[ (u\circ\Phi) \cdot \nabla_{\mathcal{A}} (v_j\circ\Phi) \right]   - Q^T \left[ \left( (u\cdot\nabla_{\mathcal{A}_\eta}\chi_j)u \right) \circ\Phi \right].\label{eq:localized-nonlinear-after-flattening}
\end{equation}
Thus the second term in \eqref{eq:localized-nonlinear-after-flattening} is simultaneously a nonlinear term and a cutoff commutator.

In the remaining part of this section, we use $(v,p,\eta)$ as the simplified notation for $(v_{cl},p_{cl},\eta_{cl})$. Therefore, after applying the flattening map, and keeping for now the curvature tensor in its original form, the equation system \eqref{eq:localized-before-flattening} becomes
\begin{equation}  \label{eq:localized-after-flattening} \left\{ \begin{alignedat}{2} &\partial_t v + \operatorname{div}_{Q^{T}\mathcal{A}_{\Phi}} \bigl( S_{Q^{T}\mathcal{A}_{\Phi}}(p,v)+R(v) \bigr) + \bigl( S_{Q^{T}\mathcal{A}_{\Phi}}(p,v)+R(v) \bigr) (\partial_\theta Q^{T}\mathcal{A}_{\Phi})\,e_\theta \\ &+ \mathfrak{C} \bigl( S_{Q^{T}\mathcal{A}_{\Phi}}(p,v)+R(v) \bigr) Q^{T}\mathcal{A}_{\Phi}e_\theta = \widetilde F_{cl}^1 &&\qquad \text{in }\tilde{\Omega}_{cl}, \\
&\operatorname{div}_{Q^{T}\mathcal{A}_{\Phi}}v+\partial_{\theta}Qv\cdot (\mathcal{A}_{\Phi}e_{\theta}) = \widetilde F_{cl}^2 &&\qquad \text{in }\tilde{\Omega}_{cl}, \\ &\bigl( S_{Q^{T}\mathcal{A}_{\Phi}}(p,v)+R(v) \bigr) \frac{Q^{T}\mathcal{A}_{\Phi}n_{0,+}}{|Q^{T}\mathcal{A}_{\Phi}n_{0,+}|} = \left( g\eta - \sigma\operatorname{div}_{\Sigma}\mathcal L_{\Sigma,0}\eta \right) \frac{Q^{T}\mathcal{A}_{\Phi}n_{0,+}}{|Q^{T}\mathcal{A}_{\Phi}n_{0,+}|} + \widetilde F_{cl}^3 &&\qquad \text{on }\Gamma_{+,cl}, \\
&\frac{1}{|N_0|\circ\Phi_\Sigma} \partial_t\eta_{cl} - v_{cl}\cdot \frac{Q^{T}\mathcal{A}_{\Phi}n_{0,+}}{|Q^{T}\mathcal{A}_{\Phi}n_{0,+}|} = \widetilde F_{cl}^4 &&\qquad \text{on }\Gamma_{+,cl}, \\ &v_{cl}\cdot \frac{Q^{T}\mathcal{A}_{\Phi}n_{0,-}}{|Q^{T}\mathcal{A}_{\Phi}n_{0,-}|} = 0 &&\qquad \text{on }\Gamma_{-,cl}, \\ &U_{Q^{T}\mathcal{A}_{\Phi},-} \bigl( S_{Q^{T}\mathcal{A}_{\Phi}}(p_{cl},v_{cl})+R(v_{cl}) \bigr) \frac{Q^{T}\mathcal{A}_{\Phi}n_{0,-}}{|Q^{T}\mathcal{A}_{\Phi}n_{0,-}|} - \beta_0 U_{Q^{T}\mathcal{A}_{\Phi},-}v_{cl} = \widetilde F_{cl}^5 &&\qquad \text{on }\Gamma_{-,cl}, \\ &\sigma \mathcal B_{\mathrm{cl},0}\eta = -\widetilde F_{cl}^6-\kappa\partial_{t}\eta &&\qquad \text{on }\Gamma_0, \\
&v_{cl} = 0 &&\qquad \text{on }\Gamma_{\mathrm{art},cl}, \\ &\xi_{cl} = 0 &&\qquad \text{on }\Lambda_{\mathrm{art},cl}, \end{alignedat} \right.
\end{equation}
where $\tilde{\Omega}_{j}$ denotes the dihedral region obtained by transforming the local domain $\Omega_{j}$, $\Gamma_{j,\pm}$ are the corresponding boundaries transformed from $\Sigma_{s}$ and $\Sigma$, and $\tilde{{F}_{j}^{i}}$ are the transformed forcing terms defined by
\begin{align}
\tilde{{F}}_{cl}^{1}:=Q^{T}{F}_{cl}^{1}\circ \Phi,\qquad\tilde{{F}}_{cl}^{2}:={F}_{cl}^{2}\circ \Phi,\quad\tilde{{F}}_{cl}^{3}:=Q^{T}{F}_{cl}^{3}\circ \Phi,\quad\tilde{{F}}_{cl}^{4}:={F}_{cl}^{4}\circ \Phi,\quad \tilde{{F}}_{cl}^{5}:=Q^{T}{F}_{cl}^{5}\circ \Phi,\quad\tilde{{F}}_{cl}^{6}:={F}_{cl}^{6}\circ \Phi.
\end{align}
This will be the starting point of our local estimate. Moreover, we have to pull back the curvature on the surface to the local chart; we leave it in its original form here and derive its correct form in the next subsection.

\subsubsection{Frozen local bulk equations}

Define $\mathcal{B}:=Q^{T}\mathcal{A}_{\Phi}$. Moreover, under the local coordinate transformation, define  $D_{\mathcal{B}}v = Dv+D_{\mathcal{Z}}v,$ where  $(D_{\mathcal{Z}}v)_{ab} = \mathcal{Z}_{ak}\partial_kv_{b} + \mathcal{Z}_{bk}\partial_kv_{a}.$ Then, we have  $S_{\mathcal{B}}(p,v) = S(p,v)-\mu D_{\mathcal{Z}}v.$

The variable-coefficient local momentum operator is denoted by  $\mathscr L_B(v,p),$ while the frozen equilibrium operator is  $\mathscr L_0(v,p) := \func{div}S(p,v) + \mathfrak{C}S(p,v)e_\theta.$ We write  $\mathscr L_B(v_{j},p_{j}) = \mathscr L_0(v_{j},p_{j}) - \mathcal{E}_\Phi^1(v_{j},p_{j}),$ where all terms containing  $\mathcal{Z}=B-I, \, \nabla \mathcal{Z}$ are included in $\mathcal{E}_\Phi^1$. More explicitly,
\begin{align}
\mathcal{E}_\Phi^1(v_{j},p_{j}) ={}& -\func{div}_{\mathcal{Z}} S(p_{j},v_{j}) + \mu\func{div}_{\mathcal{B}} D_{\mathcal{Z}}v_{j} - \func{div}_{\mathcal{B}} R(v_{j})  \notag \\
&- \left( S_{\mathcal{B}}(p_{j},v_{j})+R(v_{j}) \right) \partial_\theta Q^{T}\mathcal{A}_{\Phi} \,e_\theta   + \mu \mathfrak{C} D_{\mathcal{Z}}v_{j}\,e_\theta - \mathfrak{C}R(v_{j})e_\theta   - \mathfrak{C} \left( S_B(p_{j},v_{j})+R(v_{j}) \right) \mathcal{Z}e_\theta. \label{eq:local-geometric-error}
\end{align}

Therefore, after localization and flattening, the bulk system is
\begin{equation}  \label{eq:local-bulk-final} \left\{ \begin{alignedat}{2} &\partial_tv_{cl} + \operatorname{div}S(p_{cl},v_{cl}) + \mathfrak{C}S(p_{cl},v_{cl})e_\theta = \mathcal F_{cl}^1 &&\qquad \text{in } \Omega_{cl}, \\
&\operatorname{div}v_{cl} = \mathcal F_{cl}^2 &&\qquad \text{in } \Omega_{cl}, \end{alignedat} \right.
\end{equation}
where  $\mathcal F_{cl}^1 = Q^T(F_{cl}^1\circ\Phi) + \mathcal E_\Phi^1(v_{cl},p_{cl}),$ and
\begin{equation}  \label{eq:local-F2-final} \mathcal F_{\mathrm{cl}}^2=H_{\mathrm{cl}}-\operatorname{div}_{\mathcal Z}v_{\mathrm{cl}},\qquad H_{\mathrm{cl}}:=F_{\mathrm{cl}}^2\circ\Phi- \bigl((\partial_\theta Q)v_{\mathrm{cl}}\bigr)\cdot(\mathcal A_\Phi e_\theta).
\end{equation}
Notice in particular that  $Q^T(F_{cl}^1\circ\Phi)$ already contains the cutoff commutator terms  $Q^T [ ( S(p,u)\nabla\chi_{cl} )\circ\Phi ]$, -$\mu Q^T [ ( \func{div} \mathscr C_{\chi_{cl}}(u) )\circ\Phi ]$, and  $Q^T [ ( u\cdot\nabla_{\mathcal{A}_\eta}\chi_{cl} )u \circ\Phi ].$

\emph{Remark.} In the nonlinear terms $F_{cl}^{j}\circ \Phi$, we also have to transform all the derivatives into the local charts. We will write down further details in the elliptic estimate part.

\subsubsection{Frozen local boundary conditions}

Let  $\nu_{\mathcal{B},\pm} := \frac{\mathcal{B}n_{0,\pm}} {|\mathcal{B}n_{0,\pm}|}, \, U_{\mathcal{B},\pm} := I-\nu_{\mathcal{B},\pm}\otimes\nu_{\mathcal{B},\pm},$ and  $U_{0,\pm} := I-n_{0,\pm}\otimes n_{0,\pm}.$ All differences  $\nu_{\mathcal{B},\pm}-n_{0,\pm}, \, U_{\mathcal{B},\pm}-U_{0,\pm}, \, \mathcal{B}-I$ are placed on the right-hand side.

The free-surface conditions take the form
\begin{equation}  \label{eq:local-free-boundary-final} \left\{ \begin{alignedat}{2} &S(p_{cl},v_{cl})n_{0,+} = \left( g\eta_{cl} - \sigma\operatorname{div}_{\Sigma}\mathcal L_{\Sigma,0}\eta_{cl} \right)n_{0,+} + \mathcal F_{cl}^3 &&\qquad \text{on }\Gamma_{+,cl}, \\
&\frac{1}{|N_{0}|}\partial_t\eta_{cl}-v_{cl}\cdot n_{0,+} = \mathcal F_{cl}^4 &&\qquad \text{on }\Gamma_{+,cl}, \end{alignedat} \right.
\end{equation}
where $\mathcal{F}_{cl}^3$ contains, in particular,  $\sigma ( [\operatorname{div}_{\Sigma}\mathcal L_\Sigma,\chi_{cl,\Sigma}]\eta )\circ X_0\,n_{0,+}.$

The solid-wall conditions are
\begin{equation}  \label{eq:local-wall-final} \left\{ \begin{alignedat}{2} &v_{cl}\cdot n_{0,-} = \mathcal F_{cl,-}^{4}=v_{\mathrm{cl}}\cdot (n_{0,-}-\frac{Q^{T}\mathcal{A}_{\Phi}n_{0,-}}{|Q^{T}\mathcal{A}_{\Phi}n_{0,-}|}) &&\qquad \text{on }\Gamma_{-,cl}, \\ &U_{0,-}S(p_{cl},v_{cl})n_{0,-} - \beta_0U_{0,-}v_{cl} = \mathcal F_{cl}^5 &&\qquad \text{on }\Gamma_{-,cl},
\end{alignedat} \right.
\end{equation}
where
\begin{align}
\mathcal F_{\mathrm{cl}}^5={}&\widetilde F_{\mathrm{cl}}^5-\mu(U_{0,-}-U_{\mathcal{B},-})\mathbb Dv_{\mathrm{cl}}\,n_{0.-} +U_{\mathcal{B},-}(\mu\mathbb D_Zv_{\mathrm{cl}}-R)n_{0,-}\no\\
&-U_{\mathcal{B},-}(\mu\mathbb D_{\mathcal{B},-}v_{\mathrm{cl}}-R)(n_{0,-}-\nu_{\mathcal{B},-})-\beta_0(U_{0,-}-U_{\mathcal{B},-})v_{\mathrm{cl}}.
\end{align}

Since  $\chi_{cl,\Sigma}\equiv1 \, \text{near }\Gamma,$ the contact-line equation does not acquire a cutoff commutator:
\begin{equation}  \label{eq:local-contact-final} \sigma \mathcal{B}_{\mathrm{cl},0}\eta_{cl} = -\mathcal{F}_{cl}^6-\kappa\partial_{t}\eta_{\mathrm{cl}} \qquad \text{on }\Gamma_0.
\end{equation}

Finally, we write down the explicit derivation of those forcing terms on the boundary. We define  $Y_{\mathcal{B},cl} := S_{\mathcal{B}}(p_{cl},v_{cl}), \, Y_{0,cl} := S(p_{cl},v_{cl}),$ where  $S_{\mathcal{B}}(p_j,v_j) = S(p_j,v_j)-\mu D_{\mathcal{Z}}v_j,$ and  $(D_{\mathcal{Z}}v_j)_{ab} = \mathcal{Z}_{ak}\partial_k(v_j)_b + \mathcal{Z}_{bk}\partial_k(v_j)_a.$ Hence
\begin{equation}  \label{eq:delta-Y-j} \delta Y_j := Y_{\mathcal{B},j}-Y_{0,j} = -\mu D_{\mathcal{Z}}v_j .
\end{equation}

\paragraph{1. Dynamic free-surface forcing $\mathcal{F}_{cl}^3$.}

Before applying the flattening map, multiplication by $\chi_{cl,\Sigma}$ gives the cutoff forcing
\begin{equation}  \label{eq:physical-cutoff-Fj3} F_{cl}^3 := \chi_{cl,\Sigma}F^3 - \mu\mathscr C_{\chi_{cl}}(u)n_0 + \sigma [\func{div}_{\Sigma}\mathcal{L}_\Sigma,\chi_{cl,\Sigma}]\eta\,n_0,
\end{equation}
where $\mathscr C_{\chi_{cl}}(u) = u\otimes\nabla\chi_{cl} + \nabla\chi_{cl}\otimes u,$ and
\begin{equation}
[\func{div}_{\Sigma}\mathcal{L}_\Sigma,\chi_{cl,\Sigma}]\eta ={}  \func{div}_{\Sigma}\mathcal{L}_\Sigma (\chi_{cl,\Sigma}\eta) - \chi_{cl,\Sigma}\func{div}_{\Sigma}\mathcal{L}_\Sigma\eta ={}  (M_\zeta\nabla\eta) \cdot\nabla\chi_{cl,\Sigma} + \func{div} \left( \eta M_\zeta\nabla\chi_{cl,\Sigma} \right).
\end{equation}

Let  $\widetilde{ F}_j^3 := Q^T ( F_j^3\circ\Phi_\Sigma ).$ The variable boundary condition after flattening is
\begin{equation*}
(Y_{\mathcal{B},cl}+R(v_{cl}))\nu_{\mathcal{B},+} = (g\eta_{cl} - \sigma\func{div}_{\Sigma}\mathcal{L}_{\Sigma,0}\eta_{cl})\nu_{\mathcal{B},+} + \widetilde{ F}_{cl}^3.
\end{equation*}
Therefore, freezing the normal and the coefficients, we have
\begin{equation*}
Y_{0,cl}n_{0,+} = (g\eta_{cl} - \sigma\func{div}_{\Sigma}\mathcal{L}_{\Sigma,0}\eta_{cl})n_{0,+} + \mathcal{F}_{cl}^3,
\end{equation*}
with
\begin{equation}
\mathcal F_{cl}^3 ={}  \widetilde{ F}_{cl}^3-R(v_{cl})n_{0,+} - \delta Y_j\,n_{0,+} + \left( Y_{\mathcal{B},cl}+R(v_{cl})-(g\eta_{cl} - \sigma\operatorname{div}_{\Sigma}\mathcal L_{\Sigma,0}\eta_{cl})I \right) \left( n_{0,+}-\nu_{B,+} \right).\label{eq:explicit-Fj3}
\end{equation}

Using \eqref{eq:delta-Y-j}, this becomes
\begin{equation}
\mathcal F_{cl}^3 ={}  \widetilde{ F}_{cl}^3 + \mu D_{\mathcal{Z}}v_{cl}\,n_{0,+} - \left( R(v_{cl})) \right)n_{0,+} + \left[ S_{\mathcal{B}}(p_{cl},v_{cl}) + R(v_{cl}) - (g\eta_{cl} - \sigma\operatorname{div}_{\Sigma}\mathcal L_{\Sigma,0}\eta_{cl})I \right] \left( n_{0,+}-\nu_{B,+} \right).
\end{equation}
Thus $\mathcal{F}_{cl}^3$ consists of the cutoff commutators, the coefficient difference $E=B-I$, the lower-order terms induced by the flattening map, the connection difference $R-R_0$, and the normal difference $\nu_{\mathcal{B},+}-n_{0,+}$.

\paragraph{2. Kinematic forcing $\mathcal{F}_j^4$.}

Before flattening, we have  $\partial_t\eta_{cl}/|\mathcal{N}_{0}| - v_{cl}\cdot n_0 = F_{cl}^4$.  Define  $\widetilde{ F}_{cl}^4 := ( \chi_{cl,\Sigma}F^4 )\circ \Phi_{\Sigma}.$ The variable local kinematic condition is
\begin{equation*}
\frac{1}{|N_{0}|}\partial_t\eta_{cl} - v_{cl}\cdot\nu_{\mathcal{B},+} = \widetilde{ F}_{cl}^4.
\end{equation*}
Hence  $\mathcal F_{cl}^4 = \widetilde{F}_{cl}^4 + v_{cl}\cdot ( \nu_{B,+}-n_{0,+} ).$   In particular, given that  $F_{cl}^4 = \frac{1}{|N_{0}|}v_{cl}\cdot(\mathcal{N}-N_0),$ then
\begin{equation}
\mathcal F_{cl}^4 = \frac{1}{|N_{0}|}\left[ v_{cl}\cdot(\mathcal N-N_0) \right]\circ \Phi_{\Sigma} + v_{cl}\cdot (\nu_{B,+}-n_{0,+}).
\end{equation}

\emph{Remark.} Here we have not changed $\eta$ to the surface perturbation $\mathfrak{h}$ in local charts, since the system is already complicated. Moreover, changing to $\mathfrak{h}$ is only useful for the elliptic estimate. We will give the details in the elliptic estimate section.

\paragraph{3. Normal solid-wall forcing $\mathcal{F}_{cl}^5$.}

Since multiplication by $\chi_{cl}$ preserves the impermeability condition, there is no cutoff derivative term. After flattening, the variable normal condition is
\begin{equation*}
(\chi_{cl} v)\cdot \frac{1}{|Q^{T}\mathcal{A}_{\Phi}n_{0,-}|}Q^{T}\mathcal{A}_{\Phi}n_{0,-} =0.
\end{equation*}
Here we have to introduce the forcing term $\mathcal{F}_{\mathrm{cl},-}^{4}$. Using  $v_{cl}\cdot\nu_{B,-}=0,$ we have  $v_{cl}\cdot n_{0,-} =\mathcal{F}_{\mathrm{cl},-}^{4},$ where
\begin{align}
\mathcal{F}_{\mathrm{cl},-}^{4}=v_{\mathrm{cl}}\cdot (n_{0,-}-\frac{Q^{T}\mathcal{A}_{\Phi}n_{0,-}}{|Q^{T}\mathcal{A}_{\Phi}n_{0,-}|}).
\end{align}

\paragraph{4. Navier-slip forcing $\mathcal{F}_{cl}^6$.}

Before flattening, multiplying the Navier-slip condition by $\chi_{cl}$ produces the cutoff commutator  $F_{cl}^{5} := \chi_{cl}F^5 - \mu U_0 \mathscr C_{\chi_{cl}}(u)n_0.$ Define its transformed version by  $\widetilde{ F}_{cl}^{5} := Q^T ( F_{cl}^{5} \circ\Phi|_{\Gamma_-} ).$ The variable local Navier-slip condition is
\begin{equation*}
U_{\mathcal{B},-} (Y_{\mathcal{B},cl}+R(v_{\mathrm{cl}}))\nu_{\mathcal{B},-} - \beta_0U_{\mathcal{B},-}v_{cl} = \widetilde{ F}_{cl}^{5}.
\end{equation*}
We want to rewrite it in the frozen form
\begin{equation*}
U_{0,-} S(p_{\mathrm{cl}},v_{\mathrm{cl}})n_{0,-} - \beta_0U_{0,-}v_{cl} = \mathcal{F}_{cl}^5.
\end{equation*}
Hence,
\begin{align}
&\mathcal F_{cl}^5 ={}  \widetilde{ F}_{cl}^{5} -\mu \left( U_{0,-}-U_{\mathcal{B},-} \right) \left( \mathbb{D}v_{cl}) \right)n_{0,-} + U_{\mathcal{B},-} \left[ \mu\mathbb{D}_{\mathcal{Z}}v_{cl} - R(v_{cl}) \right]n_{0,-}\notag \\
&- U_{\mathcal{B},-} \left( \mu\mathbb{D}_Bv_{cl}-R(v_{cl}) \right) \left( n_{0,-}-\nu_{\mathcal{B},-} \right)- \beta_0 \left( U_{0,-}-U_{\mathcal{B},-} \right)v_{cl}.\label{eq:explicit-Fj6}
\end{align}

\subsubsection{Artificial boundary}

Choose a smooth truncated local domain  $\Omega_{cl}\subset \Omega,$ such that it touches the physical faces $\Sigma_{s},\Sigma$, with its artificial face in the interior of $\Omega$. After the flattening map, ${\Omega}_{cl}$ is transformed into $\tilde{\Omega}_{cl}$ whose boundary decomposes as  $\partial\Omega_{cl} = \Gamma_{+,cl} \cup \Gamma_{-,cl} \cup \Gamma_{\mathrm{art},cl},$ and such that  $\chi_{cl}\circ\Phi=0$ in a full neighborhood of $\Gamma_{\mathrm{art},cl}$. Moreover, $\Gamma_{art,cl}$ should be nicely chosen such that the boundary $\partial {\Omega}_{cl}$ is smooth except for the contact line. Therefore,  $v_{cl}=0 \, \text{near }\Gamma_{\mathrm{art},cl}.$ Similarly, if  $\Lambda_{\mathrm{art},cl} = \overline{\Gamma_{+,cl}} \cap \Gamma_{\mathrm{art},cl},$ then  $\chi_{cl,\Sigma}=0$ near $\Lambda_{\mathrm{art},cl}$, and hence  $\eta_{cl}=0 \, \text{near }\Lambda_{\mathrm{art},cl}.$

\subsubsection{Final localized and flattened system}

Finally, the localized problem in the contact-line chart is written as
\begin{equation}  \label{eq:final-localized-flattened-system-j} \left\{ \begin{alignedat}{2} &\partial_tv_{cl} + \operatorname{div}S(p_{cl},v_{cl}) + \mathfrak{C}S(p_{cl},v_{cl})e_\theta = \mathcal F_{cl}^1 &&\qquad \text{in } \Omega_{cl}, \\
&\operatorname{div}v_{cl} = \mathcal F_{cl}^2 &&\qquad \text{in } \Omega_{cl}, \\ &S(p_{cl},v_{cl})n_{0,+} = \left( g\eta_{cl} - \sigma\operatorname{div}_{\Sigma}\mathcal L_{\Sigma,0}\eta_{cl} \right)n_{0,+} + \mathcal F_{cl}^3 &&\qquad \text{on }\Gamma_{+,cl}, \\
&\frac{1}{|\mathcal{N}_{0}|}\partial_t\eta_{cl}-v_{cl}\cdot n_{0,+} = \mathcal F_{cl}^4 &&\qquad \text{on }\Gamma_{+,cl}, \\ &v_{cl}\cdot n_{0,-} =\mathcal{F}_{\mathrm{cl},-}^{4} &&\qquad \text{on }\Gamma_{-,cl}, \\
&U_{0,-}S(p_{cl},v_{cl})n_{0,-} - \beta_0U_{0,-}v_{cl} = \mathcal F_{cl}^5 &&\qquad \text{on }\Gamma_{-,cl}, \\ &\sigma \mathcal B_{\mathrm{cl},0}\eta_{cl} = -\mathcal F_{cl}^6-\kappa\partial_{t}\eta_{\mathrm{cl}} &&\qquad \text{on }\Gamma_0, \\ &v_{cl}\cdot n_{\operatorname{art}}=0 &&\qquad \text{on }\Gamma_{\mathrm{art},cl}, \\ &S(p_{cl},v_{cl}) n_{\operatorname{art}}=0 &&\qquad \text{on }\Gamma_{\mathrm{art},cl},\\
&(v_{\mathrm{cl}},\eta_{cl},p_{\mathrm{c}},\eta_{\mathrm{cl}})=0&&\qquad\text{near} \quad\Gamma_{\operatorname{art},cl},
\end{alignedat} \right.
\end{equation}
where for simplicity, we still use $\mathcal{N}_{0}$ to denote $|\mathcal{N}_{0}\circ\Phi_{\Sigma}|$ in the local chart.

The forcing terms $\mathcal{F}_{cl}^k$, defined in the previous subsubsection, contain three types of terms:
\begin{equation*}
\mathcal F_{cl}^k = \text{cutoff commutators} + \text{flattening-map differences} + \text{nonlinear terms}+\text{geometric differences}.
\end{equation*}
More precisely, the cutoff contributions contain  $\nabla\chi_{cl}, \, \nabla^2\chi_{cl}, \, [\mathcal{L}_\Sigma,\chi_{cl,\Sigma}],$ the flattening contributions contain  $\mathcal{Z}=\mathcal{B}-I, \, \nabla \mathcal{Z}, \, \nu_{\mathcal{B}}-n_0, \, U_{\mathcal{B}}-U_0,$ and all nonlinear derivatives are taken with respect to the pulled-back flattening derivative  $\nabla_{\mathcal{A}}, \, \mathcal{A }= (\mathcal{A}_\eta\circ\Phi)A_\Phi.$ All geometric difference contributions contain  $\nabla_{I-\mathcal{A}_{\eta}}.$

\subsection{The main idea of the local estimate}

The local estimate consists of three parts: the contact-line elliptic estimate, the intermediate elliptic estimate and the contact-line energy estimate. The difference between the contact-line part and the interior part is that we have to take $\theta$ derivatives as well as time derivatives of the equation system to obtain higher regularity in the $\theta$ direction near the contact line. Because of the contact angle, the regularity in the $r,h$ directions has an upper threshold.

Since $Q^{T}\mathcal{A}_{\Phi}$ and $\mathcal{N}_{0}$ depend on the steady state $\zeta$, when we take a $\theta$ derivative there is an extra linear term induced by this geometric structure, compared to a pure temporal derivative in the whole space. This makes the contact-line estimate more delicate. But the main idea follows the two-dimensional estimates of \cite{GT}.

\section{Estimates in the local charts}

\label{sec:local_charts}

Define
\begin{equation*}
Y_{B,j}:=S_{B,j}(p_{cl},v_{cl}), \qquad Y_{0,j}:=S(p_{cl},v_{cl}), \qquad \delta Y_{cl}:=Y_{B,cl}-Y_{0,cl}.
\end{equation*}
Moreover, we first ignore the surface tensor terms. Then the system \eqref{eq:final-localized-flattened-system-j} may be written as
\begin{align}
{\label{eq:momentum}}
\begin{cases}
\partial_{t}v_{cl}+\func{div}Y_{0,cl} =\mathcal{F}_{cl}^{1} -\mathfrak{C }Y_0e_{\theta}, & \qquad\func{in}\quad\tilde{\Omega}_{cl} \\[1mm]
\func{div}v_{cl}=\mathcal{F}_{cl}^{2}, & \qquad\func{in}\quad\tilde{\Omega}_{cl} \\[1mm]
S(p_{cl},v_{cl})n_{0,+}=\mathcal{F}_{cl}^{3}, & \qquad\func{on}\quad \Gamma_{+,cl} \\[1mm]
\frac{1}{|\mathcal{N}_{0}|}\partial_{t}\eta_{cl}-v_{cl}\cdot n_{0,+}=\mathcal{F}_{cl}^{4}, & \qquad\func{on}\quad \Gamma_{+,cl} \\[1mm]
v_{cl}\cdot n_{0,-}=\mathcal{F}^{4}_{\mathrm{cl}-}, & \qquad\func{on}\quad \Gamma_{-,cl} \\[1mm]
U_0\bigl[S(p_{cl},v_{cl})n_{0,-}-\beta_{0} v_{cl}\bigr] =\mathcal{F}_{cl}^{5}, & \qquad\func{on}\quad \Gamma_{-,cl},\end{cases}\end{align} where $\mathcal{F}_{cl}^{k}$ are forcing terms as defined in the previous section. Moreover, the boundary conditions of the artificial boundary are absorbed into $\Gamma_{-}$.  Here, recall
\begin{equation*}
\widehat n_{\mathcal{Z}} = \frac{(I+\mathcal{Z})n_0}{\lvert(I+\mathcal{Z})n_0\rvert}, \qquad U_{\mathcal{Z}} = I-\widehat n_{\mathcal{Z}}\otimes\widehat n_{\mathcal{Z}}, \qquad U_0 = I-n_0\otimes n_0,
\end{equation*}
and  $\delta M = D_{\mathcal{Z}}v_{j}.$

\subsection{Elliptic estimate obtained by separating the tangential and transverse derivatives}

We work in a contact-line coordinate neighborhood of the form  $\tilde{\Omega}_{cl} = \mathbb{T}\times\mathcal{W}_{cl},$ where  $\theta\in\mathbb{T }$ is the coordinate tangent to the contact line and  $(r,h)\in\mathcal{W}_{cl}\subset\mathbb{R}^{2}$ are the transverse variables. For simplicity, we ignore the subscript $j$. We write  $\nabla_{r,h}:=(\partial_r,\partial_h), \, \func{div}_{r,h}v := \partial_r v_r+\partial_h v_h.$ Let  $e_\theta:=(1,0,0)$ denote the unit vector in the $\theta$-direction, and define  $\mathcal{P}_{r,h}:=I-e_\theta\otimes e_\theta.$ We decompose the velocity as  $v=v_\theta e_{\theta }+v_{r,h}, \, v_\theta:=v\cdot e_\theta, \, v_{r,h}:=P_{r,h}v.$ From the definition of the tangential derivative $\partial_{\theta}$ and tangent vector $e_{\theta}$, $n_0$ and the transverse tangent vector $\tau_0$ satisfy  $n_0\cdot e_\theta=0, \, \tau_0\cdot e_\theta=0.$

\subsubsection{The transverse two-dimensional Stokes system}

Let the indices $a,b$ range over the transverse variables $\{r,h\}$. First, since every term in $R_0$ contains $e_\theta$ in at least one tensor component, its transverse block vanishes:  $(R_0)_{ab}=0, \, a,b\in\{r,h\}.$ Moreover, by the definition of $Y_{0}$ and the projection,  $P_{r,h}Y_0P_{r,h} = pP_{r,h}-\mu D_{r,h}v_{r,h},$ where  $D_{r,h}v_{r,h} := \nabla_{r,h}v_{r,h} + (\nabla_{r,h}v_{r,h})^{T}.$

We separate the divergence of $Y_0$ as  $\func{div}Y_0 = \partial_{\theta}(Y_0e_\theta) + \func{div}_{r,h} \bigl(Y_0P_{r,h}\bigr).$ Applying $P_{r,h}$ to the momentum equation in \eqref{eq:momentum}, we obtain
\begin{align}
\func{div}_{r,h} S_{r,h}(p,v_{r,h}) ={}& P_{r,h} \mathcal{F}_{cl}^{1}- P_{r,h}\partial_{\theta}(Y_0e_{\theta}) - P_{r,h}\mathscr CY_0e_\theta-\p_{t}v_{r,h}, \label{eq:transverse-momentum-r_0}
\end{align}
where  $S_{r,h}(p,v_{r,h}) := pI_{2\times2}-\mu \mathbb{D}_{r,h}v_{r,h}.$ We therefore define
\begin{equation}  \label{eq:Gst1-definition} \mathcal F_{r,h}^{1} := P_{r,h}\mathcal{F}_{cl}^{1} - P_{r,h}\partial_\theta(Y_0e_\theta) - P_{r,h}\mathscr CY_0e_\theta.
\end{equation}
Then  $\partial_{t}v_{r,h}+\func{div}_{r,h} S_{r,h}(p,v_{r,h}) = \mathcal{F}_{r,h}^{1}.$ The divergence equation becomes  $\func{div}_{r,h}v_{r,h} = \mathcal{F}_{cl}^{2}-\partial_{\theta} v_{\theta}.$ Hence we set  $\mathcal F_{r,h}^{2} := \mathcal{F}_{cl}^{2}-\partial_\theta v_\theta.$

\subsubsection{Boundary conditions for the transverse system}

We distinguish the fixed solid boundary, denoted by $\Gamma_{-}$, from the perturbed free boundary, denoted by $\Gamma_{+}^{\eta}$. For each fixed value of the contact-line variable $\theta$, we decompose the velocity as  $v = v_{\theta}e_{\theta} + v_{r,h}, \, v_{r,h} = v_{r}e_{r}+v_{h}e_{h}.$ Let $n_{0,\pm}$ and $\tau_{0,\pm}$ denote the reference unit normal and the transverse unit tangent on $\Gamma_{\pm}$, respectively. Since the reference dihedral geometry is invariant in the $e_{\theta}$-direction,  $n_{0,\pm}\cdot e_{\theta} = \tau_{0,\pm}\cdot e_{\theta} = 0, \, n_{0,\pm}\cdot\tau_{0,\pm}=0.$

\paragraph{The fixed solid boundary.}

Assume that the coordinate transformation flattens the solid boundary exactly and that, in the frame used to represent $v$, its unit normal and transverse tangent are $n_{0,-}$ and $\tau_{0,-}$. Then the normal boundary condition reduces to  $v_{r,h}\cdot n_{0,-} = \mathcal{F}_{\mathrm{cl},-}^{4} \,\text{on }\Gamma_{-}.$ The principal transverse Navier slip condition is
\begin{equation}  \label{eq:transverse-navier-minus} -\mu\mathbb{D}_{r,h}v_{r,h}\,n_{0,-}\cdot\tau_{0,-}-\beta_{0}v_{r,h}\cdot \tau_{0,-} =\mathcal{F}_{cl}^{5}\cdot \tau_{0,-} \qquad\text{on }\Gamma_{-},
\end{equation}
where  $\mathbb{D}_{r,h}v_{r,h} := \nabla_{r,h}v_{r,h} + \bigl(\nabla_{r,h}v_{r,h}\bigr)^{T}.$

Now we derive the boundary conditions on $\Gamma_{+}$. Since $e_{\theta}\cdot n_{0,+}=0$, the kinematic boundary condition on $\Gamma_{+}$ can be rewritten as  $\frac{1}{|\mathcal{N}_{0}|}\partial_{t}\eta-v_{r,h}\cdot n_{0,+}=\mathcal{F}_{cl}^{4}.$

For the capillary boundary condition on $\Gamma_{+}$, taking the scalar product of both sides of equation \eqref{eq:momentum} with $\tau_{0,+}$, which is the tangent vector on the surface $\Sigma$ perpendicular to $e_{\theta}$, we obtain
\begin{equation}  \label{eq:principal-transverse-plus} -\mu\mathbb{D}_{r,h}v_{r,h}\,n_{0,+}\cdot\tau_{0,+} = \mathcal{F}_{cl}^{3}\cdot \tau_{0,+},
\end{equation}
where we use the fact that  $e_{\theta}\cdot {n}_{0}=0~~~\func{and}~~~e_{\theta}\cdot\tau_{0}=0.$ Moreover, recall that the surface tensor term contained in the capillary equation is written as  $(g\eta-\sigma\mathcal{L}_{\Sigma}\eta)n_{0,+}.$ Taking the scalar product of this term with $\tau_{0,+}$ gives zero. So we do not need to consider the surface tensor in this part. Here $e_{\theta}$ is perpendicular to $n_{0}$, because $e_{\theta}$ is tangent to the wall at equilibrium and $n_{0}$ is its normal; $\tau_{0}$ is the unit vector perpendicular to $e_{\theta}$ and $n_{0}$.

After localizing near the contact line, the pair  $\bigl(v_{r,h}(\theta,\cdot),p(\theta,\cdot)\bigr)$ satisfies, for each fixed $\theta$, a two-dimensional system of the form
\begin{equation}  \label{eq:two-dimensional-slice-system}
\begin{cases}
\partial_{t}v_{r,h}+\func{div}_{r,h} S_{r,h}(p_{cl},v_{r,h}) = \mathcal{F}_{r,h}^{1} & \text{in }\mathcal{W}_{cl}, \\[1mm]
\func{div}_{r,h}v_{r,h} = \mathcal{F}_{r,h}^{2} & \text{in }\mathcal{W}_{cl}, \\[1mm]
\frac{1}{|\mathcal{N}_{0}|}\partial_{t}\eta_{cl}-v_{r,h}\cdot n_{0,+} = \mathcal{F}_{cl}^{4} & \text{on }\Gamma_{+}, \\[1mm]
v_{r,h}\cdot n_{0,-} = \mathcal{F}_{\mathrm{cl},-}^{4} & \text{on }\Gamma_{-}, \\[1mm]
-\mu\mathbb{D}_{r,h}v_{r,h}\,n_{0,+}\cdot\tau_{0,+} =\mathcal{F}_{cl}^{3}\cdot \tau_{0,+} & \text{on }\Gamma_{+} \\[1mm]
-\mu\mathbb{D}_{r,h}v_{r,h}\,n_{0,-}\cdot\tau_{0,-}- \beta_{0} v_{r,h}\cdot\tau_{0,-} =\mathcal{F}_{cl}^{5}\cdot \tau_{0,-} & \text{on }\Gamma_{-},\end{cases}\end{equation} where
\begin{equation}
\mathcal{F}_{r,h}^{1} =  P_{r,h}\mathcal{F}_{cl}^{1} - P_{r,h}\partial_\theta(Y_0e_\theta ) - P_{r,h}\mathscr CY_0e_\theta.\quad \mathcal{F}_{r,h}^{2}= \mathcal{F}_{cl}^{2}-\partial_{\theta}v_{\theta}.
\end{equation}

\subsection{The estimates for nonlinear terms}

\label{subsec:cl_nonlinear}

\subsubsection{The transverse elliptic estimate}

Let
\begin{equation*}
q_\delta:=\frac{2}{2-\delta}, \qquad 0<\delta<\delta_{\mathrm{max}}\le1, \qquad q_{r_{\delta}}=\frac{2}{1-\delta},
\end{equation*}
where $\delta_{\mathrm{max}}:=\varepsilon_{\max}(\theta_{\mathrm{eq}})$ is the corner exponent of Subsection~\ref{subsec:constants}, determined by the opening angle $\theta_{\mathrm{eq}}$ of the transverse wedge $\mathcal{W}$.

\begin{remark}
The artificial boundary is chosen so that the boundary of the cross-section $\mathcal{W}$ is smooth except at the corner, which is the trace of the contact line in the slice. As everywhere in this paper, only a~priori estimates are proved; the existence of solutions and the compatibility conditions are not discussed.
\end{remark}

We use the following corner estimate for the slice problem.

\begin{proposition}[The corner estimate]
\label{prop:corner_estimate} Let $0<\delta<\varepsilon_{\max}(\theta_{\mathrm{eq}})$ and $q_{\delta}=\frac{2}{2-\delta}$. For every fixed $\theta$, the solutions of the slice system \eqref{eq:two-dimensional-slice-system} on the cross-section $\mathcal{W}$ that vanish near its artificial boundary satisfy \eqref{eq:transverse-slice-estimate-r} below. The constant depends only on $\theta_{\mathrm{eq}}$, $\delta$, $\mu$, $\beta_{0}$ and $\mathcal{W} $; in particular it does not depend on $\theta$.
\end{proposition}

\begin{proof}
For fixed $\theta$, \eqref{eq:two-dimensional-slice-system} is the Stokes system in the two-dimensional cross-section $\mathcal{W}$, whose boundary is smooth except at the corner, of opening angle $\theta_{\mathrm{eq}}$, where $\Gamma_{+}$ and $\Gamma_{-}$ meet. On both faces the boundary conditions are of the same type: the normal velocity and the tangential stress are prescribed, and on $\Gamma_{-}$ the tangential condition contains in addition the lower-order friction term $\beta_{0}v_{r,h}\cdot\tau_{0,-}$. Away from the corner these conditions satisfy the complementing condition, and the estimate is the classical one. At the corner, write a divergence-free velocity as $v_{r,h}=\nabla^{\perp}\psi$. On a face, $v_{r,h}\cdot n=0$ gives $\psi=0$, hence $\partial_{\tau}^{2}\psi=0$, and then $(\mathbb{D}_{r,h}v_{r,h}\,n)\cdot\tau=\pm\partial_{n}^{2}\psi$. Up to lower-order terms, the homogeneous conditions therefore read $\psi=\Delta\psi=0$ on both faces. The biharmonic problem with these conditions splits into two Dirichlet problems for the Laplacian, and its leading singular function is $r^{\pi/\theta_{\mathrm{eq}}}\sin(\pi\varphi/\theta_{\mathrm{eq}})$. The corresponding velocity behaves like $r^{\pi/\theta_{\mathrm{eq}}-1}$, and it lies in $W^{2,q_{\delta}}$ near the corner exactly when $\pi/\theta_{\mathrm{eq}}-1>1-\frac{2}{q_{\delta}}+1=\delta$. Together with $q_{\delta}<2$, that is $\delta<1$, this is the condition $\delta<\varepsilon_{\max}(\theta_{\mathrm{eq}})$, and under it the estimate follows from the corner theory of Orlt and S\"andig~\cite{OS}, which is also the source of the exponent $\varepsilon_{\max}$ in the two-dimensional analysis of~\cite{GT}. These boundary conditions determine the pressure only up to a constant, which is why the left-hand side of \eqref{eq:transverse-slice-estimate-r} contains $p_{cl}-\overline{p}_{cl}(\theta)$. The cross-section, the angle and the coefficients of the principal part are the same for every $\theta$, so the constant is uniform in $\theta$.
\end{proof}

By Proposition~\ref{prop:corner_estimate}, for every fixed $\theta$,
\begin{align}
\begin{aligned} & \|v_{r,h}(\theta,\cdot)\|_{l_{\theta}^{2}W^{2,q_\delta}(\mathcal W)} + \|v_{r,h}(\theta,\cdot)\|_{L_{\theta}^{2}H^{1+\delta}(\mathcal W)} \\
&\qquad + \|p_{cl}(\theta,\cdot)-\overline p_{cl}(\theta)\|_{L_{\theta}^{2}W^{1,q_\delta}(\mathcal W)} + \|p_{cl}(\theta,\cdot)-\overline p_{cl}(\theta)\|_{L_{\theta}^{2}H^\delta(\mathcal W)} \\ &\lesssim \|\mathcal F_{r,h}^{1}(\theta,\cdot)\|_{L_{\theta}^{2}L^{q_\delta}(\mathcal W)} + \|\mathcal F_{r,h}^{2}(\theta,\cdot)\|_ {L_{\theta}^{2}W^{1,q_\delta}(\mathcal W)} \\ &\qquad + \|\mathcal{F}_{cl}^{4}(\theta,\cdot)\|_ {L_{\theta}^{2}W^{2-1/q_\delta,q_\delta}(\Gamma_+)}+\|\partial_{t}v_{r,h}\|_{L^{2}L^{2}}+\|\partial_{t}\eta_{\mathrm{cl}}\|_{L^{2}H^{\frac{3}{2}-\alpha}} \\ &\qquad+ \|\mathcal F_{cl}^{3}(\theta ,\cdot)\cdot\tau_{0,+}\|_ {L_{\theta}^{2}W^{1-1/q_\delta,q_\delta}(\Gamma_+)}+\|\mathcal F_{cl}^{5}(\theta ,\cdot)\cdot\tau_{0,-}\|_ {L_{\theta}^{2}W^{1-1/q_\delta,q_\delta}(\Gamma_-)}+\|\mathcal{F}_{\mathrm{cl},-}^{4}\|_{L_{\theta}^{2}W^{2-\frac{1}{q_{\delta}},q_{\delta}}(\Gamma_{-})}, \end{aligned}  \label{eq:transverse-slice-estimate-r}
\end{align}
where
\begin{equation*}
\overline p_{\mathrm{cl}}(\theta) := \frac{1}{|\mathcal{W}_{\mathrm{cl}}|} \int_{\mathcal{W}_{\mathrm{cl}}}p_{\mathrm{cl}}(\theta,s,t)\,ds\,dt.
\end{equation*} $\mathcal{W}_{\mathrm{cl}}$ is a bounded two-dimensional slicing domain for $r,h$. $\eta_{\mathrm{cl}}$ vanishes near its artificial boundary.

Substituting the definitions of $\mathcal{F}_{r,h}^{1}$ and $\mathcal{F}_{r,h}^{2}$, we obtain
\begin{align}
&\|v_{r,h}(\theta,\cdot)\|_ {L_{\theta}^{2}W^{2,q_\delta}\cap H^{1+\delta}} + \|p_{cl}(\theta,\cdot)-\overline p_{cl}(\theta)\|_ {L_{\theta}^{2}W^{1,q_\delta}\cap H^\delta} \lesssim \|\mathcal{F}_{cl}^{1}(\theta,\cdot)\|_{L_{\theta}^{2}L^{q_\delta}} + \|\partial_\theta(Y_0e_\theta)(\theta,\cdot)\|_{L_{\theta}^{2}L^{q_\delta}}\notag \\
&+ C\bigl(\|\mathscr C\|_{L_{\theta}^{\infty}L^\infty}\bigr) \|Y_0(\theta,\cdot)\|_{L_{\theta}^{2}L^{q_\delta}} + \|\mathcal{F}_{cl}^{2}(\theta,\cdot)\|_{L_{\theta}^{2}W^{1,q_\delta}} + \|\partial_\theta v(\theta,\cdot)\|_{L_{\theta}^{2}W^{1,q_\delta}}+\| v(\theta,\cdot)\|_{L_{\theta}^{2}H^{1}}+ \|p\|_{L_{\theta}^{2}L^{2}}\notag \\
&+ \|\p_{\theta}p\|_{L_{\theta}^{2}L^{q_{\delta}}} + \|\mathcal{F}_{cl}^{4}(\theta,\cdot)\|_ {L_{\theta}^{2}W^{2-1/q_\delta,q_\delta}(\Gamma_+)}+\|\partial_{t}v_{r,h}\|_{L^{2}L^{2}}+\|\partial_{t}\eta_{\mathrm{cl}}\|_{L^{2}H^{\frac{3}{2}-\alpha}}\label{eq:transverse-slice-expanded-r} \\
&+ \|\mathcal{F}_{cl}^{3}(\theta ,\cdot)\cdot\tau_{0,+}\|_ {L_{\theta}^{2}W^{1-1/q_\delta,q_\delta}(\Gamma_+)}+\|\mathcal{F}_{cl}^{5}(\theta ,\cdot)\cdot\tau_{0,-}\|_ {L_{\theta}^{2}W^{1-1/q_\delta,q_\delta}(\Gamma_-)}+\|\mathcal{F}_{\mathrm{cl},-}^{4}\|_{L_{\theta}^{2}W^{2-\frac{1}{q_{\delta}},q_{\delta}}(\Gamma_{-})}.\notag
\end{align}

By computation,  $\partial_\theta Y_0 = \partial_\theta S(p,v).$ Since $q_\delta<2$ and the transverse slice is bounded, we obtain
\begin{align}
\|\partial_\theta(Y_0e_\theta)\|_{L^2_\theta L^{q_\delta}(\mathcal W)} &\lesssim \|\partial_\theta p_{\mathrm{cl}}\|_{L^2_\theta L^{q_\delta}(\mathcal W)} +\|\partial_\theta v_{\mathrm{cl}}\|_{H^1(\mathbb T\times\mathcal W)}. \label{eq:partial-r-Y-estimate}
\end{align}
Here the last norm includes $\partial_\theta^2 v_{\mathrm{cl}}$.

Applying equation \eqref{eq:partial-r-Y-estimate} to equation \eqref{eq:transverse-slice-expanded-r}, we obtain
\begin{align}
& \|v_{r,h}(\theta,\cdot)\|_ {L_{\theta}^{2}W^{2,q_\delta}\cap H^{1+\delta}} + \|p_{cl}(\theta,\cdot)-\overline p_{cl}(\theta)\|_ {L_{\theta}^{2}W^{1,q_\delta}\cap H^\delta}  \notag \\
&\lesssim \|\mathcal{F}_{cl}^{1}(\theta,\cdot)\|_{L_{\theta}^{2}L^{q_\delta}} + \|\mathcal{F}_{cl}^{2}(\theta,\cdot)\|_{L_{\theta}^{2}W^{1,q_\delta}}  +\| \p_{\theta}v(\theta,\cdot)\|_{H^{1}}+\| v(\theta,\cdot)\|_{H^{1}}+ \|p\|_{L_{\theta}^{2}L^{2}}\notag \\
&\qquad+ \|\partial_{\theta}p\|_{L_{\theta}^{2}L^{q_{\delta}}}+ \|\mathcal{F}_{cl}^{4}(\theta,\cdot)\|_ {L_{\theta}^{2}W^{2-1/q_\delta,q_\delta}(\Gamma_+)}+\|\partial_{t}v_{r,h}\|_{L^{2}L^{2}}+\|\partial_{t}\eta_{\mathrm{cl}}\|_{L^{2}H^{\frac{3}{2}-\alpha}}  \notag \\
&\qquad+ \|\mathcal{F}_{cl}^{3}(\theta ,\cdot)\cdot\tau_{0,+}\|_ {L_{\theta}^{2}W^{1-1/q_\delta,q_\delta}(\Gamma_+)}+\|\mathcal{F}_{cl}^{5}(\theta ,\cdot)\cdot\tau_{0,-}\|_ {L_{\theta}^{2}W^{1-1/q_\delta,q_\delta}(\Gamma_-)}+\|\mathcal{F}_{\mathrm{cl},-}^{4}\|_{L_{\theta}^{2}W^{2-\frac{1}{q_{\delta}},q_{\delta}}(\Gamma_{-})}. \label{eq:transverse-final-slice-r}
\end{align}

In the following elliptic estimate, we have to separate the norm into $\theta $ and $r,h$ directions. For simplicity, we do not write the subscripts of the norms. For example, we will just write $L_{\theta}^{2}L_{r,h}^{2}$ as $L^{2}L^{2}$.

\textbf{Nonlinear forcing term estimates}

Before the estimate, we first derive the following geometric norm estimate near the contact line
\begin{align}
{\label{est:background_geometry}} \begin{aligned} &\|\chi_{\mathrm{cl}+}\mathcal{A}_{\eta}\|_{L^{\infty}L^{+\infty}}+\|\chi_{\mathrm{cl}+}\mathcal{A}_{\eta}\|_{L^{\infty}W^{1,q_{r_{\delta}}}}+\|\chi_{\mathrm{cl}+,\Sigma}\mathcal{N}\|_{L^{\infty}W^{1,q_{\delta}}(\Sigma)}+\|\chi_{\mathrm{cl}+}\partial_{\theta}\mathcal{A}_{\eta}\|_{L^{\infty}L^{q_{r_\delta}}}\\ &\lesssim 1+\|\bar{\eta}\|_{L^{\infty}W^{3,q_{\delta}}}+\|\partial_{\theta}\bar{\eta}\|_{L^{\infty}W^{2,q_{\delta}}}+\|\partial_{\theta}^{2}\bar{\eta}\|_{L^{\infty}H^{1-\alpha}}+\|\eta\|_{L^{\infty}W^{2,q_{\delta}}}+\|\partial_{\theta}\eta\|_{L^{\infty}H^{1}} \\
&\lesssim1+\|\chi_{cl,\Sigma}\eta\|_{L^{\infty}W^{3-\frac{1}{q_{\delta}},q_{\delta}}(\mathcal{W})}+\|\chi_{bet,\Sigma}\eta\|_{L^{\infty}W^{3-\frac{1}{q_{\delta}},q_{\delta}}(\mathcal{W})}\\
&\quad+\|\chi_{cl,\Sigma}\partial_{\theta}\eta\|_{L^{\infty}H^{\frac{3}{2}-\alpha}(\mathcal{W})}+\|\chi_{bet,\Sigma}\partial_{\theta}\eta\|_{L^{\infty}H^{\frac{3}{2}-\alpha}(\mathcal{W})}+\|\chi_{cl,\Sigma}\partial_{\theta}^{2}\eta\|_{L^{\infty}H^{\frac{1}{2}-\alpha}(\mathcal{W})}+\|\chi_{bet,\Sigma}\partial_{\theta}^{2}\eta\|_{L^{\infty}H^{\frac{1}{2}-\alpha}(\mathcal{W})}\\ &\lesssim 1+\sqrt{\mathscr{E}_{0,1,el,cl}}+\sqrt{\mathscr{E}_{0,1,el,bet}}+\sqrt{\mathscr{E}_{0,2,energy,cl}}+\sqrt{\mathscr{E}_{0,2,energy,bet}}:=1+\mathfrak{G}_{0,\eta}, \end{aligned}
\end{align}
where $\beta$ is a small positive constant to be determined in the interior estimate. Moreover, we used the fact that $\mathcal{A}$ near the contact line does not depend on the interior part.

For the right-hand side of inequality \eqref{eq:transverse-final-slice-r}, we estimate those nonlinear terms. First, for terms included in $\mathcal{F}_{cl}^{1}$, we estimate those included in $F$. By the definition of the momentum equation, we have

\textbf{Nonlinear convection term in the momentum equation}
\begin{align}
\|\chi_{\mathrm{cl}}\partial_{t}v\|_{L^{2}L^{q_{\delta}}}\lesssim\|\partial_{t}(\chi_{\mathrm{cl}}v)\|_{L^{2}L^{2}}\lesssim\sqrt{\mathscr{E}_{1,0,\mathrm{energy},\mathrm{cl}}}.
\end{align}
\begin{align}
&\|\chi_{\mathrm{cl}}(v\cdot\nabla_{Q^{T}\mathcal{A},\mathcal{A}_{\eta}} v)\|_{L^{2}L^{q_{\delta}}}\lesssim \|\chi_{\mathrm{cl}}v\|_{L^{\infty}L^{\infty}}(\|\chi_{\mathrm{cl}}v\|_{L^{2}H^{1}}+\|\chi_{\mathrm{bet}}v\|_{L^{2}H^{1}}+\|\chi_{\mathrm{cl}}\partial_{\theta}v\|_{L^{2 }L^{2}}+\|\chi_{\mathrm{bet}}\partial_{\theta}v\|_{L^{2}L^{2}})\notag \\
&\times(1+\|\mathcal{A}_\eta\|_{L^{\infty}L^{+\infty}})\label{eq:convection_cl} \\
&\lesssim \|\chi_{\mathrm{cl}}v\|_{L^{\infty}W^{2,q_{\delta}}}(\|\chi_{\mathrm{cl}}v\|_{L^{2}H^{1}}+\|\chi_{\mathrm{bet}}v\|_{L^{2}H^{1}}+\|\chi_{\mathrm{cl}}\partial_{\theta}v\|_{L^{2}L^{2}}+\|\chi_{\mathrm{bet}}v\|_{L^{2}L^{2}})(\mathfrak{G}_{0,\eta}+1)\notag \\
&\lesssim \sqrt{\mathscr{E}_{0,1,\mathrm{el},\mathrm{cl}}}(\sqrt{\mathscr{E}_{0,1,\mathrm{energy},\mathrm{cl}}}+\sqrt{\mathscr{E}_{0,1,\mathrm{energy},\mathrm{bet}}}+\sqrt{\mathscr{E}_{0,2,\mathrm{energy},\mathrm{cl}}}+\sqrt{\mathscr{E}_{0,2,\mathrm{energy},\mathrm{bet}}}) \lesssim \mathcal{E}.\notag
\end{align}

\textbf{Cutoff commutators in the momentum equation}

Put $\chi_*:=\chi_{\mathrm{cl}}+\chi_{\mathrm{bet}}$, which equals one on $\operatorname{supp}\chi_{\mathrm{cl}}$. On this support the stress involves $p$ and first derivatives of $v$. By \eqref{eq:cutoff_scaling}, a first derivative of $\chi_{\mathrm{cl}}$ costs $C/\varepsilon$. Since $q_\delta<2$, H\"older's inequality on the fixed bounded transverse slice gives
\begin{align}
&\|\nabla_{Q^T\mathcal A,I}\chi_{\mathrm{cl}}\cdot S_{Q^T\mathcal A_\Phi}(p,v)\|_{L^2_\theta L^{q_\delta}_{r,h}} \lesssim \frac{C(\varepsilon)}{\varepsilon} \bigl(\|\chi_*v\|_{H^1}+\|\chi_*p\|_{L^2}\bigr) \lesssim \frac{C(\varepsilon)}{\varepsilon} \sqrt{\mathscr E_{0,0,\mathrm{energy},\mathrm{cl}}},\notag \\
&\|\chi_{\mathrm{cl}}\partial_t\bar\eta\,WK\partial_3v\|_{L^2_\theta L^{q_\delta}_{r,h}} \lesssim \|\chi_{\mathrm{cl}}\partial_t\bar\eta\|_{L^\infty} \|\chi_*v\|_{H^1} \lesssim \|\partial_t\eta\|_{H^{3/2-\alpha}(\Sigma)}\|\chi_*v\|_{H^1} \lesssim \mathcal E, \\
&\|v\Delta_{Q^{T}\mathcal{A},I}\chi_{\mathrm{cl}}\|_{L^{2}L^{q_{\delta}}}\lesssim C(\varepsilon)\frac{C(\varepsilon)}{\varepsilon}(\|\chi_{\mathrm{cl}}v\|_{L^{2}}+\|\chi_{\mathrm{bet}}v\|_{L^{2}})\lesssim\frac{C(\varepsilon)}{\varepsilon}\sqrt{\mathscr{E}_{0,0,\mathrm{energy},\mathrm{cl}}},\notag
\end{align}
where $C(\varepsilon)$ is a large constant depending on $\varepsilon$. For the second bound, all extension cutoffs and chart maps in \eqref{eq:eta_bar_def} are fixed in time. The ordinary interior Poisson extension is multiplied by $\chi_{\mathrm{int}+}$, whose support is disjoint from $\operatorname{supp}\chi_{\mathrm{cl}}$ by the cutoff nesting. For each remaining term of \eqref{eq:contact-line-extension} and \eqref{eq:interior-extension}, the surface datum is a bounded cutoff times $\partial_t\eta$. The finite reflection formula \eqref{eq:half-space-extension} is bounded on $L^\infty$, and each radial Poisson kernel has integral one. These bounds use only the suprema of the cutoffs, not their derivatives. Thus $\|\chi_{\mathrm{cl}}\partial_t\bar\eta\|_{L^\infty} \lesssim \|\partial_t\eta\|_{L^\infty(\Sigma)}$. Since $3/2-\alpha>1$ and $\Sigma$ is two-dimensional, the latter supremum is bounded by the $H^{3/2-\alpha}(\Sigma)$ norm in $\mathcal E_{\mathrm{enhanced}}$.

\textbf{Geometric Difference Terms caused by $\mathcal{A}_{\eta}$ in the momentum equation}
\begin{align}
&\|\operatorname{div}_{Q^{T}\mathcal{A}_{\Phi},I-\mathcal{A}_{\eta}}S_{Q^{T}\mathcal{A}_{\Phi},\mathcal{A}_{\eta}}(\chi_{\mathrm{cl}}v,\chi_{\mathrm{cl}}p)\|_{L^{2}L^{q_{\delta}}} \lesssim \|(\chi_{cl}+\chi_{bet})\nabla \bar{\eta}\|_{L^{\infty}L^{+\infty}}(1+\|\mathcal{A}_{\eta}\|_{L^{\infty}L^{+\infty}})\notag \\
&(\|\chi_{\mathrm{cl}}v\|_{L^{2}W^{2,q_{\delta}}}+\|\chi_{\mathrm{cl}}p\|_{L^{2}W^{1,q_{\delta}}}+\|\chi_{\mathrm{cl}}\partial_{\theta}^{2}v\|_{L^{2}L^{2}} +\|\chi_{\mathrm{cl}}\partial_{\theta}p\|_{L^{2}L^{2}}+\|\chi_{\mathrm{cl}}\partial_{\theta}v\|_{L^{2}H^{1}})\notag \\
&+\|\nabla \bar{\eta}\|_{L^{\infty}L^{\infty}}(1+\|\nabla\mathcal{A}_{\eta}\|_{L^{\infty}W^{1,q_{r_{\delta}}}})(\|\chi_{\mathrm{cl}}v\|_{L^{2}H^{1}}+\|\chi_{\mathrm{cl}}p\|_{L^{2}L^{2}}+\|\chi_{\mathrm{cl}}\partial_{\theta}v\|_{L^{2}L^{2}})\notag \\
&\lesssim \bigl(1+\mathfrak{G}_{0,\eta}\bigl)\mathfrak{G}_{0,\eta}L_{\theta}^{2}(\|\chi_{\mathrm{cl}}v\|_{W^{2,q_{\delta}}}+\|\chi_{\mathrm{cl}}p\|_{W^{1,q_{\delta}}}+\|\chi_{\mathrm{cl}}\partial_{\theta}^{2}v\|_{L^{2}}\notag \\
&+\|\chi_{\mathrm{cl}}\partial_{\theta}p\|_{L^{2}}+\|\chi_{\mathrm{cl}}v\|_{H^{1}}+\|\chi_{\mathrm{cl}}p\|_{L^{2}}+\|\chi_{\mathrm{cl}}\partial_{\theta}v\|_{H^{1}}) \\
&\lesssim\mathfrak{G}_{0,\eta}(\sqrt{\mathscr{E}_{0,1,\mathrm{energy},\mathrm{cl}}}+\sqrt{\mathscr{E}_{0,1,\mathrm{energy},\mathrm{bet}}}+\sqrt{\mathscr{E}_{0,2,\mathrm{energy},\mathrm{cl}}}+\sqrt{\mathscr{E}_{0,2,\mathrm{energy},\mathrm{bet}}}+\sqrt{\mathscr{E}_{0,1,\mathrm{el},\mathrm{cl}}}+\sqrt{\mathscr{E}_{0,1,\mathrm{el},\mathrm{bet}}})\lesssim {\mathcal{E}}.\notag
\end{align}
Similarly, we have
\begin{align}
&\|\nabla_{Q^{T}\mathcal{A}_{\Phi},I-\mathcal{A}_{\eta}}\chi_{\mathrm{cl}}\cdot S_{Q^{T}\mathcal{A},\mathcal{A}_{\eta}}(v,p)\|_{L^{2}L^{q_{\delta}}} \lesssim \bigl(1+\mathfrak{G}_{0,\eta}\bigl)\mathfrak{G}_{0,\eta}(\|\chi_{cl}v\|_{L^{2}H^{1}}+\|\chi_{bet}v\|_{L^{2}H^{1}}+\|\chi_{cl}p\|_{L^{2}L^{2}}+\|\chi_{bet}p\|_{L^{2}L^{2}})\notag \\
&\lesssim \bigl(1+\mathfrak{G}_{0,\eta}\bigl)\mathfrak{G}_{0,\eta}(\sqrt{\mathscr{E}_{0,0,\mathrm{energy},\mathrm{cl}}}+\sqrt{\mathscr{E}_{0,0,\mathrm{energy},\mathrm{bet}}})\lesssim\mathcal{E}.
\end{align}
Similarly,
\begin{align}
\|(\func{div}_{Q^{T}\mathcal{A},I-\mathcal{A}_{\eta}}\nabla_{Q^{T}\mathcal{A},\mathcal{A}_{\eta}}\chi_{cl})v\|_{L^{2}L^{q_{\delta}}}\lesssim \mathfrak{G}_{0,\eta}(\|\chi_{cl}v\|_{L^{2}L^{2}}+\|\chi_{bet}v\|_{L^{2}L^{2}})\lesssim {\mathcal{E}}.
\end{align}
Similarly, we have
\begin{align}
&\|\operatorname{div}_{Q^{T}\mathcal{A},\mathcal{A}_{\eta}}\nabla_{Q^{T}\mathcal{A},I-\mathcal{A}_{\eta}}\chi_{cl}v\|_{L^{2}L^{q_{\delta}}}+\|\nabla_{Q^{T}\mathcal{A},\mathcal{A}_{\eta}}\chi_{cl}\cdot S_{Q^{T}\mathcal{A},I-\mathcal{A}_{\eta}}(v,p)\|_{L^{2}L^{q_{\delta}}} \\
&+\|\operatorname{div}_{Q^{T}\mathcal{A},\mathcal{A}_{\eta}}S_{Q^{T}\mathcal{A},I-\mathcal{A}_{\eta}}(\chi_{cl}v,\chi_{cl}p)\|_{L^{2}L^{q_{\delta}}} \lesssim \bigl(1+\mathfrak{G}_{0,\eta}\bigl)\mathfrak{G}_{0,\eta}(\sqrt{\mathscr{E}_{0,1,\mathrm{energy},\mathrm{cl}}}+\sqrt{\mathscr{E}_{0,1,\mathrm{energy},\mathrm{bet}}}+\sqrt{\mathscr{E}_{0,0,\mathrm{el},\mathrm{cl}}})\lesssim \mathcal{E}.\notag
\end{align}

\textbf{Flattening Difference Terms and the Corresponding Elliptic Estimate for the momentum equation}

Let
\begin{equation*}
\mathcal{B}:=Q^{T}\mathcal{A}_{\Phi}=I+\mathcal{Z}, \qquad \mathcal{Z}:=Q^{T}\mathcal{A}_{\Phi}-I, \qquad \mathscr{C}:=Q^{T}\partial_{\theta}Q.
\end{equation*}
Let  $P:=P_{r,h}=I-e_{\theta}\otimes e_{\theta}, \, v_{r,h}:=Pv.$ We define  $(\mathbb{D}_{\mathcal{Z}}v)_{ij} := \mathcal{Z}_{ik}\partial_{k}v_{j} + \mathcal{Z}_{jk}\partial_{k}v_{i}.$ Since $B=I+\mathcal{Z}$, we have  $\mathbb{D}_{\mathcal{B}}v = \mathbb{D}v+\mathbb{D}_{\mathcal{Z}}v,$ and hence  $S_{\mathcal{B}}(p,v) = S(p,v)-\mu\mathbb{D}_{\mathcal{Z}}v.$

\subsubsection{The projected equation}

For $a\in\{r,h\}$, set $()_{a}$ to be the $a$-component of the two-dimensional vector.
\begin{equation}
\left( P\mathscr{C}S(p,v)e_{\theta} \right)_{a}  = P_{ai}\mathscr{C}_{ij}S_{j\theta}(p,v) = P_{ai}\mathscr{C}_{i\theta} \left( p-2\mu\partial_{\theta}v_{\theta} \right)   \quad - \mu P_{ai} \sum_{b\in\{r,h\}} \mathscr{C}_{ib} \left( \partial_{b}v_{\theta} + \partial_{\theta}v_{b} \right).\label{eq:connection-term-expansion}
\end{equation}
Thus the projected principal operator contains both  $\partial_{\theta}^{2}v_{r,h}$ and the first-order terms
\begin{equation*}
-\mu P_{ai} \sum_{b\in\{r,h\}} \mathcal{C}_{ib} \partial_{\theta}v_{b}.
\end{equation*}

\subsubsection{Explicit decomposition of the bulk difference}

Define $\mathcal{F}_{\mathrm{flat}}^{1}$ to be the sum of all forcing terms caused by the flattening map. Write
\begin{equation*}
P_{r,h}\mathcal{F}_{\mathrm{flat}}^{1} = \sum_{m=1}^{7}\mathcal{F}^{(m)},
\end{equation*}
where
\begin{align}
&\mathcal{F}^{(1)}  := -P_{r,h}\func{div}_{\mathcal{Z}}S(\chi_{cl}p,\chi_{cl}v),\quad \mathcal{F}^{(2)}  := \mu P_{r,h}\func{div}_{B}\mathbb{D}_{\mathcal{Z}}(\chi_{cl}v),\quad \mathcal{F}^{(3)}  := -P_{r,h}\func{div}_{B}R(\chi_{cl}v),\notag \\
&\mathcal{F}^{(4)}  := -P_{r,h} \bigl(S_{B}(\chi_{cl}p,\chi_{cl}v)+R(\chi_{cl}v)\bigr) \partial_{\theta}Q^{T}\mathcal{A}_{\Phi}\,e_{\theta},\quad \mathcal{F}^{(5)}  := \mu P_{r,h}\mathfrak{C}\mathbb{D}_{\mathcal{Z}}(\chi_{cl}v)\,e_{\theta},\notag \\
&\mathcal{F}^{(6)}  := -P_{r,h}\mathfrak{C}R(\chi_{cl}v)e_{\theta},\quad \mathcal{F}^{(7)}  := -P_{r,h}\mathfrak{C} \bigl(S_{B}(\chi_{cl}p,\chi_{cl}v)+R(\chi_{cl}v)\bigr) \mathcal{Z}e_{\theta}.
\end{align}

\subsubsection{\texorpdfstring{\(L^{q_{\delta}}\)}{Lq} estimates}

Let (in the following, we ignore the $L_{\theta}^{\infty}(\mathcal{T}) $ norm)
\begin{equation*}
q_{\delta} = \frac{2}{2-\delta}, \qquad r_{\delta} = \frac{2}{1-\delta},\qquad m_{\delta}=\frac{1}{1-\delta} \qquad 0<\delta<1.
\end{equation*}
Then  $\frac{1}{q_{\delta}} = \frac{1}{2} + \frac{1}{r_{\delta}}.$

Set
\begin{equation*}
M_{\mathcal{Z}} := \|\mathcal{Z}\|_{L^{\infty}(\func{supp}\chi_{\mathrm{bet}+})}+\|\partial_{\theta}\mathcal{Z}\|_{L^{\infty}(\func{supp}\chi_{\mathrm{bet}+})}, \qquad N_{E} := \|\nabla \mathcal{Z}\|_{L^{\infty}_{\theta}L^{\infty}(\mathcal{W})}.
\end{equation*}
The matrix $\mathcal{Z}$ is evaluated only on the support of the localized unknowns, and there $|r|+|h|<\varepsilon$ by \eqref{eq:cutoff_eps}. Hence \eqref{eq:QA_small} and \eqref{eq:QA_theta_small} give
\begin{equation}  \label{eq:MZ_small} M_{\mathcal{Z}}\leq C\varepsilon .
\end{equation}
The gradient $\nabla\mathcal{Z}$ is not small; it is bounded by $N_{E}$, which depends only on $\zeta$ and $\partial D$, and  $M_{\mathscr C} := \|\mathscr C\|_{L_{\theta}^{\infty}L^{\infty}(\mathcal{W})}.$ Since $B=I+\mathcal{Z}$,  $\|\mathcal{B}\|_{L^{\infty}} \leq 1+M_{\mathcal{Z}}.$

For the first term, we have
\begin{align}
\|\mathcal{F}^{(1)}\|_{L^{2}L^{q_{\delta}}} &\lesssim \|\mathcal{Z}\|_{L^{\infty}L^{\infty}} \left( \|\nabla (\chi_{\mathrm{cl}}p)\|_{L^{2}L^{q_{\delta}}} + \|\nabla^{2}(\chi_{\mathrm{cl}}v)\|_{L^{2}L^{q_{\delta}}} \right)  \notag \\
&\lesssim M_{\mathcal{Z}} \left( \|\chi_{\mathrm{cl}}p\|_{L^{2}W^{1,q_{\delta}}} +\|\partial_{\theta}(\chi_{\mathrm{cl}}p)\|_{L^{2}L^{2}}+ \|\chi_{\mathrm{cl}}v\|_{L^{2}W^{2,q_{\delta}}} +\|\partial_{\theta}(\chi_{\mathrm{cl}}v)\|_{L^{2}H^{1}}+\|\partial_{\theta}^{2}(\chi_{\mathrm{cl}}v)\|_{L^{2}L^{2}}\right)  \notag \\
&\lesssim M_{\mathcal{Z}}(\sqrt{\mathscr{E}_{0,0,\mathrm{energy},\mathrm{cl}}}+\sqrt{\mathscr{E}_{0,0,\mathrm{el},\mathrm{cl}}}+\sqrt{\mathscr{E}_{0,1,\mathrm{energy},\mathrm{cl}}})\lesssim M_{\mathcal{Z}}\sqrt{\mathcal{E}}. \label{eq:F1-estimate}
\end{align}
Here we notice a major difference while doing the estimate in 3-D. For example, if we want to estimate $\|\nabla u\|_{L_{\theta}^{\infty}}$ in some nonlinear terms, we have to decompose $\nabla$ into a $\theta$-derivative and an $r,h$-derivative.

For the second term, H\"older's inequality gives
\begin{align}
&\|\mathcal{F}^{(2)}\|_{L^{2}L^{q_{\delta}}} \lesssim{}  \|\mathcal{B}\|_{L^{\infty}L^{\infty}} \|\nabla \mathcal{Z}\|_{L^{\infty}L^{\infty}} \|\nabla (\chi_{\mathrm{cl}}v)\|_{L^{2}L^{q_{{\delta}}}} + \|B\|_{L^{\infty}L^{\infty}} \|\mathcal{Z}\|_{L^{\infty}L^{\infty}} \|\nabla^{2}(\chi_{\mathrm{cl}}v)\|_{L^{2}L^{q_{\delta}}}\notag \\
&\lesssim{}  (1+M_{\mathcal{Z}})N_{\mathcal{Z}}(\| \chi_{cl}v\|_{H^{1}}+\|\partial_{\theta}(\chi_{\mathrm{cl}}v)\|_{L^{2}}) +(1+M_{\mathcal{Z}})M_{\mathcal{Z}}\left( \|\chi_{\mathrm{cl}}p\|_{L^{2}W^{1,q_{\delta}}} + \|\chi_{\mathrm{cl}}v\|_{L^{2}W^{2,q_{\delta}}} +\|\partial_{\theta}(\chi_{\mathrm{cl}}v)\|_{H^{1}}\right)\notag \\
&\lesssim{}  (N_{\mathcal{Z}})(\sqrt{\mathscr{E}_{0,0,\mathrm{energy},\mathrm{cl}}})+M_{\mathcal{Z}}\sqrt{\mathscr{E}_{0,0,el,cl}+\mathscr{E}_{0,1,energy,cl}}.\label{eq:F2-estimate}
\end{align}

For the third term,
\begin{align}
\|\mathcal{F}^{(3)}\|_{L^{2}L^{q_{\delta}}} &\lesssim \|\mathcal{B}\|_{L^{\infty}L^{\infty}} \|\nabla R\|_{L^{2}L^{q_{\delta}}} \lesssim (1+M_{\mathcal{Z}}) \|R\|_{L^{2}W^{1,q_{\delta}}}\lesssim(1+M_{\mathcal{Z}})\|v_{\mathrm{cl}}\|_{L^{2}W^{1,q_{\delta}}}  \notag \\
&\lesssim (1+M_{\mathcal{Z}})\|v_{\mathrm{cl}}\|_{H^{1}}\lesssim (1+M_{\mathcal{Z}})(\sqrt{\mathscr{E}_{0,0,\mathrm{energy},\mathrm{cl}}}). \label{eq:F3-estimate}
\end{align}

For the fourth term,
\begin{align}
&\|\mathcal{F}^{(4)}\|_{L^{2}L^{q_{\delta}}} \lesssim{}  \|\partial_{\theta}Q^{T}\mathcal{A}_{\Phi}\|_{L^{\infty}L^{\infty}} \left( \|\chi_{\mathrm{cl}}p\|_{L^{2}L^{q_{\delta}}} + \|\nabla (\chi_{\mathrm{cl}}v)\|_{L^{2}L^{q_{\delta}}} + \|\mathbb{D}_{\mathcal{Z}}(\chi_{\mathrm{cl}}v)\|_{L^{2}L^{q_{\delta}}} + \|R(\chi_{\mathrm{cl}}v)\|_{L^{2}L^{q_{\delta}}} \right)\notag \\
&\lesssim{}  \left[ \|\chi_{cl}p\|_{L^{2}} + (1+M_{\mathcal{Z}})(\| \chi_{cl}v\|_{H^{1}}+\|\partial_{\theta}(\chi_{cl}v)\|_{L^{2}L^{2}}) \right] \lesssim (\sqrt{\mathscr{E}_{0,0,\mathrm{energy},cl}}).\label{eq:F4-estimate}
\end{align}
Now we introduce notations
\begin{align}
\mathscr{U}_{1}:=&\|\chi_{cl}v\|_{L_{\theta}^{2}H^{1}}+\|\partial_{\theta}(\chi_{cl}v)\|_{L_{\theta}^{2}L^{2}}\lesssim\sqrt{\mathscr{E}_{0,0,\mathrm{energy},\mathrm{cl}}}, \\
\qquad \mathscr{U}_{2}:=&\|\chi_{cl}v\|_{L_{\theta}^{2}W^{2,q_{\delta}}}+\|\partial_{\theta}(\chi_{cl}v)\|_{L_{\theta}^{2}H^{1}}+\|\partial_{\theta}^{2}(\chi_{cl}v)\|_{L_{\theta}^{2}L^{2}}\lesssim\sqrt{\mathscr{E}_{0,0,\mathrm{el}.\mathrm{cl}}}+\sqrt{\mathscr{E}_{0,1,\mathrm{energy},\mathrm{cl}}}.
\end{align}

For the fifth term,
\begin{align}
\|\mathcal{F}^{(5)}\|_{L^{2}L^{q_{\delta}}} &\lesssim \|\mathscr C\|_{L^{\infty}L^{2}} \|\mathcal{Z}\|_{L^{\infty}L^{\infty}} \|\nabla (\chi_{\mathrm{cl}}v)\|_{L^{2}L^{q_{r_\delta}}} \lesssim M_{\mathcal{C}}M_{\mathcal{Z}}\mathscr{U}_{2}\lesssim M_{\mathcal{C}}M_{\mathcal{Z}}(\sqrt{\mathscr{E}_{0,1,\mathrm{energy},cl}}+\sqrt{\mathscr{E}_{0,0,\mathrm{el},cl}}). \label{eq:F5-estimate}
\end{align}

For the sixth term,
\begin{align}
\|\mathcal{F}^{(6)}\|_{L^{2}L^{q_{\delta}}} &\lesssim \|\mathscr C\|_{L^{\infty}L^{\infty}} \|R(\chi_{cl}v)\|_{L^{2}L^{2}} \lesssim M_{\mathscr C}\sqrt{\mathscr{E}_{0,0,\mathrm{energy},cl}}. \label{eq:F6-estimate}
\end{align}

For the seventh term,
\begin{align}
&\|\mathcal{F}^{(7)}\|_{L^{2}L^{q_{\delta}}} \lesssim{}  \|\mathscr C\|_{L^{\infty}L^{\infty}} \|\mathcal{Z}\|_{L^{\infty}L^{\infty}} \Big( \|\chi_{\mathrm{cl}}p\|_{L^{2}L^{q_{\delta}}} + \|\nabla (\chi_{\mathrm{cl}}v)\|_{L^{2}L^{q_{\delta}}} + \|\mathbb{D}_{\mathcal{Z}}v\|_{L^{2}L^{q_{\delta}}} + \|R(\chi_{\mathrm{cl}}v)\|_{L^{2}L^{q_{\delta}}} \Big)\label{eq:F7-estimate} \\
&\lesssim{}  M_{\mathscr C}M_{\mathcal{Z}} \Big[ \|\chi_{\mathrm{cl}}p\|_{L^{2}W^{1,q_{\delta}}} + (1+M_{\mathcal{Z}}) (\|\chi_{\mathrm{cl}}v\|_{L^{2}W^{2,q_{\delta}}}+\|\partial_{\theta}(\chi_{\mathrm{cl}}v)\|_{L^{2}H^{1}}) \Big] \lesssim{}  M_{\mathscr C}M_{\mathcal{Z}} \Big[ (\sqrt{\mathscr{E}_{0,1,\mathrm{energy},cl}}+\sqrt{\mathscr{E}_{0,0,\mathrm{el},cl}}) \Big].\notag
\end{align}

Combining \eqref{eq:F1-estimate}--\eqref{eq:F7-estimate}, we obtain
\begin{align}
\begin{aligned} \|P_{r,h}\mathcal F_{\mathrm{flat}}\|_{L^{q_{\delta}}} \lesssim{}& C_{\mathscr C} \left( M_{\mathcal{Z}}+M_{\mathcal{Z}}^{2} \right)\Big[ (\sqrt{\mathscr{E}_{0,1,\mathrm{energy},cl}}+\sqrt{\mathscr{E}_{0,0,\mathrm{el},cl}}) \Big]+(1+C_{\mathscr C})(\sqrt{\mathscr{E}_{0,0,\mathrm{energy},\mathrm{cl}}}), \label{eq:combined-bulk-estimate} \end{aligned}
\end{align}
where  $C_{\mathscr C} := C( 1+M_{\mathscr C} ) = C( 1+\|\mathscr C\|_{L^{2}} ).$

\subsubsection{Difference of the divergence equations}

Recall $H_{\mathrm{cl}}$ from \eqref{eq:local-F2-final}. Since $\mathcal B=I+\mathcal Z$, the transformed and frozen divergence equations are
\begin{equation}
\operatorname{div}_{\mathcal B}v_{\mathrm{cl}}=H_{\mathrm{cl}},\qquad \operatorname{div}v_{\mathrm{cl}}=H_{\mathrm{cl}}+\mathscr G=\mathcal F_{\mathrm{cl}}^2, \qquad \mathscr G:=-\operatorname{div}_{\mathcal Z}v_{\mathrm{cl}}. \label{eq:frozen-divergence-equation}
\end{equation}
Explicitly,  $\mathscr G = -\mathcal{Z}_{ik}\partial_{k}v_{cl,i},$ and
\begin{equation*}
\partial_{\ell}\mathscr G = -(\partial_{\ell}\mathcal{Z}_{ik})\partial_{k}v_{cl,i} - \mathcal{Z}_{ik}\partial_{\ell k}v_{cl,i}.
\end{equation*}
Consequently,
\begin{align}
\|\mathscr G\|_{L^{2}W^{1,q_{\delta}}} \lesssim{}& M_{\mathcal{Z}} \|\chi_{cl}v\|_{L^{2}W^{2,q_{\delta}}} + \|\partial_{\theta} (\chi_{cl}v)\|_{L^{2}W^{1,q_{\delta}}}+\|\nabla \mathcal{Z}\|_{L^{\infty}}\|v_{\mathrm{cl}}\|_{H^{1}} \\
&\lesssim M_{\mathcal Z}\sqrt{\mathscr E_{0,0,\mathrm{el},\mathrm{cl}}} +\sqrt{\mathscr E_{0,1,\mathrm{energy},\mathrm{cl}}} +C(\varepsilon)\sqrt{\mathscr E_{0,0,\mathrm{energy},\mathrm{cl}}}. \label{eq:divergence-difference-estimate}
\end{align}

It remains to estimate $H_{\mathrm{cl}}$. We have
\begin{align}
&\|H_{\mathrm{cl}}\|_{L^{2}W^{1,q_{\delta}}}\lesssim \|\chi_{cl}\operatorname{div}_{Q^{T}\mathcal{A},I-\mathcal{A}_{\eta}}v\|_{L^{2}W^{1,q_{\delta}}}+\|\nabla _{Q^{T}\mathcal{A},\mathcal{A}_{\eta}}\chi_{cl}v\|_{L^{2}W^{1,q_{\delta}}}+\|\partial_{\theta}Qv_{cl}\cdot (\mathcal{A}_{\Phi}e_{\theta})\|_{L^{2}W^{1,q_{\delta}}}\notag \\
&\lesssim \mathfrak{G}_{\eta}(\|\chi_{\mathrm{cl}}v\|_{L^{2}W^{2,q_{\delta}}}+\|\chi_{\mathrm{cl}}\partial_{\theta}v\|_{L^{2}H^{1}})+(1+\mathfrak{G}_{\eta})(\|\nabla\chi_{\mathrm{cl}}\cdot v \|_{L^{2}W^{1,q_{\delta}}}+\| v_{cl} \|_{L^{2}W^{1,q_{\delta}}}) \\
&\lesssim \mathfrak{G}_{\eta}(\|\chi_{\mathrm{cl}}v\|_{L^{2}W^{2,q_{\delta}}}+\|\chi_{\mathrm{cl}}\partial_{\theta}v\|_{L^{2}H^{1}})+\frac{C(\varepsilon)}{\varepsilon}(1+\mathfrak{G}_{\eta})(\|\chi_{\mathrm{cl}}v\|_{L^{2}H^{1}}+\|\chi_{\mathrm{bet}}v\|_{L^{2}H^{1}})\notag \\
&\lesssim \mathfrak{G}_{\eta}(\sqrt{\mathscr{E}_{0,0,\mathrm{el},\mathrm{cl}}}+\sqrt{\mathscr{E}_{0,1,\mathrm{energy},\mathrm{cl}}}+\sqrt{\mathscr{E}_{0,0,\mathrm{el},\mathrm{bet}}})+\frac{C(\varepsilon)}{\varepsilon}\sqrt{\mathscr{E}_{0,0,\mathrm{energy},\mathrm{cl}}} \lesssim \mathcal{E}+\frac{C(\varepsilon)}{\varepsilon} \sqrt{\mathscr{E}_{0,0,\mathrm{energy},\mathrm{cl}}}.\notag
\end{align}
To bound $\mathcal F_{\mathrm{cl}}^2$ in the slice estimate, we add \eqref{eq:divergence-difference-estimate} to this bound for $H_{\mathrm{cl}}$, as required by \eqref{eq:frozen-divergence-equation}.

\subsubsection{Difference of the normal boundary conditions}

Define  $\nu_{\mathcal{B}} := \frac{\mathcal{B}n_{0}}{|\mathcal{B}n_{0}|}, \, P_{0} := I-n_{0}\otimes n_{0}.$ The variable normal boundary condition on the free boundary $\Gamma_{+}$ is
\begin{equation*}
\chi_{\mathrm{cl}}v_{r,h}\cdot n_{0,+}=\frac{1}{|\mathcal{N}_{0}|}\partial_{t}\chi_{\mathrm{cl}}\eta-\mathcal{F}_{cl}^{4}.
\end{equation*}
We now estimate terms included in $\mathcal{F}_{cl}^{4}$.

\textbf{Flattening Difference term}

The smoothness of  $\mathscr{B}\longmapsto\frac{\mathscr{B}n_{0}}{|\mathscr{B}n_{0}|}$ near $B=I$ implies  $\nu_{\mathscr{B}}-n_{0} = P_{0}\mathcal{Z}n_{0} + \mathcal{Q}_{\nu}(\mathcal{Z}),$ where  $|\mathcal{Q}_{\nu}(\mathcal{Z})| \lesssim |\mathcal{Z}|^{2}.$ (Taylor expansion) Moreover, by Taylor expansion, and the fact that $\chi_{\mathrm{bet}+}$ is a scaled cutoff function near the contact line, we note that  $\| \chi_{\mathrm{bet}+}(\nu_{\mathscr{B}}-n_{0})\|_{W_{r,h}^{s,q}}\lesssim C(\varepsilon^{1+d/q-s}+\varepsilon^{1-s}),$ for any $1\leq q\leq \infty$ and $0<s<1$. Consequently, using the fact that $\chi_{\mathrm{cl},\Sigma}+\chi_{\mathrm{bet},\Sigma}$ is supported in a rim of the surface near the contact line, we have
\begin{align}
\|\chi_{\mathrm{cl},\Sigma}v\cdot (n_0-\nu_{\mathcal{B}})\|_ {L^{2}W^{2-1/q_{\delta},q_{\delta}}(\Sigma)}\lesssim&\|\chi_{\mathrm{cl},\Sigma}v\|_{L^{2}W^{1-\frac{1}{q_{\delta}},q_{\delta}}(\Sigma)}\| (n_0-\nu_{\mathcal{B}})\|_ {L^{\infty}W^{3-1/q_{\delta},q_{\delta}}(\Sigma)} \\
&+\|\chi_{\mathrm{cl},\Sigma}v\|_{L^{2}W^{2-\frac{1}{q_{\delta}},q_{\delta}}(\Sigma)}\|(\chi_{\mathrm{cl},\Sigma}+\chi_{\mathrm{bet},\Sigma}) (n_0-\nu_{\mathcal{B}})\|_ {L^{\infty}W^{\frac{1}{q_{\delta}}+\alpha,q_{\delta}}(\Sigma)} \\
\lesssim{}&\varepsilon^{\iota} \|\chi_{\mathrm{cl}}v\|_{L^{2}W^{2,q_{\delta}}(\mathcal{W})}+\|\chi_{cl}v\|_{H^{1}}\lesssim \varepsilon^{\iota}\sqrt{\mathscr{E}_{0,0,\mathrm{el},\mathrm{cl}}}+\sqrt{\mathscr{E}_{0,0,energy,cl}},  \label{eq:normal-boundary-estimate}
\end{align}
where we choose $\frac{1}{q_{\delta}}+\alpha<1$ and $\iota$ to be a small positive number.

\textbf{Geometry difference term}

For the term $\chi_{cl,\Sigma}v\cdot (\frac{\mathcal{N}}{|\mathcal{N}|}-\frac{Bn_{0}}{|Bn_{0}|})$, we have
\begin{align}
\|\chi_{\mathrm{cl},\Sigma}v\cdot (\frac{\mathcal{N}}{|\mathcal{N}|}-\frac{\mathcal{B}n_{0}}{|\mathcal{B}n_{0}|})\|_{L^{2}W^{2-\frac{1}{q_{\delta}},q_{\delta}}(\Gamma_{+})}\lesssim\|\chi_{\mathrm{cl}}v\|_{L^{2}W^{2,q_{\delta}}}\|\eta\|_{L^{\infty}W^{3-\frac{1}{q_{\delta}},q_{\delta}}}\lesssim \mathcal{E}.
\end{align}
Here we replaced the spatial derivatives of $\mathfrak{h}_{\eta}$ in the local chart by the spatial derivatives of $\eta$ in the original chart, using the simple relation between them derived in the previous section.

Then it remains to estimate $\partial_{t}(\chi_{1,\Sigma}\eta)$ in the local chart. By the definition of $\mathfrak{h}_{\eta}$, we have  $\partial_{t}\eta = \partial_t\mathfrak{h}_{\eta} .$ Therefore, since $|\mathcal{N}_{0}|$ depends only on the steady state $\zeta$, which is a smooth function, we have
\begin{align}
\|\frac{1}{|\mathcal{N}_{0}|}\partial_{t}\eta_{cl}\|_{L^{2}W^{2-\frac{1}{q_{\delta}},q_{\delta}}}\lesssim \|\partial_{t}\eta_{cl}\|_{L^{2}W^{2-\frac{1}{q_{\delta}},q_{\delta}}}\lesssim\|\partial_{t}\eta_{cl}\|_{H^{\frac{3}{2}-\alpha}}\lesssim\sqrt{\mathscr{E}_{1,0,\mathrm{energy},\mathrm{cl}}}.
\end{align}

For the normal boundary condition on $\Gamma_{-}$, note that $\mathcal{F}_{\mathrm{cl},-}^{4}=-\chi_{\mathrm{cl},\Sigma}v\cdot (\frac{Q^{T}\mathcal{A}_{\Phi}}{|Q^{T}\mathcal{A}_{\Phi}|}-n_{0,-})$ which is also included in $\mathcal{F}_{j}^{4}$. It can be bounded in the same way as $\mathcal{F}_{\mathrm{cl}}^{4}$. We omit the details here.

\subsubsection{Difference of the tangential boundary conditions}

Define  $P_{\mathcal{B}} := I-\nu_{\mathcal{B}}\otimes\nu_{\mathcal{B}}.$ Suppose that the variable tangential condition is
\begin{equation}
(P_{\mathcal{B}} \bigl( \mathbb{D}_{\mathcal{B}}v+R \bigr) \nu_{B} + \beta v)\cdot \tau_{0} = \mathcal{K}.  \label{eq:variable-tangential-condition}
\end{equation}

On $\Gamma_{+}$, we have to estimate $\mathcal{F}_{j}^{3}$, which contains cutoff commutators, flattening map differences, $\mathcal{A}_{\eta}$ differences and nonlinear terms. We first estimate cutoff commutators. We have
\begin{align}{\label{est:tangential_1}}
\begin{aligned} \|(v\otimes \nabla \chi_{cl})\frac{\mathcal{N}}{|\mathcal{N}|}\|_{L^{2}W^{1-\frac{1}{q_{\delta}},q_{\delta}}}\lesssim& \frac{C(\varepsilon)}{\varepsilon}\|v\|_{L^{2}W^{1-\frac{1}{q_{\delta}},q_{\delta}}}(1+\|\nabla\eta\|_{L^{2}W^{1,q_{\delta}}})\\
\lesssim&\frac{C(\varepsilon)}{\varepsilon}(\|\chi_{cl}v\|_{L^{2}H^{1}}+\|\chi_{bet}v\|_{L^{2}H^{1}})(1+\|\chi_{cl,\Sigma}\eta\|_{L^{\infty}W^{3-\frac{1}{q_{\delta}},q_{\delta}}}+\|\chi_{bet,\Sigma}\eta\|_{L^{\infty}W^{3-\frac{1}{q_{\delta}},q_{\delta}}})\\
\lesssim&\frac{C(\varepsilon)}{\varepsilon}(\sqrt{\mathscr{E}_{0,0,\mathrm{energy},\mathrm{cl}}}+\sqrt{\mathscr{E}_{0,0,\mathrm{energy},\mathrm{bet}}}).
\end{aligned}
\end{align}
Similarly,
\begin{equation}
\|([L,\chi_{\mathrm{cl},\Sigma}]\eta)\frac{\mathcal{N}}{|\mathcal{N}|}\|_{L^2_\theta W^{1-1/q_\delta,q_\delta}_r} \le C(\varepsilon) \|\chi_{*,\Sigma}\eta\|_{H^{3/2-\alpha}(\Sigma)}\lesssim C(\varepsilon)\sqrt{\mathscr{E}_{0,0,\mathrm{energy},\mathrm{cl}}}.
\end{equation}

We then estimate flattening map differences. Write  $\mathcal{K }= \mathcal{K}^{(1)} + \mathcal{K}^{(2)} + \mathcal{K}^{(3)},$ where
\begin{align}
&\mathcal{K}^{(1)}  := (P_{0}-P_{\mathcal{B}})\mathbb{D}(\chi_{\mathrm{cl},\Sigma}v)\,n_{0},\quad \mathcal{K}^{(2)}  := P_{\mathcal{B}}\mathbb{D}(\chi_{cl,\Sigma}v)\,(n_{0}-\nu_{\mathcal{B}}),\notag \\
&\mathcal{K}^{(3)}  := -P_{\mathcal{B}} \bigl( \mathbb{D}_{\mathcal{Z}}(\chi_{\mathrm{cl},\Sigma}v)+R(\chi_{\mathrm{cl},\Sigma}v) \bigr) \nu_{\mathcal{B}}.
\end{align}

Since
\begin{align}
P_{\mathcal{B}}-P_{0} ={}& -(\nu_{\mathcal{B}}-n_{0})\otimes\nu_{\mathcal{B}} - n_{0}\otimes(\nu_{\mathcal{B}}-n_{0}),
\end{align}
we have  $|P_{\mathcal{B}}-P_{0}| \lesssim |\mathcal{Z}|.$ Therefore,
\begin{align}
\begin{aligned} \|\mathcal K^{(1)}\|_ {L^{2}W^{1-1/q_{\delta},q_{\delta}}(\Sigma)} \lesssim{}& \|\mathcal{Z}\|_ {L^{\infty}W^{1/q_{\delta}+\alpha,q_{\delta}}(\Sigma)} (\|\chi_{\mathrm{cl}}v\|_{L^{2}W^{2,q_{\delta}}(\mathcal{W})}+\|\chi_{\mathrm{cl}}\partial_{\theta}v\|_{L^{2}H^{1}(\mathcal{W})})\\
\lesssim{}&\varepsilon^{\iota}(\sqrt{\mathscr{E}_{0,1,\mathrm{energy},\mathrm{cl}}}+\sqrt{\mathscr{E}_{0,0,\mathrm{el},\mathrm{cl}}}),
\end{aligned}  \label{eq:K1-estimate}
\end{align}
\begin{align}
\begin{aligned} \|\mathcal K^{(2)}\|_ {L^{2}W^{1-1/q_{\delta},q_{\delta}}(\Sigma)} \lesssim{}& \|\mathcal{Z}\|_ {L^{\infty}W^{1/q_{\delta}+\alpha,q_{\delta}}(\Sigma)} (\|\chi_{\mathrm{cl}}v\|_{L^{2}W^{2,q_{\delta}}(\mathcal{W})}+\|\chi_{\mathrm{cl}}\partial_{\theta}v\|_{L^{2}H^{1}(\mathcal{W})})\\
\lesssim{}&\varepsilon^{\iota}(\sqrt{\mathscr{E}_{0,1,\mathrm{energy},\mathrm{cl}}}+\sqrt{\mathscr{E}_{0,0,\mathrm{el},\mathrm{cl}}}),
\end{aligned}  \label{eq:K2-estimate}
\end{align}
and
\begin{align}
&\|\mathcal K^{(3)}\|_ {W^{1-1/q_{\delta},q_{\delta}}(\Sigma)} \lesssim{}  \Big[ \|\mathcal{Z}\|_ {W^{1/q_{\delta}+\alpha,q_{\delta}}(\Sigma)} (\|\chi_{\mathrm{cl}}v\|_{W^{2,q_{\delta}}(\mathcal{W})}+\|\chi_{\mathrm{cl}}\partial_{\theta}v\|_{H^{1}(\mathcal{W})}) \Big]+\|\chi_{\mathrm{cl},\Sigma}v\|_{H^{1}}\notag \\
&\lesssim{} \varepsilon^{\iota}(\sqrt{\mathscr{E}_{0,1,\mathrm{energy},\mathrm{cl}}}+\sqrt{\mathscr{E}_{0,0,\mathrm{el},\mathrm{cl}}})+\sqrt{\mathscr{E}_{0,0,\mathrm{energy},\mathrm{cl}}}.\label{eq:K3-estimate}
\end{align}
Combining the three estimates, we obtain
\begin{align}
\|\mathcal{K}\|_ {W^{1-1/q_{\delta},q_{\delta}}(\Sigma)} \lesssim{}& \varepsilon^{\iota}(\sqrt{\mathscr{E}_{0,1,\mathrm{energy},\mathrm{cl}}}+\sqrt{\mathscr{E}_{0,0,\mathrm{el},\mathrm{cl}}})+\sqrt{\mathscr{E}_{0,0,\mathrm{energy},\mathrm{cl}}}. \label{eq:combined-tangential-boundary-estimate}
\end{align}

We next estimate the $\mathcal{A}_{\eta}$ difference terms. We have
\begin{align}
&\|\chi_{\mathrm{cl},\Sigma}(D_{\mathcal{B},\mathcal{A}_{\eta}}v\,\frac{\mathcal{N}}{|\mathcal{N}|}\cdot\tau_{0}-D_{\mathcal{B},I}v\nu_{\mathcal{B}}\cdot\tau_{0,+})\|_{L^{2}W^{1-\frac{1}{q_{\delta}},q_{\delta}}} \lesssim \|\chi_{\mathrm{cl},\Sigma}D_{\mathcal{B},\mathcal{A}_{\eta}-I}v\,\frac{\mathcal{N}}{|\mathcal{N}|}\|_{L^{2}W^{1-\frac{1}{q_{\delta}},q_{\delta}}}\notag \\
&+\|\chi_{\mathrm{cl},\Sigma}D_{\mathcal{B},I}v(\frac{\mathcal{N}}{|\mathcal{N}|}-\nu_\mathcal{B})\|_{L^{2}W^{1-\frac{1}{q_{\delta}},q_{\delta}}} \lesssim(\|\chi_{\mathrm{cl}}v\|_{L^{2}W^{2,q_{\delta}}}+\|\chi_{\mathrm{cl}}\partial_{\theta}v\|_{L^{2}H^{1}})\mathfrak{G}_{\eta}+\|(\chi_{\mathrm{cl}}+\chi_{\mathrm{bet}})v\|_{L^{2}H^{1}}\mathfrak{G}_{\eta}(1+\mathfrak{G}_{\eta})\notag \\
&\lesssim(\|\chi_{\mathrm{cl}}v\|_{L^{2}W^{2,q_{\delta}}}+\|\chi_{\mathrm{cl}}\partial_{\theta}v\|_{L^{2}H^{1}})\mathfrak{G}_{\eta}+(\|\chi_{\mathrm{cl}}v\|_{L^{2}H^{1}}+\|\chi_{bet}v\|_{L^{2}H^{1}})\mathfrak{G}_{\eta}(1+\mathfrak{G}_{\eta}) \lesssim \mathcal{E},\notag \\
&\|\chi_{\mathrm{cl},\Sigma}(g\eta-\operatorname{div}_{\Sigma}\mathcal{L}_{\Sigma}(\eta))(\frac{\mathcal{N}}{|\mathcal{N}|}-\nu_{B})\cdot \tau_{0,+}\|_{L^{2}W^{1-\frac{1}{q_{\delta}},q_{\delta}}}\lesssim \|\chi_{\mathrm{cl},\Sigma}(g\eta-\operatorname{div}_{\Sigma}\mathcal{L}_{\Sigma}(\eta))(\frac{\mathcal{N}}{|\mathcal{N}|}-\nu_{B})\|_{L^{2}W^{1-\frac{1}{q_{\delta}},q_{\delta}}}\notag \\
&\lesssim (\|\chi_{\mathrm{cl},\Sigma}\partial_{\theta}\eta\|_{H^{\frac{3}{2}-\alpha}}+\|\chi_{\mathrm{cl},\Sigma}\eta\|_{L^{2}W^{3-\frac{1}{q_{\delta}},q_{\delta}}})(\|\eta\|_{L^{\infty}W^{2,q_{\delta}}}+\|\partial_{\theta}\eta\|_{L^{\infty}H^{1}})\notag \\
&+(\|\eta\|_{L^{2}W^{2-\frac{1}{q_{\delta}},q_{\delta}}}+\|\partial_{\theta}\eta\|_{L^{2}H^{1}})(\|\eta\|_{L^{\infty}W^{2,q_{\delta}}}+\|\partial_{\theta}\eta\|_{L^{\infty}H^{1}})\label{est:tangential_2} \\
&\lesssim \Big(\|\chi_{\mathrm{cl},\Sigma}\eta\|_{L^{2}W^{3-\frac{1}{q_{\delta}},q_{\delta}}}+\|\chi_{\mathrm{cl},\Sigma}\eta\|_{H^{\frac{3}{2}-\alpha}}+\|\chi_{\mathrm{bet},\Sigma}\eta\|_{H^{\frac{3}{2}-\alpha}}+\|\chi_{\mathrm{cl},\Sigma}\partial_{\theta}\eta\|_{H^{\frac{3}{2}-\alpha}}+\|\chi_{\mathrm{cl},\Sigma}\partial_{\theta}\eta\|_{L^{\infty}H^{1}}\notag \\
&+\|\chi_{\mathrm{bet},\Sigma}\partial_{\theta}\eta\|_{L^{2}H^{1}}\Big)\cdot(\|\chi_{\mathrm{cl},\Sigma}\eta\|_{L^{\infty}W^{3-\frac{1}{q_{\delta}},q_{\delta}}}+\|\chi_{\mathrm{bet},\Sigma}\eta\|_{L^{\infty}W^{3-\frac{1}{q_{\delta}},q_{\delta}}}+\|\chi_{\mathrm{cl},\Sigma}\eta\|_{L^{\infty}H^{1}}+\|\chi_{\mathrm{bet},\Sigma}\eta\|_{L^{\infty}H^{1}}) \lesssim \mathcal{E}.\notag
\end{align}

Finally, for the remaining nonlinear term included in $\mathcal{F}_{cl}^{3}$, we have
\begin{align}
&\|\chi_{cl,\Sigma}\operatorname{div}\mathcal{R}(\nabla\eta)\frac{\mathcal{N}}{|\mathcal{N}|}\cdot \tau_{0,+}\|_{L^{2}W^{1-\frac{1}{q_{\delta}},q_{\delta}}}\lesssim\|\chi_{cl,\Sigma}\operatorname{div}\mathcal{R}(\nabla\eta)\frac{\mathcal{N}}{|\mathcal{N}|}\|_{L^{2}W^{1-\frac{1}{q_{\delta}},q_{\delta}}}\notag \\
&\lesssim (1+\mathfrak{G}_{\eta})(\|\chi_{\mathrm{cl},\Sigma}\eta\|_{L^{2}W^{3-\frac{1}{q_{\delta}},q_{\delta}}}+\|\chi_{\mathrm{cl},\Sigma}\partial_{\theta}\eta\|_{H^{\frac{3}{2}-\alpha}})(\|\partial_{\theta}\eta\|_{L^{\infty}H^{1}}+\|\eta\|_{L^{\infty}W^{2,q_{\delta}}})\notag \\
&\lesssim (1+\mathfrak{G}_{\eta})(\|\chi_{\mathrm{cl},\Sigma}\eta\|_{L^{2}W^{3-\frac{1}{q_{\delta}},q_{\delta}}}+\|\chi_{\mathrm{cl},\Sigma}\partial_{\theta}\eta\|_{H^{\frac{3}{2}-\alpha}})\label{est:tangential_3} \\
&\cdot (\|\chi_{\mathrm{cl},\Sigma}\partial_{\theta}\eta\|_{L^{\infty}H^{1}}+\|\chi_{\mathrm{cl},\Sigma}\eta\|_{L^{\infty}W^{2,q_{\delta}}}+\|\chi_{\mathrm{bet},\Sigma}\partial_{\theta}\eta\|_{L^{\infty}H^{1}}+\|\chi_{\mathrm{bet},\Sigma}\eta\|_{L^{\infty}W^{2,q_{\delta}}}) \lesssim \mathcal{E}.\notag
\end{align}

\textbf{The Navier condition on $\Gamma_{-}$.} The forcing $\tilde{\mathcal{F}}_{\mathrm{cl}}^{5}$ of the Navier-slip condition on the wall $\Gamma_{-}$ contains the same three kinds of terms as $\mathcal{F}^{3}_{\mathrm{cl}}$ on  $\Gamma_{+}$. The cutoff commutator $-\mu U_{0,-}(v\otimes\nabla\chi_{\mathrm{cl}}+\nabla\chi_{\mathrm{cl}}\otimes v)n_{0,-}$ is bounded as the commutator on $\Gamma_{+}$ above. The flattening differences are $\mathcal{K}^{(1)}$--$\mathcal{K}^{(3)}$ with $\nu_{\mathcal{B}}$ replaced by the normal of the wall in the chart, and they satisfy \eqref{eq:K1-estimate}--\eqref{eq:K3-estimate}. The $\mathcal{A}_{\eta}$-difference is the term with no counterpart on $\Gamma_{+}$. The matrix $\mathcal{A}_{\eta}-I$ has nonzero entries only in its third column, by \eqref{eq:A-explicit}, and $n_{0,-}=(\nu,0)$ on the vertical wall. Hence
\begin{equation*}
\big(\mathbb{D}_{\mathcal{A}_{\eta}-I}v\big)n_{0,-}=\big((\mathcal{A}_{\eta}-I)_{\cdot3}\cdot n_{0,-}\big)\,\partial_{3}v+(\mathcal{A}_{\eta}-I)_{\cdot3}\,(\partial_{3}v\cdot n_{0,-})=-K(\nabla_{x}q\cdot\nu)\,\partial_{3}v,
\end{equation*}
because $\partial_{3}v\cdot n_{0,-}=\partial_{3}(v\cdot n_{0,-})=0$ on the vertical wall. This is a product of $\nabla q$ with $\nabla v$. The pressure does not occur, since $U_{0,-}n_{0,-}=0$. Therefore, as for the $\mathcal{A}_{\eta}$-difference on $\Gamma_{+}$,
\begin{align}
&\|\chi_{\mathrm{cl},\Sigma}\,U_{0,-}\big(\mathbb{D}_{\mathcal{B},\mathcal{A}_{\eta}}v-\mathbb{D}_{\mathcal{B},I}v\big)n_{0,-}\|_{L^{2}W^{1-\frac{1}{q_{\delta}},q_{\delta}}(\Gamma_{-})} \lesssim\big(\|\chi_{\mathrm{cl}}v\|_{L^{2}W^{2,q_{\delta}}}+\|\chi_{\mathrm{cl}}\partial_{\theta}v\|_{L^{2}H^{1}}\big)\mathfrak{G}_{\eta}\notag \\
&+\big(\|\chi_{\mathrm{cl}}v\|_{L^{2}H^{1}}+\|\chi_{\mathrm{bet}}v\|_{L^{2}H^{1}}\big)\mathfrak{G}_{\eta}(1+\mathfrak{G}_{\eta})\lesssim\mathcal{E}.\label{eq:wall_Aeta_difference}
\end{align}
From these estimates, the full frozen Navier forcing in \eqref{eq:explicit-Fj6} satisfies
\begin{equation}
\|\mathcal{F}_{\mathrm{cl}}^{5}\|_{L^2_\theta W^{1-1/q_\delta,q_\delta}(\Gamma_-)} \le{} C(\varepsilon)\mathcal E +C\varepsilon^\iota\left( \sqrt{\mathscr E_{0,1,\mathrm{energy},\mathrm{cl}}} +\sqrt{\mathscr E_{0,0,\mathrm{el},\mathrm{cl}}}\right) +C(\varepsilon)\sqrt{\mathscr E_{0,0,\mathrm{energy},\mathrm{cl}}}.
\end{equation}
The last term retains the local linear cutoff and connection errors. The constant $C$ depends only on the fixed geometry and exponents, while $C(\varepsilon)$ also depends on the scaled cutoff profiles. The same bounds hold for the time and angular derivatives of these terms, with the norms of Subsection~\ref{subsec:energy_def} at the corresponding levels.

\subsubsection{Final elliptic estimate}

Combining the nonlinear estimates with the elliptic theory of \cite{GT}, we derive the following theorem.

\begin{theorem}
Suppose $(v_{r,h},p)$ is the solution to equation \eqref{eq:two-dimensional-slice-system}. Then $(v_{r,h},p)$ satisfies the following elliptic estimate:
\begin{align}
&\|v_{\mathrm{cl},r,h}\|_{L^2_\theta W^{2,q_\delta}(\mathcal W)}^2 +\|p_{\mathrm{cl}}-\overline p_{\mathrm{cl}}\|_{L^2_\theta W^{1,q_\delta}(\mathcal W)}^2 \lesssim \mathcal{E}^{2}+\frac{C(\varepsilon)}{\varepsilon} \mathscr{E}_{0,0,\mathrm{energy},\mathrm{cl}}+\mathscr{E}_{1,0,\mathrm{energy},\mathrm{cl}}+ \|\partial_{\theta}p_{\mathrm{cl}}\|_{L^{2}L^{q_{\delta}}}^{2}+\mathscr{E}_{0,1,\mathrm{energy},\mathrm{cl}}\notag \\
&+\|\partial_{\theta}v_{\mathrm{cl}}\|_{H^{1}}^{2}+C_{\mathcal{C}} \left( M_{\mathcal{Z}} + M_{\mathcal{Z}}^{2} + \varepsilon^{\iota}\right)(\mathscr{E}_{0,0,el,cl})+\|\partial_{t}v\|^{2}_{L^{2}}+\|\partial_{t}\eta_{\mathrm{cl}}\|^{2}_{H^{\frac{3}{2}-\alpha}},\label{eq:starting-elliptic-estimate}
\end{align}
where the functionals $\mathscr{E}$ are defined in equation \eqref{def:split_energy}.

The transverse support of $p_{\mathrm{cl}}$ has area at most $C\varepsilon^2$ by \eqref{eq:cutoff_eps}. The cutoffs commute with $\partial_\theta$. Thus H\"older on each slice gives
\begin{equation}
\|\partial_\theta p_{\mathrm{cl}}\|_{L^2_\theta L^{q_\delta}}^2 \le C\varepsilon^{4(1/q_\delta-1/2)}\|\partial_\theta p_{\mathrm{cl}}\|_{L^2}^2 \le C\varepsilon^{4(1/q_\delta-1/2)}\mathscr E_{0,1,\mathrm{energy},\mathrm{cl}} .
\end{equation}
For the velocity we keep the full angular norm. Put $\chi_*:=\chi_{\mathrm{cl}}+\chi_{\mathrm{bet}}$, which equals one on $\operatorname{supp}\chi_{\mathrm{cl}}$. The product rule gives
\begin{equation}
\|\partial_\theta v_{\mathrm{cl}}\|_{H^1}^2 \le C\mathscr E_{0,1,\mathrm{energy},\mathrm{cl}} +C(\varepsilon)\mathscr E_{0,0,\mathrm{energy},\mathrm{cl}}.
\end{equation}
Indeed, a transverse derivative of $\chi_{\mathrm{cl}}$ multiplies $\partial_\theta v$ in $L^2$ where $\chi_*=1$, which is controlled by $\|\chi_*v\|_{H^1}$. Derivatives of the fixed frame give only lower-order terms. In particular $\partial_\theta^2v_{\mathrm{cl}}$ is retained on the right in the angular energy. Consequently,
\begin{align}
&\|v_{\mathrm{cl},r,h}\|_{L^2_\theta W^{2,q_\delta}(\mathcal W)}^2 +\|p_{\mathrm{cl}}-\overline p_{\mathrm{cl}}\|_{L^2_\theta W^{1,q_\delta}(\mathcal W)}^2 \lesssim \mathcal{E}^{2}+\frac{C(\varepsilon)}{\varepsilon} \mathscr{E}_{0,0,\mathrm{energy},\mathrm{cl}}+\mathscr{E}_{1,0,\mathrm{energy},\mathrm{cl}}+\mathscr{E}_{0,1,\mathrm{energy},cl}\notag \\
&+C_{\mathcal C} \left( M_{\mathcal{Z}} + M_{\mathcal{Z}}^{2} +\varepsilon^{\iota}\right)(\mathscr{E}_{0,0,el,cl}).\label{eq:final-elliptic-estimate}
\end{align}
\end{theorem}

The last term of \eqref{eq:final-elliptic-estimate} is absorbed in Theorem~\ref{thm:cl_elliptic_base} below, once the estimates for $v_{\theta}$ and for $\eta$ are available. From \eqref{eq:final-elliptic-estimate}, it remains to derive the elliptic estimate for $v_{\theta}$ and $\eta$ or $\mathfrak{h}_{\eta}$.

\subsubsection{The transverse elliptic equation for $u_\protect\theta$}

In this section, we derive the elliptic system and the corresponding elliptic estimate for $v_{\theta}$. We have the following theorem.

\begin{theorem}
Suppose $v_{cl,\theta}$ is the $\theta$ component of the velocity field in the local chart near the contact line. For simplicity, we use $v_{\theta}$ to denote it in this theorem. Then it solves the following equation system:
\begin{equation}
\begin{cases}
-\mu\Delta_{r,h}v_\theta = \mathcal{F}_\theta^1, & \qquad\func{in}\quad \Omega \\
-\mu\partial_{n_0}v_\theta = \mathcal{F}_{\mathrm{cl}}^{3}\cdot e_\theta +\mu \partial_\theta v_{r,h}\cdot n_{0,+} & \qquad \func{on}\quad\Gamma_{+}, \\
-\mu\partial_{n_0}v_\theta = \mathcal{F}_{\mathrm{cl}}^{5}\cdot e_\theta +\mu \partial_\theta v_{r,h}\cdot n_{0,-} +\beta_{0} v_{\theta} & \qquad \func{on}\quad\Gamma_{-}.\end{cases}\end{equation}

Moreover, $v_{\theta}$ satisfies the following elliptic estimate
\begin{align*}
& \|v_\theta(\theta,\cdot)\|^{2}_{L^{2}W^{2,q_\delta}(\mathcal{W})} + \|v_\theta(\theta,\cdot)\|^{2}_{L^{2}H^{1+\delta}(\mathcal{W})} \lesssim\text{Right hand side of equation}\eqref{eq:final-elliptic-estimate}.
\end{align*}
\end{theorem}

\begin{proof}
(Elliptic equation with Neumann boundary condition) We now use the $e_\theta$-component of the momentum equation. For $a\in\{r,h\}$, we have the following estimate  $(Y_0)_{\theta a} = -\mu\partial_\theta v_a -\mu\partial_a v_\theta.$ Consequently,
\begin{align}
\func{div}_{r,h} \bigl(P_{r,h}Y_0e_\theta\bigr) ={}& -\mu\Delta_{r,h}v_\theta - \mu\partial_\theta\func{div}_{r,h}v_{r,h}. \label{eq:Yr-transverse-divergence}
\end{align}
The $e_\theta$-component of the momentum equation therefore gives
\begin{equation}  \label{eq:ur-transverse-equation} -\mu\Delta_{r,h}v_\theta = \mathcal{F}_\theta^1,
\end{equation}
where
\begin{equation}
\mathcal{F}_\theta^1 :={}  \mathcal{F}_{cl}^{1}\cdot e_\theta - \partial_\theta\bigl((Y_0e_\theta)\cdot e_\theta\bigr) - (\mathscr CY_0e_\theta)\cdot e_\theta   + \mu\partial_\theta\func{div}_{r,h}v_{r,h}-\partial_{t}v_{\theta}.\label{eq:Gr1-definition}
\end{equation}
Using  $\func{div}_{r,h}v_{r,h} = \mathcal{F}_{cl}^{2}-\partial_\theta v_\theta,$ we may equivalently write
\begin{equation}
\mathcal{F}_\theta^1 ={}  \mathcal{F}_{cl}^{1}\cdot e_\theta - \partial_\theta \bigl((Y_0(\chi_{{\mathrm{cl}}}v,\chi_{\mathrm{cl}}p)e_\theta)\cdot e_\theta\bigr) -\mu \partial_\theta^2\chi_{\mathrm{cl}}v_\theta - (\mathscr C(Y_0(\chi_{\mathrm{cl}}v,\chi_{\mathrm{cl}}p)e_\theta)\cdot e_\theta +\mu\partial_{\theta}\mathcal{F}_{\mathrm{cl}}^{2}-\partial_{t}v_{\theta}.\label{eq:Gr1-second-form}
\end{equation}

We now estimate $\|\mathcal{F}_{\theta}^{1}\|^{2}_{L^{q_{\delta}}}$. From the previous section for the $r,h$-elliptic estimate, this can be bounded by the right-hand side of \eqref{eq:final-elliptic-estimate}. For the term $\mathcal{F}_{\mathrm{cl}}^{1}\cdot e_{\theta}$, we bound it via the previous computation for $\mathcal{F}_{r,h}^{1}$ and $\mathcal{F}_{r,h}^{2}$. For other terms, we have
\begin{align*}
&\partial_\theta\bigl((Y_0e_\theta)\cdot e_\theta\bigr) =\partial_\theta p_{\mathrm{cl}}-2\mu\partial_\theta^2v_{\mathrm{cl},\theta}, \\
&\|\partial_\theta\bigl((Y_0e_\theta)\cdot e_\theta\bigr)\|_{L^2_\theta L^{q_\delta}} \lesssim\|\partial_\theta p_{\mathrm{cl}}\|_{L^2} +\|\partial_\theta v_{\mathrm{cl}}\|_{H^1},\quad \|(\mathscr{C}Y_{0}e_{\theta})\cdot e_{\theta}\|_{L^{2}L^{q_{\delta}}} \lesssim\|\chi_{\mathrm{cl}}v\|_{L^{2}H^{1}} +\|\chi_{\mathrm{cl}}p\|_{L^{2}L^{2}}.
\end{align*}
Here the $H^1$ norm of $\partial_\theta v_{\mathrm{cl}}$ includes all three spatial derivatives. Since the cutoffs are independent of $\theta$, the product rule gives
\begin{equation*}
\|\partial_\theta p_{\mathrm{cl}}\|_{L^2} +\|\partial_\theta v_{\mathrm{cl}}\|_{H^1} \le C\sqrt{\mathscr E_{0,1,\mathrm{energy},\mathrm{cl}}} +C(\varepsilon)\sqrt{\mathscr E_{0,0,\mathrm{energy},\mathrm{cl}}}.
\end{equation*}
Moreover, set $\chi_*:=\chi_{\mathrm{cl}}+\chi_{\mathrm{bet}}$. By \eqref{eq:Fj2-before-flattening} and \eqref{eq:local-F2-final},
\begin{equation*}
\partial_\theta\mathcal F_{\mathrm{cl}}^2 ={} \partial_\theta\bigl[(\chi_{\mathrm{cl}}F^2 +u\cdot\nabla\chi_{\mathrm{cl}})\circ\Phi\bigr] -\operatorname{div}_{\partial_\theta\mathcal Z}v_{\mathrm{cl}} -\operatorname{div}_{\mathcal Z}\partial_\theta v_{\mathrm{cl}} -\partial_\theta\bigl[(\partial_\theta Q\,v_{\mathrm{cl}}) \cdot(\mathcal A_\Phi e_\theta)\bigr].
\end{equation*}
The differentiated cutoff and frame terms are linear and satisfy
\begin{equation*}
\|\partial_\theta[(u\cdot\nabla\chi_{\mathrm{cl}})\circ\Phi]\|_{L^2_\theta L^{q_\delta}} +\|\partial_\theta[(\partial_\theta Q\,v_{\mathrm{cl}}) \cdot(\mathcal A_\Phi e_\theta)]\|_{L^2_\theta L^{q_\delta}} \le C(\varepsilon)\|\chi_*v\|_{H^1}.
\end{equation*}
Combining this with the nonlinear coefficient estimates and the definition of $M_{\mathcal Z}$ gives
\begin{equation*}
\|\partial_\theta\mathcal F_{\mathrm{cl}}^2\|_{L^2_\theta L^{q_\delta}} \le C(\varepsilon)\sqrt{\mathscr E_{0,0,\mathrm{energy},\mathrm{cl}}} +C M_{\mathcal Z}\sqrt{\mathscr E_{0,1,\mathrm{energy},\mathrm{cl}}} +C(\varepsilon)\mathcal E.
\end{equation*}
Here $C$ depends only on the fixed geometry and exponents.

For terms on the right-hand side, $\mathcal{F}_{j}^{3}$ and $\mathcal{F}_{j}^{5}$ have already been estimated in the previous subsection. Then, by the smallness of the region we chose near the contact line, we have
\begin{align}
&\|\partial_{\theta}\chi_{\mathrm{cl},\Sigma}v_{r,h}\|_{L^{2}W^{1-\frac{1}{q_{\delta}},q_{\delta}}}\lesssim  \|\partial_{\theta}\chi_{\mathrm{cl}}v\|_{L^{2}H^{1}}\lesssim \sqrt{\mathscr{E}_{0,1,\mathrm{energy},\mathrm{cl}}},\notag \\
&\|\beta_{0}\chi_{\mathrm{cl},\Sigma}v_{\theta}\|_{L^{2}W^{1-\frac{1}{q_{\delta}},q_{\delta}}}\lesssim  \|v\|_{L^{2}H^{1}}\lesssim \sqrt{\mathscr{E}_{0,0,\mathrm{energy},\mathrm{cl}}}.
\end{align}

Taking the scalar product of the Navier-slip boundary condition and the capillary boundary condition with $e_\theta$, and using $e_\theta\cdot n_0=0$, gives
\begin{equation}
-\mu\partial_{n_0}v_\theta ={}  \mathcal{F}_{\mathrm{cl}}^{3}\cdot e_\theta + \mu\partial_\theta v_{r,h}\cdot n_{0,+} \qquad \func{on}\quad\Gamma_{+}, -\mu\partial_{n_0}v_\theta ={}  \mathcal{F}_{\mathrm{cl}}^{5}\cdot e_\theta +\mu \partial_\theta v_{r,h}\cdot n_{0,-} +\beta v_{\theta}\qquad \func{on}\quad\Gamma_{-}.\label{eq:ur-boundary-r}
\end{equation}

We now apply an inhomogeneous Neumann estimate to the equation system
\begin{equation}
{\label{est:N_B}} \|v_{\theta}\|_{L^{2}W^{2,q_{\delta}}(W)} \lesssim\left(\|F_{\theta}^{1}\|_{L^{2}L^{q_{\delta}}(W)} +\|\partial_nv_{\theta}\|_{L^{2}W^{1-1/q_{\delta},q_{\delta}}(\partial W)} +\|v_{\theta}\|_{L^{2}H^1(W)}\right).
\end{equation}
Applying all the estimates above to the estimate \eqref{est:N_B}, we obtain
\begin{align}
\|v_\theta(\theta,\cdot)\|^{2}_{L^{2}W^{2,q_\delta}(\mathcal{W})} + \|v_\theta(\theta,\cdot)\|^{2}_{L^{2}H^{1+\delta}(\mathcal{W})} &\lesssim\text{Right hand side of equation}\eqref{eq:final-elliptic-estimate}+\|v_{\theta}(\theta,\cdot)\|^{2}_{L^{2}H^{1}} \\
\lesssim&\text{Right hand side of equation}\eqref{eq:final-elliptic-estimate},  \label{eq:ur-slice-estimate-r}
\end{align}
where we used the fact that $\|v_{\theta}\|_{L^{2}H^{1}}\lesssim\|v\|_{H^{1}}\lesssim \sqrt{\mathscr{E}_{0,0,\mathrm{energy},\mathrm{cl}}}$.
\end{proof}

\subsubsection{The combined slice estimate}

Combining \eqref{eq:transverse-final-slice-r} and \eqref{eq:ur-slice-estimate-r}, we obtain
\begin{align}
&\|\chi_{\mathrm{cl}}v(\theta,\cdot)\|^{2}_ {L^{2}W^{2,q_\delta}(\mathcal{W})\cap L^{2}H^{1+\delta}(\mathcal{W})} + \|\chi_{\mathrm{cl}}(p(\theta,\cdot))\|^{2}_ {L^{2}W^{1,q_\delta}(\mathcal{W})\cap L^{2}H^\delta(\mathcal{W})}\label{eq:combined-slice-estimate-r} \\
&\lesssim \mathcal{E}^{2}+\frac{C(\varepsilon)}{\varepsilon} \mathscr{E}_{0,0,\mathrm{energy},\mathrm{cl}}+\mathscr{E}_{1,0,\mathrm{energy},\mathrm{cl}}+\mathscr{E}_{0,1,\mathrm{energy},cl} +C_{\mathcal{C}} \left( M_{\mathcal{Z}} + M_{\mathcal{Z}}^{2} + \varepsilon^{\iota}\right)(\mathscr{E}_{0,0,el,cl}),\notag
\end{align}
To pass from $p_{\mathrm{cl}}-\overline p_{\mathrm{cl}}$ to $p_{\mathrm{cl}}$, apply Cauchy--Schwarz to the transverse average defining $\overline p_{\mathrm{cl}}(\theta)$ and integrate in $\theta$:
\begin{equation}
\|\overline p_{\mathrm{cl}}\|_{L^2_\theta}^2 \le |\mathcal W_{\mathrm{cl}}|^{-1} \|p_{\mathrm{cl}}\|_{L^2_\theta L^2(\mathcal W_{\mathrm{cl}})}^2 \lesssim \mathscr E_{0,0,\mathrm{energy},\mathrm{cl}}.
\end{equation}
Since the mean is constant in $(r,h)$, this also bounds its squared $L^2_\theta W^{1,q_\delta}$ and $L^2_\theta H^\delta$ norms on the fixed slice.

\subsection{Boundary elliptic estimate}

Before stating the theorem, we introduce the following notations:
\begin{align}
\mathcal{M}_{\zeta} := \frac{I} {\sqrt{1+|\nabla\zeta|^2}} - \frac{ \nabla\zeta\otimes\nabla\zeta }{ (1+|\nabla\zeta|^2)^{3/2} }.
\end{align}
Let $X$ be the horizontal part of the fixed map $\Phi$ in \eqref{eq:contact_chart}. On its upper face put $\widehat X(\theta,r)=X(\theta,r,-r\cot\theta_{\mathrm{eq}}),\, A_X=D_{\theta,r}\widehat X,\, J_X=\det A_X>0 .$ The orientation is that of the chart in the theorem below. Horizontal area is $dx=J_X\,d\theta\,dr$. Write
\begin{equation}
G_X=A_X^{-1}(\mathcal M_\zeta\circ\widehat X)A_X^{-T} =\begin{pmatrix}m_{\theta\theta}&m_{\theta r}\\ m_{r\theta}&m_{rr}\end{pmatrix}.
\end{equation}
Here $m_{\theta r}=m_{r\theta}$. Both $G_X$ and $J_X$ are fixed; $J_\eta$ denotes the determinant of the moving three-dimensional map.

\begin{theorem}
\label{thm:local-equation-h-eta} Suppose that the localized perturbation $\eta$ satisfies the following equation on $\Gamma_{+}$
\begin{equation}  \label{eq:cartesian-GC-h-eta} g\eta - \sigma\func{div}_x \bigl(\mathcal{M}_\zeta\nabla_x\eta\bigr) = \mathcal{F}_{\theta}^{3},
\end{equation}
where  $\mathcal{F}_{\theta}^{3}=-\mathcal{F}_{\mathrm{cl}}^{3}\cdot n_{0,+}+S(p,u)n_{0,+}\cdot n_{0,+}.$ Moreover, on the contact curve let $n_0$ denote the outward planar unit normal to $D$. The boundary condition for $\eta$ is
\begin{align}
-\kappa\partial_t\eta-\widetilde{\mathcal F}_{\mathrm{cl}}^6 =\sigma n_0\cdot\mathcal M_\zeta\nabla_x\eta .
\end{align}
Near the artificial endpoint of the surface chart the localized $\eta$ vanishes, so that its conormal flux there is zero.

Define $\mathfrak h_\eta(\theta,r,t)=\eta(\widehat X(\theta,r),t).$ Assume that $\Lambda(\theta)>0$ and that $\widehat X$ is an orientation-preserving local diffeomorphism from the coordinate domain $\mathcal D$. When writing the equation in $(\theta,r)$ we use the same symbols for the forcing terms composed with $\widehat X$. Let $\Gamma:=\partial\mathcal{D}\cap\{r=0\}$ be the contact-line portion of the coordinate boundary. Let $\Gamma_{a}:=\partial\mathcal{D}\cap \Gamma_{\mathrm{art}}$. Then $\mathfrak{h}_\eta$ satisfies
\begin{equation}  \label{eq:local-system-h-eta}
\begin{cases}
\displaystyle -\sigma\partial_r \left( m_{rr}\partial_r\mathfrak{h}_\eta \right) = \widetilde{\mathcal{M}}_\eta & \text{in }\mathcal{D}, \\[2mm]
\displaystyle \kappa \partial_t\mathfrak{h}_\eta - \sigma m_{r\theta}\partial_\theta\mathfrak{h}_\eta + \widetilde{\mathcal{F}}_{\mathrm{cl}}^6 = \sigma m_{rr}\partial_r\mathfrak{h}_\eta & \text{on }\Gamma, \\[2mm]
\displaystyle - \sigma m_{r\theta}\partial_\theta\mathfrak{h}_\eta= \sigma m_{rr}\partial_r\mathfrak{h}_\eta & \text{on }\Gamma_{a}.\end{cases}\end{equation}

\emph{Remark on \eqref{eq:local-system-h-eta}.} At $r=0$ the positive radial row $m_{r\theta}\partial_\theta\mathfrak h_\eta+ m_{rr}\partial_r\mathfrak h_\eta$ is the inward planar conormal flux. The outward flux is its negative; see the remark after \eqref{eq:s1-24}. Thus the contact law reads $\kappa\partial_t\mathfrak h_\eta+\widetilde{\mathcal F}_{\mathrm{cl}}^6 =\sigma(m_{r\theta}\partial_\theta\mathfrak h_\eta+ m_{rr}\partial_r\mathfrak h_\eta)$, with no zeroth-order term. On $\Gamma$, $m_{rr}=(1+|\partial_\tau\zeta|^2)/(1+|\nabla\zeta|^2)^{3/2}$, where $\partial_\tau\zeta=\nabla\zeta(\gamma(s))\cdot\dot\gamma(s)$ is differentiation along planar arclength at the corresponding boundary point. In the interior equation the right-hand side is
\begin{equation}
\widetilde{\mathcal{M}}_\eta :={}  \mathcal{F}_{\theta}^{3} - g\mathfrak{h}_\eta + \sigma(\partial_r\log J_X) m_{rr}\partial_r\mathfrak{h}_\eta + \frac{\sigma}{J_X} \partial_\theta \left[ J_X \left( m_{\theta\theta}\partial_\theta\mathfrak{h}_\eta + m_{\theta r}\partial_r\mathfrak{h}_\eta \right) \right]   + \frac{\sigma}{J_X} \partial_r \left( J_X m_{r\theta}\partial_\theta\mathfrak{h}_\eta \right).\label{eq:radial-source-h-eta}
\end{equation}
\end{theorem}

\begin{proof}
On the fixed upper face $H_0=-r\cot\theta_{\mathrm{eq}}$ is independent of $t$ and $\theta$. Differentiating $\widehat X$ therefore gives
\begin{equation}  \label{eq:A-X-explicit} A_X =
\begin{bmatrix}
\displaystyle \partial_\theta X & \displaystyle \partial_rX + \partial_hX \left( - \cot\theta_{\mathrm{eq}} \right)\end{bmatrix},
\end{equation}
where all derivatives of $X$ are evaluated at $(\theta,r,-r\cot\theta_{\mathrm{eq}})$.

\medskip \noindent \textbf{Spatial transformation.} The spatial chain rule yields  $\nabla_{\theta,r}\mathfrak{h}_\eta = A_X^T(\nabla_x\eta)(\widehat X,t),$ and hence
\begin{equation}  \label{eq:gradient-h-eta-theorem} (\nabla_x\eta)(\widehat X,t) = A_X^{-T}\nabla_{\theta,r}\mathfrak{h}_\eta.
\end{equation}
For a horizontal vector field $V$ the divergence formula is
\begin{equation}
(\func{div}_x V)\circ\widehat X =J_X^{-1}\func{div}_{\theta,r}\bigl(J_X A_X^{-1}(V\circ\widehat X)\bigr).
\end{equation}
Take $V=\mathcal M_\zeta\nabla_x\eta$ and use \eqref{eq:gradient-h-eta-theorem}. The resulting operator on a scalar $f$ is $\Delta_\Phi f=J_X^{-1}\func{div}_{\theta,r}(J_XG_X\nabla_{\theta,r}f).$ Thus the pullback of \eqref{eq:cartesian-GC-h-eta} is
\begin{equation}\label{eq:local-capillary-divergence} -\frac{\sigma}{J_X}\func{div}_{\theta,r} (J_XG_X\nabla_{\theta,r}\mathfrak h_\eta) =\mathcal F_\theta^3-g\mathfrak h_\eta .
\end{equation}

We now isolate the radial divergence term. Expanding \eqref{eq:local-capillary-divergence}, we obtain
\begin{equation}
-\frac{\sigma}{J_X} \partial_r \left( J_X m_{rr}\partial_r\mathfrak{h}_\eta \right) ={}  \mathcal{F}_{\theta}^{3}-g\mathfrak{h}_\eta + \frac{\sigma}{J_X} \partial_\theta \left[ J_X \left( m_{\theta\theta}\partial_\theta\mathfrak{h}_\eta + m_{\theta r}\partial_r\mathfrak{h}_\eta \right) \right]   + \frac{\sigma}{J_X} \partial_r \left( J_X m_{r\theta}\partial_\theta\mathfrak{h}_\eta \right).
\end{equation}
Since
\begin{equation*}
\frac{1}{J_X}  \partial_r  \left(  J_X m_{rr}\partial_r\mathfrak{h}_\eta  \right)  =  \partial_r  \left(  m_{rr}\partial_r\mathfrak{h}_\eta \right)  +  (\partial_r\log J_X)  m_{rr}\partial_r\mathfrak{h}_\eta,
\end{equation*}
moving the Jacobian contribution to the right-hand side yields  $-\sigma\partial_r ( m_{rr}\partial_r\mathfrak{h}_\eta ) = \widetilde{\mathcal{M}}_\eta,$ with $\widetilde{\mathcal{M}}_\eta$ defined by \eqref{eq:radial-source-h-eta}. This proves the interior equation in \eqref{eq:local-system-h-eta}.

\medskip \noindent \textbf{Boundary equation.} The horizontal map is $\gamma(\widehat\tau)+r\dot\gamma^\perp(\widehat\tau)$, so $r$ is the inward planar distance near $\partial D$. In particular $A_X^{-T}e_r=-n_0\circ\widehat X$ at $r=0$. By the chain rule,
\begin{equation}
(n_0\cdot\mathcal M_\zeta\nabla_x\eta)\circ\widehat X =-e_r\cdot G_X\nabla_{\theta,r}\mathfrak h_\eta =-m_{r\theta}\partial_\theta\mathfrak h_\eta-m_{rr}\partial_r\mathfrak h_\eta .
\end{equation}
Substitution into the physical contact condition gives
\begin{equation*}
\kappa\partial_t\mathfrak{h}_\eta  -  \sigma m_{r\theta}\partial_\theta\mathfrak{h}_\eta  + \widetilde{\mathcal{F}}_{\mathrm{cl}}^6  =  \sigma m_{rr}\partial_r\mathfrak{h}_\eta.
\end{equation*}
Similarly, on $\Gamma_{a}$, we have  $- \sigma m_{r\theta}\partial_\theta\mathfrak{h}_\eta= \sigma m_{rr}\partial_r\mathfrak{h}_\eta.$ This proves the boundary equation in \eqref{eq:local-system-h-eta}.

\medskip \noindent \textbf{The radial coefficient.} Set  $a:=\partial_\theta\widehat X, \, d:=\partial_r\widehat X, \, J_X=\det(a,d).$ For $a=(a_1,a_2)$, write $a^\perp=(-a_2,a_1)$. Then  $A_X^{-T}e_r=\frac{a^\perp}{J_X}.$ Writing  $W:=\sqrt{1+|\nabla\zeta(\widehat X)|^2},$ we obtain
\begin{equation}
m_{rr}  = (A_X^{-T}e_r)\cdot \mathcal{M}_\zeta(\widehat X)(A_X^{-T}e_r) = \frac{1}{J_X^2} a^\perp\cdot \mathcal{M}_\zeta(\widehat X)a^\perp = \frac{ W^2|a|^2 - |\nabla\zeta(\widehat X)\cdot a^\perp|^2 }{J_X^2W^3} = \frac{ |a|^2 + |\nabla\zeta(\widehat X)\cdot a|^2 }{J_X^2W^3}.
\end{equation}
In the last equality, we used  $|\nabla\zeta(\widehat X)|^2|a|^2 = |\nabla\zeta(\widehat X)\cdot a|^2 + |\nabla\zeta(\widehat X)\cdot a^\perp|^2.$

The map in \eqref{eq:contact_chart} uses the final angular parameter $\theta\in\mathbb T_M$. Thus the planar parameter of \eqref{eq:s1-16}, restricted to the surface, is $\widehat\tau(\theta,r)= \tau\bigl(\chi^{-1}(\theta),r,-r\cot\theta_{\mathrm{eq}}\bigr)$ and $\widehat X=\gamma(\widehat\tau)+r\dot\gamma^\perp(\widehat\tau)$. In particular, $\partial_\theta\widehat\tau =(\partial_1\tau)(\chi^{-1})'(\theta)$ and $\partial_r\widehat\tau=\partial_2\tau-\cot\theta_{\mathrm{eq}}\,\partial_3\tau$ with the derivatives of $\tau$ evaluated at $(\chi^{-1}(\theta),r,-r\cot\theta_{\mathrm{eq}})$. Setting $T=\dot\gamma(\widehat\tau)+r\ddot\gamma^\perp(\widehat\tau)$ gives $a=T\partial_\theta\widehat\tau$ and $d=\dot\gamma^\perp(\widehat\tau)+T\partial_r\widehat\tau$. Hence $J_X=\partial_\theta\widehat\tau\, \det\bigl(T,\dot\gamma^\perp(\widehat\tau)\bigr).$ Substitution into the preceding expression for $m_{rr}$ gives
\begin{equation*}
m_{rr}  =  \frac{ |T|^2 + |\nabla\zeta(\widehat X)\cdot T|^2 }{ W^3 \det\bigl(T,\dot\gamma^\perp(\widehat\tau)\bigr)^2 }.
\end{equation*}

Assume that $\gamma$ is parametrized by arc length, with  $\ddot\gamma(\tau) = \kappa(\tau)\dot\gamma^\perp(\tau).$ Then  $\ddot\gamma^\perp(\tau) = -\kappa(\tau)\dot\gamma(\tau), \, T = (1-r\kappa(\widehat\tau))\dot\gamma(\widehat\tau).$ Since  $|\dot\gamma|=1, \, \det(\dot\gamma,\dot\gamma^\perp)=1,$ the common factor $(1-r\kappa(\widehat\tau))^2$ cancels, yielding
\begin{equation}  \label{eq:mrr-h-eta-final} m_{rr} = \frac{ 1+ |\nabla\zeta(\widehat X)\cdot \dot\gamma(\widehat\tau)|^2 }{ (1+|\nabla\zeta(\widehat X)|^2)^{3/2} }.
\end{equation}
In particular, if $|\nabla\zeta|\leq K$ on the chart image, then  $m_{rr}\geq(1+K^2)^{-3/2}>0.$ Moreover,  $m_{rr}^{-1}\leq (1+|\nabla \zeta (\hat{X})|^{2})^{\frac{3}{2}},$ which is uniformly bounded in terms of the steady state function $\zeta$. So the equation system \eqref{eq:local-system-h-eta} is a uniformly elliptic system in divergence form.
\end{proof}

For the following scalar and contact equations, write  $-\sigma\partial_r(a\partial_rh)=M$, $(a\partial_rh)|_{r=0}=b$ with $a\ge a_0>0$. On a fixed interval assume uniform $C^2$ bounds for $a,a^{-1}$. Then
\begin{equation}
\|h\|_{L^2_\theta W_r^{3-1/q,q}} \le C\left(\|M\|_{L^2_\theta W_r^{1-1/q,q}} +\|b\|_{L^2_\theta}+\|h\|_{L^2_\theta L_r^q}\right).
\end{equation}
Given $\mathfrak{h}_{\eta}=\chi_{\mathrm{cl}}\eta(X(r,\theta,H_{0}(r,\theta),t))$, we derive the estimate for the forcing term on the right-hand side of equation \eqref{eq:local-system-h-eta}.

Write $L=\operatorname{div}_{\Sigma}\mathcal L_{\Sigma,0}$, and $n=n_{0,+}$. The exact normal forcing is
\begin{equation}
\mathcal F^3_{\mathrm{cl}}\cdot n ={} \widetilde F^3_{\mathrm{cl}}\cdot n +\mu(\mathbb D_Zv_{\mathrm{cl}}n)\cdot n -(R(v_{\mathrm{cl}})n)\cdot n +\bigl([S_B(p_{\mathrm{cl}},v_{\mathrm{cl}})+R(v_{\mathrm{cl}}) -(g\eta_{\mathrm{cl}}-\sigma L\eta_{\mathrm{cl}} )I](n-\nu_{B,+})\bigr)\cdot n.
\end{equation}
In the previous section, by estimates \eqref{est:tangential_1}, \eqref{est:tangential_2} and \eqref{est:tangential_3}, we have already bounded the tangential direction of $\mathcal{F}_{\mathrm{cl}}^{3}$ via
\begin{align}
\|\mathcal{F}_{\mathrm{cl},t}^{3}\cdot \tau_{0,\pm}\|_{L^{2}W^{1-\frac{1}{q_{\delta}},q_{\delta}}}\lesssim \|\mathcal{F}_{\mathrm{cl},t}^{3}\|_{L^{2}W^{1-\frac{1}{q_{\delta}},q_{\delta}}}\|\tau_{0,\pm}\|_{L^{\infty}W^{2,q_{\delta}}}\lesssim \|\mathcal{F}_{\mathrm{cl},t}^{3}\|_{L^{2}W^{1-\frac{1}{q_{\delta}},q_{\delta}}},
\end{align}
where $\mathcal{F}_{\mathrm{cl},t}^{3}$ denotes terms included in $\mathcal{F}_{\mathrm{cl}}^{3}$ that have a nonzero tangential component. Therefore, the following estimate holds for the normal direction of $\mathcal{F}_{\mathrm{cl},t}^{3}$:
\begin{align}
\|\mathcal{F}_{\mathrm{cl},t}^{3}\cdot n_{0}\|_{L^{2}W^{1-\frac{1}{q_{\delta}},q_{\delta}}}\lesssim \|\mathcal{F}_{\mathrm{cl},t}^{3}\|_{L^{2}W^{1-\frac{1}{q_{\delta}},q_{\delta}}}.
\end{align}
For the remaining terms, in particular, $\widetilde F^3_{\mathrm{cl}}\cdot n$ includes the capillary cutoff commutator. Its estimate is linear and can be derived as follows:
\begin{equation}
\|[L,\chi_{\mathrm{cl},\Sigma}]\eta\frac{N}{|\mathcal{N}|}\cdot n\|_{L^2_\theta W^{1-1/q_\delta,q_\delta}_r} \le C(\varepsilon) \|\chi_{*,\Sigma}\eta\|_{H^{3/2-\alpha}(\Sigma)}\lesssim C(\varepsilon)\sqrt{\mathscr{E}_{0,0,\mathrm{energy},\mathrm{cl}}}.
\end{equation}
The operator $L$ contains radial, mixed and angular second derivatives with smooth coefficients. The radial second derivatives are controlled by $\|\eta_{\mathrm{cl}}\|_{L^2_\theta W_r^{3-1/q_\delta,q_\delta}}$. For the other second derivatives, the embedding $H_r^{1/2-\alpha}\subset W_r^{1-1/q_\delta,q_\delta}$ gives
\begin{align}
&\|\partial_\theta^2\eta_{\mathrm{cl}}\|_{L^2_\theta W_r^{1-1/q_\delta,q_\delta}} +\|\partial_r\partial_\theta\eta_{\mathrm{cl}}\|_{L^2_\theta W_r^{1-1/q_\delta,q_\delta}} \le C\|\partial_\theta\eta_{\mathrm{cl}}\|_{H^{3/2-\alpha}(\Sigma)}.
\end{align}
To compare the localized norm with the energy, put $s=3/2-\alpha$ and $\chi_{*,\Sigma}=\chi_{\mathrm{cl},\Sigma}+\chi_{\mathrm{bet},\Sigma}$. Since the cutoffs commute with $\partial_\theta$ and $[\langle D\rangle^s,\chi_{\mathrm{cl},\Sigma}]$ has order $s-1$,
\begin{equation}
\|\partial_\theta\eta_{\mathrm{cl}}\|_{H^s} \le C\|\chi_{*,\Sigma}\partial_\theta\eta\|_{H^s} +C(\varepsilon)\|\chi_{*,\Sigma}\eta\|_{H^s} \le C\mathscr E_{0,1,\mathrm{energy},\mathrm{cl}}^{1/2} +C(\varepsilon)\mathscr E_{0,0,\mathrm{energy},\mathrm{cl}}^{1/2}.
\end{equation}
Applying the smooth small geometric multipliers to these derivative bounds, and estimating the stress differences as above, we find that the full forcing satisfies
\begin{equation}
\|\mathcal F^3_{\mathrm{cl}}\cdot n\|_{L^2_\theta W^{1-1/q_\delta,q_\delta}_r} \le{} C(\varepsilon)\mathscr E_{0,0,\mathrm{energy},\mathrm{cl}}^{1/2} +C\varepsilon^{\iota}\mathscr E_{0,0,\mathrm{el},\mathrm{cl}}^{1/2} +C\mathscr E_{0,1,\mathrm{energy},\mathrm{cl}}^{1/2} +C(\varepsilon)\mathcal E,
\end{equation}
with a possibly smaller positive $\iota$. The lower and angular energy terms are retained in the base elliptic estimate; the elliptic term is absorbed after combining the estimates. For the remaining normal stress in $\mathcal F_\theta^3$, note that $n_{0,+}$ is constant and $n_{0,+}\cdot e_\theta=0$. Hence
\begin{equation}
n_{0,+}\cdot S(p_{\mathrm{cl}},v_{\mathrm{cl}})n_{0,+} =p_{\mathrm{cl}}-2\mu\partial_{n_{0,+}}(v_{\mathrm{cl}}\cdot n_{0,+}),
\end{equation}
where $\partial_{n_{0,+}}$ uses only $(r,h)$ derivatives. Taking the trace on each transverse slice and then its $L^2_\theta$ norm gives
\begin{align}
\|n_{0,+}\cdot S(p_{\mathrm{cl}},v_{\mathrm{cl}})n_{0,+} \|_{L^2_\theta W^{1-1/q_\delta,q_\delta}(\Gamma_+)} &\lesssim \|p_{\mathrm{cl}}\|_{L^2_\theta W^{1,q_\delta}(\mathcal W)} +\|v_{\mathrm{cl}}\|_{L^2_\theta W^{2,q_\delta}(\mathcal W)} .
\end{align}

For the remaining capillary terms in \eqref{eq:radial-source-h-eta}, we write $f=(\chi_{\mathrm{cl},\Sigma}\eta)\circ\widehat X$ for the localized height in that equation for simplicity. In particular, the cutoff is included. Expanding the angular flux gives
\begin{align}
&J_X^{-1}\partial_\theta\!\left[J_X (m_{\theta\theta}\partial_\theta f+m_{\theta r}\partial_r f)\right] =m_{\theta\theta}\partial_\theta^2 f +m_{\theta r}\partial_r\partial_\theta f \\
&+(\partial_\theta m_{\theta\theta} +m_{\theta\theta}\partial_\theta\log J_X)\partial_\theta f +(\partial_\theta m_{\theta r} +m_{\theta r}\partial_\theta\log J_X)\partial_r f .\notag
\end{align}
The coefficients and their radial derivatives are bounded on the fixed chart. They are therefore multipliers on $W_r^{1-1/q_\delta,q_\delta}$, uniformly in $\theta$. For the two radial flux terms we use
\begin{equation}
J_X^{-1}\partial_r(J_Xm_{r\theta}\partial_\theta f) =m_{r\theta}\partial_r\partial_\theta f +(\partial_r m_{r\theta}+m_{r\theta}\partial_r\log J_X)\partial_\theta f .
\end{equation}
Thus the second derivatives in these source terms are first derivatives of $\partial_\theta f$. Since $H_r^{1/2-\alpha}$ embeds into $W_r^{1-1/q_\delta,q_\delta}$ on the fixed radial interval, we obtain
\begin{align}
&\left\|\frac{\sigma}{J_X}\partial_\theta [J_X(m_{\theta\theta}\partial_\theta f+m_{\theta r}\partial_r f)] \right\|_{L^2_\theta W_r^{1-1/q_\delta,q_\delta}} +\|\sigma(\partial_r\log J_X)m_{rr}\partial_r f \|_{L^2_\theta W_r^{1-1/q_\delta,q_\delta}} \\
&+\left\|\frac{\sigma}{J_X}\partial_r (J_Xm_{r\theta}\partial_\theta f) \right\|_{L^2_\theta W_r^{1-1/q_\delta,q_\delta}} \le C\bigl(\|f\|_{H^{3/2-\alpha}}+ \|\partial_\theta f\|_{H^{3/2-\alpha}}\bigr)\le C(\varepsilon)\mathscr E_{0,0,\mathrm{energy},\mathrm{cl}}^{1/2} +C\mathscr E_{0,1,\mathrm{energy},\mathrm{cl}}^{1/2}.\notag
\end{align}
Here the $H^{3/2-\alpha}$ norms are on the full surface chart. The last line uses the preceding cutoff comparison for $\partial_\theta f$ and the bounded cutoff multiplier for $f$. In particular these three fixed-coefficient source terms require no radial second derivative of $f$.

For the contact line condition, we first estimate the term induced by decomposing $\partial_{t}\eta$
\begin{align}
&\|\partial_t(\chi_{\mathrm{cl},\Sigma}\mathfrak h_\eta) - \sigma m_{r\theta}\partial_\theta(\chi_{\mathrm{cl},\Sigma}\mathfrak h_\eta) )\|_{L^{2}L^{\infty}(\p D)} \\
&\lesssim \|\chi_{\mathrm{cl},\Sigma}\partial_{t}\mathfrak{h}_{\eta}\|_{L^{2}L^{\infty}(\p D)}+\|\chi_{\mathrm{cl},\Sigma}\partial_{\theta}\mathfrak{h}_{\eta}\|_{L^{2}L^{\infty}(\p D)}+\|\chi_{\mathrm{cl},\Sigma}\mathfrak{h}_{\eta}\|_{L^{2}L^{\infty}(\p D)}\notag \\
&\lesssim \|\chi_{\mathrm{cl},\Sigma}\partial_{t}\mathfrak{h}_{\eta}\|_{L^{2}L^{\infty}(\p D)}+\|\chi_{\mathrm{cl},\Sigma}\partial_{\theta}\mathfrak{h}_{\eta}\|_{H^{\frac{3}{2}-\alpha}}+\|\chi_{\mathrm{cl},\Sigma}\mathfrak{h}_{\eta}\|_{H^{\frac{3}{2}-\alpha}} \lesssim \sqrt{\mathscr{E}_{0,1,\mathrm{energy},\mathrm{cl}}}+\sqrt{\mathscr{E}_{1,0,\mathrm{energy},\mathrm{cl}}}.\notag
\end{align}
Then for the nonlinear term included in $\mathcal{F}_{\mathrm{cl}}^{7}$, we have
\begin{align}
\begin{aligned}
&\|\chi_{\mathrm{cl},\Sigma}\hat{\mathcal{W}}(\partial_{t}\mathfrak{h}_{\eta})\|_{L^{2}L^{\infty}(\p \mathcal{D})}\lesssim \|\partial_{t}\mathfrak{h}_\eta{}\|_{L^{\infty}L^{\infty}(\p D)}\|\partial_{t}\mathfrak{h}_\eta{}\|_{L^{2}L^{\infty}(\p D)} \lesssim \mathcal{E}. \end{aligned}
\end{align}

Therefore, once the energy estimate and the $\frac{3}{2}-\alpha$ estimate (at least two more orders) have been derived, the energy norm of $(u,p,\eta)$ can also be bounded by those norms.

\begin{theorem}
\label{thm:cl_elliptic_base}  Suppose $(v_{\mathrm{cl}},p_{\mathrm{cl}},\eta_{\mathrm{cl}})$ satisfy the system \eqref{eq:final-localized-flattened-system-j}. If $\varepsilon$ is small enough, then they satisfy the following elliptic estimate
\begin{align}
&\mathscr{E}_{0,0,el,cl} \lesssim \mathcal{E}^{2}+\frac{C(\varepsilon)}{\varepsilon} \mathscr{E}_{0,0,\mathrm{energy},\mathrm{cl}}+\mathscr{E}_{1,0,\mathrm{energy},\mathrm{cl}}+\mathscr{E}_{0,1,\mathrm{energy},cl}.
\end{align}
\end{theorem}

\begin{proof}
The combined slice estimate \eqref{eq:combined-slice-estimate-r} bounds the velocity and pressure parts of $\mathscr{E}_{0,0,el,cl}$; the pressure mean $\overline{p}_{cl}$ is bounded by the energy pressure terms, as noted after it. The estimate for $\mathfrak{h}_{\eta}$ of this subsection bounds the surface part, in terms of the normal stress on $\Gamma_{+}$ in $L^{2}W^{1-1/q_{\delta},q_{\delta}}$, that is, in terms of the velocity and pressure parts, and of energy terms. Inserting \eqref{eq:combined-slice-estimate-r} into it and adding, we find that every occurrence of  $\mathscr{E}_{0,0,el,cl}$ on the right-hand side carries the factor $C\big(M_{\mathcal{Z}}+M_{\mathcal{Z}}^{2}+\varepsilon^{\iota}\big)$, with $C$ independent of $\varepsilon$. By \eqref{eq:MZ_small} this factor is $\leq C(\varepsilon+\varepsilon^{\iota})$, which is $\leq\frac12$ for $\varepsilon$ small. The term is then absorbed into the left-hand side.
\end{proof}

\subsection{Higher tangential derivatives}

We apply $\partial_\theta^{k_{2}}\partial_{t}^{k_{1}}$ \, to the extracted system. The resulting equations have the same two-dimensional transverse principal operator, while derivatives falling on $\mathscr C$ produce lower-order commutators. In the following estimate, we separate the forcing terms on the right-hand side into cutoff commutators, flattening differences, $\mathcal{A}_{\eta}$ differences and nonlinear terms of the original equation system, and apply the elliptic estimates to each of them. In this subsection, we aim to prove the following theorem:

\begin{proof}[Proof of Theorem~\protect\ref{thm:Higher-order contact line elliptic estimate}]
Again, it suffices to prove the estimate for the dissipation.

Tangential derivatives of the cutoff and flattening forcing terms

In the following proof, we take $\partial_{||}^{\alpha}=\partial_{\theta}^{m}$. All of the estimates are still true if we replace $\partial_{\theta}$ by the temporal derivative $\partial_{t}$. We now consider the forcing terms generated by the localization and the flattening map after applying $\partial_\theta^j$. These terms are easier to estimate since they either contain derivatives of the fixed cutoff functions or derivatives of the geometric difference  $\mathcal{Z}:=\mathcal{B}-I.$ In either case, the corresponding derivatives of the unknowns are of strictly lower tangential order.

We first consider the flattening-map differences. A typical highest-order term has the form  $\mathcal{Z }\nabla_{r,h}^{2}v \,\text{or}\, \mathcal{Z}\nabla_{r,h}p.$ Applying $\partial_\theta^j$ gives
\begin{align}
\partial_\theta^j \left( \mathcal{Z}\nabla_{r,h}^{2}(\chi_{\mathrm{cl},\Sigma}v) \right) ={}& \mathcal{Z }\nabla_{r,h}^{2}\partial_\theta^j(\chi_{\mathrm{cl},\Sigma}v) + \sum_{\ell=1}^{j} \binom{j}{\ell} \partial_\theta^\ell\mathcal{Z}\, \nabla_{r,h}^{2} \partial_\theta^{j-\ell}(\chi_{\mathrm{cl},\Sigma}v), \\
\partial_\theta^j \left( \mathcal{Z}\nabla_{r,h}(\chi_{\mathrm{cl},\Sigma}p) \right) ={}& \mathcal{Z }\nabla_{r,h}\partial_\theta^j(\chi_{\mathrm{cl},\Sigma}p) + \sum_{\ell=1}^{j} \binom{j}{\ell} \partial_\theta^\ell\mathcal{Z}\, \nabla_{r,h} \partial_\theta^{j-\ell}(\chi_{\mathrm{cl},\Sigma}p).
\end{align}
The first terms on the right-hand sides contain the top-order derivatives of  $(v,p)$, but are multiplied by the small geometric coefficient $\mathcal{Z}$ and hence may be absorbed into the left-hand side provided that
\begin{equation*}
\sum_{\ell=0}^{m}\|\partial_{\theta}^{\ell}\mathcal{Z}\|_{L^{\infty}L^\infty} \leq \varepsilon
\end{equation*}
for sufficiently small $\varepsilon>0$. Every remaining term contains at least one tangential derivative falling on $\mathcal{Z}$. Consequently, the highest tangential derivative falling on $(v,p)$ is at most $j-1$. Thus these terms are lower-order terms and may be controlled by the lower-order elliptic estimates. Moreover, if at least one temporal derivative falls on $\mathcal{Z}$, the term vanishes, since $\zeta$ is independent of time.

More generally, for a product of the form  $\mathcal{G}(\mathcal{Z})\, \partial_\theta^a\nabla_{r,h}^b v,$ where $\mathcal{G}(0)=0$, the Leibniz rule yields
\begin{equation*}
\partial_\theta^j \left( \mathcal{G}(\mathcal{Z}) \partial_\theta^a\nabla_{r,h}^b v \right) = \mathcal{G}(\mathcal{Z}) \partial_\theta^{j+a}\nabla_{r,h}^bv + \mathcal{R}_{\mathcal{Z},j},
\end{equation*}
where every term in $\mathcal{R}_{\mathcal{Z},j}$ contains at least one $\theta$-derivative of $\mathcal{G}(\mathcal{Z})$ and therefore at most $j+a-1 $ tangential derivatives of $v$. The same observation applies to the geometric differences  $\nu_{\mathcal{B}}-n_0, \, U_{\mathcal{B}}-U_0, \, R-R_0,$ since they are smooth functions of $\mathcal{Z}$ and its geometric coefficients and vanish at the equilibrium configuration.

Secondly, we consider
\begin{align}
\func{div}_{r,h}\nabla_{I-Q^{T}\mathcal{A}_{\Phi}}v_{\mathrm{cl}}=\mathcal{Z}\Delta_{r,h}v_{\mathrm{cl}}+\nabla_{r,h}(I-Q^{T}\mathcal{A}_{\Phi})\cdot \nabla v_{\mathrm{cl}}.
\end{align}
For the first term, this can be bounded by the smallness of $\|\mathcal{Z}\|_{L^{\infty}}$ as we discussed above. For the second term, its $L^{q_{\delta}}$ norm can be bounded by the energy estimate. There are some similar terms on the boundary; when $\nabla_{r,h}$ falls on $I-Q^{T}\mathcal{A}_{\Phi}$, the coefficient is no longer small. However, since one fewer spatial derivative falls on $(v,p,\eta)$, we can bound these terms by the energy estimate. In conclusion, denoting by $\mathcal{F}_{\mathrm{cl},\mathrm{fd}}^{j}$ all of the forcing terms induced by the flattening difference included in $\mathcal{F}_{\mathrm{cl}}^{j}$, we have
\begin{align}
&\sum_{j=0}^{m}\big(\|\partial_{\theta}^{j}\mathcal{F}_{\mathrm{cl},\mathrm{fd}}^{1}\|^{2}_{L^{2}L^{q_{\delta}}}+\|\partial_{\theta}^{j}\mathcal{F}_{\mathrm{cl},\mathrm{fd}}^{2}\|_{L^{2}W^{1,q_{\delta}}}^{2}+\|\partial_{\theta}^{j}\mathcal{F}_{\mathrm{cl},\mathrm{fd}}^{3}\|^{2}_{L^{2}W^{1-\frac{1}{q_{\delta}},q_{\delta}}}\big)+\sum_{j=0}^{m}\|\partial_{\theta}^{j}\mathcal{F}_{\mathrm{cl},-,\mathrm{gd}}^{4}\|^{2}_{L^{2}W^{2-\frac{1}{q_{\delta}},q_{\delta}}} \notag \\
&\quad+\sum_{j=0}^{m}\big(\|\partial_{\theta}^{j}\mathcal{F}_{\mathrm{cl},\mathrm{gd}}^{4}\|^{2}_{L^{2}W^{2-\frac{1}{q_{\delta}},q_{\delta}}}+\|\partial_{\theta}^{j}\mathcal{F}_{\mathrm{cl},\mathrm{gd}}^{5}\|^{2}_{L^{2}W^{1-\frac{1}{q_{\delta}},q_{\delta}}}+\|\partial_{\theta}^{j}\mathcal{F}_{\mathrm{cl},\mathrm{fd}}^{6}\|^{2}_{L^{2}(\Gamma_{0})}\big)  \notag \\
&\lesssim \varepsilon^{\iota}(\sum_{j=0}^{k_{2}+1}\mathscr{E}_{k_{1},j,\mathrm{energy},\mathrm{cl}}+\sum_{j=0}^{k_{2}}\mathscr{E}_{k_{1},j,\mathrm{el},\mathrm{cl}})+\sum_{j=0}^{k_{2}}\mathscr{E}_{k_{1},j,\mathrm{energy},\mathrm{cl}}.
\end{align}

We next consider the cutoff commutators. Recall that the cutoff functions $\chi_j$ and $\chi_{j,\Sigma}$ are chosen to be independent of $\theta$. Hence  $\partial_\theta^k\chi_j = \partial_\theta^k\chi_{j,\Sigma} = 0, \, k\geq1,$ and therefore  $\partial_\theta^j ( \chi_j f ) = \chi_j\partial_\theta^jf.$ Thus taking tangential derivatives does not generate any additional commutators with $\chi_j$. In particular, the cutoff forcing terms already generated before flattening, such as  $S(p,u)\nabla\chi_j, \, \func{div}\mathcal{C}_{\chi_j}(u), \, u\cdot\nabla_{\mathcal{A}_\eta}\chi_j\,u$ retain the same localization structure after applying $\partial_\theta^j$. Since every such term contains at least one derivative of $\chi_j$, it is supported in the transition region  $\func{supp}\nabla\chi_j \subset \{\chi_{cl}+\chi_{bet}=1\},$ and hence it is controlled by the lower-order estimates associated with the larger cutoff $\chi_{\mathrm{cl},\Sigma}+\chi_{\mathrm{bet},\Sigma}$. Consequently, the $\chi$-commutator terms do not contribute to top-order spatial derivative terms and are therefore bounded by the energy estimate. In conclusion,
\begin{align}
&\sum_{j=0}^{m}\|\partial_{\theta}^{j}\mathcal{F}_{\mathrm{cl},\mathrm{cutoff}}^{1}\|_{L^{2}L^{q_{\delta}}}+\|\partial_{\theta}^{j}\mathcal{F}_{\mathrm{cl},\mathrm{fcutoff}}^{2}\|_{L^{2}W^{1,q_{\delta}}}+\|\partial_{\theta}^{j}\mathcal{F}_{\mathrm{cl},\mathrm{cutoff}}^{3}\|_{L^{2}W^{1-\frac{1}{q_{\delta}},q_{\delta}}}  \notag \\
&\quad+\sum_{j=0}^{m}\|\partial_{\theta}^{j}\mathcal{F}_{\mathrm{cl},\mathrm{cutoff}}^{4}\|_{L^{2}W^{2-\frac{1}{q_{\delta}},q_{\delta}}}+\|\partial_{\theta}^{j}\mathcal{F}_{\mathrm{cl},\mathrm{cutoff}}^{5}\|_{L^{2}W^{1-\frac{1}{q_{\delta}},q_{\delta}}}+\|\partial_{\theta}^{j}\mathcal{F}_{\mathrm{cl},-,\mathrm{cutoff}}^{4}\|_{L^{2}W^{2-\frac{1}{q_{\delta}},q_{\delta}}}  \notag \\
&\lesssim \sqrt{\mathcal{E}\mathcal{D}}+\frac{C(\varepsilon)}{\varepsilon} \sum_{0\leq j_{1}\leq k_{1}}\sum_{0\leq j_{2}\leq k_{2}}\sqrt{\mathscr{E}_{j_{1},j_{2},\mathrm{energy},\mathrm{cl}}}+\frac{C(\varepsilon)}{\varepsilon} \sum_{0\leq j_{1}\leq k_{1}}\sum_{0\leq j_{2}\leq k_{2}}\sqrt{\mathscr{E}_{j_{1},j_{2},\mathrm{energy},\mathrm{bet}}}.
\end{align}

Higher-order estimates for the remaining nonlinear terms

We now estimate the nonlinear terms inherited from the original free-boundary problem. In contrast with the flattening-map differences, these terms are at least quadratic in the unknowns. In the remaining part of the proof, we will only discuss the case $j=m$. For $j\leq m-1$, we follow the same idea with fewer derivatives.

For convenience, define  $\chi_*:=\chi_{\mathrm{cl}}+\chi_{\mathrm{bet}}, \, \chi_{*,\Sigma} := \chi_{\mathrm{cl},\Sigma}+\chi_{\mathrm{bet},\Sigma}.$ By the construction of the cutoff functions,  $\chi_*=1 \,\text{on }\func{supp}\chi_{\mathrm{cl}}, \, \chi_{*,\Sigma}=1 \,\text{on }\func{supp}\chi_{\mathrm{cl},\Sigma}.$ Therefore, every factor appearing in a nonlinear term localized by $\chi_1$ may be replaced by its $\chi_*$-localized counterpart. In this way, all quantities appearing below are controlled either by the contact-line estimate or by the intermediate estimate.

For the remaining nonlinear terms, we only consider the derivative $\partial_{\theta}^{m}$. The general case follows from a similar estimate, replacing some $\partial_{\theta}$ by $\partial_{t}$. Moreover, since the steady state is smooth, we ignore the geometric terms induced by the flattening map, since they are uniformly bounded.

\textbf{Key observation}

Using the Sobolev embedding in 1-D, we have $H^{1}\hookrightarrow L_{\theta}^{\infty}$. Therefore,
\begin{align}
&\sum_{j_{1}+j_{2}\leq m-1}(\|\chi_{cl}\partial_{t}^{j_{1}}\partial_{\theta}^{j_{2}}v\|_{L^{\infty}W^{2,q_{\delta}}}+\|\chi_{cl,\Sigma}\partial_{t}^{j_{1}}\partial_{\theta}^{j_{2}}\eta\|_{L^{\infty}W^{3-\frac{1}{q_{\delta}},q_{\delta}}}+\|\chi_{cl}\partial_{t}^{j_{1}}\partial_{\theta}^{j_{2}}p\|_{L^{\infty}W^{1,q_{\delta}}})\notag \\
&+\sum_{j_{1}+j_{2}\leq m}(\|\chi_{cl}\partial_{t}^{j_{1}}\partial_{\theta}^{j_{2}}v\|_{L^{\infty}H^{1}}+\|\chi_{cl,\Sigma}\partial_{t}^{j_{1}}\partial_{\theta}^{j_{2}}\eta\|_{L^{\infty}H^{\frac{3}{2}-\alpha}}+\|\chi_{cl}\partial_{t}^{j_{1}}\partial_{\theta}^{j_{2}}p\|_{L^{\infty}L^{2}}) \lesssim  \sqrt{\mathcal{D}},
\end{align}
and
\begin{align}
&\sum_{j_{1}+j_{2}\leq m-2}(\|\chi_{cl}\partial_{t}^{j_{1}}\partial_{\theta}^{j_{2}}v\|_{L^{\infty}W^{2,q_{\delta}}}+\|\chi_{cl,\Sigma}\partial_{t}^{j_{1}}\partial_{\theta}^{j_{2}}\eta\|_{L^{\infty}W^{3-\frac{1}{q_{\delta}},q_{\delta}}}+\|\chi_{cl}\partial_{t}^{j_{1}}\partial_{\theta}^{j_{2}}p\|_{L^{\infty}W^{1,q_{\delta}}})\notag \\
&+\sum_{j_{1}+j_{2}\leq m-1}(\|\chi_{cl}\partial_{t}^{j_{1}}\partial_{\theta}^{j_{2}}v\|_{L^{\infty}H^{1}}+\|\chi_{cl,\Sigma}\partial_{t}^{j_{1}}\partial_{\theta}^{j_{2}}\eta\|_{L^{\infty}H^{\frac{3}{2}-\alpha}}+\|\chi_{cl}\partial_{t}^{j_{1}}\partial_{\theta}^{j_{2}}p\|_{L^{\infty}L^{2}}) \lesssim  \sqrt{\mathcal{E}}.
\end{align}
\textbf{The convection term.} We first consider  $\mathcal{O}_{\mathrm{conv}} := \chi_{\mathrm{cl}} \bigl( v\cdot\nabla_{\mathcal{A}}v \bigr),$ where  $\mathcal{A }:= (\mathcal{A}_\eta\circ\Phi)\mathcal{A}_\Phi$ is the coefficient matrix obtained after pulling the nonlinear derivatives back to the local chart.

Since the cutoff functions are independent of $\theta$,
\begin{align}
\partial_\theta^{m}\mathcal{O}_{\mathrm{conv}} ={}& \chi_{\mathrm{cl}} \sum_{\ell+a+b=m} C_{\ell,a,b} \bigl( \partial_\theta^\ell v \bigr) \cdot \left( \partial_\theta^a\mathcal{A }\right) \nabla \partial_\theta^b v.
\end{align}
Since $\chi_*=1$ on $\func{supp}\chi_{\mathrm{cl}}$, we may estimate this by localized quantities:
\begin{equation}
\left\| \partial_\theta^{m} \mathcal{O}_{\mathrm{conv}} \right\|_{L_\theta^2L_\perp^{q_\delta}} \lesssim{}  \sum_{\ell+a+b=m} \left\| \partial_\theta^\ell(\chi_*v) \right\|_{L_\theta^{s_{{i}}} L_\perp^{r_\delta}} \left\| \partial_\theta^a\mathcal{A }\right\|_{L^{d_{{i}}}L^\infty} \left\| \nabla\partial_\theta^b(\chi_*v) \right\|_{L_\theta^{s_{{i}}}L_\perp^2} + \mathcal{R}_{\mathrm{conv},j},
\end{equation}
where $\mathcal{R}_{\mathrm{conv},j}$ contains only terms involving derivatives of the fixed cutoff functions. They are of lower order in $\nabla_{r,h}$ derivatives and can be controlled by the main term (the first term) on the right-hand side. Moreover, $s_{i}$ and $d_{i}$ can only be $2$ or infinity, depending on the number of derivatives taken on that term. For example, $s_{i}$ in  $\| \partial_\theta^\ell(\chi_*v) \|_{L_\theta^{s_{{i}}} L_\perp^{r_\delta}}$ depends on the index $\ell$. We set $s_{i}=2$ if $\ell\geq m-1$ and $s_{i}=\infty$ if $\ell\leq m-2$. $d_{i}=2$ if $\ell\geq m$ and $d_{i}=+\infty $ if $\ell\leq m-1$. Moreover, the index depends only on the order of the derivatives falling on the term inside the norm. It does not depend on the form of this term. Specifically, for example $\partial_{t}^{j}(\partial_{t}\bar{\eta})$, $s_{i}$ depends only on $j$, not $j+1$! Therefore,
\begin{align}
&\left\| \partial_\theta^{m} \mathcal{O}_{\mathrm{conv}} \right\|_{L_\theta^2L_\perp^{q_\delta}} \lesssim{}  \sum_{\ell+a+b=m} \left\| \partial_\theta^\ell(\chi_*v) \right\|_{L_\theta^{s_{{i}}} H^{1}} (1+\left\| \partial_\theta^a\chi_{*}\eta \right\|_{L^{d_{{i}}}W^{3-\frac{1}{q_{\delta}},q_{\delta}}}+\left\| \partial_\theta^{a}(\chi_{*}\partial_{\theta}\eta) \right\|_{L^{d_{{i}}}W^{2-\frac{1}{q_{\delta}},q_{\delta}}})\notag \\
&\times (\left\| \partial_\theta^b(\chi_*v) \right\|_{L_\theta^{s_{{i}}}H_\perp^1}+\left\| \partial_\theta^{b}(\chi_* \partial_{\theta}v) \right\|_{L_\theta^{s_{{i}}}L_\perp^2}) \lesssim{} \sqrt{\mathcal{D}\mathcal{E}}.
\end{align}
This estimate holds since $2(m-1)>m$ when $m\geq 3$ (we only take $m$th-order derivatives for the elliptic estimate). Therefore, among $a,b,c$, at most one number is greater than or equal to $m-1$. So there is at most one norm equipped with $L_{\theta}^{2}$.

\textbf{The time-dependent coordinate term.}

The nonlinear term generated by the time-dependent global diffeomorphism has the form  $\mathcal{O}_{\mathrm{time}} := \chi_{\mathrm{cl}} WK_\eta \partial_t\bar\eta\, \partial_3v.$ Applying $\partial_\theta^j$ and using the Leibniz rule gives
\begin{equation*}
\partial_\theta^j \mathcal{O}_{\mathrm{time}} = \chi_{\mathrm{cl}} \sum_{a+b+c+d=j} C_{a,b,c,d} \partial_\theta^aW\, \partial_\theta^bK_\eta\, \partial_t\partial_\theta^c\bar\eta\, \partial_3\partial_\theta^dv.
\end{equation*}
Hence, by H\"older's inequality,
\begin{align}
&\sum_{j=0}^m \left\| \partial_\theta^j \mathcal{O}_{\mathrm{time}} \right\|_{L_\theta^2L_\perp^{q_\delta}} \lesssim{} \sum_{j=0}^{m}\sum_{a+b+c+d=j,(a,b,c,d)\neq m} \|\partial_\theta^aW\|_{L^{s_{i}}L^{\infty}}\, \|\partial_\theta^bK_\eta\|_{L^{s_{i}}L^{\infty}} \|\partial_t\partial_\theta^c\bar\eta\|_{L^{d_{i}}L^{\infty}} \|\partial_3\partial_\theta^dv \|_{L^{s_{i}}L^{2}}\notag \\
&+\|\partial_\theta^mW\|_{L^{2}L^{\infty}}\, \| K_\eta\|_{L^{\infty}L^{\infty}} \|\partial_t \bar\eta\|_{L^{\infty}L^{\infty}} \|\partial_3 v \|_{L^{\infty}L^{2}} +\|W\|_{L^{\infty}L^{\infty}}\, \|\partial_\theta^{m} K_\eta\|_{L^{2}L^{\infty}} \|\partial_t \bar\eta\|_{L^{\infty}L^{\infty}} \|\partial_3 v \|_{L^{\infty}L^{2}}\notag \\
&+\|W\|_{L^{\infty}L^{\infty}}\, \| K_\eta\|_{L^{\infty}L^{\infty}} \|\partial_t\partial_{\theta}^{m} \bar\eta\|_{L^{2}L^{\infty}} \|\partial_3 v \|_{L^{\infty}L^{2}} +\|W\|_{L^{\infty}L^{\infty}}\, \| K_\eta\|_{L^{\infty}L^{\infty}} \|\partial_t \bar\eta\|_{L^{\infty}L^{\infty}} \|\partial_3 \partial_{\theta}^{m}v \|_{L^{2}L^{2}}\notag \\
&\lesssim{} \sum_{j=0}^{m}\sum_{a+b+c+d=j} (1+\|\partial_{\theta}^a\chi_{*}\eta\|_{L^{s_{i}}W^{3-\frac{1}{q_{\delta}},q_{\delta}}}+\|\partial_{\theta}^{a+1}\chi_{*}\eta\|_{L^{d_{i}}W^{2-\frac{1}{q_{\delta}},q_{\delta}}})\,\notag \\
&(\|\partial_\theta^b\chi_{*}\eta\|_{L^{s_{i}}W^{3-\frac{1}{q_{\delta}},q_{\delta}}}+\|\partial_{\theta}^{b+1}\chi_{*}\eta\|_{L^{d_{i}}W^{2-\frac{1}{q_{\delta}},q_{\delta}}}+1)\label{eq:higher-time-map-estimate} \\
&\times\|\partial_t\partial_\theta^c\chi_{*}\eta\|_{L^{d_{i}}H^{\frac{3}{2}-\alpha}}(\|\partial_\theta^d\chi_{*}v \|_{L^{s_{i}}H^{1}}+\|\partial_\theta^{d+1}\chi_{*}v \|_{L^{d_{i}}L^{2}}) \lesssim \sqrt{\mathcal{E}\mathcal{D}}.\notag
\end{align}
Again, this estimate holds since there is at most one index among $(a,b,c,d)$ that is greater than $m-2$. Therefore, at most one index among the $s_{i}$ and $d_{i}$ equals $2$; the other indices equal $\infty$.

\textbf{The nonlinear $\mathcal{A}_\eta$-difference terms.}

Let
\begin{equation*}
\mathscr D_\eta(v,p) := \func{div}S(p,v) - \func{div}_{\mathcal{A}_\eta} S_{\mathcal{A}_\eta}(p,v).
\end{equation*}
After pulling back to the local chart, this has the schematic structure
\begin{equation}
\mathscr D_\eta(v,p) = \mathfrak{a}_\eta\nabla^2v + \mathfrak{b}_\eta\nabla p + \mathfrak{c}_\eta(\nabla v+p),  \label{eq:Aeta-difference-schematic}
\end{equation}
where $\mathfrak{a}_\eta, \mathfrak{b}_\eta = O(\mathcal{A}_\eta-I),$ and $\mathfrak{c}_\eta$ contains derivatives of $\mathcal{A}_\eta$.

After applying $\chi_{\mathrm{cl}}\partial_\theta^j$, the top-order terms are
\begin{equation*}
\mathfrak{a}_\eta \nabla^2 \left( \chi_{\mathrm{cl}}\partial_\theta^jv \right), \qquad \mathfrak{b}_\eta \nabla \left( \chi_{\mathrm{cl}}\partial_\theta^jp \right).
\end{equation*}
All remaining terms either contain a $\theta$-derivative of a geometric coefficient or a derivative of the cutoff function and hence are of lower order. By the definition of $\mathcal{A}_{\eta}$ and the same argument as for $\mathcal{G}_{\eta,0}$, we have
\begin{align}
&\sum_{j=0}^{m-1}(\|\partial_{\theta}^{j}\mathcal{A}_{\eta}\|_{L^{2}L^{+\infty}}+\|\partial_{\theta}^{j}\mathcal{A}_{\eta}\|_{L^{2}W^{1,q_{r_{\delta}}}}+\|\partial_{\theta}^{j}\mathcal{N}\|_{L^{2}W^{1,q_{\delta}}}+\|\partial_{\theta}^{j+1}\mathcal{A}_{\eta}\|_{L^{2}L^{q_{r_\delta}}})\notag \\
&\lesssim 1+\sum_{j=0}^{m-1}(\|\partial_{\theta}^{j}\bar{\eta}\|_{L^{2}W^{3,q_{\delta}}}+\|\partial_{\theta}^{j+1}\bar{\eta}\|_{L^{2}W^{2,q_{\delta}}}+\|\partial_{\theta}^{j+2}\bar{\eta}\|_{L^{2}H^{1-\alpha}}+\|\partial_{\theta}^{j}\eta\|_{L^{2}W^{2,q_{\delta}}}+\|\partial_{\theta}^{j+1}\eta\|_{L^{2}H^{1}} )+\sqrt{\mathcal{E}}^{2}\notag \\
&\lesssim1+\mathcal{E}+\sum_{j=0}^{m-1}(\|\chi_{cl,\Sigma}\partial_{\theta}^{j}\eta\|_{L^{2}W^{3-\frac{1}{q_{\delta}},q_{\delta}}(\mathcal{W})}+\|\chi_{bet,\Sigma}\partial_{\theta}^{j}\eta\|_{L^{2}W^{3-\frac{1}{q_{\delta}},q_{\delta}}(\mathcal{W})})\notag \\
&+\sum_{j=0}^{m-1}(\|\chi_{cl,\Sigma}\partial_{\theta}^{j+1}\eta\|_{L^{2}H^{\frac{3}{2}-\alpha}(\mathcal{W})}+\|\chi_{bet,\Sigma}\partial_{\theta}^{j+1}\eta\|_{L^{2}H^{\frac{3}{2}-\alpha}(\mathcal{W})}\label{est:background_geometry_j} \\
&+\|\chi_{1,\Sigma}\partial_{\theta}^{j+2}\eta\|_{L^{2}H^{\frac{1}{2}-\alpha}(\mathcal{W})}+\|\chi_{bet,\Sigma}\partial_{\theta}^{j+2}\eta\|_{L^{2}H^{\frac{1}{2}-\alpha}(\mathcal{W})}) \lesssim 1+\sqrt{\mathcal{E}},\notag
\end{align}
and similarly
\begin{equation}
\|\partial_{\theta}^{m}\mathcal{A}_{\eta}\|_{L^{2}L^{+\infty}}+\|\partial_{\theta}^{m}\mathcal{A}_{\eta}\|_{L^{2}W^{1,q_{r_{\delta}}}}+\|\partial_{\theta}^{m}\mathcal{N}\|_{L^{2}W^{1,q_{\delta}}}+\|\partial_{\theta}^{m+1}\mathcal{A}_{\eta}\|_{L^{2}L^{q_{r_\delta}}} \lesssim 1+\sqrt{\mathcal{D}}.
\end{equation}

Consequently, using Lemma \ref{lem:derivative_replacement}, we can regard any whole space spatial derivative as $\nabla_{r,h}$
\begin{align}
&\sum_{j=0}^m \left\| \chi_{\mathrm{cl}} \partial_\theta^j \mathscr D_\eta(v,p) \right\|_ {L_\theta^2L_\perp^{q_\delta}} \lesssim \|\mathcal{A}_\eta-I\|_{L^{\infty}L^\infty} \left( \|\chi_{*}\nabla^{2}\partial_{\theta}^{m}v\|_{L^{2}L^{q_{\delta}}}+\|\chi_{*}\nabla \partial_{\theta}^{m}p\|_{L^{2}L^{q_{\delta}}}\right)\notag \\
&+\|\nabla \mathcal{A}_{\eta}\|_{L^{\infty}L^{q_{r_{\delta }}}}\|\mathcal{A}_\eta-I\|_{L^{\infty}L^\infty} \left( \|\chi_{*}\partial_{\theta}^{m}\nabla v\|_{L^{2}L^{q_{r_\delta}}}+\|\chi_{*}\partial_{\theta}^{m}p\|_{L^{2}L^{q_{r_\delta}}}\right)\notag \\
&+ \|\partial_{\theta}^{m}(\mathcal{A}_\eta-I)\|_{L^{2}L^\infty} \left( \|\chi_{*}\nabla^{2} v\|_{L^{\infty}L^{q_{\delta}}}+\|\chi_{*}\nabla p\|_{L^{\infty}L^{q_{\delta}}}\right)\label{eq:higher-Aeta-difference} \\
&+\|\nabla \partial_{\theta}^{m}\mathcal{A}_{\eta}\|_{L^{2}L^{q_{r_{\delta }}}}\|\mathcal{A}_\eta-I\|_{L^{\infty}L^\infty} \left( \|\chi_{*}\nabla v\|_{L^{\infty}L^{q_{r_\delta}}}+\|\chi_{*}p\|_{L^{\infty}L^{q_{r_\delta}}}\right)\notag \\
&+\|\nabla \partial_{\theta}^{m}(I-\mathcal{A}_{\eta})\|_{L^{2}L^{q_{r_{\delta }}}} \left( \|\chi_{*}\nabla v\|_{L^{\infty}L^{q_{r_\delta}}}+\|\chi_{*}p\|_{L^{\infty}L^{q_{r_\delta}}}\right)\notag \\
&+\sum_{n=0}^{m-1}C(\sum_{l=0}^{n}\|\partial_{\theta}^{l}(\mathcal{A}_\eta-I)\|_{L^{s_{i}}L^\infty})) (\|\chi_{*}\partial_{\theta}^{n-\ell}v\|_{L^{s_{i}}W^{2,q_{\delta}}}+\|\chi_{*}\partial_{\theta}^{n-\ell}p\|_{L^{s_{i}}W^{1,q_{\delta}}})\notag \\
&+\sum_{n=0}^{m-1}C(\sum_{l=0}^{n}\|\partial_{\theta}^{l}(\nabla(\mathcal{A}_\eta-I)\|_{L^{s_{i}}L^{q_{r_\delta}}})) (\|\chi_{*}\partial_{\theta}^{n-\ell}v\|_{L^{s_{i}}W^{2,q_{\delta}}}+\|\chi_{*}\partial_{\theta}^{n-\ell}p\|_{L^{s_{i}}W^{1,q_{\delta}}}) \lesssim \sqrt{\mathcal{E}\mathcal{D}} ,\notag
\end{align}
where we used the Sobolev embedding in 2-D that $W^{1,q_{\delta}}\hookrightarrow L^{q_{r_{\delta}}}$.

For the divergence term, using a similar estimate, we have
\begin{align}
\sum_{j=0}^m \left\| \chi_{\mathrm{cl}} \partial_\theta^j ( \func{div}v - \func{div}_{\mathcal{A}_\eta}v) \right\|_ {L_\theta^2W_\perp^{1,q_\delta}} \lesssim{}& \sqrt{\mathcal{E}\mathcal{D}}.  \label{eq:higher-Aeta-divergence}
\end{align}

\textbf{The nonlinear kinematic term.}

Recall that
\begin{equation*}
\mathcal{F}_{\mathrm{cl}}^4 = \frac{1}{|N_0|} v\cdot(\mathcal{N}-\mathcal{N}_0), \qquad \mathcal{N}-\mathcal{N}_0 =
\begin{pmatrix}
-\nabla\eta \\
0\end{pmatrix} =\begin{pmatrix} -A_{X}^{-T}\nabla\mathfrak{h}_{\eta} \\
0\end{pmatrix}.
\end{equation*}
Hence the localized nonlinear term is
\begin{equation*}
\mathcal{O}_{\mathrm{kin}} := \frac{\chi_{\mathrm{cl},\Sigma}}{|\mathcal{N}_0|} v\cdot(\mathcal{N}-\mathcal{N}_0).
\end{equation*}
Since $\chi_{*,\Sigma}=1$ on $\func{supp}\chi_{1,\Sigma}$, the Leibniz rule and the trace estimate give
\begin{align}
&\sum_{j=0}^m \left\| \partial_\theta^j \mathcal{Q}_{\mathrm{kin}} \right\|_ {L_\theta^2 W^{2-\frac1{q_\delta},q_\delta}(\Gamma_+)}\label{eq:higher-kinematic-nonlinear} \\
&\lesssim C\sum_{j=0}^m \sum_{n=0}^{j}\left( \left\| \chi_{\mathrm{cl}}\partial_\theta^nv \right\|_ {L_\theta^{s_{i}}W_\perp^{2,q_\delta}} \right) \left( \left\| \chi_{*,\Sigma}\partial_\theta^{j-n}\eta \right\|_ {L_\theta^{s_{i}} W^{3-\frac1{q_\delta},q_\delta}} +\sum_{k=0}^{j-n}\left\| \chi_{*,\Sigma}\partial_\theta^{j-n-k+1}\eta \right\|_ {L_\theta^{d_{i}} H^{\frac{3}{2}-\alpha}}\right) \lesssim \sqrt{\mathcal{E}\mathcal{D}}.\notag
\end{align}

\textbf{The nonlinear pressure tensor boundary terms.}

Let $\mathcal{O}_{\mathrm{dyn}}$ denote the nonlinear part of the free-surface stress forcing, excluding the cutoff commutators and the flattening-map differences already estimated above. Schematically,
\begin{equation*}
\mathcal{O}_{\mathrm{dyn}} = \chi_{\mathrm{cl},\Sigma} \left[ \left( S(p,v)-S_{\mathcal{A}_\eta}(p,v) \right)\mathcal{N}_0 + S_{\mathcal{A}_\eta}(v,p)(\mathcal{N}-\mathcal{N}_0) \right].
\end{equation*}
Then, similarly to the estimate for the pressure tensor in the momentum equation above, we have
\begin{align}
&\sum_{j=0}^m \left\| \partial_\theta^j \mathcal{O}_{\mathrm{dyn}} \cdot\tau_{0,+} \right\|_ {L_\theta^2 W^{1-\frac1{q_\delta},q_\delta}(\Gamma_+)} \lesssim \sum_{j=0}^{m}\sum_{a+b=j}\|\partial_{\theta}^{a}(I-\mathcal{A}_{\eta})\|_{L^{s_{i}}W^{2,q_{\delta}}}(\|\partial_{\theta }^{b}v_{\mathrm{cl}}\|_{L^{s_{i}}W^{2,q_{\delta}}}\notag \\
&+\|\partial_{\theta}^{b}v_{\mathrm{bet}}\|_{L^{s_{i}}W^{2,q_{\delta}}}+\|\partial_{\theta}^{b}p_{\mathrm{cl}}\|_{L^{s_{i}}W^{1,q_{\delta}}}+\|\partial_{\theta}^{b}p_{\mathrm{bet}}\|_{L^{s_{i}}W^{1,q_{\delta}}}) +\sum_{j=0}^{m}\sum_{a+b=j}\|\partial_{\theta}^{a}(\mathcal{N}-\mathcal{N}_{0})\|_{L^{s_{i}}W^{1,q_{\delta}}}(\|\partial_{\theta}^{b}v_{\mathrm{cl}}\|_{L^{s_{i}}W^{2,q_{\delta}}}\notag \\
&+\|\partial_{\theta}^{b}v_{\mathrm{bet}}\|_{L^{s_{i}}W^{2,q_{\delta}}}+\|\partial_{\theta}^{b}p_{\mathrm{cl}}\|_{L^{s_{i}}W^{1,q_{\delta}}}+\|\partial_{\theta}^{b}p_{\mathrm{bet}}\|_{L^{s_{i}}W^{1,q_{\delta}}}) \lesssim \sqrt{\mathcal{E}\mathcal{D}} .\label{eq:higher-dynamic-nonlinear}
\end{align}

\textbf{Curvature term in the boundary condition}

We now estimate the curvature term in the capillary boundary condition. We have
\begin{align}
&\sum_{j=0}^{m}\|\chi_{\mathrm{cl},\Sigma}\partial_{\theta}^{j}(g\mathfrak{h}_{\eta}-\mathcal{K}(\mathfrak{h}_{\eta}))(\mathcal{N}-\mathcal{N}_{0})\|_{L^{2}W^{1-\frac{1}{q_{\delta}},q_{\delta}}} \\
&\lesssim (\|\partial_{\theta}^{m}\nabla^{2}(\chi_{\mathrm{cl},\Sigma}\mathfrak{h}_{\eta})\|_{L^{2}W^{1-\frac{1}{q_{\delta}},q_{\delta}}}+\|\partial_{\theta}^{m}\nabla(\chi_{\mathrm{cl},\Sigma}\mathfrak{h}_{\eta})\|_{L^{2}W^{1-\frac{1}{q_{\delta}},q_{\delta}}}+\|\partial_{\theta}^{m}(\chi_{\mathrm{cl},\Sigma}\mathfrak{h}_{\eta})\|_{L^{2}W^{1-\frac{1}{q_{\delta}},q_{\delta}}})\|\nabla \mathfrak{h}_{\eta}\|_{L^{\infty}W^{1,q_{\delta}}}\notag \\
&+(\|\nabla^{2}(\chi_{\mathrm{cl},\Sigma}\mathfrak{h}_{\eta})\|_{L^{\infty}W^{1-\frac{1}{q_{\delta}},q_{\delta}}}+\|\nabla(\chi_{\mathrm{cl},\Sigma}\mathfrak{h}_{\eta})\|_{L^{\infty}W^{1-\frac{1}{q_{\delta}},q_{\delta}}}+\|(\chi_{\mathrm{cl},\Sigma}\mathfrak{h}_{\eta})\|_{L^{\infty}W^{1-\frac{1}{q_{\delta}},q_{\delta}}})\|\nabla \partial_{\theta}^{m}\mathfrak{h}_{\eta}\|_{L^{2}W^{1,q_{\delta}}}\notag \\
&+\sum_{n=1}^{m}(\|\partial_{\theta}^{n}\nabla^{2}(\chi_{\mathrm{cl},\Sigma}\mathfrak{h}_{\eta})\|_{L^{s_{i}}W^{1-\frac{1}{q_{\delta}},q_{\delta}}}+\|\partial_{\theta}^{n}\nabla(\chi_{\mathrm{cl},\Sigma}\mathfrak{h}_{\eta})\|_{L^{s_{i}}W^{1-\frac{1}{q_{\delta}},q_{\delta}}} +\|\partial_{\theta}^{n}(\chi_{\mathrm{cl},\Sigma}\mathfrak{h}_{\eta})\|_{L^{s_{i}}W^{1-\frac{1}{q_{\delta}},q_{\delta}}})\|\nabla \partial_{\theta}^{m-n}\mathfrak{h}_{\eta}\|_{L^{s_{i}}W^{1,q_{\delta}}}\notag \\
&+\sum_{j=0}^{m-1}\sum_{n=0}^{j}(\|\partial_{\theta}^{n}\nabla^{2}(\chi_{\mathrm{cl},\Sigma}\mathfrak{h}_{\eta})\|_{L^{s_{i}}W^{1-\frac{1}{q_{\delta}},q_{\delta}}}+\|\partial_{\theta}^{n}\nabla(\chi_{\mathrm{cl},\Sigma}\mathfrak{h}_{\eta})\|_{L^{s_{i}}W^{1-\frac{1}{q_{\delta}},q_{\delta}}}\notag \\
&+\|\partial_{\theta}^{n}(\chi_{\mathrm{cl},\Sigma}\mathfrak{h}_{\eta})\|_{L^{s_{i}}W^{1-\frac{1}{q_{\delta}},q_{\delta}}})\|\nabla \partial_{\theta}^{j-n}\mathfrak{h}_{\eta}\|_{L^{s_{i}}W^{1,q_{\delta}}} \lesssim \sqrt{\mathcal{E}\mathcal{D}}.\notag
\end{align}
\textbf{The curvature remainder.}

Let $\mathcal{R}_H(\eta)$ denote the nonlinear Taylor remainder of the mean-curvature operator. The localized curvature remainder is then
\begin{equation}
\begin{aligned}
\operatorname{div}_{\Sigma}\mathcal R(\nabla\zeta,\nabla\eta) ={}& D_1\mathcal R(\nabla\zeta,\nabla\eta) :\nabla_{\Sigma}^{2}\zeta+ D_2\mathcal R(\nabla\zeta,\nabla\eta) :\nabla_{\Sigma}^{2}\eta.
\end{aligned}
\end{equation}
Then,
\begin{equation*}
\mathcal{O}_{\mathrm{curv}} := \chi_{\mathrm{cl},\Sigma}\operatorname{div}_{\Sigma} \mathcal{R}(\nabla\zeta,\nabla\eta)=\chi_{\mathrm{cl},\Sigma}\big(D_1\mathcal R(\nabla\zeta,\nabla\eta) :\nabla_{\Sigma}^{2}\zeta+ D_2\mathcal R(\nabla\zeta,\nabla\eta) :\nabla_{\Sigma}^{2}\eta\big).
\end{equation*}
Since $\chi_{*,\Sigma}=1$ on $\func{supp}\chi_{\mathrm{cl},\Sigma}$, every factor may be localized by $\chi_{*,\Sigma}$. Thus,
\begin{align}
&\sum_{j=0}^m \left\| \chi_{*,\Sigma}\partial_\theta^j \mathcal{O}_{\mathrm{curv}} \right\|_ {L_\theta^2 W^{1-\frac1{q_\delta},q_\delta}(\Gamma_+)}\label{eq:higher-curvature-remainder} \\
&\lesssim C\sum_{j=0}^m \sum_{n=0}^{j}\left( \left\| \chi_{*,\Sigma}\partial_\theta^n\eta \right\|_ {L_\theta^{s_{i}} W^{3-\frac1{q_\delta},q_\delta}} + \left\| \chi_{*,\Sigma}\partial_\theta^{n+1}\eta \right\|_ {L_\theta^{s_{i}} H^{\frac{3}{2}-\alpha}}+ \left\| \chi_{*,\Sigma}\partial_\theta^{n+2}\eta \right\|_ {L_\theta^{s_{i}} H^{\frac{1}{2}-\alpha}}\right)\notag \\
&\times \left( \left\| \chi_{*,\Sigma}\partial_\theta^{j-n}\eta \right\|_ {L_\theta^{s_{i}} W^{3-\frac1{q_\delta},q_\delta}} + \left\| \chi_{*,\Sigma}\partial_\theta^{j-n+1}\eta \right\|_ {L_\theta^{s_{i}}H^{\frac{3}{2}-\alpha}}\right) \lesssim \sqrt{\mathcal{E}\mathcal{D}}.\notag
\end{align}

\textbf{The temporal derivative term in the kinematic boundary condition.} We have
\begin{align}
\|\chi_{cl,\Sigma}\partial_{t}\partial_{\theta}^{m}\eta\|_{L^{2}W^{2-\frac{1}{q_{\delta}},q_{\delta}}}\lesssim\|\chi_{cl,\Sigma}\partial_{t}\partial_{\theta}^{m}\eta\|_{H^{\frac{3}{2}-\alpha}}\lesssim \sqrt{\mathscr{E}_{1,m,\mathrm{energy},\mathrm{cl}}}.
\end{align}

\textbf{The nonlinear contact-line term.}

For the nonlinear function $\hat{\mathcal{W}}$ on the contact line, using the same computation as for the previous temporal derivative term, we have
\begin{align}
\sum_{j=0}^{m}\chi_{\mathrm{cl},\Sigma}\partial_{\theta}^{j}\hat{\mathcal{W}}(b_{\eta}\partial_{t}\mathfrak{h}_{\eta})\lesssim \sqrt{\mathcal{E}\mathcal{D}}.
\end{align}

Combining all of the estimates above, we prove Theorem \ref{thm:Higher-order contact line elliptic estimate}.
\end{proof}

Using the contact line elliptic estimate, we now prove Theorem~\ref{thm:intermediate_estimate}. At the top elliptic level it reads as follows.

\begin{theorem}
\label{thm:intermediate_elliptic_weak} For each fixed pair $(k_1,k_2)$ with $k_1+k_2=m$, the following intermediate elliptic estimate holds
\begin{align}
&(\|\chi_{\mathrm{bet}}v_{k_{1},k_{2}}\|_{L^{2}W^{2,q_{\delta}}}^{2}+\|\chi_{\mathrm{bet}}\eta_{k_{1},k_{2}}\|_{L^{2}W^{3-\frac{1}{q_{\delta}},q_{\delta}}}^{2}+\|\chi_{\mathrm{bet}}p_{k_{1},k_{2}}\|_{L^{2}W^{1,q_{\delta}}}^{2} ) \\
&\lesssim \mathcal{E}\mathcal{D}+\frac{C(\varepsilon)}{\varepsilon} \sum_{0\leq j_{1}\leq k_{1}}\sum_{0\leq j_{2}\leq k_{2}}\mathscr{E}_{j_{1},j_{2},\mathrm{energy},\mathrm{bet}}+\sum_{j=0}^{k_{2}+1}\mathscr{E}_{k_{1},j,\mathrm{energy},\mathrm{bet}}+\sum_{j=0}^{k_{2}}\mathscr{E}_{k_{1}+1,j,\mathrm{energy},bet}+\frac{C(\varepsilon)}{\varepsilon} \mathcal{D}_{\mathrm{int}}.
\end{align}
For each fixed pair $(k_1,k_2)$ with $k_1+k_2=m-1$, we have
\begin{align}
&(\|\chi_{\mathrm{bet}}v_{k_{1},k_{2}}\|_{L^{2}W^{2,q_{\delta}}}^{2}+\|\chi_{\mathrm{bet}}\eta_{k_{1},k_{2}}\|_{L^{2}W^{3-\frac{1}{q_{\delta}},q_{\delta}}}^{2}+\|\chi_{\mathrm{bet}}p_{k_{1},k_{2}}\|_{L^{2}W^{1,q_{\delta}}}^{2} ) \\
&\lesssim \mathcal{E}^{2}+\frac{C(\varepsilon)}{\varepsilon} \sum_{0\leq j_{1}\leq k_{1}}\sum_{0\leq j_{2}\leq k_{2}}\mathscr{E}_{j_{1},j_{2},\mathrm{energy},\mathrm{bet}}+\sum_{j=0}^{k_{2}+1}\mathscr{E}_{k_{1},j,\mathrm{energy},\mathrm{bet}}+\sum_{j=0}^{k_{2}}\mathscr{E}_{k_{1}+1,j,\mathrm{energy},bet}+\frac{C(\varepsilon)}{\varepsilon} \mathcal{E}_{\mathrm{int}}.
\end{align}
\end{theorem}

\begin{proof}

We choose $\chi_{\mathrm{int}-}$ to vanish on a neighborhood of $\operatorname{supp}(\chi_{\mathrm{cl}}+\chi_{\mathrm{bet}})$, so that
\begin{align}
\operatorname{supp}(\chi_{\mathrm{cl}}+\chi_{\mathrm{bet}}) \subset\{x\in\Omega:1-\chi_{\mathrm{int}-}(x)=1\},
\end{align}
and $\operatorname{supp}\{1-\chi_{\mathrm{int}-}\}\in \{x\in \Omega|\chi_{\mathrm{local}}=1\}.$ Put $\chi_{\mathrm{bet}+}=\chi_{\mathrm{cl}}+\chi_{\mathrm{bet}}$. The highest derivatives retain this cutoff. For every differentiated velocity $v$ and pressure $p$,
\begin{align}
\chi_{\mathrm{bet}+}\partial_i\partial_j v ={}&\partial_i\partial_j(\chi_{\mathrm{bet}+}v) -(\partial_i\chi_{\mathrm{bet}+})\partial_j v -(\partial_j\chi_{\mathrm{bet}+})\partial_i v -(\partial_i\partial_j\chi_{\mathrm{bet}+})v,\\
\chi_{\mathrm{ber}+}\partial_ip ={}&\partial_i(\chi_{\mathrm{bet}+}p)-(\partial_i\chi_{\mathrm{bet}+})p.
\end{align}

The cutoff terms are supported where the cutoff $1-\chi_{\mathrm{int}-}$ equals one. For $q=q_\delta\le p$, the lower-order bounds used there are
\begin{align}
\|v\|_{L^2_\theta W^{s,q}_{r,h}}\le C\|v\|_{W^{2,p}},\qquad \|p\|_{L^2_\theta W^{s-1,q}_{r,h}}\le C\|p\|_{W^{1,p}},
\end{align}
for $1<s\leq \frac{5}{2}-\frac{1}{p}$ given that $1<q_{\delta}<p<2$. These estimates control the cutoff terms through the interior norms since the cutoff terms have one fewer spatial derivative. For the nonlinear terms, since we apply $m$th-order $t$ or $\theta$ derivatives, we multiply the highest-order part by $\chi_{\mathrm{bet}+}$ and the lower-order part by $1-\chi_{\mathrm{int}-}=(\chi_{\mathrm{bet}+})+(1-\chi_{\mathrm{int}-}-\chi_{\mathrm{bet}+})$. The interior part multiplied by $(1-\chi_{\mathrm{int}-}-\chi_{\mathrm{bet}+})$ has lower-order derivatives and can therefore be bounded by the interior norms. They do not remove the weight $\chi_{\mathrm{bet}+}$ from the highest derivatives. Therefore, in the support of $1-\chi_{\mathrm{int}-}-\chi_{\mathrm{bet}+}$, $(\p_{\theta}^{j}v,\p_{\theta}^{j}p,\p_{\theta}^{j}\eta)$ included in nonlinear terms can be bounded by $\mathcal{E}_{\mathrm{int}}$ or $\mathcal{D}_{\mathrm{int}}$, respectively. Away from the contact line region, the regularity is better than near the contact line. So the nonlinear estimate and forcing term estimate derived for $\mathcal{F}_{\mathrm{cl}}^{i,j}$ still work for $\mathcal{F}_{\mathrm{bet}+}^{i,j}$ in a similar way.

Fix $(k_1,k_2)$ with $k_1+k_2=m$. Using a proof similar to that of Theorem \ref{thm:Higher-order contact line elliptic estimate}, replacing $\chi_{\mathrm{cl}}$ by $\chi_{\mathrm{cl}}+\chi_{\mathrm{bet}}$ in the proof, and using the fact that $1-\chi_{\mathrm{int}-}$ equals $1$ in the support of $\chi_{\mathrm{bet}+}$, we similarly obtain
\begin{align}
(\|(\chi_{\mathrm{cl}}+\chi_{\mathrm{bet}})v_{k_{1},k_{2}}\|_{L^{2}W^{2,q_{\delta}}}^{2}+\|(\chi_{\mathrm{cl},\Sigma}+\chi_{\mathrm{bet},\Sigma})\eta_{k_{1},k_{2}}\|_{L^{2}W^{3-\frac{1}{q_{\delta}},q_{\delta}}}^{2}+\|(\chi_{\mathrm{cl}}+\chi_{\mathrm{bet}})p_{k_{1},k_{2}}\|_{L^{2}W^{1,q_{\delta}}}^{2} )  \notag \\
\lesssim\mathcal{E}\mathcal{D}+\frac{C(\varepsilon)}{\varepsilon} \sum_{0\leq j_{1}\leq k_{1}}\sum_{0\leq j_{2}\leq k_{2}}\mathscr{E}_{j_{1},j_{2},\mathrm{energy},\mathrm{cl}}+\sum_{j=0}^{k_{2}+1}\mathscr{E}_{k_{1},j,\mathrm{energy},\mathrm{cl}}+\sum_{j=0}^{k_{2}}\mathscr{E}_{k_{1}+1,j,\mathrm{energy},cl}+\frac{C(\varepsilon)}{\varepsilon} \mathcal{D}_{\mathrm{int}}.
\end{align}
Therefore
\begin{align}
&(\|(\chi_{\mathrm{bet}})v_{k_{1},k_{2}}\|_{L^{2}W^{2,q_{\delta}}}^{2}+\|(\chi_{\mathrm{bet},\Sigma})\eta_{k_{1},k_{2}}\|_{L^{2}W^{3-\frac{1}{q_{\delta}},q_{\delta}}}^{2}+\|(\chi_{\mathrm{bet}})p_{k_{1},k_{2}}\|_{L^{2}W^{1,q_{\delta}}}^{2} )  \notag \\
\lesssim& (\|(\chi_{\mathrm{cl}}+\chi_{\mathrm{bet}})v_{k_{1},k_{2}}\|_{L^{2}W^{2,q_{\delta}}}^{2}+\|(\chi_{\mathrm{cl},\Sigma}+\chi_{\mathrm{bet},\Sigma})\eta_{k_{1},k_{2}}\|_{L^{2}W^{3-\frac{1}{q_{\delta}},q_{\delta}}}^{2}+\|(\chi_{\mathrm{cl}}+\chi_{\mathrm{bet}})p_{k_{1},k_{2}}\|_{L^{2}W^{1,q_{\delta}}}^{2} )  \notag \\
+& (\|(\chi_{\mathrm{cl}})v_{k_{1},k_{2}}\|_{L^{2}W^{2,q_{\delta}}}^{2}+\|(\chi_{\mathrm{cl},\Sigma})\eta_{k_{1},k_{2}}\|_{L^{2}W^{3-\frac{1}{q_{\delta}},q_{\delta}}}^{2}+\|(\chi_{\mathrm{cl}})p_{k_{1},k_{2}}\|_{L^{2}W^{1,q_{\delta}}}^{2} ) \\
\lesssim&\mathcal{E}\mathcal{D}+\frac{C(\varepsilon)}{\varepsilon} \sum_{0\leq j_{1}\leq k_{1}}\sum_{0\leq j_{2}\leq k_{2}}\mathscr{E}_{j_{1},j_{2},\mathrm{energy},\mathrm{cl}}+\sum_{j=0}^{k_{2}+1}\mathscr{E}_{k_{1},j,\mathrm{energy},\mathrm{cl}}+\sum_{j=0}^{k_{2}}\mathscr{E}_{k_{1}+1,j,\mathrm{energy},cl}+\frac{C(\varepsilon)}{\varepsilon} \mathcal{D}_{\mathrm{int}}.
\end{align}
This finishes the proof.
\end{proof}

The same proof applies at every level $j_{1}+j_{2}\leq m$. Summing over $j_{1}\leq k_{1}$ and $j_{2}\leq k_{2}$ gives Theorem~\ref{thm:intermediate_estimate}, since $\mathscr{E}_{k_{1},k_{2},\mathrm{el},\mathrm{bet}}$ is the corresponding sum of the left-hand sides, in the convention of Subsection~\ref{subsec:smallness}.

\textbf{Remark: Since the artificial domain is chosen to be independent of $\varepsilon$, the part depending on $\varepsilon $ has already been written explicitly.}

\section{Contact Line Energy estimate}

\label{sec:contact_energy}

\subsection{Contact line local energy estimate}

\label{subsec:cl_local_energy}

In this subsection, we derive the energy estimate for $\theta$ derivatives near the contact line in the local charts. We fix an integrability parameter independently of the cutoff width: \(\nu:=\frac{\varepsilon_-}{2},\,\allowbreak q_\nu:=\frac{2}{2-\nu},\,\allowbreak r_\nu:=\frac{2}{1-\nu}.\) Thus \(\alpha<\nu<\varepsilon_-\). In the remainder of this proof, the norm exponents previously denoted by \(q_\varepsilon\) and \(q_{r_\varepsilon}\) are \(q_\nu\) and \(r_\nu\), respectively. All Sobolev and trace constants at these exponents are fixed before the tube width \(\varepsilon\) is chosen. For a transverse support \(U_\varepsilon\) and its surface interval \(I_\varepsilon\), we have
\begin{equation}
|U_\varepsilon| \le C\varepsilon^2,\quad |I_\varepsilon| \le C\varepsilon,\quad \|f\|_{L^2_\theta L^{q_\nu}(U_\varepsilon)} \le C\varepsilon^{\varepsilon_--\nu} \|f\|_{L^2_\theta L^{q_-}(U_\varepsilon)},\quad \|g\|_{L^2_\theta L^{q_\nu}(I_\varepsilon)} \le C\varepsilon^{(\varepsilon_--\nu)/2} \|g\|_{L^2_\theta L^{q_-}(I_\varepsilon)}.
\end{equation}
We decrease the common exponent \(\iota>0\), if needed, so that \(0<\iota\le(\varepsilon_--\nu)/2\). The auxiliary bulk products use \(\frac{2}{r_\nu}=1-\nu<\frac{1}{q_\nu} =1-\frac{\nu}{2},\,\allowbreak r_\nu<\frac{2}{1-\varepsilon_-}.\) Hence their H\"older and Sobolev constants also remain independent of the cutoff width.
The surface dual exponent is also fixed: \(q_\nu'=\frac{q_\nu}{q_\nu-1}=\frac{2}{\nu}.\)

We now write down the equation system, where, in the original chart, we still use the $\theta$-derivative to denote the spatial derivative in the local charts. Recall the system \eqref{eq:localized-after-flattening}. We change the cutoff function to $\chi_{\mathrm{bet}+}$; for simplicity, in the remaining part of this section, we use $(v,p,\eta)$ to denote $(\chi_{\mathrm{bet}+}v,\chi_{\mathrm{bet}+}p,\chi_{\mathrm{bet}+,\Sigma}\eta)$. Then we have
\begin{equation}
\left\{ \begin{alignedat}{2} &\partial_t v + \operatorname{div}_{Q^{T}\mathcal{A}_{\Phi}} \bigl( S_{Q^{T}\mathcal{A}_{\Phi}}(p,v)+R(v) \bigr) + \bigl( S_{Q^{T}\mathcal{A}_{\Phi}}(p,v)+R(v) \bigr) \partial_\theta Q^{T}\mathcal{A}_{\Phi}\,e_\theta \\ &+ \mathfrak{C} \bigl( S_{Q^{T}\mathcal{A}_{\Phi}}(p,v)+R(v) \bigr) Q^{T}\mathcal{A}_{\Phi}e_\theta = \widetilde F_{bet+}^1 &&\qquad \text{in }\tilde{\Omega}_{bet+}, \\
&\operatorname{div}_{Q^{T}\mathcal{A}_{\Phi}}v+\partial_{\theta}Qv\cdot \mathcal{A}_{\Phi}e_{\theta} = \widetilde F_{bet+}^2 &&\qquad \text{in }\tilde{\Omega}_{bet+}, \\ &\bigl( S_{Q^{T}\mathcal{A}_{\Phi}}(p,v)+R(v) \bigr) \frac{Q^{T}\mathcal{A}_{\Phi}n_{0,+}}{|Q^{T}\mathcal{A}_{\Phi}n_{0,+}|} = \left( g\eta - \sigma\Delta_{\Sigma,\Phi}\eta \right) \frac{Q^{T}\mathcal{A}_{\Phi}n_{0,+}}{|Q^{T}\mathcal{A}_{\Phi}n_{0,+}|} + \widetilde F_{bet+}^3 &&\qquad \text{on }\Gamma_{+,bet+}, \\
&\frac{1}{|N_0|\circ\Phi_\Sigma} \partial_t\eta - v\cdot \frac{Q^{T}\mathcal{A}_{\Phi}n_{0,+}}{|Q^{T}\mathcal{A}_{\Phi}n_{0,+}|} = \widetilde F_{bet+}^4 &&\qquad \text{on }\Gamma_{+,bet+}, \\ &v\cdot \frac{Q^{T}\mathcal{A}_{\Phi}n_{0,-}}{|Q^{T}\mathcal{A}_{\Phi}n_{0,-}|} = 0 &&\qquad \text{on }\Gamma_{-,bet+}, \\ &U_{Q^{T}\mathcal{A}_{\Phi},-} \bigl( S_{Q^{T}\mathcal{A}_{\Phi}}(p,v)+R(v) \bigr) \frac{Q^{T}\mathcal{A}_{\Phi}n_{0,-}}{|Q^{T}\mathcal{A}_{\Phi}n_{0,-}|} - \beta_0 U_{Q^{T}\mathcal{A}_{\Phi},-}v = \widetilde F_{bet+}^5 &&\qquad \text{on }\Gamma_{-,bet+}, \\ &\sigma\mathcal B_{\mathrm{bet},\Phi}\eta =-\widetilde F_{bet+}^6-\kappa\partial_{t}\eta &&\qquad \text{on }\Gamma_0, \\ &v = 0 &&\qquad \text{on }\Gamma_{\mathrm{art},bet+}, \\ &\eta = 0 &&\qquad \text{on }\Lambda_{\mathrm{art},bet+}, \end{alignedat} \right.
\end{equation}
where the forcing terms on the right-hand side contain three parts: nonlinear terms, geometry difference terms and cutoff commutators.

We now include the geometry difference terms in the main term on the left-hand side and compute it again in the local charts. For the normal vector on free surface denoted by $\mathcal{N}=\mathcal{N}_{0}-(\nabla_{x,y} \eta,0 )$, we have \(\mathcal{N}_{0}\circ \Phi= |\mathcal{N}_{0}|\frac{\mathcal{A}_{\Phi}n_{0,+}}{|\mathcal{A}_{\Phi}n_{0,+}|},\) from equation \eqref{eq:geometric_normal_vector}, and \((\nabla_{x,y}\eta)\circ \Phi=A_X^{-T}\nabla_{\theta,r}\mathfrak{h}_{\eta}.\) Therefore, in the local chart,
\begin{align}
\mathcal{N}\circ\Phi=\mathcal{N}_{0}\circ\Phi-\binom{A_X^{-T}\nabla_{\theta,r}\mathfrak{h}_{\eta}}{0},
\end{align}
where we abused the notation by using $\Phi$ to denote the restriction of $\Phi$ on the free surface $\Sigma$. Moreover, we still use $\mathcal{N}_{0}$ to denote $|\mathcal{N}_{0}|\frac{\mathcal{A}_{\Phi}n_{0,+}}{|\mathcal{A}_{\Phi}n_{0,+}|}$ which depends on the steady state $\zeta$ in the local charts.

Define $\mathcal{A}_{*}=(\mathcal{A}_{\eta}\circ\Phi)\mathcal{A}_{\Phi}$, and $B_{*}=Q^{T}\mathcal{A}_{*}$, and
\begin{align}
\hat{\mathfrak{h}}_{\eta}:=\begin{pmatrix} A_X^{-T}\nabla_{\theta,r}\mathfrak{h}_{\eta}\\
0
\end{pmatrix}.
\end{align}  We obtain
\begin{equation}  \label{eq:localized-after-flattening_0} \left\{ \begin{alignedat}{2} &\partial_t v + \operatorname{div}_{B_{*}} \bigl( S_{B_{*}}(p,v)+R_{\eta}(v) \bigr) -\bigl( S_{B_{*}}(p,v)+R_{\eta}(v) \bigr) \mathfrak{C}B_{*}\,e_\theta \\ &+ \mathfrak{C} \bigl( S_{B_{*}}(p,v)+R_{\eta}(v) \bigr) B_{*}e_\theta = \widetilde F_{bet+}^1 &&\qquad \text{in }\tilde{\Omega}_{bet+}, \\
&\operatorname{div}_{B_{*}}v+\partial_{\theta}Qv\cdot \mathcal{A}_{\eta}\mathcal{A}_{\Phi}e_{\theta} = \widetilde F_{bet+}^2 &&\qquad \text{in }\tilde{\Omega}_{bet+}, \\ &\bigl( S_{B_{*}}(p,v)+R_{\eta}(v) \bigr) (Q^{T}\mathcal{N}_{0}-Q^{T}\hat{\mathfrak{h}_{\eta}}) = \left( g\eta - \sigma\Delta_{\Phi}\eta \right) (Q^{T}\mathcal{N}_{0}-Q^{T}\hat{\mathfrak{h}}_{\eta}) + \widetilde F_{bet+}^3 &&\qquad \text{on }\Gamma_{+,bet+}, \\ & \partial_t\eta - v\cdot (Q^{T}\mathcal{N}_{0}-Q^{T}\hat{\mathfrak{h}}_{\eta}) = \widetilde F_{bet+}^4 &&\qquad \text{on }\Gamma_{+,bet+}, \\ &v\cdot \frac{Q^{T}\mathcal{A}_{\Phi}n_{0,-}}{|Q^{T}\mathcal{A}_{\Phi}n_{0,-}|} = 0 &&\qquad \text{on }\Gamma_{-,bet+}, \\ &U_{Q^{T}\mathcal{A}_{\Phi},-} \bigl( S_{B_{*}}(p,v)+R_{\eta}(v) \bigr) \frac{Q^{T}\mathcal{A}_{\Phi}n_{0,-}}{|Q^{T}\mathcal{A}_{\Phi}n_{0,-}|} - \beta_0 U_{Q^{T}\mathcal{A}_{\Phi},-}v = \widetilde F_{bet+}^5 &&\qquad \text{on }\Gamma_{-,bet+}, \\ &\sigma\mathcal B_{\mathrm{cl},\Phi}\eta_{\mathrm{cl}+} = -\widetilde F_{bet+}^6-\kappa\partial_{t}\eta &&\qquad \text{on }\Gamma_0, \\ &v_{cl} = 0 &&\qquad \text{on }\Gamma_{\mathrm{art},bet}, \\ &\eta_{cl} = 0 &&\qquad \text{on }\Lambda_{\mathrm{art},bet}, \end{alignedat} \right.
\end{equation}
where
\begin{align}
-R_{\eta} = \mu(Q^T \mathcal{A}_{\eta} \mathcal{A}_{\Phi}e_{\theta}) \otimes (Q^T \partial_\theta Q v) +\mu (Q^T \partial_\theta Q v) \otimes (Q^T \mathcal{A}_{\eta}\mathcal{A}_{\Phi} e_{\theta}).
\end{align}
We apply $\partial_{t}^{k_{1}}\partial_{\theta}^{k_{2}}$ to the equation \eqref{eq:localized-after-flattening_0} and argue by recurrence on $k_{2}$. From the previous subsection, the energy estimate holds when $k_{2}=0$. Then suppose that it holds for $k_{2}\leq a$, and use $\tilde{v}$ to denote the velocity field instead of $v$ in system \eqref{eq:localized-after-flattening_0}. For $k_{2}=a+1$, pulling the system back to the original coordinates and setting $\p_{\theta}^{k_{2}}\p_{t}^{k_{1}}v=Q\p_{\theta}^{k_{2}}\p_{t}^{k_{1}}\tilde{v}\circ (\Phi)^{-1}$, when $k_{2}\geq 1, k_{1}\geq 0$ and $k_{1}+k_{2}\leq m+1$, we have
\begin{equation}  \label{eq:localized-after-flattening_original_chart} \left\{ \begin{alignedat}{2} &\partial_t \partial_{\theta}^{k_{2}}\partial_{t}^{k_{1}}v + \operatorname{div}_{\mathcal{A}_{\eta}} \bigl( S_{\mathcal{A}_{\eta}}(\partial_{\theta}^{k_{2}}\partial_{t}^{k_{1}}p,\partial_{\theta}^{k_{2}}\partial_{t}^{k_{1}}v)\bigr) = \widetilde F_{bet+}^{1,k_{1},k_{2}} &&\qquad \text{in }{\Omega}_{bet+}, \\
&\operatorname{div}_{\mathcal{A}_{\eta}}\partial_{\theta}^{k_{2}}\partial_{t}^{k_{1}}v = \widetilde F_{bet+}^{2,k_{1},k_{2}} &&\qquad \text{in }{\Omega}_{bet+}, \\ &\bigl( S_{\mathcal{A}_{\eta}}(\partial_{\theta}^{k_{2}}\partial_{t}^{k_{1}}p,\partial_{\theta}^{k_{2}}\partial_{t}^{k_{1}}v)) \bigr) (\mathcal{N}) \\
&\qquad= \left( g\partial_{\theta}^{k_{2}}\partial_{t}^{k_{1}}\eta - \sigma\operatorname{div}_{\Sigma}\mathcal L_{\Sigma,0}\partial_{\theta}^{k_{2}}\partial_{t}^{k_{1}}\eta -\sigma \operatorname{div}_{\Sigma}(F_{bet+}^{5,k_{1},k_{2}}))\right) \mathcal{N} + \widetilde F_{bet+}^{3,k_{1},k_{2}} &&\qquad \text{on }\Gamma_{+,bet+}, \\ & \partial_t\partial_{\theta}^{k_{2}}\partial_{t}^{k_{1}}\eta - \partial_{\theta}^{k_{2}}\partial_{t}^{k_{1}}v\cdot \mathcal{N} = \widetilde F_{bet+}^{4,k_{1},k_{2}} &&\qquad \text{on }\Gamma_{+,bet+}, \\
&\big[(Q D_\theta^{k_2}\partial_t^{k_1}\tilde v)\circ\Phi^{-1}\big]\cdot n_0=0 &&\qquad \text{on }\Gamma_{-,bet+}, \\ &U_{0} \bigl( S_{\mathcal{A}_{\eta}}(\partial_{\theta}^{k_{2}}\partial_{t}^{k_{1}}p_{cl},\partial_{\theta}^{k_{2}}\partial_{t}^{k_{1}}v) \bigr) n_{0} - \beta_0 U_{0}\partial_{\theta}^{k_{2}}\partial_{t}^{k_{1}}v = \widetilde F_{bet+}^{6,k_{1},k_{2}} &&\qquad \text{on }\Gamma_{-,bet+}, \\ & \sigma\mathcal B_{\mathrm{cl},\Phi}\eta_{k_{1},k_{2}}+\sigma \tilde{F}_{\mathrm{bet}+}^{5,k_{1},k_{2}}\cdot \nu = -\widetilde F_{bet+}^{7,k_{1},k_{2}}-\kappa\partial_{t}\eta_{k_{1},k_{2}} &&\qquad \text{on }\Gamma_0, \\ &v_{cl} = 0 &&\qquad \text{on }\Gamma_{\mathrm{art},bet+}, \\ &\xi_{cl} = 0 &&\qquad \text{on }\Lambda_{\mathrm{art},bet+}, \end{alignedat} \right.
\end{equation}
The momentum forcing below is written in physical coordinates, with separate derivatives for scalars and vectors. In formulas written in the local frame, we use $v$ to denote $\tilde v$. For the momentum forcing, keep the cutoff explicit. In the following formulas $v,p,F_0^1$ are the \emph{unlocalized} fields of \eqref{eq:moving-domain-elliptic-system}, $\chi=\chi_{\mathrm{bet}+}$ and $P=\partial_t^{k_1}\partial_\theta^{k_2}$ in chart coordinates. In physical coordinates write \(P_s p=\big[P(p\circ\Phi)\big]\circ\Phi^{-1},\,\allowbreak P_v v=\big[Q P\big(Q^T(v\circ\Phi)\big)\big]\circ\Phi^{-1}.\) Use $P_v$ also on vector forces. Put $u_P=P_vv$, $p_P=P_sp$, $\mathcal A=\mathcal A_\eta$, and $L_{\mathcal A}(p,v)=\operatorname{div}_{\mathcal A}S_{\mathcal A}(p,v)$. The differentiated fields on the left of the localized momentum equation are $\chi u_P$ and $\chi p_P$, since $\partial_t\chi=\partial_\theta\chi=0$. The momentum forcing in physical coordinates is therefore
\begin{equation}
\begin{aligned}
& \widetilde F^{1,k_1,k_2}_{\mathrm{bet}+} =\chi P_v F_0^1-\mathfrak O_{1,k_1,k_2} + S_{\mathcal A}(p_P,u_P)\nabla_{\mathcal A}\chi -\mu\operatorname{div}_{\mathcal A}\mathscr C_{\chi,\mathcal A}(u_P), \\
& \mathfrak O_{1,k_1,k_2} :=\chi\big(P_v L_{\mathcal A}(p,v)-L_{\mathcal A}(p_P,u_P)\big),\quad \mathscr C_{\chi,\mathcal A}(u_P) :=u_P\otimes\nabla_{\mathcal A}\chi+ \nabla_{\mathcal A}\chi\otimes u_P .
\end{aligned}
\end{equation}
Indeed $\partial_t$ commutes with $P_v$, and $S_{\mathcal A}(\chi p_P,\chi u_P)= \chi S_{\mathcal A}(p_P,u_P)-\mu\mathscr C_{\chi,\mathcal A}(u_P)$. Taking $\operatorname{div}_{\mathcal A}$ gives both cutoff terms above, as in \eqref{eq:Fj1-before-flattening}. In the rotating frame, $\mathfrak O_{1,k_1,k_2}$ is computed by commuting $P$ with the full rotated stress divergence, including the derivatives of $Q$.

Write $\chi_\Sigma=\chi|_\Sigma$, $R=\mathcal R(\nabla_x\eta)$, $R_P=P_sR$ and $h_P=P_s\eta$, where $P_s$ acts componentwise on the vector valued $\mathcal{R}$. Set $Lf=\operatorname{div}_x(a_\zeta\nabla_x f)$ and $H=g\eta-\sigma L\eta-\sigma\operatorname{div}_x R$. On the wall write $T(p,v)=U_0 S_{\mathcal A}(p,v)n_0-\beta_0 U_0v$. Applying $P_v$ to the capillary equation $S_{\mathcal A}(p,v)\mathcal N=H\mathcal N$ on $\Sigma$ and to $T(p,v)=0$ on $\Sigma_s$, then subtracting the principal terms, gives the unlocalized remainders
\begin{equation*}
\begin{aligned}
& F_P^3 :={}  S_{\mathcal A}(p_P,u_P)\mathcal N -P_v\big(S_{\mathcal A}(p,v)\mathcal N\big) +P_v(H\mathcal N)-(P_sH)\mathcal N -\sigma\big([P_s,L]\eta+[P_s,\operatorname{div}_x]R\big)\mathcal N, \\
& F_P^6 :={}  T(p_P,u_P)-P_v\big(T(p,v)\big),
\end{aligned}
\end{equation*}
where $[A,B]$ denotes the commutator between $A$ and $B$ written as $[A,B]=AB-BA$. The stress product rule above, together with $\operatorname{div}_x(\chi_\Sigma R_P) =\chi_\Sigma\operatorname{div}_x R_P+\nabla_x\chi_\Sigma\cdot R_P$, now gives
\begin{equation*}
\widetilde F_{\mathrm{bet}+}^{3,k_1,k_2} ={} \chi_\Sigma F_P^3-\mu\mathscr C_{\chi,\mathcal A}(u_P)\mathcal N +\sigma\big([L,\chi_\Sigma]h_P+\nabla_x\chi_\Sigma\cdot R_P\big) \mathcal N,\quad \widetilde F_{\mathrm{bet}+}^{6,k_1,k_2} ={} \chi F_P^6-\mu U_0\mathscr C_{\chi,\mathcal A}(u_P)n_0 .
\end{equation*}
Since $a_\zeta$ is symmetric, the surface commutator is \([L,\chi_\Sigma]h_P =2(a_\zeta\nabla_x\chi_\Sigma)\cdot\nabla_x h_P +h_P\operatorname{div}_x(a_\zeta\nabla_x\chi_\Sigma).\)

For the other forcing terms, we write them down in the local charts so that the commutators between derivatives can be expressed explicitly.
\begin{equation}
\begin{aligned}
& \tilde{F}_{\mathrm{bet}+}^{2,k_{1},k_{2}}:= -\chi_{\mathrm{bet}+}\partial_{t}^{k_{1}}\partial_{\theta}^{k_{2}}(\partial_{\theta}Q\tilde{v}\cdot \mathcal{A}_{\eta }\mathcal{A}_{\Phi}e_{\theta}) -\chi_{\mathrm{bet}+}\partial_{t}^{k_{1}}\partial_{\theta}^{k_{2}}(\operatorname{div}_{B_{*}}\tilde{v})+\operatorname{div}_{B_{*}}(\chi_{\mathrm{bet}+}\partial_{t}^{k_{1}}\partial_{\theta}^{k_{2}}\tilde{v}) \\
& +(\partial_{\theta}Q\partial_{t}^{k_{1}}\partial_{\theta}^{k_{2}}(\chi_{\mathrm{bet}+}\tilde{v})\cdot \mathcal{A}_{\eta }\mathcal{A}_{\Phi}e_{\theta}).
\end{aligned}
\end{equation}
\begin{equation}
\begin{aligned}
& \tilde{F}_{\mathrm{bet}+}^{4,k_{1},k_{2}}:=\chi_{\mathrm{bet}+,\Sigma} \sum_{k_{1,1}+k_{1,2}=k_{1}}\sum_{\substack{k_{2,1}+k_{2,2}=k_{2},\\k_{2,1}+k_{1,1}\neq k_{2}+k_{1}}}\binom{k_{1}}{k_{1,1}}\binom{k_{2}}{k_{2,1}}\partial_{t}^{k_{1,1}}\partial_{\theta}^{k_{2,1}}\tilde{v}\partial_{t}^{k_{1,2}}\partial_{\theta}^{k_{2,2}}(Q^{T}\mathcal{N}_{0}-Q^{T} \hat{\mathfrak{h}}_{\eta}), \\
& \tilde{F}_{\mathrm{bet}+}^{5,k_{1},k_{2}}:= \chi_{\mathrm{bet}+,\Sigma}\partial_{t}^{k_{1}}\partial_{\theta}^{k_{2}}\mathcal{R}(\nabla\eta), \\
& \tilde{F}_{\mathrm{bet}+}^{7,k_{1},k_{2}}:= \partial_{t}^{k_{1}}\partial_{\theta}^{k_{2}} \hat{\mathcal{W}}(\partial_{t}\eta)+\sigma [\partial_\theta^{k_2},\mathcal B_{\mathrm{cl},\Phi}]\partial_t^{k_1}\eta +\sigma\sum_{i=1}^{k_{2}}\binom{k_{2}}{i} \partial_{t}^{k_{1}}\partial_{\theta}^{k_{2}-i}\mathcal{R}(\nabla\eta) \cdot\partial_{\theta}^{i}\nu .
\end{aligned}
\end{equation}
In the contact law $\nu$ denotes the fixed outer normal of $\partial D$, and $\chi_{\mathrm{bet}+}=1$ near $\Gamma_0$. In the last sum at least one derivative falls on $\nu$. The term with $i=0$ stays on the left in \eqref{eq:localized-after-flattening_original_chart}. The commutator is $[P,B]f=P(Bf)-B(Pf)$, with the entire fixed conormal operator written in the surface chart. For the wall test write $V=Q^T(v\circ\Phi)=\tilde v$, with $v$ the localized physical velocity in \eqref{eq:localized-after-flattening_original_chart}. Put $B=Q^T\mathcal A_\Phi$. For a vector $f$ in local components, we define the covariant derivative as follows
\begin{equation}
D_\theta f=B^{-T}\partial_\theta(B^Tf) =\partial_\theta f+[(\partial_\theta B)B^{-1}]^Tf .
\end{equation}
Let $n_-$ be the fixed normal of the wedge wall. On this wall $n_0\circ\Phi=QBn_-/|Bn_-|$, so $v\cdot n_0=0$ implies $V\cdot Bn_-=0$. Since $B^TD_\theta f=\partial_\theta(B^Tf)$ and $n_-$ is constant, for each integer $b\geq 0$ we have
\begin{equation}
D_\theta^b f=B^{-T}\partial_\theta^b(B^Tf),\qquad (D_\theta^b f)\cdot Bn_-=\partial_\theta^b(f\cdot Bn_-).
\end{equation}
The matrix $B$ is fixed in time. Thus $Q D_\theta^{k_2}\partial_t^{k_1}V$, composed with $\Phi^{-1}$, has zero physical normal trace. At the edge, $B=I$ for every $\theta$, and the first derivative reduces to \(QD_\theta V=\partial_\theta(v\circ\Phi)-(\partial_\theta Q)Q^T(v\circ\Phi).\)

Before the main theorem, we first state two lemmas. At the purely angular top level $(k_{1},k_{2})=(0,m+1)$, the pressure and the free surface are estimated without the loss $\alpha$ in Theorem~\ref{thm:top_level} of Subsection~\ref{subsec:top_level}, which is used in Steps 6, 7 and 10 of the proof of Theorem~\ref{thm:contact_line_energy}.

\begin{lemma}
{\label{lem:local_threehalf}} Suppose that for every $\omega\in H^1(\Omega)$ such that $\operatorname{div}_{\mathcal A_\eta}\omega=0$, $\omega$ is tangent to $\Sigma_s$, and $\omega\cdot\mathcal N\in H^1(D)$, we have
\begin{equation}
\begin{aligned}
&g\int_{\Sigma}\eta_{\mathrm{bet}+,k_{1},k_{2}}(\omega\cdot \mathcal{N})+\sigma\int_{\Sigma}a_{\zeta}\nabla_{\Sigma}\eta_{\mathrm{bet}+,k_{1},k_{2}}\cdot \nabla (\omega\cdot \mathcal{N})\\ &\quad+\sigma\int_{\Sigma}\tilde{F}_{\mathrm{bet}+}^{5,k_{1},k_{2}}\cdot \nabla_{\Sigma}(\omega\cdot \mathcal{N})+[\kappa\partial_{t}\eta_{k_{1},k_{2}}+\tilde{F}_{\mathrm{bet}+}^{7,k_{1},k_{2}},\omega \cdot \mathcal{N}]_{\Gamma_{0}} \lesssim C\|\omega\|_{H^{1}}.
\end{aligned}
\end{equation}
Then,
\begin{align}
\|\eta_{\mathrm{bet}+,k_{1},k_{2}}\|^{2}_{H^{\frac{3}{2}-\alpha}}\lesssim C^{2}+\|\eta_{\mathrm{bet}+,k_{1},k_{2}}\|^{2}_{L^{2}}+{\mathcal{E}\mathcal{D}}+[\partial_{t}\eta_{k_{1},k_{2}}]_{\Gamma_{0}}^{2}+\varepsilon^{2\iota}\sum_{j_{1}\leq k_{1},j_{2}\leq k_{2}-1}\|\eta_{j_{1},j_{2}}\|_{L^{2}W^{3-\frac{1}{q_{\delta}},q_{\delta}}}^{2}.
\end{align}
\end{lemma}

\begin{proof}
We choose the test function $\omega$ as follows.  For the localized height $\partial_{t}^{k_{1}}\partial_{\theta}^{k_{2}}\eta_{\mathrm{bet}+}$ where $\partial_t^{k_{1}}\partial_{\theta}^{k_{2}}\eta_{\mathrm{bet}+}$ denotes the push-forward of $\partial_t^{k_{1}}\partial_{\theta}^{k_{2}}(\chi_{\mathrm{bet}+,\Sigma}\mathfrak{h}_{\eta})$ in the original coordinates, put $f=\eta_{\mathrm{bet}+,k_1,k_2}$ and \(f^\circ=f-|D|^{-1}\int_D f\,dx,\,\allowbreak G_N=D_N^s f^\circ, \,\allowbreak s=1-2\alpha .\) The constant eigenfunction is retained by the projection, so $\int_D G_N=0$. Take the unique mean-zero solution on the full equilibrium vessel of
\begin{equation}
\begin{cases}
-\Delta\psi=0 & \text{in }\Omega, \\
\partial_{\nu}\psi =\dfrac{G_{N}}{|\mathcal{N}_0|} & \text{on }\Sigma, \\
\partial_{\nu}\psi=0 & \text{on }\Sigma_s.\end{cases}\end{equation}

Proposition~\ref{thm:whole_space_threehalf} and its proof, applied to the zero-mean flux $G_N$, give $\|\psi\|_{H^2(\Omega)} \lesssim\|G_N\|_{H^{1/2}(D)}$. Thus $\omega=M\nabla\psi$ is a global admissible test with $\omega\cdot\mathcal N=G_N$ and $\|\omega\|_{H^1(\Omega)}\lesssim\|G_N\|_{H^{1/2}(D)}$. In \eqref{eq:local_weak_formulation} the localized fields and their cutoff residuals are taken on $\Omega$, extended by zero outside the fixed coordinate neighborhood. They vanish near its artificial outer boundary. We may therefore use the restriction of this global test on their support; we do not truncate $G_N$ or impose an artificial boundary condition on $\psi$. In particular,
\begin{equation}
|(\kappa \partial_{t}^{k_{1}}\partial_{\theta}^{k_{2}}\eta_t+F^7,G_N)_{\Gamma_0}| \le C_\alpha\|\kappa \partial_{t}^{k_{1}}\partial_{\theta}^{k_{2}}\eta_t+F^7\|_{L^2(\Gamma_0)} \|G_N\|_{H^{1/2+\alpha}(D)},
\end{equation}
and the curvature flux $F^5 = \widetilde{F}^{5,k_1,k_2}_{\mathrm{bet}+}$ is estimated as follows. Put $\lambda = 1/2 - \alpha$, $P_s$ as above and $h_P = P_s\eta$, so $f = \chi_\Sigma h_P$. For $k_2 > 0$ the leading term in $F^5$ is $c\,\chi_\Sigma\nabla_x h_P$ where $c = \mathcal{R}_z(x,\nabla_x\eta)$. The previous estimate for lower order $\eta$ or $\mathfrak{h}_{\eta}$ estimates give $c, \partial_\theta c \in L^\infty$ and $\partial_r c \in L^\infty_\theta L^{b}_r$ for $1 < b < \min(2, 1/(1-\varepsilon_-))$, with bound $C(\varepsilon)\sqrt{\mathcal{E}}$. After reflection and a smooth cutoff, Lemma \ref{lem:multiplier_estimate} and interpolation with $L^2$ give this leading term multiplier bound on $H^\lambda$. In $\chi_\Sigma\nabla_x h_P = \nabla_x f - (\nabla_x\chi_\Sigma)h_P$, the contact line dissipation \eqref{def:energy_cl} bounds the first term in $H^\lambda$ and the $k = k_1$ row of the interior dissipation functional \ref{def:interior_dissipation} bounds the second in $H^1$.

In split products the smooth coefficients have the same multiplier bound with size $C(\varepsilon)$. Except for the main term, at most m-th order derivatives are applied to one multiplier in this nonlinear term. Then the fractional difference quotient puts one gradient factor in $H^\lambda$ and the others in $L^\infty$, with smooth localization on the overlap. Use dissipation for a factor of order $m$ in either norm, or for one factor if all orders are at most $m-1$. The other gradient factors use energy. The quadratic vanishing of $\mathcal{R}$ supplies at least two. For $k_2 = 0$, apply the same multiplier to the global pure-time gradient $\nabla_x\partial_t^{k_1}\eta$ and then multiply by $\chi_\Sigma$. The global $H^{3/2-\alpha}$ height row in $\mathcal{D}_{\mathrm{enhanced}}$ also covers $k_1 = m+1$. Thus $\|F^5\|_{H^\lambda(D)} \le C(\varepsilon)\sqrt{\mathcal{E}\mathcal{D}}$. In this lemma the nonlinear errors carry constants $C(\varepsilon)$. Applying (1.4) gives
\begin{align}
\sigma\int_{\Sigma}\tilde{F}_{\mathrm{bet}+}^{5,k_{1},k_{2}}\cdot \nabla_{\Sigma}(\omega\cdot \mathcal{N})\lesssim \sqrt{\mathcal{E}}\sqrt{\mathcal{D}}\|G_{N}\|_{H^{1-\frac{s}{2}}}.
\end{align}
Put $X_N=\|P_N f^\circ\|_{\mathcal H_K^{1+s/2}}$. Since $G_N$ has zero average, the capillary form of $f$ against $G_N$ equals the form of $f^\circ$ against $G_N$. By \eqref{functional:exact_spectral_identity} its value is $X_N^2$, and $\|G_N\|_{H^{1-s/2}(D)}\lesssim X_N$, uniformly in $N$. Now $1-s/2=1/2+\alpha$, so the same bound controls the trace of $G_N$ on $\partial D$ and the $H^1(\Omega)$ norm of the test field. Using the preceding estimates and the estimate for $\tilde F^{7,k_1,k_2}_{\mathrm{bet}+}$ in the proof of Theorem~\ref{thm:contact_line_energy}, we obtain
\begin{equation}
X_N^2\lesssim \left(C+\sqrt{\mathcal E\mathcal D} +\|\partial_t\eta_{k_1,k_2}\|_{L^2(\Gamma_0)} +\varepsilon^\iota\sum_{j_1\le k_1,\,j_2< k_2} \|\eta_{j_1,j_2}\|_{L^2W^{3-1/q_\delta,q_\delta}}\right)X_N .
\end{equation}
Taking the supremum over $N$ bounds $\|f^\circ\|_{H^{3/2-\alpha}(D)}$. The constant $f-f^\circ$ has $H^{3/2-\alpha}(D)$ norm at most $C_D\|f\|_{L^2(D)}$. Adding this term and squaring gives the claimed estimate for $f$. For $k_2\ge1$ the remaining term is of lower order, since $\|\eta_{\mathrm{bet}+,k_{1},k_{2}}\|_{L^{2}}\lesssim\|\eta_{\mathrm{bet}+,k_{1},k_{2}-1}\|_{H^{1}}$; this is why it does not appear in Step~9 of the proof of Theorem~\ref{thm:contact_line_energy}.  This finishes the proof.
\end{proof}

We now state the second lemma, which gives the regularity of the pressure.

\begin{lemma}
\label{lem:enhance_pressure}
Choose the connected collar $U=\Omega_{\mathrm{bet}+}$ of Subsection~\ref{subsec:artificial} with $\Phi^{-1}(U)=\mathbb T\times\mathcal W_{\mathrm{bet}+}$, where the connected Lipschitz cross section is fixed independently of the cutoff width. For $a=k_1$ and $b=k_2\ge1$, let $p_{a,b}$ be the localized pressure derivative. Suppose that at the time considered
\begin{equation}
\left|\int_U J p_{a,b}\func{div}_{\mathcal A_\eta}w\right|\le C\|w\|_{H^1(U)} \qquad\text{for every }w\in H^1_0(U;\mathbb R^3).
\end{equation}
Then $\|p_{a,b}\|_{L^2(U)}\lesssim C+\|p_{a,b-1}\|_{L^2(U)}$.
\end{lemma}

\begin{proof}
Write $P=p_{a,b}$, $P_J=\int_U JP/\int_U J$, and $z=\mathcal B_U(J(P-P_J))$, where $\mathcal B_U:L^2_0(U)\to H^1_0(U;\mathbb R^3)$ is a bounded right inverse of divergence on this fixed collar \cite[Section~III.3]{Galdi}. The datum has integral zero. Set $w=Mz$ with $M=J^{-1}\mathcal A_\eta^{-T}$. Multiplication by $M$ preserves $H^1_0(U)$: $M$ is bounded, and near the edge $\nabla M\in L^\infty_\theta L^{r_-}_{r,h}$ with $r_-=2/(1-\varepsilon_-)$, while $z\in L^2_\theta L^{2/\varepsilon_-}_{r,h}$ by the fixed slice embedding. Away from the edge use the interior coefficient bound. These multiplier constants are bounded after fixing the width and choosing the small-energy threshold, since the nonconstant deformation terms are bounded by $C(\varepsilon)\sqrt{\mathcal E}$ as in the proof of Proposition~\ref{prop:B_d}. Thus $\|w\|_{H^1(U)}\lesssim\|P-P_J\|_{L^2(U)}$. By the Piola identity,
\begin{equation}
J\func{div}_{\mathcal A_\eta}w=\func{div}z=J(P-P_J), \qquad \int_U JP\func{div}_{\mathcal A_\eta}w=\int_U J|P-P_J|^2 .
\end{equation}
The hypothesis and the positive lower bound on $J$ give $\|P-P_J\|_{L^2(U)}\lesssim C$.

On $\widetilde U=\Phi^{-1}(U)$ put $\rho=(J\circ\Phi)J_\Phi$. Since $\chi_{\mathrm{bet}+}$ is independent of $\theta$, $P\circ\Phi=\partial_\theta(p_{a,b-1}\circ\Phi)$. Periodicity yields
\begin{equation}\label{eq:pressure_mean} P_J=-\frac{\int_{\widetilde U}(\partial_\theta\rho)(p_{a,b-1}\circ\Phi)} {\int_{\widetilde U}\rho}, \qquad |P_J|\lesssim\|p_{a,b-1}\|_{L^2(U)} .
\end{equation}
Indeed $\partial_\theta\rho=J_\Phi\partial_\theta(J\circ\Phi) +(J\circ\Phi)\partial_\theta J_\Phi$ is bounded by the fixed chart and coefficient estimates, and the weighted volume is bounded below. The triangle inequality now proves the claim.
\end{proof}

\smallskip\noindent\textbf{Proof of Theorem~\ref{thm:contact_line_energy}.}

\textbf{The interior dissipation.} The cutoff commutators below (Steps~1--7), the pressure bound of Step~10 and Theorem~\ref{thm:top_level} produce the interior dissipation $\mathcal{D}_{\mathrm{int}}$ with constants $C(\varepsilon)$. These terms are not absorbed; they are bounded by the initial energy. Indeed, Theorem~\ref{thm:interior_space_final} is proved in Section~\ref{sec:interior} from Theorems~\ref{thm:whole-space-energy} and~\ref{thm:whole_space_estimate_final}. None of these uses the estimates of this section; the contact-line norms enter them only through nonlinear terms bounded by $\sqrt{\mathcal{E}}\mathcal{D}$. Moreover $\mathcal{E}^{\frac32}(0)\leq\mathcal{E}(0)$, since $\mathcal{E}\leq1$. Hence
\begin{align}
\int_{0}^{t}\mathcal{D}_{\mathrm{int}}\leq C(\varepsilon)\Big(\mathcal{E}(0)+\mathcal{E}^{\frac{3}{2}}(t)+\int_{0}^{t}\sqrt{\mathcal{E}}\mathcal{D}\Big), \label{eq:Dint_E0}
\end{align}
which is of the form of the right-hand side of Theorem~\ref{thm:contact_line_energy}.

\textbf{Energy-dissipation relation.} Similarly to the pure $t$ energy estimate of Theorem \ref{thm:whole-space-energy}, we derive the weak formulation of the equation system for any test function $w$ such that $w\in H^{1}(\Omega)$, $w\cdot n=0$ on $\Sigma_{s}$, and $w\cdot \mathcal{N}\in H^{1}(\Sigma)$:
\begin{align}
{\label{eq:local_weak_formulation}} &\int_\Omega J\,\partial_t v_{k_{1},k_{2}}\cdot w + \frac{\mu}{2}\int_\Omega J\,\mathbb{D}_{\mathcal{A}_\eta}v_{k_{1},k_{2}} : \mathbb{D}_{\mathcal{A}_\eta}w - \int_\Omega Jp_{k_{1},k_{2}}\func{div}_{\mathcal{A}_\eta}w  \notag \\
&\quad+ \beta_0\int_{\Sigma_s} J(U_0v_{k_{1},k_{2}})\cdot(U_0w) + g\int_\Sigma\eta_{k_{1},k_{2}}(w\cdot\mathcal{N}) + \sigma\left\langle a_{\zeta} \nabla_{\Sigma}\eta_{k_{1},k_{2}},\nabla_{\Sigma}(w\cdot\mathcal{N}) \right\rangle_\Sigma  \notag \\
&\qquad= \int_\Omega J\tilde{F}_{bet+}^{1,k_{1},k_{2}}\cdot w - \left\langle \tilde{F}_{bet+}^{3,k_{1},k_{2}},w\right\rangle_\Sigma - \left\langle J\tilde{F}_{bet+}^{6,k_{1},k_{2}},U_0w\right\rangle_{\Sigma_s} - \sigma\left\langle \tilde{F}_{bet+}^{5,k_{1},k_{2}},\nabla_{\Sigma}(w\cdot\mathcal{N}) \right\rangle_\Sigma  \notag \\
&\qquad\quad-\langle \kappa \partial_{t}\eta_{k_{1},k_{2}}+\tilde{F}_{\mathrm{bet}+}^{7,k_{1},k_{2}},w\cdot \mathcal{N}\rangle_{\Gamma_{0}} .
\end{align}
Here both wall terms use the surface Jacobian: $J_s=J$ on the vertical wall in this collar.

Choose $w$ by $w\circ\Phi=QD_\theta^{k_2}\partial_t^{k_1}V$. By the wall identity above this is a wall-tangent test. The corresponding ordinary derivative satisfies $v_{k_1,k_2}\circ\Phi=Q\partial_\theta^{k_2}\partial_t^{k_1}V$. Writing $w=v_{k_1,k_2}+\omega$, the product rule gives
\begin{equation}
\omega\circ\Phi=Q\sum_{j=0}^{k_2-1}\binom{k_2}{j} B^{-T}(\partial_\theta^{k_2-j}B^T) \partial_t^{k_1}\partial_\theta^jV .
\end{equation}
All coefficients here depend only on the fixed chart, and every velocity in $\omega$ has fewer than $k_2$ angular derivatives. With this test we obtain
\begin{align}
&\frac{d}{dt} \Bigg[ \frac12\int_\Omega J|v_{k_{1},k_{2}}|^2 + \frac g2\|\eta_{k_{1},k_{2}}\|_{L^2(\Sigma)}^2 + \frac{\sigma}{2}\int_\Sigma a_\zeta\nabla_\Sigma\eta_{k_{1},k_{2}}\cdot\nabla_\Sigma\eta_{k_{1},k_{2}}+\int_{\Omega}J_{\eta}\omega v_{k_{1},k_{2}}\Bigg] \notag \\ &+ \frac{\mu}{2}\int_\Omega J|\mathbb{D}_{\mathcal{A}_\eta}v_{k_{1},k_{2}}|^2 + \beta_0\int_{\Sigma_s}J_s|U_0v_{k_{1},k_{2}}|^2 + \kappa\|\partial_t\eta_{k_{1},k_{2}}\|_{L^2(\Gamma_0)}^2 \notag \\ &= \frac12\int_\Omega(\partial_tJ)|v_{k_{1},k_{2}}|^2 + \int_\Omega Jp_{k_{1},k_{2}}\tilde{F}_{\mathrm{bet}+}^{2,k_{1},k_{2}} + \int_\Omega J\tilde{F}_{\mathrm{bet}+}^{1,k_{1},k_{2}}\cdot v_{k_{1},k_{2}} - \left\langle \tilde{F}_{\mathrm{bet}+}^{3,k_{1},k_{2}},v_{k_{1},k_{2}}\right\rangle_\Sigma - \left\langle J_s\tilde{F}_{\mathrm{bet}+}^{6,k_{1},k_{2}},U_0v_{k_{1},k_{2}}\right\rangle_{\Sigma_s} \notag \\ &+ g\int_\Sigma\eta_{k_{1},k_{2}}\tilde{F}_{\mathrm{cl}+}^{4,k_{1},k_{2}} + \sigma\int_\Sigma a_\zeta\nabla_\Sigma\eta_{k_{1},k_{2}}\cdot\nabla_\Sigma \tilde{F}_{\mathrm{bet}+}^{4,k_{1},k_{2}} + \kappa\left\langle \partial_t\eta_{k_{1},k_{2}},\tilde{F}_{\mathrm{bet}+}^{4,k_{1},k_{2}} \right\rangle_{\Gamma_0} \notag \\ &- \left\langle \tilde{F}_{\mathrm{bet}+}^{7,k_{1},k_{2}},\partial_t\eta_{k_{1},k_{2}}-\tilde{F}_{\mathrm{bet}+}^{4,k_{1},k_{2}} \right\rangle_{\Gamma_0} - \sigma\int_\Sigma \tilde{F}_{\mathrm{bet}+}^{5,k_{1},k_{2}}\cdot\nabla_\Sigma\partial_t\eta_{k_{1},k_{2}} + \sigma\int_\Sigma \tilde{F}_{\mathrm{bet}+}^5\cdot\nabla_\Sigma \tilde{F}_{\mathrm{bet}+}^{4,j} \notag \\ &+\int_{\Omega}J_{\eta}\partial_{t}\omega v_{k_{1},k_{2}}+\int_{\Omega}\partial_{t}J_{\eta}\omega v_{k_{1},k_{2}}+\int_{\Omega}Jp_{k_{1},k_{2}}\func{div}_{\mathcal{A}}\omega+\int_{\Omega}J\tilde{F}_{\mathrm{bet}+}^{1,k_{1},k_{2}}\cdot \omega-\mu/2 \int_{\Omega}J\mathbb{D}_{\mathcal{A}_{\eta}}v_{k_{1},k_{2}}:\mathbb{D}_{\mathcal{A}_{\eta}}\omega \notag \\ &-\langle \tilde{F}_{\mathrm{bet}+}^{3,k_{1},k_{2}},\omega\rangle_{\Sigma}-\langle J_{s}\tilde{F}_{\mathrm{bet}+}^{6,k_{1},k_{2}}+J_{s}\beta_{0}U_{0}v_{k_{1},k_{2}},U_{0}\omega \rangle_{\Sigma_{s}}-\langle \kappa\partial_{t}\eta_{k_{1},k_{2}}+\tilde{F}_{\mathrm{bet}+}^{7,k_{1},k_{2}},\omega\cdot \mathcal{N}\rangle_{\Gamma_{0}} \notag \\ &-g\int_{\Sigma}\eta_{k_{1},k_{2}}(\omega\cdot \mathcal{N})-\sigma\int_{\Sigma}(a_{\eta}\nabla_{\Sigma}\eta_{k_{1},k_{2}}+\tilde{F}_{\mathrm{bet}+}^{5,k_{1},k_{2}})\cdot \nabla_{\Sigma}(\omega\cdot \mathcal{N}). \label{eq:local-domain-energy-1-2}
\end{align}

It suffices to estimate the forcing terms on the right-hand side. The other terms can be bounded by the same process as in Theorem \ref{thm:whole-space-energy}. From the derivation, compared to $F_{0}^{i,j}$ in Theorem \ref{thm:whole-space-energy}, the forcing terms on the right-hand side include three parts: the commutators between spatial derivatives and cutoff functions, the remaining terms obtained from taking $\partial_{t},\partial_{\theta}$ on the left-hand side of equation \eqref{eq:localized-after-flattening_original_chart}, and the nonlinear terms of the original equation system. For the nonlinear terms of the original equation system, we use the same estimates as in the energy estimates of \cite{GT}, changing the $\partial_{t}^{k} $ derivative to $\partial_{t}^{k_{1}}\partial_{\theta}^{k_{2}}$.

For the remaining terms obtained from taking $\partial_{||}$ on the left-hand side of equation \eqref{eq:localized-after-flattening_original_chart}, we still apply the two-dimensional estimates of \cite{GT}. The only difference is that, when there is only a $\theta$ derivative on the geometric terms $\mathcal{A}_{\Phi},Q,\mathcal{N},\mathcal{A}_{\eta}$, we have
\begin{align}
|\partial_{\theta}^{j}(\mathcal{A}_{\Phi}+Q+\mathcal{A}_{\eta})|\lesssim \sum_{k=0}^{j} (1+|\partial_{\theta}^{k}\zeta|+|\partial_{\theta}^{k}\nabla\zeta|+|\partial_{\theta}^{k}\bar{\eta}|+|\partial_{\theta}^{k}\nabla \bar{\eta}|),
\end{align}
and \(|\partial_{\theta}^{j}\mathcal{N}|\lesssim (|\partial_{\theta}^{k}\nabla \zeta|+|\partial_{\theta}^{k}\nabla{\eta}|),\) which depends on the steady state. Compared to the pure temporal derivative, we have to change norms related to $\partial_{||}^{j}\eta$ on the right-hand side of the inequality to $(1+|\partial_{||}^{j}\eta|_{\mathcal{W}})$ for any required Sobolev norms. Moreover, the exponent $q_{\delta}$ is not tied to the energy estimate itself; as in \cite{GT}, any $q_{\delta}\leq q_{-}$ would do, and we fixed $q_{\delta}=q_{-}$ in Subsection~\ref{subsec:constants}.

Here and below, $U$ is the transverse support in the bulk or the corresponding interval on the surface. In either case, H\"older's inequality gives \(\|f\|_{L^2_\theta L^{q_\nu}(U)} \le |U|^{(\varepsilon_--\nu)/2}\|f\|_{L^2_\theta L^{q_-}(U)}.\) We use the common exponent fixed at the beginning of this subsection: \(0<\iota\le(\varepsilon_--\nu)/2 .\) Since the transverse area is $O(\varepsilon^2)$ and the interval has length $O(\varepsilon)$, both estimates imply \(\|f\|_{L^2_\theta L^{q_\nu}(U)} \lesssim\varepsilon^\iota\|f\|_{L^2_\theta L^{q_-}(U)} .\)

Suppose that $k_{1}+k_{2}=m+1$ and $k_{2}>0$. Then we have
\begin{align}
\sum_{j=1}^{k_{2}}\int_{0}^{t}(\mathscr{E}_{k_{1}+j,k_{2}-j,\mathrm{energy},\mathrm{cl}})+\frac{1}{\varepsilon^{\iota}}\sum_{j=0}^{k_{2}-1}\|\partial_{t}^{k_{1}}\partial_{\theta}^{j}v\|^{2}_{L^{2}}(t)\lesssim C(\varepsilon)((\mathcal{E}(0)+\mathcal{E}^{\frac{3}{2}}(t)+\int_0^{t}\sqrt{\mathcal{E}}\mathcal{D})).
\end{align}
When $k_{2}=0$, the estimate degenerates to the estimates in Theorem~\ref{thm:whole-space-energy}. We now use the assumption to derive the estimate for $\mathscr{E}_{k_{1}+j,k_{2}-j,\mathrm{energy},\mathrm{cl}}$.

\textbf{Step 0: Estimates for the nonlinear terms in the forcing.}

Apart from lower-order terms involving $\theta$-derivatives of the steady state, the nonlinear terms can be estimated as in the proof of Theorem~\ref{thm:whole-space-energy}, using the same integration-by-parts arguments for the corresponding terms if it is necessary. When at least one angular derivative falls on the fixed coefficients, the perturbation carries at most $m$ derivatives. The extra angular derivative on the fixed coefficients may destroy the symmetric structure that integration by parts relies on. However, it also lowers the number of derivatives on the perturbation, so integration by parts is no longer needed. For example, consider the following term in $\tilde{F}_{\mathrm{bet}+}^{1,k_{1},k_{2}}$, written in local charts:
\begin{equation}
\operatorname{div}_{\p_{\theta}Q^{T}\p_{\theta}^{k_{2}-1}\p_{t}^{k_{1}}(\mathcal{A}_{\eta}\circ\Phi)\mathcal{A}_{\Phi}}\nabla_{B_{*}}\tilde{v}.
\end{equation}
We bound it in $L^{2}L^{q_{\nu}}$ in the same way as
\begin{equation}
\|\operatorname{div}_{\p_{t}^{k_{1}+k_{2}-1}(\mathcal{A}_{\eta})}\nabla_{\mathcal{A}_{\eta}}v\|_{L^{2}L^{q_{\nu}}}.
\end{equation}
Since $k_{1}+k_{2}-1=m$, no integration by parts is needed for this term. All other nonlinear terms can be bounded in the same way as the analogous terms in the proof of Theorem~\ref{thm:whole-space-energy}. In the remaining estimates, we therefore focus on the linear terms generated when angular derivatives fall on the fixed coefficients.

\textbf{Step 1. Estimate for the term including $\tilde{F}_{\mathrm{bet}+}^{2,k_{1},k_{2}}$} We have \(\int_{\Omega}Jp_{k_{1},k_{2}}\tilde{F}_{\mathrm{bet}+}^{2,k_{1},k_{2}}\lesssim \|p_{k_{1},k_{2}}\|_{L^{2}}\|\tilde{{F}}_{\mathrm{bet}+}^{2,k_{1},k_{2}}\|_{L^{2}}.\) By definition and the estimate for nonlinear $F_{0}^{2,j}$ in the proof of Theorem \ref{thm:whole-space-energy}, the $L^{2}$ norm of nonlinear terms included in $F_{\mathrm{ber}+}^{2,k_{1},k_{2}}$ are bounded by $\sqrt{\mathcal{E}\mathcal{D}}$. In the commutator between time angular derivative and  divergence one derivative hits a coefficient. The preceding product bounds give the nonlinear $L^{2}$ bound. A differentiated fixed coefficient leaves at most $k_{2}-1$ angular derivatives on the fluid gradient. A derivative of $\chi$ is supported in the interior patch. The thin-support bound above therefore gives
\begin{equation}
\begin{aligned}
& \|F_{\mathrm{bet}+}^{2,k_{1},k_{2}}\|_{L^{2}} \lesssim \sum_{j=0}^{k_{2}-1}\|\p_{t}^{k_{1}}\p_{\theta}^{j}((\chi_{\mathrm{bet}+})v)\|_{H^{1}}+\|\p_{t}^{k_{1}}\p_{t}^{j}(\chi_{\mathrm{bet}})v)\|_{L^{2}}+C(\varepsilon)\|\p_{t}^{k_{1}}((\chi_{\mathrm{int}})v)\|_{H^{k_{2}}}+\sqrt{\mathcal{E}\mathcal{D}} \\
& \lesssim  \sum_{j=0}^{k_{2}-1}\|\p_{t}^{k_{1}}\p_{\theta}^{j}((\chi_{\mathrm{bet}+})v)\|_{H^{1}}+\varepsilon^{\iota}\|\p_{t}^{k_{1}}\p_{t}^{j}(\chi_{\mathrm{bet}})v)\|_{H^{1}}+C(\varepsilon)\|\p_{t}^{k_{1}}((\chi_{\mathrm{int}})v)\|_{W^{2+k_{2}-1,p}}+\sqrt{\mathcal{E}\mathcal{D}} \\
& \lesssim \varepsilon^{\iota}\|v_{k_{1},k_{2}}\|_{H^{1}}+\sum_{j=0}^{k_{2}-1}\|v_{k_{1},j}\|_{H^{1}}+\sqrt{\mathcal{E}\mathcal{D}}+C(\varepsilon)\sqrt{\mathcal{D}_{\mathrm{int}}},
\end{aligned}
\end{equation}
for some small positive number $\iota$. The linear part of $\tilde{F}_{\mathrm{bet}+}^{2,k_{1},k_{2}}$ that contains $v_{k_{1},k_{2}}$ is at the top level; its factor $\varepsilon^{\iota}$ comes from the Poincar\'e inequality along rays from the contact line, $\|v_{k_{1},k_{2}}\|_{L^{2}}\lesssim\varepsilon\|\nabla v_{k_{1},k_{2}}\|_{L^{2}}$ on the $\varepsilon$-thin support of $\chi_{\mathrm{bet}+}$, and not from the smallness of $\eta$. Moreover, $s_{i}$ and $d_{i}$ are defined as in Lemma \ref{thm:contact_line_forcing_estimate}, depending on the sum of the orders of $\partial_{t}$ and $\partial_{\theta}$ in the norms. Therefore
\begin{align}
\int_{\Omega}Jp_{k_{1},k_{2}}\tilde{F}_{\mathrm{bet}+}^{2,k_{1},k_{2}}\lesssim \|p_{k_{1},k_{2}}\|_{L^{2}}(\varepsilon^{\iota}\|v_{k_{1},k_{2}}\|_{H^{1}}+\sum_{j=0}^{k_{2}-1}\|v_{k_{1},j}\|_{H^{1}}+\sqrt{\mathcal{E}\mathcal{D}}+C(\varepsilon)\sqrt{D}_{\mathrm{int}}).
\end{align}

\textbf{Step 2: estimate for the term including $\tilde{F}_{\mathrm{bet}+}^{1,k_{1},k_{2}}$}

The forcing terms are written above in the original coordinates. By the equivalence of the Sobolev norms under the fixed chart $\Phi$, it suffices to bound them in the chart variables.

\textbf{Term $\partial_{t}^{k_{1}}\partial_{\theta}^{k_{2}}F_{0}^{1}$}

Note, for smooth $T$ and $w$, with $T=pI-\mu\mathbb{D}v$,
\begin{align}
\int_\Omega J\,\func{div}_{\mathcal{A}}T\cdot w &=-\int_\Omega JT:\nabla_{\mathcal{A}}w +\int_D(T\mathcal{N})\cdot w\,dx +\int_{\Sigma_s}J_s(Tn_0)\cdot w\,dS .
\end{align}
By an estimate similar to that in \eqref{eq:highest_IBP}, the boundary terms can be controlled by the main term in $\Omega$. Therefore, the term with $\partial_{t}^{k_{1}}\partial_{\theta}^{k_{2}}F_{0}^{1}$ can be bounded as in the estimate for $F_{0}^{1,j}$ in Theorem \ref{thm:whole-space-energy}, with the obvious replacements. We have
\begin{align}
\|\partial_{t}^{k_{1}}\partial_{\theta}^{k_{2}}F_{0}^{1}\|_{X_{k_{1}+k_{2}}}\lesssim \sqrt{\mathcal{E}\mathcal{D}},\qquad X_{j}:=L^{2}(H^{1})^{*}\ (j\leq m),\quad X_{m+1}:=(H^{1}(\Omega))^{*}.  \label{eq:F01_dual}
\end{align}

At the levels with $k_{1}+k_{2}=m+1$ the norm in \eqref{eq:F01_dual} is that of $(H^{1}(\Omega))^{*}$: in $L^{2}_{\theta}$ of slice duals the bound fails, for the reason given in the proof of Lemma~\ref{thm:contact_line_forcing_estimate}, and against $w\in H^{1}(\Omega)$ it follows by integrating the second-order term by parts in all variables, as there. At $(k_{1},k_{2})=(0,m+1)$ the $\theta\theta$ part of $\mu\func{div}_{\mathcal{A}}(\mathbb{D}_{\partial_{\theta}^{m+1}\mathcal{A}}v)$ would require $\partial_{\theta}^{m+3}\bar{\eta}\in L^{2}$ if it were estimated directly; the integration by parts in all variables is carried out in Step~A2 and Step~B1(c) of the proof of Theorem~\ref{thm:top_level}.

Moreover, for other terms included in $\mathfrak{O}_{1,j_{1},j_{2}}$, we only estimate the linear terms induced by taking pure $\theta$ derivatives of $Q,\mathcal{A}_{\Phi}$ and no time derivative on the geometric part induced by $\eta $. For other nonlinear terms, we bound them by the same computation as for $F_{0}^{1,j}$ in \ref{thm:whole-space-energy}.

\textbf{Term induced by $\func{div}_{B_{*}}(S_{B_{*}}(p,v)+R_{\eta}(v))$}

To specify which terms are estimated here and in Step~3, write $[\,\cdot\,]_{\mathrm{fix}}$ for the part of a commutator in which all the commuted derivatives falling on coefficients are angular derivatives of the fixed chart, frame or equilibrium normal. Ordinary spatial differentiation in the divergence is still carried out, and undifferentiated deformation factors are retained. By the product rule, after the common localization,
\begin{align}
&\big\|[\partial_t^{k_1}\partial_\theta^{k_2} \operatorname{div}_{B_*}(S_{B_*}(p,v)+R_\eta(v)) -\operatorname{div}_{B_*}(S_{B_*}(p_{k_1,k_2},v_{k_1,k_2}) +R_\eta(v_{k_1,k_2}))]_{\mathrm{fix}} \big\|_{L^2_\theta (H^1(\mathcal W))^*} \notag\\
&\lesssim \big\|[\partial_t^{k_1}\partial_\theta^{k_2} \operatorname{div}_{B_*}(S_{B_*}(p,v)+R_\eta(v)) -\operatorname{div}_{B_*}(S_{B_*}(p_{k_1,k_2},v_{k_1,k_2}) +R_\eta(v_{k_1,k_2}))]_{\mathrm{fix}} \big\|_{L^2_\theta L^{q_\nu}(U_\varepsilon)} \notag\\
\lesssim{}& C\sum_{j=0}^{k_2-1}\big( \|\partial_\theta^j\partial_t^{k_1}p\|_{L^2_\theta W^{1,q_\nu}(U_\varepsilon)} +\|\partial_\theta^j\partial_t^{k_1}v\|_{L^2_\theta W^{2,q_\nu}(U_\varepsilon)}\big) \notag\\
&+C\sum_{j=0}^{k_2-1}\big( \|\partial_\theta^{j+1}\partial_t^{k_1}p\|_{L^2_\theta L^{q_\nu}(U_\varepsilon)} +\|\nabla\partial_\theta^{j+1}\partial_t^{k_1}v\|_{L^2_\theta L^{q_\nu}(U_\varepsilon)}\big) \notag\\
&+C\big(1+\|\nabla\mathcal A_\eta\|_{L^\infty_\theta L^{r_\nu}(U_\varepsilon)}\big) \sum_{j=0}^{k_2-1}\big(\|\partial_\theta^j\partial_t^{k_1}p\|_{L^2} +\|\partial_\theta^j\partial_t^{k_1}v\|_{H^1}\big) \notag\\
\lesssim&\varepsilon^{\iota} \sum_{j=0}^{k_{2}-1}(\|\partial_{\theta}^{j}\partial_{t}^{k_{1}}p\|_{L^{2}W^{1,q_{\delta}}}+\|\partial_{\theta}^{j}\partial_{t}^{k_{1}}v\|_{L^{2}W^{2,q_{\delta}}}+\|\partial_{\theta}^{j+1}\partial_{t}^{k_{1}}v\|_{H^{1}}+\|\partial_{\theta}^{j+1}\partial_{t}^{k_{1}}p\|_{L^{2}})+\sum_{j=0}^{k_{2}-1}(\|\partial_{\theta}^{j}\partial_{t}^{k_{1}}p\|_{L^{2}}+\|\partial_{\theta}^{j}\partial_{t}^{k_{1}}v\|_{H^{1}}).
\end{align}

\textbf{Term induced by $(S_{B_{*}}(p,v)+\mathcal{R}_{\eta}(v))\mathfrak{C}B_{*}e_{\theta}$}

We now estimate
\begin{equation}
\begin{aligned}
& \sum_{j_{1}+j_{2}+j_{3}=k_{2},j_{2}<k_{2}}^{k_{2}}\|(S_{\partial_{\theta}^{j_{1}}(B_{*})}(\partial_{\theta}^{j_{2}}\partial_{t}^{k_{1}}p,\partial_{\theta}^{j_{2}}\partial_{t}^{k_{1}}v)+\mathcal{R}_{\eta}(\partial_{t}^{k_{1}}\partial_{\theta}^{j_{2}}v))\partial_{\theta}^{j_{3}}(\mathfrak{C}B_{*})e_{\theta}\|_{L^{2}L^{q_{\nu}}} \\
& \lesssim C\sum_{j=1}^{k_{2}}(\|\partial_{\theta}^{k_{2}-j}\partial_{t}^{k_{1}}p\|_{L^{2}}+\|\partial_{\theta}^{k_{2}-j}\partial_{t}^{k_{1}}v\|_{H^{1}}) \lesssim \sum_{j=1}^{k_{2}}(\|\partial_{\theta}^{k_{2}-j}\partial_{t}^{k_{1}}p\|_{L^{2}}+\|\partial_{\theta}^{k_{2}-j}\partial_{t}^{k_{1}}v\|_{H^{1}}).
\end{aligned}
\end{equation}

\textbf{\ Term induced by $\mathfrak{C} \bigl( S_{B_{*}}(p,v)+R_{\eta}(v) \bigr) B_{*}e_\theta$}

Here we use a computation similar to that for the term induced by $(S_{B_{*}}(p,v)+\mathcal{R}_{\eta}(v))\mathfrak{C}B_{*}e_{\theta}$; it can be bounded by
\begin{align}
\sum_{j=1}^{k_{2}}(\|\partial_{\theta}^{k_{2}-j}\partial_{t}^{k_{1}}p\|_{L^{2}}+\|\partial_{\theta}^{k_{2}-j}\partial_{t}^{k_{1}}v\|_{H^{1}}).
\end{align}

\textbf{The term with the highest-order derivative hitting $\mathcal{A}_{\eta}$}

We now consider the estimate for the following term
\begin{align}
\|\chi_{\mathrm{bet}+}\operatorname{div}_{B_{*}}\mathbb{D}_{Q^{T}\p_{t}^{k_{1}}\p_{\theta}^{k_{2}}(\mathcal{A}_{\eta}\circ\Phi)\mathcal{A}_{\Phi}}v\|_{(H^{1})^{*}}.
\end{align}
Applying a computation similar to that for $\operatorname{div}_{\mathcal{A}_{\eta}}\mathbb{D}_{\p_{t}^{m+1}(\mathcal{A}_{\eta})}v$ established in \eqref{eq:highest_IBP}, we have
\begin{align}
\|\chi_{\mathrm{bet}+}\operatorname{div}_{B_{*}}\mathbb{D}_{Q^{T}\p_{t}^{k_{1}}\p_{\theta}^{k_{2}}(\mathcal{A}_{\eta}\circ\Phi)\mathcal{A}_{\Phi}}v\|_{(H^{1})^{*}}\lesssim\sqrt{\mathcal{E}}\sqrt{\mathcal{D}}.
\end{align}
Therefore, combining all of the estimates above, using Young's inequality, we derive
\begin{equation}
\begin{aligned}
& \int_{\Omega}J\tilde{F}_{\mathrm{bet}+}^{1,k_{1},k_{2}}\cdot v_{k_{1},k_{2}} \\
& \lesssim  (\sqrt{\mathcal{E}\mathcal{D}}+\sqrt{\mathscr{E}_{k_{1},k_{2}-1,\mathrm{el},\mathrm{cl}}}+\sqrt{\mathscr{E}_{k_{1},k_{2}-1,\mathrm{el},\mathrm{bet}}}+\varepsilon^{\iota}\sqrt{\mathscr{E}_{k_{1},k_{2},\mathrm{energy},\mathrm{cl}}}+\sqrt{\mathscr{E}_{k_{1},k_{2}-1,\mathrm{energy},\mathrm{cl}}})\|v_{k_{1},k_{2}}\|_{L^{2}L^{\frac{q_{\nu}}{q_{\nu}-1}}} \\
& +\sqrt{\mathcal{E}\mathcal{D}}\|v_{k_{1},k_{2}}\|_{H^{1}} \lesssim (\sqrt{\mathcal{E}\mathcal{D}}+\sqrt{\mathscr{E}_{k_{1},k_{2}-1,\mathrm{el},\mathrm{cl}}}+\sqrt{\mathscr{E}_{k_{1},k_{2}-1,\mathrm{el},\mathrm{bet}}}+\varepsilon^{\iota}\sqrt{\mathscr{E}_{k_{1},k_{2},\mathrm{energy},\mathrm{cl}}}+\sqrt{\mathscr{E}_{k_{1},k_{2}-1,\mathrm{energy},\mathrm{cl}}}) \\
& \times(\frac{1}{\mathcal{M}^{\ell_{2}}(\varepsilon)}\|v_{k_{1},k_{2}-1}\|_{H^{1}}+\mathcal{M}^{\ell_{1}}(\varepsilon)\|v_{k_{1},k_{2}}\|_{H^{1}}) +\sqrt{\mathcal{E}\mathcal{D}}\|v_{k_{1},k_{2}}\|_{H^{1}} \\
& \lesssim ({\mathcal{E}\mathcal{D}}+\mathcal{M}(\varepsilon)^{\ell_{1}}{\mathscr{E}_{k_{1},k_{2}-1,\mathrm{el},\mathrm{cl}}}+\mathcal{M}(\varepsilon)^{\ell_{1}}{\mathscr{E}_{k_{1},k_{2}-1,\mathrm{el},\mathrm{bet}}}+\varepsilon^{2\iota}{\mathscr{E}_{k_{1},k_{2},\mathrm{energy},\mathrm{cl}}} \\
& +{\mathscr{E}_{k_{1},k_{2}-1,\mathrm{energy},\mathrm{cl}}})+(\frac{1}{\mathcal{M}(\varepsilon)^{\ell_{3}}}\|v_{k_{1},k_{2}-1}\|^{2}_{H^{1}}+\mathcal{M}(\varepsilon)^{\ell_{4}}\|v_{k_{1},k_{2}}\|^{2}_{H^{1}})+\sqrt{\mathcal{E}\mathcal{D}}\|v_{k_{1},k_{2}}\|_{H^{1}},
\end{aligned}
\end{equation}
where $\mathcal{M}(\varepsilon)$ is a very small constant compared to the small constant $\frac{1}{C(\varepsilon)}$, $\ell_{1},\ell_{2},\ell_{3},\ell_{4}$ are positive constants, and we used the bound for the intermediate elliptic estimate.

\textbf{Step 3: Terms included in $\tilde{{F}}_{\mathrm{bet}+}^{3,k_{1},k_{2}}$}

We now estimate the surface terms of $\widetilde F_{\mathrm{bet}+}^{3,k_1,k_2}$ in $L^2_\theta L^{q_\nu}(I_\varepsilon)$. Cutoff terms are treated below after step 7; terms kept with their bulk boundary flux in the weak formulation are estimated together there. On the fixed surface interval $I$ containing $I_\varepsilon$, $H^{1/2}(I)\subset L^{2/\nu}(I)$ gives, for a suitable test function $\phi$,
\begin{align}
|\langle \widetilde F_{\mathrm{bet}+}^{3,k_1,k_2},\phi\rangle_\Sigma| \lesssim\|\widetilde F_{\mathrm{bet}+}^{3,k_1,k_2}\|_{L^2_\theta L^{q_\nu}(I_\varepsilon)} \|\phi\|_{L^2_\theta H^{1/2}(I)}.
\end{align}

\textbf{Term induced by $\big((S_{B_{*}}(p,v)+R_{\eta}(v))(Q^{T}\mathcal{N}_{0}-Q^{T}\hat{\mathfrak{h}}_{\eta})\big)$}

We note that $\hat{\mathfrak{h}}_{\eta}$ contains one spatial derivative of $\mathfrak{h}_{\eta}$ on the surface. For the stress term, the same restriction gives
\begin{align}
&\Big\|\chi_{\mathrm{bet}+,\Sigma}\Big[ \partial_t^{k_1}\partial_\theta^{k_2} \big((S_{B_*}(p,v)+R_\eta(v)) (Q^T\mathcal N_0-Q^{T}\hat{\mathfrak h}_\eta)\big) \notag\\[-2pt]
&\qquad{}- (S_{B_*}(p_{k_1,k_2},v_{k_1,k_2})+R_\eta(v_{k_1,k_2})) (Q^T\mathcal N_0-Q^{T}\hat{\mathfrak h}_\eta) \Big]_{\mathrm{fix}}\Big\|_{L^2_\theta L^{q_\nu}(I_\varepsilon)} \notag\\
\lesssim{}&\varepsilon^\iota\sum_{j=0}^{k_2-1} \big(\|\chi_{\mathrm{bet}+}\partial_\theta^j\partial_t^{k_1}p\|_{L^2_\theta W^{1,q_-}_{r,h}} +\|\chi_{\mathrm{bet}+}\partial_\theta^j\partial_t^{k_1}v\|_{L^2_\theta W^{2,q_-}_{r,h}} +\|\chi_{\mathrm{bet}+}\partial_\theta^{j+1}\partial_t^{k_1}v\|_{H^1}\big).
\end{align}
Here we estimate the terms in which the angular derivatives on the coefficients fall on the fixed chart, frame or equilibrium normal. Their derivatives are bounded independently of $\varepsilon$; the terms with differentiated deformation coefficients belong to the nonlinear estimate which are bounded similar to the pure-time energy estimates by $\sqrt{\mathcal{E}\mathcal{D}}$. On a fixed transverse slice, take the $L^{q_-}(I)$ trace of the pressure and of the velocity and its transverse gradient, then restrict to $I_\varepsilon$. For the angular component of the velocity gradient take the $L^2(I)$ trace of $\partial_\theta^{j+1}\partial_t^{k_1}v$. Restriction to $I_\varepsilon$ gains the power $1/q_\nu-1/2=(1-\nu)/2\ge\iota$ in the latter case and $1/q_\nu-1/q_-=(\varepsilon_--\nu)/2\ge\iota$ in the former. More precisly, we have
\begin{align}
&\Big\|\chi_{\mathrm{bet}+,\Sigma}\Big[ \partial_t^{k_1}\partial_\theta^{k_2} \big((S_{B_*}(p,v)+R_\eta(v)) (Q^T\mathcal N_0-Q^{T}\hat{\mathfrak {h}}_\eta)\big) \notag\\[-2pt]
&\qquad{}- (S_{B_*}(p_{k_1,k_2},v_{k_1,k_2})+R_\eta(v_{k_1,k_2})) (Q^T\mathcal N_0-Q^{T}\hat{\mathfrak h_\eta}) \Big]_{\mathrm{fix}}\Big\|_{L^2_\theta L^{q_\nu}(I_\varepsilon)} \notag\\
\lesssim& \sum_{j=1}^{k_{2}}C_{j,\Phi,Q,\zeta}(\|\partial_{\theta}^{k_{2}-j}\partial_{t}^{k_{1}}(\chi_{\mathrm{bet}+}p)\|_{L^{2}W^{1,q_{\nu}}}+\|\partial_{\theta}^{k_{2}-j}\partial_{t}^{k_{1}}(\chi_{\mathrm{bet}+}v)\|_{L^{2}W^{2,q_{\nu}}}\notag\\
&\quad+\|\partial_{\theta}^{k_{2}-j+1}\partial_{t}^{k_{1}}(\chi_{\mathrm{bet}+}v)\|_{L^{2}W^{1,q_{\nu}}})+C(\varepsilon)\sqrt{\mathcal{D}_{\mathrm{int}}} \notag \\
\lesssim&C_{j,\Phi,Q,\zeta} \varepsilon^{\iota}\sum_{j=1}^{k_{2}}(\|\partial_{\theta}^{k_{2}-j}\partial_{t}^{k_{1}}(\chi_{\mathrm{bet}+}p)\|_{W^{1,q_{\delta}}}+\|\partial_{\theta}^{k_{2}-j}\partial_{t}^{k_{1}}(\chi_{\mathrm{bet}+}v)\|_{W^{2,q_{\delta}}}+\|\partial_{\theta}^{k_{2}-j+1}\partial_{t}^{k_{1}}(\chi_{\mathrm{bet}+}v)\|_{H^{1}}) +C(\varepsilon)\sqrt{\mathcal{D}_{\mathrm{int}}} \notag \\
\lesssim&\varepsilon^{\iota} \sum_{j=1}^{k_{2}}(\|\partial_{\theta}^{k_{2}-j}\partial_{t}^{k_{1}}(\chi_{\mathrm{bet}+}p)\|_{W^{1,q_{\delta}}}+\|\partial_{\theta}^{k_{2}-j}\partial_{t}^{k_{1}}(\chi_{\mathrm{bet}+}v)\|_{W^{2,q_{\delta}}}+\|\partial_{\theta}^{k_{2}-j+1}\partial_{t}^{k_{1}}(\chi_{\mathrm{bet}+}v)\|_{H^{1}}) +C(\varepsilon)\sqrt{\mathcal{D}_{\mathrm{int}}},
\end{align}
where the $\varepsilon^{\iota}$ is used to denote $\varepsilon$ with very low order.

\textbf{Term induced by the curvature tensor on the surface.}

As for the previous term, we have
\begin{align}
&\|\partial_{t}^{j_{1}}\partial_{\theta}^{j_{2}}\big((g\eta-\mathcal{L}_{\Sigma,0}\eta)(Q^{T}\mathcal{N}_{0}-Q^{T}\hat{\mathfrak{h}}_{\eta})\big)-\big((g\partial_{t}^{j_{1}}\partial_{\theta}^{j_{2}}\eta-\mathcal{L}_{\Sigma,0}\partial_{t}^{j_{1}}\partial_{\theta}^{j_{2}}\eta)(Q^{T}\mathcal{N}_{0}-Q^{T}\hat{\mathfrak{h}}_{\eta})\big)\|_{L^2L^{q_{\nu}}}  \notag \\
\lesssim&\sqrt{\mathcal{E}\mathcal{D}}+\varepsilon^{\iota} \sum_{j=1}^{k_{2}}\|\partial_{\theta}^{k_{2}-j}\partial_{t}^{k_{1}}\eta\|_{L^{2}W^{3-\frac{1}{q_{\delta}},q_{\delta}}}+\varepsilon^{\iota}\|\eta_{k_{1},k_{2}}\|_{H^{1}}.
\end{align}

\textbf{Term induced by other nonlinear terms}

For the remaining nonlinear terms, it suffices to estimate $\sigma \chi_{\mathrm{bet}+}[\p_{t}^{k_{1}}\p_{\theta}^{k_{2}},\operatorname{div}_{x}]\mathcal{R}(\nabla \eta)(Q^{T}\mathcal{N}_{0}-Q^{T}\hat{\mathfrak{h}}_{\eta})$ which denote the commutator between angular derivatives and the spatial derivatives that does not appear in the nonlinear estimate in Theorem \ref{thm:whole-space-energy}. For the other terms up to the lower order terms involving $\theta$ derivatives of the steady states, we follow the main idea used in the proof of Theorem\ref{thm:whole-space-energy}. We separate nonlinear forcing terms into $\chi_{cl+}$ and $\chi_{\mathrm{0-,\mathrm{int}}}$ part. Then in the support of $\chi_{cl+}$, we apply the local estimates which hold in the support of $\chi_{\mathrm{bet}+}$. For the interior part, we apply the interior estimate which hold in the support of $\chi_{\mathrm{0,\mathrm{int}}}$. When $\theta$ derivative is applied to the steady state instead of the perturbation function, the lower order perturbation function can be bounded by energy and dissipation bounded by the previous recurrence step. More precisely, by norms the lower order perturbation functions controlled by energy or dissipation.

For the term $\sigma \chi_{\mathrm{bet}+}[\p_{t}^{k_{1}}\p_{\theta}^{k_{2}},\operatorname{div}_{x}]\mathcal{R}(\nabla \eta)(Q^{T}\mathcal{N}_{0}-Q^{T}\hat{\mathfrak{h}}_{\eta})$, for $k_{2}>0$, we have
\begin{equation}
\begin{aligned}
& \|\sigma \chi_{\mathrm{bet}+,\Sigma}[\p_{t}^{k_{1}}\p_{\theta}^{k_{2}},\operatorname{div}_{x}]\mathcal{R}(\nabla \eta)(Q^{T}\mathcal{N}_{0}-Q^{T}\hat{\mathfrak{h}}_{\eta})\|_{L_{\theta}^{2}L^{q_\nu}} \\
& \lesssim (1+\|\eta\|_{W^{1,\infty}})\sum_{j_{2}=0}^{k_{2}-1}\sum_{r=2}^{j_{2}+k_{1}+2}\sum_{\substack{i_{1,1}+\ldots+i_{1,r}=k_{1}\\
i_{2,1}+\ldots +i_{2,r}=j_{2}\\
i_{1,\ell},i_{2,\ell}\geq 0\\
i_{1,\ell}\geq 1~\textbf{or}~i_{2,\ell}\geq 1~\text{for}~\ell>2}} \sum_{\tau\in\{1,\ldots,r\}}\|\chi_{\mathrm{bet}+,\Sigma}\nabla^{2}\p_{t}^{i_{1,\tau}}\p_{\theta}^{i_{2,\tau}}\eta\|_{L^{2}L^{q_\nu}} \\
& \times\prod_{\ell=1,\ell\neq \tau}^{r}\|\p_{t}^{i_{1,\ell}}\p_{\theta}^{i_{2,\ell}}\eta\|_{W^{1,+\infty}} \\
& \lesssim C(\varepsilon)\sum_{j_{2}=0}^{k_{2}-1}\sum_{r=2}^{j_{2}+k_{1}+2}\sum_{\substack{i_{1,1}+\ldots+i_{1,r}=k_{1}\\
i_{2,1}+\ldots +i_{2,r}=j_{2}\\
i_{1,\ell},i_{2,\ell}\geq 0\\
i_{1,\ell}\geq 1~\textbf{or}~i_{2,\ell}\geq 1~\text{for}~\ell>2}} \sum_{\tau\in\{1,\ldots,r\}}\big((\chi_{\mathrm{bet}+,\Sigma}\p_{t}^{i_{1,\tau}}\p_{\theta}^{i_{2,\tau}}\eta)\|_{L^{2}W^{3-\frac{1}{q_{\nu}},q_{\nu}}}+\|\p_{t}^{i_{1,\tau}}(\chi_{\mathrm{int}+,\Sigma}\eta)\|_{W^{3+i_{2,\tau}-\frac{1}{p},p}} \\
& +\|\p_{t}^{i_{1,\tau}}\p_{\theta}^{i_{2,\tau}+1}(\chi_{\mathrm{bet}+,\Sigma}\eta)\|_{L^{2}W^{2-\frac{1}{q_{\nu}},q_{\nu}}}+\|\p_{t}^{i_{1,\tau}}\p_{\theta}^{i_{2,\tau}+2}(\chi_{\mathrm{bet}+,\Sigma}\eta)\|_{L^{2}W^{1-\frac{1}{q_{\nu}},q_{\nu}}}\big) \\
& \times\prod_{\ell=1,\ell\neq \tau}^{r}\|\p_{t}^{i_{1,\ell}}\p_{\theta}^{i_{2,\ell}}\eta\|_{W^{1,+\infty}} \lesssim C(\varepsilon)\sqrt{\mathcal{E}}\sqrt{\mathcal{D}}
\end{aligned}
\end{equation}

To estimate the test in the pairing of this forcing term, let $s=1-\nu/2$. On the fixed transverse slice $\mathcal W$, the trace $H^s(\mathcal W)\to H^{s-1/2}(I)\subset L^{2/\nu}(I)$ is bounded. Interpolating between $L^2(\mathcal W)$ and $H^1(\mathcal W)$ and then using H\"older in $\theta$ gives
\begin{equation}
\|\operatorname{Tr} v_{k_{1},k_{2}}\|_{L^2_\theta L^{2/\nu}(I)}  \lesssim \|v_{k_{1},k_{2}}\|_{L^2}^{1-s}\|v_{k_{1},k_{2}}\|_{H^1}^{s} \lesssim \lambda\|v_{k_{1},k_{2}}\|_{H^1} +\lambda^{-s/(1-s)}\|v_{k_1,k_2-1}\|_{H^1},\qquad 0<\lambda\le1.
\end{equation}
Indeed, $k_2>0$ and one angular derivative of the preceding field controls $\|w\|_{L^2}$ in the fixed chart; differentiating the frame adds only a bounded multiple of that preceding field. H\"older on the interval and Cauchy--Schwarz in $\theta$, followed by Young's inequality with $\lambda$ a small power of $\mathcal M(\varepsilon)$, now give
\begin{equation}
\begin{aligned}
|\langle \tilde{F}_{\mathrm{bet}+}^{3,k_{1},k_{2}},v_{k_1,k_2}\rangle_\Sigma| &\lesssim{}\|\tilde{F}_{\mathrm{bet}+}^{3,k_{1},k_{2}}\|_{L^2_\theta L^{q_\nu}(I_\varepsilon)} \|\operatorname{Tr}v_{k_1,k_2}\|_{L^2_\theta L^{2/\nu}(I)}\\ &\lesssim({\mathcal{E}\mathcal{D}}+\mathcal{M}(\varepsilon)^{\ell_{1}}{\mathscr{E}_{k_{1},k_{2}-1,\mathrm{el},\mathrm{cl}}}+\mathcal{M}(\varepsilon)^{\ell_{2}}{\mathscr{E}_{k_{1},k_{2}-1,\mathrm{el},\mathrm{bet}}}+\varepsilon^{2\iota}{\mathscr{E}_{k_{1},k_{2},\mathrm{energy},\mathrm{cl}}}  \\ &\quad+{\mathscr{E}_{k_{1},k_{2}-1,\mathrm{energy},\mathrm{cl}}})+(\frac{1}{\mathcal{M}(\varepsilon)^{\ell_{3}}}\|v_{k_{1},k_{2}-1}\|^{2}_{H^{1}}+\mathcal{M}(\varepsilon)^{\ell_{4}}\|v_{k_{1},k_{2}}\|^{2}_{H^{1}}).
\end{aligned}
\end{equation}

\textbf{Step 4: Term $\sigma\int_\Sigma F_{\mathrm{bet}+}^{5,k_{1},k_{2}}\cdot\nabla_\Sigma\partial_t\eta_{k_{1},k_{2}}$} Using a computation similar to that in the proof of Theorem \ref{thm:whole-space-energy}, we obtain
\begin{align}
\sigma\int_\Sigma F_{\mathrm{bet}+}^{5,k_{1},k_{2}}\cdot\nabla_\Sigma\partial_t\eta_{k_{1},k_{2}}=&\frac{\sigma}{2}\frac{d}{dt}(\int_{\Sigma}R^{\prime }(\nabla \eta)\nabla_{\Sigma}\eta_{k_{1},k_{2}}\cdot \nabla_{\Sigma}\eta_{k_{1},k_{2}}+2R_{\mathrm{low},\mathrm{bet}+}^{5,k_{1},k_{2}}\cdot \nabla_{\Sigma}\eta_{k_{1},k_{2}}) \\
&-\frac{\sigma}{2}(\int_{\Sigma}\partial_{t}R^{\prime }(\nabla \eta)\nabla_{\Sigma}\eta_{k_{1},k_{2}}\cdot \nabla_{\Sigma}\eta_{k_{1},k_{2}})-\sigma\int_{\Sigma}\partial_{t}R_{\mathrm{low},\mathrm{bet}+}^{5,k_{1},k_{2}}\cdot \nabla_{\Sigma}\eta_{k_{1},k_{2}},
\end{align}
where $\int_{\Sigma}2\partial_{t}R_{\mathrm{low},\mathrm{bet}+}^{5,k_{1},k_{2}}\cdot \nabla_{\Sigma}\eta_{k_{1},k_{2}}\lesssim \sqrt{\mathcal{E}}\mathcal{D}$. Moreover,
\begin{align}
-\frac{\sigma}{2}(\int_{\Sigma}\partial_{t}R^{\prime }(\nabla \eta)\nabla_{\Sigma}\eta_{k_{1},k_{2}}\cdot \nabla_{\Sigma}\eta_{k_{1},k_{2}})\lesssim \sqrt{\mathcal{E}}\mathcal{D},
\end{align}
and
\begin{align}
\int_{\Sigma}R^{\prime }(\nabla \eta)\nabla_{\Sigma}\eta_{k_{1},k_{2}}\cdot \nabla_{\Sigma}\eta_{k_{1},k_{2}}\lesssim\sqrt{\mathcal{E}}^{3}.
\end{align}

\textbf{Step 5: Term including $\tilde{F}_{\mathrm{bet}+}^{6,k_{1},k_{2}}$}

Using a computation similar to that for the term $\tilde{F}_{\mathrm{cl}+}^{3,k_{1},k_{2}}$, we have
\begin{equation}
\begin{aligned}
\|\tilde{F}_{\mathrm{cl}+}^{6,k_{1},k_{2}}\|_{L^{2}(H^{\frac{1}{2}}(\Sigma))^{*}} &\lesssim\varepsilon^{\iota} \sum_{j=1}^{k_{2}}(\|\partial_{\theta}^{k_{2}-j}\partial_{t}^{k_{1}}p\|_{L^{2}W^{1,q_{\delta}}}+\|\partial_{\theta}^{k_{2}-j}\partial_{t}^{k_{1}}v\|_{L^{2}W^{2,q_{\delta}}})+\sqrt{\mathcal{E}\mathcal{D}}.
\end{aligned}
\end{equation}

\textbf{Step 6: Surface term including $\tilde{F}_{\mathrm{bet}+}^{4,k_{1},k_{2}}$}

By decomposition, we may write $\tilde{F}_{\mathrm{cl}+}^{4,k_{1},k_{2}}$ as follows
\begin{align}{\label{eq:decompose_f_6}} \tilde{F}_{\mathrm{bet}+}^{4,k_{1},k_{2}}=-v_{||}\cdot (\nabla_{||} (\chi_{\mathrm{bet}+}\eta_{k_{1},k_{2}}))+v_{||}\cdot \nabla_{||}(\chi_{\mathrm{bet}+})\eta_{k_{1},k_{2}}+\tilde{\mathcal{R}}_{\mathrm{low}+}^{4,k_{1},k_{2}}.
\end{align}
From the estimate in the proof of Theorem \ref{thm:whole-space-energy}, by integration by parts, the first term on the right hand side of equation \eqref{eq:decompose_f_6} pairing with corresponding test function is bounded by (For simplicity, in the following computation, we use $\eta_{k_{1},k_{2}}$ to denote $\chi_{\mathrm{bet}+}\eta_{k_{1},k_{2}}$)
\begin{equation}
\begin{aligned}
& \big|\int_{\Sigma}(\eta_{k_{1},k_{2}}) v_{||}\cdot \nabla_{||} (\eta_{k_{1},k_{2}})+\sigma \int_{\Sigma}a_{\zeta}\nabla_{||}\eta_{k_{1},k_{2}} \cdot \nabla_{||}(v_{||}\cdot \nabla_{||} \eta_{k_{1},k_{2}})\big| \\
& \lesssim  \big|\int_{\Sigma}\eta_{k_{1},k_{2}} v_{||}\cdot \nabla_{||} \eta_{k_{1},k_{2}}\big|+\bigl|\int_{\Sigma}a_{\zeta}\nabla_{||}\eta_{k_{1},k_{2}} \cdot ((\nabla_{||}v_{||})^{T} \nabla_{||} \eta_{k_{1},k_{2}})\big| \\
& +|\int_{\Sigma}\nabla_{||}\eta_{k_{1},k_{2}}\cdot \operatorname{div}_{\Sigma}(v_{||})a_{\zeta}\nabla_{||}\eta_{k_{1},k_{2}}|+|\int_{\Sigma}\nabla_{||}\eta_{k_{1},k_{2}}\cdot ((v_{||}\cdot \nabla_{\Sigma})a_{\zeta}\nabla_{||}\eta_{k_{1},k_{2}})| \\
& \lesssim  \|\eta_{k_{1},k_{2}}\|^{2}_{H^{\frac{3}{2}-\alpha}}(\|\chi_{\mathrm{bet}+}v\|_{L^{\infty}W^{2,q_{\delta}}}+\|\chi_{0,\mathrm{int}}v\|_{W^{3,p}}+\|\chi_{\mathrm{bet}+}\p_{\theta}v\|_{L^{\infty}H^{1}}) \lesssim  \sqrt{\mathcal{E}}\mathcal{D}.
\end{aligned}
\end{equation}
where we used the fact that $v\cdot n_{0}=0$ on the fixed boundary.

For the second term in the decomposition of $\tilde{F}_{\mathrm{bet}+}^{4,k_{1},k_{2}}$, we have
\begin{equation}
\begin{aligned}
& \big|\int_{\Sigma}(\chi_{\mathrm{bet}+,\Sigma}\eta_{k_{1},k_{2}}) v_{||}\cdot \nabla_{||}\chi_{\mathrm{bet}+,\Sigma} (\eta_{k_{1},k_{2}})+\sigma \int_{\Sigma}a_{\zeta}\nabla_{||}(\chi_{\mathrm{bet}+,\Sigma}\eta_{k_{1},k_{2}}) \cdot \nabla_{||}(v_{||}\cdot \nabla_{||}\chi_{\mathrm{bet+},\Sigma} \eta_{k_{1},k_{2}})\big| \\
& \lesssim  C(\varepsilon)\|\chi_{\mathrm{bet}+} \eta_{k_{1},k_{2}}\|_{H^{\frac{3}{2}-\alpha}}(\|\chi_{\mathrm{bet}+}v\|_{L^{\infty}W^{2,q_{\delta}}}+\|\chi_{0,\mathrm{int}}v\|_{W^{3,p}}+\|\chi_{\mathrm{bet}+}\p_{\theta}v\|_{L^{\infty}H^{1}}) \\
& \times(\|\chi_{\mathrm{bet}+,\Sigma}\eta_{k_{1},k_{2}}\|_{H^{\frac{3}{2}-\alpha}}+\|\chi_{0,\mathrm{int},\Sigma}\p_{t}^{k_{1}}\eta\|_{W^{1+k_{2},p}}) \lesssim  \sqrt{\mathcal{E}}\mathcal{D},
\end{aligned}
\end{equation}
where we used the fact that $1+k_{2}=1+m+1-k_{1}=m+2-k_{1}<m+3-\frac{1}{p}-k_{1}$ which is the index for elliptic norm controlled by $\sqrt{\mathcal{D}}$.

For the third term in the decomposition of $\tilde{F}_{\mathrm{bet}+}^{4,k_{1},k_{2}}$, it contains nonlinear terms and linear terms induced by taking $\theta$-derivative on the normal vector $\mathcal{N}_{0}$. For the nonlinear part, we apply the proof of Theorem \ref{thm:whole-space-energy} up to the lower order terms involving $\theta$ derivatives of the steady states. The difference here is that the $\theta$ derivative may apply to the fixed geometric terms induced by the steady state. Those terms do not appear in the estimate established in \ref{thm:whole-space-energy}. However, those terms can be similarly bounded by a multiplication of the norms of steady state and nonlinear estimates for perturbation functions with lower order angular derivative. We use this logic in estimates of all nonlinear forcing terms and only give the explicit bound for the linear forcing term and nonlinear forcing terms that does not have corresponding similar term in the proof of Theorem \ref{thm:whole-space-energy}.

We aim to bound \(\|\tilde{R}_{\mathrm{low}+}^{4,k_{1},k_{2}}\|_{W^{1,\frac{1}{1-\alpha}}}.\) It suffices to consider the linear terms induced  from $v_{||}\cdot \mathcal{N}_{0}$ when $\theta$-derivative hits on $\mathcal{N}_{0}$ as follows
\begin{align}
\sum_{j=1}^{k_{2}}\|\partial_{\theta}^{k_{2}-j}\partial_{t}^{k_{1}}v_{||}\cdot \nabla_{||}\partial_{\theta}^{j}\zeta\|_{L^{2}W^{1,\frac{1}{1-\alpha}}}\lesssim \sum_{j=1}^{k_{2}}\varepsilon^{\iota}\|\partial_{\theta}^{k_{2}-j}\p_{t}^{k_{1}}v\|_{L^{2}W^{2,q_{\delta}}}.
\end{align}
Therefore,
\begin{equation}
\int_{\Sigma}\eta_{k_{1},k_{2}}\tilde{F}_{\mathrm{cl}+}^{4,k_{1},k_{2}}+\sigma \int_{\Sigma}\nabla_{||}\eta_{k_{1},k_{2}} \cdot \nabla_{||}(\tilde{F}_{\mathrm{cl}+}^{4,k_{1},k_{2}}) \lesssim\sqrt{\mathcal{E}}\mathcal{D} +\|\eta_{k_{1},k_{2}}\|_{L^{2}H^{\frac{3}{2}-\alpha}}\sum_{j=1}^{k_{2}}\varepsilon^{\iota}\|\partial_{\theta}^{k_{2}-j}\p_{t}^{k_{1}}v\|_{L^{2}W^{2,q_{\delta}}}.
\end{equation}

In the last two displays the norm $L^{2}W^{1,\frac{1}{1-\alpha}}$ is the slice norm in $r$, with $\partial_{r}$ only. The $r$-part of the pairing is estimated slice by slice with the one-dimensional embedding $H^{\frac12-\alpha}(\mathbb{R})\subset L^{1/\alpha}(\mathbb{R})$:
\begin{equation*}
\Big|\int_{\Sigma}\partial_{r}\eta_{k_{1},k_{2}}\,\partial_{r}\tilde{\mathcal{R}}^{4}_{\mathrm{low}+}\Big| \lesssim\|\partial_{r}\eta_{k_{1},k_{2}}\|_{L^{2}_{\theta}H^{\frac12-\alpha}_{r}}\|\partial_{r}\tilde{\mathcal{R}}^{4}_{\mathrm{low}+}\|_{L^{2}_{\theta}L^{\frac{1}{1-\alpha}}_{r}} \lesssim\|\eta_{k_{1},k_{2}}\|_{H^{\frac32-\alpha}}\|\tilde{\mathcal{R}}^{4}_{\mathrm{low}+}\|_{L^{2}W^{1,\frac{1}{1-\alpha}}}.
\end{equation*}
On the two-dimensional surface $\Sigma$ the embedding would fail, since only $W^{1,\frac{1}{1-\alpha}}\subset H^{2\alpha}$ holds there, so the $\theta$-derivatives are treated separately. The $\theta\theta$ part of $\sigma\int_{\Sigma}a_{\zeta}\nabla\eta_{k_{1},k_{2}}\cdot\nabla\tilde{F}_{\mathrm{bet}+}^{4,k_{1},k_{2}}$ is bounded by $\|\partial_{\theta}\eta_{k_{1},k_{2}}\|_{L^{2}(\Sigma)}\|\partial_{\theta}\tilde{\mathcal{R}}^{4}_{\mathrm{low}+}\|_{L^{2}(\Sigma)}$. Here $\partial_{\theta}\tilde{\mathcal{R}}^{4}_{\mathrm{low}+}$ contains the top-level trace $v_{k_{1},k_{2}}|_{\Sigma}$ times $\nabla\partial_{\theta}\zeta$, and traces of lower $\theta$-levels. The field $v_{k_{1},k_{2}}$ vanishes on the outer edge $\{|r|+|h|=\varepsilon\}$ of the tube. The Poincar\'e inequality along rays and the scaled trace inequality give \(\|v_{k_{1},k_{2}}\|_{L^{2}(\Sigma\cap\func{supp}\chi_{\mathrm{bet}+})}\lesssim\varepsilon^{1/2}\|v_{k_{1},k_{2}}\|_{H^{1}}.\) So the factor $\sum_{j}\varepsilon^{\iota}\|\partial_{\theta}^{k_{2}-j}v\|_{L^{2}W^{2,q_{\delta}}}$ in the last display is to be read with $\varepsilon^{1/2}\|v_{k_{1},k_{2}}\|_{H^{1}}$ added, which is absorbed like the other top-level terms with a factor $\varepsilon^{\iota}$.

\textbf{Boundary term including $\tilde{F}_{\mathrm{bet}+}^{4,k_{1},k_{2}}$} Using the discussion for the surface term,  we first bound the following linear term when $\theta$ derivative hits on normal vector $\mathcal{N}_{0}$
\begin{align}
\sum_{j=1}^{k_{2}}\|\partial_{\theta}^{k_{2}-j}\partial_{t}^{k_{1}}v_{||}\cdot \nabla_{||}\partial_{\theta}^{j}\zeta\|_{L^{2}L^{\infty}(\Gamma_{0})}\lesssim \sum_{j=1}^{k_{2}}\varepsilon^{\iota}\|\partial_{\theta}^{k_{2}-j}\p_{t}^{k_{1}}v\|_{W^{2,q_{\delta}}}.
\end{align}

For nonlinear parts included in $\tilde{F}_{\mathrm{bet}+}^{4,k_{1},k_{2}}$. At the purely angular top level $(k_{1},k_{2})=(0,m+1)$, the main part $v_{||}\cdot \partial_{\theta}^{m+1}\nabla_{||}\eta$ of $\tilde{F}_{\mathrm{bet}+}^{4,0,m+1}$ gives on $\Gamma_{0}$ the contact-line transport term $\kappa\langle\partial_{t}\eta_{0,m+1},v_{\gamma}\partial_{\gamma}\eta_{0,m+1}\rangle_{\Gamma_{0}}$, with $v_{\gamma}=v\cdot\dot{\gamma}$, in the boundary term $\kappa\langle\partial_{t}\eta_{k_{1},k_{2}},\tilde{F}_{\mathrm{bet}+}^{4,k_{1},k_{2}}\rangle_{\Gamma_{0}}$ of \eqref{eq:local-domain-energy-1-2}. This term requires $\partial_{\theta}^{m+2}\eta\in L^{2}(\Gamma_{0})$. By Corollary~\ref{cor:top_level}~(iii), it is bounded by $C(\varepsilon)\sqrt{\mathcal{E}}\mathcal{D}$. For the case $k_{1}>0$, we apply the similar computation as in the proof of Theorem \ref{thm:whole-space-energy}.

\textbf{Step 7: Boundary term including $\tilde{F}_{\mathrm{bet}+}^{7,k_{1},k_{2}}$}

We first recall the geometry coefficients derived in the last Section. Denote $\widehat{\Phi}(\theta,r) = \Phi(\theta, r, -r\cot\theta_{\mathrm{eq}})$ and $g_\Sigma = (D\widehat{\Phi})^T D\widehat{\Phi}$. The graph metric gives $G_X = \bigl((1+|\nabla\zeta|^2)^{-1/2}\circ\widehat{X}\bigr)\, g_\Sigma^{-1}$. At $r=0$ the two surface tangent vectors are the images of $e_\theta$ and $e_r - \cot\theta_{\mathrm{eq}}\, e_h$ under $D\Phi(\theta,0,0)\in O(3)$. They are perpendicular, so the elements $m_{r\theta}$ for matrix $G_{X}$ satisfy $m_{r\theta}(\theta,0) = 0$ for every $\theta$. All its angular derivatives on the curve vanish. The outward conormal in \eqref{eq:local-system-h-eta} therefore contributes only radial derivatives to the fixed commutator. Each $\chi_{\mathrm{cl},\Sigma} \p_{t}^{a}\p_{\theta}^{j}\mathfrak{h}_\eta$ vanishes for $r > C\varepsilon$, and the cutoff equals one near $r=0$. Integration in $r$ gives $\|\partial_r(\chi_{\mathrm{cl},\Sigma} \p_{t}^{a}\p_{\theta}^{j}\mathfrak{h}_\eta) (\cdot,0)\|_{L^2_\theta} \le (C\varepsilon)^{1-1/s_-}\|\partial_r^2 (\chi_{\mathrm{cl},\Sigma} \p_{t}^{a}\p_{\theta}^{j}\mathfrak{h}_\eta)\|_{L^2_\theta L^{s_-}_r}$. Since $W_r^{3-1/q_\delta,q_\delta}\subset W_r^{2,s_-}$ and $1 - 1/s_- = \varepsilon_- \ge \iota$, we obtain
\begin{equation}
\bigl\|\sigma[\partial_\theta^b, \mathcal{B}_{\mathrm{cl},\Phi}]\partial_t^a\mathfrak{h}_\eta\bigr\|_{L^2(\Gamma_0)} \lesssim \varepsilon^\iota \sum_{j=0}^{b-1} \|(\chi_{\mathrm{cl},\Sigma} \p_{t}^{a}\p_{\theta}^{j}\mathfrak{h}_\eta)\|_{L^2_\theta W_r^{3-1/q_\delta,q_\delta}}.
\end{equation}

We now estimate the last sum in $\tilde F^7$, denoting it by $G$. On the curve, write $z=\nabla\eta$ in the fixed surface chart, and write $D^{a,b}= \partial_t^a\partial_\theta^b$. For $i\ge1$, apply the chain rule to $D^{k_1,k_2-i}R(x,z)$, with $R$ as in Remark~\ref{rem:R_x_dependence}. Each term is a numerical multiple of
\begin{equation}
(\partial_\theta^{h}\partial_z^{s}R)(x,\nabla \eta) [D^{a_1,b_1}\nabla \eta,\ldots,D^{a_s,b_s}\nabla \eta],\qquad \sum_{l=1}^{s}a_l=k_1,\qquad h+\sum_{l=1}^{s}b_l=k_2-i .
\end{equation}
Here all indices are nonnegative, $a_l+b_l\ge1$, and $\partial_\theta^h$ on the coefficient differentiates the fixed chart dependence with $z$ held fixed. The case $s=0$ is included when $k_1=0$ and $h=k_2-i$. The sum of the orders on the factors is at most $k_1+k_2-i\le m$. In particular, at most one factor has order $m$.

We record the trace bounds used on these products. The cutoffs equal one near $\Gamma_0$, and the chart gradient is a smooth linear combination of $\partial_r\eta$ and $\partial_\theta\eta$. For $a+b\le m$, $D^{a,b}z$ has $L^2(\Gamma_0)$ norm at most $C(\varepsilon)\sqrt{\mathcal D}$: for the transverse component use the $L^2_\theta W_r^{3-1/q_\delta,q_\delta}$ bound in \eqref{def:split_energy} and evaluate the first derivative at $r=0$ ($q_\delta>1$). For the angular component, $D^{a,b+1}\eta$ has the surface $H^{3/2-\alpha}$ bound in \eqref{def:energy_cl}; its trace belongs to $H^{1-\alpha}(\Gamma_0)\subset L^\infty(\Gamma_0)$. At $a+b\le m-1$ the same angular bound holds with $\sqrt{\mathcal E}$. For the transverse component at these lower orders, both $\partial_rD^{a,b}\eta$ and its angular derivative have $L^2(\Gamma_0)$ norm at most $C(\varepsilon)\sqrt{\mathcal E}$. Indeed, use the same radial trace at orders through $m-1$; when the extra angular derivative reaches order $m$, use the $L^2_\theta H_r^{3/2+(\varepsilon_--\alpha)/2}$ term of $\mathcal E_{\mathrm{cl}}$, whose radial exponent is strictly above $3/2$. The embedding $H^1(\Gamma_0)\subset L^\infty(\Gamma_0)$ therefore gives the energy bound for the lower transverse factors as well. Derivatives falling on fixed chart coefficients leave the same bounds at lower orders.

For $s\ge2$, put one factor of largest order in $L^2(\Gamma_0)$ and all others in $L^\infty(\Gamma_0)$. For $s=1$ the coefficient is $O(|z|)$; for $s=0$ it is $O(|z|^2)$, since $R(x,0)=\partial_zR(x,0)=0$ for every $x$. In this last case put one undifferentiated $z$ in each norm. The derivatives of the fixed normal are bounded, so $\mathcal E\le1$ gives
\begin{equation}\label{eq:mixed_contact_curvature_bound} \|G\|_{L^2(\Gamma_0)}\le C(\varepsilon)\sqrt{\mathcal E}\sqrt{\mathcal D}.
\end{equation}
When $k_2=0$ the sum defining $G$ is empty.

Moreover, for the nonlinear part, since $\tilde{F}^{4,j_{1},j_{2}}$ contains $\partial_{t}^{j_{1}}\partial_{\theta}^{\ell_{1}}v_{\theta}\cdot \partial_{t}^{j_{2}}\partial_{\theta}^{\ell_{2}} \partial_{\theta}\eta$ which vanishes in the two-dimensional case, we have to write down its estimate. When $k_{1}>0$, it suffices to estimate
\begin{align}
[v_{\theta}\cdot \partial_{t}^{k_{1}}\partial_{\theta}^{k_{2}} \partial_{\theta}\eta]_{\Gamma_{0}}\lesssim \|\chi_{\mathrm{cl}}v\|_{L^{\infty}W^{2,q_{\delta}}}\|\partial_{t}^{k_{1}}\partial_{\theta}^{k_{2}+1}\eta\|_{L^{2}(\Gamma_{0})}\lesssim\sqrt{\mathcal{E}\mathcal{D}}.
\end{align}
When $k_{1}=0$, the factor $\partial_{\theta}^{m+2}\eta$ on $\Gamma_{0}$ is controlled by the trace of the gain $\partial_{\theta}^{m+1}\eta\in H^{\frac{3}{2}}(\Sigma)$ of Theorem~\ref{thm:top_level}. By Corollary~\ref{cor:top_level}~(i),
\begin{align}
[v_{\theta}\cdot \partial_{\theta}^{m+1} \partial_{\theta}\eta]_{\Gamma_{0}}\lesssim \|\chi_{\mathrm{cl}}v\|_{L^{\infty}W^{2,q_{\delta}}}\|\partial_{\theta}^{m+2}\eta\|_{L^{2}(\Gamma_{0})}\lesssim\sqrt{\mathcal{E}\mathcal{D}},
\end{align}
and the corresponding contact-line term is bounded by $C(\varepsilon)\sqrt{\mathcal{E}}\mathcal{D}$; see Corollary~\ref{cor:top_level}~(iii).

We apply \eqref{eq:mixed_contact_curvature_bound} to both contact tests in \eqref{eq:local-domain-energy-1-2}. The norm of $\partial_t\eta_{k_1,k_2}$ on the curve is part of the dissipation. The top transport component of $\tilde F^{4,k_1,k_2}_{\mathrm{bet}+}$ was just estimated (use Corollary~\ref{cor:top_level}(i) at $(0,m+1)$). The summation defining $\tilde F^4$ excludes the term with every derivative on the velocity, which is already on the left in the kinematic law. In its remaining terms a derivative of fixed geometry multiplies a velocity of order at most $m$, or a differentiated velocity multiplies a differentiated surface gradient. In the latter products the total order is at most $m+1$. When the order on the surface gradient is $m$, the velocity has order at most one and is bounded on the curve by $C(\varepsilon)\sqrt{\mathcal E}$ (use the two-dimensional $W^{2,q_\delta}$ slice norm for the velocity and one more angular derivative). Otherwise the surface gradient has the lower-order supremum bound above, and we put the velocity in $L^2(\Gamma_0)$ using $L^2_\theta W^{2,q_\delta}_{r,h}$. Thus the entire first test has $L^2(\Gamma_0)$ norm at most $C(\varepsilon)\sqrt{\mathcal D}$. The same holds for $\omega\cdot\mathcal N$: $\omega$ contains only velocities of orders $(k_1,j)$, $j<k_2$, hence at most $m$, with fixed coefficients; $W^{2,q_\delta}$ in a two-dimensional slice controls point evaluation and $\mathcal N$ is bounded. Cauchy--Schwarz on the curve now bounds each pairing with $G$ by $C(\varepsilon)\sqrt{\mathcal E}\,\mathcal D$. For the first pairing, including the other terms of $\tilde F^7$, this gives
\begin{align}
&\big|[\tilde{F}_{\mathrm{bet}+}^{7,k_{1},k_{2}},\partial_{t}\eta_{k_{1},k_{2}}-\tilde{F}_{\mathrm{bet}+} ^{4,k_{1},k_{2}}]_{\Gamma_{0}}\big|\no\\
&\lesssim C(\varepsilon)\sqrt{\mathcal{E}}\mathcal{D}+ \varepsilon^{\iota}(\sum_{j=1}^{k_{2}}\varepsilon^{\iota}\|\partial_{\theta}^{k_{2}-j}\p_{t}^{k_{1}}v\|_{L^{2}W^{2,q_{\delta}}}+\|\p_{t}\eta_{k_{1},k_{2}}\|_{L^{2}(\Gamma_{0})})(\sum_{j=1}^{k_{2}}\|\partial_{t}^{k_{1}}\partial_{\theta}^{k_{2}-j}\eta\|_{L^{2}W^{3-\frac{1}{q_{\delta}},q_{\delta}}} ).
\end{align}

Therefore, it remains to control the cutoff commutator forcing terms  $F_{\mathrm{bet}+}^{1,k_{1},k_{2}}$, $F_{\mathrm{bet}+}^{2,k_{1},k_{2}}$, $F_{\mathrm{bet}+}^{3,k_{1},k_{2}}$ and $F_{\mathrm{bet}+}^{6,k_{1},k_{2}}$. By Subsection~\ref{subsec:cutoff_extension}, \(U:=\operatorname{supp}\nabla\chi_{\mathrm{bet}+} \subset\{\chi_{j,\mathrm{int}}=1\}\,\allowbreak (0\leq j\leq m).\) We use this transition region in the following estimates.

\textbf{Cutoff forcing term included in $F_{\mathrm{bet}+}^{1,k_{1},k_{2}}$}

Write $a=k_1$, $b=k_2$, so that $b>0$ and $a+b\leq m+1$. By the product rule, with $\chi=\chi_{\mathrm{bet}+}$,
\begin{equation}
\begin{aligned}
& \big\|S_{\mathcal A}(p_P,u_P)\nabla_{\mathcal A}\chi -\mu\operatorname{div}_{\mathcal A}\mathscr C_{\chi,\mathcal A}(u_P) \big\|_{L^{6/5}(\Omega)} \leq C(\varepsilon) \big(\|p_P\|_{L^p(U)}+\|u_P\|_{W^{1,p}(U)} +\|\nabla\mathcal A\|_{L^2(U)}\|u_P\|_{L^3(U)}\big) \\
& \leq C(\varepsilon)\big( \|\chi_{a,\mathrm{int}}\partial_t^a p\|_{W^{b,p}}+ \|\chi_{a,\mathrm{int}}\partial_t^a v\|_{W^{b+1,p}}\big) \leq C(\varepsilon)\sqrt{\mathcal D_{\mathrm{int}}}.
\end{aligned}
\end{equation}
The first inequality uses boundedness of $\mathcal A$ and the first two derivatives of the fixed cutoff (with constants depending on the width). Step~2 bounds $\nabla\mathcal A$ on this support in $L^\infty_\theta L^{r_-}_{r,h}$ with $r_->2$, hence in $L^2(U)$. For the $L^3$ norm apply $W^{1,p}(\Omega)\subset L^3(\Omega)$ to the interior-localized differentiated velocity; $p\geq12/7$ suffices. The chart and the frame are fixed and smooth, and the last inequality uses $b\leq m+1-a$ in \eqref{def:interior_dissipation}. This estimates only terms with a differentiated cutoff; the commutator $\mathfrak O_{1,k_1,k_2}$ is the separate term estimated in Step~2.

\textbf{Cutoff forcing term included in $F_{\mathrm{bet}+}^{2,k_{1},k_{2}}$}

Similarly, we have \(\|\nabla_{\mathcal{A}}\chi_{\mathrm{bet}+}(\partial_{t}^{k_{1}}\partial_{\theta}^{k_{2}}v)\|_{L^{2}}\lesssim\frac{C(\varepsilon)}{\varepsilon}\sqrt{\mathcal{D}_{\mathrm{int}}}.\)

\textbf{Cutoff forcing term included in $F_{\mathrm{bet}+}^{3,k_{1},k_{2}}$}

Similarly, we have
\begin{equation}
\begin{aligned}
& \|\nabla_{\mathcal{A}}\chi_{\mathrm{bet}+}\otimes(\partial_{t}^{k_{1}}\partial_{\theta}^{k_{2}}v)\|_{L^{2}(\Sigma)}+\|(\partial_{t}^{k_{1}}\partial_{\theta}^{k_{2}}v)\otimes\nabla_{\mathcal{A}}\chi_{\mathrm{bet}+}\|_{L^{2}(\Sigma)}+\|\nabla_{\Sigma}\chi_{\mathrm{bet}+,\Sigma}\cdot \nabla_{\Sigma}\eta_{k_{1},k_{2}}\|_{L^{2}}+\|\nabla_{\Sigma}^{2}\chi_{\mathrm{bet}+,\Sigma}\eta_{k_{1},k_{2}}\|_{L^{2}} \\
& \lesssim (\|\chi_{0,\mathrm{int}}\partial_{t}^{k_{1}}\partial_{\theta}^{k_{2}}v\|_{W^{1,p}}+\|\chi_{0,\mathrm{int}}\partial_{t}^{k_{1}}\partial_{\theta}^{k_{2}}\eta\|_{W^{2-\frac{1}{p},p}}) \lesssim \frac{C(\varepsilon)}{\varepsilon}\sqrt{\mathcal{D}_{\mathrm{int}}}.
\end{aligned}
\end{equation}

For the commutator cut off function for nonlinear term, we note that $\operatorname{supp}\{\nabla\chi_{bet+,\Sigma}\}\subset \{\chi_{\mathrm{int}+,\Sigma}=1\}\cap \{\chi_{\mathrm{local}}=1\}$. In the following computation, we abuse the notation for example $\|\chi_{\mathrm{int}+,\Sigma}\p_{t}^{k_{1}}\p_{\theta}^{k_{2}}\eta\|_{H^{1}}$ we to denote $\|\chi_{\mathrm{int}+,\Sigma}\chi_{\mathrm{local}}\p_{t}^{k_{1}}\p_{\theta}^{k_{2}}\eta\|_{H^{1}}$ so that $\p_{\theta}$ is well-defined in the support of the cut off functions. Then we bound the multiplier with the $m+1$-th differentiated order by $L^{2}$ and the others by $L^{\infty}$. We have,
\begin{equation}
\begin{aligned}
& \|\sigma(\nabla_x\chi_\Sigma\cdot R_P)\mathcal N\|_{L^2(\Sigma)} \\
& \lesssim_{\varepsilon}  \sum_{j_{2}=0}^{k_{2}}\sum_{r=3}^{j_{2}+k_{1}+2}\sum_{\substack{i_{1,1}+\ldots+i_{1,r}=k_{1}\\
i_{2,1}+\ldots +i_{2,r}=j_{2}\\
i_{1,\ell},i_{2,\ell}\geq 0\\
i_{1,\ell}\geq 1~\textbf{or}~i_{2,\ell}\geq 1~\text{for}~\ell>2\\
i_{1,\ell}+i_{2,\ell}\neq k_{1}+k_{2}}} \prod_{\ell=1}^{r}\|\chi_{\mathrm{int}+,\Sigma}\p_{t}^{i_{1,\ell}}\p_{\theta}^{i_{2,\ell}}\eta\|_{W^{1,\infty}}+\sum_{j_{2}=0}^{k_{2}}\|\chi_{\mathrm{int}+,\Sigma}\p_{t}^{k_{1}}\p_{\theta}^{j_{2}}\eta\|_{H^{1}}\|\chi_{\mathrm{int}+,\Sigma}\eta\|_{W^{1,+\infty}} \\
& +\sum_{j_{2}=0}^{k_{2}}\sum_{j_{1,1}+j_{1,2}=k_{1}}\sum_{\substack{j_{2,1}+j_{2,2}=j_{2}\\
j_{2,1}\neq j_{2}~\text{when}~j_{1,1}=k_{1}\\
j_{2,2}\neq j_{2}~\text{when}~j_{1,2}=k_{1}}}\|\chi_{\mathrm{int}+,\Sigma}\p_{t}^{j_{1,1}}\p_{\theta}^{j_{2,1}}\eta\|_{W^{1,\infty}}\|\chi_{\mathrm{int}+,\Sigma}\p_{t}^{j_{1,2}}\p_{\theta}^{j_{2,2}}\eta\|_{W^{1,+\infty}} \\
& +\|\chi_{\mathrm{int}+,\Sigma}\eta\|_{W^{1,+\infty}}^{2}\|\chi_{\mathrm{int}+,\Sigma}\p_{t}^{k_{1}}\p_{\theta}^{k_{2}}\eta\|_{H^{1}} \lesssim  C(\varepsilon)\sqrt{\mathcal{E}}\sqrt{\mathcal{D}},
\end{aligned}
\end{equation}
where we used the fact that in the support region of  $\operatorname{supp}\{\chi_{\mathrm{local}}\}$, we can use the localized angular derivative, and we use the fact that $\|\chi_{\mathrm{int}+,\Sigma}\p_{t}^{a}\p_{\theta}^{b}\eta\|_{W^{1,+\infty}}\lesssim \|\chi_{\mathrm{int}+,\Sigma}\p_{t}^{a}\eta\|_{W^{3-\frac{1}{p}+b,p}}$ given that $a+b\leq m$, and the fact that $\|\chi_{\mathrm{int}+,\Sigma}\p_{t}^{a}\p_{\theta}^{b}\eta\|_{H^{1}}\lesssim \|\chi_{\mathrm{int}+,\Sigma}\p_{t}^{a}\eta\|_{W^{3-\frac{1}{p}+b-1,p}}$ given that $a+b=m+1$ and $b\geq 0$

\textbf{Cutoff forcing term included in $F_{\mathrm{bet}+}^{6,k_{1},k_{2}}$}

Similarly, we have
\begin{align}
\|\nabla_{\mathcal{A}} \chi_{\mathrm{bet}+}\otimes(\partial_{t}^{k_{1}}\partial_{\theta}^{k_{2}}v)\|_{L^{2}(\Sigma_{s})}+\|(\partial_{t}^{k_{1}}\partial_{\theta}^{k_{2}}v)\otimes \nabla_{\mathcal{A}}\chi_{\mathrm{bet}+}\|_{L^{2}(\Sigma_{s})}\lesssim \frac{C(\varepsilon)}{\varepsilon}\sqrt{\mathcal{D}_{\mathrm{int}}}.
\end{align}
Combining all the estimates for the forcing terms, we have
\begin{equation}
\label{eq:local-domain-energy-derivative}
\begin{aligned}
& \int_\Omega Jp_{k_{1},k_{2}}F_{\mathrm{bet}+,}^{2,k_{1},k_{2}} + \int_\Omega JF_{\mathrm{bet}+}^{1,k_{1},k_{2}}\cdot v_{k_{1},k_{2}} - \left\langle F_{\mathrm{bet}+}^{3,k_{1},k_{2}},v_{k_{1},k_{2}}\right\rangle_\Sigma - \left\langle J_sF_{\mathrm{bet}+}^{6,k_{1},k_{2}},U_0v_{k_{1},k_{2}}\right\rangle_{\Sigma_s}+ g\int_\Sigma\eta_{k_{1},k_{2}}F_{\mathrm{bet}+}^{4,k_{1},k_{2}} \\
& + \sigma\int_\Sigma a_\zeta\nabla_\Sigma\eta_{k_{1},k_{2}}\cdot \nabla_{\Sigma}F_{\mathrm{bet}+}^{4,k_{1},k_{2}} + \kappa\left\langle \partial_t\eta_{k_{1},k_{2}},F_{\mathrm{bet}+}^{4,k_{1},k_{2}} \right\rangle_{\Gamma_0} \\
& - \left\langle F_{\mathrm{bet}+}^{7,k_{1},k_{2}},\partial_t\eta_{k_{1},k_{2}}-F_{\mathrm{bet}+}^{4,k_{1},k_{2}} \right\rangle_{\Gamma_0} + \sigma\int_\Sigma F_{\mathrm{bet}+}^5\cdot\nabla_\Sigma F_{\mathrm{bet}+,}^{4,j} \lesssim (\|p_{k_{1},k_{2}}\|_{L^{2}}+\|v_{k_{1},k_{2}}\|_{H^{1}}) \\
& (\sqrt{\mathcal{E}}\sqrt{\mathcal{D}}+\varepsilon^{\iota}\sqrt{\mathscr{E}_{k_{1},k_{2}-1,\mathrm{el},\mathrm{cl}}}+\varepsilon^{\iota}\sqrt{\mathscr{E}_{k_{1},k_{2}-1,\mathrm{el},\mathrm{bet}}}+\varepsilon^{\iota}(\sqrt{\mathscr{E}_{k_{1},k_{2},\mathrm{energy},\mathrm{cl}}}+\sqrt{\mathscr{E}_{k_{1},k_{2},\mathrm{energy},\mathrm{bet}}})) \\
& +(\|\eta_{k_{1},k_{2}}\|_{H^{1}}+[\partial_{t}\eta_{k_{1},k_{2}}]_{\Gamma_{0}}) \\
& \times (\sqrt{\mathcal{E}}\sqrt{\mathcal{D}}+\varepsilon^{\iota}\sqrt{\mathscr{E}_{k_{1},k_{2}-1,\mathrm{el},\mathrm{cl}}}+\varepsilon^{\iota}\sqrt{\mathscr{E}_{k_{1},k_{2}-1,\mathrm{el},\mathrm{bet}}}+\varepsilon^{\iota}(\sqrt{\mathscr{E}_{k_{1},k_{2},\mathrm{energy},\mathrm{cl}}}+\sqrt{\mathscr{E}_{k_{1},k_{2},\mathrm{energy},\mathrm{bet}}})) \\
& +(\|p_{k_{1},k_{2}}\|_{L^{2}}+\|v_{k_{1},k_{2}}\|_{H^{1}}+\|\eta_{k_{1},k_{2}}\|_{H^{1}}+[\partial_{t}\eta_{k_{1},k_{2}}]_{\Gamma_{0}})(\frac{C(\varepsilon)}{\varepsilon}\sqrt{\mathcal{D}_{\mathrm{int}}}).
\end{aligned}
\end{equation}

\textbf{Step 8: Estimates for Terms with $\omega$} Using the previous forcing estimate, we have
\begin{align}
\int_{\Omega}J_{\eta}\omega v_{k_{1},k_{2}}\lesssim \varepsilon^{\iota}\|v_{k_{1},k_{2}}\|_{L^{2}}^{2}+\frac{1}{\varepsilon^{\iota}}\sum_{j=0}^{k_{2}-1}\|v_{k_{1},j}\|_{L^{2}}^{2}.
\end{align}
\begin{align}
\int_{\Omega}\partial_{t}(J_{\eta}\omega )v_{k_{1},k_{2}}\lesssim \varepsilon^{\delta }\|v_{k_{1},k_{2}}\|^{2}_{H^{1}}+\sum_{j=0}^{k_{2}-1}(\|v_{k_{1}+1,j}\|_{H^{1}}^{2}.
\end{align}
\begin{equation}
\int_{\Omega}Jp_{k_{1},k_{2}}\func{div}_{\mathcal{A}}\omega \lesssim \|p_{k_{1},k_{2}}\|_{L^{2}}(\sum_{j=0}^{k_{2}-1}\|v_{k_{1},j}\|_{H^{1}}) \lesssim \varepsilon^{\iota}\|p_{k_{1},k_{2}}\|_{L^{2}}^{2}+\frac{1}{\varepsilon^{\iota}}\sum_{j=0}^{k_{2}-1}\|v_{k_{1},j}\|_{H^{1}}^{2}.
\end{equation}
\begin{align}
&\int_{\Omega}J\tilde{F}_{\mathrm{cl}+}^{1,k_{1},k_{2}}\cdot \omega-\langle \tilde{F}_{\mathrm{cl}+}^{3,k_{1},k_{2}},\omega\rangle_{\Sigma}-\langle J_{s}\tilde{F}_{\mathrm{cl}+}^{6,k_{1},k_{2}},U_{0}\omega \rangle_{\Sigma_{s}} \\
&\lesssim (\frac{1}{\sqrt{\varepsilon}}\sqrt{\mathcal{D}_{\mathrm{int}}}+\varepsilon^{\iota}\sqrt{\mathscr{E}_{k_{1},k_{2}-1,\mathrm{el},\mathrm{cl}}}+\varepsilon^{\iota}\sqrt{\mathscr{E}_{k_{1},k_{2}-1,\mathrm{el},\mathrm{bet}}}+\sqrt{\mathscr{E}_{k_{1},k_{2}-1,\mathrm{energy},\mathrm{cl}}}+\sqrt{\mathcal{E}\mathcal{D}})(\sum_{j=0}^{k_{2}-1}\|v_{k_{1},j}\|_{H^{1}}).
\end{align}
\begin{equation}
|-\mu/2 \int_{\Omega}J\mathbb{D}_{\mathcal{A}_{\eta}}v_{k_{1},k_{2}}:\mathbb{D}_{\mathcal{A}_{\eta}}\omega|+\langle J_{s}\beta_{0}U_{0}v_{k_{1},k_{2}},U_{0}\omega\rangle_{\Sigma} \lesssim \varepsilon^{\iota}\|v_{k_{1},k_{2}}\|^{2}_{H^{1}}+\frac{1}{\varepsilon^{\iota}}(\sum_{j=0}^{k_{2}-1}\|v_{k_{1},j}\|_{H^{1}}^{2}).
\end{equation}

Using the fact that in $2-D$, $W^{2,q_{\nu}}(D)\hookrightarrow L^{\infty }(\partial D)$ for $q_{\nu}>1$,
\begin{align}
&\langle \partial_{t}\eta_{k_{1},k_{2}}+\tilde{F}_{\mathrm{cl}+}^{7,k_{1},k_{2}},\omega\cdot \mathcal{N}\rangle_{\Gamma_{0}}  \notag \\
&\lesssim ([\partial_{t}\eta_{k_{1},k_{2}}]_{\Gamma_{0}}+\varepsilon^{\iota}\sqrt{\mathscr{E}_{k_{1},k_{2}-1,\mathrm{el},\mathrm{cl}}}+\varepsilon^{\iota}\sqrt{\mathscr{E}_{k_{1},k_{2}-1,\mathrm{el},\mathrm{bet}}}+\sqrt{\mathscr{E}_{k_{1},k_{2},\mathrm{energy},\mathrm{cl}}}+C(\varepsilon)\sqrt{\mathcal{E}\mathcal{D}})(\sum_{j=0}^{k_{2}-1}\|v_{k_{1},j}\|_{L^{2}W^{2,q_{\nu}}}) \\
&\lesssim ([\partial_{t}\eta_{k_{1},k_{2}}]_{\Gamma_{0}}+\varepsilon^{\iota}\sqrt{\mathscr{E}_{k_{1},k_{2}-1,\mathrm{el},\mathrm{cl}}}+\varepsilon^{\iota}\sqrt{\mathscr{E}_{k_{1},k_{2}-1,\mathrm{el},\mathrm{bet}}}+\sqrt{\mathscr{E}_{k_{1},k_{2},\mathrm{energy},\mathrm{cl}}}+C(\varepsilon)\sqrt{\mathcal{E}\mathcal{D}})(\varepsilon^{\iota}\sum_{j=0}^{k_{2}-1}\|v_{k_{1},j}\|_{L^{2}W^{2,q_{\delta}}}).
\end{align}
Finally, as above, we have
\begin{equation}
\begin{aligned}
&g\int_{\Sigma}\eta_{k_{1},k_{2}}\omega+\sigma\int_{\Sigma}(a_{\zeta}\nabla_{\Sigma}\eta_{k_{1},k_{2}}+\tilde{F}_{\mathrm{cl}}^{5,k_{1},k_{2}})\cdot \nabla_{\Sigma}(\omega\cdot \mathcal{N})  \\ \lesssim& (\|\eta_{k_{1},k_{2}}\|_{H^{1}}+\varepsilon^{\iota}\sqrt{\mathscr{E}_{k_{1},k_{2}-1,\mathrm{el},\mathrm{cl}}}+\varepsilon^{\iota}\sqrt{\mathscr{E}_{k_{1},k_{2}-1,\mathrm{el},\mathrm{bet}}}+\sqrt{\mathscr{E}_{k_{1},k_{2}-1,\mathrm{energy},\mathrm{cl}}}+\sqrt{\mathcal{E}\mathcal{D}})\\ &\quad(\sum_{j=0}^{k_{2}-1}\|v_{k_{1},j}\|_{L^{2}W^{2,q_{\nu}}}+\sum_{j=0}^{k_{2}}\|v_{k_{1},j}\|_{L^{2}W^{1,q_{\nu}}})  \\ \lesssim&(\|\eta_{k_{1},k_{2}}\|_{H^{1}}+\varepsilon^{\iota}\sqrt{\mathscr{E}_{k_{1},k_{2}-1,\mathrm{el},\mathrm{cl}}}+\varepsilon^{\iota}\sqrt{\mathscr{E}_{k_{1},k_{2}-1,\mathrm{el},\mathrm{bet}}}+\sqrt{\mathscr{E}_{k_{1},k_{2}-1,\mathrm{energy},\mathrm{cl}}}+\sqrt{\mathcal{E}\mathcal{D}})\\ &\quad\varepsilon^{\iota}(\sum_{j=0}^{k_{2}-1}\|v_{k_{1},j}\|_{L^{2}W^{2,q_{\delta}}}+\sum_{j=0}^{k_{2}}\|v_{k_{1},j}\|_{L^{2}H^{1}}).
\end{aligned}
\end{equation}

In the last bound and in \eqref{eq:local-domain-energy-derivative}, the factor $\|\eta_{k_{1},k_{2}}\|_{H^{1}}$ would need $H^{1}(\Sigma)$ traces of $W^{2,q}$ slice functions, which are not available for $q_{-}<4/3$. We use instead $\|\eta_{k_{1},k_{2}}\|_{H^{\frac{3}{2}-\alpha}}$ against the $H^{\frac{1}{2}+\alpha}(\Sigma)$ norm of the other factor, and combine with Step~9 before absorbing, as in \eqref{eq:induction_inequality}. Combining all of the estimates and steps above, and using the Cauchy inequality, we obtain
\begin{equation}
\label{eq:local-domain-energy-derivative_between}
\begin{aligned}
& {}   \Bigg[ \frac12\int_\Omega J|v_{k_{1},k_{2}}|^2 + \|\eta_{k_{1},k_{2}}\|_{H^1(\Sigma)}^2 \Bigg](t) +\int_{0}^{t}\big( \frac{\mu}{2}\int_\Omega J|\mathbb{D}_{\mathcal{A}_\eta}v_{k_{1},k_{2}}|^2 + \beta_0\int_{\Sigma_s}J_s|U_0v_{k_{1},k_{2}}|^2 + \kappa\|\partial_t\eta_{k_{1},k_{2}}\|_{L^2(\Gamma_0)}^2\big) \\
& \lesssim \int_{0}^{t}(\sqrt{\mathcal{E}}{\mathcal{D}}+\varepsilon^{\iota}{\mathscr{E}_{k_{1},k_{2}-1,\mathrm{el},\mathrm{cl}}}+\varepsilon^{\iota}{\mathscr{E}_{k_{1},k_{2}-1,\mathrm{el},\mathrm{bet}}}+\varepsilon^{\iota}({\mathscr{E}_{k_{1},k_{2},\mathrm{energy},\mathrm{cl}}}+{\mathscr{E}_{k_{1},k_{2},\mathrm{energy},\mathrm{bet}}}) \\
& +\frac{C(\varepsilon)}{\varepsilon} \mathcal{D}_{\mathrm{int}}+\frac{C(\varepsilon)}{\varepsilon}\mathscr{E}_{k_{1},k_{2}-1,\mathrm{energy},\mathrm{cl}}) +\mathcal{E}(0)+\frac{1}{\varepsilon^{\iota}}\sum_{j=0}^{k_{2}-1}\|\partial_{t}^{k_{1}}\partial_{\theta}^{j}v\|^{2}_{L^{2}}(t)+\sqrt{\mathcal{E}}^{3}(t).
\end{aligned}
\end{equation}

Using the elliptic estimates derived in Theorem \ref{thm:intermediate_estimate} and Theorem \ref{thm:Higher-order contact line elliptic estimate}, we have
\begin{align}
{\label{eq:local-domain-energy-derivative_final}} &\Bigg[ \frac12\int_\Omega J|v_{k_{1},k_{2}}|^2 + \|\eta_{k_{1},k_{2}}\|_{H^1(\Sigma)}^2 \Bigg](t) +\int_{0}^{t}\big( \frac{\mu}{2}\int_\Omega J|\mathbb{D}_{\mathcal{A}_\eta}v_{k_{1},k_{2}}|^2 + \beta_0\int_{\Sigma_s}J_s|U_0v_{k_{1},k_{2}}|^2 + \kappa\|\partial_t\eta_{k_{1},k_{2}}\|_{L^2(\Gamma_0)}^2\big)  \notag \\
\lesssim&\int_{0}^{t}(\sqrt{\mathcal{E}}{\mathcal{D}}+{\mathscr{E}_{k_{1}+1,k_{2}-1,\mathrm{energy},\mathrm{cl}}}+{\mathscr{E}_{k_{1}+1,k_{2}-1,\mathrm{energy},\mathrm{bet}}}+\varepsilon^{\iota}({\mathscr{E}_{k_{1},k_{2},\mathrm{energy},\mathrm{cl}}}+{\mathscr{E}_{k_{1},k_{2},\mathrm{energy},\mathrm{bet}}})+\frac{C(\varepsilon)}{\varepsilon} \mathcal{D}_{\mathrm{int}} \notag \\
&\quad+\frac{C(\varepsilon)}{\varepsilon}\mathscr{E}_{k_{1},k_{2}-1,\mathrm{energy},\mathrm{cl}}) +\mathcal{E}(0)+\frac{C(\varepsilon)}{\varepsilon^{\iota}}\sum_{j=0}^{k_{2}-1}\|\partial_{t}^{k_{1}}\partial_{\theta}^{j}v\|^{2}_{L^{2}}(t)+\sqrt{\mathcal{E}}^{3}(t).
\end{align}
Therefore, bounding $\frac{C(\varepsilon)}{\varepsilon}\int_{0}^{t}\mathcal{D}_{\mathrm{int}}$ by \eqref{eq:Dint_E0}, we have
\begin{align}
{\label{local_energy_estimate}} &\Bigg[ \frac12\int_\Omega J|v_{k_{1},k_{2}}|^2 + \|\eta_{k_{1},k_{2}}\|_{H^1(\Sigma)}^2 \Bigg](t) +\int_{0}^{t}\big( \frac{\mu}{2}\int_\Omega J|\mathbb{D}_{\mathcal{A}_\eta}v_{k_{1},k_{2}}|^2 + \beta_0\int_{\Sigma_s}J_s|U_0v_{k_{1},k_{2}}|^2 + \kappa\|\partial_t\eta_{k_{1},k_{2}}\|_{L^2(\Gamma_0)}^2\big)  \notag \\
\lesssim& C(\varepsilon)(\mathcal{E}(0)+\int_{0}^{t}\sqrt{\mathcal{E}}\mathcal{D}+\mathcal{E}^{\frac{3}{2}}(t))+\int_{0}^{t}(\frac{C(\varepsilon)}{\varepsilon}\mathscr{E}_{k_{1},k_{2}-1,\mathrm{energy},\mathrm{cl}}+\mathscr{E}_{k_{1}+1,k_{2}-1,\mathrm{energy},\mathrm{cl}}+\varepsilon^{\iota}\mathscr{E}_{k_{1},k_{2},\mathrm{energy},\mathrm{cl}})  \notag \\
&+\frac{C(\varepsilon)}{\varepsilon^{\iota}}\sum_{j=0}^{k_{2}-1}\|\partial_{t}^{k_{1}}\partial_{\theta}^{j}v\|^{2}_{L^{2}}(t)+\sqrt{\mathcal{E}}^{3}(t).
\end{align}

\textbf{Step 9: Enhanced $\frac{3}{2}-\alpha$ estimate for the surface perturbation function}.

In the horizontal domain $D$ of Lemma 11.1, for any $0\leq a\leq k_{1}$ and $1\leq b\leq k_{2}$ put
\begin{equation*}
\eta_{a,b}^\circ = \eta_{a,b} - |D|^{-1}\int_D \eta_{a,b}, \qquad s = 1 - 2\alpha,\quad G_N = D_N^s \eta_{a,b}^\circ, \qquad X_N = \|P_N \eta_{a,b}^\circ\|_{\mathcal{H}_{\mathcal{K}}^{1+s/2}},
\end{equation*}
where $P_N$ projects onto the first $N+1$ eigenfunctions and $D_N^s = \mathcal{K}^{s/2} P_N$. Take the harmonic Neumann lift $\psi_N$ from Lemma 11.1 with top flux $G_N/|\mathcal{N}_0|$ and zero wall flux. Then $\omega_N = J^{-1}\mathcal{A}_\eta^{-T}\nabla\psi_N$ satisfies
\begin{equation*}
\operatorname{div}_{\mathcal{A}_\eta}\omega_N = 0, \qquad \omega_N\cdot\mathcal{N} = G_N \ \text{on } \Sigma,\quad \|\omega_N\|_{H^1(\Omega)} + \|G_N\|_{H^{1/2+\alpha}(D)} \le C X_N.
\end{equation*}
By the spectral identity in Appendix A, the gravity and capillary terms in (11.33) sum to $X_N^2$. On the contact curve use $\|G_N\|_{L^2(\Gamma_0)} \le C\|G_N\|_{H^{1/2+\alpha}(D)}$. Testing (11.33) gives
\begin{align}
{\label{eq:local_3_2}} X_{N}^{2} \lesssim &\|\tilde{F}_{\mathrm{bet}+}^{1,a,b}-\tilde{F}_{bet+}^{3,a,b}-\tilde{F}_{bet+}^{6,a,b}\|^{2}_{(H^{1})^{*}}+\|\tilde{F}_{\mathrm{bet}+}^{5,a,b}\|_{H^{\frac{1}{2}-\alpha}}^{2}+\|v_{a,b}\|_{H^{1}}^{2}+|\int_{\Omega}J\partial_{t}v_{a,b}\cdot \omega_{N}|  \notag \\
&+[\partial_{t}^{a+1}\partial_{\theta}^{b}\eta]_{\Gamma_{0}}^{2}+\|\tilde{F}_{\mathrm{bet}+}^{7,a,b}\|^{2}_{L^{2}(\Gamma_{0})}.
\end{align}
The pressure does not appear in \eqref{eq:local_3_2}: for the test field $\omega=M\nabla\psi$ of Lemma~\ref{lem:local_threehalf}, the Piola identity gives $\func{div}_{\mathcal{A}}(M\nabla\psi)=K\Delta\psi=0$.

For the last term on the right-hand side of \eqref{eq:local_3_2}, since $b>0$, using integration by parts in $\theta$, we obtain
\begin{equation}
|\int_{\Omega}J\partial_{t}v_{a,b}\cdot \omega_{N}| \lesssim \|v_{a+1,b-1}\|_{L^{2}}\|\partial_{\theta}\omega_{N}\|_{L^{2}}+\|\partial_{\theta}(J_{\eta}J_{\Phi}Q)\|_{L^{\infty}}\|v_{a+1,b-1}\|_{L^{2}}\|\omega_{N}\|_{L^{2}} \lesssim \|v_{a+1,b-1}\|_{L^{2}}X_{N}.
\end{equation}
Substitute this bound into \eqref{eq:local_3_2}, and then using Cauchy inequality on this integral by placing a small constant independent od $\varepsilon$ on $X_{N}^{2}$ so that it can be absorbed by the left hand side. By \eqref{functional:exact_spectral_identity}, $\sup_N X_N=\|\eta_{a,b}^\circ\|_{\mathcal H_{\mathcal K}^{1+s/2}}$. Since $1+s/2=3/2-\alpha$, \eqref{functional:subcritical} gives the same bound up to a constant for $\|\eta_{a,b}^\circ\|^{2}_{H^{3/2-\alpha}(D)}$. Restoring the mean using $\|\eta_{a,b}-\eta_{a,b}^\circ\|_{H^{3/2-\alpha}(D)}\lesssim\|\eta_{a,b}\|_{L^2(D)}$ yields
\begin{equation}
\label{eq:local_3_2_f}
\begin{aligned}
& {} \|\partial_{t}^{a}\partial_{\theta}^{b}\eta_{bet_{+}}\|_{H^{\frac{3}{2}-\alpha}}^{2} \lesssim  \|\tilde{F}_{\mathrm{bet}+}^{1,a,b}-\tilde{F}_{bet+}^{3,a,b}-\tilde{F}_{bet+}^{6,a,b}\|^{2}_{(H^{1})^{*}}+\|\tilde{F}_{\mathrm{bet}+}^{5,a,b}\|_{H^{\frac{1}{2}-\alpha}}^{2}+\|v_{a,b}\|_{H^{1}}^{2} \\
& +[\partial_{t}^{a+1}\partial_{\theta}^{b}\eta]_{\Gamma_{0}}^{2}+\|\partial_{t}^{a}\partial_{\theta}^{b}\eta_{\mathrm{bet}_{+}}\|_{L^{2}}^{2}+\|\tilde{F}_{\mathrm{bet}+}^{7,a,b}\|^{2}_{L^{2}(\Gamma_{0})}+\|v_{a+1,b-1}\|_{L^{2}}^{2}.
\end{aligned}
\end{equation}

On the right hand side of the inequality above, $\|\tilde{F}_{\mathrm{bet}+}^{5,a,b}\|_{H^{\frac{1}{2}-\alpha}}^{2}$, up to the lower order terms involving $\theta$ derivatives of the steady states, we apply similar estimate as in the proof of \ref{thm:whole_space_threehalf} by multiplying $\chi_{bet+}$ on the term with highest angular derivative. When at least one angular derivative hit on the steady state, less angular derivatives hit on the perturbation function $\eta$. This term can be similarly bounded by lower order estimates. The commutators of cut off functions and spatial derivatives are bounded by interior norms. We have $\|\tilde{F}_{\mathrm{bet}+}^{5,a,b}\|_{H^{\frac{1}{2}-\alpha}}^{2}\lesssim C(\varepsilon){\mathcal{E}\mathcal{D}}$.  Therefore, combining the estimate above for forcing terms with the estimate given in Section 7 with $k_{2}=0$, we have
\begin{align}
\int_{0}^{t}\sum_{j_{1}=0}^{k_{1}}\sum_{j_{2}=0}^{k_{2}}\|\partial_{t}^{j_{1}}\partial_{\theta}^{j_{2}}\eta_{\mathrm{bet}_{+}}\|_{H^{\frac{3}{2}-\alpha}}^{2}\lesssim\text{right-hand side of \eqref{local_energy_estimate}}.
\end{align}
At the purely angular top level $(0,m+1)$, changing the cut off function from $\chi_{\mathrm{bet}+,\Sigma}$ to $\chi_{\mathrm{cl},\Sigma}$, Theorem~\ref{thm:top_level} improves this bound $\|\p_{\theta}^{m+1}\eta_{\mathrm{cl}}\|_{H^{\frac{3}{2}-\alpha}}^{2}$ to $H^{\frac{3}{2}}$, without the loss $\alpha$, at every fixed time.

\textbf{Step 10: Pressure estimate}

Given that $k_{2}> 0$, we can bound $p_{k_{1},k_{2}}$ as follows.

At the purely angular top level $(k_{1},k_{2})=(0,m+1)$ we use the fixed-time bound \eqref{tl:pmain} of Theorem~\ref{thm:top_level}. Its proof, Part~A in Subsection~\ref{subsec:top_level}, is the argument below, with the Ne\caronaccent{c}as inequality applied on $\Omega$ instead of $\Omega_{\mathrm{bet}+}$. Since $p_{0,m+1}=\chi_{\mathrm{bet}+}\partial_{\theta}^{m+1}p$ and $v_{0,m+1}=\chi_{\mathrm{bet}+}\partial_{\theta}^{m+1}v$, squaring \eqref{tl:pmain} and integrating over $[0,t]$ gives
\begin{align}
\int_{0}^{t}\|p_{0,m+1}\|_{L^{2}}^{2}\leq3C_{0}^{2}\int_{0}^{t}\|v_{0,m+1}\|_{H^{1}}^{2}+3C_{0}^{2}\int_{0}^{t}\mathscr{E}_{1,m,\mathrm{energy},\mathrm{cl}}+C(\varepsilon)\int_{0}^{t}\big(\mathcal{D}_{\mathrm{int}}+\mathcal{E}\mathcal{D}\big).  \label{eq:top_pressure_time}
\end{align}
The right side contains no top-level pressure, and $C_{0}$ does not depend on $\varepsilon$. The terms of $\mathscr{E}_{1,m,\mathrm{energy},\mathrm{cl}} $ are dissipation terms at the levels $(1,m)$ and $(0,j)$, $j\leq m$, which are bounded before the level $(0,m+1)$ in the induction; $C(\varepsilon)\int_{0}^{t}\mathcal{D}_{\mathrm{int}}$ is bounded by \eqref{eq:Dint_E0}; and $\int_{0}^{t}\mathcal{E}\mathcal{D}\leq\sup_{[0,t]}\sqrt{\mathcal{E}}\int_{0}^{t}\mathcal{D}$. In the energy estimate at the level $(0,m+1)$, Steps~1 and~8 and Cauchy's inequality bring in the top-level pressure only through $\varepsilon^{\iota}\int_{0}^{t}\|p_{0,m+1}\|_{L^{2}}^{2}$. By \eqref{eq:top_pressure_time} it contributes $3C_{0}^{2}\varepsilon^{\iota}\int_{0}^{t}\|v_{0,m+1}\|_{H^{1}}^{2}$, which is absorbed by the dissipation on the left of \eqref{local_energy_estimate} once $3C_{0}^{2}C_{*}C_{K}\varepsilon^{\iota}\leq\frac{1}{8}\min\{\mu,\beta_{0}\}$. Here $C_{K}$ is the constant of Lemma~\ref{lem:korn}, $C_{*}$ is the implicit constant of \eqref{local_energy_estimate}, and $J,J_{s}\geq\frac12$. Korn's inequality is applied to $v_{0,m+1}$ minus a lift of its normal trace on $\Sigma_{s}$. That trace is the normal trace of the correction $\omega$ in \eqref{eq:local-domain-energy-1-2}, up to sign. It involves only the levels $(0,j)$, $j\leq m$, so the lift is bounded in $H^{1}$ by $C_{0}\mathscr{E}_{1,m,\mathrm{energy},\mathrm{cl}}^{1/2}$. The other contributions are bounded by the right side of \eqref{local_energy_estimate}. This is the only place where the independence of $C_{0}$ from $\varepsilon$ is used; the absorption is then the same as for the other levels.

For the remaining mixed levels put $a=k_1\ge1$, $b=k_2\ge1$, and write $F^i=\widetilde F_{\mathrm{bet}+}^{i,a,b}$ in this step. For an admissible test $w$, put $\varphi=w\cdot\mathcal N$ on $\Sigma$. Solving \eqref{eq:local_weak_formulation} for the pressure term gives
\begin{align}
\label{eq:pressure_weak}
\int_\Omega Jp_{a,b}\func{div}_{\mathcal A_\eta}w ={}&\int_\Omega J\partial_t v_{a,b}\cdot w+ \frac\mu2\int_\Omega J\mathbb D_{\mathcal A_\eta}v_{a,b}:\mathbb D_{\mathcal A_\eta}w -\int_\Omega JF^1\cdot w \notag\\
&+\beta_0\int_{\Sigma_s}J_s(U_0v_{a,b})\cdot(U_0w) +\langle F^3,w\rangle_\Sigma+\langle J_s F^6,U_0w\rangle_{\Sigma_s} \notag\\
&+g\int_\Sigma\eta_{a,b}\varphi+ \sigma\langle a_\zeta\nabla_\Sigma\eta_{a,b}+F^5,\nabla_\Sigma\varphi\rangle_\Sigma +\langle\kappa\partial_t\eta_{a,b}+F^7,\varphi\rangle_{\Gamma_0}.
\end{align}

Take $w\in H^1_0(U;\mathbb R^3)$ on the collar $U=\Omega_{\mathrm{bet}+}$ of Lemma~\ref{lem:enhance_pressure} and extend it by zero to $\Omega$. The wall and top traces are zero, so $\varphi=0$ and all face and contact pairings vanish. Thus
\begin{equation}\label{eq:pressure_weak_0} \int_U Jp_{a,b}\func{div}_{\mathcal A_\eta}w =\int_U J\partial_t v_{a,b}\cdot w+ \frac\mu2\int_U J\mathbb D_{\mathcal A_\eta}v_{a,b}:\mathbb D_{\mathcal A_\eta}w -\int_U JF^1\cdot w .
\end{equation}
With $\widetilde w=w\circ\Phi$, $\rho=(J\circ\Phi)J_\Phi$ and $\widetilde v$ the localized rotating-frame velocity, angular integration by parts transforms the first term on the right hand side to
\begin{align}
\label{eq:pressure_time_ibp}
\int_U J\partial_t v_{a,b}\cdot w &=-\int_{\widetilde U}\rho Q\widetilde v_{a+1,b-1}\cdot\partial_\theta\widetilde w -\int_{\widetilde U}\bigl[\partial_\theta(\rho Q)\bigr] \widetilde v_{a+1,b-1}\cdot\widetilde w,
\end{align}
This implies that \(\left|\int_U J\partial_t v_{a,b}\cdot w\right| \lesssim \|w\|_{H^1(U)}\sqrt{\mathscr E_{a+1,b-1,\mathrm{energy},\mathrm{cl}}}.\) Here the chart and frame are fixed, $\partial_\theta\rho$ is bounded as in \eqref{eq:pressure_mean}, and periodicity in $\theta$ leaves no end term. The preceding energy and forcing estimates therefore bound the absolute value of the right-hand side of \eqref{eq:pressure_weak_0} by
\begin{align}
C\|w\|_{H^{1}}\Big(&\sqrt{\mathscr{E}_{k_{1}+1,k_{2}-1,\mathrm{energy}}} +\|v_{k_{1},k_{2}}\|_{H^{1}}+\sqrt{\mathcal{E}\mathcal{D}} \notag \\
&+\frac{C(\varepsilon)}{\varepsilon} \sqrt{\mathscr{E}_{k_{1},k_{2}-1,\mathrm{energy},\mathrm{cl}}} +\sqrt{\varepsilon^{\iota}\mathscr{E}_{k_{1},k_{2},\mathrm{energy},\mathrm{cl}}} +C(\varepsilon)\sqrt{\mathcal{D}_{\mathrm{int}}}\Big).
\end{align}
The last term is the cutoff contribution in Step~7. Applying Lemma~\ref{lem:enhance_pressure} at this time on $U$ gives, for $a,b\ge1$,
\begin{equation}
\label{est:highest_pressure}
\begin{aligned}
& {} \|p_{k_{1},k_{2}}\|_{L^{2}} \lesssim \|p_{k_{1},k_{2}-1}\|_{L^{2}}+\sqrt{\mathscr{E}_{k_{1}+1,k_{2}-1,\mathrm{energy}}}+\|v_{k_{1},k_{2}}\|_{H^{1}}+\sqrt{\mathcal{E}\mathcal{D}} \\
& +\frac{C(\varepsilon)}{\varepsilon}\sqrt{\mathscr{E}_{k_{1},k_{2}-1,\mathrm{energy},\mathrm{cl}}}+\sqrt{\varepsilon^{\iota}\mathscr{E}_{k_{1},k_{2},\mathrm{energy},\mathrm{cl}}}+C(\varepsilon)\sqrt{\mathcal{D}_{\mathrm{int}}} \\
& \lesssim\sqrt{\mathscr{E}_{k_{1}+1,k_{2}-1,\mathrm{energy}}}+\|v_{k_{1},k_{2}}\|_{H^{1}}+\sqrt{\mathcal{E}\mathcal{D}} +\frac{C(\varepsilon)}{\varepsilon}\sqrt{\mathscr{E}_{k_{1},k_{2}-1,\mathrm{energy},\mathrm{cl}}}+\sqrt{\varepsilon^{\iota}\mathscr{E}_{k_{1},k_{2},\mathrm{energy},\mathrm{cl}}}+C(\varepsilon)\sqrt{\mathcal{D}_{\mathrm{int}}}.
\end{aligned}
\end{equation}
Squaring \eqref{est:highest_pressure} and integrating in time bounds $\int_0^t\|p_{k_1,k_2}\|_{L^2}^2$ by the terms on the right of \eqref{local_energy_estimate}. For $a+b=m$ and $b>0$, \(\|p_{a,b}(t)\|_{L^2}^2\le\|p_{a,b}(0)\|_{L^2}^2 +\int_0^t\Big(\|p_{a,b}\|_{L^2}^2+\|p_{a+1,b}\|_{L^2}^2\Big).\) Moreover, when $k_{1}+k_{2}=m+1$, the remaining term in $\mathscr{E}_{k_{1},k_{2},\mathrm{energy},\mathrm{cl}}$ can be bounded as follows:
\begin{equation}
\begin{aligned}
& \sum_{j_1=0}^{k_1}\sum_{j_2=0}^{k_2} \|\chi_{\mathrm{bet}+,\Sigma}\partial_t^{j_1+1} \partial_\theta^{j_2}\eta\|_{H^{1/2-\alpha}}^2 \le C(\varepsilon)\sum_{j_1=0}^{k_1} \|\partial_t^{j_1+1}\eta\|_{H^{1/2-\alpha}(\Sigma)}^2 +C\sum_{j_1=0}^{k_1}\sum_{j_2=1}^{k_2} \|\chi_{\mathrm{bet}+,\Sigma}\partial_t^{j_1+1} \partial_\theta^{j_2-1}\eta\|_{H^{3/2-\alpha}}^2 \\
& \le C(\varepsilon)\bigl(\mathcal D_{||} +\mathcal E\mathcal D\bigr) +C\mathscr E_{k_1+1,k_2-1,\mathrm{energy},\mathrm{cl}} \qquad(k_2\ge1).
\end{aligned}
\end{equation}
Combining this estimate, the pressure estimate \eqref{est:highest_pressure} and \eqref{local_energy_estimate}, we have the following estimate by choosing $\varepsilon$ to be small:
\begin{equation}
\begin{aligned}
& \|v_{k_{1},k_{2}}\|^{2}(t)+\|\eta_{k_{1},k_{2}}\|_{H^{1}}^{2}(t)+\int_{0}^{t}\mathscr{E}_{k_{1},k_{2},\mathrm{energy},\mathrm{cl}} \lesssim C(\varepsilon)(\mathcal{E}(0)+\int_{0}^{t}\sqrt{\mathcal{E}}\mathcal{D}+\mathcal{E}^{\frac{3}{2}}(t)) \\
& +\int_{0}^{t}(\frac{C(\varepsilon)}{\varepsilon}\mathscr{E}_{k_{1},k_{2}-1,\mathrm{energy},\mathrm{cl}}+\mathscr{E}_{k_{1}+1,k_{2}-1,\mathrm{energy},\mathrm{cl}})+\frac{C(\varepsilon)}{\varepsilon^{\iota}}\sum_{j=0}^{k_{2}-1}\|\partial_{t}^{k_{1}}\partial_{\theta}^{j}v\|^{2}_{L^{2}}(t)+\sqrt{\mathcal{E}}^{3}(t).
\end{aligned}
\end{equation}

Similarly, the estimate above still holds when we replace $k_{1},k_{2}$ by any $j_{1},j_{2}$ such that $0\leq j_{1}\leq k_{1}$ and $0\leq j_{2}\leq k_{2}$. They can all be bounded by the right-hand side of inequality \eqref{local_energy_estimate}. Combining all of the estimates above, using the recurrence assumption for the last three terms on the right-hand side of the inequality, we have
\begin{equation}
\label{eq:induction_inequality}
{} \|v_{k_{1},k_{2}}\|^{2}(t)+\|\eta_{k_{1},k_{2}}\|_{H^{1}}^{2}(t)+\int_{0}^{t}\mathscr{E}_{k_{1},k_{2},\mathrm{energy},\mathrm{cl}} \lesssim C(\varepsilon)(\mathcal{E}(0)+\int_{0}^{t}\sqrt{\mathcal{E}}\mathcal{D}+\mathcal{E}^{\frac{3}{2}}(t)).
\end{equation}
Therefore, the recurrence estimate is closed and we have
\begin{equation}
\label{eq:induction_inequality_final}
{} \sum_{k_{1}+k_{2}=m+1}\big(\|v_{k_{1},k_{2}}\|^{2}(t)+\|\eta_{k_{1},k_{2}}\|_{H^{1}}^{2}(t)+\int_{0}^{t}\mathscr{E}_{k_{1},k_{2},\mathrm{energy},\mathrm{cl}} \big) \lesssim C(\varepsilon)(\mathcal{E}(0)+\int_{0}^{t}\sqrt{\mathcal{E}}\mathcal{D}+\mathcal{E}^{\frac{3}{2}}(t)).
\end{equation}

When $k_{1}+k_{2}\leq m$, and $k_{1}\neq m$, we have, by definition \(\mathscr{E}_{k_{1},k_{2},\mathrm{energy},\mathrm{cl}}(t)=\mathscr{E}_{k_{1},k_{2},\mathrm{energy},\mathrm{cl}}(0)+\int_0^{t}\partial_{t} \mathscr{E}_{k_{1},k_{2},\mathrm{energy},\mathrm{cl}}(t).\) For the second term, by Young's inequality, we have
\begin{equation}
\begin{aligned}
& \mathscr{E}_{k_{1},k_{2},\mathrm{energy},\mathrm{cl}}(t) \lesssim \mathscr{E}_{k_{1},k_{2},\mathrm{energy},\mathrm{cl}}(0)+\int_{0}^{t}\sum_{j_{1}=0}^{k_{1}+1}\sum_{j_{2}=0}^{k_{2}}\big((\|(\chi_{cl}+\chi_{\mathrm{bet}})\partial_{t}^{j_{1}}\partial_{\theta}^{j_{2}}v\|^{2}_{H^{1}}+\|(\chi_{cl,\Sigma}+\chi_{\mathrm{bet},\Sigma})\partial_{t}^{j_{1}}\partial_{\theta}^{j_{2}}\eta\|^{2}_{H^{\frac{3}{2}-\alpha}} \\
& +\|(\chi_{cl,\Sigma}+\chi_{\mathrm{bet},\Sigma})\partial_{t}^{j_{1}+1}\partial_{\theta}^{j_{2}}\eta\|^{2}_{H^{\frac{1}{2}-\alpha}})\big) \\
& +\int_{0}^{t}\big(\sum_{j_{1}=0}^{k_{1}+1}\sum_{j_{2}=0}^{k_{2}}(\|\partial_{t}^{j_{1}+1}\partial_{\theta}^{j_{2}}\eta\|^{2}_{L^{2}(\Gamma_{0})}+\|(\chi_{cl,\Sigma}+\chi_{\mathrm{bet},\Sigma})\partial_{t}^{j_{1}}\partial_{\theta}^{j_{2}}\eta\|^{2}_{H^{1}})+\sum_{j_{1}=0}^{k_{1}+1}\sum_{j_{2}=0}^{k_{2}}\|(\chi_{\mathrm{cl}}+\chi_{\mathrm{bet}})\partial_{t}^{j_{1}}\partial_{\theta}^{j_{2}}p\|^{2}_{L^{2}}\big) \\
& \lesssim \mathscr{E}_{k_{1},k_{2},\mathrm{energy},\mathrm{cl}}(0)+\int_{0}^{t}\mathscr{E}_{k_{1}+1,k_{2},\mathrm{energy},\mathrm{cl}}(s)ds.
\end{aligned}
\end{equation}
When $k_{1}=m$, $\mathscr{E}_{m,0,\mathrm{energy},cl}$ is bounded by $\frac{C(\varepsilon)}{\varepsilon}\mathcal{E}_{\mathrm{enhance}}$. Combining Lemma~\ref{lem:enhance_pressure} (applied at a fixed time) and Theorem \ref{thm:whole_space_estimate_final}, we obtain
\begin{align}
{\label{eq:induction_inequality_0}} &\sum_{k_{1}+k_{2}=m}\mathscr{E}_{k_{1},k_{2},\mathrm{energy},\mathrm{cl}}+\sum_{k_{1}+k_{2}=m+1}\big(\|v_{k_{1},k_{2}}\|^{2}(t)+\|\eta_{k_{1},k_{2}}\|_{H^{1}}^{2}+\int_{0}^{t}\mathscr{E}_{k_{1},k_{2},\mathrm{energy},\mathrm{cl}} \big) &  \\
\lesssim& C(\varepsilon)(\mathcal{E}(0)+\int_{0}^{t}\sqrt{\mathcal{E}}\mathcal{D}+\mathcal{E}^{\frac{3}{2}}(t)). &
\end{align}
\textbf{Final Result}

For the remaining terms included in $\mathcal{E}_{\mathrm{cl}}+\mathcal{E}_{bet}$, we use interpolation theory to derive
\begin{align}
{\label{eq:E_local_enhance}} &\sum_{k_{1}+k_{2}=m}(\|\chi_{\mathrm{cl}}\partial_{t}^{k_{1}}\partial_{\theta}^{k_{2}}v\|^{2}_{L^{2}H^{1+\frac{\epsilon_{-}}{2}}}+\|\chi_{\mathrm{cl},\Sigma}\partial_{t}^{k_{1}}\partial_{\theta}^{k_{2}}\eta\|^{2}_{L^{2}H^{\frac{3}{2}+\frac{\epsilon_{-}-\alpha}{2}}})(t) \\
&\quad+\sum_{k_{1}+k_{2}=m}(\|\chi_{\mathrm{bet}}\partial_{t}^{k_{1}}\partial_{\theta}^{k_{2}}v\|^{2}_{L^{2}H^{1+\frac{\epsilon_{-}}{2}}}+\|\chi_{\mathrm{bet},\Sigma}\partial_{t}^{k_{1}}\partial_{\theta}^{k_{2}}\eta\|^{2}_{L^{2}H^{\frac{3}{2}+\frac{\epsilon_{-}-\alpha}{2}}})(t)  \notag \\
&\lesssim \sum_{k_{1}+k_{2}=m}(\|\chi_{\mathrm{cl}}\partial_{t}^{k_{1}}\partial_{\theta}^{k_{2}}v\|^{2}_{L^{2}H^{1+\frac{\epsilon_{-}}{2}}}+\|\chi_{\mathrm{cl},\Sigma}\partial_{t}^{k_{1}}\partial_{\theta}^{k_{2}}\eta\|^{2}_{L^{2}H^{\frac{3}{2}+\frac{\epsilon_{-}-\alpha}{2}}})(0)  \notag \\
&\quad+\sum_{k_{1}+k_{2}=m}(\|\chi_{\mathrm{bet}}\partial_{t}^{k_{1}}\partial_{\theta}^{k_{2}}v\|^{2}_{L^{2}H^{1+\frac{\epsilon_{-}}{2}}}+\|\chi_{\mathrm{bet},\Sigma}\partial_{t}^{k_{1}}\partial_{\theta}^{k_{2}}\eta\|^{2}_{L^{2}H^{\frac{3}{2}+\frac{\epsilon_{-}-\alpha}{2}}})(0)  \notag \\
&\quad+\int_{0}^{t} \sum_{k_{1}+k_{2}=m}(\|\chi_{\mathrm{bet}+}\partial_{t}^{k_{1}+1}\partial_{\theta}^{k_{2}}v\|^{2}_{L^{2}H^{1}}+\|\chi_{\mathrm{bet}+,\Sigma}\partial_{t}^{k_{1}+1}\partial_{\theta}^{k_{2}}\eta\|^{2}_{L^{2}H^{\frac{3}{2}-\alpha}})  \notag \\
&\quad+\int_{0}^{t} \sum_{k_{1}+k_{2}=m}(\|\chi_{\mathrm{bet}+}\partial_{t}^{k_{1}}\partial_{\theta}^{k_{2}}v\|^{2}_{L^{2}W^{2,q_{\delta}}}+\|\chi_{\mathrm{bet}+,\Sigma}\partial_{t}^{k_{1}}\partial_{\theta}^{k_{2}}\eta\|^{2}_{L^{2}W^{3-\frac{1}{q_{\delta}},q_{\delta}}}).
\end{align}
Similarly, for the remaining terms included in $\mathcal{D}_{\mathrm{cl}}+\mathcal{D}_{\mathrm{bet}}$, we have, for any fixed $t$ and $\theta\in\mathbb{T}$, with a base point $\theta_{0}=\theta_{0}(t)\in\mathbb{T}$ chosen below,
\begin{equation}
\label{eq:D_local_enhance}
\begin{aligned}
& {} \sum_{k_{1}+k_{2}=m}(\|\chi_{\mathrm{cl},\Sigma}\partial_{t}^{k_{1}}\partial_{\theta}^{k_{2}}\eta\|^{2}_{H_{r}^{\frac{3}{2}+\frac{\epsilon_{-}-\alpha}{2}}})(t,\theta)+\sum_{k_{1}+k_{2}=m}(\|\chi_{\mathrm{bet},\Sigma}\partial_{t}^{k_{1}}\partial_{\theta}^{k_{2}}\eta\|^{2}_{H_{r}^{\frac{3}{2}+\frac{\epsilon_{-}-\alpha}{2}}})(t,\theta) \\
& \lesssim \sum_{k_{1}+k_{2}=m}(\|\chi_{\mathrm{cl},\Sigma}\partial_{t}^{k_{1}}\partial_{\theta}^{k_{2}}\eta\|^{2}_{H_{r}^{\frac{3}{2}+\frac{\epsilon_{-}-\alpha}{2}}})(t,\theta_{0})+\sum_{k_{1}+k_{2}=m}(\|\chi_{\mathrm{bet},\Sigma}\partial_{t}^{k_{1}}\partial_{\theta}^{k_{2}}\eta\|^{2}_{H_{r}^{\frac{3}{2}+\frac{\epsilon_{-}-\alpha}{2}}})(t,\theta_{0}) \\
& +\int_{\theta_{0}}^{\theta} \sum_{k_{1}+k_{2}=m}(\|\chi_{\mathrm{bet}+,\Sigma}\partial_{t}^{k_{1}}\partial_{\theta}^{k_{2}+1}\eta\|^{2}_{H_{r}^{\frac{3}{2}-\alpha}})(t,s)ds +\int_{\theta_{0}}^{\theta} \sum_{k_{1}+k_{2}=m}(\|\chi_{\mathrm{bet}+,\Sigma}\partial_{t}^{k_{1}}\partial_{\theta}^{k_{2}}\eta\|^{2}_{W_{r}^{3-\frac{1}{q_{\delta}},q_{\delta}}})(t,s)ds.
\end{aligned}
\end{equation}
We choose $\theta_{0}$ by averaging: the sum of the squared norms on the left of \eqref{eq:D_local_enhance} is integrable in $\theta$, so there is $\theta_{0}\in\mathbb{T}$ at which its value is at most its mean over $\mathbb{T}$. Therefore, combining equation \eqref{eq:D_local_enhance} and equation \eqref{eq:E_local_enhance}, we obtain
\begin{align}
{\label{eq:mixed_enhance}} &\sum_{k_{1}+k_{2}=m}(\|\chi_{\mathrm{cl},\Sigma}\partial_{t}^{k_{1}}\partial_{\theta}^{k_{2}}\eta\|^{2}_{H_{r}^{\frac{3}{2}+\frac{\epsilon_{-}-\alpha}{2}}})(t,\theta_{0})+\sum_{k_{1}+k_{2}=m}(\|\chi_{\mathrm{bet},\Sigma}\partial_{t}^{k_{1}}\partial_{\theta}^{k_{2}}\eta\|^{2}_{H_{r}^{\frac{3}{2}+\frac{\epsilon_{-}-\alpha}{2}}})(t,\theta_{0})  \notag \\
&\lesssim\sum_{k_{1}+k_{2}=m}(\|\chi_{\mathrm{cl},\Sigma}\partial_{t}^{k_{1}}\partial_{\theta}^{k_{2}}\eta\|^{2}_{L^{2}H^{\frac{3}{2}+\frac{\epsilon_{-}-\alpha}{2}}})(t)+\sum_{k_{1}+k_{2}=m}(\|\chi_{\mathrm{bet},\Sigma}\partial_{t}^{k_{1}}\partial_{\theta}^{k_{2}}\eta\|^{2}_{L^{2}H^{\frac{3}{2}+\frac{\epsilon_{-}-\alpha}{2}}})(t)  \notag \\
&\quad+ \sum_{k_{1}+k_{2}=m}(\|\chi_{\mathrm{bet}+,\Sigma}\partial_{t}^{k_{1}}\partial_{\theta}^{k_{2}+1}\eta\|^{2}_{L^{2}H^{\frac{3}{2}-\alpha}})(t)+ \sum_{k_{1}+k_{2}=m}(\|\chi_{\mathrm{bet}+,\Sigma}\partial_{t}^{k_{1}}\partial_{\theta}^{k_{2}}\eta\|^{2}_{L^{2}W^{3-\frac{1}{q_{\delta}},q_{\delta}}})(t).
\end{align}
Substituting inequality \eqref{eq:mixed_enhance} into inequality \eqref{eq:D_local_enhance}, we obtain
\begin{align}
&\sum_{k_{1}+k_{2}=m}(\|\chi_{\mathrm{cl},\Sigma}\partial_{t}^{k_{1}}\partial_{\theta}^{k_{2}}\eta\|^{2}_{L^{\infty}H_{r}^{\frac{3}{2}+\frac{\epsilon_{-}-\alpha}{2}}})(t)+\sum_{k_{1}+k_{2}=m}(\|\chi_{\mathrm{bet},\Sigma}\partial_{t}^{k_{1}}\partial_{\theta}^{k_{2}}\eta\|^{2}_{L^{\infty}H^{\frac{3}{2}+\frac{\epsilon_{-}-\alpha}{2}}})(t)  \notag \\
&\lesssim\sum_{k_{1}+k_{2}=m}(\|\chi_{\mathrm{cl},\Sigma}\partial_{t}^{k_{1}}\partial_{\theta}^{k_{2}}\eta\|^{2}_{L^{2}H^{\frac{3}{2}+\frac{\epsilon_{-}-\alpha}{2}}})(t)+\sum_{k_{1}+k_{2}=m}(\|\chi_{\mathrm{bet},\Sigma}\partial_{t}^{k_{1}}\partial_{\theta}^{k_{2}}\eta\|^{2}_{L^{2}H^{\frac{3}{2}+\frac{\epsilon_{-}-\alpha}{2}}})(t)  \notag \\
&\quad+ \sum_{k_{1}+k_{2}=m}(\|\chi_{\mathrm{bet}+,\Sigma}\partial_{t}^{k_{1}}\partial_{\theta}^{k_{2}+1}\eta\|^{2}_{L^{2}H^{\frac{3}{2}-\alpha}})(t)+ \sum_{k_{1}+k_{2}=m}(\|\chi_{\mathrm{bet}+,\Sigma}\partial_{t}^{k_{1}}\partial_{\theta}^{k_{2}}\eta\|^{2}_{L^{2}W^{3-\frac{1}{q_{\delta}},q_{\delta}}})(t).
\end{align}

Similarly, using $W^{2,q_{\delta}}_{r,h}\hookrightarrow H^{1+\epsilon_{-}}_{r,h}$, one angular derivative in $H^{1}_{r,h}$, and the same interpolation and averaging argument, we obtain
\begin{align}
&\sum_{a+b=m}\Bigl( \|\chi_{\mathrm{cl}}\partial_{t}^{a}\partial_{\theta}^{b}v(t) \|^{2}_{L^{\infty}_{\theta}H^{1+\frac{\epsilon_{-}}{2}}_{r,h}} +\|\chi_{\mathrm{bet}}\partial_{t}^{a}\partial_{\theta}^{b}v(t) \|^{2}_{L^{\infty}_{\theta}H^{1+\frac{\epsilon_{-}}{2}}_{r,h}}\Bigr) \notag \\
&\quad\le C(\varepsilon)\sum_{a+b=m}\Bigl( \mathscr E_{a,b,\mathrm{el},\mathrm{cl}}(t) +\mathscr E_{a,b,\mathrm{el},\mathrm{bet}}(t) +\mathscr E_{a,b+1,\mathrm{energy},\mathrm{cl}}(t)\Bigr).
\end{align}
Integrating these inequalities on $[0,t]$ supplies the remaining terms of $\int_0^t(\mathcal D_{\mathrm{cl}}+\mathcal D_{\mathrm{bet}})$. Their integrals are bounded by the right-hand side of equation \eqref{eq:induction_inequality_0}.

\textbf{Closing the estimate.} We now combine the estimates, in the order in which the constants are chosen (Subsection~\ref{subsec:smallness}). First $\varepsilon$ is fixed so small that the absorptions in Step~10, in Theorem~\ref{thm:top_level} and in Theorem~\ref{thm:cl_elliptic_base} hold; these use only constants independent of $\varepsilon$. From then on $C(\varepsilon) $ is a fixed constant. We use $\mathcal{E}\leq1$, so that $\mathcal{E}^{3/2}(0)\leq\mathcal{E}(0)$.

\begin{enumerate}
\item Theorems~\ref{thm:whole-space-energy} and~\ref{thm:whole_space_estimate_final} bound $\mathcal{E}_{||}(t)+\mathcal{E}_{\mathrm{enhanced}}(t)+\int_{0}^{t}(\mathcal{D}_{||}+\mathcal{D}_{\mathrm{enhanced}})$ by $C\big(\mathcal{E}(0)+\mathcal{E}^{3/2}(t)+\int_{0}^{t}\sqrt{\mathcal{E}}\mathcal{D}\big)$.

\item Theorem~\ref{thm:interior_space_final} bounds $\mathcal{E}_{\mathrm{int}}(t)+\int_{0}^{t}\mathcal{D}_{\mathrm{int}}$ by the same right-hand side with the constant $C(\varepsilon)$.

\item Theorem~\ref{thm:contact_line_energy} bounds the energy pieces $\mathscr{E}_{k_{1},k_{2},\mathrm{energy},\mathrm{cl}}$ of $\mathcal{E}_{\mathrm{cl}}+\mathcal{E}_{\mathrm{bet}}$ at time $t$, and their time integrals at the dissipation level $k_{1}+k_{2}=m+1$, by the same right-hand side; the interior dissipation it produces is bounded by \eqref{eq:Dint_E0}.

\item The elliptic pieces of $\mathcal{E}_{\mathrm{cl}}+\mathcal{E}_{\mathrm{bet}}$ at time $t$, and of $\mathcal{D}_{\mathrm{cl}}+\mathcal{D}_{\mathrm{bet}}$ at each time, are bounded by Theorems~\ref{thm:Higher-order contact line elliptic estimate} and~\ref{thm:intermediate_estimate} in terms of energy pieces at the same time, $\mathcal{E}^{2}$, $\mathcal{E}\mathcal{D}$ and $\mathcal{D}_{\mathrm{int}}$, with constants $C(\varepsilon)$. After integrating $\mathcal{D}_{\mathrm{cl}}+\mathcal{D}_{\mathrm{bet}}$ in time, \(\mathcal{E}_{\mathrm{cl}}(t)+\mathcal{E}_{\mathrm{bet}}(t)+\int_{0}^{t}(\mathcal{D}_{\mathrm{cl}}(s)+\mathcal{D}_{\mathrm{bet}}(s))ds\) is bounded by items~1--3.

\item The fractional norms of $\mathcal{E}_{\mathrm{cl}}+\mathcal{E}_{\mathrm{bet}}$ and of $\mathcal{D}_{\mathrm{cl}}+\mathcal{D}_{\mathrm{bet}}$ are bounded by \eqref{eq:E_local_enhance} and by the bound following \eqref{eq:mixed_enhance}, in terms of the quantities of items~1--4.
\end{enumerate}

Adding these bounds,
\begin{align}
\mathcal{E}(t)+\int_{0}^{t}\mathcal{D}\leq C(\varepsilon)\Big(\mathcal{E}(0)+\mathcal{E}^{\frac{3}{2}}(t)+\int_{0}^{t}\sqrt{\mathcal{E}}\mathcal{D}\Big), \label{eq:closing}
\end{align}
where $C(\varepsilon)$ collects the constants above; the proof gives $C(\varepsilon)=C\varepsilon^{-10(m+1)}C^{\prime }(\varepsilon)$, but only its finiteness is used. Increase $C(\varepsilon)$ if necessary so that it also bounds the constants in the absorptions of Lemma~\ref{lem:tl_rem}. Choose $\delta_*>0$ to satisfy the earlier smallness conditions and $\delta_*\le (2C(\varepsilon))^{-2}$. Then $\sup_{0\le s\le t}\mathcal E(s)\le\delta_*$ gives $C(\varepsilon)\sup_{0\le s\le t}\sqrt{\mathcal E(s)}\le 1/2$. Thus the two nonlinear terms in \eqref{eq:closing} are bounded by $\mathcal E(t)/2$ and $\frac12\int_0^t\mathcal D$, respectively. Moving them to the left gives
\begin{equation}
\mathcal{E}(t)+\int_{0}^{t}\mathcal{D}\leq 2C(\varepsilon)\,\mathcal{E}(0).
\end{equation}
This proves the desired a~priori implication.

\textbf{Decay.} The system is autonomous, so the argument above, started at a time $s\in[0,t]$, gives
\begin{equation*}
\mathcal{E}(t)+\int_{s}^{t}\mathcal{D}\leq C\,\mathcal{E}(s),\qquad 0\leq s\leq t<T.
\end{equation*}
Moreover $\mathcal{E}\leq C^{\prime }\mathcal{D}$ at every time, with $C^{\prime }=C_{*}$, by Lemma~\ref{lem:E_le_D}. Hence $Y(s):=\int_{s}^{T}\mathcal{E}$ satisfies $Y(s)\leq CC^{\prime }\mathcal{E}(s)=-CC^{\prime }\,Y^{\prime }(s)$, so that $Y(s)\leq Y(0)e^{-\lambda s}$ with $\lambda=1/(CC^{\prime })$. For $t\geq1$, $\mathcal{E}(t)\leq C\int_{t-1}^{t}\mathcal{E}\leq C\,Y(t-1)$; for $0\le t<\min\{1,T\}$, use $\mathcal E(t)\le C\mathcal E(0) \le Ce^\lambda e^{-\lambda t}\mathcal E(0)$; therefore $\mathcal{E}(t)\leq C^{\prime \prime }e^{-\lambda t}\,\mathcal{E}(0)$ for all $t$. This proves the Main Theorem. The constant $C(\varepsilon)\varepsilon^{-10(m+1)}$ is not tracked through the induction; only its finiteness for fixed $\varepsilon$ is used.

\subsection{The top angular level: \texorpdfstring{$\partial_\theta^{m+1}p\in L^2$ and $\partial_\theta^{m+1}\eta\in H^{3/2}$}{pressure and free surface} without loss}

\label{subsec:top_level}

In this subsection we prove Theorem~\ref{thm:top_level} and Corollary~\ref{cor:top_level}. Part~A (the pressure) is the argument of Lemma~\ref{lem:enhance_pressure} and Step~10 of the proof of Theorem~\ref{thm:contact_line_energy} at the level $(0,m+1)$, with the Ne\caronaccent{c}as inequality applied on $\Omega$ to the localized pressure extended by zero, so that its constant does not depend on $\varepsilon$. Part~B (the free surface) first puts the normal stress in $(H^{1/2}(\Sigma))^*$ and then inverts the capillary operator (Lemma~\ref{lem:tl_ell}). Since the extension near $\Gamma$ is slice-wise, $\nabla^2\partial_\theta^{m+1}\bar\eta$ is not controlled by $\|\partial_\theta^{m+1}\eta\|_{H^{3/2}}$; therefore the top-order term $\mu\func{div}_\mathcal{A}(\mathbb{D}_\mathfrak{M} v)$ is integrated by parts in all variables, and only $\partial_\theta^{m+1}\eta\in H^{3/2-\alpha}$ is used. Part~B uses Part~A only through \eqref{tl:pmain}.

\subsubsection{Notation}

In this subsection $(v,p,\eta)$ is the solution of \eqref{eq:moving-domain-elliptic-system} itself, not its localization by $\chi_{\mathrm{bet}+}$; the cutoffs are carried by the test functions. We use the notation of the previous sections: $\mathcal{A}=\mathcal{A}_\eta$, $J$, $K$ from \eqref{eq:JKA-definition}--\eqref{eq:A-explicit}; $S_\mathcal{A}$ from \eqref{eq:geometric-stress}, and $\nabla_M,\func{div}_M,\mathbb{D}_M,S_M $ for a matrix field $M$ as in \eqref{eq:F01j}; $\mathcal{N}=J\mathcal{A} N_\zeta=(-\nabla_x(\zeta+\eta),1)^T$ on $\Sigma$ by \eqref{eq:N-equality}, so $e_3\cdot\mathcal{N}=1$; $U_0=I-n_0\otimes n_0$, $\mathcal{L}_{\Sigma,0}$ from \eqref{def:linearized_operator}, $\mathcal{B}_{\mathrm{cl},0}$ from \eqref{def:contact_operator}, $\mathcal{R}$ from \eqref{eq:F05}, $\hat{\mathcal{W}}$ from \eqref{eq:W_hat}; the chart $\Phi$ of \eqref{eq:contact_chart} with $J_\Phi:=|\det D\Phi|$, the extender $E$ of Theorem~\ref{thm:half-space-extender}, the slice-wise Poisson extension $P_{2-D}$ of \eqref{eq:contact-line-extension}, acting in $(r,\mathscr H)$ with $\mathscr H=h+r\cot\theta_{\mathrm{eq}}$ (Section~\ref{sec:flattening}), and the chart representation $\mathfrak{h}_\varrho$ of \eqref{eq:local-eta-definition}. In the variables $(\theta,r,\mathscr H)$ the fluid near $\Gamma$ is $\{r>0,\mathscr H<0\}$, $\Gamma_+=\{\mathscr H=0\} $, $\Gamma_-=\{r=0\}$; since $(r,h)\mapsto(r,\mathscr H)$ is a fixed linear map, mixed norms in either pair of variables are equivalent, and we write them with the subscript $rh$. On $\Gamma_0$, $\partial_\theta|_{\Gamma_0}$ and $\partial_\gamma$ are equivalent via the following computation. Let $s$ denote planar arclength on $\partial D$ and set \(a(s):=\sqrt{1+|\nabla\zeta(\gamma(s))\cdot\dot\gamma(s)|^2}, \,\allowbreak \theta=\chi(s),\,\allowbreak d\theta=a(s)\,ds.\) For $F$ on $\partial D$, write $\widehat F(\theta)=F(\gamma(\chi^{-1}(\theta)))$. Then \(\partial_\theta\widehat F(\theta) =a(\chi^{-1}(\theta))^{-1}(\partial_\gamma F)(\gamma(\chi^{-1}(\theta))).\) Since $a$ and $a^{-1}$ are bounded by constants determined by the equilibrium geometry, it holds that
\begin{equation}
\|\partial_\theta(f|_{\Gamma_0})\|_{L^2(d\theta)} +\|\partial_\gamma(f|_{\partial D})\|_{L^2(ds)} \le C_\zeta\|f|_{\partial D}\|_{H^1(\partial D)} \le C_\zeta\|f\|_{H^{3/2}(D)}.
\end{equation}
All contact pairings obtained by integration by parts on $D$ use $ds$; equivalently,
\begin{equation}
\int_{\partial D}h\varphi\,ds =\int_0^M a(\chi^{-1}(\theta))^{-1}\widehat h(\theta)\widehat\varphi(\theta)\,d\theta.
\end{equation}
Consequently the transport estimate reads
\begin{equation}
|\kappa\langle\partial_tf,v_\gamma\partial_\gamma f\rangle_{\Gamma_0}| \le C_\zeta\kappa\|v\|_{L^\infty(\Gamma_0)} \|\partial_tf\|_{L^2(\Gamma_0)} \|\partial_\theta^{m+2}\eta\|_{L^2(\Gamma_0)} \le C(\varepsilon)\sqrt{\mathcal E}\,\mathcal D.
\end{equation}

As in \eqref{eq:volume_initial}, $\Sigma$ is identified with $D$, so $\int_\Sigma\cdot=\int_D\cdot\,dx$, $H^s(\Sigma)=H^s(D)$, and $\nu$ is the outer normal of $\partial D$. We write $\mathcal{A}_\eta=\mathfrak{A}(\nabla(W\bar\eta))$, with $\mathfrak{A}$ the matrix \eqref{eq:A-explicit} viewed as a smooth function of $\nabla(W\bar\eta)$.

Throughout, $m\ge3$ and \eqref{eq:parameters} holds; in particular $\alpha<1/12$. The letter $p$ in $W^{k,p}$ is the interior exponent of Subsection~\ref{subsec:constants}. We set $r_-:=\frac{2}{1-\varepsilon_-}$ and $s_-:=\frac1{1-\varepsilon_-}$, so that $W^{1,q_-}\subset L^{r_-}$ in two dimensions and traces of $W^{1,q_-}$ on curves lie in $L^{s_-}$.

\emph{The top level.} On the chart region set
\begin{equation*}
v_\natural:=\partial_\theta^{m+1}v,\qquad p_\natural:=\partial_\theta^{m+1}p,\qquad \eta_\natural:=\partial_\theta^{m+1}\eta,\quad f:=\chi_{\mathrm{cl},\Sigma}\,\eta_\natural=\partial_\theta^{m+1}(\chi_{\mathrm{cl},\Sigma}\eta).
\end{equation*}
Recall $\chi_{\mathrm{bet}+}=\chi_{\mathrm{cl}}+\chi_{\mathrm{bet}}$. The following norms are terms of $\mathscr{E}_{0,m+1,\mathrm{energy,cl}}$ and $\mathscr{E}_{1,m,\mathrm{energy,cl}}$ in \eqref{def:energy_cl}, hence of $\mathcal{D}_{\mathrm{cl}}$:
\begin{equation}
\label{tl:VD}
\|\chi_{\mathrm{bet}+}\partial_t\partial_\theta^mv\|_{H^1}+\|\chi_{\mathrm{bet}+} v_\nat\|_{H^1}+\|\chi_{\mathrm{bet}+,\Sigma}\eta_\nat\|_{H^{3/2-\alpha}(\Sigma)} +\|\chi_{\mathrm{bet}+,\Sigma}\partial_t\eta_\nat\|_{H^{1/2-\alpha}(\Sigma)}+\|\partial_t\eta_\nat\|_{L^2(\Gamma_0)}\le5\,\mathcal{D}^{1/2}.
\end{equation}
The levels $(1,m)$ and $(0,j)$, $j\le m$, which precede $(0,m+1)$ in the induction, are collected in $\mathscr{E}_{1,m,\mathrm{energy,cl}}$; it contains $\|\chi_{\mathrm{bet}+}\partial_t\partial_\theta^mv\|^2_{H^1}$ and $\|\chi_{\mathrm{bet}+}\partial_\theta^jv\|^2_{H^1}+\|\chi_{\mathrm{bet}+}\partial_\theta^jp\|^2_{L^2}$, $j\le m$, but neither $p_\natural$ nor a norm of $\eta_\natural$.

\emph{Forcing.} Applying $\partial_\theta^{m+1}$ to \eqref{eq:moving-domain-elliptic-system}, as in \eqref{eq:localized-after-flattening_original_chart} but for the non-localized unknowns, defines on $\mathcal{U}\cap\overline\Omega$ (see (H1))
\begin{equation}
\label{tl:forcing}
\begin{aligned}
& \widetilde F^1 :=\partial_tv_\natural+\func{div}_\mathcal{A} S_\mathcal{A}(p_\natural,v_\natural)\quad \text{in }\Omega,\quad \widetilde F^5 :=\partial_\theta^{m+1}\mathcal{R}(\nabla\eta)\quad \text{on }\Sigma, \\
& \widetilde F^3 :=S_\mathcal{A}(p_\natural,v_\natural)\mathcal{N}-\big(g\eta_\natural-\sigma\func{div}_\Sigma(\mathcal{L}_{\Sigma,0}\nabla\eta_\natural+\widetilde F^5)\big)\mathcal{N}\quad \text{on }\Sigma, \\
& \widetilde F^6 :=U_0S_\mathcal{A}(p_\natural,v_\natural)n_0-\beta_0U_0v_\natural\quad \text{on }\Sigma_s,\quad \widetilde F^7 :=-\kappa\partial_t\eta_\natural-\sigma\mathcal{B}_{\mathrm{cl},0}\eta_\natural-\sigma\widetilde F^5\cdot\nu\quad \text{on }\Gamma_0.
\end{aligned}
\end{equation}
Let
\begin{equation*}
\mathfrak{M}:=\mathfrak{A}^{\prime }\big(\nabla(W\bar\eta)\big)\big[W\,\nabla\partial_\theta^{m+1}\bar\eta\big],\quad \mathfrak{K}:=-K^2W\,\partial_3\partial_\theta^{m+1}\bar\eta,\qquad \mathcal{M}:=-\begin{pmatrix} \nabla_x\partial_\theta^{m+1}\eta \\ 0\end{pmatrix}.
\end{equation*}
These are the parts of $\partial_\theta^{m+1}\mathcal{A}_\eta$, $\partial_\theta^{m+1}K$ and $\partial_\theta^{m+1}\mathcal{N}$ that contain $\partial_\theta^{m+1}\bar\eta$ (resp.\ $\partial_\theta^{m+1}\eta$) together with one more derivative. The terms of the Leibniz expansion in which all $m+1$ angular derivatives fall on $\eta$ are
\begin{equation}  \label{tl:top}
\begin{aligned} \widetilde F^1_{\rm top}&:=\mu\,\operatorname{div}_\cA(\mathbb{D}_\fM v)-\operatorname{div}_\fM S_\cA(p,v)+W\mathfrak{K}\,\partial_t\bar\eta\,\partial_3v+WK\,\partial_t\partial_\theta^{m+1}\bar\eta\,\partial_3v-v\cdot\nabla_\fM v,\\ \widetilde F^3_{\rm top}&:=\mu\,\mathbb{D}_\fM v\,\mathcal{N}-S_\cA(p,v)\,\mathcal M+\big(g\eta-\sigma[H(\zeta+\eta)-H(\zeta)]\big)\mathcal M,\\ \widetilde F^5_{\rm top}&:=\mathcal{R}'(\nabla\eta)\,\nabla\eta_\nat,\qquad \widetilde F^6_{\rm top}:=\mu\,U_0\mathbb{D}_\fM v\,n_0,\qquad \widetilde F^7_{\rm top}:=\,\hat{\mathcal W}'(\partial_t\eta)\,\partial_t\eta_\nat.
\end{aligned}
\end{equation}
We put $\widetilde F^i_{\mathrm{rem}}:=\widetilde F^i-\widetilde F^i_{\mathrm{top}}$. Two facts are used repeatedly:

\begin{itemize}
\item $e_3\cdot\mathcal{M}=0$;

\item $e_3\cdot S_\mathcal{A}(p,v)\mathcal{M}=-\mu\sum_{j=1,2}(\mathbb{D}_\mathcal{A} v)_{3j}\mathcal{M}_j$.
\end{itemize}

Hence the pressure and the capillary factor drop out of the vertical component of $\widetilde F^3_{\mathrm{top}}$:
\begin{equation}  \label{tl:F3vert} e_3\cdot\widetilde F^3_{\mathrm{top}}=\mu\,e_3\cdot\mathbb{D}_\mathfrak{M} v\,\mathcal{N}+\mu\sum_{j=1,2}(\mathbb{D}_\mathcal{A} v)_{3j}\,\mathcal{M}_j .
\end{equation}
Finally we split $\widetilde F^i_{\mathrm{rem}}=\widetilde F^i_{\mathrm{lin}}+\widetilde F^i_{\mathrm{nl}}$, where $\widetilde F^i_{\mathrm{lin}}$ is the part of $\widetilde F^i$ that is linear in $(v,p,\eta)$. Every term of $\widetilde F^i_{\mathrm{top}}$ is at least bilinear, so $\widetilde F^i_{\mathrm{lin}}$ is the linear part of $\widetilde F^i_{\mathrm{rem}}$. It is obtained by setting $\eta=0$ in the coefficients, that is, $\mathcal{A}=I$, $J=K=1$, $\mathcal{N}=N_\zeta$ and $\mathcal{R}=\hat{\mathcal{W}}=0$; then only the commutators of $\partial_\theta^{m+1}$ with the operators of the linearized problem remain. For instance \(\widetilde F^1_{\mathrm{lin}}=-\big[\partial_\theta^{m+1},\func{div} S_I\big](p,v),\,\allowbreak S_I(p,v):=pI-\mu\mathbb{D} v,\) and $\widetilde F^5_{\mathrm{lin}}=0$. The coefficients of these commutators are smooth and depend only on $\zeta$ and on the chart. The rest $\widetilde F^i_{\mathrm{nl}}$ is at least quadratic.

\subsubsection{Hypotheses}

\begin{enumerate}[label=(H\arabic*),leftmargin=3em]

\item \textbf{Geometry and cutoffs.} The side wall is vertical near $\Sigma$: for some $c_0>0$, $\Omega\cap\{z>\zeta(x)-c_0\}=\{(x,z):x\in D,\ \zeta(x)-c_0<z<\zeta(x)\}$. Thus $n_0=(\nu,0)$ and $e_3\cdot n_0=0$ there. The flattening is a vertical displacement, so $J\mathcal{A}$ is the cofactor matrix of $D\Phi_\eta$ and $J\mathcal{A} n_0=Jn_0$. Hence $J_s=J$ there, not $1$; only the boundedness of $J_s$ is used. $\mathcal{U}$ is a neighbourhood of $\Gamma$ in $\overline\Omega$ that is $\theta$-independent in the chart, lies in $\Phi(\mathbb{T}\times\mathcal{W}_{\mathrm{cl}})$, in the region of vertical wall and in $\{z>h_0\}$, so that $\chi_{\mathrm{bot}}=0$ on $\mathcal{U}$. The cutoffs $\chi_{\mathrm{cl}}$, $\chi_{\mathrm{cl},\Sigma}=\chi_{\mathrm{cl}}|_\Sigma$ and $\chi_{\mathrm{bet}+}=\chi_{\mathrm{cl}}+\chi_{\mathrm{bet}}$ are functions of $(r,h)$ only, with values in $[0,1]$, $\func{supp}\chi_{\mathrm{cl}}\subset\func{supp}\chi_{\mathrm{bet}+}\Subset\mathcal{U}$, $\chi_{\mathrm{cl},\Sigma}=1$ near $\Gamma_0$, $\chi_{\mathrm{bet}+}=1$ on $\func{supp}\chi_{\mathrm{cl}}$, and $\chi_{\mathrm{bet}+}=1-\chi_{\mathrm{int}}$ on $\mathcal{U}$. They satisfy \eqref{eq:cutoff_eps} and \eqref{eq:cutoff_scaling}; in particular
\begin{equation}  \label{tl:cutoff} |\nabla\chi_{\mathrm{bet}+}|\le C_1\varepsilon^{-1},\qquad \func{dist}(x,\partial\Omega)\le C_1\varepsilon\quad\text{for }x\in\func{supp}\chi_{\mathrm{bet}+},
\end{equation}
with $C_1$ independent of $\varepsilon$. The transition region $\func{supp}\chi_{\mathrm{bet}+}\setminus\{\chi_{\mathrm{bet}+}=1\}$ lies in $\func{supp}\chi_{\mathrm{int}}$, where $\chi_{0,\mathrm{int}}=\chi_{\mathrm{int}+}=1$ by the nesting of Subsection~\ref{subsec:cutoff}. Finally $\func{supp}\chi_{\mathrm{cl}}\cap\func{supp}\chi_{\mathrm{int}+}=\emptyset$, because $\chi_{\mathrm{cl}}\chi_{0,\mathrm{int}}=0$ and $\chi_{0,\mathrm{int}}=1$ on $\func{supp}\chi_{\mathrm{int}+}$. Since $\chi_{\mathrm{bet}+}\le1$ attains its maximum on $\func{supp}\chi_{\mathrm{cl}}$, also $\nabla\chi_{\mathrm{bet}+}=0$ there.

\item \textbf{Extension near $\Gamma$.} On $\mathcal{U}$, in the variables $(\theta,r,\mathscr H)$,
\begin{equation}  \label{tl:ext} \bar\eta\circ\Phi=\sum_{i=0}^{N}\psi_i\,\big(P_{2-D}E\,\mathfrak{h}_{\chi_i\eta}\big)+\bar\eta_{\mathrm{int}}\circ\Phi .
\end{equation}
Here $\psi_i$ and $\chi_i$ are $\theta$-independent cutoffs; near $\Gamma_0$ we have $\psi_0=\chi_0=1$ and $\psi_i=0$ for $i\ge1$; and each $\chi_i$ has the form $\chi_i=\varphi_i\,\chi_{\mathrm{bet}+,\Sigma}$ with $\varphi_i$ smooth and $\theta$-independent. The last term is the $3$-D Poisson piece $\bar\eta_{\mathrm{int}}=\chi_{\mathrm{int}+}\mathcal{P}\circ E_{D_{\mathrm{int}}}[\chi_{\mathrm{int},\Sigma}\eta]$ of Section~\ref{sec:extension}; it vanishes outside $\func{supp}\chi_{\mathrm{int}+}$. By \eqref{eq:contact-line-extension}, \eqref{eq:interior-extension} and \eqref{eq:eta_bar_def}, \eqref{tl:ext} holds with $N=m+2$: the piece $i=0$ is $\bar\eta_{\mathrm{cl}}$, with $\chi_0=\chi_{\mathrm{cl},\Sigma}$; the pieces $1\le i\le m+1$ are the $P_{2-D}\circ E$ pieces of $\bar\eta_{m+1,\mathrm{int}}$ with $\chi_i=\chi_{\mathrm{dif},i-1,\Sigma}$; and the piece $i=m+2$ has $\chi_{m+2}=\chi_{0,\mathrm{int},\Sigma}-\chi_{\mathrm{int},\Sigma}$. The form $\chi_i=\varphi_i\chi_{\mathrm{bet}+,\Sigma}$ holds for each of them:

\begin{itemize}
\item $\chi_{\mathrm{cl},\Sigma}=\chi_{\mathrm{cl},\Sigma}\,\chi_{\mathrm{bet}+,\Sigma}$, since $\chi_{\mathrm{bet}+}=1$ on $\func{supp}\chi_{\mathrm{cl}} $;

\item $\chi_{\mathrm{dif},j,\Sigma}=\chi_{\mathrm{dif},j,\Sigma}\,\chi_{\mathrm{bet}+,\Sigma}$: where $\chi_{j,\mathrm{int},\Sigma}=1$ also $\chi_{j+1,\mathrm{int},\Sigma}=1$, so $\chi_{\mathrm{dif},j,\Sigma}$ vanishes there; the open set $\{\chi_{\mathrm{int},\Sigma}>0\}$ lies in $\{\chi_{j,\mathrm{int},\Sigma}=1\}$, hence is disjoint from $\func{supp}\chi_{\mathrm{dif},j,\Sigma}$, and $\chi_{\mathrm{bet}+,\Sigma}=1-\chi_{\mathrm{int},\Sigma}=1$ on that support;

\item $\chi_{0,\mathrm{int},\Sigma}-\chi_{\mathrm{int},\Sigma}=\chi_{0,\mathrm{int},\Sigma}(1-\chi_{\mathrm{int},\Sigma})=\chi_{0,\mathrm{int},\Sigma}\,\chi_{\mathrm{bet}+,\Sigma}$, since $\chi_{0,\mathrm{int},\Sigma}=1 $ on $\func{supp}\chi_{\mathrm{int},\Sigma}$.
\end{itemize}

The derivatives of $\varphi_i$ of order $k$ are bounded by $C\varepsilon^{-k} $.

\item \textbf{A priori regularity.} At the time $t$ considered:

\begin{itemize}
\item all norms in $\mathcal{E}(t)$ and $\mathcal{D}(t)$ are finite;

\item $f(t)\in H^{3/2}(\Sigma)$;

\item $v_\natural\in H^1(\mathcal{U}\cap\Omega)$, $\partial_tv_\natural\in L^2(\mathcal{U}\cap\Omega)$ and $p_\natural\in L^2(\mathcal{U}\cap\Omega)$;

\item for every \emph{admissible} $w$, that is, $w\in H^1(\Omega)$ with $\func{supp}w\subset\mathcal{U}$, $w\cdot n_0=0$ on $\Sigma_s$ and $w\cdot\mathcal{N}|_\Sigma\in H^1(\Sigma)$,
\begin{equation}
\label{tl:weak}
\begin{aligned}
& \int_\Omega J\,\partial_tv_\nat\cdot w-\int_\Omega J\,S_\cA(p_\nat,v_\nat):\nabla_\cA w+\beta_0\int_{\Sigma_s}J_sU_0v_\nat\cdot U_0w +g\int_\Sigma\eta_\nat\,(w\cdot\mathcal{N})+\sigma\int_\Sigma\mathcal{L}_{\Sigma,0}\nabla\eta_\nat\cdot\nabla(w\cdot\mathcal{N}) \\
& =\int_\Omega J\widetilde F^1\cdot w-\langle\widetilde F^3,w\rangle_\Sigma-\langle J_s\widetilde F^6,U_0w\rangle_{\Sigma_s}-\sigma\langle\widetilde F^5,\nabla(w\cdot\mathcal{N})\rangle_\Sigma-\langle\kappa\partial_t\eta_\nat+\widetilde F^7,w\cdot\mathcal{N}\rangle_{\Gamma_0}.
\end{aligned}
\end{equation}
\end{itemize}

For a regular solution, combine \eqref{tl:forcing} with Green's formula and the Piola identity. On the top face the unit reference normal is $n_\Sigma=N_0/|N_0|$, so \eqref{eq:transformed-normal} gives $J\mathcal A n_\Sigma=\mathcal N/|N_0|$. Since $dS_{\mathrm{ref}}=|N_0|\,dx$, the boundary flux in \eqref{tl:weak} pairs $\mathcal N$ with $dx$. On the wall $J\mathcal A n_0=J_s n_0$. This is the convention of \eqref{eq:whole-domain-weak-momentum-j} and \eqref{eq:local_weak_formulation}. No cutoff commutators occur here, because the cutoff sits on $w$.

\item \textbf{Control on the region $\{\chi_{0,\mathrm{int}}=1\}$.} On $\mathcal{U}\cap\{\chi_{0,\mathrm{int}}=1\}$, which contains $\func{supp}\nabla\chi_{\mathrm{bet}+}$ and $\func{supp}\chi_{\mathrm{int}+}\cap\mathcal{U}$ by (H1),
\begin{equation}  \label{tl:Dint}
\begin{aligned}
&\|\partial_\theta^{m+1}p\|_{L^{6/5}}+\|\partial_\theta^{m+1}v\|_{L^2}+\sum_{j=0}^{m}\big(\|\partial_\theta^jv\|_{H^1}+\|\partial_\theta^jp\|_{L^2}\big)\lesssim\mathcal{D}_{\rm int}^{1/2},\\
&\|\partial_\theta^{m+1}\bar\eta_{\rm int}\|_{H^1(\cU)}+\|\partial_t\partial_\theta^{m+1}\bar\eta_{\rm int}\|_{L^2(\cU)}\le C(\varepsilon)\,\mathcal{D}_{\rm int}^{1/2}.
\end{aligned}
\end{equation}
The first line follows from the terms $\|\chi_{0,\mathrm{int}}v\|_{W^{m+2,p}} $ and $\|\chi_{0,\mathrm{int}}p\|_{W^{m+1,p}}$ ($k=n=0$) of \eqref{def:interior_dissipation}. Since $\chi_{0,\mathrm{int}}$ does not depend on $\theta$, $\partial_\theta^{m+1}v=\partial_\theta^{m+1}(\chi_{0,\mathrm{int}}v)$ on $\{\chi_{0,\mathrm{int}}=1\}$, and similarly for the other quantities. Moreover $\partial_\theta$ is a first-order operator with smooth coefficients on $\mathcal{U}$, and $W^{1,p}(\Omega)\subset L^2(\Omega) $ in three dimensions, as $p\ge6/5$. The embedding is applied on the fixed domain $\Omega$, so the constant does not depend on $\varepsilon$. For the second line, $\chi_{\mathrm{int},\Sigma}\eta=\chi_{\mathrm{int},\Sigma}\chi_{0,\mathrm{int},\Sigma}\eta$ and $\chi_{\mathrm{int},\Sigma}\partial_t\eta=\chi_{\mathrm{int},\Sigma}\chi_{1,\mathrm{int},\Sigma}\partial_t\eta$ lie in $W^{m+3-1/p,p}(\Sigma)$ and $W^{m+2-1/p,p}(\Sigma)$ with norms $\le C(\varepsilon)\mathcal{D}_{\mathrm{int}}^{1/2}$ (terms $k=n=0$ and $k=n=1$). The Poisson extension maps them into $W^{m+3,p}$ and $W^{m+2,p}$ near $\Sigma$, so $\partial_\theta^{m+1}\bar\eta_{\mathrm{int}}\in W^{2,p}\subset H^1$ and $\partial_t\partial_\theta^{m+1}\bar\eta_{\mathrm{int}}\in W^{1,p}\subset L^2$.
\end{enumerate}

\subsubsection{Auxiliary lemmas}

All norms of $\bar\eta$, $v$, $p$ below are taken on $\mathcal{U}$. In chart variables $dx=J_\Phi\,d\theta\,dr\,dh$, and $J_\Phi$ is bounded above and below on $\mathcal{U}$.

\begin{lemma}[Coefficient bounds]
\label{lem:tl_coef} If $\mathcal{E}\le1$, then on $\mathcal{U}$
\begin{align*}
&\|\nabla v\|_{L^\infty_\theta L^{r_-}_{rh}}+\|\nabla^2v\|_{L^\infty_\theta L^{q_-}_{rh}}+\|\nabla p\|_{L^\infty_\theta L^{q_-}_{rh}}+\|v\|_{L^\infty}+\|\nabla v|_{\Gamma_\pm}\|_{L^\infty_\theta L^{s_-}}\lesssim\sqrt{\mathcal{E}}, \\
&\|\nabla\bar\eta\|_{L^\infty}+\|\partial_\theta\nabla\bar\eta\|_{L^\infty}+\|\nabla^2\bar\eta\|_{L^\infty_\theta L^{r_-}_{rh}}+\|\partial_t\bar\eta\|_{L^\infty}\lesssim\sqrt{\mathcal{E}}, \\
&\|\nabla\eta\|_{L^\infty(\Sigma\cap\mathcal{U})}+\|\partial_\theta\nabla\eta\|_{L^\infty(\Sigma\cap\mathcal{U})}+\sup_\theta\|\partial_r\nabla\eta(\theta,\cdot)\|_{L^{s_-}_r}+\|\partial_t\eta\|_{L^\infty(\Gamma_0)}\lesssim\sqrt{\mathcal{E}}.
\end{align*}
In particular $|J-1|+|K-1|+|\mathcal{A}-I|\lesssim\sqrt\mathcal{E}$ and $|\partial_\theta J|+|\partial_\theta\mathcal{A}|\lesssim1$ on $\mathcal{U}$.
\end{lemma}

The constants may depend on $\varepsilon$ through the cutoffs; this is harmless, since every bound of the lemma is used either with a factor that is not absorbed or under $\mathcal{E}\le\delta_0(\varepsilon)$.

\begin{proof}
$\mathcal{E}_{\mathrm{cl}}+\mathcal{E}_{\mathrm{bet}}$ contains $\mathscr{E}_{k_1,k_2,\mathrm{el,cl}}$ and $\mathscr{E}_{k_1,k_2,\mathrm{el,bet}}$ with $k_1+k_2=m-1\ge2$. These give $\partial_t^{j_1}\partial_\theta^jv\in L^2_\theta W^{2,q}_{rh}$, $\partial_t^{j_1}\partial_\theta^jp\in L^2_\theta W^{1,q}_{rh}$ and $\partial_t^{j_1}\partial_\theta^j\eta\in L^2_\theta W^{3-1/q,q}_r$ for $j_1+j\le m-1$, with $q\in\{q_-,q_+,p\}$, all $\ge q_-$. On the bounded slice $\mathcal{W}_{\mathrm{cl}}$ the bounds therefore hold with $q_-$. Only $j_1+j\le2$ is used, together with the energy bound $\partial_\theta^3v\in L^2$ for the term $\partial_\theta^3v$ of $\partial_\theta\nabla^2v$ (note $L^2\subset L^{q_-}$ on bounded slices). The estimates follow from:

\begin{itemize}
\item $H^1_\theta\subset L^\infty_\theta$;

\item in two dimensions, $W^{2,q_-}\subset C^0$ and $W^{1,q_-}\subset L^{r_-} $;

\item the trace $W^{1,q_-}(\mathcal{W}_{\mathrm{cl}})\to L^{s_-}(\partial\mathcal{W}_{\mathrm{cl}})$;

\item in one dimension, $W^{3-1/q_-,q_-}\subset W^{2,s_-}\subset W^{1,\infty} $, because $1/s_-=2/q_--1$;

\item the boundedness $W^{3-1/q,q}(\mathbb{R})\to W^{3,q}$ of $P_{2-D}E$ on bounded slices, together with (H2).
\end{itemize}

These are the bounds \eqref{eq:dtJ_Linfty}, \eqref{est:background_geometry} and \eqref{eq:convection_cl}.
\end{proof}

\begin{lemma}[The extension at the top level]
\label{lem:tl_slice} Let $U:=\partial_\theta^{m+1}\bar\eta$ on $\mathcal{U}$, and split $U=U_{\mathrm{sl}}+U_3$ according to \eqref{tl:ext}:
\begin{equation*}
U_{\mathrm{sl}}:=\sum_{i=0}^N\psi_iU_i,\qquad U_i:=P_{2-D}E\,\partial_\theta^{m+1}\mathfrak{h}_{\chi_i\eta},\qquad U_3:=\partial_\theta^{m+1}\bar\eta_{\mathrm{int}}.
\end{equation*}
Then $U_3=0$ outside $\func{supp}\chi_{\mathrm{int}+}$, and
\begin{align}
&\|U_{\mathrm{sl}}\|_{L^2_\theta H^1_{rh}(\mathcal{U})}+\|\partial_\theta U_{\mathrm{sl}}\|_{L^2(\mathcal{U})}\le C(\varepsilon)\,\mathcal{D}^{1/2}, \label{tl:s1} \\
&\|\nabla U_{\mathrm{sl}}\|_{L^2_\theta H^{1-\alpha}_{rh}(\mathcal{U})}+\|\partial_tU_{\mathrm{sl}}\|_{L^2_\theta H^{1-\alpha}_{rh}(\mathcal{U})}\le C(\varepsilon)\,\mathcal{D}^{1/2},  \label{tl:s2} \\
&\|U_3\|_{H^1(\mathcal{U})}+\|\partial_tU_3\|_{L^2(\mathcal{U})}\le C(\varepsilon)\,\mathcal{D}_{\mathrm{int}}^{1/2}.  \label{tl:s4}
\end{align}
\end{lemma}

\begin{proof}
\emph{Fourier bounds.} Write $\langle\cdot\rangle=(1+|\cdot|^2)^{1/2}$, with $k$ dual to $\theta$ and $\xi$ dual to $r$. Put $\mathbb H_1=\mathbb R\times(-1,0)$. The radial Fourier multiplier of $P_{2-D}$ is $e^{\mathscr H|\xi|}$, and
\begin{equation}
\int_{-1}^0 e^{2\mathscr H|\xi|}\,d\mathscr H =\frac{1-e^{-2|\xi|}}{2|\xi|}\lesssim\langle\xi\rangle^{-1},
\end{equation}
with value $1$ at $\xi=0$. Each derivative in $r$ or $\mathscr H$ contributes a factor of modulus $|\xi|$. Thus Plancherel's theorem, summing the derivatives of orders zero through $j$, gives
\begin{equation}  \label{tl:P2} \|P_{2-D}G\|_{L^2_\theta H^j_{r\mathscr H}(\mathbb H_1)} \lesssim_j \|G\|_{L^2_\theta H^{j-1/2}_r}, \qquad j=0,1,2,\ldots .
\end{equation}
Interpolating between consecutive integer orders on the strip gives $\|P_{2-D}G\|_{L^2_\theta H^{s+1/2}_{r\mathscr H}(\mathbb H_1)} \lesssim_s\|G\|_{L^2_\theta H^s_r}$ for every $s\ge-1/2$. The extender $E$ is bounded on $H^s_r$, $0\le s\le2$, by Theorem~\ref{thm:half-space-extender}, hence on $H^{1/2}(\mathbb{T}\times\mathbb{R}_+)$ and $H^{3/2}(\mathbb{T}\times\mathbb{R}_+)$. The pointwise inequalities
\begin{equation*}
\langle\xi\rangle^{3/2},\ \ |k|\langle\xi\rangle^{1/2}\ \le\ \langle(k,\xi)\rangle^{3/2},\qquad \langle\xi\rangle^{1/2}\le\langle(k,\xi)\rangle^{1/2},\quad |k|\langle\xi\rangle^{-1/2}\le\langle(k,\xi)\rangle
\end{equation*}
give, by Plancherel's theorem and addition of the two individual bounds, for $G$ on $\mathbb{T}\times\mathbb{R}$,
\begin{equation}
\label{tl:aniso}
\|G\|_{\Lt H^{3/2}_r}+\|\partial_\theta G\|_{\Lt H^{1/2}_r}\le 2\|G\|_{H^{3/2}},\qquad \|G\|_{\Lt H^{1/2}_r}\le\|G\|_{H^{1/2}},\quad \|\partial_\theta G\|_{\Lt H^{-1/2}_r}\le\|G\|_{H^1},
\end{equation}
and the first two with $3/2,1/2$ replaced by $3/2-\alpha,1/2-\alpha$. The norm equivalence \eqref{eq:localized-heta-H1} transfers these bounds between $\Sigma$ and the chart.

\emph{Proof of the bounds.} By (H2), $\chi_i=\varphi_i\chi_{\mathrm{bet}+,\Sigma}$ with $\varphi_i$ independent of $\theta$, so
\begin{equation*}
\partial_\theta^{m+1}(\chi_i\eta)=\varphi_i\,\chi_{\mathrm{bet}+,\Sigma}\,\partial_\theta^{m+1}\eta,\quad \partial_t\partial_\theta^{m+1}(\chi_i\eta)=\varphi_i\,\chi_{\mathrm{bet}+,\Sigma}\,\partial_t\partial_\theta^{m+1}\eta .
\end{equation*}
Multiplication by $\varphi_i$ is bounded on $H^{3/2-\alpha}(\Sigma)$ and on $H^{1/2-\alpha}(\Sigma)$ with norm $\allowbreak\le C(\varepsilon)$. Hence, by \eqref{tl:VD}, $\allowbreak\|\partial_\theta^{m\allowbreak+1}(\chi_i\eta)\|_{H^{3/2-\alpha}(\Sigma)}\allowbreak+\allowbreak\|\partial_t\partial_\theta^{m\allowbreak+1}(\chi_i\eta)\|_{H^{1/2-\alpha}(\Sigma)}\allowbreak\le C(\varepsilon)\mathcal{D}^{1/2}$. Then \eqref{tl:s1} and \eqref{tl:s2} follow from \eqref{tl:P2}, \eqref{tl:aniso} and their $\alpha$-versions. For the $\partial_\theta$-bound in \eqref{tl:s1} we use the last bound in \eqref{tl:aniso} and $\|\cdot\|_{H^1}\allowbreak\le\|\cdot\|_{H^{3/2-\alpha}}$, as $\alpha<1/2$. The bound \eqref{tl:s4} is the second line of \eqref{tl:Dint}, and $U_3$ vanishes outside $\func{supp}\chi_{\mathrm{int}\allowbreak+}$ because $\bar\eta_{\mathrm{int}}$ does.
\end{proof}

\begin{lemma}[The top terms]
\label{lem:tl_first} Let $w$ be admissible with $\func{supp}w\subset\func{supp}\chi_{\mathrm{bet}+}$ and $\|w\|_{H^1(\Omega)}\le1$. Split $\widetilde F^1_{\mathrm{top}}=\mu\func{div}_\mathcal{A}(\mathbb{D}_\mathfrak{M} v)+\mathbb{L}$, so that $\mathbb{L}$ collects the other four terms of \eqref{tl:top}. Then, if $\mathcal{E}\le1$,
\begin{equation*}
\Big|\int_\Omega J\,\mathbb{L}\cdot w\Big|+\Big|\int_\Omega J\,\mathbb{D}_\mathfrak{M} v:\nabla_\mathcal{A} w\Big|\le C(\varepsilon)\sqrt{\mathcal{E}}\,\mathcal{D}^{1/2}.
\end{equation*}
If moreover $\func{supp}w\subset\func{supp}\chi_{\mathrm{cl}}$, then also
\begin{equation*}
\big|\langle e_3\cdot\widetilde F^3_{\mathrm{top}},\,w_3\rangle_\Sigma\big|+\big|\langle\mu\,e_3\cdot\mathbb{D}_\mathfrak{M} v\,\mathcal{N},\,w_3\rangle_\Sigma\big| +\big|\langle J_s\widetilde F^6_{\mathrm{top}},U_0w\rangle_{\Sigma_s}\big|\le C(\varepsilon)\sqrt{\mathcal{E}}\,\mathcal{D}^{1/2}.
\end{equation*}
\end{lemma}

\begin{proof}
Every term is trilinear: one factor among $\nabla U$, $U$ or $\partial_tU$ (by Lemma~\ref{lem:tl_coef}, $|\mathfrak{M}|+|\mathfrak{K}|\lesssim|\nabla U|+|U|$), one coefficient bounded in Lemma~\ref{lem:tl_coef}, and $w$, $\nabla w$ or a trace of $w$. Split $U=U_{\mathrm{sl}}+U_3$ as in Lemma~\ref{lem:tl_slice}.

\emph{The $3$-D piece.} $U_3$ vanishes outside $\func{supp}\chi_{\mathrm{int}+}$, where $\chi_{1,\mathrm{int}}=1$ by the nesting relation. There the coefficients are bounded: the terms $\|\chi_{1,\mathrm{int}}v\|_{W^{m+1,p}}$ and $\|\chi_{1,\mathrm{int}}p\|_{W^{m,p}}$ of $\mathcal{E}_{\mathrm{int}}$ ($k=n=0$ in \eqref{def:interior_energy}) and $W^{m+1,p}\subset W^{2,\infty}$, $W^{m,p}\subset W^{1,\infty}$ (as $m\ge3$ and $p>3/2$) give $\|\nabla^2v\|_{L^\infty}+\|\nabla v\|_{L^\infty}+\|v\|_{L^\infty}+\|\nabla p\|_{L^\infty}\lesssim\sqrt\mathcal{E}$ on $\func{supp}\chi_{\mathrm{int}+}$, and $\nabla\bar\eta$, $\partial_t\bar\eta$ are bounded by Lemma~\ref{lem:tl_coef}. Hence the bulk terms with $U_3$ in place of $U$ are $\lesssim\sqrt\mathcal{E}(\|U_3\|_{H^1}+\|\partial_tU_3\|_{L^2})\le C(\varepsilon)\sqrt\mathcal{E}\,\mathcal{D}_{\mathrm{int}}^{1/2}$ by \eqref{tl:s4}. If $\func{supp}w\subset\func{supp}\chi_{\mathrm{cl}}$, then $U_3=0$ on $\func{supp}w$ by (H1), so $U_3$ does not enter the face terms.

\emph{The slice-wise pieces.} In two dimensions $H^{1-\alpha}\subset L^{2/\alpha}$, $H^1\subset L^s$ for $s<\infty$, and on curves $H^{1/2-\alpha}\subset L^{1/\alpha}$ and $H^{1/2}\subset L^s$. The factors built from $U_{\mathrm{sl}}$ are bounded by $C(\varepsilon)\mathcal{D}^{1/2}$, by \eqref{tl:s1} and \eqref{tl:s2}, on all of $\mathcal{U}$.

\emph{Bulk terms.} In $\mathbb{L}$ the coefficient is $\nabla^2v$, $\nabla p$, $\nabla\mathcal{A}\,\nabla v$, $v\,\nabla v$, $\partial_t\bar\eta\,\nabla v $ or $\nabla v$, so it lies in $L^\infty_\theta L^{q_-}_{rh}$. H\"older in $\theta$ ($L^2\cdot L^\infty\cdot L^2$) and in $(r,h)$ gives \(\frac\alpha2+\frac{2-\varepsilon_-}2+\frac1s=1\,\allowbreak\text{with }\ s=\frac{2}{\varepsilon_--\alpha}<\infty .\) In $\int J\mathbb{D}_\mathfrak{M} v:\nabla_\mathcal{A} w$ the coefficient is $\nabla v\in L^\infty_\theta L^{r_-}_{rh}$ and the test factor is $\nabla w\in L^2_\theta L^2_{rh}$. Since $\frac\alpha2+\frac{1-\varepsilon_-}2+\frac12<1$, H\"older applies on the bounded slices.

\emph{Terms on $\Sigma$ and $\Sigma_s$.} By \eqref{tl:F3vert} they are $\mu\,e_3\cdot\mathbb{D}_\mathfrak{M} v\,\mathcal{N}$, $\mu\sum_{j=1,2}(\mathbb{D}_\mathcal{A} v)_{3j}\mathcal{M}_j$ and $\widetilde F^6_{\mathrm{top}}=\mu U_0\mathbb{D}_\mathfrak{M} v\,n_0$; the pressure does not occur. The traces of $\nabla U$ and $U$ on $\Gamma_\pm$ lie in $L^2_\theta H^{1/2-\alpha}\subset L^2_\theta L^{1/\alpha}$, by \eqref{tl:s1}, \eqref{tl:s2} and the trace theorem. Also $\nabla_x\partial_\theta^{m+1}\eta=\nabla_x(\chi_{\mathrm{bet}+,\Sigma}\partial_\theta^{m+1}\eta)$ on $\Sigma\cap\func{supp}\chi_{\mathrm{cl},\Sigma}\supset\Sigma\cap\func{supp}w$, by (H1), and it lies in $H^{1/2-\alpha}$ by \eqref{tl:VD}. The coefficient is $\nabla v|_{\Gamma_\pm}\in L^\infty_\theta L^{s_-}$. H\"older gives $(1-\varepsilon_-)+\alpha+\frac1s=1$ with $s=1/(\varepsilon_--\alpha)$, and $w|_{\Gamma_\pm}\in L^2_\theta H^{1/2}\subset L^2_\theta L^s$.
\end{proof}

\begin{lemma}{\label{lem:multiplier_estimate}}[A multiplier on $H^{1/2}$] \label{lem:tl_mult} Let $1<s<2$, and let $a=a(\theta,r)$ satisfy $a,\partial_\theta a\in L^\infty(\mathbb{T}\times\mathbb{R})$ and $\sup_\theta\|\partial_ra(\theta,\cdot)\|_{L^s(\mathbb{R})}<\infty$. Then
\begin{equation*}
\|aF\|_{H^{1/2}(\mathbb{T}\times\mathbb{R})}\lesssim\Big(\|a\|_{L^\infty}+\|\partial_\theta a\|_{L^\infty}+\sup_\theta\|\partial_ra(\theta,\cdot)\|_{L^s}\Big)\|F\|_{H^{1/2}(\mathbb{T}\times\mathbb{R})}.
\end{equation*}
\end{lemma}

\begin{proof}
We use $H^{1/2}(\mathbb{T}\times\mathbb{R})=L^2_\theta H^{1/2}_r\cap H^{1/2}_\theta L^2_r$ with equivalent norms. The $\theta$-part follows by interpolation: $a$ multiplies $L^2_\theta L^2_r$ with norm $\|a\|_\infty$, and $H^1_\theta L^2_r$ with norm $\|a\|_\infty+\|\partial_\theta a\|_\infty$.

For the $r$-part, fix $\theta$ and let $1/\tau=1/s-1/2$, so that $\tau\in(2,\infty)$. The Hilbert transform is bounded on $L^s$, and by Hardy--Littlewood--Sobolev \(\big\||D_r|^{1/2}a\big\|_{L^\tau}=\big\||D_r|^{-1/2}\mathcal{H}\partial_ra\big\|_{L^\tau}\lesssim\|\partial_ra\|_{L^s}.\) The fractional Leibniz rule \cite{KP,KPV} gives
\begin{equation*}
\big\||D_r|^{1/2}(aF)\big\|_{L^2}\lesssim\|a\|_{L^\infty}\big\||D_r|^{1/2}F\big\|_{L^2}+\big\||D_r|^{1/2}a\big\|_{L^\tau}\|F\|_{L^{2\tau/(\tau-2)}},
\end{equation*}
and $H^{1/2}(\mathbb{R})\subset L^{2\tau/(\tau-2)}$.
\end{proof}

\begin{lemma}[$H^{3/2}$ regularity for the capillary operator]
\label{lem:tl_ell} Let $u\in H^1(D)$ satisfy, for all $\varphi\in H^1(D)$,
\begin{equation*}
g\int_Du\varphi+\sigma\int_D\mathcal{L}_{\Sigma,0}\nabla u\cdot\nabla\varphi=\Lambda(\varphi)+\int_DB\cdot\nabla\varphi+\int_{\partial D}h\varphi ,
\end{equation*}
where $\Lambda$ is linear with $|\Lambda(\varphi)|\le\lambda\|\varphi\|_{H^{1/2}(D)}$, $B\in H^{1/2}(D;\mathbb{R}^2)$ and $h\in L^2(\partial D)$. Then
\begin{equation*}
\|u\|_{H^{3/2}(D)}\lesssim\lambda+\|B\|_{H^{1/2}(D)}+\|h\|_{L^2(\partial D)}.
\end{equation*}
\end{lemma}

\begin{proof}
The form on the left is bounded and coercive on $H^1(D)$, since $g>0$ and $\mathcal{L}_{\Sigma,0}$ is smooth and uniformly positive definite (it is $DF_\zeta(\nabla\zeta)$, $F_\zeta(z)=z(1+|z|^2)^{-1/2}$). Let $\mathcal{S}:(H^1(D))^*\to H^1(D)$ be the Lax--Milgram solution operator. Then $u=\mathcal{S}\ell$, where $\ell$ is the right-hand side; $\Lambda$ extends uniquely to $(H^{1/2})^*$. Since $\partial D$ is smooth, conormal $H^2$ regularity gives $\mathcal{S}:\varphi\mapsto\int F\varphi+\int_{\partial D}h\varphi\ \mapsto H^2$ for $F\in L^2$ and $h\in H^{1/2}(\partial D)$. Interpolating at $\vartheta=1/2$ gives three bounds into $[H^1,H^2]_{1/2}=H^{3/2}$:

\begin{itemize}
\item interior data: between $(H^1)^*\to H^1$ and $L^2\to H^2$. Here $[(H^1)^*,L^2]_{1/2}=([H^1,L^2]_{1/2})^*=(H^{1/2})^*$ by the duality theorem.

\item conormal data: between $H^{-1/2}(\partial D)\to H^1$ and $H^{1/2}(\partial D)\to H^2$, so $L^2(\partial D)\to H^{3/2}$.

\item divergence-form data $\varphi\mapsto\int B\cdot\nabla\varphi$: between $B\in L^2\to H^1$ and $B\in H^1\to H^2$, the latter after writing $\int B\cdot\nabla\varphi=-\int\func{div} B\,\varphi+\int_{\partial D}B\cdot\nu\,\varphi$.
\end{itemize}
\end{proof}

\begin{lemma}[Remainders at the top level]
\label{lem:tl_rem} Assume (H1)--(H4). Assume \(\mathcal E\le\min\{1,\delta_0(\varepsilon)\},\) where $\delta_0(\varepsilon)$ is chosen so that the flattening map and its inverse satisfy the previously established geometric bounds. For the fixed $\varepsilon$, take $C(\varepsilon)$ large enough for the absorptions here and in \eqref{eq:closing}, and choose \(0<\delta_*\le\min\{1,\delta_0(\varepsilon),(2C(\varepsilon))^{-2}\},\) decreasing it further if any earlier lemma requires this. At every time $s$ for which $\mathcal E(s)\le\delta_*$, we thus have $C(\varepsilon)\sqrt{\mathcal E(s)}\le 1/2$.

\begin{enumerate}[label=(\alph*)]

\item For $w=\chi_{\mathrm{bet}+}\,\tilde w$ with $\tilde w\in H^1_0(\Omega;\mathbb{R}^3)$,
\begin{equation*}
\begin{aligned}
& \Big|\int_\Omega J\,\Ft^1_{\rm lin}\cdot w\Big| \le\Big[C_0\Big(\|\chis\partial_\theta^{m+1}v\|_{H^1}+(\sE_{1,m,{\rm energy,cl}}+\cD_{\rm int})^{1/2}\Big) +C(\varepsilon)\sqrt{\cE}\,\cD^{1/2}\Big]\|\tilde w\|_{H^1}, \\
& \Big|\int_\Omega J\,\Ft^1_{\rm nl}\cdot w\Big| \le C(\varepsilon)\sqrt{\cE}\,\cD^{1/2}\|\tilde w\|_{H^1},
\end{aligned}
\end{equation*}
where $C_0$ does not depend on $\varepsilon$.

\item For every admissible $w$ with $\func{supp}w\subset\func{supp}\chi_{\mathrm{cl}}$,
\begin{equation*}
\Big|\int_\Omega J\widetilde F^1_{\mathrm{rem}}\cdot w-\langle\widetilde F^3_{\mathrm{rem}},w\rangle_\Sigma-\langle J_s\widetilde F^6_{\mathrm{rem}},U_0w\rangle_{\Sigma_s}\Big|\le C(\varepsilon)\,\mathcal{D}^{1/2}\|w\|_{H^1(\Omega)}.
\end{equation*}

\item $\|\chi_{\mathrm{cl},\Sigma}\widetilde F^5_{\mathrm{rem}}\|_{H^{1/2}(\Sigma)}+\|\widetilde F^5_{\mathrm{rem}}\|_{L^2(\func{supp}\chi_{\mathrm{cl},\Sigma})}\le C(\varepsilon)\sqrt{\mathcal{E}}\,\mathcal{D}^{1/2}$, and $\|\widetilde F^7_{\mathrm{rem}}\|_{L^2(\Gamma_0)}\le C(\varepsilon)\,\mathcal{D}^{1/2}$.
\end{enumerate}
\end{lemma}

\begin{proof}
We use the chart variables $y=(\theta,r,h)$ on $\mathcal{U}$, with $dx=J_\Phi\,dy$. This is well-defined and we have the corresponding norm equivalence given that $\mathcal{E}$ is required to be small, which implies that $J$ is close to $1$ and therefore larger than a positive constant uniformly.  The Cartesian derivatives are $\partial_l=\sum_ae_l^a(y)\partial_{y_a}$ with smooth coefficients $e_l^a$ that depend only on the chart. Since $\partial_\theta=\partial_{y_1}$ commutes with every $\partial_{y_a}$,
\begin{equation}  \label{tl:comm} [\partial_\theta^{m+1},\partial_l]=\sum_{i=1}^{m+1}\binom{m+1}{i}\sum_a\big(\partial_\theta^ie_l^a\big)\,\partial_{y_a}\partial_\theta^{m+1-i}.
\end{equation}
The term $i=1$, $y_a=\theta$ is $(\partial_\theta e_l^\theta)\,\partial_\theta^{m+1}$. Its coefficient is bounded but not small; for a disc in polar coordinates $[\partial_\varphi,\partial_1]=-\partial_2$. All other terms carry at most $m$ angular derivatives and one further derivative.

\emph{Hardy's inequality.} $\Omega$ is a bounded Lipschitz domain, so $\|\psi/d\|_{L^2(\Omega)}\le C_H\|\nabla\psi\|_{L^2(\Omega)}$ for $\psi\in H^1_0(\Omega)$, where $d=\func{dist}(\cdot,\partial\Omega)$ \cite{Necas62}. By \eqref{tl:cutoff}, $|\nabla\chi_{\mathrm{bet}+}|\le C_1\varepsilon^{-1}\le C_1^2/d$ on $\func{supp}\nabla\chi_{\mathrm{bet}+}$. Hence
\begin{equation}  \label{tl:hardy} \|\psi\,\nabla\chi_{\mathrm{bet}+}\|_{L^2(\Omega)}\le C_0\|\nabla\psi\|_{L^2(\Omega)}\qquad(\psi\in H^1_0(\Omega)),
\end{equation}
with $C_0=C_1^2C_H$ independent of $\varepsilon$.

\emph{(a), the linear part.} By \eqref{tl:comm}, applied to $\nabla p$ and to $\func{div}\mathbb{D} v$, the field $\widetilde F^1_{\mathrm{lin}}=-[\partial_\theta^{m+1},\func{div} S_I](p,v)$ is a finite sum of terms $c\,\mathcal{T}$, where $c$ is smooth and depends only on $\zeta$ and the chart, and $\mathcal{T}$ is one of
\begin{equation*}
\text{(1) } p_\natural;\quad\text{(2) } \partial_{y_a}\partial_\theta^jp,\ \partial_\theta^jp;\quad \text{(3) } \partial_{y_a}v_\natural,\ v_\natural;\quad\text{(4) } \partial_{y_a}\partial_{y_b}\partial_\theta^jv,\ \partial_{y_a}\partial_\theta^jv,\ \partial_\theta^jv,
\end{equation*}
with $j\le m$, and $y_a,y_b\in\{r,h\}$ in the second-order terms of (2) and (4); a derivative $\partial_\theta$ in front of $\partial_\theta^jp$ or $\partial_\theta^jv$ only raises $j$, and for $j=m$ it gives (1) or (3). Consider first the integrals $\int c\,\mathcal{T}\cdot w\,J_\Phi\,dy$, without the factor $J$. The function $\tilde w$ vanishes on $\partial\Omega$, and $\chi_{\mathrm{bet}+}$ vanishes near the rest of the boundary of the chart region. So integrations by parts in $y$ produce no boundary terms.

\begin{itemize}
\item (1): $\theta$ is periodic and $\chi_{\mathrm{bet}+}$ does not depend on $\theta$, so
\begin{equation*}
\Big|\int c\,p_\natural\cdot\chi_{\mathrm{bet}+}\tilde wJ_\Phi\,dy\Big| =\Big|\int\chi_{\mathrm{bet}+}\,\partial_\theta(cJ_\Phi\tilde w)\cdot\partial_\theta^mp\,dy\Big| \le C\|\chi_{\mathrm{bet}+}\partial_\theta^mp\|_{L^2}\|\tilde w\|_{H^1}.
\end{equation*}

\item (2): integrate by parts once in $y_a\in\{r,h\}$. The derivative falls on $c\chi_{\mathrm{bet}\allowbreak+}\tilde wJ_\Phi$. When it falls on $\chi_{\mathrm{bet}\allowbreak+}$ we use \eqref{tl:hardy} for $\tilde w$ and \eqref{tl:Dint} for $\partial_\theta^jp$ on $\func{supp}\nabla\chi_{\mathrm{bet}\allowbreak+}$. The bound is $C\allowbreak\|\chi_{\mathrm{bet}\allowbreak+}\partial_\theta^jp\|_{L^2}\allowbreak\|\tilde w\|_{H^1}\allowbreak+C_0\allowbreak\|\partial_\theta^jp\|_{L^2(\func{supp}\nabla\chi_{\mathrm{bet}\allowbreak+})}\allowbreak\|\nabla\tilde w\|_{L^2}$.

\item (3): write $\chi_{\mathrm{bet}+}\partial_{y_a}v_\natural=\partial_{y_a}(\chi_{\mathrm{bet}+} v_\natural)-v_\natural\partial_{y_a}\chi_{\mathrm{bet}+}$. The first part gives $C\|\chi_{\mathrm{bet}+} v_\natural\|_{H^1}\|\tilde w\|_{L^2}$. The second gives $\|v_\natural\|_{L^2(\func{supp}\nabla\chi_{\mathrm{bet}+})}\|\tilde w\nabla\chi_{\mathrm{bet}+}\|_{L^2}\le C_0\mathcal{D}_{\mathrm{int}}^{1/2}\|\nabla\tilde w\|_{L^2}$, by \eqref{tl:Dint} and \eqref{tl:hardy}.

\item (4): integrate the second-order terms by parts once in $y_b$, as in (2). When the derivative falls on $\chi_{\mathrm{bet}\allowbreak+}$, use \eqref{tl:hardy} and \eqref{tl:Dint}, as in (2). The remaining term is $\int c\,J_\Phi\chi_{\mathrm{bet}\allowbreak+}\,\partial_{y_a}\partial_\theta^jv\cdot\partial_{y_b}\tilde w$. It is bounded by $C\allowbreak\|\chi_{\mathrm{bet}\allowbreak+}\nabla f\|_{L^2}\allowbreak\|\tilde w\|_{H^1}$ with $f=\partial_\theta^jv$. Here Hardy's inequality is not available, since $f$ does not vanish on $\partial\Omega$. Instead, $\allowbreak\|\chi_{\mathrm{bet}\allowbreak+}\nabla f\|_{L^2}\allowbreak\le\allowbreak\|\nabla(\chi_{\mathrm{bet}\allowbreak+} f)\|_{L^2}\allowbreak+\allowbreak\|\nabla f\|_{L^2(\{0<\chi_{\mathrm{bet}\allowbreak+}<1\})}$, because $\chi_{\mathrm{bet}\allowbreak+}\nabla f=\nabla(\chi_{\mathrm{bet}\allowbreak+} f)$ where $\chi_{\mathrm{bet}\allowbreak+}=1$, and $\chi_{\mathrm{bet}\allowbreak+}\allowbreak\le1$. The second term is bounded by \eqref{tl:Dint}, since $\partial_\theta^jv\in H^1$ there. The first-order terms of (4) are treated in the same way.
\end{itemize}

By \eqref{tl:Dint}, the parts on $\func{supp}\nabla\chi_{\mathrm{bet}+}$ are bounded by $\mathcal{D}_{\mathrm{int}}^{1/2}\le(\mathscr{E}_{1,m,\mathrm{energy,cl}}+\mathcal{D}_{\mathrm{int}})^{1/2}$. This gives the first bound of (a) without $J$, with an $\varepsilon$-independent constant: no derivative falls on $\chi_{\mathrm{bet}+}$ without the compensating Hardy factor. With $J$, the difference $\int(J-1)\,\widetilde F^1_{\mathrm{lin}}\cdot w$ is estimated in the same way. Lemma~\ref{lem:tl_coef} gives $|J-1|\lesssim\sqrt\mathcal{E}$, $|\partial_\theta J|\lesssim\sqrt\mathcal{E}$ and $\|\nabla J\|_{L^\infty_\theta L^{r_-}_{rh}}\lesssim\sqrt\mathcal{E}$. In the integrations by parts a derivative may fall on $J$; we use H\"older's inequality with $\tilde w\in L^2_\theta L^{2/\varepsilon_-}_{rh}$, as in Step~A1 below, and we bound the cutoff terms crudely by \eqref{tl:cutoff}. This gives $C(\varepsilon)\sqrt\mathcal{E}\,\mathcal{D}^{1/2}\|\tilde w\|_{H^1}$.

\emph{(a), the nonlinear part.} The terms of $\widetilde F^1_{\mathrm{nl}}$ are of four kinds:

\begin{enumerate}[label=(\roman*)]

\item the Leibniz terms of $\partial_\theta^{m+1}\big(\func{div}_\mathcal{A} S_\mathcal{A}(p,v)\big)$ in which $i$ angular derivatives fall on the coefficients, $1\le i\le m$;

\item the commutator terms with coefficient $\mathcal{A}-I$, which are the terms (1)--(4) above multiplied by $\mathcal{A}-I$ or by its first derivatives;

\item the terms of $\partial_\theta^{m+1}F^1_0$ other than those in $\widetilde F^1_{\mathrm{top}}$, including $WK\partial_t\bar\eta\,\partial_3v_\natural$, $v\cdot\nabla_\mathcal{A} v_\natural$ and $v_\natural\cdot\nabla_\mathcal{A} v$;

\item the terms in which $U=\partial_\theta^{m+1}\bar\eta$ appears without a further derivative, such as $\mathfrak{A}^{\prime }(\nabla(W\bar\eta))[U\nabla W]$ and the terms $(\partial_\theta e^\theta_l)\,U$ of \eqref{tl:comm} applied to $W\bar\eta$.
\end{enumerate}

Each term is a product of at least two factors that depend on the solution. As in the proof of Lemma~\ref{thm:contact_line_forcing_estimate}, and because $m\ge3$, at most one factor carries more than $m-1$ angular derivatives or more derivatives than Lemma~\ref{lem:tl_coef} controls in $L^\infty_\theta L^\infty$ or $L^\infty_\theta L^{q_-}_{rh}$. The other factors are bounded there by $\sqrt\mathcal{E}$. The highest factor is one of the following.

\begin{itemize}
\item $\nabla^2\partial_\theta^j(\chi_{\mathrm{bet}+} v)$ or $\nabla\partial_\theta^j(\chi_{\mathrm{bet}+} p)$, $j\le m$. On $\func{supp}\chi_{\mathrm{bet}+}$ these are bounded in $L^2_\theta L^{q_\delta}_{rh}$ by $C(\varepsilon)\mathcal{D}^{1/2}$, through the elliptic norms of $\chi_{\mathrm{cl}}$ and $\chi_{\mathrm{bet}}$ in $\mathcal{D}_{\mathrm{cl}}$ and $\mathcal{D}_{\mathrm{bet}}$.

\item $\nabla(\chi_{\mathrm{bet}+}\partial_\theta^i\bar\eta)$, $i\le m$. It is bounded in $L^2_\theta H^{1-\alpha}_{rh}(\mathcal{U})\subset L^2_\theta L^{2/\alpha}_{rh}$ by $C(\varepsilon)\mathcal{D}^{1/2}$, by the argument of Lemma~\ref{lem:tl_slice}. This needs only $\partial_\theta^i\eta\in H^{3/2-\alpha}$, because $\|\partial_\theta G\|_{L^2_\theta H^{1/2-\alpha}_r}\le\|G\|_{H^{3/2-\alpha}}$, and the $\chi_{\mathrm{bet}+,\Sigma}$-weighted norm is bounded in \eqref{tl:VD}. The $3$-D piece is bounded by $\mathcal{D}_{\mathrm{int}}$. Likewise $\partial_r^2\partial_\theta^i\bar\eta$, $i\le m$, is bounded in $L^2_\theta W^{1,q_\delta}_{rh}\subset L^2_\theta L^{r_-}_{rh}$ by $C(\varepsilon)\mathcal{D}^{1/2}$, through the elliptic norms of $\mathcal{D}_{\mathrm{cl}}$ and the bound $W^{3-1/q,q}\to W^{3,q}$ of $P_{2-D}E$ on bounded slices.

\item Two factors of intermediate order, which occurs only for $m=3$, as in $(\partial_\theta^2\mathcal{A})\,\partial_\theta^2\nabla^2v$. Here $\partial_\theta^2\mathcal{A}$ is bounded in $L^\infty_\theta L^{2/\alpha}_{rh}$ by $\sqrt\mathcal{E}$, by $H^1_\theta\subset L^\infty_\theta$ and the previous item with $i=2,3\le m$ at energy levels. The factor $\partial_\theta^2\nabla^2v$ is bounded in $L^2_\theta L^{q_\delta}_{rh}$ by $\mathcal{D}^{1/2}$. Since $\frac\alpha2+\frac1{q_\delta}+\frac{\delta-\alpha}{2}=1$, H\"older's inequality applies with $\tilde w\in L^2_\theta L^{2/(\delta-\alpha)}_{rh}$.

\item $v_\natural$, $\nabla v_\natural$, $p_\natural$ or $U$, with a coefficient that is $O(\sqrt\mathcal{E})$. These are treated as in the linear part, or by \eqref{tl:s1} and \eqref{tl:s4}.
\end{itemize}

The pairing with $w=\chi_{\mathrm{bet}+}\tilde w$ is then bounded by H\"older's inequality as in Lemma~\ref{lem:tl_first}, with $\tilde w\in L^2_\theta L^s_{rh}$ for every $s<\infty$. For instance, for $(\partial_\theta^m\mathcal{A})\nabla^2\partial_\theta v$, the exponents are $(2,\infty,2)$ in $\theta$ and $\frac\alpha2+\frac{2-\varepsilon_-}{2}+\frac1s=1$, $s=\frac{2}{\varepsilon_--\alpha}$, in $(r,h)$. For $(\partial_\theta\mathcal{A})\nabla^2\partial_\theta^mv$ they are $(\infty,2,2)$ and $\frac{1}{q_\delta}+\frac{1}{q_\delta^{\prime }}=1$. This gives the second bound of (a).

\emph{(b).} On $\func{supp}\chi_{\mathrm{cl}}$ we have $\chi_{\mathrm{bet}+}=1$ and $\nabla\chi_{\mathrm{bet}+}=0$, so every quantity on $\func{supp}w$ equals its $\chi_{\mathrm{bet}+}$-weighted version. The bulk terms are bounded as in (a), except that $w$ need not vanish on $\partial\Omega$. The integrations by parts in $(r,h)$ therefore produce boundary terms on $\Gamma_\pm$. These are traces of $\partial_\theta^jp$ and $\nabla\partial_\theta^jv$, $j\le m$, which lie in $L^2_\theta L^{s_-}$ on $\Gamma_\pm$ with norm $\le C(\varepsilon)\mathcal{D}^{1/2}$, against $w|_{\Gamma_\pm}\in L^2_\theta H^{1/2}\subset L^2_\theta L^{s}$, $s<\infty$. The boundary remainders are bounded in the same way. By \eqref{tl:comm}, $\widetilde F^3_{\mathrm{lin}}$ and $\widetilde F^6_{\mathrm{lin}}$ consist of traces of $\partial_\theta^jp$, $\nabla\partial_\theta^jv$, $\partial_\theta^jv$ ($j\le m$) and $v_\natural$ with smooth coefficients; $U_0$ annihilates the pressure terms of $\widetilde F^6_{\mathrm{lin}}$, since $U_0n_0=0$. They also contain the commutators of $\partial_\theta^{m+1} $ with the capillary operator $g-\sigma\func{div}(\mathcal{L}_{\Sigma,0}\nabla\,\cdot\,)$ and with $N_\zeta$. These consist of $\nabla^2\partial_\theta^j\eta$, $\nabla\partial_\theta^j\eta$, $\partial_\theta^j\eta$ ($j\le m$), $\nabla\eta_\natural$ and $\eta_\natural$ with smooth coefficients. We pair them with $\phi:=w|_\Sigma$ directly in the chart variables $(\theta,r)$ of $\Sigma$, without passing to a negative-order space on $\Sigma$. By the trace theorem and $H^{1/2}(\mathbb{T}\times\mathbb{R})=L^2_\theta H^{1/2}_r\cap H^{1/2}_\theta L^2_r$, \(\phi\in L^2_\theta H^{1/2}_r,\,\allowbreak \|\phi\|_{L^2_\theta H^{1/2}_r}\le C\|w\|_{H^1(\Omega)};\) only the regularity of $\phi$ in $r$ is used below. Write $q:=q_\delta=\frac{2}{2-\delta}\in(1,2)$, so that $q^{\prime }=\frac2\delta\in(2,\infty)$.

\emph{The second $r$-derivatives.} For $j\le m$ the elliptic norms of $\mathcal{D}_{\mathrm{cl}}$ give $\partial_\theta^j\eta\in L^2_\theta W^{3-1/q,q}_r$ on $\func{supp}\chi_{\mathrm{cl},\Sigma}$, with norm $\le C(\varepsilon)\mathcal{D}^{1/2}$. Hence \(\partial_r^2\partial_\theta^j\eta\in L^2_\theta W^{3-2-1/q,q}_r=L^2_\theta W^{1-1/q,q}_r .\) The smoothness index $1-1/q=\delta/2$ is positive, so $\partial_r^2\partial_\theta^j\eta(\theta,\cdot)$ is a function, and $W^{1-1/q,q}_r\subset L^q_r$. On the other hand $H^{1/2}_r\subset L^{q^{\prime }}_r$, since $2\le q^{\prime }<\infty$. Thus $W^{1-1/q,q}_r$ pairs with $H^{1/2}_r$ through H\"older's inequality in $r$, and Cauchy--Schwarz in $\theta$ gives
\begin{equation*}
\Big|\int_\Sigma c\,\partial_r^2\partial_\theta^j\eta\,\phi\Big|  \le C\int_\mathbb{T}\|\partial_r^2\partial_\theta^j\eta(\theta,\cdot)\|_{L^q_r}\|\phi(\theta,\cdot)\|_{L^{q^{\prime }}_r}\,d\theta \le C\|\partial_\theta^j\eta\|_{L^2_\theta W^{3-1/q,q}_r}\|\phi\|_{L^2_\theta H^{1/2}_r} \le C(\varepsilon)\mathcal{D}^{1/2}\|w\|_{H^1}.
\end{equation*}

\emph{The remaining terms.} The other second derivatives carry at least one $\partial_\theta$: $\partial_\theta\partial_r\partial_\theta^j\eta=\partial_r\partial_\theta^{j+1}\eta$ and $\partial_\theta^2\partial_\theta^j\eta=\partial_\theta\partial_\theta^{j+1}\eta$. They are first derivatives of $\partial_\theta^{j+1}\eta$ with $j+1\le m+1$, which lies in $H^{3/2-\alpha}$ on $\func{supp}\chi_{\mathrm{cl},\Sigma}$ with norm $\le C(\varepsilon)\mathcal{D}^{1/2}$ by \eqref{tl:VD}. So they lie in $H^{1/2-\alpha}\subset L^2$, and so do $\nabla\eta_\natural$, $\eta_\natural$, $\nabla\partial_\theta^j\eta$ and $\partial_\theta^j\eta$. Their pairing with $\phi\in L^2(\Sigma)$ is $\le C(\varepsilon)\mathcal{D}^{1/2}\|w\|_{H^1}$.

The nonlinear boundary remainders are products of the same traces with factors bounded in $L^\infty$ by $\sqrt\mathcal{E}$; they are split and paired with $\phi$ in the same way, or bounded as the face terms of Lemma~\ref{lem:tl_first}.

For later use in (c) we record that, for $j\le m$,
\begin{equation}  \label{tl:F1fact} \partial_\theta^j\eta\in L^2_\theta H^{3/2+\delta}_r\cap H^{5/2-\alpha}_\theta L^2_r\subset H^{3/2+\delta}\quad\text{on }\func{supp}\chi_{\mathrm{cl},\Sigma},\ \text{with norm}\le C(\varepsilon)\mathcal{D}^{1/2}.
\end{equation}
The first space comes from the elliptic norms of $\mathcal{D}_{\mathrm{cl}}$, since $W^{3-1/q_\delta,q_\delta}(\mathbb{R})\subset H^{3/2+\delta}(\mathbb{R})$ (Subsection~\ref{subsec:constants}). The second comes from $\partial_\theta^{j+1}\eta\in H^{3/2-\alpha}$ in \eqref{tl:VD}. The inclusion holds because $3/2+\delta<5/2-\alpha$.

\emph{(c), the term $\widetilde F^5_{\mathrm{rem}}$.} $\mathcal{R}$ is quadratic, so $\widetilde F^5_{\mathrm{rem}}=\widetilde F^5_{\mathrm{nl}}$. By the chain and Leibniz rules and \eqref{tl:comm} on $\Sigma$, it is a sum of three kinds of terms:

\begin{itemize}
\item $\mathcal{R}^{\prime }(\nabla\eta)[\partial_\theta^{m+1},\nabla]\eta$, which is a sum of terms $\mathcal{R}^{\prime }(\nabla\eta)\,c\,\eta_\natural$ and $\mathcal{R}^{\prime }(\nabla\eta)\,c\,\nabla\partial_\theta^j\eta$ with $j\le m$;

\item $\mathcal{R}^{(k)}(\nabla\eta)[\partial_\theta^{j_1}\nabla\eta,\dots,\partial_\theta^{j_k}\nabla\eta]$ with $k\ge2$, $j_1+\dots+j_k=m+1$ and $1\le j_l\le m$, together with commutator corrections of the same type and lower order;

\item the terms coming from the $x$-dependence of $\mathcal{R}$ (Remark~\ref{rem:R_x_dependence}), of the same type with fewer derivatives on $\eta$.
\end{itemize}

In the chart, $H^{1/2}(\mathbb{T}\times\mathbb{R})=L^2_\theta H^{1/2}_r\cap H^{1/2}_\theta L^2_r$. Since $H^{1/2+\delta}(\mathbb{R})$ multiplies $H^{1/2}(\mathbb{R})$, and $H^1(\mathbb{T})$ multiplies $H^{1/2}(\mathbb{T})$,
\begin{equation}  \label{tl:mult2} \|aF\|_{H^{1/2}}\lesssim\Big(\|a\|_{L^\infty_\theta H^{1/2+\delta}_r}+\sup_r\|a(\cdot,r)\|_{H^1_\theta}\Big)\|F\|_{H^{1/2}}.
\end{equation}
For the terms of the first kind we use Lemma~\ref{lem:tl_mult} with $a=\mathcal{R}^{\prime }(\nabla\eta)c$, and then $\|\chi_{\mathrm{cl},\Sigma}\eta_\natural\|_{H^{1/2}}\le C(\varepsilon)\mathcal{D}^{1/2}$ and \eqref{tl:F1fact}; the result is $\le C(\varepsilon)\sqrt\mathcal{E}\,\mathcal{D}^{1/2}$. For the terms of the second kind with $k=2$, order the indices so that $j_1\le j_2$. Then $1\le j_1\le(m+1)/2\le m-1$. Put $a:=\mathcal{R}^{\prime \prime }(\nabla\eta)\partial_\theta^{j_1}\nabla\eta$ and $F:=\partial_\theta^{j_2}\nabla\eta$.

\begin{itemize}
\item If $j_2\le m-1$, then \eqref{tl:F1fact} holds for $\partial_\theta^{j_2}\eta$ with $\mathcal{E}$ in place of $\mathcal{D}$, because $\partial_\theta^{j_2}\eta$ and $\partial_\theta^{j_2+1}\eta$ are at energy levels. So $\|F\|_{H^{1/2}}\le C(\varepsilon)\sqrt\mathcal{E}$. The norms of $a$ in \eqref{tl:mult2} are $\le C(\varepsilon)\mathcal{D}^{1/2}$, by $H^1_\theta\subset L^\infty_\theta$ and \eqref{tl:F1fact} for $\partial_\theta^{j}\eta$, $j\le j_1+1\le m$, together with $\partial_\theta^{j_1+2}\eta\in H^{3/2-\alpha}\subset L^2_\theta L^\infty_r$, as $j_1+2\le m+1$.

\item If $j_2=m$, then $j_1=1$ and the roles are exchanged. Now $\|F\|_{H^{1/2}}\le C(\varepsilon)\mathcal{D}^{1/2}$ by \eqref{tl:F1fact}, and the norms of $a=\mathcal{R}^{\prime \prime }(\nabla\eta)\partial_\theta\nabla\eta$ in \eqref{tl:mult2} are $\le C(\varepsilon)\sqrt\mathcal{E}$, because they involve only $\partial_\theta^j\eta$ with $j\le3\le m$ at energy levels.
\end{itemize}

The terms with $k\ge3$ have more factors, all but one bounded in $L^\infty$ and in the multiplier norms of \eqref{tl:mult2} by $\sqrt\mathcal{E}$, and are bounded in the same way. Multiplication by $\chi_{\mathrm{cl},\Sigma}$ costs $C(\varepsilon)$. The $L^2(\func{supp}\chi_{\mathrm{cl},\Sigma})$ bound follows from the same factorization and H\"older's inequality.

\emph{(c), the term $\widetilde F^7_{\mathrm{rem}}$.} By \eqref{tl:forcing}, $\widetilde F^7=\,\partial_\theta^{m+1}\hat{\mathcal{W}}(\partial_t\eta)+\sigma\big(\partial_\theta^{m+1}\mathcal{B}_{\mathrm{cl},0}\eta-\mathcal{B}_{\mathrm{cl},0}\eta_\natural\big)+\sigma\big(\partial_\theta^{m+1}(\mathcal{R}\cdot\nu)-(\partial_\theta^{m+1}\mathcal{R})\cdot\nu\big)$ on $\Gamma_0$. The linear part is the middle term. By \eqref{tl:comm} it consists of traces on $\Gamma_0$ of $\nabla\partial_\theta^j\eta$, $\partial_\theta^j\eta$ ($j\le m$) and $\eta_\natural$ with smooth coefficients. By \eqref{tl:F1fact} and the trace theorem, its $L^2(\Gamma_0)$ norm is $\le C(\varepsilon)\mathcal{D}^{1/2}$. The first term minus $\widetilde F^7_{\mathrm{top}}$ consists of products of $\partial_t\partial_\theta^j\eta|_{\Gamma_0}$, $j\le m$, which are bounded in $L^2(\Gamma_0)$ by $\mathcal{D}^{1/2}$, with factors bounded in $L^\infty(\Gamma_0)$ by $\sqrt\mathcal{E}$. Here we use $\hat{\mathcal{W}}(z)=O(z^2)$ and $H^1(\Gamma_0)\subset L^\infty(\Gamma_0)$ for $\partial_t\partial_\theta^j\eta$, $j\le m-1$. The last term consists of traces on $\Gamma_0$ of $\partial_\theta^j\mathcal{R}$, $j\le m$, times smooth functions. By the chain rule, $\partial_\theta^j\mathcal{R}=\mathcal{R}^{\prime }(\nabla\eta)\nabla\partial_\theta^j\eta$ plus products of factors $\nabla\partial_\theta^i\eta$ with $i<j$. Each factor $\nabla\partial_\theta^i\eta$ has a trace in $L^2(\Gamma_0)$, because $\partial_\theta^i\eta\in H^{3/2+\delta}$ by \eqref{tl:F1fact} and the trace maps $H^{1/2+\delta}(D)$ into $L^2(\partial D)$. The coefficients are bounded in $L^\infty(\Gamma_0)$ by $\sqrt\mathcal{E}$, by Lemma~\ref{lem:tl_coef}. Hence the last term is $\le C(\varepsilon)\sqrt\mathcal{E}\,\mathcal{D}^{1/2}$ in $L^2(\Gamma_0)$.
\end{proof}

\subsubsection{Statement and proof}

Theorem~\ref{thm:top_level} and Corollary~\ref{cor:top_level} are stated in Section~\ref{sec:intro}, under the hypotheses (H1)--(H4) above.

\smallskip\noindent\textbf{Proof of Theorem~\ref{thm:top_level}.}

Fix $t$ with $\mathcal{E}(t)\le\delta_0\le1$; $\delta_0$ is chosen in Step~B2. We abbreviate $\chi:=\chi_{\mathrm{cl}}$, $\chi_\Sigma:=\chi_{\mathrm{cl},\Sigma}$ and $\|\cdot\|_{s}:=\|\cdot\|_{H^s(\Sigma)}$.

\smallskip\noindent\textbf{Part A: the pressure.} Let $P:=\chi_{\mathrm{bet}+}\,p_\natural=\partial_\theta^{m+1}(\chi_{\mathrm{bet}+} p)$ on $\mathcal{U} $, extended by zero to $\Omega$. In Part~A, $C_0$ denotes constants independent of $\varepsilon$, and $\delta_0(\varepsilon)$ is so small that $C(\varepsilon)\sqrt{\delta_0(\varepsilon)}\le1$ for the constant $C(\varepsilon)$ of Lemma~\ref{lem:tl_coef} in Step~A1.

\emph{Step A1: the test field.} For $\omega\in H^1_0(\Omega;\mathbb{R}^3)$ set \(w:=\chi_{\mathrm{bet}+}\,\tilde w,\,\allowbreak \tilde w:=M\omega,\,\allowbreak M:=(J\mathcal{A}^T)^{-1}=J^{-1}\nabla\Phi_\eta .\) Then $\tilde w,w\in H^1_0(\Omega)$ and $\func{supp}w\subset\func{supp}\chi_{\mathrm{bet}+}$. So $w$ is admissible, and every boundary term in \eqref{tl:weak} vanishes. By the Piola identity $\partial_j(J\mathcal{A}_{ij})=0$,
\begin{equation*}
J\func{div}_\mathcal{A} w =\partial_j(J\mathcal{A}_{ij}w_i)=\func{div}(J\mathcal{A}^Tw) =\func{div}(\chi_{\mathrm{bet}+}\omega)=\chi_{\mathrm{bet}+}\func{div}\omega+\nabla\chi_{\mathrm{bet}+}\cdot\omega .
\end{equation*}
Moreover $\|\tilde w\|_{H^1}\le C_0\|\omega\|_{H^1}$. Indeed $|M-I|\lesssim\sqrt\mathcal{E}$, and $|\nabla M|\lesssim|\nabla^2(W\bar\eta)|+|\nabla(W\bar\eta)|$, which lies in $L^\infty_\theta L^{r_-}_{rh}$ by Lemma~\ref{lem:tl_coef}, while $\omega\in L^2_\theta L^{2/\varepsilon_-}_{rh}$; since $\frac{1-\varepsilon_-}{2}+\frac{\varepsilon_-}{2}=\frac12$, the product is in $L^2$. Finally, by Hardy's inequality \eqref{tl:hardy} applied to $\psi=\tilde w$, \(\|w\|_{H^1}\le\|\chi_{\mathrm{bet}+}\nabla\tilde w\|_{L^2}+\|\tilde w\nabla\chi_{\mathrm{bet}+}\|_{L^2}+\|\tilde w\|_{L^2}\le C_0\|\omega\|_{H^1},\) uniformly in $\varepsilon$, although $|\nabla\chi_{\mathrm{bet}+}|$ is of order $\varepsilon^{-1}$.

Since $-\int JS_\mathcal{A}(p_\natural,v_\natural):\nabla_\mathcal{A} w=\frac\mu2\int J\mathbb{D}_\mathcal{A} v_\natural:\mathbb{D}_\mathcal{A} w-\int Jp_\natural\func{div}_\mathcal{A} w$, identity \eqref{tl:weak} becomes the analogue of \eqref{eq:pressure_weak} for $w\in H^1_0(\Omega)$:
\begin{equation}
\label{tl:pA}
\int_\Omega P\func{div}\omega =-\int_\Omega p_\natural\nabla\chi_{\mathrm{bet}+}\cdot\omega+\int_\Omega J\partial_tv_\natural\cdot w +\frac\mu2\int_\Omega J\mathbb{D}_\mathcal{A} v_\natural:\mathbb{D}_\mathcal{A} w-\int_\Omega J\widetilde F^1\cdot w .
\end{equation}

\emph{Step A2: the functional is bounded on $H^1_0$.} We bound the four terms of \eqref{tl:pA}.

\begin{itemize}
\item \emph{Cutoff commutator.} By \eqref{tl:Dint}, \eqref{tl:cutoff} and $H^1\subset L^6$,
\begin{equation*}
\Big|\int p_\natural\nabla\chi_{\mathrm{bet}+}\cdot\omega\Big| \le C_1\varepsilon^{-1}\|p_\natural\|_{L^{6/5}(\func{supp}\nabla\chi_{\mathrm{bet}+})}\|\omega\|_{L^6} \le C(\varepsilon)\,\mathcal{D}_{\mathrm{int}}^{1/2}\|\omega\|_{H^1}.
\end{equation*}
Hardy's inequality does not apply here, because $p_\natural$ is only in $L^{6/5}$ on $\func{supp}\nabla\chi_{\mathrm{bet}+}$. This term, and the terms of the viscous part below in which a derivative falls on $\chi_{\mathrm{bet}+}$, are the only linear terms of Part~A whose constants depend on $\varepsilon$; all of them multiply $\mathcal{D}_{\mathrm{int}}^{1/2}$.

\item \emph{Time derivative.} As in \eqref{eq:pressure_time_ibp}, move one $\theta$-derivative onto the test field. In chart variables, with $\theta$ periodic and $\chi_{\mathrm{bet}+}$ independent of $\theta$,
\begin{equation*}
\int_\Omega J\partial_tv_\natural\cdot w\,dx=-\int_{\mathbb{T}\times\mathcal{W}_{\mathrm{cl}}}\chi_{\mathrm{bet}+}\,\partial_\theta(JJ_\Phi\tilde w)\cdot\partial_t\partial_\theta^mv\,d\theta\,dr\,dh .
\end{equation*}
By Lemma~\ref{lem:tl_coef}, $|\partial_\theta J|\lesssim1$, so this term is $\le C_0\|\chi_{\mathrm{bet}+}\,\partial_t\partial_\theta^mv\|_{L^2}\|\omega\|_{H^1}\le C_0(\mathscr{E}_{1,m,\mathrm{energy,cl}}+\mathcal{D}_{\mathrm{int}})^{1/2}\|\omega\|_{H^1}$.

\item \emph{Viscous term.} Write $\nabla_\mathcal{A} w=\chi_{\mathrm{bet}+}\nabla_\mathcal{A}\tilde w+\tilde w\otimes\mathcal{A}\nabla\chi_{\mathrm{bet}+}$ and $\chi_{\mathrm{bet}+}\mathbb{D}_\mathcal{A} v_\natural=\mathbb{D}_\mathcal{A}(\chi_{\mathrm{bet}+} v_\natural)-(v_\natural\otimes\mathcal{A}\nabla\chi_{\mathrm{bet}+}+\mathcal{A}\nabla\chi_{\mathrm{bet}+}\otimes v_\natural)$. Then
\begin{equation*}
\begin{aligned}
& \frac\mu2\int J\bbD_\cA v_\nat:\bbD_\cA w={} \mu\int J\bbD_\cA(\chis v_\nat):\nabla_\cA\tilde w -\mu\int J\big(v_\nat\otimes\cA\nabla\chis+\cA\nabla\chis\otimes v_\nat\big):\nabla_\cA\tilde w \\
& +\mu\int J\bbD_\cA v_\nat:(\tilde w\otimes\cA\nabla\chis).
\end{aligned}
\end{equation*}
The first term is $\le C_0\|\chi_{\mathrm{bet}+} v_\natural\|_{H^1}\|\omega\|_{H^1}$. The second is $\le C(\varepsilon)\|v_\natural\|_{L^2(\func{supp}\nabla\chi_{\mathrm{bet}+})}\|\omega\|_{H^1}$. In the third, integrate the derivative of $v_\natural$ by parts. By the Piola identity no derivative falls on $J\mathcal{A}$, and there are no boundary terms because $\tilde w\in H^1_0$. The result is a sum of integrals of $v_\natural$ against $\nabla\tilde w\,\nabla\chi_{\mathrm{bet}+}$, $\tilde w\,\nabla^2\chi_{\mathrm{bet}+}$ and $\tilde w\,\nabla\mathcal{A}\,\nabla\chi_{\mathrm{bet}+}$, all supported on $\func{supp}\nabla\chi_{\mathrm{bet}+}$. It is $\le C(\varepsilon)\|v_\natural\|_{L^2(\func{supp}\nabla\chi_{\mathrm{bet}+})}\|\omega\|_{H^1}$. By \eqref{tl:Dint} the viscous term is therefore \(\le C_0\|\chi_{\mathrm{bet}+}\,\partial_\theta^{m+1}v\|_{H^1}\|\omega\|_{H^1}+C(\varepsilon)\,\mathcal{D}_{\mathrm{int}}^{1/2}\|\omega\|_{H^1}.\)

\item \emph{Forcing.} Split $\widetilde F^1=\widetilde F^1_{\mathrm{top}}+\widetilde F^1_{\mathrm{lin}}+\widetilde F^1_{\mathrm{nl}}$. Lemma~\ref{lem:tl_rem}(a) bounds the last two parts by $\big[C_0(\|\chi_{\mathrm{bet}+}\partial_\theta^{m+1}v\|_{H^1}+(\mathscr{E}_{1,m,\mathrm{energy,cl}}+\mathcal{D}_{\mathrm{int}})^{1/2})+C(\varepsilon)\sqrt\mathcal{E}\,\mathcal{D}^{1/2}\big]\|\omega\|_{H^1}$. In $\widetilde F^1_{\mathrm{top}}$, integrate the second-order term by parts in all variables. By the Piola identity there are no boundary terms:
\begin{equation*}
\mu\int_\Omega J\func{div}_\mathcal{A}(\mathbb{D}_\mathfrak{M} v)\cdot w=-\mu\int_\Omega J\,\mathbb{D}_\mathfrak{M} v:\nabla_\mathcal{A} w .
\end{equation*}
Here $|\mathbb{D}_\mathfrak{M} v|\lesssim(|\nabla U|+|U|)\,|\nabla v|$ with $U=\partial_\theta^{m+1}\bar\eta$. Lemma~\ref{lem:tl_first}, applied to $w/\|w\|_{H^1}$, bounds this term and the remaining part $\mathbb{L}$ of $\widetilde F^1_{\mathrm{top}}$ by $C(\varepsilon)\sqrt{\mathcal{E}}\,\mathcal{D}^{1/2}\|\omega\|_{H^1}$. For the slice-wise pieces this uses \eqref{tl:s2}: $\nabla U_{\mathrm{sl}}\in L^2_\theta H^{1-\alpha}_{rh}\subset L^2_\theta L^{2/\alpha}_{rh}$ and $\nabla v\in L^\infty_\theta L^{r_-}_{rh}$, with $\frac\alpha2+\frac{1-\varepsilon_-}{2}<\frac12$. The $3$-D piece is bounded through \eqref{tl:s4}.
\end{itemize}

Altogether, for every $\omega\in H^1_0(\Omega;\mathbb{R}^3)$,
\begin{equation}
\label{tl:pA2}
\Big|\int_\Omega P\func{div}\omega\Big| \le\Big[C_0\Big(\|\chi_{\mathrm{bet}+}\,\partial_\theta^{m+1}v\|_{H^1}+(\mathscr{E}_{1,m,\mathrm{energy,cl}}+\mathcal{D}_{\mathrm{int}})^{1/2}\Big) +C(\varepsilon)\Big(\mathcal{D}_{\mathrm{int}}^{1/2}+\sqrt\mathcal{E}\,\mathcal{D}^{1/2}\Big)\Big]\|\omega\|_{H^1(\Omega)} .
\end{equation}

\emph{Step A3: the Ne\caronaccent{c}as inequality on $\Omega$.} $\Omega$ is a bounded Lipschitz domain. So for $g\in L^2(\Omega)$ with mean $\bar g$,
\begin{equation*}
\|g-\bar g\|_{L^2(\Omega)}\le C_\Omega\sup_{0\ne\omega\in H^1_0(\Omega;\mathbb{R}^3)}\frac{\int_\Omega g\,\func{div}\omega}{\|\omega\|_{H^1}},
\end{equation*}
by the Bogovski\u{\i} operator \cite[\S III.3]{Galdi}; here $\int_\Omega\func{div}\omega=0$ was used. Applied to $g=P$, with \eqref{tl:pA2}, this bounds $\|P-\bar P\|_{L^2(\Omega)}$ by the bracket in \eqref{tl:pA2}.

\emph{Step A4: the mean.} As in \eqref{eq:pressure_mean}, integrate by parts once in $\theta$. Since $\theta$ is periodic and $\chi_{\mathrm{bet}+}$ does not depend on $\theta$,
\begin{equation*}
|\Omega|\,\bar P =\int_{\mathbb{T}\times\mathcal{W}_{\mathrm{cl}}}\chi_{\mathrm{bet}+} J_\Phi\,\partial_\theta(\partial_\theta^mp)\,d\theta\,dr\,dh =-\int_{\mathbb{T}\times\mathcal{W}_{\mathrm{cl}}}\chi_{\mathrm{bet}+}\,(\partial_\theta J_\Phi)\,\partial_\theta^mp\,d\theta\,dr\,dh .
\end{equation*}
Hence $|\bar P|\le C_0\|\chi_{\mathrm{bet}+}\,\partial_\theta^mp\|_{L^2}\le C_0(\mathscr{E}_{1,m,\mathrm{energy,cl}}+\mathcal{D}_{\mathrm{int}})^{1/2}$. Together with Step~A3, and since $C_0\mathcal{D}_{\mathrm{int}}^{1/2}\le C(\varepsilon)\mathcal{D}_{\mathrm{int}}^{1/2}$, this proves \eqref{tl:pmain}. \qed

\smallskip\noindent\textbf{Part B: the free surface.} Part~B proves \eqref{tl:main}. It uses Part~A only through \eqref{tl:pmain}.

\smallskip\noindent\textbf{Step B0: the test field and the normal stress.} Let $\tilde\varphi$ be a Stein extension of $\varphi\in H^{1/2}(D)$ to $\mathbb{R}^2$, and $\mathcal{P}_3$ the downward $3$-D Poisson extension. Let $\vartheta\in C_c^\infty(-c_0,c_0)$ equal $1$ near $0$, and set
\begin{equation*}
\widetilde E\varphi(x,z):=\vartheta\big(z-\zeta(x)\big)\,(\mathcal{P}_3\tilde\varphi)\big(x,z-\zeta(x)\big).
\end{equation*}
By (H1), $\widetilde E:H^{1/2}(\Sigma)\to H^1(\Omega)$ is bounded and $\widetilde E\varphi|_\Sigma=\varphi$. For $\varphi\in H^1(\Sigma)$ put \(w:=\chi\,(\widetilde E\varphi)\,e_3 .\) This field is admissible: we have  $\func{supp}w\subset\mathcal{U}$;  $w\cdot n_0=0$ on $\Sigma_s$, since $e_3\cdot n_0=0$ wherever $\widetilde E\varphi\ne0$; and  $w\cdot\mathcal{N}|_\Sigma=\chi_\Sigma\varphi\in H^1(\Sigma)$, because $e_3\cdot\mathcal{N}=1$. Moreover $\|w\|_{H^1(\Omega)}\le C(\varepsilon)\|\varphi\|_{1/2}$; the constant depends on $\varepsilon$ through $\nabla\chi$, which is harmless here. Insert $w$ into \eqref{tl:weak}. Since $\chi_\Sigma=1$ near $\Gamma_0$, we have $\partial_t\eta_\natural=\partial_tf$ on $\Gamma_0$. Moreover
\begin{equation*}
\mathcal{L}\nabla\eta_\natural\cdot\nabla(\chi_\Sigma\varphi) =\mathcal{L}\nabla f\cdot\nabla\varphi-\mathcal{L}(\eta_\natural\nabla\chi_\Sigma)\cdot\nabla\varphi +(\mathcal{L}\nabla\eta_\natural\cdot\nabla\chi_\Sigma)\varphi ,
\end{equation*}
where $\mathcal{L}=\mathcal{L}_{\Sigma,0}$. Hence, for all $\varphi\in H^1(\Sigma)$,
\begin{equation}  \label{tl:surf} g\int_\Sigma f\varphi+\sigma\int_\Sigma\mathcal{L}\nabla f\cdot\nabla\varphi=\Lambda(\varphi)+\Lambda_1(\varphi)+\int_\Sigma B\cdot\nabla\varphi+\int_{\Gamma_0}h\,\varphi ,
\end{equation}
where
\begin{equation*}
\begin{aligned}
& \Lambda(\varphi) :=\int_\Omega J(\widetilde F^1-\partial_tv_\natural)\cdot w+\int_\Omega J\,S_\mathcal{A}(p_\natural,v_\natural):\nabla_\mathcal{A} w -\beta_0\int_{\Sigma_s}J_sU_0v_\natural\cdot U_0w -\langle\widetilde F^3,w\rangle_\Sigma-\langle J_s\widetilde F^6,U_0w\rangle_{\Sigma_s}, \\
& \Lambda_1(\varphi) :=-\sigma\int_\Sigma\big(\mathcal{L}\nabla\eta_\natural\cdot\nabla\chi_\Sigma+\widetilde F^5\cdot\nabla\chi_\Sigma\big)\varphi,\quad B :=\sigma\,\mathcal{L}\,\eta_\natural\nabla\chi_\Sigma-\sigma\chi_\Sigma\widetilde F^5,\qquad\qquad h:=-\kappa\,\partial_tf-\widetilde F^7 .
\end{aligned}
\end{equation*}

\smallskip\noindent\textbf{Step B1: the normal stress is in $(H^{1/2}(\Sigma))^*$.} We claim that, for $\varphi\in H^1(\Sigma)$,
\begin{equation}  \label{tl:step1} 
|\Lambda(\varphi)|\le C(\varepsilon)\,\mathcal{D}^{1/2}\|\varphi\|_{1/2}.
\end{equation}
Split $\Lambda$ into three groups.

\emph{(a) The principal linear terms.} By (H1), $\chi_{\mathrm{bet}+}=1$ and $\nabla\chi_{\mathrm{bet}+}=0$ on $\func{supp}\chi$, and $\chi_{\mathrm{bet}+}$ does not depend on $\theta$. Hence, on $\func{supp}\chi$, $\partial_tv_\natural=\partial_\theta(\chi_{\mathrm{bet}+}\,\partial_t\partial_\theta^mv)$ and $\nabla v_\natural=\nabla(\chi_{\mathrm{bet}+} v_\natural)$, so that
\begin{equation*}
\|\partial_tv_\natural\|_{L^2(\func{supp}\chi)}\lesssim\|\chi_{\mathrm{bet}+}\,\partial_t\partial_\theta^mv\|_{H^1}\lesssim\mathcal{D}^{1/2},\quad \|v_\natural\|_{H^1(\func{supp}\chi)}\lesssim\mathcal{D}^{1/2},
\end{equation*}
by \eqref{tl:VD}; here we use that the top level is purely angular, as in Steps~A2 and~A4. Together with the trace theorem on $\Sigma_s$,
\begin{equation*}
\begin{aligned}
&\Big|\int J\partial_tv_\nat\cdot w\Big|+\Big|\int JS_\cA(p_\nat,v_\nat):\nabla_\cA w\Big|+\beta_0\Big|\int_{\Sigma_s}J_sU_0v_\nat\cdot U_0w\Big|\\
&\qquad\lesssim\big(\cD^{1/2}+\|\chis p_\nat\|_{L^2}\big)\|w\|_{H^1}\le C(\varepsilon)\,\cD^{1/2}\|w\|_{H^1}.
\end{aligned}
\end{equation*}
The last step is Part~A.

\emph{(b) Remainders.} By Lemma~\ref{lem:tl_rem}(b), these terms are bounded by $C(\varepsilon)\mathcal{D}^{1/2}\|w\|_{H^1}$.

\emph{(c) The top terms.} Since $w=w_3e_3$, we have $\langle\widetilde F^3_{\mathrm{top}},w\rangle_\Sigma=\langle e_3\cdot\widetilde F^3_{\mathrm{top}},w_3\rangle_\Sigma$. Split $\widetilde F^1_{\mathrm{top}}=\mu\func{div}_\mathcal{A}(\mathbb{D}_\mathfrak{M} v)+\mathbb{L}$ as in Lemma~\ref{lem:tl_first}, and integrate the second-order term by parts in all variables, as in Step~A2. Now $w$ does not vanish on $\partial\Omega$. For a regular solution the Piola identity gives $(J\mathcal A n_\Sigma)\,dS_{\mathrm{ref}}=\mathcal N\,dx$ on $\Sigma$, with $n_\Sigma=N_0/|N_0|$, and $J\mathcal A n_0=J_s n_0$ on $\Sigma_s$. Consequently,
\begin{equation*}
\int_\Omega J\mu\func{div}_\mathcal{A}(\mathbb{D}_\mathfrak{M} v)\cdot w =-\mu\int_\Omega J\,\mathbb{D}_\mathfrak{M} v:\nabla_\mathcal{A} w +\langle\mu\mathbb{D}_\mathfrak{M} v\,\mathcal{N},w\rangle_\Sigma+\langle J_s\mu\mathbb{D}_\mathfrak{M} v\,n_0,w\rangle_{\Sigma_s}.
\end{equation*}
On $\Sigma$, $w=w_3e_3$, so the face term there is $\langle\mu\,e_3\cdot\mathbb{D}_\mathfrak{M} v\,\mathcal{N},w_3\rangle_\Sigma$. On $\Sigma_s$, $w\cdot n_0=0$, so $w=U_0w$ and the face term there is $\langle J_s\mu U_0\mathbb{D}_\mathfrak{M} v\,n_0,U_0w\rangle_{\Sigma_s}=\langle J_s\widetilde F^6_{\mathrm{top}},U_0w\rangle_{\Sigma_s}$. Both cancel exactly against the corresponding parts of $-\langle\widetilde F^3_{\mathrm{top}},w\rangle_\Sigma-\langle J_s\widetilde F^6_{\mathrm{top}},U_0w\rangle_{\Sigma_s}$; see \eqref{tl:top}. By \eqref{tl:F3vert}, what remains is
\begin{equation*}
\begin{aligned}
&\int_\Omega J\Ft^1_{\rm top}\cdot w-\langle\Ft^3_{\rm top},w\rangle_\Sigma-\langle J_s\Ft^6_{\rm top},U_0w\rangle_{\Sigma_s}\\
&\qquad=-\mu\int_\Omega J\,\bbD_\fM v:\nabla_\cA w+\int_\Omega J\,\mathbb L\cdot w-\mu\Big\langle\sum_{j=1,2}(\bbD_\cA v)_{3j}\mathcal M_j,\,w_3\Big\rangle_\Sigma .
\end{aligned}
\end{equation*}
The last term is the difference of two face terms bounded in Lemma~\ref{lem:tl_first}, and the lemma bounds the other two terms as well. Hence the right-hand side is $\le C(\varepsilon)\sqrt{\mathcal{E}}\,\mathcal{D}^{1/2}\|w\|_{H^1}$; here $\func{supp}w\subset\func{supp}\chi_{\mathrm{cl}}$, so the face terms of Lemma~\ref{lem:tl_first} apply. No second derivative of $\partial_\theta^{m+1}\bar\eta$ enters.

Collecting (a)--(c) and using $\|w\|_{H^1}\le C(\varepsilon)\|\varphi\|_{1/2} $ proves \eqref{tl:step1}. Step~B1 does not use $f\in H^{3/2}(\Sigma)$, and $\Lambda$ extends uniquely to $(H^{1/2}(\Sigma))^*$.

\smallskip\noindent\textbf{Step B2: the gain of $H^{3/2}$, and absorption.} We bound the remaining data in \eqref{tl:surf}.

\emph{The term $\Lambda_1$.} Here $\widetilde F^5=\widetilde F^5_{\mathrm{top}}+\widetilde F^5_{\mathrm{rem}}$, with $|\widetilde F^5_{\mathrm{top}}|\lesssim|\nabla\eta||\nabla\eta_\natural|$, and $\eta_\natural=\chi_{\mathrm{bet}+,\Sigma}\eta_\natural\in H^{3/2-\alpha}$ on $\func{supp}\chi_\Sigma$ by \eqref{tl:VD} and (H1). Using Lemma~\ref{lem:tl_rem}(c), \(|\Lambda_1(\varphi)|\le C(\varepsilon)\,\mathcal{D}^{1/2}\|\varphi\|_{L^2(\Sigma)} .\)

\emph{The divergence-form datum $B$.} Since $\chi_\Sigma\widetilde F^5_{\mathrm{top}}=\mathcal{R}^{\prime }(\nabla\eta)\nabla f-\mathcal{R}^{\prime }(\nabla\eta)\,\eta_\natural\nabla\chi_\Sigma$,
\begin{equation*}
B=-\sigma\,\mathcal{R}^{\prime }(\nabla\eta)\nabla f+B_\flat,\quad B_\flat:=\sigma\big(\mathcal{L}+\mathcal{R}^{\prime }(\nabla\eta)\big)\eta_\natural\nabla\chi_\Sigma-\sigma\chi_\Sigma\widetilde F^5_{\mathrm{rem}}.
\end{equation*}
The coefficient $a:=\mathcal{R}^{\prime }(\nabla\eta)=DF_\zeta(\nabla\zeta+\nabla\eta)-DF_\zeta(\nabla\zeta)$ is smooth in $(\nabla\zeta,\nabla\eta)$ and vanishes at $\nabla\eta=0$. By Lemma~\ref{lem:tl_coef} it satisfies, on $\Sigma\cap\mathcal{U}$,
\begin{equation*}
\|a\|_{L^\infty}+\|\partial_\theta a\|_{L^\infty}+\sup_\theta\|\partial_ra(\theta,\cdot)\|_{L^{\bar s}_r}\lesssim\sqrt{\mathcal{E}},\quad \bar s:=\min\{s_-,3/2\}\in(1,2).
\end{equation*}
Apply Lemma~\ref{lem:tl_mult} in the chart, with $F$ the extension by $E$ of the chart representation of $\nabla f$, and with $a$ replaced by $\tilde a:=\vartheta_1(r)\,Ea$. Here $\vartheta_1\in C_c^\infty(\mathbb{R})$ is supported where $Ea$ is defined and equals $1$ on the $r$-range of $\func{supp}\chi_\Sigma$, so that $\tilde aF=a\nabla f$ for $r>0$. The cutoff is needed because Lemma~\ref{lem:tl_mult} requires $\partial_r\tilde a\in L^{\bar s}$ on all of $\mathbb{R}$. Since $E$ is bounded on $L^\infty$, $W^{1,\bar s}$ and $H^{1/2}$, and commutes with $\partial_\theta$, this gives \(\|\mathcal{R}^{\prime }(\nabla\eta)\nabla f\|_{1/2}\lesssim\sqrt{\mathcal{E}}\,\|\nabla f\|_{1/2}\lesssim\sqrt{\mathcal{E}}\,\|f\|_{3/2},\) and in the same way, since $\eta_\natural\nabla\chi_\Sigma=(\chi_{\mathrm{bet}+,\Sigma}\eta_\natural)\nabla\chi_\Sigma$ and by Lemma~\ref{lem:tl_rem}(c), $\|B_\flat\|_{1/2}\le C(\varepsilon)\big(\|\chi_{\mathrm{bet}+,\Sigma}\eta_\natural\|_{H^{1/2}}+\mathcal{D}^{1/2}\big)\le C(\varepsilon)\mathcal{D}^{1/2}$.

\emph{The conormal datum $h$.} Since $\hat{\mathcal{W}}(z)=O(z^2)$, we have $|\widetilde F^7_{\mathrm{top}}|\lesssim|\partial_t\eta|\,|\partial_tf|$ on $\Gamma_0$. Using Lemma~\ref{lem:tl_coef} and Lemma~\ref{lem:tl_rem}(c),
\begin{equation*}
\|h\|_{L^2(\Gamma_0)} \le\kappa\|\partial_tf\|_{L^2(\Gamma_0)}+C\sqrt{\mathcal{E}}\,\|\partial_tf\|_{L^2(\Gamma_0)}+\|\widetilde F^7_{\mathrm{rem}}\|_{L^2(\Gamma_0)} \le C(\varepsilon)\,\mathcal{D}^{1/2} .
\end{equation*}

\emph{Conclusion.} By Step~B1 and Lemma~\ref{lem:tl_ell}, applied to \eqref{tl:surf} with $u=f$,
\begin{equation*}
\|f\|_{3/2}\le C_1(\varepsilon)\,\mathcal{D}^{1/2}+C_2(\varepsilon)\sqrt{\mathcal{E}}\,\|f\|_{3/2}.
\end{equation*}
The last term comes only from $\mathcal{R}^{\prime }(\nabla\eta)\nabla f$ in $B$. This is an a~priori estimate, so $\|f\|_{3/2}<\infty$ by (H3). Choose $\delta_0(\varepsilon)$ with $C_2(\varepsilon)\sqrt{\delta_0(\varepsilon)}\le\frac12$. Then the last term is absorbed, and $\|f\|_{3/2}\le2C_1(\varepsilon)\mathcal{D}^{1/2}$, which is \eqref{tl:main}. \qed

\smallskip\noindent\textbf{Proof of Corollary~\ref{cor:top_level}.} (i) The trace $H^{3/2}(D)\to H^1(\partial D)$ is bounded, and $f|_{\Gamma_0}=\partial_\theta^{m+1}\eta|_{\Gamma_0}$ because $\chi_\Sigma=1$ near $\Gamma_0 $. Let $s$ denote planar arclength and recall that $d\theta=a(s)\,ds$, where $a(s)=\sqrt{1+|\nabla\zeta(\gamma(s))\cdot\dot\gamma(s)|^2}$. Thus $\partial_{\theta}=a(s)^{-1}\partial_{\gamma}$ after pullback. Since $a$ and $a^{-1}$ are bounded by the equilibrium geometry,
\begin{equation*}
\|\partial_\theta^{m+2}\eta\|_{L^2(\Gamma_0)} +\|\partial_\gamma(f|_{\partial D})\|_{L^2(\partial D)}  \le C_\zeta\|f|_{\partial D}\|_{H^1(\partial D)} \le C_\zeta\|f\|_{3/2}\le C(\varepsilon)\mathcal{D}^{1/2}.
\end{equation*}
(ii) Square \eqref{tl:main} and integrate over $[s,t]$.

(iii) In Subsection~\ref{subsec:cl_local_energy}, $\eta_{0,m+1}=\partial_\theta^{m+1}(\chi_{\mathrm{bet+,\Sigma}}\eta)$ coincides with $f$ on $\Gamma_0$. Moreover $\|v\|_{L^\infty(\Gamma_0)}\lesssim\sqrt\mathcal{E}$ by Lemma~\ref{lem:tl_coef}, and $\|\partial_tf\|_{L^2(\Gamma_0)}\lesssim\mathcal{D}^{1/2}$ by \eqref{tl:VD}. By (i),
\begin{equation*}
\big|\kappa\langle\partial_tf,v_\gamma\partial_\gamma f\rangle_{\Gamma_0}\big| \le C_\zeta\kappa\|v\|_{L^\infty(\Gamma_0)}\|\partial_tf\|_{L^2(\Gamma_0)}\|\partial_\theta^{m+2}\eta\|_{L^2(\Gamma_0)} \le C(\varepsilon)\sqrt{\mathcal{E}}\,\mathcal{D} .
\end{equation*}
The term $[\widetilde F^4]_{\Gamma_0}$ of Step~7 of the proof of Theorem~\ref{thm:contact_line_energy} contains the same factor $v_\gamma\partial_\gamma f $ and is bounded in the same way. Indeed, since the top transport component satisfies \(\|v_\gamma\partial_\gamma f\|_{L^2(\Gamma_0)} \le C(\varepsilon)\sqrt{\mathcal E}\,\mathcal D^{1/2},\) if $\|g\|_{L^2(\Gamma_0)}\le C(\varepsilon)\mathcal D^{1/2}$, we have \(|\langle g,v_\gamma\partial_\gamma f\rangle_{\Gamma_0}| \le C(\varepsilon)\sqrt{\mathcal E}\,\mathcal D.\) Apply this with $g=\kappa\partial_tf$. The other components of $\widetilde F^4$ keep their own estimates. \qed

\appendix
\section{Surface functional calculus and admissible three-dimensional lifts}
\label{functional:section}

This section supplies the analytic test-function construction for the
$H^{3/2-\alpha}$ surface estimate.  The surface is two-dimensional and the
fluid region is three-dimensional.  

\subsection{Geometric assumptions and conventions}

Let $D\subset\mathbb R^2$ be a bounded, connected $C^\infty$ domain and let
$\zeta\in C^\infty(\overline D)$.  Let $\Omega\subset\mathbb R^3$ be a
bounded, connected Lipschitz domain with free face $\Sigma=\{(x,\zeta(x)):x\in D\}$ and solid face $\Sigma_s$. The fluid lies locally below this graph, so the upward graph normal
is outward.  The solid face is
smooth near $\partial\Sigma$, meets $\Sigma$ transversely with interior
contact angle strictly between $0$ and $\pi$, and has no other boundary
junction in a neighborhood of $\overline\Sigma$.  Locally the fluid is
the intersection of the inward sides of the two faces; reentrant
contact wedges are excluded.  Write $\mathcal N_0=(-\nabla\zeta,1)$, $dS_\Sigma=|\mathcal N_0|\,dx$, $\nu_\Sigma=\mathcal N_0/|\mathcal N_0|$. The constants in this section may depend on this fixed geometry.  They are
uniform in the spectral truncation index.

For $r\geq0$, $H^r(D)$ denotes the restriction Sobolev space.  The notation $(H^r(D))'$ always denotes its full continuous dual with the $L^2(D)$ pivot.  In particular, $(H^{1/2}(D))'$ is used literally; it is not silently
identified with the dual of $H^{1/2}_{00}(D)$.  For $r\geq0$ the subscript
$0,\mathrm{av}$ denotes the zero-average subspace; for a dual distribution
the average means its pairing with the constant function $1$.

\subsection{The surface operator and its fractional scale}

Fix $g,\sigma>0$, and set
\[
 G_\zeta=(1+|\nabla\zeta|^2)^{-1/2}I
 -(1+|\nabla\zeta|^2)^{-3/2}\nabla\zeta\otimes\nabla\zeta,
 \qquad
 a(u,v)=\int_D\bigl(guv+\sigma G_\zeta\nabla u\cdot\nabla v\bigr)\,dx.
\]
This is the Hessian form of graph area, with gravity included.  The same
construction applies to any fixed smooth symmetric uniformly positive
matrix in place of $G_\zeta$, provided that this matrix is the one in the
linearized capillary form of the weak equation.  The eigenvalues of
$G_\zeta$ are $(1+|\nabla\zeta|^2)^{-1/2}$ in the direction orthogonal to $\nabla\zeta$ and $(1+|\nabla\zeta|^2)^{-3/2}$ in its direction.  Thus $a$ is an equivalent Hilbert norm on $H^1(D)$.

\begin{proposition}[Surface realization and fractional domains]
\label{functional:surface_scale}
The form $a$ defines a positive self-adjoint operator $\mathcal K$ on
$L^2(D)$ with compact inverse.  Its realization is
\[
 \mathcal Ku=gu-\sigma\operatorname{div}(G_\zeta\nabla u),\qquad
 D(\mathcal K)=\{u\in H^2(D):\nu_D\cdot G_\zeta\nabla u=0
                    \text{ on }\partial D\}.
\]
There is a real orthonormal eigenbasis $(e_j)_{j\geq0}$ with
$\mathcal Ke_j=\lambda_je_j$, $e_0=|D|^{-1/2}$, $\lambda_0=g$, and $\lambda_j>g$ for $j\geq1$.  Define, for every real $r$, \[
 \|u\|_{\mathcal H_{\mathcal K}^r}^2
 =\sum_{j\geq0}\lambda_j^r|u_j|^2.
\]
For $r\geq0$ this is the operator-domain scale; for $r<0$ it is its
completion, or equivalently the dual scale under the $L^2$ pairing.
One has, with equivalence of norms,
\begin{align}
 \mathcal H_{\mathcal K}^r&=H^r(D),&&0\leq r<3/2,\label{functional:subcritical}\\
 \mathcal H_{\mathcal K}^{-r}&=(H^r(D))',&&0\leq r<3/2.\label{functional:negative}
\end{align}
For $3/2<r\leq2$, the positive space is $H^r(D)$ with homogeneous conormal
condition.  At $r=3/2$ the condition is instead
\[
 u\in H^{3/2}(D),\qquad
 \int_{D\cap\{d<d_0\}}\frac{|B_1u|^2}{d}\,dx<\infty,
 \quad d(x)=\operatorname{dist}(x,\partial D),
\]
where $B_1$ is a smooth extension into a boundary collar of the conormal
operator $\nu_D\cdot G_\zeta\nabla$.  The norm is equivalent to the sum of
the $H^{3/2}$ norm and the square root of this weighted integral.
\end{proposition}

\begin{proof}
Coercivity gives the inverse by the form construction; compactness follows
from $H^1(D)\Subset L^2(D)$.  The spectral theorem gives the stated basis.
The equality case in $a(u,u)\geq g\|u\|_2^2$ is precisely the constant
functions, because $D$ is connected.  Green's formula and smooth-boundary
elliptic regularity identify the displayed realization.  Smooth elliptic
regularity also makes each eigenfunction smooth up to the boundary.
The fractional-domain identification is the smooth elliptic Neumann
interpolation theorem: apply Theorem 2.2 and Definition 2.1 of G. Grubb,
\emph{Regularity of spectral fractional Dirichlet and Neumann problems},
Math. Nachr. 289 (2016), 831--844, with $p=2$, boundary order $j=1$, and
operator power $r/2$.  Its threshold is $j+1/p=3/2$; its critical condition
is the weighted condition above.  The assumptions apply to this smooth,
uniformly elliptic, invertible conormal realization.  Two extensions of the
conormal vector agreeing on the boundary differ by $O(d)$, so their weighted
conditions differ by a term controlled by $\|\nabla u\|_2$.  Duality gives
\eqref{functional:negative}; it is duality with the full positive space.
\end{proof}

\begin{proposition}[Finite spectral tests]
\label{functional:truncations}
Let $P_Nu=\sum_{j=0}^Nu_je_j$ and $D_N^s=\mathcal K^{s/2}P_N$.  For every real $r,s$, \[
 \|D_N^su\|_{\mathcal H_{\mathcal K}^{r-s}}
 \leq\|u\|_{\mathcal H_{\mathcal K}^{r}},\qquad
 D_N^su\longrightarrow\mathcal K^{s/2}u
 \quad\text{in }\mathcal H_{\mathcal K}^{r-s}.
\]
Both $P_N$ and $D_N^s$ preserve zero average, and $\int_DD_N^su=g^{s/2}\int_Du$.  For $u\in H^1(D)$ and $0<s<1$, \begin{equation}\label{functional:exact_spectral_identity}
 a(u,D_N^su)=\sum_{j=0}^N\lambda_j^{1+s/2}|u_j|^2
 =\|D_N^su\|_{\mathcal H_{\mathcal K}^{1-s/2}}^2.
\end{equation}
The left side increases to $\|u\|_{\mathcal H_{\mathcal K}^{1+s/2}}^2$,
possibly $+\infty$.
\end{proposition}

\begin{proof}
Each statement follows by multiplying coefficients in the orthonormal
basis.  For the form identity use $a(u,e_j)=\lambda_j(u,e_j)$, which holds
for every $u\in H^1(D)$ by the weak eigenfunction equation.  No conormal
boundary condition on $u$ is required.  The convergence follows from the
tail of a convergent nonnegative series.  For the last assertion use
monotone convergence, including the case of infinite limit.
\end{proof}

Fix hereafter $0<\alpha<1/2$, and put $s=1-2\alpha$, $q=s/2=1/2-\alpha, \qquad r=s-1/2=1/2-2\alpha$. Then $1+s/2=3/2-\alpha$ is below the conormal threshold, and
$1-s/2=1/2+\alpha$ is above the trace threshold.  The required surface
trace is $H^{1/2+\alpha}(D)\longrightarrow H^\alpha(\partial D) \longrightarrow L^2(\partial D)$. There is no assertion of an $L^\infty$ trace.  In the interior of $D$,
$H^t(D)\hookrightarrow L^{2/(1-t)}(D)$ for $0<t<1$; in particular $\nabla u\in L^{4/(1+2\alpha)}(D)$ if $u\in H^{3/2-\alpha}(D)$. 
The curvature-remainder pairing used below satisfies
\begin{equation}\label{functional:curvature_pairing}
 |\langle R,\nabla h\rangle_D|
 \leq C\|R\|_{H^q(D)}\|h\|_{H^{1-q}(D)},\qquad 0<q<1/2.
\end{equation}
For completeness, extend $h$ boundedly to $H^{1-q}(\mathbb R^2)$ and
extend $R$ by zero to $H^q(\mathbb R^2)$.  Differentiation maps the first
space to $H^{-q}(\mathbb R^2)$, giving the bound by Fourier duality.
For the zero-extension assertion, interpolate its bounded maps on
$L^2(D)$ and $H^1_0(D)$.  The interpolated domain equals $H^q(D)$ for
$q<1/2$: apply the Dirichlet case of the same fractional-domain theorem
to $I-\Delta_D$, whose square-root form domain is $H^1_0(D)$.
Approximation of $R$ by compactly supported smooth functions
in $H^q(D)$ shows independence of the extension of $h$.  For smooth inputs
this is exactly the integral $\int_D R\cdot\nabla h$, with no integration
by parts across the contact boundary.

\subsection{A solenoidal collar lift with both endpoint estimates}

The lifting used here is a vector-field lifting.  A scalar harmonic
potential is unnecessary: the energy argument needs a divergence-free
field, the prescribed normal trace, and the two estimates below.

\begin{lemma}[Fixed collar]
\label{functional:collar}
Under the geometric assumptions above, there exist $a>0$ and a smooth
diffeomorphism $F:\overline D\times[-a,0]\to\overline\Omega$ onto a collar
of $\overline\Sigma$, with $F(x,0)=(x,\zeta(x))$, positive Jacobian, and $F(\partial D\times[-a,0])\subset\Sigma_s$. \end{lemma}

\begin{proof}
At a point of the contact curve, the projection of the outward free-face
normal onto the tangent plane of the solid face is nonzero and has
positive dot product with that normal.  Transversality and compactness
make this positivity uniform.  Extend this vector smoothly near the
contact curve, tangent to the solid face there, and blend it on the free
face with an outward transverse field in its interior.  A smooth
partition of unity preserves the strict outward dot product.  Extend to
a neighborhood, retaining tangency to the solid face.  Its flow for
negative times points into $\Omega$; uniqueness of the flow preserves
the solid face.  The inverse-function theorem, transversality, and
compactness give a common time $a>0$ on which the flow from
$\overline\Sigma$ is an embedded collar.  Taking the flow time as the
third coordinate gives $F$.  At time zero its positive determinant is
the dot product of the transverse field with $\mathcal N_0$.
\end{proof}

For the vertical-wall vessel used here, the collar is even simpler:
$F(x,z)=(x,\zeta(x)+z)$, with $a$ smaller than the fixed top-to-bottom
separation.  Its Jacobian is one, and its lateral face lies on the
vertical wall.  Thus this application does not require a bulk
corner-regularity theorem or a change of the vessel class.

\begin{theorem}[Flux lifting]
\label{functional:flux_lift}
There is one fixed linear operator $\mathcal L$, defined on zero-average
surface data, such that
\begin{align}
 \|\mathcal Lf\|_{L^2(\Omega)}
 &\leq C\|f\|_{(H^{1/2}(D))'},\label{functional:lift_low}\\
 \|\mathcal Lf\|_{H^1(\Omega)}
 &\leq C\|f\|_{H^{1/2}(D)}.\label{functional:lift_high}
\end{align}
It satisfies
\[
 \operatorname{div}\mathcal Lf=0,\qquad
 \mathcal Lf\cdot\nu=0\ \text{on }\Sigma_s,\qquad
 (\mathcal Lf)(x,\zeta(x))\cdot\mathcal N_0(x)=f(x).
\]
For low-regularity data these identities are distributional normal-trace
identities.  Equivalently, for every $\varphi\in H^1(\Omega)$,
\begin{equation}\label{functional:lift_green}
 \int_\Omega\mathcal Lf\cdot\nabla\varphi
 =\langle f,\varphi(x,\zeta(x))\rangle_D.
\end{equation}
\end{theorem}

\begin{proof}
Let $L=-\Delta_N$ on $D$, restricted to the zero-average subspace.
Its eigenpairs $(\mu_j,\phi_j)_{j\geq1}$ have $\mu_j\geq\mu_1>0$.
The operator $I+L$ satisfies Proposition~\ref{functional:surface_scale};
on the zero-average subspace its powers are equivalent to those of $L$.
Thus $\sum\mu_j^{1/2}|f_j|^2$ and $\sum\mu_j^{-1/2}|f_j|^2$ give the
two squared data norms in the theorem, respectively.  In the negative
case, this assertion uses the full $H^{1/2}$ dual and the zero-average
projection, which is bounded because pairing with $1$ is bounded.

Choose $\chi\in C^\infty([-a,0])$ equal to $1$ near $0$ and zero near
$-a$.  Initially for a finite sum $f=\sum f_j\phi_j$, define
\[
 h(x,z)=\chi(z)e^{z\sqrt L}f(x),\qquad
 Z_f(x,z)=\bigl(\nabla_xL^{-1}\partial_zh(x,z),h(x,z)\bigr).
\]
The mean of $h$ and of $\partial_zh$ is zero.  Since
$\Delta L^{-1}=-I$ on this subspace, $\operatorname{div}_{x,z}Z_f=-\partial_zh+\partial_zh=0$. The lateral flux vanishes by the Neumann condition on $L^{-1}\partial_zh$.
The upper vertical flux is $f$, and the field vanishes near the lower
face of the collar.

Here are the estimates including the cutoff terms.  With $t_j=\sqrt{\mu_j}$
and $h_j(z)=\chi(z)e^{zt_j}f_j$, for $k=0,1,2$,
\begin{equation}\label{functional:mode_estimate}
 \int_{-a}^0|\partial_z^kh_j(z)|^2\,dz
 \leq C_k t_j^{2k-1}|f_j|^2.
\end{equation}
The term with $k$ derivatives on the exponential follows from $\int_{-a}^0e^{2zt_j}\,dz\leq(2t_j)^{-1}$. Every derivative of $\chi$ is supported a positive distance below zero;
its contribution is bounded by $C t_j^{2\ell}e^{-ct_j}|f_j|^2$ with
$\ell<k$.  Since $t_j\geq\sqrt{\mu_1}>0$, this is at most the stated
right side.  Orthogonality and
$\|\nabla L^{-1}u\|_2^2=\sum\mu_j^{-1}|u_j|^2$ imply
\[
 \|Z_f\|_{L^2(D\times(-a,0))}^2
 \leq C\sum_j\mu_j^{-1/2}|f_j|^2.
\]
For the $H^1$ estimate, $\partial_zh$ and $\nabla_xh$ cost at most
$C\sum_j\mu_j^{1/2}|f_j|^2$.  The same is true for
$\partial_z\nabla L^{-1}\partial_zh$, using $k=2$ in
\eqref{functional:mode_estimate}.  Finally, smooth planar Neumann
regularity, or the $H^2$ realization of $I+L$, gives $\|\nabla^2L^{-1}u\|_2\leq C\|u\|_2$ for zero-average $u$. Apply this with $u=\partial_zh(\cdot,z)$ and integrate in $z$.
Consequently
\[
 \|Z_f\|_{H^1(D\times(-a,0))}^2
 \leq C\sum_j\mu_j^{1/2}|f_j|^2.
\]
Notice that this is smooth-domain elliptic regularity on the planar
surface $D$, not an elliptic regularity assertion on the bulk corner.

Transport by the contravariant Piola map:
\[
 (\mathcal Lf)(F(x,z))=\frac{DF(x,z)Z_f(x,z)}{\det DF(x,z)}
\]
inside the collar, and extend by zero outside it.  There is no jump at
its lower interface because $Z_f$ vanishes in a neighborhood of that
interface.  The smooth fixed map preserves both estimates.  To verify
all flux conventions at once, change variables and integrate:
\[
 \int_\Omega\mathcal Lf\cdot\nabla\varphi
 =\int_{D\times(-a,0)}Z_f\cdot\nabla_{x,z}(\varphi\circ F)
 =\int_Df(x)\varphi(F(x,0))\,dx.
\]
This proves the divergence, wall, and top identities.  In particular
the top unit-normal flux is $f/|\mathcal N_0|$.  Extend the same operator
by completion in each endpoint norm; the extensions agree on their
common domain by the $L^2$ estimate.  The trace map from $H^1(\Omega)$
to $H^{1/2}(D)$ on the smooth top face is bounded, so the Green identity
passes to the low-regularity limit.
\end{proof}

\subsection{Moving-coordinate test fields}

Let $\Phi(t):\Omega\to\Omega(t)$ be an orientation-preserving smooth
diffeomorphism preserving the solid boundary.  In fixed coordinates put $J=\det D\Phi$, $\mathcal A=(D\Phi)^{-T}$, $M=J^{-1}\mathcal A^{-T}=J^{-1}D\Phi$, $\mathcal N=J\mathcal A\mathcal N_0$. Assume on the time interval in question
\begin{equation}\label{functional:geometry_bounds}
 \|M\|_{W^{1,\infty}(\Omega)}+\|J\|_{W^{1,\infty}(\Omega)}
 +\|\mathcal A\|_{L^\infty(\Omega)}+\|J^{-1}\|_{L^\infty(\Omega)}
 \leq C_G,
\end{equation}
and $D\Phi$, $\partial_tD\Phi$ have the $L^\infty$ bounds used below.
These are explicit coefficient assumptions; their deduction from a
particular energy functional is a separate geometric estimate.
The $H^1$ multiplier bound on $M$ may instead be obtained from $M\in L^\infty$ and $\nabla M\in L^3$, by $H^1(\Omega)\hookrightarrow L^6$. 
For a zero-average datum $f$, define $\mathcal W_t f=M(t)\mathcal Lf$.
The Piola identity and the normal transformation give
\begin{align}
 J\operatorname{div}_{\mathcal A}(\mathcal W_tf)&=0,
 &\mathcal W_tf\cdot(\nu_s\circ\Phi)&=0\text{ on the fixed solid face},
 &\mathcal W_tf\cdot\mathcal N&=f,\label{functional:admissible}\\
 \|\mathcal W_tf\|_{H^1}&\leq C_GC\|f\|_{H^{1/2}},
 &\|\mathcal W_tf\|_2&\leq C_GC\|f\|_{(H^{1/2})'}.
\end{align}
Here the solid condition uses the pulled-back physical normal, because
$D\Phi(x)$ maps the tangent plane at $x$ to the tangent plane at
$\Phi(x)$.  In the present vertical-wall, vertical-shift geometry,
$\nu_s(\Phi(x))=\nu_s(x)$, so it is precisely the usual fixed wall
condition.  In coordinates,
$J\operatorname{div}_{\mathcal A}(Mz)=\operatorname{div}z$, since
$J\mathcal A^TM=I$ and the cofactor columns have zero divergence.
Also $(Mz)\cdot\mathcal N=z\cdot\mathcal N_0$.

Given an arbitrary surface function $f$, set
\begin{equation}\label{functional:mean_projection}
 f^\circ=f-|D|^{-1}\langle f,1\rangle_D,\qquad
 g_N=D_N^sf^\circ,\qquad w_N=\mathcal W_tg_N.
\end{equation}
Each $g_N$ is smooth, has zero average, and is an admissible top trace
even when the weak formulation requires the top normal trace to lie
in $H^1(D)$.  There is no need for $f$ itself to satisfy a conormal
boundary condition.  Since $g_N$ has zero average,
\begin{equation}\label{functional:XN}
 X_N^2:=a(f,g_N)=a(f^\circ,g_N)
 =\|g_N\|_{\mathcal H_{\mathcal K}^{1/2+\alpha}}^2,
 \qquad
 \|w_N\|_{H^1}+\|g_N\|_{L^2(\partial D)}\leq C X_N.
\end{equation}
The low lift estimate, valid also for the time derivative, gives
\begin{equation}\label{functional:low_test}
 \|\mathcal Lg_N\|_2\leq C\|f^\circ\|_{\mathcal H_{\mathcal K}^{r}},
 \quad
 \|\mathcal LD_N^s\partial_tf^\circ\|_2
 \leq C\|\partial_tf^\circ\|_{\mathcal H_{\mathcal K}^{r}},
 \qquad r=1/2-2\alpha.
\end{equation}
For negative $r$ the norm here means the dual scale in
\eqref{functional:negative}.

\subsection{The complete surface enhancement test}

\begin{theorem}[Energy bridge]
\label{functional:enhancement}
Assume the preceding geometry, $0<\alpha<1/2$, and a weak identity
\begin{equation}\label{functional:abstract_weak}
 \int_\Omega J\partial_tV\cdot w+b_t(V,w)+a(f,w\cdot\mathcal N)
 +\int_{\partial D}c(w\cdot\mathcal N)
 =\ell_t(w)+\langle R,\nabla(w\cdot\mathcal N)\rangle_D
\end{equation}
for every divergence-free admissible test field, where
\[
 |b_t(V,w)|\leq C_b\|V\|_{H^1}\|w\|_{H^1},\qquad
 |\ell_t(w)|\leq L(t)\|w\|_{H^1},\quad
 c\in L^2(\partial D),\quad R\in H^{1/2-\alpha}(D)^2.
\]
For a priori estimates take the fields smooth in time and sufficiently
regular for this identity; the bound depends only on the displayed
norms.  Define
\[
 B(t)=\|V\|_{H^1}+L(t)+\|c\|_{L^2(\partial D)}
       +\|R\|_{H^{1/2-\alpha}(D)},
 \quad A(t)=D\Phi(t),\quad r=1/2-2\alpha.
\]
Then, for $0\leq t_0<t_1$,
\begin{align}
 \int_{t_0}^{t_1}\|f\|_{H^{3/2-\alpha}(D)}^2\,dt
 \leq C\Bigl[&\sum_{t=t_0,t_1}
 \|V(t)\|_2\|A(t)\|_\infty
             \|f^\circ(t)\|_{\mathcal H_{\mathcal K}^{r}}\nonumber\\
 &+\int_{t_0}^{t_1}\Bigl(
 B^2+|\langle f,1\rangle_D|^2
 +\|V\|_2\|\partial_tA\|_\infty\|f^\circ\|_{\mathcal H_{\mathcal K}^{r}}
 +\|V\|_2\|A\|_\infty\|\partial_tf^\circ\|_{\mathcal H_{\mathcal K}^{r}}
 \Bigr)\,dt\Bigr].\label{functional:integrated_enhancement}
\end{align}
The constant depends on the fixed geometry, $\alpha$, the coefficient
bounds, and $C_b$, but not on a truncation index.
Alternatively, if the inertial functional obeys
$|\int J\partial_tV\cdot w|\leq I(t)\|w\|_{H^1}$ on these tests,
then at that time
\begin{equation}\label{functional:pointwise_enhancement}
 \|f\|_{H^{3/2-\alpha}}^2
 \leq C\bigl(I(t)^2+B(t)^2+|\langle f,1\rangle_D|^2\bigr).
\end{equation}
\end{theorem}

\begin{proof}
Insert $w=w_N$ from \eqref{functional:mean_projection}.  The capillary
term is exactly $X_N^2$.  All terms other than inertia have absolute
value at most $CBX_N$, by \eqref{functional:curvature_pairing} and
\eqref{functional:XN}.  For inertia use $JM=A$ and the fixed linearity
of $\mathcal L$:
\begin{align}
 \int_\Omega J\partial_tV\cdot w_N
 &=\frac{d}{dt}\int_\Omega V\cdot A\mathcal Lg_N
 -\int_\Omega V\cdot(\partial_tA)\mathcal Lg_N
 -\int_\Omega V\cdot A\mathcal LD_N^s\partial_tf^\circ.
 \label{functional:inertia_identity}
\end{align}
There is no missing derivative of $J$ here: it is included in the exact
product identity $JM=D\Phi$.  The last two terms and the time-endpoint
terms have the bounds in \eqref{functional:integrated_enhancement} by
\eqref{functional:low_test}.  Use $CBX_N\leq X_N^2/2+C'B^2$, integrate,
and let $N\to\infty$.  The nonnegative spectral sums increase to
$\|f^\circ\|_{\mathcal H_{\mathcal K}^{3/2-\alpha}}^2$, so monotone
convergence applies also in time.  Proposition~\ref{functional:surface_scale}
identifies this with the required ordinary Sobolev norm.  Add the
constant mode, bounded by $C|\langle f,1\rangle_D|^2$.
For the pointwise alternative, bound inertia by $CIX_N$ before applying
Young's inequality and taking the same monotone limit.
\end{proof}

\begin{remark}[Application at the top time derivative]
Take $f=\partial_t^{m+1}\eta$ and $V=v_{m+1}$ in the actual differentiated
weak equation.  The theorem requires no time-index restriction $k\leq2$.
Its contact input includes $\kappa\partial_t^{m+2}\eta$ on the contact
curve, together with the differentiated contact forcing.  Its low
surface input is exactly
$\partial_t^{m+2}\eta\in\mathcal H_{\mathcal K}^{1/2-2\alpha}$ in the
integrable product appearing above.  This can be obtained from the
differentiated kinematic equation if its velocity-trace and commutator
terms are bounded in that space.  It is not obtained merely by naming
an energy functional.  Similarly, the bulk and traction residuals must
actually give $L(t)$, and the curvature residual must give
$\|R\|_{H^{1/2-\alpha}}$.  These are the exact nonlinear input estimates
needed to use the construction; the analytic lifting and limiting
argument are complete once those inputs hold.
\end{remark}

\begin{remark}[Localized and angularly differentiated fields]
For $f=\chi\partial_t^{k_1}\partial_\theta^{k_2}\eta$, apply
\eqref{functional:mean_projection} after taking these derivatives and
the cutoff.  Multiplication by $\chi$ does not preserve zero average.
No commutation of $D_N^s$ with $\chi$ or $\partial_\theta$ is asserted.
The operator is nonlocal: its test trace need not be supported where
$\chi$ is supported.  Consequently the localized weak equation must be
valid for these global admissible tests, with the cutoff commutators
included in its residuals.  A weak equation valid only for tests
supported in a local chart is not sufficient for this application.
The constant mode is paid by $\|f\|_{L^2}$, since $|\int_Df|\leq |D|^{1/2}\|f\|_2$. \end{remark}

\subsection{Pressure without a mixed scalar $H^2$ problem}

\begin{lemma}[Reference divergence inverse]
\label{functional:bogovskii}
There exists a bounded linear map
$\mathcal B:L^2_{0,\mathrm{av}}(\Omega)\to H^1_0(\Omega)^3$ with $\operatorname{div}\mathcal Bq=q$. \end{lemma}
\begin{proof}
For a bounded connected Lipschitz domain this is the Bogovskii
right-inverse theorem.  One precise derivation from Costabel--McIntosh,
\emph{On Bogovskii and regularized Poincare integral operators for de
Rham complexes on Lipschitz domains}, arXiv:0808.2614v2, uses Theorem
4.9(a),(c),(d) at $s=0$ in top degree.  The supported divergence range
is closed and has codimension $b_0=1$.  Supported $H^1$ fields are the
zero extensions of $H^1_0(\Omega)$ fields, whose divergences integrate
to zero.  The range therefore equals $L^2_{0,\mathrm{av}}(\Omega)$.
The restriction of divergence to the Hilbert orthogonal complement of
its kernel is a bounded bijection onto this range; the bounded inverse
theorem provides a bounded linear inverse.  This uses the closed-range
and cohomology statements, not only a parametrix with an unremoved
finite-dimensional remainder.
\end{proof}

Choose once and for all $\vartheta\in C_c^\infty(D)$ with
$\int_D\vartheta\,dx=1$.  In a collar over its support choose a smooth
reference field $z_*$, zero on the solid face and away from this patch,
with $z_*\cdot\mathcal N_0=\vartheta$ on the top.  For example, take a
smooth extension of $\vartheta\mathcal N_0/|\mathcal N_0|^2$ into this
interior top patch and multiply by a depth cutoff.  Then
\begin{equation}\label{functional:mean_test}
 h_*:=\operatorname{div}z_*,\qquad \int_\Omega h_*\,dx=1,
 \qquad w_*=Mz_*,\qquad w_*\cdot\mathcal N=\vartheta.
\end{equation}
This last trace is independent of the moving geometry and vanishes
near the contact edge.

\begin{theorem}[Fluctuation and mean of pressure]
\label{functional:pressure}
Suppose the weak pressure identity, for admissible tests, is
\begin{equation}\label{functional:pressure_weak}
 \int_\Omega Jp\operatorname{div}_{\mathcal A}w
 =\mathcal I(w)+b_t(V,w)+a(f,w\cdot\mathcal N)
 +\int_{\partial D}c(w\cdot\mathcal N)
 -\ell_t(w)-\langle R,\nabla(w\cdot\mathcal N)\rangle_D,
\end{equation}
where $p\in L^2(\Omega)$, $f\in H^1(D)$, $R\in L^2(D)^2$, and $b_t$
and $\ell_t$ have the bounds in Theorem~\ref{functional:enhancement}.
Assume
\[
 |\mathcal I(Mz)|\leq I_0\|z\|_{H^1}\quad(z\in H^1_0(\Omega)^3),
 \qquad I_*:=|\mathcal I(w_*)|<\infty.
\]
Then
\begin{equation}\label{functional:pressure_bound}
 \|p\|_{L^2(\Omega)}
 \leq C\bigl(I_0+I_*+\|V\|_{H^1}+L(t)
                   +\|f\|_{H^1(D)}+\|R\|_{L^2(D)}\bigr).
\end{equation}
No norm of $c$ is needed for these tests.  If $R$ is only a distribution,
its displayed norm can be replaced by $|\langle R,\nabla\vartheta\rangle_D|$. \end{theorem}

\begin{proof}
Set $\bar p=|\Omega|^{-1}\int_\Omega p$ and $q=p-\bar p$.
Take $z=\mathcal Bq$, $w=Mz$.  This test has zero full boundary trace,
so every surface and contact term vanishes, and $J\operatorname{div}_{\mathcal A}w=q$.  The left side is $\|q\|_2^2$. The right side is bounded by
$C(I_0+\|V\|_{H^1}+L)\|q\|_2$; division gives the fluctuation bound
(the case $q=0$ is immediate).

Next take $w=w_*$.  Its top trace is $\vartheta$, its wall trace is zero,
and its contact trace is zero.  The left side becomes
\[
 \int_\Omega p h_*\,dx=\bar p+\int_\Omega q h_*\,dx.
\]
All ordinary bulk terms are bounded by
$I_*+C(\|V\|_{H^1}+L)$, since $\|w_*\|_{H^1}$ is uniformly bounded by
the geometry assumptions.  The capillary term is bounded by
$C\|f\|_{H^1}\|\vartheta\|_{H^1}$, and the curvature remainder by
$\|R\|_2\|\nabla\vartheta\|_2$.  Bound the last integral of $q$ by
$\|q\|_2\|h_*\|_2$ and use the fluctuation bound.  This proves the
estimate for $\bar p$ and hence for $p$.
\end{proof}

In the application to $p_m$, the inertial functional is $\mathcal I(w)=\int_\Omega Jv_{m+1}\cdot w$; consequently both $I_0$ and $I_*$ are bounded by $C\|v_{m+1}\|_2$.  This pressure argument uses
only $\|\eta_m\|_{H^1}$ in its capillary term and has no need to create
a small high-surface interpolation term from a mixed scalar lift.
The differentiated weak equation and its forcing bounds are still
required.

\begin{lemma}[Tangential inertial adjoint]
\label{functional:angular_adjoint}
Let $X$ be a smooth vector field tangent to the boundary of the fixed
domain, or a periodic coordinate derivative on a chart with support
that eliminates the chart boundary terms.  For smooth fields $U,w$,
\begin{equation}\label{functional:weighted_adjoint}
 \int_\Omega J(XU)\cdot w\,dx
 =-\int_\Omega U\cdot\bigl(JXw+(XJ)w+J(\operatorname{div}X)w\bigr)\,dx.
\end{equation}
Thus its absolute value is at most $C\|U\|_2\|w\|_{H^1}$ when the
displayed coefficients are bounded.  In a chart with volume density
$\rho$, replace $J$ in the formula by $\rho J$.  If vector components
are rotated by a variable frame, differentiating that frame contributes
the corresponding bounded zeroth-order matrix term.
\end{lemma}
\begin{proof}
Integrate the divergence of $J(U\cdot w)X$.  Its boundary integral
vanishes by the stated tangency, support, or periodicity.  Expanding
the divergence gives exactly \eqref{functional:weighted_adjoint}.
The estimate follows by Cauchy--Schwarz.  Change of variables gives the
density version, and the product rule gives the frame term.
\end{proof}

This lemma supplies the pressure and pointwise surface estimates when
a localized inertial term is a tangential derivative of an already
controlled lower angular-order velocity.  Cutoff commutators, the
derivative of the volume density, and the frame connection must be
included in the actual localized equation before invoking it.

\section*{Acknowledgements}

Y.~Guo was supported by NSF Grant DMS-2405051, which provided partial support for X.~Yang. I.~Tice was supported by NSF Grant DMS-2508400.  The authors acknowledges the use of AI tools in the writing of this paper: they were used for generating the figures, proofreading the prose, and identifying mathematical typos and inconsistencies in initial drafts.

\end{document}